%% file: MHD-book.tex
\documentclass[a4paper]{amsbook}
\usepackage{hyperref}% Used to make clickable links between sections. Also makes a table of contents for the pdf viewer (cannot use [draft] option).
\usepackage{amsmath, amssymb}
\usepackage{graphicx}
\usepackage{dsfont}% Use this as replacement for \mathbb, e.g.\ \mathds{1}. Much much nicer than \mathbbm.
\theoremstyle{plain}

\newtheorem{thm}{Theorem}[chapter]%[section]
\newtheorem{prop}[thm]{Proposition}
\newtheorem{lem}[thm]{Lemma}

\theoremstyle{definition}

\newtheorem*{defi}{Definition} % [chapter]

\theoremstyle{remark}

\newtheorem*{rem}{Remark} % [chapter]

\numberwithin{section}{chapter}
\numberwithin{equation}{chapter}

\def\eqdefa{\buildrel\hbox{\footnotesize def}\over =}

\def\P {\mathbb{P}}
\def\L {\mathbf{L}}

\def\K {\mathbf{K}}

\def\CD {\mathcal{D}}

\def\CE {\mathcal{E}}

\def\cK {\mathcal{K}}
\def\cL {\mathcal{L}}

\def\cQ {\mathcal{Q}}

\def\a {{\alpha}}

\def\g {{\gamma}}

\def\eps {{\epsilon}}

\def\d {{\partial}}

\newcommand{\DIV}{\operatorname{div}}
\newcommand{\ROT}{\operatorname{rot}}
\newcommand{\Div}{\operatorname{div}}
\newcommand{\rot}{\operatorname{rot}}

\newcommand{\Span}{\operatorname{span}}

\newcommand{\Ker}{\operatorname{Ker}}

\newcommand{\ba}{\begin{aligned}}
\newcommand{\ea}{\end{aligned}}

\newcommand{\be}{\begin{equation}}
\newcommand{\ee}{\end{equation}}

\def\Ma{\hbox{{\rm Ma}}}
\def\Kn{\hbox{{\rm Kn}}}
\def\St{\hbox{{\rm St}}}
\def\Re{\hbox{{\rm Re}}}

\newcommand{\set}[2]{\left\{#1\,:\,#2\right\}}
\newcommand{\ip}[1]{\left\langle#1\right\rangle}

\begin{document}

\frontmatter

\title{From the Vlasov-Maxwell-Boltzmann system to incompressible viscous electro-magneto-hydrodynamics}

%    Information for first author
\author{Diogo Ars\'enio}
%    Address of record for the research reported here
\address{CNRS \& Universit\'e Paris Diderot, Institut de Math\'ematiques de Jussieu - Paris Rive Gauche, B\^atiment Sophie Germain, 8 place Aur\'elie Nemours, 75205 Paris Cedex 13, France}
%    Current address
\email{Diogo.Arsenio@imj-prg.fr}
%    \thanks will become a 1st page footnote.
% \thanks{}

%    Information for second author
\author{Laure Saint-Raymond}
\address{Universit\'e Pierre et Marie Curie \& Ecole Normale Sup\'erieure, D\'epartement de Math\'ematiques et Applications, 45 rue d'Ulm, 75230 Paris Cedex 05, France}
\email{Laure.Saint-Raymond@ens.fr}

\date{\today}
% \subjclass[2010]{Primary ;\\Secondary }

% \begin{abstract}
% 	bla bla bla
% \end{abstract}

\maketitle

\setcounter{tocdepth}{5}% Includes subsections, etc, in the table of contents.
\tableofcontents

\input{preface}

\mainmatter

\part{Formal derivations and macroscopic weak stability}

\input{VMB}

\input{formal}

\input{stability}

\part{Conditional convergence results}\label{part 2}

\input{convergence0}

\input{weak0}

\input{constraint0}

\input{strong0}

\input{highconstraint0}

\input{conservation0}

% \input{fluctuationlemma} % erase this chapter in the end, when everything has been incorporated

\input{t-oscillations0}

\input{Grad0}
\input{rel-entropy0}

\part{Unconditional convergence results~: long-range interactions}\label{part 3}
\part{Unconditional convergence results~: short-range interactions}\label{part 4}

% %%%%%%%%%%%%%%%%%%%%%%%%%%%%%%%%%%%%%%%%%%%%%%%%%%%%%%%%%%%%%%%%%%%%%%%%%%%
% \part{Unconditional convergence results~: long-range interactions}\label{part 3}
%
% \input{assumptions1}
%
% \input{existence1}
%
% \input{convergence1}
%
%
%
% %%%%%%%%%%%%%%%%%%%%%%%%%%%%%%%%%%%%%%%%%%%%%%%%%%%%%%%%%%%%%%%%%%%%%%%%%%%
% \part{Unconditional convergence results~: short-range interactions}\label{part 4}
%
% % \input{assumptions2}
%
% \input{existence2}
%
% \input{convergence2}

% \input{x-comp} % This file contains equations and bits of proofs for Young measures on the x compactness and hypoellipticity. Uses obsolete techniques on convexity of renormalizations and defect measures. Check again later and then erase. Probably safe to erase the whole file.

%%%%%%%%%%%%%%%%%%%%%%%%%%%%%%%%%%%%%%%%%%%%%%%%%%%%%%%%%%%%%%%%%%%%%%%%%%%

\appendix

\part*{Appendices}

\input{appendixtransfer}

\input{appendixyoung}

\input{appendixhypoellipticity}

\backmatter

% \newpage

% \part*{Bibliography}

\makeatletter% Use this commands so that the pdf bookmark (use hyperref package) corresponding to the bibliography appears at the level of a \part and not as a chapter in the appendices.
\renewcommand*{\toclevel@chapter}{-1} % Put chapter depth at the same level as \part.
% \chapter{Epilogue}
% \renewcommand*{\toclevel@chapter}{0} % Put chapter depth back to its default value.
\makeatother

\bibliographystyle{plain}
\bibliography{MHDbiblio.bib}

\end{document}

%% file: preface.tex
\chapter*{Preface}

The present book aims at presenting in a systematic, painstaking and rather exhaustive way the incompressible viscous fluid limits of the Vlasov-Maxwell-Boltzmann system for one or two species. In these regimes, the evolution of the fluid is governed by equations of Navier-Stokes-Fourier type, with some electromagnetic forcing. Depending on the precise scaling, this forcing term takes on various forms~: it may be linear or nonlinear, electrostatic or governed by some hyperbolic wave equations, possibly constrained by some relation of Ohm's type.

From the mathematical point of view, these models have very different behaviors~; in particular, the existence and stability of solutions require sometimes very weak notions of solutions. The asymptotic analysis, which consists most often in retrieving the structure of the limiting system in the scaled Vlasov-Maxwell-Boltzmann system, uses therefore various mathematical methods with important technical refinements.
Thus, in order to make the reading easier, different tools will be presented in separated chapters.

\bigskip

The first part of this work is devoted to the systematic formal analysis of viscous hydrodynamic limits. Chapter \ref{presentingVMB} introduces the Vlasov-Maxwell-Boltzmann system as well as its formal properties. An important point to be noted is that the a priori bounds coming from these physical laws are not enough to prove the {\bf existence of global solutions}, even in the renormalized sense of DiPerna and Lions \cite{diperna}, which is a major difficulty for the study of fast relaxation limits. This actually explains the dividing of the three other parts of this book, of increasing difficulty, giving rigorous convergence results in more and more general settings.

Chapter \ref{formal-chap} introduces the different scaling parameters arising in the system, and details the formal steps leading to the constraint relations and the evolution equations in each regime. We therefore obtain a rather precise {\bf classification of physically relevant models} for viscous incompressible plasmas, some of which actually do not seem to have been previously described in the literature.

Chapter \ref{weak stability} presents a mathematical analysis of these different models. The most singular of them have a behavior which is actually more similar to the incompressible Euler equations than to the Navier-Stokes equations~: the lack of weak stability does not allow to prove the existence of global solutions, with the exception of very weak solutions in the spirit of dissipative solutions introduced by Lions for the Euler equations \cite{lions7}. This {\bf lack of stability for limiting systems} is the second major difficulty for the study of hydrodynamic limits.

\bigskip

The goal of the second part is to make precise and rigorous the convergence results described formally in the first part. In order to isolate the difficulties which are specific to the asymptotic analysis, we choose here to prove first conditional results, i.e.\ to consider the convergence of renormalized solutions even though their existence is not known. This of course does not imply the convergence of weaker solutions which will be studied in the sequel (renormalized solutions with defect measure, and a fortiori solutions with Young measures), but most of the proof will remain unchanged.
The important point is that the analysis is based essentially on the uniform estimates coming from the scaled entropy inequality, which holds in all situations.

Furthermore, we will focus exclusively on two typical regimes, namely leading~:
\begin{itemize}
	\item from the one species Vlasov-Maxwell-Boltzmann equations to the incompressible quasi-static Navier-Stokes-Fourier-Maxwell-Poisson system~;
	
	\item from the two species Vlasov-Maxwell-Boltzmann equations to the two-fluid incompressible Navier-Stokes-Fourier-Maxwell system with Ohm's law in the case of strong interspecies collisions, or to the two-fluid incompressible Navier-Stokes-Fourier-Maxwell system with solenoidal Ohm's law in the case of weak interspecies collisions.
\end{itemize}
These asymptotic regimes are critical, in the sense that they are the most singular ones among the formal asymptotics mentioned in Chapter \ref{formal-chap} and that all remaining regimes can be treated rigorously by similar or even simpler arguments.

We will not detail in this preface the content of all chapters of the second part, but rather insist on the main points requiring a treatment different from usual hydrodynamic limits \cite{SR}.
In the case with only one species, the main difference is due to the fact that the transport equation contains force terms involving a derivative with respect to $v$, which does not allow to transfer equi-integrability from the $v$-variable to the $x$-variable as in \cite{golse3}. This is a major complication. The new idea here consists in getting first some {\bf strong compactness in $v$} by using regularizing properties of the gain operator \cite{lions} and, then, in transferring this strong compactness onto the spatial variable by means of {\bf refined hypoelliptic arguments} developed in \cite{arsenio}. The second important difference comes from the fast temporal oscillations which couple {\bf acoustic and electromagnetic modes}. Note that we introduce here a simple method to avoid dealing with non local projections.

Overall, we are eventually able to establish through weak compactness methods a very general result (Theorem \ref{NS-WEAKCV}) asserting the convergence of renormalized solutions of the one species Vlasov-Maxwell-Boltzmann system towards weak solutions of corresponding macroscopic systems.

In the case of two species, the situation not only requires to exploit the methods for one species, it is considerably more complex~:
\begin{itemize}
\item First of all, there is {\bf an additional scaling parameter} measuring the strength of interspecies interactions (and, incidentally, the typical size of the electric current, which can be much smaller than the bulk velocities of each of the two species of particles)~: this implies that the (formal) expansions involve a larger number of terms (for instance, the constraint equations are derived at different orders).

\item Secondly, the {\bf linearized collision operator} has a more complicated vectorial structure. The inversion of fluxes and the computation of dissipation terms in the limiting energy inequalities are therefore substantially more technical.

\item In the most singular regimes, we get nonlinear constraint equations. This means that {\bf renormalization methods}, {\bf compensated compactness techniques} and {\bf controls on the conservation defects} are already required at this stage of the proof.

\item We have no sufficient uniform a priori bound on the electric current to handle nonlinear terms, which prevents from taking limits in the approximate conservation of momentum law. To avoid this difficulty we need to introduce a {\bf modified conservation law involving the Poynting vector}.

\item Even in this more suitable form, the evolution equations are not stable under weak convergence, and we have no equi-integrability in these singular regimes. We develop therefore some improved modulated entropy method, which allows to consider renormalized solutions without important restriction on the initial data. Note that this {\bf renormalized modulated entropy method} should also lead to some improvements concerning the convergence of the Boltzmann equation (without any electromagnetic field) to the Navier-Stokes equations for ill-prepared initial data.
\end{itemize}

\bigskip

The third and fourth parts (which will be published in a second volume) are more technical. They show how to adapt the arguments presented in the conditional case of the second part to take into account the state of the art about the Cauchy theory for the Vlasov-Maxwell-Boltzmann system.

In the case of long-range microscopic interactions giving rise to a collision cross-section with a singularity for grazing collisions, treated in the third part, we start by proving the existence of renormalized solutions with a defect measure in the spirit of the construction by Alexandre and Villani \cite{alexandre}. This result, which is important independently of the study of hydrodynamic limits, has been addressed in the note \cite{arsenio-sr-CRAS}. The study of hydrodynamic limits follows then essentially the lines of \cite{arsenio3} (combined with the results of the conditional part). We would like however to mention some important contributions~:

\begin{itemize}
\item The first one concerns the {\bf  estimate of the defect measure}. A refined analysis of the convergence of approximate solutions to the Vlasov-Maxwell-Boltzmann system shows that the defect measure can be controlled by the entropy dissipation. This remark allows for a simplification of the proofs from \cite{arsenio3}, especially the passage to the limit in the kinetic equation leading to the characterization of the limiting form of the dissipation, and the control of conservation defects.

\item The other simplification is related to the {\bf renormalization process}. Here we choose a decomposition of the renormalized collision operator which allows both to control the singularity due to the collision cross-section, and to preserve the good scalings for the fluctuation. In particular, the same decomposition can be used for the control of the transport and of the conservation defects (with a loop estimate).
\end{itemize}

\bigskip

In the case of general microscopic interactions (including for instance the case of hard spheres), it is not known how to prove the convergence of approximation schemes of the Vlasov-Maxwell-Boltzmann system, due to a lack of compactness produced by the electromagnetic interaction. The existence of renormalized solutions is therefore still an open problem. Nevertheless, Lions \cite{lions3} has defined a very weak notion of solution --~the measure-valued renormalized solutions~-- defined as limit of approximate solutions~: the equation to be satisfied involves indeed Young measures.

In the fourth part, we begin by refining the control of Young measures for such solutions by the entropy inequality. We then proceed by showing that the estimates obtained in the second part are very stable, so that they can be generalized with Young measures. By using convexity properties and Jensen inequalities, we can extend all the arguments, and operate both the moment method and the entropy method in more singular regimes. This extension to solutions of the Vlasov-Maxwell-Boltzmann system defined in a very weak sense shows that the methods based on the entropy inequality are extremely robust, and that the convergence is essentially determined by the limiting system.

These good asymptotic properties seem to further indicate that the measure-valued solutions defined by Lions (which have never been really studied from the qualitative point of view) are relevant in some sense.

%% Please "sign" your preface
\vspace{1cm}
\begin{flushright}\noindent
Paris, France,\hfill {\it Diogo Ars\'enio \& Laure Saint-Raymond}\\
January 2016\hfill \null\\
\end{flushright}

%% file: VMB.tex
\chapter{The Vlasov-Maxwell-Boltzmann system}\label{presentingVMB}

In the present monograph, we intend to investigate in a rather systematic way the {\bf scaling limits of the Vlasov-Maxwell-Boltzmann system}
\begin{equation}\label{VMB}
	\begin{cases}
		\begin{aligned}
			\d_t f + v \cdot \nabla_x f + \frac{q}{m} \left( E + v \wedge B \right) \cdot \nabla_v f &= Q(f,f), \\
			\text{\footnotesize(Vlasov-Boltzmann)}\hspace{-50mm} &\\
			\mu_0\epsilon_0\d_t E - \ROT B &= - \mu_0q \int_{\mathbb{R}^3} fv dv, \\
			\text{\footnotesize(Amp\`ere)}\hspace{-50mm} &\\
			\d_t B + \ROT E& = 0, \\
			\text{\footnotesize(Faraday)}\hspace{-50mm} &\\
			\DIV E &=\frac{q}{\epsilon_0}\left(\int_{\mathbb{R}^3} fdv -1\right), \\
			\text{\footnotesize(Gauss)}\hspace{-50mm} &\\
			\DIV B &=0, \\
			\text{\footnotesize(Gauss)}\hspace{-50mm} &\\
		\end{aligned}
	\end{cases}
\end{equation}
leading to viscous incompressible magnetohydrodynamics, and to justify rigorously the corresponding asymptotics.

More precisely, the Vlasov-Maxwell-Boltzmann system describes the evolution of a gas of one species of charged particles (cations and anions (or electrons), i.e.\ positively and negatively charged ions, respectively) of mass $m>0$ and charge $q\in\mathbb{R}$, subject to auto-induced electromagnetic forces. Such a gas of charged particles, under a global neutrality condition, is called a plasma. The particle number density $f(t,x,v)\geq 0$, where $t\in [0,\infty)$, $x\in \Omega\subset\mathbb{R}^3$ and $v\in\mathbb{R}^3$, represents the distribution of particles which, at time $t$, are at position $x$ and have velocity $v$.

The evolution of the density $f$ is governed by the Vlasov-Boltzmann equation, which is the first line of \eqref{VMB}. In essence, it tells that the variation of the density $f$ along the trajectories of the particles (represented by the transport term $\partial_t f+v\cdot\nabla_x f$) is subject to the influence of a Lorentz force $q\left( E + v \wedge B \right)$ (represented by the Vlasov term $\frac{q}{m} \left( E + v \wedge B \right) \cdot \nabla_v f$) and inter-particle collisions in the gas (represented by the Boltzmann collision operator $Q(f,f)$).

The Lorentz force acting on the gas is auto-induced. That is, the electric field $E(t,x)$ and the magnetic field $B(t,x)$ are generated by the motion of the particles in the plasma itself. Their evolution is governed by Maxwell's equations, which are the remaining lines of \eqref{VMB}, namely Amp\`ere's equation, Faraday's equation and Gauss' laws. Here, the physical constants $\mu_0,\epsilon_0>0$ are, respectively, the vacuum permeability (or magnetic constant) and the vacuum permittivity (or electric constant). Recall that the speed of light is determined by the formula $c=\frac{1}{\sqrt{\mu_0\epsilon_0}}$.

\bigskip

We will also consider the \textbf{two species Vlasov-Maxwell-Boltzmann system}
\begin{equation}\label{2VMB}
	\begin{cases}
		\begin{aligned}
			\d_t f^+ + v \cdot \nabla_x f^+ + \frac{q^+}{m^+} \left( E + v \wedge B \right) \cdot \nabla_v f^+ &= Q(f^+,f^+) + Q(f^+,f^-), \\
			\text{\footnotesize(Vlasov-Boltzmann for cations)}\hspace{-50mm} &\\
			\d_t f^- + v \cdot \nabla_x f^- - \frac{q^-}{m^-} \left( E + v \wedge B \right) \cdot \nabla_v f^- &= Q(f^-,f^-) + Q(f^-,f^+), \\
			\text{\footnotesize(Vlasov-Boltzmann for anions)}\hspace{-50mm} &\\
			\mu_0\epsilon_0\d_t E - \ROT B &= - \mu_0 \int_{\mathbb{R}^3} \left(q^+f^+-q^-f^-\right)v dv, \\
			\text{\footnotesize(Amp\`ere)}\hspace{-50mm} &\\
			\d_t B + \ROT E& = 0,\\
			\text{\footnotesize(Faraday)}\hspace{-50mm} &\\
			\DIV E &=\frac{1}{\epsilon_0}\int_{\mathbb{R}^3} \left(q^+f^+-q^-f^-\right)dv, \\
			\text{\footnotesize(Gauss)}\hspace{-50mm} &\\
			\DIV B &=0, \\
			\text{\footnotesize(Gauss)}\hspace{-50mm} &\\
		\end{aligned}
	\end{cases}
\end{equation}
which is more physically accurate, since it describes the evolution of a gas of two species of oppositely charged particles (cations of charge $q^+>0$ and mass $m^+>0$, and anions of charge $-q^-<0$ and $m^->0$), subject to auto-induced electromagnetic forces.

Thus, the particle number density $f^+(t,x,v)\geq 0$ represents the distribution of the positively charged ions (i.e.\ cations), while the particle number density $f^-(t,x,v)\geq 0$ represents the distribution of the negatively charged ions (i.e.\ anions). Note that the collision operators $Q(f^+,f^-)$ and $Q(f^-,f^+)$ have been added to the right-hand sides of the respective Vlasov-Boltzmann equations in \eqref{2VMB} in order to account for the variations in the densities $f^+$ and $f^-$ due to interspecies collisions.

We refer to \cite{montgomery1964plasma} for a discussion of the validity of such systems from a physical viewpoint.

Henceforth, for the mere sake of mathematical convenience, we will make the simplification that both kinds of particles have the exact same mass $m^\pm=m>0$ and charge $q^\pm=q>0$. Even though this reduction may first appear rather unphysical, it remains nevertheless a reasonable approximation since the mass of cations and anions only differs by the mass of a few electrons, which is several orders of magnitude less than that of atomic nuclei. Anyway, we believe that the essential mathematical difficulties are contained in this case, and we expect that most of the analysis contained in this work carries over to the case of distinct masses, as long as they remain of a comparable order of magnitude. We refer to \cite{jang} for an independent formal study of some hydrodynamic limits of the two species Vlasov-Maxwell-Boltzmann system including the case of unequal masses, leading in particular to a formal justification of the Hall effect, which we will not address here.

% But first, for the mere sake of simplicity, let us focus on the one fluid system \eqref{VMB}. More on the two species system \eqref{2VMB} later on.

\bigskip

The mathematical framework we shall consider is the one defined by physical a priori estimates, namely {\bf entropy and energy bounds}, which corresponds to renormalized or even weaker solutions of the Vlasov-Maxwell-Boltzmann systems. For the sake of simplicity, we will consider, throughout this work, that the spatial domain is, in fact, the whole space $\Omega=\mathbb{R}^3$, thus avoiding the complicated discussion of boundary conditions.

\bigskip

The general strategy that will be used to study magnetohydrodynamic limits is therefore based on uniform a priori bounds, weak compactness and either the {\bf moment method of Grad}, or some generalized {\bf relative entropy method},  which are the only way to deal with very weak solutions~: strong convergence  requires indeed --~at least~-- that local conservation laws are satisfied.

\section{The Boltzmann collision operator}\label{defining collision operator}

The Boltzmann collision operator, present in the right-hand side of the Vlasov-Boltzmann equations in \eqref{VMB} and \eqref{2VMB}, is the quadratic form, acting on the velocity variable, associated to the bilinear operator
\begin{equation}
	\label{boltz-operator}
	Q(f,h)=\int _{\mathbb{R}^3} \int_{\mathbb{S}^2}
	\left(f'h'_* - fh_*\right) b(v-v_*,\sigma) d\sigma dv_*,
\end{equation}
where we have used the standard abbreviations
$$
	f=f(v), \qquad f'=f(v'), \qquad h_*=h(v_*), \qquad h_*'=h(v_*'),
$$
with $(v',v'_*)$ given by
$$
	v'={v+v_*\over 2} +{|v-v_*|\over 2} \sigma,
	\qquad
	v'_*={v+v_*\over 2} -{|v-v_*|\over 2} \sigma.
$$
One can easily show that the quadruple $(v,v_*,v',v_*')$ parametrized by $\sigma\in\mathbb{S}^{2}$ provides the family of all solutions to the system of four equations
\begin{equation}\label{conservation momentum energy}
	\begin{aligned}
		v+v_* & = v'+v_*',\\
		|v|^2+|v_*|^2 & = |v'|^2+|v_*'|^2.
	\end{aligned}
\end{equation}
At the kinetic level, these relations express the fact that interparticle collisions are assumed to be elastic and thus conserve momentum and energy, where $(v,v_*)$ denote the pre-collisional velocities and $(v',v_*')$ denote the post-collisional velocities of two interacting particles. Notice that the transformation $(v,v_*,\sigma) \mapsto \left(v',v'_*,\frac{v-v_*}{\left|v-v_*\right|}\right)$ is involutional.

It is to be emphasized that the definition of the Boltzmann operator for interspecies collisions with distinct masses is more complex. Indeed, in this case, the microscopic conservations of momentum and energy are
\begin{equation*}
	\begin{aligned}
		m_+v_++m_-v_- & = m_+v_+'+m_-v_-',\\
		m_+|v_+|^2+m_-|v_-|^2 & = m_+|v_+'|^2+m_-|v_-'|^2.
	\end{aligned}
\end{equation*}
Therefore, the masses must appear in the convolution relations defining the mixed collision operators, which become highly singular whenever the mass ratio tends to infinity or to zero. Again, for mathematical convenience, we will not deal with this case and stick to equal masses.

\bigskip

The Boltzmann collision operator can therefore be split, at least formally, into a gain term and a loss term
$$
	\begin{aligned}
		Q(f,h) & =Q^+(f,h)-Q^-(f,h) \\
		& = \int _{\mathbb{R}^3 \times \mathbb{S}^2}
		f'h'_* b dv_*d\sigma
		-\int _{\mathbb{R}^3 \times \mathbb{S}^2}
		fh_* b dv_* d\sigma.
	\end{aligned}
$$
The loss term counts all collisions in which a given particle of velocity $v$ will encounter another particle, of velocity $v_*$, and thus will change its velocity leading to a loss of particles of velocity $v$, whereas the gain term measures the number of particles of velocity $v$ which are created due to some collision between particles of velocities $v'$ and $v'_*$.

The cross-section $b=b(z,\sigma)$, or collision kernel, where $(z,\sigma)\in\mathbb{R}^3\times\mathbb{S}^{2}$, present in the integrand of \eqref{boltz-operator}, is a measurable function positive almost everywhere, which somehow measures the statistical repartition of post-collisional velocities $(v',v_*')$ given the pre-collisional velocities $(v,v_*)$. Its precise form depends crucially on the nature of the microscopic interactions, thus, it is determined by the intermolecular forces that are being considered. However, due to the Galilean invariance of collisions, it only depends on the magnitude of the relative velocity $|z|$ and on the deviation angle $\theta$, or deflection (scattering) angle, defined by $\cos\theta=k\cdot\sigma$ where $k=\frac{z}{|z|}$. We will therefore sometimes abuse notation and write $b(z,\sigma)=b(|z|,\cos\theta)$ without any confusion since the arguments of $b$ are then either vectors or scalars.

It is a common mathematical simplification, called the \textbf{cutoff} assumption, to suppose that the cross-section is at least locally integrable, i.e.\ $b(z,\sigma)\in L^1_\mathrm{loc}\left(\mathbb{R}^3\times\mathbb{S}^2\right)$. However, this hypothesis fails to hold when long-range interaction forces are present between the particles in the gas. Thus, in this \textbf{non-cutoff} case, the collision kernel is non-integrable. This is due to a strong singularity of the kernel in the angular variable created by the enormous amount of grazing collisions in the gas, i.e.\ collisions whose deflection angle is almost null.

\bigskip

For instance, if the particles are assumed to interact via a given repulsive potential $\Phi(r)$, where $r>0$ denotes the distance between two interacting particles, then the post-collisional velocities and especially the deviation angle $\theta$ can be computed in terms of the impact parameter $\beta$, i.e.\ the distance of closest approach if the particles were not to interact, and the relative velocity $z=v-v_*$ as the result of a classical scattering problem (see \cite{cercignani} for instance)~:
$$
	\theta(\beta,z)
	= \pi-2\int _0^{\frac\beta{s_0}} {du\over \sqrt{1-u^2-{4\over |z|^2} \Phi\left({\beta \over u}\right)}},
$$
where $s_0$ is the positive root of
$$
	1-{\beta^2\over s_0^2} -4 {\Phi(s_0)\over |z|^2}=0.
$$
Then the cross-section $b$ is implicitly defined by
$$
	b(|z|,\cos \theta)={\beta\over \sin \theta} {\partial \beta \over \partial \theta}|z|.
$$
It can be made fully explicit in the case of hard spheres
$$
	b(|z|,\cos \theta) =a^2 |z|,
$$
where $a>0$ is the (scaled) radius of the spheres.%, and in the case of Coulomb interactions where $b$ is given by Rutherford's formula.

As shown by Maxwell, it is possible to obtain a rather explicit expression for a wide class of physically relevant collision kernels (see \cite{villani} and references therein), namely the so-called inverse power kernels. This terminology stems from the fact that these kernels model a gas whose particles interact according to an {inverse power potential} $\Phi(r)=\frac{1}{r^{s-1}}$, where $r>0$ represents the distance between two particles and $s>2$. Maxwell's calculations show that in such a case one has
\begin{equation*}% \label{cross section 1}
	b\left(|z|,\cos\theta\right)=|z|^\gamma b_0(\cos\theta),\qquad \gamma=\frac{s-5}{s-1},
\end{equation*}
where the angular cross-section $b_0(\cos\theta)$ is smooth on $\theta\in (0,\pi)$ and has a non-integrable singularity at $\theta=0$ behaving as
\begin{equation*}% \label{cross section 2}
	b(\cos\theta)\sin\theta\sim\frac{1}{\theta^{1+\nu}},\qquad \nu=\frac{2}{s-1},
\end{equation*}
where the factor $\sin\theta$ accounts for the Jacobian determinant of spherical coordinates. Notice that, in this particular situation, $b(z,\sigma)$ is thus not locally integrable, which is not due to the specific form of inverse power potential. In fact, one can show (see \cite{villani}) that a non-integrable singularity arises if and only if forces of infinite range are present in the gas.% This, in turn, is a consequence of the enormous amount of grazing collisions, i.e.\ collisions with such large impact parameter $\beta$ that the deviation angle $\theta$ is almost null, existent as soon as long-range interactions between particles are considered.

The  case of Maxwellian molecules $s=5$ corresponds to $\gamma=0$, which is not physically relevant but enables one to perform many explicit calculations in agreement with physical observations. It is customary to loosely classify cross-sections into two categories~: hard and soft, respectively corresponding to the super-Maxwellian ($s>5$) and the sub-Maxwellian cases ($s<5$). We will however not employ this dichotomy since our hypotheses will allow us to treat all hard and soft kernels in a single unified theory.

It turns out that the limiting case $s=2$, which corresponds to Coulombian interactions, is not well suited for Boltzmann's equation as the Boltzmann collision operator should be replaced by the Landau operator in order to handle that situation (see \cite{villani}). The other limiting case $s=\infty$ corresponds formally to the hard spheres case.

\section{Formal macroscopic properties}\label{macro properties}

Using the well-known facts (see \cite{cercignani2}) that transforming $(v,v_*,\sigma)\mapsto(v_*,v,-\sigma)$ and $(v,v_*,\sigma)\mapsto \left(v',v_*',\frac{v-v_*}{\left|v-v_*\right|}\right)$ merely induces mappings with unit Jacobian determinants, known as the pre-post-collisional changes of variables or simply collisional symmetries, one can show that
\begin{equation}\label{sym integral}
	\begin{aligned}
		& \int_{\mathbb{R}^3} Q(f,f)(v) \varphi(v) dv \\
		& = \frac14 \int _{\mathbb{R}^3 \times \mathbb{R}^3 \times \mathbb{S}^2}
		\left(f'f'_* - ff_*\right) b(v-v_*,\sigma)
		\left(\varphi+\varphi_*-\varphi'-\varphi_*'\right)
		dvdv_* d\sigma,
	\end{aligned}
\end{equation}
for all $f(v)$ and $\varphi(v)$ regular enough. It then follows from \eqref{conservation momentum energy} that the above integral vanishes if and only if $\varphi(v)$ is a collision invariant, i.e.\ any linear combination of $\left\{1,v_1,v_2,v_3,|v|^2\right\}$.

Thus, successively multiplying the Vlasov-Boltzmann equation in \eqref{VMB} by the collision invariants and then integrating in velocity yields formally the local conservation laws
\begin{equation}\label{macroscopic conservation laws}
	\partial_t\int_{\mathbb{R}^3}f
	\begin{pmatrix} 1 \\ v \\ \frac{|v|^2}{2} \end{pmatrix}
	dv
	+
	\nabla_x\cdot\int_{\mathbb{R}^3}f
	\begin{pmatrix} v \\ v\otimes v \\ \frac{|v|^2}{2}v \end{pmatrix}
	dv
	=\frac{q}{m}
	\int_{\mathbb{R}^3}f
	\begin{pmatrix} 0 \\ E+v\wedge B \\ E\cdot v \end{pmatrix}
	dv,
\end{equation}
which provide the link to a macroscopic description of the gas.

In the case of two species \eqref{2VMB}, we obtain (recall that we are assuming equal masses $m^\pm=m$ and charges $q^\pm=q$)
\begin{equation}\label{macroscopic conservation mass two species}
	\partial_t\int_{\mathbb{R}^3}f^\pm
	dv
	+
	\nabla_x\cdot\int_{\mathbb{R}^3}f^\pm v
	dv
	=0
\end{equation}
and
\begin{equation}\label{macroscopic conservation laws two species}
	\begin{aligned}
		\partial_t\int_{\mathbb{R}^3}\left(f^++f^-\right)
		\begin{pmatrix} v \\ \frac{|v|^2}{2} \end{pmatrix}
		dv +
		\nabla_x & \cdot \int_{\mathbb{R}^3}\left(f^++f^-\right)
		\begin{pmatrix} v\otimes v \\ \frac{|v|^2}{2}v \end{pmatrix}
		dv \\
		& =\frac{q}{m}
		\int_{\mathbb{R}^3}\left(f^+-f^-\right)
		\begin{pmatrix} E+v\wedge B \\ E\cdot v \end{pmatrix}
		dv.
	\end{aligned}
\end{equation}

\bigskip

On the other hand, the standard energy estimates for Maxwell's system in \eqref{VMB} and \eqref{2VMB} (we refer to \cite{jackson} for more details on Maxwell's equations) are obtained, first, by taking the scalar product of the Amp\`ere and Faraday equations with $E$ and $B$, respectively, and summing the resulting quantities, which yields the conservation laws for one species
\begin{equation}\label{maxwell energy}
	\partial_t\left(\frac{\mu_0\epsilon_0|E|^2+|B|^2}{2}\right)+\nabla_x\cdot \left(E\wedge B\right)
	=-\mu_0 q \int_{\mathbb{R}^3} f E\cdot v dv,
\end{equation}
and for two species
\begin{equation}\label{maxwell energy two species}
	\partial_t\left(\frac{\mu_0\epsilon_0|E|^2+|B|^2}{2}\right)+\nabla_x\cdot \left(E\wedge B\right)
	=-\mu_0 q \int_{\mathbb{R}^3} \left( f^+-f^- \right) E\cdot v dv.
\end{equation}
Second, by taking the vector product of the Amp\`ere and Faraday equations with $B$ and $E$, respectively, employing Gauss' laws when necessary and summing the resulting quantities, which yields the conservation laws for one species
\begin{equation}\label{maxwell poynting}
	\begin{aligned}
		\mu_0\epsilon_0\partial_t \left(E\wedge B\right) + \nabla_x
		\left(\frac{\mu_0\epsilon_0|E|^2+|B|^2}{2}\right)
		& -\nabla_x\cdot\left(\mu_0\epsilon_0 E\otimes E+B\otimes B\right) \\
		& =-\mu_0 q \int_{\mathbb{R}^3}f \left(E+v\wedge B\right)dv + \mu_0 q E,
	\end{aligned}
\end{equation}
and for two species
\begin{equation}\label{maxwell poynting two species}
	\begin{aligned}
		\mu_0\epsilon_0\partial_t \left(E\wedge B\right) + \nabla_x
		\left(\frac{\mu_0\epsilon_0|E|^2+|B|^2}{2}\right)
		& -\nabla_x\cdot\left(\mu_0\epsilon_0 E\otimes E+B\otimes B\right) \\
		& =- \mu_0 q \int_{\mathbb{R}^3}\left(f^+-f^-\right) \left(E+v\wedge B\right)dv,
	\end{aligned}
\end{equation}
Notice the similitude of the source terms in \eqref{macroscopic conservation laws}, \eqref{maxwell energy}, \eqref{maxwell poynting}, and in \eqref{macroscopic conservation laws two species}, \eqref{maxwell energy two species}, \eqref{maxwell poynting two species}.

\bigskip

The other very important feature of the Boltzmann equation comes also from the symmetries of the collision operator. Without caring about integrability issues, we plug $\varphi=\log f$ into the symmetrized integral \eqref{sym integral} and use the properties of the logarithm to find
\begin{equation}\label{def D}
	\begin{aligned}
		D(f) & \eqdefa - \int_{\mathbb{R}^3} Q(f,f)\log f dv \\
		& = \frac14 \int _{\mathbb{R}^3 \times \mathbb{R}^3 \times \mathbb{S}^2}
		\left(f'f'_* - ff_*\right)  \log \left({f'f'_*\over ff_*}\right) b(v-v_*,\sigma)
		dvdv_* d\sigma \geq 0.
	\end{aligned}
\end{equation}
The so defined \textbf{entropy dissipation} $\int_{\mathbb{R}^3}D(f)(t,x)dx$ is non-negative and the functional $\int_0^t\int_{\mathbb{R}^3}D(f)(s,x)dxds$ is therefore nondecreasing on $t> 0$.

This leads to Boltzmann's $H$-theorem, also known as the second principle of thermodynamics, stating that the \textbf{entropy}
$$
	\int_{\mathbb{R}^3} f\log f dv
$$
is (at least formally) a Lyapunov functional for the Boltzmann equation. Indeed, formally multiplying the Vlasov-Boltzmann equation in \eqref{VMB} by $\log f$ and then integrating in space and velocity clearly leads to
\begin{equation}\label{Hthm}
	\frac{d}{dt}\int_{\mathbb{R}^3} f\log f(t,x,v) dv
	+ \nabla_x\cdot \int_{\mathbb{R}^3} f\log f(t,x,v) v dv
	+ D(f)(t,x) = 0.
\end{equation}
A similar procedure on the two species Vlasov-Boltzmann equations in \eqref{2VMB} yields
\begin{equation}\label{Hthm 2 species}
	\begin{aligned}
		\frac{d}{dt}\int_{\mathbb{R}^3} & \left(f^+\log f^+ + f^-\log f^-\right)(t,x,v) dv \\
		& + \nabla_x\cdot \int_{\mathbb{R}^3} \left(f^+\log f^+ + f^-\log f^-\right)(t,x,v) vdv \\
		& + \left(D\left(f^+\right)+D\left(f^-\right) + D\left(f^+,f^-\right)\right)(t,x)
		= 0,
	\end{aligned}
\end{equation}
where we have denoted the mixed entropy dissipation
\begin{equation}\label{def D mixed}
	\begin{aligned}
		D\left(f,h\right) & \eqdefa - \int_{\mathbb{R}^3} Q(f,h)\log f + Q(h,f)\log h dv \\
		& = \frac12 \int _{\mathbb{R}^3 \times \mathbb{R}^3 \times \mathbb{S}^2}
		\left(f'h'_* - fh_*\right)  \log \left({f'h'_*\over fh_*}\right) b(v-v_*,\sigma)
		dvdv_* d\sigma \geq 0.
	\end{aligned}
\end{equation}

\bigskip

As for the equation $Q(f,f)=0$, it is possible to show, since necessarily $D(f)=0$ in this case, that it is only satisfied by the so-called Maxwellian distributions $M_{R,U,T}$ defined by
\begin{equation*}
	M_{R,U,T}(v)=\frac{R}{(2\pi T)^{\frac{3}{2}}}e^{-\frac{\left|v-U\right|^2}{2 T}},
\end{equation*}
where $R\in\mathbb{R}_+$, $U\in\mathbb{R}^3$ and $T\in\mathbb{R}_+$ are respectively the macroscopic density, bulk velocity and temperature, under some appropriate choice of units. The relation $Q(f,f)=0$ expresses the fact that collisions are no longer responsible for any variation in the density and so, that the gas has reached statistical equilibrium. In fact, it is possible to show that if the density $f$ is a Maxwellian distribution for some $R(t,x)$, $U(t,x)$ and $T(t,x)$, then the macroscopic conservation laws \eqref{macroscopic conservation laws} turn out to constitute a compressible Euler system with electromagnetic forcing terms.

Similarly, for two species of particles, if the plasma reaches thermodynamic equilibrium so that the equations $Q\left(f^+,f^+\right)+Q\left(f^+,f^-\right)=0$ and $Q\left(f^-,f^-\right)+Q\left(f^-,f^+\right)=0$ are solved simultaneously, then necessarily $D\left(f^+\right)+D\left(f^-\right) + D\left(f^+,f^-\right)=0$, which implies that $f^+=M_{R^+,U^+,T^+}$ and $f^-=M_{R^-,U^-,T^-}$ with $U^+=U^-$ and $T^+=T^-$, but not necessarily equal masses. In this case, it is possible to show that the macroscopic system of conservation laws \eqref{macroscopic conservation mass two species}-\eqref{macroscopic conservation laws two species} constitute a compressible Euler system with electromagnetic forcing terms.

\bigskip

Finally, we define the (global) relative entropy, for any particle number density $f\geq 0$ and any Maxwellian distribution $M_{R,U,T}$, by
\begin{equation}\label{def H}
	H\left(f|M_{R,U,T}\right)(t)=
	\int_{\mathbb{R}^3\times\mathbb{R}^3}
	\left(f\log\frac{f}{M_{R,U,T}}-f+M_{R,U,T}\right)(t)dxdv \geq 0.
\end{equation}
We will more simply denote the relative entropy by $H(f)$, whenever the relative Maxwellian distribution is clearly implied. The global control of the relative entropies follows then from the non-negativity of the entropies dissipations.  Indeed, combining the $H$-theorem \eqref{Hthm} with the global conservation of mass and energy from \eqref{macroscopic conservation laws} and Maxwell's energy conservation \eqref{maxwell energy}, it is in general possible to establish for one species (see \cite{diperna3}, for instance), further integrating in time and space, by virtue of the convexity properties of the entropies and the entropy dissipations, the following weaker \textbf{relative entropy inequality}, for any $t>0$,
\begin{equation}\label{ent-ineq}
	\begin{aligned}
		\int_{\mathbb{R}^3\times\mathbb{R}^3} & \left(f \log {f \over M} - f +M\right) (t) dxdv \\
		& + \frac 1{2} \int_{\mathbb{R}^3} \left( \frac{\eps_0}{m} |E|^2+ {1\over m\mu_0} |B|^2 \right)(t) dx
		+ \int_0^t\int_{\mathbb{R}^3}D(f)(s) dx ds
		\\
		& \leq \int_{\mathbb{R}^3\times\mathbb{R}^3} \left(f^{\mathrm{in}}\log {f^{\mathrm{in}} \over M}
		- f^{\mathrm{in}}+M\right) dxdv \\
		& + \frac1{2}\int_{\mathbb{R}^3} \left( \frac{\eps_0}{m} |E^{\rm in}|^2+ {1\over m\mu_0} |B^{\rm in}|^2 \right) dx,
	\end{aligned}
\end{equation}
where $\left(f^{\rm in},E^{\rm in}, B^{\rm in}\right)$ denotes the initial data and $M$ denotes a global normalized Maxwellian distribution
\begin{equation*}
	M=M_{1,0,1}=\frac{1}{(2\pi)^{\frac{3}{2}}}e^{-\frac{\left|v\right|^2}{2}}.
\end{equation*}
Similarly, for two species, combining the $H$-theorem \eqref{Hthm 2 species} with the global conservation of mass and energy from \eqref{macroscopic conservation mass two species}-\eqref{macroscopic conservation laws two species} and Maxwell's energy conservation \eqref{maxwell energy two species}, we get the entropy inequality, for all $t>0$,
\begin{equation}\label{ent-ineq 2 species}
	\begin{aligned}
		\int_{\mathbb{R}^3\times\mathbb{R}^3} & \left(f^+ \log {f^+ \over M} - f^+ +M\right) (t)
		+ \left(f^- \log {f^- \over M} - f^- +M\right) (t) dxdv \\
		& + \frac 1{2} \int_{\mathbb{R}^3} \left( \frac{\eps_0}{m} |E|^2+ {1\over m\mu_0} |B|^2 \right)(t) dx \\
		& + \int_0^t\int_{\mathbb{R}^3}\left(D\left(f^+\right)+D\left(f^-\right) + D\left(f^+,f^-\right)\right)(s) dx ds
		\\
		& \leq \int_{\mathbb{R}^3\times\mathbb{R}^3} \left(f^{+\mathrm{in}}\log {f^{+\mathrm{in}} \over M}
		- f^{+\mathrm{in}}+M\right) dxdv \\
		& + \int_{\mathbb{R}^3\times\mathbb{R}^3} \left(f^{-\mathrm{in}}\log {f^{-\mathrm{in}} \over M}
		- f^{-\mathrm{in}}+M\right) dxdv \\
		& + \frac1{2}\int_{\mathbb{R}^3} \left( \frac{\eps_0}{m} |E^{\rm in}|^2+ {1\over m \mu_0} |B^{\rm in}|^2 \right) dx,
	\end{aligned}
\end{equation}
where $\left(f^{+ {\rm in}},f^{- {\rm in}},E^{\rm in}, B^{\rm in}\right)$ denotes the initial data.

Generally speaking, the $H$-theorem and the entropy inequalities \eqref{ent-ineq} and \eqref{ent-ineq 2 species} together with the conservation laws \eqref{macroscopic conservation laws} and \eqref{macroscopic conservation mass two species}-\eqref{macroscopic conservation laws two species} constitute key elements in the study of hydrodynamic limits.% In particular, these properties are the only requirement to derive formal limits and obtain a classification of magnetohydrodynamic regimes, which will be done in the next chapter.

\section{The mathematical framework}

The construction of suitable global solutions to the Vlasov-Maxwell-Boltzmann system \eqref{VMB}
$$
	\begin{cases}
		\begin{aligned}
			\d_t f + v \cdot \nabla_x f + \left( E + v \wedge B \right) \cdot \nabla_v f &= Q(f,f), \\
			\d_t E - \ROT B &= - \int_{\mathbb{R}^3} fv dv, \\
			\d_t B + \ROT E& = 0, \\
			\DIV E &=\int_{\mathbb{R}^3} fdv -1, \\
			\DIV B &=0,
		\end{aligned}
	\end{cases}
$$
or to the two species Vlasov-Maxwell-Boltzmann system \eqref{2VMB}
\begin{equation*}
	\begin{cases}
		\begin{aligned}
			\d_t f^\pm + v \cdot \nabla_x f^\pm \pm \left( E + v \wedge B \right) \cdot \nabla_v f^\pm &= Q(f^\pm,f^\pm) + Q(f^\pm,f^\mp), \\
			\d_t E - \ROT B &= -  \int_{\mathbb{R}^3} \left(f^+-f^-\right)v dv, \\
			\d_t B + \ROT E& = 0, \\
			\DIV E &=\int_{\mathbb{R}^3} \left(f^+-f^-\right)dv, \\
			\DIV B &=0,
		\end{aligned}
	\end{cases}
\end{equation*}
for large initial data is considered of outstanding difficulty, due to the lack of dissipative phenomena in Maxwell's equations, which are hyperbolic. Here, for the sake of simplicity, we have discarded all free parameters, since these are irrelevant for the existence theory. Thus, so far, the only known answer to this problem is due to Lions in \cite{lions3}, where a rather weak notion of solutions was derived~: the so-called measure-valued renormalized solutions. However, these solutions failed to reach mathematical consensus on their usefulness due to their very weak aspect.

It should be mentioned that an alternative approach yielding strong solutions, provided smallness and regularity assumptions on the initial data are satisfied, was obtained more recently by Guo in \cite{guo}. But such solutions fall out of the scope of our derivation of hydrodynamic limits since they are not based on the physical entropy and energy estimates. Anyway, were we to consider such strong solution, our approach and strategy would remain strictly the same, for, as we are about to see in Chapter \ref{formal-chap} below, the only uniform bounds valid in the hydrodynamic limit are precisely the physical entropy and energy estimates.

\bigskip

This poor understanding of the mathematical theory of the Vlasov-Maxwell-Boltzmann system is the reason why getting rigorous convergence results is so complex. For the sake of readability, we have therefore decided to separate the different kinds of difficulties.

\bigskip

\noindent$\bullet$ In a first part, we will prove {\bf conditional convergence results} restricting our attention to the case of Maxwellian cross-sections, i.e.\ $b \equiv 1$, for mere technical simplicity, and assuming the existence of \textbf{renormalized solutions} to \eqref{VMB} and \eqref{2VMB}, which is actually not known. It is to be emphasized that, even if this notion of solution is relatively rough, the convergence proof in this weak case has no purely technical difficulty specific to this roughness. Indeed, were we to deal with stronger solutions, the strategy of proof would not be any different or easier because we are considering here only the uniform bounds which come from physical estimates.

In this framework, we can focus on the key arguments of the convergence proof, which are not so different from the ones used for hydrodynamic limits of neutral gases. A crucial point is to understand how to get strong compactness on macroscopic fields, which cannot be dealt with using $L^1$ mixing lemma such as in \cite{golse3} because of the electromagnetic forcing terms. We will thus first prove strong compactness with respect to velocity, and then use refined hypoelliptic estimates established in \cite{arsenio} in order to transfer the strong compactness to the spatial variable (Chapter \ref{hypoellipticity}).

The other key point which requires a specific treatment is the study of fast time oscillations insofar as they possibly couple weak compressibility  with strong electromagnetic effects (Chapter \ref{oscillations}).

\bigskip

The second and third part will be then devoted to the understanding of additional technical difficulties related to the fact that we are not able to build renormalized solutions to the Vlasov-Maxwell-Boltzmann systems, but only even weaker solutions.

\bigskip

\noindent$\bullet$ In the case of {\bf singular collision kernels}, using the regularizing properties of the collision operator with respect to $v$, we will actually show the existence of \textbf{renormalized solutions with a defect measure} in the sense of Alexandre and Villani.
The major change is the fact that the renormalized kinetic equation is replaced by an inequality (the consistency coming from the conservation of mass). This leads to the introduction of a defect measure.

The important new step of the convergence proof is then to establish that this defect measure vanishes in the fast relaxation limit, which comes from refined entropy dissipation estimates. 

There are also many additional technical steps due to the singularity of the collision kernel, which makes the control of the conservation defects and the hypoelliptic transfer of compactness more difficult.

\bigskip

\noindent$\bullet$ In the apparently simpler case of {\bf cutoff collision kernels}, because of the lack of strong compactness estimates, we are not able to prove that approximate solutions $f_N$ to the Vlasov-Maxwell-Boltzmann systems \eqref{VMB} and \eqref{2VMB} converge to actual renormalized solutions. Indeed, without strong compactness properties, it is not possible to establish that $\beta(f_N) \to \beta(f)$ for any renormalization $\beta$, which accounts for the introduction of Young measures and the definition of a very rough notion of solution, namely the \textbf{measure-valued renormalized solutions}. Of course the physical meaning of such weak solutions is unclear, which probably explains why they have not been studied so far.

Nevertheless, we will establish here that --~in the fast relaxation limit~-- they exhibit the expected behavior, converging to the relevant magnetohydrodynamic model, which can be considered as an indication of their physical relevance.

The key point of the proof will be to obtain integrated versions of all estimates with respect to the Young measures, and to prove that asymptotically the Young measures are not seen by the limiting equation, even though they do not converge to Dirac masses due to lack of uniqueness of solutions in the limiting systems.

%% file: formal.tex
\chapter{Scalings and formal limits}\label{formal-chap}

In view of what is known on hydrodynamic limits of the Boltzmann equation (see \cite{SR3} and the references therein), which corresponds to the particular case where particles are not charged, i.e.\ $q=0$ in \eqref{VMB}, we will focus on {\bf incompressible diffusive regimes}, since we do not expect to be able to obtain a complete mathematical derivation for other choices of scalings.

% Other references include \cite{bardos4, bardos2, golse4, golse2, golse, golse5, levermore, lions4}

\section{Incompressible viscous regimes}

In the absence of electromagnetic field, the Boltzmann equation can be rewritten in non-dimensional variables
$$
	\St \d_t f +v\cdot \nabla_x f= {1\over \Kn} Q(f,f),
$$
where we have introduced the following parameters~:
\begin{itemize}
	
	\item the Knudsen number $\Kn=\frac{\lambda_0}{l_0}$, measuring the ratio of the mean free path $\lambda_0$ to the observation length scale $l_0$~;
	
	\item the Strouhal number $\St=\frac{l_0}{c_0t_0}$, measuring the ratio of the observation length scale $l_0$ to the typical length $c_0t_0$ run by a particle during a unit of time $t_0$, where $c_0$ is the speed of sound (or thermal speed)~;
	
	\item choosing the length $l_0$, time $t_0$ and velocity scales $u_0$ in such a way that we observe a macroscopic motion, i.e.\ $u_0=\frac{l_0}{t_0}$, we have the identity $\St=\Ma$ where the Mach number $\Ma=\frac{u_0}{c_0}$ is defined as the ratio of the bulk velocity $u_0$ to the thermal speed.
	
\end{itemize}

\bigskip

Hydrodynamic approximations are obtained in the fast relaxation  limit $\Kn \to 0$, which precisely corresponds to the asymptotic regime where the fluid under consideration satisfies the continuum hypothesis, for the mean free path becomes infinitesimally small. Because of the von K\'arm\'an relation for perfect gases, we then expect the flow to be dissipative when the Reynolds number
$$
	\Re \approx \frac{\Ma}{\Kn},
$$
measuring the inverse kinematic viscosity of the gas, is of order $1$, i.e.\ when the Mach number also tends to $0$.

In order to ensure the consistency of these scaling assumptions, we will consider --~as usual~-- data which are fluctuations $g$ of order $\Ma$
$$ f=M(1+\Ma\ g),$$
around a global normalized Maxwellian equilibrium
$$
	M(v)=\frac{1}{\left(2\pi\right)^\frac{3}{2}}e^{-\frac{|v|^2}{2}},
$$
of density $1$, bulk velocity $0$ and temperature $1$.

\bigskip

Thus, as is well-known since the works of Bardos, Golse and Levermore \cite{BGL2, bardos3}, the viscous incompressible hydrodynamic regimes of collisional kinetic systems are obtained in the fast relaxation limit when the above-mentioned dimensionless numbers $\Kn$, $\St$ and $\Ma$, are all of the same order $\epsilon>0$, say. In the sequel, we will therefore restrict our attention to the scaled Vlasov-Maxwell-Boltzmann system
\begin{equation*}
	\begin{cases}
		\begin{aligned}
			\epsilon\d_t f + v \cdot \nabla_x f + \frac{ql_0}{mc_0^2} \left( E + c_0v \wedge B \right) \cdot \nabla_v f &=
			\frac{1}{\epsilon}Q(f,f),
			\\
			f & =M\left(1+\epsilon g\right),
			\\
			\eps c_0 \mu_0\epsilon_0 \d_t E - \ROT B &= - \mu_0 q c_0 l_0 \int_{\mathbb{R}^3} fv dv,
			\\
			\eps c_0 \d_t B + \ROT E& = 0,
			\\
			\DIV E &=\frac{q l_0 }{\epsilon_0}\left(\int_{\mathbb{R}^3} fdv -1\right),
			\\
			\DIV B &=0,
			\\
		\end{aligned}
	\end{cases}
\end{equation*}
and to the scaled two species Vlasov-Maxwell-Boltzmann system
\begin{equation*}
	\begin{cases}
		\begin{aligned}
			\eps \d_t f^\pm + v \cdot \nabla_x f^\pm \pm \frac{ql_0}{mc_0^2} \left( E + c_0v \wedge B \right) & \cdot \nabla_v f^\pm
			\\
			&= \frac 1\eps Q(f^\pm,f^\pm)
			+ \frac{\delta^2}{\eps} Q(f^\pm,f^\mp),
			\\
			f^\pm & =M\left(1+\epsilon g^\pm\right),
			\\
			\eps c_0 \mu_0\epsilon_0 \d_t E - \ROT B &= - \mu_0 qc_0l_0 \int_{\mathbb{R}^3} \left(f^+-f^-\right)v dv,
			\\
			\eps c_0 \d_t B + \ROT E& = 0,
			\\
			\DIV E &=\frac{q l_0}{\epsilon_0}\int_{\mathbb{R}^3} \left(f^+-f^-\right)dv,
			\\
			\DIV B &=0,
			\\
		\end{aligned}
	\end{cases}
\end{equation*}
where we have introduced another bounded parameter $\delta>0$ in front of the interspecies collision operator to differentiate the strength of interactions. The size of the parameter $\delta$ will be compared to the Knudsen number $\Kn=\eps$ and we will distinguish three cases, due to their distinct asymptotic behavior~:
\begin{itemize}
	\item $\delta\sim 1$, strong interspecies interactions~;
	\item $\delta=o(1)$ and $\frac\delta\eps$ unbounded, weak interspecies interactions~;
	\item $\delta=O(\eps)$, very weak interspecies interactions.
\end{itemize}
Notice also that we have performed the same nondimensionalization on the whole Vlasov-Maxwell-Boltzmann systems, which explains the presence of the parameters $\eps$, $c_0$ and $l_0$ in Maxwell's equations.

\section{Scalings for the electromagnetic field}

First, from \eqref{ent-ineq}, we get the scaled entropy inequality for one species, for all $t>0$,
\begin{equation}\label{pre scaled entropy}
	\begin{aligned}
		\frac1{\eps^2} \int_{\mathbb{R}^3\times\mathbb{R}^3} & \left(f \log {f \over M} - f +M\right) (t) dxdv \\
		& + \frac 1{2c_0^2\eps^2} \int_{\mathbb{R}^3} \left( \frac{\eps_0}{m} |E|^2+ {1\over m\mu_0} |B|^2 \right)(t) dx
		+\frac{1}{\epsilon^4}\int_0^t\int_{\mathbb{R}^3}D(f)(s) dx ds
		\\
		& \leq  \frac1{\eps^2} \int_{\mathbb{R}^3\times\mathbb{R}^3} \left(f^{\mathrm{in}}\log {f^{\mathrm{in}} \over M}
		- f^{\mathrm{in}}+M\right) dxdv \\
		& + \frac1{2c_0^2\eps^2}\int_{\mathbb{R}^3} \left( \frac{\eps_0}{m} |E^{\rm in}|^2+ {1\over m\mu_0} |B^{\rm in}|^2 \right) dx,
	\end{aligned}
\end{equation}
where $\left(f^{\rm in},E^{\rm in}, B^{\rm in}\right)$ denotes the initial data.

As for the two species case, from \eqref{ent-ineq 2 species}, we get the scaled entropy inequality, for all $t>0$,
\begin{equation}\label{pre scaled entropy two species}
	\begin{aligned}
		\frac1{\eps^2} \int_{\mathbb{R}^3\times\mathbb{R}^3} & \left(f^+ \log {f^+ \over M} - f^+ +M\right) (t)
		+ \left(f^- \log {f^- \over M} - f^- +M\right) (t) dxdv \\
		& + \frac 1{2c_0^2\eps^2} \int_{\mathbb{R}^3} \left( \frac{\eps_0}{m} |E|^2+ {1\over m\mu_0} |B|^2 \right)(t) dx \\
		& +\frac{1}{\epsilon^4}\int_0^t\int_{\mathbb{R}^3}\left(D\left(f^+\right)+D\left(f^-\right) + \delta^2 D\left(f^+,f^-\right)\right)(s) dx ds
		\\
		& \leq  \frac1{\eps^2} \int_{\mathbb{R}^3\times\mathbb{R}^3} \left(f^{+\mathrm{in}}\log {f^{+\mathrm{in}} \over M}
		- f^{+\mathrm{in}}+M\right) dxdv \\
		& + \frac1{\eps^2} \int_{\mathbb{R}^3\times\mathbb{R}^3} \left(f^{-\mathrm{in}}\log {f^{-\mathrm{in}} \over M}
		- f^{-\mathrm{in}}+M\right) dxdv \\
		& + \frac1{2c_0^2\eps^2}\int_{\mathbb{R}^3} \left( \frac{\eps_0}{m} |E^{\rm in}|^2+ {1\over m \mu_0} |B^{\rm in}|^2 \right) dx,
	\end{aligned}
\end{equation}
where $\left(f^{+ {\rm in}},f^{- {\rm in}},E^{\rm in}, B^{\rm in}\right)$ denotes the initial data.

Note that the entropy inequalities \eqref{pre scaled entropy} and \eqref{pre scaled entropy two species} are the only uniform controls we have on the particle number densities and on the electric and magnetic fields, meaning that whatever the repartition of the free energy at the initial time, all the contributions are expected to be of the same order.

\bigskip

Thus, up to a change of units in $E$ and $B$, namely setting
$$
	\tilde E = \frac{1}{c_0\eps} \sqrt{\frac{\eps_0}{m}} E,\qquad \tilde B = \frac 1{c_0\eps\sqrt{m\mu_0}} B,
$$
so that $\tilde E$ and $\tilde B$ are uniformly controlled by the scaled entropy inequalities \eqref{pre scaled entropy} or \eqref{pre scaled entropy two species}, we have (dropping the tildes for the sake of readability), for one species,
\begin{equation}\label{pre scaled VMB}
	\begin{cases}
		\begin{aligned}
			\epsilon\d_t f + v \cdot \nabla_x f + \left( \alpha E + \beta v \wedge B \right) \cdot \nabla_v f &=
			\frac{1}{\epsilon}Q(f,f),
			\\
			f & =M\left(1+\epsilon g\right),
			\\
			\gamma \d_t E - \ROT B &= - \frac{\beta}{\epsilon^2} \int_{\mathbb{R}^3} fv dv,
			\\
			\gamma \d_t B + \ROT E& = 0,
			\\
			\DIV E &=\frac{\alpha}{\epsilon^2}\left(\int_{\mathbb{R}^3} fdv -1\right),
			\\
			\DIV B &=0,
		\end{aligned}
	\end{cases}
\end{equation}
and, for two species,
\begin{equation}\label{pre scaled VMB two species}
	\begin{cases}
		\begin{aligned}
			\eps \d_t f^\pm + v \cdot \nabla_x f^\pm \pm \left( \alpha E + \beta v \wedge B \right) \cdot \nabla_v f^\pm &= \frac 1\eps Q(f^\pm,f^\pm) + \frac{\delta^2}{\eps} Q(f^\pm,f^\mp),
			\\
			f^\pm & =M\left(1+\epsilon g^\pm\right),
			\\
			\gamma\d_t E - \ROT B &= - \frac{\beta}{\eps^2} \int_{\mathbb{R}^3} \left(f^+-f^-\right)v dv,
			\\
			\gamma\d_t B + \ROT E& = 0,
			\\
			\DIV E &=\frac{\alpha}{\epsilon^2}\int_{\mathbb{R}^3} \left(f^+-f^-\right)dv,
			\\
			\DIV B &=0,
			\\
		\end{aligned}
	\end{cases}
\end{equation}
where there are only three free parameters left (else that $\eps$ and $\delta$) to describe the qualitative behaviors of the systems, namely~:
\begin{itemize}
	\item $\alpha = \eps \frac{ql_0}{c_0\sqrt{m\eps_0}} $ measuring the electric repulsion according to Gauss' law~;
	\item $\beta= \eps ql_0 \sqrt{\frac{\mu_0}{m}} $ measuring the magnetic induction according to Amp\`ere's law~;
	\item $\gamma = \eps c_0 \sqrt{\eps_0 \mu_0} = u_0 \sqrt{\eps_0 \mu_0} $ which is nothing else than the ratio of the bulk velocity to the speed of light.
\end{itemize}
Notice that these parameters are naturally constrained to the relation
\begin{equation*}
	\beta=\frac{\alpha\gamma}{\eps}.
\end{equation*}

We will impose some natural restrictions on the size of $\alpha$, $\beta$ and $\gamma$. First of all, we will require that $\gamma=O(1)$. Note, however, that an unbounded $\frac{\gamma}{\eps}=c_0\sqrt{\eps_0\mu_0}$ may seem physically unrealistic since it corresponds to a regime where the thermal speed (i.e.\ the speed of sound) exceeds the speed of light. As usual, such situations should only be interpreted as asymptotic regimes where appropriate physical approximations are valid.

% Moreover, in the one species case, we will demand that $\alpha \gamma =O(\eps^2)$, so that dynamics are non-trivial in the limit $\eps \to 0$. Else we would indeed have from Amp\`ere's equation \eqref{ampere} that the bulk velocity $\frac{1}{\epsilon}\int_{\mathbb{R}^3} fv dv$ tends to zero. We further assume, for one species, that $\alpha$ and $\beta$ are of order $O(\epsilon)$, so that electric and magnetic forces create bounded acceleration terms in the Vlasov-Boltzmann equation in \eqref{pre scaled VMB}. Situations where one of these parameters is large compared to $\eps$ are much more complicated. Indeed, we expect the Lorentz force to strongly penalize the system, leading asymptotically to some nonlinear macroscopic constraint that we are not able to deal with. Actually, as far as we know, there is no mathematical method to investigate such problems of nonlinear singular perturbation. For instance, understanding the dynamo effect is a related question which remains challenging.

Moreover, in the one species case, we will demand that $\alpha$ and $\beta$ are of order $O(\epsilon)$, so that electric and magnetic forces create bounded acceleration terms in the Vlasov-Boltzmann equation in \eqref{pre scaled VMB}.

Situations where one of these parameters is large compared to $\eps$ are much more complicated. Indeed, we expect the Lorentz force to strongly penalize the system, leading asymptotically to some nonlinear macroscopic constraint that we are not able to deal with in the one species case. Actually, as far as we know, there is no systematic mathematical method to investigate such problems of nonlinear singular perturbation. For instance, understanding the dynamo effect is a related question which remains challenging.

Thus, on the whole, for one species, we will consider bounded parameters $\alpha$, $\beta$ and $\gamma$ satisfying
\begin{equation*}
	\alpha=O(\eps),\qquad\beta=O(\eps),\qquad\gamma=O(1)\qquad\text{and}\qquad \beta=\frac{\alpha\gamma}{\eps}.
\end{equation*}
We will then distinguish two critical cases, namely
\begin{enumerate}
	\item $\alpha=\eps$, $\beta=\eps$, $\gamma=\eps$,
	\item $\alpha=\eps^2$, $\beta=\eps$, $\gamma=1$,
\end{enumerate}
and will explain how all other cases can be easily deduced from the above, just eliminating lower order terms which are too small. The full range of parameters will be described later on by the Figure \ref{figure 1} on page \pageref{figure 1}. For the moment, we merely emphasize that the above-mentioned critical cases correspond exactly to the vertices of the domain represented in Figure \ref{figure 1}.

% For two species, the restrictions on the size of the parameters $\alpha$, $\beta$ and $\gamma$ are not so explicitly deduced by contemplation of the system \eqref{pre scaled VMB two species}. As a matter of fact, the need of asymptotically bounded acceleration terms in the macroscopic laws associated with the Vlasov-Boltzmann equations in \eqref{pre scaled VMB two species} leads us to require that $\alpha=O(\eps)=O(\delta)$, $\beta=O(\delta)$ and $\alpha\gamma=O(\eps\delta)$. Note that $\beta=O(\eps)$ is not required.

For two species, the restrictions on the size of the parameters $\alpha$, $\beta$ and $\gamma$ are not so explicitly deduced by inspection of the system \eqref{pre scaled VMB two species}, except in the case $\delta=O(\eps)$, which lowers the order of the interspecies collision term $\frac{\delta^2}{\eps} Q(f^\pm,f^\mp)$ in \eqref{pre scaled VMB two species} and is thus analog to the one species case. However, when $\frac\delta\eps$ is unbounded, the interspecies collision term $\frac{\delta^2}{\eps} Q(f^\pm,f^\mp)$ becomes a singular perturbation and, as a matter of fact, the need of asymptotically bounded acceleration terms in the macroscopic laws associated with the Vlasov-Boltzmann equations in \eqref{pre scaled VMB two species} leads us to require that $\alpha=O(\eps)$ and $\beta=O(\delta)$. Note that $\beta=O(\eps)$ is not required in this case, which is in sharp contrast with the one species case.

Thus, on the whole, for two species, we will consider bounded parameters $\alpha$, $\beta$, $\gamma$ and $\delta$ satisfying either
\begin{equation*}
	\alpha=O(\eps), \quad\beta=O(\eps),\quad\gamma=O(1) \quad\text{and}\quad \beta=\frac{\alpha\gamma}{\eps},
\end{equation*}
when $\delta=O(\eps)$, or
\begin{equation*}
	\alpha=O(\eps),\quad\beta=O(\delta),\quad\gamma=O(1)\quad\text{and}\quad \beta=\frac{\alpha\gamma}{\eps},
\end{equation*}
otherwise.

We will then distinguish two critical cases, namely
\begin{enumerate}
	\item $\alpha=\eps$, $\beta=\eps$, $\gamma=\eps$,
	\item $\alpha=\eps^2$, $\beta=\eps$, $\gamma=1$,
\end{enumerate}
when $\delta=O(\eps)$, and
\begin{enumerate}
	\item $\alpha=\eps$, $\beta=\delta$, $\gamma=\delta$,
	\item $\alpha=\delta\eps$, $\beta=\delta$, $\gamma=1$,
\end{enumerate}
when $\frac\delta\eps$ is unbounded (note that the latter two cases coincide when $\delta\sim 1$), and will explain how all other cases can be easily deduced from the above, just eliminating lower order terms which are too small. Thus, as for one species, when $\delta=O(\eps)$, the full range of parameters will be described by the Figure \ref{figure 1} on page \pageref{figure 1}. Furthermore, when $\delta\sim 1$, the range of parameters will be represented by the Figure \ref{figure 2} on page \pageref{figure 2}, while the case $\delta=o(1)$ with $\frac{\delta}{\eps}$ unbounded will be described by the Figure \ref{figure 3} on page \pageref{figure 3}. Again, we merely emphasize, for the moment, that the above-mentioned critical cases correspond exactly to the vertices of the respective domains represented in Figures \ref{figure 1}, \ref{figure 2} and \ref{figure 3}.

\section{Formal analysis of the one species asymptotics}\label{formal one}

Thus, for a plasma of one species of particles, our starting point is the scaled system
\begin{equation}\label{scaledVMB}
	\begin{cases}
		\begin{aligned}
			\epsilon\d_t f_\eps + v \cdot \nabla_x f_\eps + \left( \alpha E_\eps + \beta v \wedge B_\eps \right) \cdot \nabla_v f_\eps &=
			\frac{1}{\epsilon}Q(f_\eps,f_\eps),
			\\
			f_\eps & =M\left(1+\epsilon g_\eps\right),
			\\
			\gamma \d_t E_\eps - \ROT B_\eps &= - \frac{\beta}{\epsilon} \int_{\mathbb{R}^3} g_\eps v M dv,
			\\
			\gamma \d_t B_\eps + \ROT E_\eps & = 0,
			\\
			\DIV E_\eps &=\frac{\alpha}{\epsilon}\int_{\mathbb{R}^3} g_\eps M dv,
			\\
			\DIV B_\eps &=0,
		\end{aligned}
	\end{cases}
\end{equation}
supplemented with some initial data satisfying
\begin{equation*}
	\frac1{\eps^2}H\left(f_\eps^{\mathrm{in}}\right)
	+ \frac1{2}\int_{\mathbb{R}^3} |E_\eps^{\rm in}|^2+|B_\eps^{\rm in}|^2  dx < \infty,
\end{equation*}
where $H\left(f_\eps^{\mathrm{in}}\right)=H\left(f_\eps^{\mathrm{in}}|M\right)$. In particular, the corresponding scaled entropy inequality, where $t>0$,
\begin{equation}\label{scaled entropy}
	\begin{aligned}
		\frac1{\eps^2} H\left(f_\eps \right)
		+ \frac 1{2} \int_{\mathbb{R}^3} |E_\eps|^2+ |B_\eps|^2 dx
		& +\frac{1}{\epsilon^4}\int_0^t\int_{\mathbb{R}^3}D(f_\eps)(s) dx ds
		\\
		& \leq  \frac1{\eps^2}H\left(f_\eps^{\mathrm{in}}\right)
		+ \frac1{2}\int_{\mathbb{R}^3} |E_\eps^{\rm in}|^2+ |B_\eps^{\rm in}|^2 dx,
	\end{aligned}
\end{equation}
guarantees that the solution will remain --~for all non-negative times~-- a fluctuation of order $\eps$ around the global equilibrium $M$~:
$$f_\eps=M(1+\eps g_\eps).$$
Note that the kinetic equation in \eqref{scaledVMB} can then be rewritten, in terms of the fluctuation $g_\eps$, as
\begin{equation}\label{boltz-lin}
	\eps\d_t g_\eps +v\cdot \nabla_x g_\eps
	+(\alpha E_\eps+\beta v\wedge B_\eps) \cdot \nabla_v g_\eps-{\alpha\over \eps}  E_\eps \cdot v\left(1+\epsilon g_\epsilon\right)
	=-\frac1\eps \cL g_\eps +\cQ(g_\eps,g_\eps),
\end{equation}
where we denote
\begin{equation}\label{def L and Q}
	\cL g =-\frac1M \left(Q(Mg,M)+Q(M,Mg)\right) \quad \text{and} \quad \cQ(g,g) =\frac1M Q(Mg,Mg).
\end{equation}

\subsection{Thermodynamic equilibrium}\label{equilibrium conv}

The entropy inequality \eqref{scaled entropy} provides uniform bounds on $E_\eps$, $B_\eps$ and $g_\eps$. Therefore, assuming some formal compactness, up to extraction of subsequences, one has
$$
	E_\eps \rightharpoonup E, \qquad B_\eps \rightharpoonup B, \qquad g_\eps \rightharpoonup g,
$$
in a weak sense to be rigorously detailed in a subsequent chapter.

Then, multiplying \eqref{boltz-lin} by $\eps$, and taking formal limits as $\eps \to 0$ shows that $\cL g=0$. It can be shown (see Proposition \ref{coercivity}, below) that the kernel of the linearized Boltzmann operator $\mathcal{L}$ coincides exactly with the vector space spanned by the collision invariants $\left\{1,v_1,v_2,v_3,|v|^2\right\}$. Thus, we conclude that $g$ is an infinitesimal Maxwellian, that is a linear combination of collision invariants
\begin{equation}\label{infinitesimal maxwellian}
	g=\rho+u\cdot v + \theta \left(\frac{|v|^2}{2}-\frac 32\right),
\end{equation}
where $\rho\in\mathbb{R}$, $u\in\mathbb{R}^3$ and $\theta\in\mathbb{R}$ only depend on $t$ and $x$, and are respectively the fluctuations of density, bulk velocity and temperature.

The fact that the fluctuations assume the infinitesimal Maxwellian form describes that the gas reaches thermodynamic (or statistical) equilibrium, in the fast relaxation limit.

\bigskip

We define now the macroscopic fluctuations of density $\rho_\eps$, bulk velocity $u_\eps$ and temperature $\theta_\eps$ by
\begin{equation*}
	\begin{aligned}
		\rho_\eps & =\int_{\mathbb{R}^3}g_\eps Mdv,\\
		u_\eps & =\int_{\mathbb{R}^3}g_\eps v Mdv,\\
		\theta_\eps & =\int_{\mathbb{R}^3}g_\eps\left(\frac{|v|^2}{3}-1\right) Mdv,
	\end{aligned}
\end{equation*}
and the hydrodynamic projection $\Pi g_\eps$ of $g_\eps$ by
\begin{equation*}
	\Pi g_\eps = \rho_\eps+u_\eps\cdot v + \theta_\eps \left(\frac{|v|^2}{2}-\frac 32\right),
\end{equation*}
which is nothing but the orthogonal projection of $g_\eps$ onto the kernel of $\mathcal{L}$ in $L^2\left(Mdv\right)$.

Note that the previous step establishing the convergence of $g_\eps$ towards thermodynamic equilibrium yields, in fact, the uniform boundedness of $\frac 1\eps \mathcal{L}g_\eps$, which implies, at least formally, that
\begin{equation}\label{hydro projection}
	g_\eps-\Pi g_\eps = O(\eps).
\end{equation}
This convergence may also be derived directly from the uniform control of the entropy dissipation $\frac{1}{\eps^4}D(f_\eps)$ in the entropy inequality \eqref{scaled entropy}, provided we can control the large values of the fluctuations. Indeed, according to \eqref{def D}, we write
\begin{equation*}
	D(f_\eps) =
	\frac14 \int _{\mathbb{R}^3 \times \mathbb{R}^3 \times \mathbb{S}^2}
	\left(\frac{f_\eps'f'_{\eps *} - f_\eps f_{\eps *}}{f_\eps f_{\eps *}}\right)
	\log \left(1+\frac{f_\eps'f'_{\eps *} - f_\eps f_{\eps *}}{f_\eps f_{\eps *}}\right)
	f_\eps f_{\eps *}b
	dvdv_* d\sigma .
\end{equation*}
Therefore, since the non-negative function $z\log(1+z)$ behaves essentially as $z^2$, for small values of $|z|$, we deduce a formal control on
\begin{equation*}
	\frac1{4\eps^4} \int _{\mathbb{R}^3 \times \mathbb{R}^3 \times \mathbb{S}^2}
	\left(\frac{f_\eps'f'_{\eps *} - f_\eps f_{\eps *}}{f_\eps f_{\eps *}}\right)^2
	f_\eps f_{\eps *}b
	dvdv_* d\sigma .
\end{equation*}
Then, since
\begin{equation*}
	f_\eps'f'_{\eps *} - f_\eps f_{\eps *}= \eps \left(g_\eps' + g'_{\eps *} - g_\eps - g_{\eps *}\right) + \eps^2\left(g_\eps'g'_{\eps *} - g_\eps g_{\eps *}\right),
\end{equation*}
we infer that $\frac 1 {\eps}\left(g_\eps' + g'_{\eps *} - g_\eps - g_{\eps *}\right)$ is uniformly bounded, which, in other words, amounts to a control on $\frac 1\eps \mathcal{L}g_\eps$.

\bigskip

The asymptotic dynamics of $(\rho_\eps,u_\eps,\theta_\eps)$ is then governed by fluid equations, to be obtained from the moment equations associated with \eqref{boltz-lin}. Thus, successively multiplying \eqref{boltz-lin} by the collision invariants $1$, $v$ and $\frac{|v|^2}{2}$, and integrating in $Mdv$, yields
\begin{equation}\label{moment-eps0}
	\begin{cases}
		\begin{aligned}
			\d_t \rho_\eps +\frac1\eps \DIV u _\eps & = 0,\\
			\d_t u_\eps +\frac1\eps \nabla_x (\rho_\eps+\theta_\eps )
			- \frac{\alpha}{\eps^2} E_\eps
			& =
			\left(\frac\alpha\eps \rho_\eps E_\eps + \frac\beta\eps u_\eps \wedge B_\eps\right)
			-\frac1\eps \DIV \int_{\mathbb{R}^3} g_\eps \phi M dv ,\\
			\frac 32 \d_t (\rho_\eps+\theta_\eps) +\frac 5{2\eps} \DIV u_\eps
			& =
			\frac{\alpha}{\eps} u_\eps \cdot E_\eps
			- \frac1\eps \DIV \int_{\mathbb{R}^3}g_\eps \psi M dv,
		\end{aligned}
	\end{cases}
\end{equation}
where
\begin{equation}\label{phi-psi-def}
	\phi(v)=v\otimes v - \frac{|v|^2}{3}\operatorname{Id},
	\qquad
	\psi(v)=\left(\frac{|v|^2}{2}-\frac 52\right)v.
\end{equation}

Recall that we are assuming $\alpha=O(\eps)$ and $\beta=O(\eps)$. Hence, the nonlinear terms in the right-hand side of \eqref{moment-eps0} containing the electromagnetic fields are expected to be bounded. Furthermore, notice that the polynomials $\phi(v)$ and $\psi(v)$ are orthogonal to the collision invariants in the $L^2(Mdv)$ inner-product. That is to say $\int_{\mathbb{R}^3} \varphi \phi M dv=0$ and $\int_{\mathbb{R}^3} \varphi \psi M dv=0$, for all collision invariants $\varphi(v)$. Since, according to \eqref{hydro projection}, $g_\eps$ converges towards an infinitesimal Maxwellian with a rate $O(\eps)$, it is therefore natural to expect, at least formally, that the terms
\begin{equation*}
	\begin{aligned}
		\frac1\eps \int_{\mathbb{R}^3} g_\eps \phi M dv & = \frac1\eps \int_{\mathbb{R}^3} \left(g_\eps-\Pi g_\eps\right) \phi M dv, \\
		\frac1\eps \int_{\mathbb{R}^3} g_\eps \psi M dv & = \frac1\eps \int_{\mathbb{R}^3} \left(g_\eps-\Pi g_\eps\right) \psi M dv,
	\end{aligned}
\end{equation*}
in \eqref{moment-eps0} are bounded and have a limit.

More precisely, it can be shown that, in general, the linearized Boltzmann operator $\mathcal{L}$ is self-adjoint and Fredholm of index zero on $L^2(Mdv)$ (or a variant of it depending on the cross-section). Therefore, its range is exactly the orthogonal complement of its kernel. It follows that $\phi\in L^2(Mdv)$ and $\psi\in L^2(Mdv)$ belong to the range of $\mathcal{L}$ and, thus, that there are inverses $\tilde\phi\in L^2(Mdv)$ and $\tilde\psi\in L^2(Mdv)$ such that
\begin{equation}\label{phi-psi-def inverses}
	\phi =\cL \tilde \phi \qquad\text{and}\qquad \psi=\cL \tilde \psi,
\end{equation}
which can be uniquely determined by the fact that they are orthogonal to the kernel of $\mathcal{L}$ (i.e.\ to the collision invariants).

Consequently, the macroscopic system \eqref{moment-eps0} can be recast as
\begin{equation}\label{moment-eps}
	\begin{cases}
		\begin{aligned}
			\d_t \rho_\eps +\frac1\eps \DIV u _\eps & = 0,\\
			\d_t u_\eps +\frac1\eps \nabla_x (\rho_\eps+\theta_\eps )
			- \frac{\alpha}{\eps^2} E_\eps
			& =
			\left(\frac\alpha\eps\rho_\eps E_\eps +\frac\beta\eps u_\eps \wedge B_\eps\right)
			-\frac1\eps \DIV \int_{\mathbb{R}^3} \mathcal{L}g_\eps \tilde\phi M dv ,\\
			\frac 32 \d_t \theta_\eps +\frac 1{\eps} \DIV u_\eps
			& =
			\frac{\alpha}{\eps} u_\eps \cdot E_\eps
			- \frac1\eps \DIV \int_{\mathbb{R}^3}\mathcal{L}g_\eps \tilde\psi M dv,
		\end{aligned}
	\end{cases}
\end{equation}
where the terms $\frac 1\eps \mathcal{L}g_\eps$ will be expressed employing the Vlasov-Boltzmann equation \eqref{boltz-lin}. The above macroscopic system \eqref{moment-eps} is coupled with Maxwell's equations on $E_\eps$ and $B_\eps$~:
\begin{equation}\label{Maxwell-eps}
	\begin{cases}
		\begin{aligned}
			\gamma \d_t E_\eps - \ROT B_\eps &= - \frac{\beta}{\epsilon} u_\eps,
			\\
			\gamma \d_t B_\eps + \ROT E_\eps & = 0,
			\\
			\DIV E_\eps &=\frac{\alpha}{\epsilon}\rho_\eps,
			\\
			\DIV B_\eps &=0.
		\end{aligned}
	\end{cases}
\end{equation}

A careful formal analysis of the whole coupled macroscopic system \eqref{moment-eps}-\eqref{Maxwell-eps} will yield the asymptotic dynamics of $\left(\rho,u,\theta,E,B\right)$.

\subsection{Macroscopic constraints}\label{macro constraint one species}

At leading order, the system \eqref{moment-eps}-\eqref{Maxwell-eps} describes the propagation of acoustic $(\rho_\eps,u_\eps,\sqrt{\frac32}\theta_\eps)$ and electromagnetic $(E_\eps,B_\eps)$ waves~:
\begin{equation}\label{wave}
	\d_t
	\begin{pmatrix}
		\rho_\eps\\ u_\eps\\ \sqrt{3\over 2}\theta_\eps \\ E_\eps \\ B_\eps
	\end{pmatrix}
	+W_\eps
	\begin{pmatrix}
		\rho_\eps\\ u_\eps\\ \sqrt{3\over 2}\theta_\eps \\ E_\eps \\ B_\eps
	\end{pmatrix}
	=O(1),
\end{equation}
where the wave operator $W_\eps$, containing the singular terms from \eqref{moment-eps}-\eqref{Maxwell-eps} and defined explicitly below, is antisymmetric (with respect to the $L^2(dx)$ inner-product) and, therefore, can only have purely imaginary eigenvalues. The semi-group generated by this operator may thus produce fast time oscillations, which we are about to discuss briefly.
\begin{enumerate}

	\item When $\gamma \sim 1$ (so that $\alpha=O(\eps^2)$), we have
	\begin{equation}\label{wave operator 1}
		W_\eps=
		\begin{pmatrix}
			0 &\frac1\eps \DIV&0&0&0\\
			\frac1\eps \nabla_x&0&\frac 1\eps\sqrt{2\over 3} \nabla_x& 0 & 0\\
			0 & \frac 1\eps\sqrt{2\over 3} \DIV &0&0&0\\
			0&0&0&0&0\\
			0&0&0&0&0
		\end{pmatrix}.
	\end{equation}
	Thus, the singular perturbation creates only high frequency acoustic waves. Consequently, averaging over fast time oscillations as $\epsilon\to 0$, we get the macroscopic constraints
	\begin{equation}\label{constraint1}
		\DIV u =0, \qquad \rho+\theta =0,
	\end{equation}
	respectively referred to as incompressibility and Boussinesq relations. These are supplemented by the asymptotic constraints coming from Gauss' laws in \eqref{Maxwell-eps}
	\begin{equation*}
		\DIV E=0,\qquad \DIV B=0.
	\end{equation*}

	\item When $\gamma =o(1)$ and $\alpha=O(\eps^2)$, we have
	\begin{equation}\label{wave operator 2}
		W_\eps=
		\begin{pmatrix}
			0 &\frac1\eps \DIV&0&0&0\\
			\frac1\eps \nabla_x&0&\frac 1\eps\sqrt{2\over 3} \nabla_x& 0 & 0\\
			0 & \frac 1\eps\sqrt{2\over 3} \DIV &0&0&0\\
			0&0&0&0&-{1 \over \gamma} \ROT\\
			0&0&0&{1 \over \gamma} \ROT &0
		\end{pmatrix}.
	\end{equation}
	Thus, the singular perturbation creates both high frequency acoustic and electromagnetic waves. However, these waves remain decoupled and have a comparable frequency of oscillation if and only if $\gamma\sim\eps$. By averaging these fast time oscillations as $\epsilon\to 0$, we get the macroscopic constraints
	\begin{equation}\label{constraint2}
		\begin{aligned}
			\DIV u & =0, & \rho+\theta & = 0,\\
			\ROT B & = 0, & \ROT E & = 0.
		\end{aligned}
	\end{equation}
	These are supplemented by the asymptotic constraints coming from Gauss' laws in \eqref{Maxwell-eps}
	\begin{equation*}
		\DIV E=0,\qquad \DIV B=0.
	\end{equation*}
	Hence,
	\begin{equation*}
		E=0,\qquad B=0.
	\end{equation*}
	
	\item When $\gamma =o(1)$ and $\frac\alpha{\eps^2}$ is unbounded, we have
	\begin{equation}\label{wave operator 3}
		W_\eps=
		\begin{pmatrix}
			0 &\frac1\eps \DIV&0&0&0\\
			\frac1\eps \nabla_x&0&\frac 1\eps\sqrt{2\over 3} \nabla_x& -{\alpha \over \eps^2} \operatorname{Id} & 0\\
			0 & \frac 1\eps\sqrt{2\over 3} \DIV &0&0&0\\
			0&{\alpha \over \eps^2} \operatorname{Id}&0&0&-{1 \over \gamma} \ROT\\
			0&0&0&{1 \over \gamma} \ROT &0
		\end{pmatrix}.
	\end{equation}
	Thus, the singular perturbation creates both high frequency acoustic and electromagnetic waves, which are coupled. These waves may or may not have comparable frequency of oscillation. By averaging these fast time oscillations as $\epsilon\to 0$, we get the macroscopic constraints
	\begin{equation}\label{constraint3}
		\begin{aligned}
			\DIV u & =0, & \nabla_x (\rho+\theta) & = \left[\frac{\alpha}{\eps}\right] E,\\
			\ROT B & = \left[\frac{\beta}{\eps}\right] u, & \ROT E & = 0,
		\end{aligned}
	\end{equation}
	where we have denoted by $\left[\frac{\alpha}{\eps}\right]$ and $\left[\frac{\beta}{\eps}\right]$ the respective limits of $\frac \alpha\eps$ and $\frac{\beta}{\eps}$ as $\eps\to 0$. As usual, when $\alpha=o(\eps)$, the weak Boussinesq relation $\nabla_x (\rho+\theta)=0$ can be improved to the strong Boussinesq relation $\rho+\theta=0$, assuming $\rho$ and $\theta$ enjoy enough integrability. These are supplemented by the asymptotic constraints coming from Gauss' laws in \eqref{Maxwell-eps}
	\begin{equation*}
		\DIV E=\left[\frac{\alpha}{\eps}\right]\rho,\qquad \DIV B=0.
	\end{equation*}
	
\end{enumerate}

The exact nature of time oscillations produced by the system \eqref{wave}, in the limit $\eps\to 0$, will be rigorously discussed, with greater detail, later on in Chapter \ref{oscillations}.

\subsection{Evolution equations}\label{evolution}

The previous step shows that, since $W_\eps$ is singular, the asymptotic dynamics of $\left(\rho_\eps,u_\eps,\sqrt{\frac 32}\theta_\eps,E_\eps,B_\eps\right)$ becomes constrained to  the kernel $\Ker W_\eps$ as $\eps\to 0$. Moreover, since $W_\eps$ is antisymmetric, its range is necessarily orthogonal to its kernel. Therefore, in order to get the asymptotic evolution equations for $\left(\rho,u,\sqrt{\frac 32}\theta, E,B\right)$, it is natural to project the system \eqref{wave} onto $\Ker W_\eps$, which will rid us of all the singular terms in \eqref{wave} and allow us to pass to the limit. In other words, we will obtain the limiting dynamics of the system \eqref{wave} by testing it against functions in $\Ker W_\eps$.

We will denote by $P:L^2(dx)\to L^2(dx)$ the Leray projector onto solenoidal vector fields and $P^\perp=\operatorname{Id}-P$ the projector onto the orthogonal complement, that is $P= -\Delta^{-1}\rot\rot$ and $P^\perp = \Delta^{-1}\nabla\DIV$.
\begin{enumerate}

	\item When $\gamma \sim 1$, the kernel of $W_\eps$, defined in \eqref{wave operator 1}, is obviously determined by all $\left(\rho_\eps^0,u_\eps^0,\sqrt{\frac 32}\theta_\eps^0,E_\eps^0,B_\eps^0\right)$ which satisfy
	\begin{equation*}
			\DIV u_\eps^0 = 0 \qquad\text{and}\qquad \rho_\eps^0+\theta_\eps^0=0.
	\end{equation*}
	It is then readily seen that its orthogonal complement $\Ker W_\eps^\perp$ is determined by all $\left(\tilde \rho_\eps,\tilde u_\eps,\sqrt{\frac 32}\tilde \theta_\eps,\tilde E_\eps,\tilde B_\eps\right)$ such that
	\begin{equation*}
			P \tilde u_\eps = 0 \qquad\text{and}\qquad \frac 32\tilde \theta_\eps-\tilde \rho_\eps=0.
	\end{equation*}
	Hence, projecting the system \eqref{moment-eps}-\eqref{Maxwell-eps} onto $\Ker W_\eps^\perp$ yields
	\begin{equation}\label{moment-eps2}
		\begin{cases}
			\begin{aligned}
				\d_t Pu_\eps
				+
				\frac1\eps P \DIV \int_{\mathbb{R}^3} \mathcal{L}g_\eps \tilde\phi M dv
				& =
				P \left(\frac{\alpha}{\eps^2}E_\eps + \frac\alpha\eps\rho_\eps E_\eps +\frac\beta\eps u_\eps \wedge B_\eps\right) ,\\
				\d_t \left(\frac32\theta_\eps-\rho_\eps\right)
				+
				\frac1\eps \DIV \int_{\mathbb{R}^3}\mathcal{L}g_\eps \tilde\psi M dv
				& =
				\frac{\alpha}{\eps} u_\eps \cdot E_\eps.
			\end{aligned}
		\end{cases}
	\end{equation}

	\item When $\gamma =o(1)$ and $\alpha=O(\eps^2)$, the kernel of $W_\eps$, defined in \eqref{wave operator 2}, is obviously determined by all $\left(\rho_\eps^0,u_\eps^0,\sqrt{\frac 32}\theta_\eps^0,E_\eps^0,B_\eps^0\right)$ which satisfy
	\begin{equation*}
		\begin{aligned}
			\DIV u_\eps^0 & = 0 , & \rho_\eps^0+\theta_\eps^0 & =0,\\
			\rot E_\eps^0 & = 0, & \rot B_\eps^0 & = 0.
		\end{aligned}
	\end{equation*}
	It is then readily seen that its orthogonal complement $\Ker W_\eps^\perp$ is determined by all $\left(\tilde \rho_\eps,\tilde u_\eps,\sqrt{\frac 32}\tilde \theta_\eps,\tilde E_\eps,\tilde B_\eps\right)$ such that
	\begin{equation*}
		\begin{aligned}
			P\tilde u_\eps & = 0, & \frac 32\tilde \theta_\eps-\tilde \rho_\eps & =0, \\
			P^\perp \tilde E_\eps & =0, & P^\perp\tilde  B_\eps & = 0.
		\end{aligned}
	\end{equation*}
	Hence, projecting the system \eqref{moment-eps} onto $\Ker W_\eps^\perp$ also yields the system \eqref{moment-eps2}. Moreover, in view of Gauss' laws, the projection of Maxwell's equations \eqref{Maxwell-eps} onto $\Ker W_\eps^\perp$ yields no useful information.
	
	\item When $\gamma =o(1)$ and $\frac\alpha{\eps^2}$ is unbounded, the wave operator $W_\eps$ is defined by \eqref{wave operator 3}. Notice then that Gauss' laws from \eqref{Maxwell-eps} are invariant under the action of the wave operator $W_\eps$. Consequently, it is enough to consider the restriction of $W_\eps$ to electromagnetic fields which verify Gauss' laws. It follows that the kernel of $W_\eps$ is obviously determined by all $\left(\rho_\eps^0,u_\eps^0,\sqrt{\frac 32}\theta_\eps^0,E_\eps^0,B_\eps^0\right)$ which satisfy
	\begin{equation*}
		\begin{aligned}
			\ROT B_\eps^0 & =\frac{\beta}{\eps}u_\eps^0, &
			\nabla_x\left(\rho_\eps^0+\theta_\eps^0\right) & =\frac{\alpha}{\epsilon}E_\eps^0, \\
			\DIV E_\eps^0 & =\frac\alpha\eps\rho_\eps^0, & \DIV B_\eps^0 & = 0 .
		\end{aligned}
	\end{equation*}
	It is then readily seen that its orthogonal complement $\Ker W_\eps^\perp$ is determined by all $\left(\tilde \rho_\eps,\tilde u_\eps,\sqrt{\frac 32}\tilde \theta_\eps,\tilde E_\eps,\tilde B_\eps\right)$ such that
	\begin{equation*}
		\begin{aligned}
			\frac 32\tilde \theta_\eps-\tilde \rho_\eps & =0 , &
			\ROT \tilde u_\eps + \frac{\beta}{\eps}\tilde B_\eps & =0, \\
			\DIV \tilde E_\eps & =\frac\alpha\eps\tilde \rho_\eps, & \DIV \tilde B_\eps & = 0 .
		\end{aligned}
	\end{equation*}
	Considering the magnetic potential $\tilde B_\eps=\rot\tilde  A_\eps$, uniquely determined if $\Div \tilde A_\eps=0$ (i.e.\ fixing the Coulomb gauge), the above set of constraints can be rephrased as
	\begin{equation*}
		\begin{aligned}
			\frac 32\tilde \theta_\eps-\tilde \rho_\eps & =0 , &
			P\tilde u_\eps + \frac{\beta}{\eps} \tilde A_\eps & =0, \\
			\DIV \tilde E_\eps & =\frac\alpha\eps\tilde \rho_\eps, & \DIV\tilde A_\eps & = 0 .
		\end{aligned}
	\end{equation*}
	Hence, projecting the system \eqref{moment-eps} onto $\Ker W_\eps^\perp$ yields
	\begin{equation}\label{moment-eps3}
		\begin{cases}
			\begin{aligned}
				\d_t \left(Pu_\eps + \frac{\beta}{\eps}A_\eps \right)
				+
				\frac1\eps P \DIV \int_{\mathbb{R}^3} \mathcal{L}g_\eps \tilde\phi M dv
				& =
				P \left(\frac\alpha\eps\rho_\eps E_\eps +\frac\beta\eps u_\eps \wedge B_\eps\right) ,\\
				\d_t \left(\frac32\theta_\eps-\rho_\eps\right)
				+
				\frac1\eps \DIV \int_{\mathbb{R}^3}\mathcal{L}g_\eps \tilde\psi M dv
				& =
				\frac{\alpha}{\eps} u_\eps \cdot E_\eps,
			\end{aligned}
		\end{cases}
	\end{equation}
	where $B_\eps=\rot A_\eps$ and $\Div A_\eps=0$, and where we have used that Faraday's equation from \eqref{Maxwell-eps} implies
	\begin{equation*}
		\frac{\beta}{\eps}\partial_t A_\eps
+\frac{\alpha}{\eps^2}PE_\eps =0.
	\end{equation*}
	
\end{enumerate}

\bigskip

There only remains to evaluate the flux terms $\frac1\eps \int_{\mathbb{R}^3} \mathcal{L}g_\eps \tilde\phi M dv$ and $\frac1\eps \int_{\mathbb{R}^3}\mathcal{L}g_\eps \tilde\psi M dv$ in \eqref{moment-eps2} and \eqref{moment-eps3}. Following \cite{BGL2, bardos3}, this is done by employing \eqref{boltz-lin} to evaluate that
\begin{equation}\label{fluxes sub}
	\frac1\eps \cL g_\eps =
	\cQ(g_\eps,g_\eps) - v\cdot \nabla_x g_\eps
	+{\alpha\over \eps}  E_\eps \cdot v + O(\eps),
\end{equation}
which yields formally in the limit, by virtue of the infinitesimal Maxwellian form \eqref{infinitesimal maxwellian},
\begin{equation*}
	\begin{aligned}
		\lim_{\eps\to 0} \frac1\eps \cL g_\eps & =
		\cQ(g,g) - v\cdot \nabla_x g
		+\left[{\alpha\over \eps}\right]  E \cdot v \\
		& =\frac 12 \mathcal{L}\left(g^2\right) - v\cdot \nabla_x g
		+\left[{\alpha\over \eps}\right]  E \cdot v \\
		& =
		\frac 12 u^t\mathcal{L}(\phi)u
		+\theta u\cdot \mathcal{L}(\psi)
		+\frac 12 \theta^2\mathcal{L}\left(\frac{|v|^4}{4}\right)\\
		& -\Div\left((\rho+\theta)v+\frac{|v|^2}{3}u+\phi u+\theta\psi\right)
		+\left[{\alpha\over \eps}\right]  E \cdot v \\
		& =
		\frac 12 u^t\mathcal{L}(\phi)u
		+\theta u\cdot \mathcal{L}(\psi)
		+\frac 12 \theta^2\mathcal{L}\left(\frac{|v|^4}{4}\right)
		-\Div\left(\phi u+\theta\psi\right),
	\end{aligned}
\end{equation*}
where we have used, in the last line, that $\DIV u =0$ and $\nabla_x (\rho+\theta) = \left[\frac{\alpha}{\eps}\right] E$, whatever the asymptotic regime.

Next, we use that $\tilde\phi$ and $\tilde\psi$ have similar symmetry properties as $\phi$ and $\psi$, thanks to the rotational invariance of $\mathcal{L}$. More precisely, following \cite{desvillettes}, it can be shown (see also \cite[Section 2.2.3]{golse0}) that there exist two scalar valued functions $\alpha,\beta:[0,\infty)\to\mathbb{R}$ such that
\begin{equation*}
	\tilde\phi(v) = \alpha\left(\left|v\right|\right)\phi(v)
	\qquad\text{and}\qquad
	\tilde\psi(v) = \beta\left(\left|v\right|\right)\psi(v),
\end{equation*}
which implies (see \cite[Lemma 4.4]{BGL2}) that
\begin{equation}\label{delta identities}
	\begin{aligned}
		\int_{\mathbb{R}^3}\phi_{ij} \tilde\phi_{kl} Mdv &
		=\mu \left(\delta_{ik}\delta_{jl}+\delta_{il}\delta_{jk} - \frac 23 \delta_{ij}\delta_{kl} \right),
		\\
		\int_{\mathbb{R}^3}\psi_{i} \tilde\psi_{j} Mdv & = \frac 52 \kappa \delta_{ij},
	\end{aligned}
\end{equation}
where
\begin{equation}\label{mu kappa}
	\mu = \frac{1}{10}\int_{\mathbb{R}^3}\phi : \tilde\phi M dv
	\qquad\text{and}\qquad
	\kappa = \frac 2{15} \int_{\mathbb{R}^3}\psi\cdot\tilde\psi M dv.
\end{equation}

Hence, we conclude through tedious but straightforward calculations that
\begin{equation*}
	\begin{aligned}
		\lim_{\eps\to 0} \frac1\eps \int_{\mathbb{R}^3} \mathcal{L}g_\eps \tilde\phi M dv 
		& =
		\int_{\mathbb{R}^3} \frac 12 \left(u^t\phi u\right)\phi Mdv
		- \int_{\mathbb{R}^3} \Div\left(\phi u\right)\tilde\phi Mdv
		\\
		& =
		u\otimes u -\frac{|u|^2}{3}\operatorname{Id}
		- \mu \left(\nabla_x u+\nabla_x^t u\right),
		\\
		\lim_{\eps\to 0} \frac1\eps \int_{\mathbb{R}^3}\mathcal{L}g_\eps \tilde\psi M dv
		& =
		\int_{\mathbb{R}^3} \theta u\cdot \psi \psi M dv
		- \int_{\mathbb{R}^3} \Div\left(\theta\psi\right)\tilde\psi M dv\\
		& =
		\frac 52 \theta u
		- \frac 52 \kappa \nabla_x\theta.
	\end{aligned}
\end{equation*}
We finally identify the advection and diffusion terms
\begin{equation*}
	\begin{aligned}
		\lim_{\eps\to 0} \frac1\eps P\Div \int_{\mathbb{R}^3} \mathcal{L}g_\eps \tilde\phi M dv
		& =
		P\left(u\cdot\nabla_x u\right) - \mu\Delta_x u,
		\\
		\lim_{\eps\to 0} \frac1\eps \Div\int_{\mathbb{R}^3}\mathcal{L}g_\eps \tilde\psi M dv
		& =
		\frac 52 u\cdot\nabla_x\theta - \frac 52 \kappa \Delta_x\theta .
	\end{aligned}
\end{equation*}

\bigskip

On the whole, provided nonlinear terms remain stable in the limiting process, we obtain the following asymptotic systems~:
\begin{enumerate}

	\item When $\gamma \sim 1$, letting $\epsilon$ tend to zero in the system \eqref{moment-eps2} coupled with Maxwell's equations \eqref{Maxwell-eps} yields
	\begin{equation*}
		\begin{cases}
			\begin{aligned}
				\d_t u +
				u\cdot\nabla_x u - \mu\Delta_x u
				& = -\nabla_x p+
				\left[\frac{\alpha}{\eps^2}\right]E + \left[\frac\beta\eps\right] u \wedge B ,\\
				\d_t \theta
				+
				u\cdot\nabla_x\theta - \kappa \Delta_x\theta
				& = 0, \\
				\left[\gamma\right] \d_t E - \ROT B &= - \left[\frac{\beta}{\epsilon}\right] u,
				\\
				\left[\gamma\right] \d_t B + \ROT E & = 0,
			\end{aligned}
		\end{cases}
	\end{equation*}
	with the constraints from \eqref{constraint1}
	\begin{equation*}
		\begin{aligned}
			\Div u & =0, & \rho+\theta & = 0,\\
			\Div E & = 0 , & \Div B & = 0.
		\end{aligned}
	\end{equation*}

	\item When $\gamma =o(1)$ and $\alpha=O(\eps^2)$, letting $\epsilon$ tend to zero in the system \eqref{moment-eps2} coupled with Maxwell's equations \eqref{Maxwell-eps} yields
	\begin{equation*}
		\begin{cases}
			\begin{aligned}
				\d_t u +
				u\cdot\nabla_x u - \mu\Delta_x u
				& = -\nabla_x p ,\\
				\d_t \theta
				+
				u\cdot\nabla_x\theta - \kappa \Delta_x\theta
				& = 0,
			\end{aligned}
		\end{cases}
	\end{equation*}
	with the constraints from \eqref{constraint2}
	\begin{equation*}
		\begin{aligned}
			\Div u & =0, & \rho+\theta & = 0,\\
			E & = 0 , & B & = 0.
		\end{aligned}
	\end{equation*}
	
	\item When $\gamma =o(1)$ and $\frac\alpha{\eps^2}$ is unbounded, letting $\epsilon$ tend to zero in the system \eqref{moment-eps3} coupled with Maxwell's equations \eqref{Maxwell-eps} yields
	\begin{equation*}
		\begin{cases}
			\begin{aligned}
				\d_t \left(u + \left[\frac{\beta}{\eps}\right]A \right)
				+
				u\cdot\nabla_x u - \mu\Delta_x u
				& = -\nabla_x p +
				\left[\frac\alpha\eps\right]\rho E +\left[\frac\beta\eps\right] u \wedge B ,\\
				\d_t \left(\frac32\theta-\rho\right)
				+
				\frac 52 u\cdot\nabla_x\theta - \frac 52 \kappa \Delta_x\theta
				& =
				\left[\frac{\alpha}{\eps}\right] u \cdot E,
			\end{aligned}
		\end{cases}
	\end{equation*}
	with the constraints from \eqref{constraint3}
	\begin{equation*}
		\begin{aligned}
			\DIV u & =0, & \nabla_x (\rho+\theta) & = \left[\frac{\alpha}{\eps}\right] E,\\
			\ROT B & = \left[\frac{\beta}{\eps}\right] u, & \ROT E & = 0,\\
			\DIV E & =\left[\frac{\alpha}{\eps}\right]\rho, & \DIV B & =0, \\
			\rot A & = B, & \Div A & = 0.
		\end{aligned}
	\end{equation*}
	The above system can be rewritten more explicitly by defining the adjusted electric field $\tilde E = - \partial_t A + E$. It then holds that
	\begin{equation*}
		\begin{cases}
			\begin{aligned}
				\d_t u
				+
				u\cdot\nabla_x u - \mu\Delta_x u
				& = -\nabla_x p + \left[\frac{\beta}{\eps}\right]\tilde E
				+ \rho \nabla_x\theta +\left[\frac\beta\eps\right] u \wedge B ,\\
				\d_t \left(\frac32\theta-\rho\right)
				+
				u\cdot\nabla_x\left(\frac32\theta-\rho\right) - \frac 52 \kappa \Delta_x\theta
				& = 0,\\
				\partial_t B + \rot \tilde E  & = 0,
			\end{aligned}
		\end{cases}
	\end{equation*}
	with the constraints
	\begin{equation*}
		\begin{aligned}
			\DIV u & =0, & \Delta_x (\rho+\theta) & = \left[\frac{\alpha}{\eps}\right]^2\rho,\\
			\ROT B & = \left[\frac{\beta}{\eps}\right] u, & \Div B & = 0,\\
			&& \DIV \tilde E & =\left[\frac{\alpha}{\eps}\right]\rho.
		\end{aligned}
	\end{equation*}
	Notice, finally, that if further $\alpha=o(\eps)$, then the above system is greatly simplified and becomes
	\begin{equation*}
		\begin{cases}
			\begin{aligned}
				\d_t u
				+
				u\cdot\nabla_x u - \mu\Delta_x u
				& = -\nabla_x p + \left[\frac{\beta}{\eps}\right]\tilde E +
				\left[\frac\beta\eps\right] u \wedge B ,\\
				\d_t \theta
				+
				u\cdot\nabla_x\theta - \kappa \Delta_x\theta
				& =
				0, \\
				\partial_t B + \rot \tilde E  & = 0,
			\end{aligned}
		\end{cases}
	\end{equation*}
	with the constraints
	\begin{equation*}
		\begin{aligned}
			\DIV u & =0, & \rho+\theta & = 0,\\
			\Div \tilde E & = 0, &  \DIV B & =0, \\
			\ROT B & = \left[\frac{\beta}{\eps}\right] u, & E & = 0.
		\end{aligned}
	\end{equation*}	
	
\end{enumerate}

% \begin{enumerate}
% \item when $\gamma =1$ (so that $\alpha=O(\eps^2)$), 
% $$\begin{aligned}
% \d_t (u +c A)+P(\nabla\cdot (u\otimes u)-b u\wedge B) -\mu \Delta u =0\\
% \d_t (3\theta-2\rho)\theta+5\nabla \cdot (u\theta)  -5\kappa \Delta \theta =0\\
%  \d_t E-\ROT B=-c u\\
% \d_t B +\ROT E =0, \quad B=\ROT A.
% \end{aligned}
% $$
% which can be rewritten in the more standard way
% \begin{equation}
% \label{slow1}
% \begin{aligned}
% \d_t u +P(\nabla\cdot (u\otimes u)-cE - b u\wedge B) -\mu \Delta u =0\\
% \d_t \theta+\nabla \cdot (u\theta) -\frac25 u\cdot E -\kappa \Delta \theta =0\\
%  \d_t E-\ROT B=-c u,\quad 
% \d_t B +\ROT E =0.
% \end{aligned}
% \end{equation}
% where $b$ and  $c$ denote the respective  limits of $\beta/\eps$ and $\alpha\gamma/\eps^2$, 
% 
% 
% \item when $\gamma =o(1)$,
% $$\begin{aligned}
% \d_t (u +c A)+P(\nabla\cdot (u\otimes u)-a\rho E -b u\wedge B) -\mu \Delta u =0\\
% \d_t (3\theta-2\rho)+5\nabla \cdot (u\theta) -2 a u\cdot E -5\kappa \Delta \theta =0\\
% \end{aligned}
% $$
% which can be rewritten in the more standard way
% \begin{equation}
% \label{slow2}
% \begin{aligned}
% \d_t u +P(\nabla\cdot (u\otimes u)-c\tilde E - \rho \nabla(\rho+\theta) - b u\wedge B) -\mu \Delta u =0\\
% \d_t (3\theta-2\rho)+\nabla \cdot (u(3\theta-2\rho))  -5\kappa \Delta \theta =0\\
% \d_t B +\ROT \tilde E =0.
% \end{aligned}
% \end{equation}
% where $b$ and  $c$ denote the respective  limits of $\beta/\eps$ and $\alpha\gamma/\eps^2$, 
% \end{enumerate}

\subsection{Summary}

\begin{figure}[!ht]
	\centering
	\includegraphics[scale=0.54]{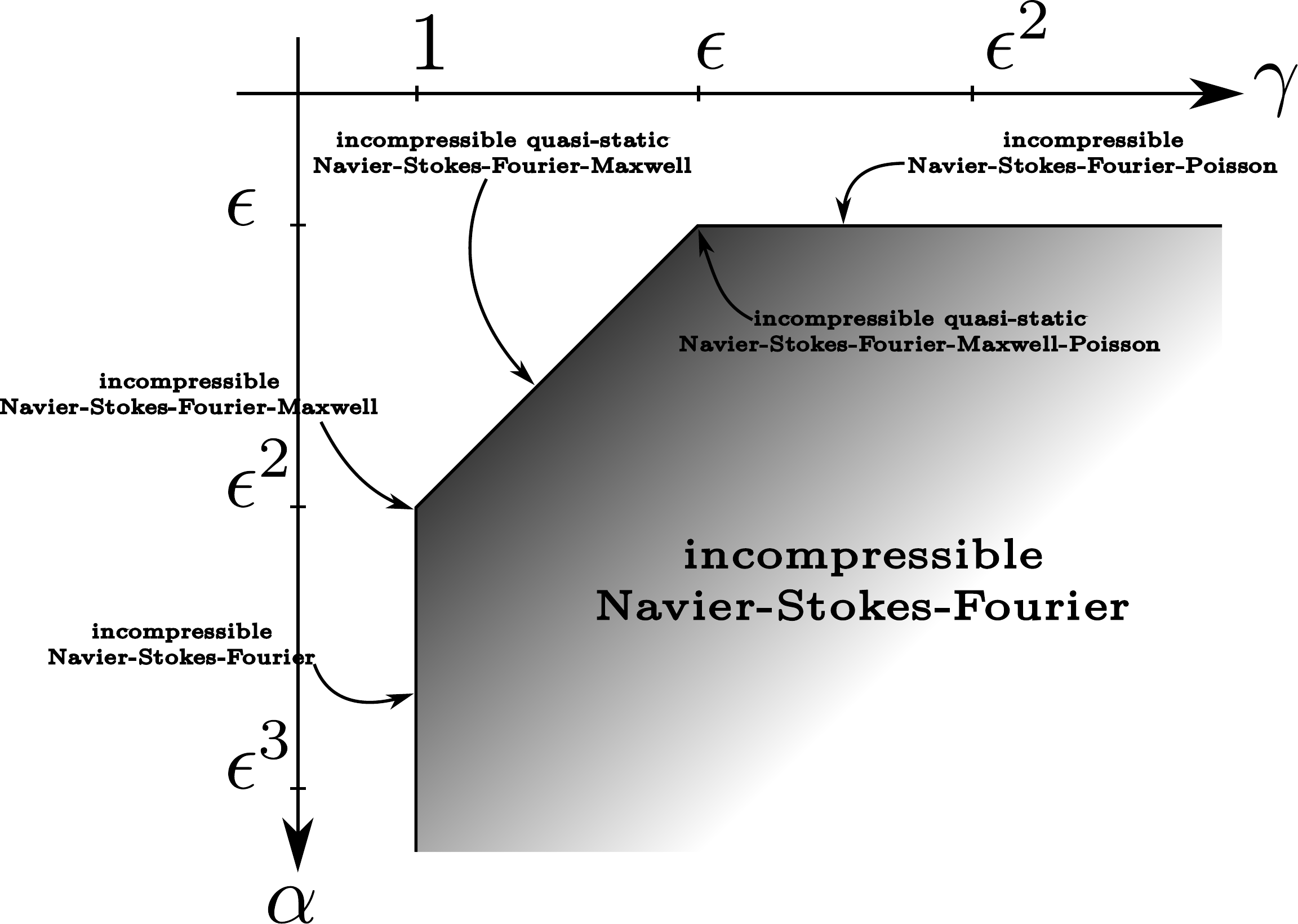}
	\caption{Asymptotic regimes of the one species Vlasov-Maxwell-Boltzmann system \eqref{scaledVMB}.}
	\label{figure 1}
\end{figure}

At last, we see that the asymptotics of the Vlasov-Maxwell-Boltzmann system \eqref{scaledVMB} can be depicted in terms of the limits of the following parameters~:
\begin{itemize}
	\item the strength of the electric induction $\alpha$,
	\item the strength of the magnetic induction $\beta=\frac{\alpha\gamma}{\eps}$,
	\item the ratio of the bulk velocity to the speed of light $\gamma$.
	% \item and the self-induction $\frac{\alpha\gamma}{\eps^2}$.
\end{itemize}
% Clearly, the nature of the systems that we have derived in Section \ref{evolution} is determined by the size of $\alpha$ and $\gamma$, rather than by the size of $\beta$. Indeed, the limit $\left[\frac\beta\eps\right]$ only appears as a coefficient in the Lorentz force in the momentum equation, while the limits $\left[\gamma\right]$, $\left[\frac{\alpha\gamma}{\eps^2}\right]$, $\left[\frac{\alpha}{\eps^2}\right]$ and $\left[\frac{\alpha}{\eps}\right]$ truly determine asymptotically the coupling between the Navier-Stokes-Fourier equations and Maxwell's equations. Furthermore, recall that the physical range is constrained to the relation $\beta=\frac{\alpha\gamma}{\eps}$, which support the viewpoint that the limiting regimes are described by only the two parameters $\alpha$ and $\gamma$.
Figure \ref{figure 1} summarizes the different asymptotic regimes, on a logarithmic scale, of the Vlasov-Maxwell-Boltzmann system \eqref{scaledVMB}.% in the physical range $\beta=\frac{\alpha\gamma}{\eps}$, so that $\left[\frac\beta\eps\right]=\left[\frac{\alpha\gamma}{\eps^2}\right]$. Again, one should keep in mind that truly physical situations further require that $\gamma=O(\eps)$, which corresponds to a non-vanishing light velocity.

% \begin{figure*}[h] %  figure placement: here, top, bottom, or page
% 	\centering
% 	\includegraphics[width=4.5in]{MHD.pdf} 
% 	\caption{Asymptotic regimes of the Vlasov-Maxwell-Boltzmann system \eqref{scaledVMB}.}
% 	% \label{fig1}
% \end{figure*}

Thus, up to multiplicative constants, we reach the following asymptotic systems of equations~:
\begin{enumerate}
	
	\item If $\alpha=o(\eps)$ and $\beta=o(\eps)$, we obtain the incompressible Navier-Stokes-Fourier system~:
	\begin{equation*}
		\begin{cases}
			\begin{aligned}
				\d_t u +
				u\cdot\nabla_x u - \mu\Delta_x u
				& = -\nabla_x p , & \Div u & =0, \\
				\d_t \theta
				+
				u\cdot\nabla_x\theta - \kappa \Delta_x\theta
				& = 0, & \rho+\theta & =0.
			\end{aligned}
		\end{cases}
	\end{equation*}
	This system satisfies the following formal energy conservation laws~:
	\begin{equation*}
		\begin{aligned}
			\frac 12\frac{d}{dt}\left\|u\right\|_{L^2_x}^2 + \mu \left\|\nabla_x u\right\|_{L^2_x}^2 & = 0, \\
			\frac 12\frac{d}{dt}\left\|\theta\right\|_{L^2_x}^2 + \kappa \left\|\nabla_x \theta\right\|_{L^2_x}^2 & = 0.
		\end{aligned}
	\end{equation*}
	
	\item If $\alpha=\eps^2$ and $\gamma=1$, we obtain the incompressible Navier-Stokes-Fourier-Maxwell system~:
	\begin{equation*}
		\begin{cases}
			\begin{aligned}
				\d_t u +
				u\cdot\nabla_x u - \mu\Delta_x u
				& = -\nabla_x p+
				E + u \wedge B , & \Div u & = 0,\\
				\d_t \theta
				+
				u\cdot\nabla_x\theta - \kappa \Delta_x\theta
				& = 0, & \rho+\theta & = 0, \\
				\d_t E - \ROT B &= -  u, & \Div E & = 0,
				\\
				\d_t B + \ROT E & = 0, & \Div B & = 0.
			\end{aligned}
		\end{cases}
	\end{equation*}
	This system satisfies the following formal energy conservation laws~:
	\begin{equation*}
		\begin{aligned}
			\frac 12\frac{d}{dt}\left(\left\|u\right\|_{L^2_x}^2 + \left\|E\right\|_{L^2_x}^2 + \left\|B\right\|_{L^2_x}^2 \right) + \mu \left\|\nabla_x u\right\|_{L^2_x}^2 & = 0, \\
			\frac 12\frac{d}{dt}\left\|\theta\right\|_{L^2_x}^2 + \kappa \left\|\nabla_x \theta\right\|_{L^2_x}^2 & = 0.
		\end{aligned}
	\end{equation*}
	
	\item If $\alpha=o(\eps)$, $\beta=\eps$ and $\gamma=o(1)$, we obtain the incompressible quasi-static Navier-Stokes-Fourier-Maxwell system~:
	\begin{equation}\label{IQSNSFM}
		\begin{cases}
			\begin{aligned}
				\d_t u
				+
				u\cdot\nabla_x u - \mu\Delta_x u
				& = -\nabla_x p + E + u \wedge B , & \Div u & = 0,\\
				\d_t \theta
				+
				u\cdot\nabla_x\theta - \kappa \Delta_x\theta
				& = 0 , & \rho+\theta & =0 , \\
				\rot B & = u, & \Div E & = 0 , \\
				\partial_t B + \rot E  & = 0, & \Div B & = 0.
			\end{aligned}
		\end{cases}
	\end{equation}
	This system satisfies the following formal energy conservation laws~:
	\begin{equation*}
		\begin{aligned}
			\frac 12\frac{d}{dt}\left(\left\|u\right\|_{L^2_x}^2 + \left\|B\right\|_{L^2_x}^2 \right) + \mu \left\|\nabla_x u\right\|_{L^2_x}^2 & = 0, \\
			\frac 12\frac{d}{dt}\left\|\theta\right\|_{L^2_x}^2 + \kappa \left\|\nabla_x \theta\right\|_{L^2_x}^2 & = 0.
		\end{aligned}
	\end{equation*}
	Here, the electric field is defined indirectly as a mere distribution, through Faraday's equation, by
	\begin{equation*}
		E=-\partial_t A,
	\end{equation*}
	where $B=\rot A$ and $\Div A=0$.
	
	\item If $\alpha=\eps$ and $\gamma=\eps$, we obtain the incompressible quasi-static Navier-Stokes-Fourier-Maxwell-Poisson system~:
	\begin{equation}\label{NSFMP}
		\begin{cases}
			\begin{aligned}
				\d_t u
				+
				u\cdot\nabla_x u - \mu\Delta_x u
				& = -\nabla_x p + E + \rho \nabla_x\theta + u \wedge B , \hspace{-20mm}&& \\
				&& \Div u & = 0,\\
				\d_t \left(\frac32\theta-\rho\right)
				+
				u\cdot\nabla_x\left(\frac32\theta-\rho\right)
				- \frac 52 \kappa \Delta_x\theta
				& = 0,
				& \Delta_x(\rho+\theta) & =\rho, \\
				\ROT B & = u, & \Div E & = \rho , \\
				\partial_t B + \rot E  & = 0, & \Div B & = 0.
			\end{aligned}
		\end{cases}
	\end{equation}
	This system satisfies the following formal energy conservation law (see Proposition \ref{energy estimate} for an explicit computation of the energy)~:
	\begin{equation*}
		\begin{aligned}
			\frac 12\frac{d}{dt}\left(\left\|\rho\right\|_{L^2_x}^2 + \left\|u\right\|_{L^2_x}^2 + \frac 32 \left\|\theta\right\|_{L^2_x}^2
			+ \left\|P^\perp E\right\|_{L^2_x}^2 + \left\|B\right\|_{L^2_x}^2 \right) & \\
			+ \mu \left\|\nabla_x u\right\|_{L^2_x}^2 & + \frac 52 \kappa \left\|\nabla_x \theta\right\|_{L^2_x}^2 = 0.
		\end{aligned}
	\end{equation*}
	Here, the solenoidal component of the electric field is defined indirectly as a mere distribution, through Faraday's equation, by
	\begin{equation*}
		PE=-\partial_t A,
	\end{equation*}
	where $B=\rot A$ and $\Div A=0$, while its irrotational component is determined, through Gauss' law, by
	\begin{equation*}
		P^\perp E = \nabla_x\left(\rho+\theta\right).
	\end{equation*}
	Notice that the equations in this system are all coupled.
	
	\item If $\alpha=\eps$ and $\gamma=o(\eps)$, we obtain the incompressible Navier-Stokes-Fourier-Poisson system~:
	\begin{equation}\label{INSFP}
		\begin{cases}
			\begin{aligned}
				\d_t u
				+
				u\cdot\nabla_x u - \mu\Delta_x u
				& = -\nabla_x p + \rho \nabla_x\theta , & \Div u & = 0,\\
				\d_t \left(\frac32\theta-\rho\right)
				+
				u\cdot\nabla_x\left(\frac32\theta-\rho\right) - \frac 52 \kappa \Delta_x\theta
				& = 0, & \Delta_x(\rho+\theta) & =\rho.
			\end{aligned}
		\end{cases}
	\end{equation}
	This system satisfies the following formal energy conservation law~:
	\begin{equation*}
		\begin{aligned}
			\frac 12\frac{d}{dt}\left(\left\|\rho\right\|_{L^2_x}^2 + \left\|u\right\|_{L^2_x}^2 + \frac 32 \left\|\theta\right\|_{L^2_x}^2
			+ \left\|\nabla_x\left(\rho+\theta\right)\right\|_{L^2_x}^2 \right) & \\
			+ \mu \left\|\nabla_x u\right\|_{L^2_x}^2 & + \frac 52 \kappa \left\|\nabla_x \theta\right\|_{L^2_x}^2 & = 0.
		\end{aligned}
	\end{equation*}
	Physically, in this system, the fluid is subject to a self-induced static electric field $E$ determined by
	\begin{equation*}
		\rot E = 0,\qquad \Div E=\rho,
	\end{equation*}
	hence
	\begin{equation*}
		E=\nabla_x(\rho+\theta).
	\end{equation*}
	Notice that the equations in this system are all coupled.
	% \begin{equation*}
	% 	\begin{cases}
	% 		\begin{aligned}
	% 			\d_t u
	% 			+
	% 			u\cdot\nabla_x u - \mu\Delta_x u
	% 			& = -\nabla_x p + \rho E , & \Div u & = 0,\\
	% 			\d_t \left(\frac32\theta-\rho\right)
	% 			+
	% 			\frac 52 u\cdot\nabla_x\theta - \kappa \Delta_x\theta
	% 			& = u \cdot E, & \Div E & = \rho,\\
	% 			\Delta_x\theta & =\left(\operatorname{Id}-\Delta_x\right)\rho, & \rot E & = 0.
	% 		\end{aligned}
	% 	\end{cases}
	% \end{equation*}

\end{enumerate}

% ============================
% = Vlasov-Poisson-Boltzmann =
% ============================

\subsection{The Vlasov-Poisson-Boltzmann system}

The Vlasov-Poisson-Boltzmann system describes the evolution of a gas of one species of charged particles (ions or electrons) subject to an auto-induced electrostatic force. This system is obtained formally from the Vlasov-Maxwell-Boltzmann system by letting the speed of light tend to infinity while all other parameters remain fixed. Accordingly, setting $\gamma = 0$ in \eqref{scaledVMB} yields the scaled Vlasov-Poisson-Boltzmann system~:
\begin{equation}\label{scaledVMP}
	\begin{cases}
		\begin{aligned}
			\epsilon\d_t f_\eps + v \cdot \nabla_x f_\eps + \alpha \nabla_x\phi_\eps \cdot \nabla_v f_\eps &=
			\frac{1}{\epsilon}Q(f_\eps,f_\eps),
			\\
			f_\eps & =M\left(1+\epsilon g_\eps\right),
			\\
			\Delta_x \phi_\eps &=\frac{\alpha}{\epsilon}\int_{\mathbb{R}^3} g_\eps M dv.
		\end{aligned}
	\end{cases}
\end{equation}
Here, the plasma is subject to a self-induced electrostatic field $E_\eps$ determined by
\begin{equation*}
	\rot E_\eps = 0,\qquad \Div E_\eps=\frac{\alpha}{\epsilon}\int_{\mathbb{R}^3} g_\eps M dv,
\end{equation*}
hence
\begin{equation*}
	E_\eps=\nabla_x\phi_\eps.
\end{equation*}
The above system is supplemented with some initial data satisfying
\begin{equation*}
	\frac1{\eps^2} H\left(f_\eps^\mathrm{in}\right)
	+ \frac1{2}\int_{\mathbb{R}^3} |E_\eps^{\rm in}|^2 dx < \infty.
\end{equation*}
In particular, solutions of \eqref{scaledVMP} satisfy the corresponding scaled entropy inequality, where $t>0$,
\begin{equation*}
	\begin{aligned}
		\frac1{\eps^2} H\left(f_\eps \right)
		+ \frac 1{2} \int_{\mathbb{R}^3} |E_\eps|^2 dx
		& +\frac{1}{\epsilon^4}\int_0^t\int_{\mathbb{R}^3}D(f_\eps)(s) dx ds
		\\
		& \leq  \frac1{\eps^2}H\left(f_\eps^{\mathrm{in}}\right)
		+ \frac1{2}\int_{\mathbb{R}^3} |E_\eps^{\rm in}|^2 dx.
	\end{aligned}
\end{equation*}

Thus, the formal asymptotic analysis of \eqref{scaledVMP} is contained in our analysis of the Vlasov-Maxwell-Boltzmann system \eqref{scaledVMB}. Specifically, setting $\gamma=\beta=0$ in the limiting systems obtained in Section \ref{evolution}, we see that the Vlasov-Poisson-Boltzmann system \eqref{scaledVMP} converges, when $\alpha=o(\eps)$, towards the incompressible Navier-Stokes-Fourier system in a Boussinesq regime, with $E=0$~:
\begin{equation*}
	\begin{cases}
		\begin{aligned}
			\d_t u +
			u\cdot\nabla_x u - \mu\Delta_x u
			& = -\nabla_x p , & \Div u & =0\\
			\d_t \theta
			+
			u\cdot\nabla_x\theta - \kappa \Delta_x\theta
			& = 0, & \rho+\theta & =0.
		\end{aligned}
	\end{cases}
\end{equation*}
While, when $\left[\frac\alpha\eps\right]\neq 0$, we find the convergence towards the incompressible Navier-Stokes-Fourier-Poisson system~:
\begin{equation*}
	\begin{cases}
		\begin{aligned}
			\d_t u
			+
			u\cdot\nabla_x u - \mu\Delta_x u
			& = -\nabla_x p && \\
			& + \rho \nabla_x\theta , & \Div u & = 0,\\
			\d_t \left(\frac32\theta-\rho\right)
			+
			u\cdot\nabla_x\left(\frac32\theta-\rho\right) - \frac 52 \kappa \Delta_x\theta
			& = 0, & \Delta_x(\rho+\theta) & = \left[\frac{\alpha}{\eps}\right]^2 \rho,
		\end{aligned}
	\end{cases}
\end{equation*}
where the electrostatic field is determined by $\left[\frac{\alpha}{\eps}\right]E=\nabla_x(\rho+\theta)$.

In fact, the Vlasov-Poisson-Boltzmann system is inherently simpler than the Vlasov-Maxwell-Boltzmann system, because it couples the Vlasov-Boltzmann equation with a simple elliptic equation, namely Poisson's equation, while the Vlasov-Maxwell-Boltzmann system couples the Vlasov-Boltzmann equation with an hyperbolic system, namely Maxwell's system of equations. Thus, the rigorous mathematical analysis on the Vlasov-Maxwell-Boltzmann system, presented in the remainder of this work, will also apply to the Vlasov-Poisson-Boltzmann system and, therefore, analog results will hold.

% ======================================
% = 2 species Vlasov-Maxwell-Boltzmann =
% ======================================

\section{Formal analysis of the two species asymptotics}\label{formal two}

We turn now to the formal asymptotic study of the incompressible viscous regimes of the two species Vlasov-Maxwell-Boltzmann system \eqref{2VMB}. Recall that we are only considering the case of equal masses and opposite charges. The analysis follows exactly the same steps as in the one species case \eqref{VMB}. However, the situation obviously becomes now more complex and general.

% First of all, there are additional free parameters~:
% \begin{itemize}
% \item the mass ratio $\frac{m_+}{m_-}$,
% \item the charge ratio $\frac{q_+}{q_-}$.
% \end{itemize}
% Note that, because of the global neutrality, the charge ratio is related to the ratio of macroscopic densities.
% Obviously the number of such parameters grows linearly in terms of the number of species.
% 
% These additional parameters arise at different places in the equations, for instance as the mass scales the thermal speed for each specie, or because the mean free path is related to the macroscopic density.
% 
% Even worse, the masses appear in the convolution relations defining the mixed collision operators
% $$m_+ v_+' +m_- v'_- = m_+ v _+ +m_- v_-, \quad m_+| v_+' |^2+m_- |v'_- |^2= m_+ |v _+|^2 +m_- |v_-|^2$$
% so that these collision operators become highly singular whenever the mass ratio tends to infinity or to zero.

For a plasma of two species of particles, our starting point is the scaled system
\begin{equation}\label{scaled VMB two species}
	\begin{cases}
		\begin{aligned}
			\eps \d_t f_\eps^\pm + v \cdot \nabla_x f_\eps^\pm \pm \left( \alpha E_\eps + \beta v \wedge B_\eps \right) \cdot \nabla_v f_\eps^\pm &= \frac 1\eps Q(f_\eps^\pm,f_\eps^\pm) + \frac{\delta^2}{\eps} Q(f_\eps^\pm,f_\eps^\mp),
			\\
			f_\eps^\pm & =M\left(1+\epsilon g_\eps^\pm\right),
			\\
			\gamma\d_t E_\eps - \ROT B_\eps &= - \frac{\beta}{\eps} \int_{\mathbb{R}^3} \left(g_\eps^+-g_\eps^-\right)v M dv,
			\\
			\gamma\d_t B_\eps + \ROT E_\eps& = 0,
			\\
			\DIV E_\eps &=\frac{\alpha}{\epsilon}\int_{\mathbb{R}^3} \left(g_\eps^+-g_\eps^-\right) M dv,
			\\
			\DIV B_\eps &=0,
			\\
		\end{aligned}
	\end{cases}
\end{equation}
supplemented with some initial data satisfying
\begin{equation*}
	\frac1{\eps^2} H\left(f_\eps^{+\mathrm{in}}\right)
	+ \frac1{\eps^2} H\left(f_\eps^{-\mathrm{in}}\right)
	+ \frac1{2}\int_{\mathbb{R}^3} |E_\eps^{\rm in}|^2+|B_\eps^{\rm in}|^2  dx < \infty,
\end{equation*}
where $H\left(f_\eps^{\pm\mathrm{in}}\right)=H\left(f_\eps^{\pm\mathrm{in}}|M\right)$. In particular, the corresponding scaled entropy inequality, where $t>0$,
\begin{equation}\label{scaled entropy two species}
	\begin{aligned}
	 	\frac1{\eps^2} H\left(f_\eps^{+}\right)
		+ \frac1{\eps^2} H\left(f_\eps^{-}\right)
		& + \frac 1{2} \int_{\mathbb{R}^3} |E_\eps|^2+ |B_\eps|^2 dx \\
		& +\frac{1}{\epsilon^4}\int_0^t\int_{\mathbb{R}^3}\left(D\left(f_\eps^+\right)+D\left(f_\eps^-\right)
		+ \delta^2 D\left(f_\eps^+,f_\eps^-\right)\right)(s) dx ds
		\\
		& \leq
		\frac1{\eps^2} H\left(f_\eps^{+\mathrm{in}}\right)
		+ \frac1{\eps^2} H\left(f_\eps^{-\mathrm{in}}\right)
		+ \frac1{2}\int_{\mathbb{R}^3} |E_\eps^{\rm in}|^2+ |B_\eps^{\rm in}|^2 dx,
	\end{aligned}
\end{equation}
guarantees that the solution will remain --~for all non-negative times~-- a fluctuation of order $\eps$ around the global equilibrium $M$~:
$$f_\eps^\pm=M(1+\eps g_\eps^\pm).$$
Note that the kinetic equations in \eqref{scaled VMB two species} can then be rewritten, in terms of the fluctuation $g_\eps$, as
\begin{equation}\label{boltz-lin two species}
	\begin{aligned}
		\eps\d_t \begin{pmatrix} g_\eps ^+\\ g_\eps^- \end{pmatrix} + v & \cdot \nabla_x \begin{pmatrix} g_\eps ^+\\ g_\eps^- \end{pmatrix}
		+ (\alpha E_\eps+\beta v\wedge B_\eps) \cdot \nabla_v \begin{pmatrix} g_\eps ^+\\ - g_\eps^- \end{pmatrix}
		- {\alpha\over \eps}  E_\eps \cdot v \begin{pmatrix} 1+\eps g_\eps ^+\\ -1-\eps g_\eps^- \end{pmatrix} \\
		% & = \frac 1\eps
		% \begin{pmatrix} - \cL g_\eps^+ +\eps\cQ(g_\eps^+,g_\eps^+) \\ - \cL g_\eps^- +\eps\cQ(g_\eps^-,g_\eps^-) \end{pmatrix}
		% +
		% \frac{\delta^2}{\eps}
		% \begin{pmatrix} - \cL \left(g_\eps^+,g_\eps^-\right) + \eps \cQ\left(g_\eps^+,g_\eps^-\right) \\ - \cL \left(g_\eps^-,g_\eps^+\right) + \eps \cQ\left(g_\eps^-,g_\eps^+\right) \end{pmatrix}\\
		& = -\frac 1\eps
		\begin{pmatrix} \cL g_\eps^+ + \delta^2\cL \left(g_\eps^+,g_\eps^-\right) \\
		\cL g_\eps^- + \delta^2 \cL \left(g_\eps^-,g_\eps^+\right) \end{pmatrix}
		+
		\begin{pmatrix} \cQ(g_\eps^+,g_\eps^+) + \delta^2 \cQ\left(g_\eps^+,g_\eps^-\right) \\ \cQ(g_\eps^-,g_\eps^-) + \delta^2 \cQ\left(g_\eps^-,g_\eps^+\right) \end{pmatrix},
	\end{aligned}
\end{equation}
where we denote
\begin{equation}\label{def L and Q two species}
	\cL (g,h) =-\frac1M \left(Q(Mg,M)+Q(M,Mh)\right) \quad \text{and} \quad \cQ(g,h) =\frac1M Q(Mg,Mh).
\end{equation}
% \begin{equation}
% \begin{aligned}
% \eps \d_t \begin{pmatrix} g_\eps ^+\\ g_\eps^- \end{pmatrix} \!+&v\cdot \nabla_x \begin{pmatrix} g_\eps ^+\\ g_\eps^- \end{pmatrix} \!+ \eps^{\frac{1-\alpha}2}\!\! \left( ( v\wedge B_\eps)\cdot \nabla_v \begin{pmatrix} g_\eps ^+\\ - g_\eps^- \end{pmatrix}\! - E_\eps\cdot v\begin{pmatrix} 1\\ -1\end{pmatrix}\right)\\
% &=-\frac1\eps \begin{pmatrix}  \cL g_\eps ^+ \\ \cL g_\eps ^-\end{pmatrix} -\frac1{\eps^\alpha}  \begin{pmatrix}  \cL (g_\eps^-,g_\eps ^+) \\ \cL (g_\eps ^+,g_\eps^-) \end{pmatrix}+ \begin{pmatrix} \cQ( g_\eps ^+, g_\eps^+)\\  \cQ( g_\eps ^-, g_\eps^-) \end{pmatrix}+O(\eps^{1-\alpha})\end{aligned}
% \end{equation}

% Gathering all these conditions together, we get the scaled Vlasov-Maxwell Boltzmann system
% \begin{equation}
% \label{VMB2}
% \begin{aligned}
% \ds { \eps \d_t f ^\pm+v\cdot \nabla_x f^\pm}
% \pm{  \eps^{\frac{1-\alpha}2} (\eps E+v\wedge B) 
% \cdot \nabla_v f^\pm }
% = 
% { {1\over\eps}Q(f^\pm,f^\pm)+ {1\over\eps^\alpha} Q(f^\mp,f^\pm)} \\
% \d_t B+\ROT E =0 \\
%  \d_t E-\ROT B =-\eps^{} (\int f^+ vdv+  \int f^- v dv) 
%  \end{aligned}
% \end{equation}
%  The  scaled entropy inequality is therefore
%  $$
% \begin{aligned}
% \frac1{\eps^2} H(f_\eps^+|M) &+\frac1{\eps^2} H(f_\eps^-|M) + \frac12 \int (|E_\eps|^2+|B_\eps|^2)dx  \\
% &+ \int_0^t \int \frac1{\eps^4}\left( D(f_\eps^+)+D(f_\eps^-)\right)  +\frac1{\eps^{3+\alpha}}  D(f_\eps^+, f_\eps^-)  dxds \leq C_0\,.
% \end{aligned}$$
% It gives uniform bounds on $E_\eps$, $B_\eps$ and $g_\eps^\pm$.

\bigskip

It turns out that, in the limit $\eps \to 0$, we will have now three types of constraints~: 
\begin{itemize}
	
	\item conditions on the velocity profiles coming from the fast relaxation towards thermodynamic equilibrium (i.e.\ small Knudsen regime, see Section \ref{constraint-section 1})~;
	
	\item linear macroscopic hydrodynamic constraints due to the weak compressibility (i.e.\ small Mach regime, see Section \ref{constraint-section 2})~;
	
	\item nonlinear macroscopic electrodynamic constraints coming from momentum and energy exchange between species due to interspecies collisions (see Section \ref{ohm section}).
	
\end{itemize}
As in the case of one species of charged particles, we expect the first two types of constraints to be weakly stable,
and thus to be derived from simple uniform a priori estimates. The procedure leading to the last couple of electrodynamic constraint equations (including Ohm's law) is a little bit more complex and will depend on the strength of interspecies collisional interactions, that is to say, on the size of $\delta>0$ compared to $\eps>0$. In fact, the nature of the whole asymptotic systems obtained in the limit $\eps\to 0$ will be conditioned by the size of $\delta$, and we will therefore distinguish three different asymptotic regimes~:
\begin{itemize}
	
	\item Very weak interspecies collisional interactions, $\delta=O(\eps)$~; in this regime, the interspecies collision operators $\frac{\delta^2}{\eps} Q(f_\eps^\pm,f_\eps^\mp)$ in \eqref{scaled VMB two species} are a regular perturbation. Therefore, the corresponding limiting systems will be composed of two hydrodynamic systems --~one for each species~-- coupled mainly through the mean field interactions of the electromagnetic forces. The derivation of these regimes will be easily deduced from the asymptotic analysis for one species from Section \ref{formal one} and will therefore be treated first in Section \ref{very weak section}.
	
	\item Weak interspecies collisional interactions, $\delta=o(1)$ and $\frac\delta\eps$ unbounded~; in this regime, the interspecies collision operators $\frac{\delta^2}{\eps} Q(f_\eps^\pm,f_\eps^\mp)$ in \eqref{scaled VMB two species} are a singular perturbation, whose order may vary from the other singular perturbations present in the system \eqref{scaled VMB two species}. In particular, it is not the most singular perturbation of \eqref{scaled VMB two species}.
	
	\item Strong interspecies collisional interactions, $\delta\sim 1$~; in this regime, the interspecies collision operators $\frac{\delta^2}{\eps} Q(f_\eps^\pm,f_\eps^\mp)$ in \eqref{scaled VMB two species} are a singular perturbation of the most singular order present in the system \eqref{scaled VMB two species}.
	
\end{itemize}

\subsection{Thermodynamic equilibrium}\label{constraint-section 1}

The entropy inequality \eqref{scaled entropy two species} provides uniform bounds on $E_\eps$, $B_\eps$, $g^+_\eps$ and $g^-_\eps$. Therefore, assuming some formal compactness, up to extraction of subsequences, one has
$$
	E_\eps \rightharpoonup E, \qquad B_\eps \rightharpoonup B, \qquad g_\eps^\pm \rightharpoonup g^\pm,
$$
in a weak sense to be rigorously detailed in a subsequent chapter.

Then, multiplying \eqref{boltz-lin two species} by $\eps$, and taking formal limits as $\eps \to 0$ shows that
$$
	\lim_{\eps\to 0}\L_\delta
	\begin{pmatrix} g^+_\eps \\ g^-_\eps \end{pmatrix} =
	\L_{[\delta]}
	\begin{pmatrix} g^+ \\ g^- \end{pmatrix} =
	\begin{pmatrix}
		0\\0
	\end{pmatrix},
$$
where
\begin{equation*}% \label{vector L 2}
	\L_\delta
	\begin{pmatrix} g \\ h \end{pmatrix}
	=
	\begin{pmatrix}
		\cL g + \delta^2\cL \left(g,h\right) \\
		\cL h + \delta^2 \cL \left(h,g\right)
	\end{pmatrix}
	=
	\begin{pmatrix}
		\cL g + \delta^2 \int_{\mathbb{R}^3\times\mathbb{S}^2} \left(g+h_*-g'-h'_*\right) b M_*dv_*d\sigma\\
		\cL h + \delta^2 \int_{\mathbb{R}^3\times\mathbb{S}^2} \left(h+g_*-h'-g_*'\right)  b M_*dv_*d\sigma
	\end{pmatrix},
\end{equation*}
and
\begin{equation*}% \label{vector L}
	\L_{[\delta]}
	\begin{pmatrix} g \\ h \end{pmatrix}
	=
	\begin{pmatrix}
		\cL g + \left[\delta\right]^2\cL \left(g,h\right) \\
		\cL h + \left[\delta\right]^2 \cL \left(h,g\right)
	\end{pmatrix}
	=
	\begin{pmatrix}
		\cL g + \left[\delta\right]^2 \int_{\mathbb{R}^3\times\mathbb{S}^2} \left(g+h_*-g'-h'_*\right) b M_*dv_*d\sigma\\
		\cL h + \left[\delta\right]^2 \int_{\mathbb{R}^3\times\mathbb{S}^2} \left(h+g_*-h'-g_*'\right)  b M_*dv_*d\sigma
	\end{pmatrix}.
\end{equation*}
It can be shown (see Proposition \ref{coercivity 2}) that, when $\left[\delta\right]\neq 0$, the kernel of the vectorial linearized Boltzmann operators $\mathbf{L}_\delta$ and $\mathbf{L}_{[\delta]}$ coincide exactly with the vector space spanned by the set
\begin{equation}\label{collision invariant two species}
	\left\{
	\begin{pmatrix}
		1\\0
	\end{pmatrix},
	\begin{pmatrix}
		0\\1
	\end{pmatrix},
	\begin{pmatrix}
		v_1\\v_1
	\end{pmatrix},
	\begin{pmatrix}
		v_2\\v_2
	\end{pmatrix},
	\begin{pmatrix}
		v_3\\v_3
	\end{pmatrix},
	\begin{pmatrix}
		|v|^2\\|v|^2
	\end{pmatrix}
	\right\}.
\end{equation}
However, when $\left[\delta\right]=0$, the kernel of $\mathbf{L}_{[\delta]}$ is larger and is composed of all vectors $\begin{pmatrix}\varphi_1(v)\\ \varphi_2(v)\end{pmatrix}$ such that $\varphi_1(v)$ and $\varphi_2(v)$ are collision invariants whose coefficients are independent.

Thus, we conclude, if $\left[\delta\right]\neq 0$, that $g^\pm$ is an infinitesimal Maxwellian of the form
\begin{equation}\label{infinitesimal maxwellian 1}
	\begin{pmatrix}
		g^+\\g^-
	\end{pmatrix}
	=
	\begin{pmatrix}
		\rho^++u\cdot v + \theta \left(\frac{|v|^2}{2}-\frac 32\right)
		\\
		\rho^-+u\cdot v + \theta \left(\frac{|v|^2}{2}-\frac 32\right)
	\end{pmatrix},
\end{equation}
while, if $\left[\delta\right]= 0$,
\begin{equation}\label{infinitesimal maxwellian 2}
	\begin{pmatrix}
		g^+\\g^-
	\end{pmatrix}
	=
	\begin{pmatrix}
		\rho^++u^+\cdot v + \theta^+ \left(\frac{|v|^2}{2}-\frac 32\right)
		\\
		\rho^-+u^-\cdot v + \theta^- \left(\frac{|v|^2}{2}-\frac 32\right)
	\end{pmatrix},
\end{equation}
where $\rho^+,\rho^-\in\mathbb{R}$, $u,u^+,u^-\in\mathbb{R}^3$ and $\theta,\theta^+,\theta^-\in\mathbb{R}$ only depend on $t$ and $x$, and are respectively the fluctuations of density, bulk velocity and temperature.

In fact, whenever $\frac\delta\eps$ is unbounded, we show below that necessarily $u^+=u^-$ and $\theta^+=\theta^-$, as well, because of higher order singular limiting constraints. Therefore, the infinitesimal Maxwellian form \eqref{infinitesimal maxwellian 2} will be assumed by the limiting fluctuations in the case $\delta=O(\eps)$ only, that is in the case of very weak interspecies collisions.

The fact that the fluctuations assume the infinitesimal Maxwellian form describes that the gas reaches thermodynamic (or statistical) equilibrium, in the fast relaxation limit.

\bigskip

We define now the macroscopic fluctuations of density $\rho_\eps^\pm$, bulk velocity $u_\eps^\pm$ and temperature $\theta_\eps^\pm$ by
\begin{equation*}
	\begin{aligned}
		\rho_\eps^\pm & =\int_{\mathbb{R}^3}g_\eps^\pm Mdv,\\
		u_\eps^\pm & =\int_{\mathbb{R}^3}g_\eps^\pm v Mdv,\\
		\theta_\eps^\pm & =\int_{\mathbb{R}^3}g_\eps^\pm\left(\frac{|v|^2}{3}-1\right) Mdv,
	\end{aligned}
\end{equation*}
and the hydrodynamic projection $\Pi g_\eps^\pm$ of $g_\eps^\pm$ by
\begin{equation*}
	\Pi g_\eps^\pm = \rho_\eps^\pm+u_\eps^\pm\cdot v + \theta_\eps^\pm \left(\frac{|v|^2}{2}-\frac 32\right),
\end{equation*}
which is nothing but the orthogonal projection of $g_\eps^\pm$ onto the kernel of $\mathcal{L}$ in $L^2\left(Mdv\right)$.

Note that the previous step establishing the convergence of $g_\eps^\pm$ towards thermodynamic equilibrium  yields, in fact, the uniform boundedness of $\frac 1\eps\L_\delta
\begin{pmatrix} g_\eps^+ \\ g_\eps^- \end{pmatrix}$. Therefore, if $\left[\delta\right]\neq 0$, we deduce, at least formally, that
\begin{equation*}
	\begin{pmatrix} g_\eps^+ \\ g_\eps^- \end{pmatrix}
	-
	\begin{pmatrix} \frac{\rho_\eps^+-\rho_\eps^-}{2} + \Pi\frac{g_\eps^++g_\eps^-}{2} \\ \frac{\rho_\eps^--\rho_\eps^+}{2} + \Pi\frac{g_\eps^++g_\eps^-}{2} \end{pmatrix}
	=O(\eps),
\end{equation*}
where $\begin{pmatrix} \frac{\rho_\eps^+-\rho_\eps^-}{2} + \Pi\frac{g_\eps^++g_\eps^-}{2} \\ \frac{\rho_\eps^--\rho_\eps^+}{2} + \Pi\frac{g_\eps^++g_\eps^-}{2} \end{pmatrix}$ clearly defines the orthogonal projection of $\begin{pmatrix} g_\eps^+ \\ g_\eps^- \end{pmatrix}$ onto the kernel of $\L_\delta $, which is spanned by \eqref{collision invariant two species}. This bound implies, in particular, that
\begin{equation*}
	(g_\eps^+-\rho_\eps^+)-(g_\eps^--\rho_\eps^-)=O(\eps)
	\qquad\text{and}\qquad g_\eps^\pm-\Pi g_\eps^\pm = O(\eps).
\end{equation*}
However, if $\left[\delta\right] = 0$, we can only formally deduce, for the moment, that
\begin{equation*}
	g_\eps^\pm-\Pi g_\eps^\pm = O(\eps).
\end{equation*}

Just as in the one species case (see Section \ref{equilibrium conv}), the convergence of $g_\eps^+$ and $g_\eps^-$ with a rate $O(\eps)$ towards their hydrodynamic projections $\Pi g_\eps^+$ and $\Pi g_\eps^-$ can also be inferred, at least formally, from the uniform control of the entropies dissipations $\frac{1}{\eps^4} D\left(f_\eps^+\right)$ and $\frac{1}{\eps^4} D\left(f_\eps^-\right)$ in \eqref{scaled entropy two species}. We are now going to show how the exact same formal reasoning applied to the control of the mixed entropy dissipation $\frac{\delta^2}{\eps^4} D\left(f_\eps^+,f_\eps^-\right)$ in \eqref{scaled entropy two species} yields formally that
\begin{equation}\label{hydro electro projection}
	(g_\eps^+-\rho_\eps^+)-(g_\eps^--\rho_\eps^-)=O\left(\frac\eps\delta\right),
\end{equation}
which is not so readily deduced by direct inspection of \eqref{boltz-lin two species}. Note that this control is relevant in the cases of weak or strong interspecies interactions only, that is when $\frac\delta\eps$ is unbounded.

Thus, as in Section \ref{equilibrium conv}, formally approximating $z\log(1+z)$ by $z^2$, which is valid in a neighborhood of $z=0$, in the definition \eqref{def D mixed} of the mixed entropy dissipation, we deduce a control of
\begin{equation*}
	\frac{\delta^2}{2\eps^4} \int _{\mathbb{R}^3 \times \mathbb{R}^3 \times \mathbb{S}^2}
	\left(\frac{{f_\eps^+}'{f_{\eps *}^-}' - f_\eps^+ f_{\eps *}^-}{f_\eps^+ f_{\eps *}^-}\right)^2
	f_\eps^+ f_{\eps *}^-b
	dvdv_* d\sigma .
\end{equation*}
Then, since
\begin{equation*}
	\begin{aligned}
		{f_\eps^+}' & {f_{\eps *}^-}' - f_\eps^+ f_{\eps *}^- \\
		& =
		\eps \left(\left(\Pi g_\eps^+\right)' + \left(\Pi g_\eps^-\right)_*'
		- \left(\Pi g_\eps^+\right) - \left(\Pi g_\eps^-\right)_*\right) \\
		& + \eps^2 \left(\left(\frac{g_\eps^+-\Pi g_\eps^+}{\eps}\right)' + \left(\frac{g_\eps^--\Pi g_\eps^-}{\eps}\right)_*'
		- \left(\frac{g_\eps^+-\Pi g_\eps^+}{\eps}\right) - \left(\frac{g_\eps^--\Pi g_\eps^-}{\eps}\right)_*\right) \\
		& + \eps^2\left({g_\eps^+}'{g_{\eps *}^-}' - g_\eps^+ g_{\eps *}^-\right),
	\end{aligned}
\end{equation*}
we infer that $\frac \delta {\eps}\left(\left(\Pi g_\eps^+\right)' + \left(\Pi g_\eps^-\right)_*' - \left(\Pi g_\eps^+\right) - \left(\Pi g_\eps^-\right)_*\right)$ is uniformly bounded. Finally, a direct computation of the integral
\begin{equation*}
	\frac{\delta^2}{\eps^2} \int _{\mathbb{R}^3 \times \mathbb{R}^3 \times \mathbb{S}^2}
	\left(\left(\Pi g_\eps^+\right)' + \left(\Pi g_\eps^-\right)_*' - \left(\Pi g_\eps^+\right) - \left(\Pi g_\eps^-\right)_*\right)^2b
	dvdv_* d\sigma,
\end{equation*}
shows that
\begin{equation*}
	u^+-u^-=O\left(\frac\eps\delta\right)\qquad\text{and}\qquad
	\theta^+-\theta^-=O\left(\frac\eps\delta\right),
\end{equation*}
which incidentally establishes \eqref{hydro electro projection}. Of course, the rigorous demonstration of such bounds, later on in Section \ref{bulk}, will necessitate the control of the large values of the fluctuations in order to justify the formal approximation of $z\log (1+z)$ by $z^2$.

On the whole, we have shown that, for all cases of strong, weak and very weak interspecies interactions, it holds
\begin{equation}\label{hydro projection two species}
	g_\eps^\pm-\Pi g_\eps^\pm = O(\eps)
	\qquad\text{and}\qquad
	(g_\eps^+-\rho_\eps^+)-(g_\eps^--\rho_\eps^-)=O\left(\frac\eps\delta\right).
\end{equation}
Note that this implies that
\begin{equation*}
	\begin{aligned}
		\begin{pmatrix} g_\eps^+ \\ g_\eps^- \end{pmatrix}
		& -
		\begin{pmatrix} \frac{\rho_\eps^+-\rho_\eps^-}{2} + \Pi\frac{g_\eps^++g_\eps^-}{2} \\ \frac{\rho_\eps^--\rho_\eps^+}{2} + \Pi\frac{g_\eps^++g_\eps^-}{2} \end{pmatrix} \\
		& =
		\begin{pmatrix} g_\eps^+ \\ g_\eps^- \end{pmatrix}
		-
		\begin{pmatrix} \Pi g_\eps^+ \\ \Pi g_\eps^- \end{pmatrix}
		+
		\begin{pmatrix} \Pi\frac{(g_\eps^+-\rho_\eps^+)-(g_\eps^--\rho_\eps^-)}{2} \\ - \Pi\frac{(g_\eps^+-\rho_\eps^+)-(g_\eps^--\rho_\eps^-)}{2} \end{pmatrix}
		=O\left(\frac\eps\delta\right).
	\end{aligned}
\end{equation*}
We will therefore henceforth denote, when considering weak or strong interspecies collisions,
\begin{equation*}
	\begin{pmatrix} h_\eps^+ \\ h_\eps^- \end{pmatrix}
	=
	\frac{\delta}{\eps}\left[\begin{pmatrix} g_\eps^+ \\ g_\eps^- \end{pmatrix}
	-
	\begin{pmatrix} \frac{\rho_\eps^+-\rho_\eps^-}{2} + \Pi\frac{g_\eps^++g_\eps^-}{2} \\ \frac{\rho_\eps^--\rho_\eps^+}{2} + \Pi\frac{g_\eps^++g_\eps^-}{2} \end{pmatrix} \right].
\end{equation*}
In particular, further note that, for weak interspecies collisions, that is whenever $\delta=o(1)$ and $\frac\delta\eps$ is unbounded,
\begin{equation}\label{h equilibrium}
		\lim_{\eps\rightarrow 0}\begin{pmatrix} h_\eps^+ \\ h_\eps^- \end{pmatrix}
		-
		\begin{pmatrix} \Pi h_\eps^+ \\ \Pi h_\eps^-
		\end{pmatrix}=
		\lim_{\eps\rightarrow 0}\begin{pmatrix} h_\eps^+ \\ h_\eps^- \end{pmatrix}
		-
		\frac{\delta}{\eps} \begin{pmatrix} \Pi\frac{(g_\eps^+-\rho_\eps^+)-(g_\eps^--\rho_\eps^-)}{2} \\ - \Pi\frac{(g_\eps^+-\rho_\eps^+)-(g_\eps^--\rho_\eps^-)}{2}
		\end{pmatrix}=0,
\end{equation}
so that the weak limits $h_\eps^\pm\rightharpoonup h^\pm$ are necessarily infinitesimal Maxwellians. But this does not seem to hold for strong interspecies interactions, that is $\delta\sim 1$.

\bigskip

In light of the above formal controls, we define new macroscopic hydrodynamic variables
\begin{equation*}
	\rho_\eps = \frac{\rho^+_\eps+\rho^-_\eps}{2}, \qquad
	u_\eps = \frac{u^+_\eps + u^-_\eps}{2}, \qquad
	\theta_\eps = \frac{\theta^+_\eps+\theta_\eps^-}{2},
\end{equation*}
and electrodynamic variables (irrelevant for very weak interspecies collisions because $\frac\delta\eps$ is bounded in this case)
\begin{equation*}
	n_\eps = \rho^+_\eps - \rho^-_\eps, \qquad
	j_\eps = \frac\delta\eps\left(u^+_\eps - u^-_\eps\right), \qquad
	w_\eps = \frac\delta\eps\left(\theta^+_\eps -\theta^-_\eps\right),
\end{equation*}
namely the electric charge $n_\eps$, the electric current $j_\eps$ and the internal electric energy $w_\eps$. We will also consider their formal weak limits
$$
	\rho_\eps \rightharpoonup \rho, \qquad u_\eps \rightharpoonup u, \qquad \theta_\eps \rightharpoonup \theta,
	\qquad n_\eps \rightharpoonup n, \qquad j_\eps \rightharpoonup j, \qquad w_\eps \rightharpoonup w.
$$
Notice that
\begin{equation*}
	\Pi h_\eps^\pm=\pm \frac 12\left(j_\eps\cdot v + w_\eps\left(\frac{|v|^2}{2}-\frac 32\right)\right),
\end{equation*}
hence, for weak interspecies collisions,
\begin{equation*}% \label{h limit 1}
		\lim_{\eps\rightarrow 0}h_\eps^\pm =h^\pm
		=\pm \frac 12\left(j\cdot v + w\left(\frac{|v|^2}{2}-\frac 32\right)\right),
\end{equation*}
whereas, for strong interspecies collisions, we only have that
\begin{equation}\label{h limit 2}
		\lim_{\eps\rightarrow 0}\Pi h_\eps^\pm =\Pi h^\pm
		=\pm \frac 12\left(j\cdot v + w\left(\frac{|v|^2}{2}-\frac 32\right)\right).
\end{equation}

\bigskip

The asymptotic dynamics of $(\rho_\eps^+,u_\eps^+,\theta_\eps^+,\rho_\eps^-,u_\eps^-, \theta_\eps^-)$, or equivalently $(\rho_\eps,u_\eps,\theta_\eps,q_\eps,j_\eps, w_\eps)$, is then governed by fluid equations, to be obtained from the moments equations associated with \eqref{boltz-lin two species}. Thus, successively multiplying \eqref{boltz-lin two species} by the collision invariants $1$, $v$ and $\frac{|v|^2}{2}$, and integrating in $Mdv$, yields
\begin{equation}\label{two fluid system eps}
	\begin{cases}
		\begin{aligned}
			\d_t \rho_\eps^\pm +\frac1\eps \DIV u_\eps^\pm & = 0,\\
			\d_t u_\eps^\pm +\frac1\eps \nabla_x \left( \rho_\eps^\pm+\theta_\eps^\pm \right)
			\mp \frac{\alpha}{\eps^2} E_\eps
			+\frac{\delta^2}{\eps^2}\int_{\mathbb{R}^3}\mathcal{L}\left(g_\eps^\pm,g_\eps^\mp\right)vMdv
			& = \\
			&  \hspace{-50mm} \pm \left(\frac\alpha\eps\rho_\eps^\pm E_\eps +\frac\beta\eps u_\eps^\pm \wedge B_\eps\right) \\
			& \hspace{-50mm} -\frac1\eps \DIV \int_{\mathbb{R}^3} g_\eps^\pm \phi M dv
			+\frac{\delta^2}{\eps}\int_{\mathbb{R}^3}\mathcal{Q}\left(g_\eps^\pm,g_\eps^\mp\right) vMdv ,\\
			\frac 32 \d_t \left(\rho_\eps^\pm+\theta_\eps^\pm\right) +\frac 5{2\eps} \DIV u_\eps^\pm
			+\frac{\delta^2}{\eps^2}\int_{\mathbb{R}^3}\mathcal{L}\left(g_\eps^\pm,g_\eps^\mp\right)\frac{|v|^2}{2}Mdv
			& =
			\pm \frac{\alpha}{\eps} u_\eps^\pm \cdot E_\eps \\
			&  \hspace{-50mm} - \frac1\eps \DIV \int_{\mathbb{R}^3}g_\eps^\pm \psi M dv
			+\frac{\delta^2}{\eps}\int_{\mathbb{R}^3}\mathcal{Q}\left(g_\eps^\pm,g_\eps^\mp\right)\frac{|v|^2}{2}Mdv,
		\end{aligned}
	\end{cases}
\end{equation}
where $\phi(v)$ and $\psi(v)$ have already been defined in \eqref{phi-psi-def}. The above system will be used in the case of very weak interspecies interactions only. For weak and strong interspecies interactions, the evolution equations can then be recast, in terms of the new hydrodynamic and electrodynamic variables, as
\begin{equation}\label{moment-eps0 two species}
	\begin{cases}
		\begin{aligned}
			\d_t \rho_\eps +\frac1\eps \DIV u_\eps & = 0,\\
			\d_t u_\eps +\frac1\eps \nabla_x \left( \rho_\eps+\theta_\eps \right)
			& =
			\left(\frac\alpha{2\eps} n_\eps E_\eps +\frac\beta{2\delta} j_\eps \wedge B_\eps\right) \\
			& -\frac1\eps \DIV \int_{\mathbb{R}^3} \frac{g_\eps^++g_\eps^-}{2} \phi M dv,\\
			\frac 32 \d_t \theta_\eps +\frac 1{\eps} \DIV u_\eps
			& =
			\frac{\alpha}{2\delta} j_\eps \cdot E_\eps
			- \frac1\eps \DIV \int_{\mathbb{R}^3}\frac{g_\eps^++g_\eps^-}{2}\psi M dv.
		\end{aligned}
	\end{cases}
\end{equation}

Recall that we are assuming $\alpha=O(\eps)$ and $\beta=O(\eps)$ for very weak interspecies collisions, i.e.\ when $\delta=O(\eps)$, and $\alpha=O(\eps)$ and $\beta=O(\delta)$ for weak and strong interspecies collisions, i.e.\ when $\frac\delta\eps$ is unbounded. Hence, the nonlinear terms in the right-hand side of \eqref{two fluid system eps}, for very weak interspecies collisions, and \eqref{moment-eps0 two species}, for weak and strong interspecies collisions, containing the electromagnetic fields are expected to be bounded. Furthermore, just as in the case of one species, the polynomials $\phi(v)$ and $\psi(v)$ are orthogonal to the collision invariants in the $L^2(Mdv)$ inner-product. That is to say $\int_{\mathbb{R}^3} \varphi \phi M dv=0$ and $\int_{\mathbb{R}^3} \varphi \psi M dv=0$, for all collision invariants $\varphi(v)$. Since, according to \eqref{hydro projection two species}, the fluctuations $g_\eps^+$ and $g_\eps^-$ converge towards infinitesimal Maxwellians with a rate $O(\eps)$, it is therefore natural to expect, at least formally, that the terms
\begin{equation*}
	\begin{aligned}
		\frac1\eps \int_{\mathbb{R}^3} g_\eps^\pm \phi M dv & = \frac1\eps \int_{\mathbb{R}^3} \left(g_\eps^\pm
		-\Pi g_\eps^\pm \right) \phi M dv, \\
		\frac1\eps \int_{\mathbb{R}^3} g_\eps^\pm \psi M dv & = \frac1\eps \int_{\mathbb{R}^3} \left( g_\eps^\pm
		-\Pi g_\eps^\pm \right) \psi M dv,
	\end{aligned}
\end{equation*}
in \eqref{two fluid system eps} and \eqref{moment-eps0 two species} are bounded and have a limit.

Thus, following the strategy for one species in Section \ref{equilibrium conv}, we rewrite \eqref{two fluid system eps} as
\begin{equation}\label{moment-eps two species very weak}
	\begin{cases}
		\begin{aligned}
			\d_t \rho_\eps^\pm +\frac1\eps \DIV u_\eps^\pm & = 0,\\
			\d_t u_\eps^\pm +\frac1\eps \nabla_x \left( \rho_\eps^\pm+\theta_\eps^\pm \right)
			\mp \frac{\alpha}{\eps^2} E_\eps
			+\frac{\delta^2}{\eps^2}\int_{\mathbb{R}^3}\mathcal{L}\left(g_\eps^\pm,g_\eps^\mp\right)vMdv
			& = \\
			&  \hspace{-52mm} \pm \left(\frac\alpha\eps\rho_\eps^\pm E_\eps +\frac\beta\eps u_\eps^\pm \wedge B_\eps\right) \\
			& \hspace{-52mm} -\frac1\eps \DIV \int_{\mathbb{R}^3} \mathcal{L}g_\eps^\pm \tilde\phi M dv
			+\frac{\delta^2}{\eps}\int_{\mathbb{R}^3}\mathcal{Q}\left(g_\eps^\pm,g_\eps^\mp\right) vMdv ,\\
			\frac 32 \d_t \theta_\eps^\pm +\frac 1{\eps} \DIV u_\eps^\pm
			+\frac{\delta^2}{\eps^2}\int_{\mathbb{R}^3}\mathcal{L}\left(g_\eps^\pm,g_\eps^\mp\right)\frac{|v|^2}{2}Mdv
			& =
			\pm \frac{\alpha}{\eps} u_\eps^\pm \cdot E_\eps \\
			&  \hspace{-52mm} - \frac1\eps \DIV \int_{\mathbb{R}^3}\mathcal{L}g_\eps^\pm \tilde\psi M dv
			+\frac{\delta^2}{\eps}\int_{\mathbb{R}^3}\mathcal{Q}\left(g_\eps^\pm,g_\eps^\mp\right)\frac{|v|^2}{2}Mdv,
		\end{aligned}
	\end{cases}
\end{equation}
and \eqref{moment-eps0 two species} as
\begin{equation}\label{moment-eps two species}
	\begin{cases}
		\begin{aligned}
			\d_t \rho_\eps +\frac1\eps \DIV u_\eps & = 0,\\
			\d_t u_\eps +\frac1\eps \nabla_x \left( \rho_\eps+\theta_\eps \right)
			& =
			\left(\frac\alpha{2\eps} n_\eps E_\eps +\frac\beta{2\delta} j_\eps \wedge B_\eps\right) \\
			& -\frac1{2\eps} \DIV \int_{\mathbb{R}^3} \mathcal{L}\left(g_\eps^++g_\eps^-\right) \tilde\phi M dv,\\
			\frac 32 \d_t \theta_\eps +\frac 1{\eps} \DIV u_\eps
			& =
			\frac{\alpha}{2\delta} j_\eps \cdot E_\eps
			- \frac1{2\eps} \DIV \int_{\mathbb{R}^3}\mathcal{L}\left(g_\eps^++g_\eps^-\right)\tilde\psi M dv,
		\end{aligned}
	\end{cases}
\end{equation}
where $\tilde \phi$ and $\tilde \psi$ are the pseudo-inverses of $\phi$ and $\psi$, respectively, defined in \eqref{phi-psi-def inverses}, and where the terms $\frac 1\eps \mathcal{L}g_\eps^\pm$ will be expressed employing the Vlasov-Boltzmann equations \eqref{boltz-lin two species}. Each of the above macroscopic systems \eqref{moment-eps two species very weak} and \eqref{moment-eps two species} is coupled with Maxwell's equations on $E_\eps$ and $B_\eps$~:
\begin{equation}\label{Maxwell-eps two species}
	\begin{cases}
		\begin{aligned}
			\gamma\d_t E_\eps - \ROT B_\eps &= - \frac{\beta}{\eps}\left(u_\eps^+ -u_\eps^-\right)
			= - \frac{\beta}{\delta}j_\eps,
			\\
			\gamma\d_t B_\eps + \ROT E_\eps& = 0,
			\\
			\DIV E_\eps &=\frac{\alpha}{\epsilon} \left(\rho_\eps^+-\rho_\eps^-\right)
			=\frac{\alpha}{\epsilon} n_\eps,
			\\
			\DIV B_\eps &=0.
		\end{aligned}
	\end{cases}
\end{equation}

A careful formal analysis of the whole coupled macroscopic systems \eqref{moment-eps two species very weak}-\eqref{Maxwell-eps two species}, for very weak interspecies collisions, and \eqref{moment-eps two species}-\eqref{Maxwell-eps two species}, for weak and strong interspecies collisions, will yield the asymptotic dynamics of $(\rho^\pm,u^\pm,\theta^\pm,E,B)$ and $(\rho,u,\theta,E,B)$, respectively. However, note that, in the case of weak or strong interspecies collisions only, the above coupled system \eqref{moment-eps two species}-\eqref{Maxwell-eps two species} remains underdetermined, as the evolution for $n_\eps$, $j_\eps$ and $w_\eps$ is missing. It turns out that the electrodynamic variables will be determined by nonlinear constraint equations. In particular, $j_\eps$ will be asymptotically determined by the so-called Ohm's law, which we derive below in Section \ref{ohm section}.

\subsection{The case of very weak interspecies collisions}\label{very weak section}

The reader should, at this point, take some time to compare the two species system \eqref{moment-eps two species very weak}-\eqref{Maxwell-eps two species} with the one species system \eqref{moment-eps}-\eqref{Maxwell-eps}. When $\delta=O(\eps)$, the coupling between cations and anions in the two species system \eqref{moment-eps two species very weak}-\eqref{Maxwell-eps two species} is caused only by the mean field interaction of the electromagnetic field $(E,B)$ and by the low order interspecies collision terms
\begin{equation*}
	\frac{\delta^2}{\eps^2}\int_{\mathbb{R}^3}\mathcal{L}\left(g_\eps^\pm,g_\eps^\mp\right)
	\begin{pmatrix}
		v \\ \frac{|v|^2}{2}
	\end{pmatrix}
	Mdv
	\qquad\text{and}\qquad
	\frac{\delta^2}{\eps}\int_{\mathbb{R}^3}\mathcal{Q}\left(g_\eps^\pm,g_\eps^\mp\right)
	\begin{pmatrix}
		v \\ \frac{|v|^2}{2}
	\end{pmatrix}
	Mdv.
\end{equation*}
As we are about to discuss, the system \eqref{moment-eps two species very weak}-\eqref{Maxwell-eps two species} essentially behaves, in the limit $\eps\to 0$, as two coupled one species systems of the kind \eqref{moment-eps}-\eqref{Maxwell-eps}.

Indeed, when compared with \eqref{moment-eps}, the only additional terms that one finds in \eqref{moment-eps two species very weak} are~:
\begin{itemize}
	
	\item The linear interspecies collision terms $\frac{\delta^2}{\eps^2}\int_{\mathbb{R}^3}\mathcal{L}\left(g_\eps^\pm,g_\eps^\mp\right)
	\begin{pmatrix}
		v \\ \frac{|v|^2}{2}
	\end{pmatrix}
	Mdv$, which converge, as $\eps\to 0$, towards
	\begin{equation*}
		\left[\frac{\delta}{\eps}\right]^2\int_{\mathbb{R}^3}\mathcal{L}\left(g^\pm,g^\mp\right)
		\begin{pmatrix}
			v \\ \frac{|v|^2}{2}
		\end{pmatrix}
		Mdv =
		\pm\left[\frac{\delta}{\eps}\right]^2
		\begin{pmatrix}
			\sigma^{-1}\left(u^+-u^-\right) \\ \frac 52\lambda^{-1}\left(\theta^+ - \theta^-\right)
		\end{pmatrix},
	\end{equation*}
	where the electrical conductivity $\sigma>0$ and the energy conductivity $\lambda>0$ are constants defined by
	\begin{equation*}
		\begin{aligned}
			\frac 1\sigma & =\frac 12 \int_{\mathbb{R}^3}v\cdot\mathcal{L}\left(v , - v \right)Mdv \\
			& =\frac 12\int_{\mathbb{R}^3\times\mathbb{R}^3\times\mathbb{S}^2}\left|v-v'\right|^2 b(v-v_*,\sigma) MM_* dvdv_*d\sigma \\
			& =\frac 1{2}\int_{\mathbb{R}^3\times\mathbb{R}^3}\left|v-v_*\right|^2 
			m(v-v_*) MM_* dvdv_*,
		\end{aligned}
	\end{equation*}
	and
	\begin{equation*}
		\begin{aligned}
			\frac 1\lambda & =\frac 1{20} \int_{\mathbb{R}^3}|v|^2\mathcal{L}\left(|v|^2,
			-|v|^2\right)Mdv \\
			& =\frac 1{20} \int_{\mathbb{R}^3\times\mathbb{R}^3\times\mathbb{S}^2}\left(|v|^2-|v'|^2\right)^2 b(v-v_*,\sigma) MM_* dvdv_*d\sigma \\
			& =\frac 1{20}\int_{\mathbb{R}^3\times\mathbb{R}^3}\left(|v|^2-|v_*|^2\right)^2 
			m(v-v_*) MM_* dvdv_*,
		\end{aligned}
	\end{equation*}
	where the cross-section for momentum and energy transfer $m(v-v_*)$ is defined in Proposition \ref{cross section transfer}.
	
	\item The nonlinear interspecies collisions terms $\frac{\delta^2}{\eps}\int_{\mathbb{R}^3}\mathcal{Q}\left(g_\eps^\pm,g_\eps^\mp\right)
	\begin{pmatrix}
		v \\ \frac{|v|^2}{2}
	\end{pmatrix}
	Mdv$, which are at least of formal order $O(\eps)$ and thus, vanish in the limit $\eps\to 0$.
	
\end{itemize}

Thus, the remainder of the formal asymptotic analysis of the two species system \eqref{moment-eps two species very weak}-\eqref{Maxwell-eps two species} follows exactly the same steps as the analysis of the one species system \eqref{moment-eps}-\eqref{Maxwell-eps} performed in Sections \ref{macro constraint one species} and \ref{evolution}, which we somewhat detail now.

Note first that the system \eqref{moment-eps two species very weak}-\eqref{Maxwell-eps two species} can be rewritten as a singular perturbation
\begin{equation*}
	\d_t
	\begin{pmatrix}
		\rho_\eps^+\\ u_\eps^+ \\ \sqrt{3\over 2}\theta_\eps^+ \\ \rho_\eps^- \\ u_\eps^-\\ \sqrt{3\over 2}\theta_\eps^- \\ E_\eps \\ B_\eps
	\end{pmatrix}
	+W_\eps
	\begin{pmatrix}
		\rho_\eps^+\\ u_\eps^+ \\ \sqrt{3\over 2}\theta_\eps^+ \\ \rho_\eps^- \\ u_\eps^-\\ \sqrt{3\over 2}\theta_\eps^- \\ E_\eps \\ B_\eps
	\end{pmatrix}
	=O(1),
\end{equation*}
which describes the propagation of waves in the system, where the wave operator is given by
\begin{equation*}
	W_\eps=
	\begin{pmatrix}
		0 &\frac1\eps \DIV&0&0&0&0&0&0\\
		\frac1\eps \nabla_x&0& \frac 1\eps\sqrt{2\over 3} \nabla_x&0&0&0& -{\alpha \over \eps^2} \operatorname{Id} & 0\\
		0 & \frac 1\eps\sqrt{2\over 3} \DIV &0&0&0&0&0&0\\
		0&0&0&0 &\frac1\eps \DIV&0&0&0\\
		0&0&0&\frac1\eps \nabla_x&0&\frac 1\eps\sqrt{2\over 3} \nabla_x& {\alpha \over \eps^2} \operatorname{Id} & 0\\
		0&0&0&0 & \frac 1\eps\sqrt{2\over 3} \DIV &0&0&0\\
		0&{\alpha \over \eps^2} \operatorname{Id}&0&0&-{\alpha \over \eps^2} \operatorname{Id}&0&0&-{1 \over \gamma} \ROT\\
		0&0&0&0&0&0&{1 \over \gamma} \ROT &0
	\end{pmatrix}.
\end{equation*}
We derive then the macroscopic constraint equations on $\left(\rho^\pm,u^\pm,\theta^\pm,E,B\right)$ reproducing the reasoning from Section \ref{macro constraint one species}.
\begin{enumerate}
	
	\item When $\gamma\sim 1$ (so that $\alpha=O(\eps^2)$), averaging over fast time oscillations as $\eps\to 0$, we get the macroscopic constraints
	\begin{equation*}
		\DIV u^\pm =0, \qquad \rho^\pm+\theta^\pm =0,
	\end{equation*}
	respectively referred to as incompressibility and Boussinesq relations. These are supplemented by the asymptotic constraints coming from Gauss' laws in \eqref{Maxwell-eps two species}
	\begin{equation*}
		\DIV E=0,\qquad \DIV B=0.
	\end{equation*}
	
	\item When $\gamma =o(1)$ and $\alpha=O(\eps^2)$, averaging over fast time oscillations as $\epsilon\to 0$, we get the macroscopic constraints
	\begin{equation*}
		\begin{aligned}
			\DIV u^\pm & =0, & \rho^\pm+\theta^\pm & = 0,\\
			\ROT B & = 0, & \ROT E & = 0.
		\end{aligned}
	\end{equation*}
	These are supplemented by the asymptotic constraints coming from Gauss' laws in \eqref{Maxwell-eps two species}
	\begin{equation*}
		\DIV E=0,\qquad \DIV B=0.
	\end{equation*}
	Hence,
	\begin{equation*}
		E=0,\qquad B=0.
	\end{equation*}
	
	\item When $\gamma =o(1)$ and $\frac\alpha{\eps^2}$ is unbounded, averaging over fast time oscillations as $\epsilon\to 0$, we get the macroscopic constraints
	\begin{equation*}
		\begin{aligned}
			\DIV u^\pm & =0, & \nabla_x \left(\rho^\pm+\theta^\pm\right) & = \pm\left[\frac{\alpha}{\eps}\right] E,\\
			\ROT B & = \left[\frac{\beta}{\eps}\right] \left(u^+-u^-\right), & \ROT E & = 0.
		\end{aligned}
	\end{equation*}
	As usual, when $\alpha=o(\eps)$, the weak Boussinesq relation $\nabla_x \left(\rho^\pm+\theta^\pm\right)=0$ can be improved to the strong Boussinesq relation $\rho^\pm+\theta^\pm=0$, assuming $\rho^\pm$ and $\theta^\pm$ enjoy enough integrability. These are supplemented by the asymptotic constraints coming from Gauss' laws in \eqref{Maxwell-eps two species}
	\begin{equation*}
		\DIV E=\left[\frac{\alpha}{\eps}\right]\left(\rho^+-\rho^-\right),\qquad \DIV B=0.
	\end{equation*}
	
\end{enumerate}

Next, following the reasoning from Section \ref{evolution}, we derive the asymptotic evolution equations associated with the two species system \eqref{moment-eps two species very weak}-\eqref{Maxwell-eps two species}. To this end, notice that \eqref{boltz-lin two species} implies, in particular, that
\begin{equation*}
	\frac1\eps \cL g_\eps^\pm =
	\cQ(g_\eps^\pm,g_\eps^\pm) - v\cdot \nabla_x g_\eps^\pm
	\pm{\alpha\over \eps}  E_\eps \cdot v + O(\eps),
\end{equation*}
which is analog to \eqref{fluxes sub} in the one species case. Hence, we obtain the advection and diffusion terms, as in Section \ref{evolution},
\begin{equation*}
	\begin{aligned}
		\lim_{\eps\to 0} \frac1\eps P\Div \int_{\mathbb{R}^3} \mathcal{L}g_\eps^\pm \tilde\phi M dv
		& =
		P\left(u^\pm\cdot\nabla_x u^\pm\right) - \mu\Delta_x u^\pm,
		\\
		\lim_{\eps\to 0} \frac1\eps \Div\int_{\mathbb{R}^3}\mathcal{L}g_\eps^\pm \tilde\psi M dv
		& =
		\frac 52 u^\pm\cdot\nabla_x\theta^\pm - \frac 52 \kappa \Delta_x\theta^\pm,
	\end{aligned}
\end{equation*}
where
\begin{equation*}
	\mu = \frac{1}{10}\int_{\mathbb{R}^3}\phi : \tilde\phi M dv
	\qquad\text{and}\qquad
	\kappa = \frac 2{15} \int_{\mathbb{R}^3}\psi\cdot\tilde\psi M dv.
\end{equation*}

We are now in a position to obtain the limiting evolution of the system \eqref{moment-eps two species very weak}-\eqref{Maxwell-eps two species}.
\begin{enumerate}
	
	\item When $\gamma \sim 1$, letting $\epsilon$ tend to zero in the system \eqref{moment-eps two species very weak}-\eqref{Maxwell-eps two species} yields
	\begin{equation*}
		\begin{cases}
			\begin{aligned}
				\d_t u^\pm +
				u^\pm\cdot\nabla_x u^\pm - \mu\Delta_x u^\pm \hspace{20mm} & \\
				\pm \left[\frac\delta\eps\right]^2\frac1\sigma\left(u^+-u^-\right)
				& = -\nabla_x p^\pm
				\pm \left[\frac{\alpha}{\eps^2}\right]E \pm \left[\frac\beta\eps\right] u^\pm \wedge B ,\\
				\d_t \theta^\pm
				+
				u^\pm\cdot\nabla_x\theta^\pm - \kappa \Delta_x\theta^\pm \hspace{20mm} & \\
				\pm \left[\frac\delta\eps\right]^2\frac1\lambda\left(\theta^+-\theta^-\right)
				& = 0, \\
				\left[\gamma\right] \d_t E - \ROT B &= - \left[\frac{\beta}{\epsilon}\right] \left(u^+-u^-\right),
				\\
				\left[\gamma\right] \d_t B + \ROT E & = 0,
			\end{aligned}
		\end{cases}
	\end{equation*}
	with the constraints
	\begin{equation*}
		\begin{aligned}
			\Div u^\pm & =0, & \rho^\pm+\theta^\pm & = 0,\\
			\Div E & = 0 , & \Div B & = 0.
		\end{aligned}
	\end{equation*}
	
	\item When $\gamma =o(1)$ and $\alpha=O(\eps^2)$, letting $\epsilon$ tend to zero in the system \eqref{moment-eps two species very weak}-\eqref{Maxwell-eps two species} yields
	\begin{equation*}
		\begin{cases}
			\begin{aligned}
				\d_t u^\pm +
				u^\pm\cdot\nabla_x u^\pm - \mu\Delta_x u^\pm \pm \left[\frac\delta\eps\right]^2\frac1\sigma\left(u^+-u^-\right)
				& = -\nabla_x p^\pm ,\\
				\d_t \theta^\pm
				+
				u^\pm\cdot\nabla_x\theta^\pm - \kappa \Delta_x\theta^\pm
				\pm \left[\frac\delta\eps\right]^2\frac1\lambda\left(\theta^+-\theta^-\right)
				& = 0,
			\end{aligned}
		\end{cases}
	\end{equation*}
	with the constraints
	\begin{equation*}
		\begin{aligned}
			\Div u^\pm & =0, & \rho^\pm+\theta^\pm & = 0,\\
			E & = 0 , & B & = 0.
		\end{aligned}
	\end{equation*}
	
	\item When $\gamma =o(1)$ and $\frac\alpha{\eps^2}$ is unbounded, letting $\epsilon$ tend to zero in the system \eqref{moment-eps two species very weak}-\eqref{Maxwell-eps two species} yields
	\begin{equation*}
		\begin{cases}
			\begin{aligned}
				\d_t \left(u^\pm \pm \left[\frac{\beta}{\eps}\right]A \right)
				+
				u^\pm\cdot\nabla_x u^\pm - \mu\Delta_x u^\pm \hspace{-5mm} & \\
				\pm \left[\frac\delta\eps\right]^2\frac1\sigma\left(u^+-u^-\right)
				& = -\nabla_x p^\pm
				\pm \left[\frac\alpha\eps\right]\rho^\pm E \pm \left[\frac\beta\eps\right] u^\pm \wedge B ,\\
				\d_t \left(\frac32\theta^\pm-\rho^\pm\right)
				+
				\frac 52 u^\pm\cdot\nabla_x\theta^\pm - \frac 52 \kappa \Delta_x\theta^\pm \hspace{-5mm} & \\
				\pm\frac 52 \left[\frac\delta\eps\right]^2\frac1\lambda\left(\theta^+-\theta^-\right)
				& =
				\pm\left[\frac{\alpha}{\eps}\right] u^\pm \cdot E,
			\end{aligned}
		\end{cases}
	\end{equation*}
	with the constraints
	\begin{equation*}
		\begin{aligned}
			\DIV u^\pm & =0, & \nabla_x \left(\rho^\pm+\theta^\pm\right) & = \pm\left[\frac{\alpha}{\eps}\right] E,\\
			\ROT B & = \left[\frac{\beta}{\eps}\right] \left(u^+-u^-\right), & \ROT E & = 0,\\
			\DIV E & =\left[\frac{\alpha}{\eps}\right]\left(\rho^+-\rho^-\right), & \DIV B & =0, \\
			\rot A & = B, & \Div A & = 0.
		\end{aligned}
	\end{equation*}
	The above system can be rewritten more explicitly by defining the adjusted electric field $\tilde E = - \partial_t A + E$. It then holds that
	\begin{equation*}
		\begin{cases}
			\begin{aligned}
				\d_t u^\pm
				+
				u^\pm\cdot\nabla_x u^\pm - \mu\Delta_x u^\pm
				\pm \left[\frac\delta\eps\right]^2\frac1\sigma\left(u^+-u^-\right)
				& = -\nabla_x p^\pm \pm \left[\frac{\beta}{\eps}\right]\tilde E \\
				& +\rho^\pm\nabla_x\theta^\pm \pm \left[\frac\beta\eps\right] u^\pm \wedge B ,\\
				\d_t \left(\frac32\theta^\pm-\rho^\pm\right)
				+
				u^\pm\cdot\nabla_x\left(\frac32\theta^\pm-\rho^\pm\right) - \frac 52 \kappa \Delta_x\theta^\pm \hspace{-1mm} & \\
				\pm\frac 52 \left[\frac\delta\eps\right]^2\frac1\lambda\left(\theta^+-\theta^-\right)
				& = 0, \\
				\partial_t B + \rot \tilde E  & = 0,
			\end{aligned}
		\end{cases}
	\end{equation*}
	with the constraints
	\begin{equation*}
		\begin{aligned}
			\DIV u^\pm & =0, & \Delta_x \left(\rho^\pm+\theta^\pm\right) & = \pm\left[\frac{\alpha}{\eps}\right]^2\left(\rho^+-\rho^-\right),\\
			\ROT B & = \left[\frac{\beta}{\eps}\right] \left(u^+-u^-\right), & \DIV B & = 0,\\
			&&\DIV\tilde E & =\left[\frac{\alpha}{\eps}\right]\left(\rho^+-\rho^-\right).
		\end{aligned}
	\end{equation*}

	Notice, finally, that if further $\alpha=o(\eps)$, then the above system is greatly simplified and becomes
	\begin{equation*}
		\begin{cases}
			\begin{aligned}
				\d_t u^\pm
				+
				u^\pm\cdot\nabla_x u^\pm - \mu\Delta_x u^\pm \hspace{20mm} & \\
				\pm \left[\frac\delta\eps\right]^2\frac1\sigma\left(u^+-u^-\right)
				& = -\nabla_x p^\pm \pm \left[\frac{\beta}{\eps}\right]\tilde E
				\pm \left[\frac\beta\eps\right] u^\pm \wedge B ,\\
				\d_t \theta^\pm
				+
				u^\pm\cdot\nabla_x \theta^\pm - \kappa \Delta_x\theta^\pm \hspace{20mm} & \\
				\pm \left[\frac\delta\eps\right]^2\frac1\lambda\left(\theta^+-\theta^-\right)
				& = 0, \\
				\partial_t B + \rot \tilde E  & = 0,
			\end{aligned}
		\end{cases}
	\end{equation*}
	with the constraints
	\begin{equation*}
		\begin{aligned}
			\DIV u^\pm & =0, & \rho^\pm+\theta^\pm & = 0,\\
			\DIV \tilde E & = 0, & \DIV B & = 0, \\
			\ROT B & = \left[\frac{\beta}{\eps}\right] \left(u^+-u^-\right), & E & = 0.
		\end{aligned}
	\end{equation*}
	
\end{enumerate}

On the whole, we conclude that, in the case of very weak interspecies collisions $\delta=O(\eps)$, the parameters $\alpha$, $\beta$ and $\gamma$ determine the asymptotics of the two species Vlasov-Maxwell-Boltzmann system \eqref{scaled VMB two species} exactly as they do determine the asymptotics of the one species Vlasov-Maxwell-Boltzmann system treated in Section \ref{formal one}. More precisely, the limiting two fluid macroscopic systems we obtain here can always be interpreted as two systems for one species --~one for cations and one for anions~-- coupled through their mean field interaction with the electromagnetic field $(E,B)$ and, whenever $\delta\sim\eps$, by an interspecies exchange of momentum and energy expressed by the linear terms $\frac1\sigma\left(u^+-u^-\right)$ and $\frac1\lambda\left(\theta^+-\theta^-\right)$. Therefore, the different asymptotic regimes for two species are also described by Figure \ref{figure 1} on page \pageref{figure 1}.

Thus, when $\delta\sim\eps$, up to multiplicative constants, we reach the following asymptotic systems of equations~:
\begin{enumerate}
	
	\item If $\alpha=o(\eps)$ and $\beta=o(\eps)$, we obtain the two fluid incompressible Navier-Stokes-Fourier system~:
	\begin{equation*}
		\begin{cases}
			\begin{aligned}
				\d_t u^\pm +
				u^\pm\cdot\nabla_x u^\pm - \mu\Delta_x u^\pm \pm \frac 1\sigma\left(u^+-u^-\right)
				& = -\nabla_x p^\pm , & \Div u^\pm & =0, \\
				\d_t \theta^\pm
				+
				u^\pm\cdot\nabla_x\theta^\pm - \kappa \Delta_x\theta^\pm \pm \frac 1\lambda\left(\theta^+-\theta^-\right)
				& = 0, & \rho^\pm+\theta^\pm & =0.
			\end{aligned}
		\end{cases}
	\end{equation*}
	This system satisfies the following formal energy conservation laws~:
	\begin{equation*}
		\begin{aligned}
			\frac 12\frac{d}{dt}\left(\left\|u^+\right\|_{L^2_x}^2 + \left\|u^-\right\|_{L^2_x}^2 \right) + \mu \left( \left\|\nabla_x u^+\right\|_{L^2_x}^2 + \left\|\nabla_x u^-\right\|_{L^2_x}^2 \right) + \frac 1\sigma \left\|u^+-u^-\right\|_{L^2_x}^2 & = 0, \\
			\frac 12\frac{d}{dt}\left(\left\|\theta^+\right\|_{L^2_x}^2 + \left\|\theta^-\right\|_{L^2_x}^2 \right) + \kappa \left( \left\|\nabla_x \theta^+\right\|_{L^2_x}^2 + \left\|\nabla_x \theta^-\right\|_{L^2_x}^2 \right) + \frac 1\lambda \left\|\theta^+-\theta^-\right\|_{L^2_x}^2 & = 0.
		\end{aligned}
	\end{equation*}
	
	\item If $\alpha=\eps^2$ and $\gamma=1$, we obtain the two fluid incompressible Navier-Stokes-Fourier-Maxwell system~:
	\begin{equation}\label{two fluid}
		\begin{cases}
			\begin{aligned}
				\d_t u^\pm +
				u^\pm\cdot\nabla_x u^\pm - \mu\Delta_x u^\pm \hspace{11mm} & && \\
				\pm \frac 1\sigma\left(u^+-u^-\right)
				& = -\nabla_x p^\pm
				\pm E \pm u^\pm \wedge B , & \Div u^\pm & = 0,\\
				\d_t \theta^\pm
				+
				u^\pm\cdot\nabla_x\theta^\pm - \kappa \Delta_x\theta^\pm \hspace{11mm} & && \\
				\pm \frac 1\lambda\left(\theta^+-\theta^-\right)
				& = 0, & \rho^\pm +\theta^\pm  & = 0, \\
				\d_t E - \ROT B &= -  \left(u^+-u^-\right), & \Div E & = 0,
				\\
				\d_t B + \ROT E & = 0, & \Div B & = 0.
			\end{aligned}
		\end{cases}
	\end{equation}
	This system satisfies the following formal energy conservation laws~:
	\begin{equation*}
		\begin{aligned}
			\frac 12\frac{d}{dt}\left(\left\|u^+\right\|_{L^2_x}^2 + \left\|u^-\right\|_{L^2_x}^2
			+ \left\|E\right\|_{L^2_x}^2 + \left\|B\right\|_{L^2_x}^2 \right) \hspace{50mm} & \\
			+ \mu \left( \left\|\nabla_x u^+\right\|_{L^2_x}^2 + \left\|\nabla_x u^-\right\|_{L^2_x}^2 \right) + \frac 1\sigma \left\|u^+-u^-\right\|_{L^2_x}^2 & = 0, \\
			\frac 12\frac{d}{dt}\left(\left\|\theta^+\right\|_{L^2_x}^2 + \left\|\theta^-\right\|_{L^2_x}^2 \right) + \kappa \left( \left\|\nabla_x \theta^+\right\|_{L^2_x}^2 + \left\|\nabla_x \theta^-\right\|_{L^2_x}^2 \right) + \frac 1\lambda \left\|\theta^+-\theta^-\right\|_{L^2_x}^2 & = 0.
		\end{aligned}
	\end{equation*}
	
	\item If $\alpha=o(\eps)$, $\beta=\eps$ and $\gamma=o(1)$, we obtain the two-fluid incompressible quasi-static Navier-Stokes-Fourier-Maxwell system~:
	\begin{equation}\label{TFIQSNSFM}
		\begin{cases}
			\begin{aligned}
				\d_t u^\pm +
				u^\pm\cdot\nabla_x u^\pm - \mu\Delta_x u^\pm \hspace{11mm} & && \\
				\pm \frac 1\sigma\left(u^+-u^-\right)
				& = -\nabla_x p^\pm \pm E \pm u^\pm \wedge B , & \Div u^\pm & = 0,\\
				\d_t \theta^\pm
				+
				u^\pm\cdot\nabla_x\theta^\pm - \kappa \Delta_x\theta^\pm \hspace{11mm} & && \\
				\pm \frac 1\lambda\left(\theta^+-\theta^-\right)
				& = 0, & \rho^\pm+\theta^\pm & =0 , \\
				\rot B & = u^+-u^-, & \Div E & = 0 , \\
				\partial_t B + \rot E  & = 0, & \Div B & = 0.
			\end{aligned}
		\end{cases}
	\end{equation}
	This system satisfies the following formal energy conservation laws~:
	\begin{equation*}
		\begin{aligned}
			\frac 12\frac{d}{dt}\left(\left\|u^+\right\|_{L^2_x}^2 + \left\|u^-\right\|_{L^2_x}^2 + \left\|B\right\|_{L^2_x}^2 \right) \hspace{60mm} & \\
			+ \mu \left( \left\|\nabla_x u^+\right\|_{L^2_x}^2 + \left\|\nabla_x u^-\right\|_{L^2_x}^2 \right) + \frac 1\sigma \left\|u^+-u^-\right\|_{L^2_x}^2 & = 0, \\
			\frac 12\frac{d}{dt}\left(\left\|\theta^+\right\|_{L^2_x}^2 + \left\|\theta^-\right\|_{L^2_x}^2 \right) + \kappa \left( \left\|\nabla_x \theta^+\right\|_{L^2_x}^2 + \left\|\nabla_x \theta^-\right\|_{L^2_x}^2 \right) + \frac 1\lambda \left\|\theta^+-\theta^-\right\|_{L^2_x}^2 & = 0.
		\end{aligned}
	\end{equation*}
	Here, the electric field is defined indirectly as a mere distribution, through Faraday's equation, by
	\begin{equation*}
		E=-\partial_t A,
	\end{equation*}
	where $B=\rot A$ and $\Div A=0$.
	
	\item If $\alpha=\eps$ and $\gamma=\eps$, we obtain the two fluid incompressible quasi-static Navier-Stokes-Fourier-Maxwell-Poisson system~:
	\begin{equation}\label{TFIQSNSFMP}
		\begin{cases}
			\begin{aligned}
				\d_t u^\pm +
				u^\pm\cdot\nabla_x u^\pm - \mu\Delta_x u^\pm & \\
				\pm \frac 1\sigma\left(u^+-u^-\right)
				& = -\nabla_x p^\pm \pm E +
				\rho^\pm \nabla_x\theta^\pm \pm u^\pm \wedge B , \\
				\Div u^\pm & = 0 , \\
				\d_t \left(\frac32\theta^\pm-\rho^\pm\right)
				+
				u^\pm\cdot\nabla_x\left(\frac32\theta^\pm-\rho^\pm\right)
				& - \frac 52 \kappa \Delta_x\theta^\pm \\
				\pm \frac5{2\lambda} \left(\theta^+-\theta^-\right)
				& = 0, \\
				\Delta_x\left(\rho^\pm+\theta^\pm\right) & =\pm\left(\rho^+-\rho^-\right), \\
				\ROT B & = u^+-u^-, \\
				\partial_t B + \rot E  & = 0, \\
				 \Div E & = \rho^+-\rho^- , \\
				 \Div B & = 0.
			\end{aligned}
		\end{cases}
	\end{equation}
	This system satisfies the following formal energy conservation law~:
	\begin{equation*}
		\begin{aligned}
			\frac 12\frac{d}{dt}\bigg(\left\|\rho^+\right\|_{L^2_x}^2 & + \left\|\rho^-\right\|_{L^2_x}^2
			+ \left\|u^+\right\|_{L^2_x}^2 + \left\|u^-\right\|_{L^2_x}^2 \\
			& + \frac 32 \left\|\theta^+\right\|_{L^2_x}^2 + \frac 32 \left\|\theta^-\right\|_{L^2_x}^2
			+ \left\|P^\perp E\right\|_{L^2_x}^2 + \left\|B\right\|_{L^2_x}^2 \bigg) & \\
			& + \mu \left( \left\|\nabla_x u^+\right\|_{L^2_x}^2 + \left\|\nabla_x u^-\right\|_{L^2_x}^2 \right)
			+ \frac52 \kappa \left( \left\|\nabla_x \theta^+\right\|_{L^2_x}^2 + \left\|\nabla_x \theta^-\right\|_{L^2_x}^2 \right) \\
			& + \frac 1\sigma \left\|u^+-u^-\right\|_{L^2_x}^2
			+ \frac 5{2\lambda} \left\|\theta^+-\theta^-\right\|_{L^2_x}^2
			= 0.
		\end{aligned}
	\end{equation*}
	Here, the solenoidal component of the electric field is defined indirectly as a mere distribution, through Faraday's equation, by
	\begin{equation*}
		PE=-\partial_t A,
	\end{equation*}
	where $B=\rot A$ and $\Div A=0$, while its irrotational component is determined, through Gauss' law, by
	\begin{equation*}
		P^\perp E = \pm\nabla_x\left(\rho^\pm+\theta^\pm\right).
	\end{equation*}
	Notice that the equations in this system are all coupled.
	
	\item If $\alpha=\eps$ and $\gamma=o(\eps)$, we obtain the two fluid incompressible Navier-Stokes-Fourier-Poisson system~:
	\begin{equation}\label{TFINSFP}
		\begin{cases}
			\begin{aligned}
				\d_t u^\pm +
				u^\pm\cdot\nabla_x u^\pm - \mu\Delta_x u^\pm
				\pm \frac 1\sigma\left(u^+-u^-\right)
				& = -\nabla_x p^\pm +
				\rho^\pm \nabla_x\theta^\pm , \\ \Div u^\pm & = 0,\\
				\d_t \left(\frac32\theta^\pm-\rho^\pm\right)
				+
				u^\pm\cdot\nabla_x\left(\frac32\theta^\pm-\rho^\pm\right)
				- \frac 52 & \kappa \Delta_x\theta^\pm \\
				\pm \frac5{2\lambda} \left(\theta^+-\theta^-\right)
				& = 0, \\
				\Delta_x\left(\rho^\pm+\theta^\pm\right) & =\pm\left(\rho^+-\rho^-\right).
			\end{aligned}
		\end{cases}
	\end{equation}
	This system satisfies the following formal energy conservation law~:
	\begin{equation*}
		\begin{aligned}
			\frac 12\frac{d}{dt}\bigg(\left\|\rho^+\right\|_{L^2_x}^2 & + \left\|\rho^-\right\|_{L^2_x}^2
			+ \left\|u^+\right\|_{L^2_x}^2 + \left\|u^-\right\|_{L^2_x}^2 \\
			& + \frac 32 \left\|\theta^+\right\|_{L^2_x}^2 + \frac 32 \left\|\theta^-\right\|_{L^2_x}^2
			+ \left\|\nabla_x\left(\rho^\pm+\theta^\pm\right)\right\|_{L^2_x}^2 \bigg) & \\
			& + \mu \left( \left\|\nabla_x u^+\right\|_{L^2_x}^2 + \left\|\nabla_x u^-\right\|_{L^2_x}^2 \right)
			+ \frac52 \kappa \left( \left\|\nabla_x \theta^+\right\|_{L^2_x}^2 + \left\|\nabla_x \theta^-\right\|_{L^2_x}^2 \right) \\
			& + \frac 1\sigma \left\|u^+-u^-\right\|_{L^2_x}^2
			+ \frac 5{2\lambda} \left\|\theta^+-\theta^-\right\|_{L^2_x}^2
			= 0.
		\end{aligned}
	\end{equation*}
	Physically, in this system, the fluid is subject to a self-induced static electric field $E$ determined by
	\begin{equation*}
		\rot E = 0,\qquad \Div E=\rho^+-\rho^-,
	\end{equation*}
	hence
	\begin{equation*}
		E=\pm\nabla_x\left(\rho^\pm+\theta^\pm\right).
	\end{equation*}
	Notice that the equations in this system are all coupled.
	
\end{enumerate}
When $\delta=o(\eps)$, one obtains the corresponding asymptotic systems by simply discarding the linear terms $\pm \frac 1\sigma\left(u^+-u^-\right)$ and $\pm \frac 1\lambda\left(\theta^+-\theta^-\right)$ in the preceding systems.

The above interpretation of two fluid systems as a coupling of one fluid systems will no longer hold for the more singular case of weak and strong interactions, i.e.\ when $\frac\delta\eps$ is unbounded, which we treat next.

\subsection{Macroscopic hydrodynamic constraints}\label{constraint-section 2}

Let us focus now on the analysis of the weak and strong interspecies collisional interactions. Contrary to the one species case, in the two species case, when $\frac\delta\eps$ is unbounded, the acoustic waves are always decoupled from the electromagnetic waves, which we treat below in Section \ref{ohm section}. We deal now with the acoustic waves.

At leading order, the system \eqref{moment-eps two species} describes the propagation of acoustic $(\rho_\eps,u_\eps,\sqrt{\frac32}\theta_\eps)$ waves~:
\begin{equation}\label{wave two species}
	\d_t
	\begin{pmatrix}
		\rho_\eps\\ u_\eps\\ \sqrt{3\over 2}\theta_\eps
	\end{pmatrix}
	+W_\eps
	\begin{pmatrix}
		\rho_\eps\\ u_\eps\\ \sqrt{3\over 2}\theta_\eps
	\end{pmatrix}
	=O(1),
\end{equation}
where the wave operator $W_\eps$, containing the singular terms from \eqref{moment-eps two species}, is defined by
\begin{equation}\label{wave operator 1 two species}
	W_\eps=
	\begin{pmatrix}
		0 &\frac1\eps \DIV&0\\
		\frac1\eps \nabla_x&0&\frac 1\eps\sqrt{2\over 3} \nabla_x\\
		0 & \frac 1\eps\sqrt{2\over 3} \DIV &0
	\end{pmatrix}.
\end{equation}

The wave operator $W_\eps$ is antisymmetric (with respect to the $L^2(dx)$ inner-product) and, therefore, can only have purely imaginary eigenvalues. The semi-group generated by this operator may thus produce fast time oscillations. Consequently, averaging over fast time oscillations as $\epsilon\to 0$, we get the macroscopic constraints
\begin{equation}\label{constraint1 two species}
	\DIV u =0, \qquad \rho+\theta =0,
\end{equation}
respectively referred to as incompressibility and Boussinesq relations. The exact nature of time oscillations produced by the system \eqref{wave two species}, in the limit $\eps\to 0$, will be rigorously discussed, with greater detail, later on in Chapter \ref{oscillations}.

\subsection{Hydrodynamic evolution equations}\label{evolution two species}

The previous step shows that, since $W_\eps$ is singular, the asymptotic dynamics of $\left(\rho_\eps,u_\eps,\sqrt{\frac 32}\theta_\eps\right)$ becomes constrained to  the kernel $\Ker W_\eps$ as $\eps\to 0$. Moreover, since $W_\eps$ is antisymmetric, its range is necessarily orthogonal to its kernel. Therefore, in order to get the asymptotic evolution equations for $\left(\rho,u,\sqrt{\frac 32}\theta\right)$, it is natural to project the system \eqref{wave two species} onto $\Ker W_\eps$, which will rid us of all the singular terms in \eqref{wave two species} and allow us to pass to the limit. In other words, we will obtain the limiting dynamics of the system \eqref{wave two species} by testing it against functions in $\Ker W_\eps$.

The kernel of $W_\eps$, defined in \eqref{wave operator 1 two species}, is obviously determined by all $\left(\rho_\eps^0,u_\eps^0,\sqrt{\frac 32}\theta_\eps^0\right)$ which satisfy
\begin{equation*}
		\DIV u_\eps^0 = 0 \qquad\text{and}\qquad \rho_\eps^0+\theta_\eps^0=0.
\end{equation*}
It is then readily seen that its orthogonal complement $\Ker W_\eps^\perp$ is determined by all $\left(\tilde \rho_\eps,\tilde u_\eps,\sqrt{\frac 32}\tilde \theta_\eps\right)$ such that
\begin{equation*}
		P \tilde u_\eps = 0 \qquad\text{and}\qquad \frac 32\tilde \theta_\eps-\tilde \rho_\eps=0.
\end{equation*}
Hence, projecting the system \eqref{moment-eps two species} onto $\Ker W_\eps^\perp$ yields
\begin{equation}\label{moment-eps2 two species}
	\begin{cases}
		\begin{aligned}
			\d_t Pu_\eps
			+\frac1{\eps} P \DIV \int_{\mathbb{R}^3} \mathcal{L}\left(\frac{g_\eps^++g_\eps^-}{2}\right) \tilde\phi M dv
			& =
			P\left(\frac\alpha{2\eps} n_\eps E_\eps +\frac\beta{2\delta} j_\eps \wedge B_\eps\right), \\
			\d_t \left(\frac32 \theta_\eps-\rho_\eps\right)
			+\frac1{\eps} \DIV \int_{\mathbb{R}^3}\mathcal{L}\left(\frac{g_\eps^++g_\eps^-}{2}\right)\tilde\psi M dv
			& =
			\frac{\alpha}{2\delta} j_\eps \cdot E_\eps.
		\end{aligned}
	\end{cases}
\end{equation}

\bigskip

There only remains to evaluate the flux terms $\frac1\eps \int_{\mathbb{R}^3} \mathcal{L}\left(\frac{g_\eps^++g_\eps^-}{2}\right) \tilde\phi M dv$ and $\frac1\eps \int_{\mathbb{R}^3}\mathcal{L}\left(\frac{g_\eps^++g_\eps^-}{2}\right) \tilde\psi M dv$ in \eqref{moment-eps2 two species}. Just as for one species in Section \ref{evolution}, following \cite{BGL2, bardos3}, this is done by employing \eqref{boltz-lin two species} to evaluate that
\begin{equation*}
	\begin{aligned}
		\frac {1+\delta^2}\eps \mathcal{L}\left(g_\eps^++g_\eps^-\right)
		& =
		\left[\cQ(g_\eps^+,g_\eps^+) + \cQ(g_\eps^-,g_\eps^-)\right]
		+ \delta^2\left[ \cQ\left(g_\eps^+,g_\eps^-\right) + \cQ\left(g_\eps^-,g_\eps^+\right) \right] \\
		& -v \cdot \nabla_x \left(g_\eps^++g_\eps^-\right)
		+O(\eps),
	\end{aligned}
\end{equation*}
% \begin{equation*}
% 	\begin{aligned}
% 		& \eps\d_t \left(g_\eps^++g_\eps^-\right) + v \cdot \nabla_x \left(g_\eps^++g_\eps^-\right)
% 		+ (\alpha E_\eps+\beta v\wedge B_\eps) \cdot \nabla_v \left(g_\eps^+-g_\eps^-\right)
% 		- \alpha E_\eps \cdot v \left(g_\eps^+-g_\eps^-\right) \\
% 		& = -\frac {1+\delta^2}\eps \mathcal{L}\left(g_\eps^++g_\eps^-\right)
% 		+
% 		\left[\cQ(g_\eps^+,g_\eps^+) + \cQ(g_\eps^-,g_\eps^-)\right]
% 		+ \delta^2\left[ \cQ\left(g_\eps^+,g_\eps^-\right) + \cQ\left(g_\eps^-,g_\eps^+\right) \right] ,
% 	\end{aligned}
% \end{equation*}
which yields formally in the limit, by virtue of the infinitesimal Maxwellian form \eqref{infinitesimal maxwellian 1},
\begin{equation*}
	\begin{aligned}
		\lim_{\eps\rightarrow 0}\frac {1}\eps \mathcal{L}\left(\frac{g_\eps^++g_\eps^-}{2}\right)
		& =
		\frac 1{2(1+[\delta]^2)}\left[\cQ(g^+,g^+) + \cQ(g^-,g^-)\right] \\
		& + \frac{[\delta]^2}{2(1+[\delta]^2)}\left[ \cQ\left(g^+,g^-\right) + \cQ\left(g^-,g^+\right) \right] \\
		& - \frac 1{1+[\delta]^2}v \cdot \nabla_x \left(\frac{g^++g^-}{2}\right) \\
		& =
		\frac 1{4}\left[\cL\left({g^+}^2\right) + \cL\left({g^-}^2\right)\right]
		- \frac 1{1+[\delta]^2}v \cdot \nabla_x \left(\frac{g^++g^-}{2}\right) \\
		& =
		\frac 12u^t\mathcal{L}(\phi)u
		+\theta u\cdot \mathcal{L}(\psi)
		+\frac 12\theta^2\mathcal{L}\left(\frac{|v|^4}{4}\right)\\
		& -\frac 1{1+[\delta]^2}\Div\left((\rho+\theta)v+\frac{|v|^2}{3}u+\phi u+\theta\psi\right) \\
		& =
		\frac 12u^t\mathcal{L}(\phi)u
		+\theta u\cdot \mathcal{L}(\psi)
		+\frac 12\theta^2\mathcal{L}\left(\frac{|v|^4}{4}\right)\\
		& -\frac 1{1+[\delta]^2}\Div\left(\phi u+\theta\psi\right),
	\end{aligned}
\end{equation*}
where we have used, in the last line, that $\DIV u =0$ and $\nabla_x (\rho+\theta) = 0$.

Next, we use that $\tilde\phi$ and $\tilde\psi$ have similar symmetry properties as $\phi$ and $\psi$, thanks to the rotational invariance of $\mathcal{L}$. More precisely, following \cite{desvillettes}, it can be shown (see also \cite[Section 2.2.3]{golse0}) that there exist two scalar valued functions $\alpha,\beta:[0,\infty)\to\mathbb{R}$ such that
\begin{equation*}
	\tilde\phi(v) = \alpha\left(\left|v\right|\right)\phi(v)
	\qquad\text{and}\qquad
	\tilde\psi(v) = \beta\left(\left|v\right|\right)\psi(v),
\end{equation*}
which implies (see \cite[Lemma 4.4]{BGL2}) that
\begin{equation*}% \label{delta identities 2}
	\begin{aligned}
		\int_{\mathbb{R}^3}\phi_{ij} \tilde\phi_{kl} Mdv &
		=(1+[\delta]^2)\mu \left(\delta_{ik}\delta_{jl}+\delta_{il}\delta_{jk} - \frac 23 \delta_{ij}\delta_{kl} \right),
		\\
		\int_{\mathbb{R}^3}\psi_{i} \tilde\psi_{j} Mdv & =(1+[\delta]^2) \frac 52 \kappa \delta_{ij},
	\end{aligned}
\end{equation*}
where
\begin{equation}\label{mu kappa 2}
	\mu = \frac{1}{10(1+[\delta]^2)}\int_{\mathbb{R}^3}\phi : \tilde\phi M dv
	\qquad\text{and}\qquad
	\kappa = \frac 2{15(1+[\delta]^2)} \int_{\mathbb{R}^3}\psi\cdot\tilde\psi M dv.
\end{equation}

Hence, we conclude through tedious but straightforward calculations that
\begin{equation*}
	\begin{aligned}
		\lim_{\eps\to 0} \frac1\eps \int_{\mathbb{R}^3} \mathcal{L}\left(\frac{g_\eps^++g_\eps^-}{2}\right) \tilde\phi M dv 
		& =
		\int_{\mathbb{R}^3} \frac 12 \left(u^t\phi u\right)\phi Mdv \\
		& - \frac 1{1+[\delta]^2} \int_{\mathbb{R}^3} \Div\left(\phi u\right)\tilde\phi Mdv
		\\
		& =
		u\otimes u -\frac{|u|^2}{3}\operatorname{Id}
		- \mu \left(\nabla_x u+\nabla_x^t u\right),
		\\
		\lim_{\eps\to 0} \frac1\eps \int_{\mathbb{R}^3}\mathcal{L}\left(\frac{g_\eps^++g_\eps^-}{2}\right) \tilde\psi M dv
		& =
		\int_{\mathbb{R}^3} \theta u\cdot \psi \psi M dv \\
		& - \frac 1{1+[\delta]^2} \int_{\mathbb{R}^3} \Div\left(\theta\psi\right)\tilde\psi M dv\\
		& =
		\frac 52 \theta u
		- \frac 52 \kappa \nabla_x\theta.
	\end{aligned}
\end{equation*}
We finally identify the advection and diffusion terms
\begin{equation*}
	\begin{aligned}
		\lim_{\eps\to 0} \frac1\eps P\Div \int_{\mathbb{R}^3} \mathcal{L}\left(\frac{g_\eps^++g_\eps^-}{2}\right) \tilde\phi M dv
		& =
		P\left(u\cdot\nabla_x u\right) - \mu\Delta_x u,
		\\
		\lim_{\eps\to 0} \frac1\eps \Div\int_{\mathbb{R}^3}\mathcal{L}\left(\frac{g_\eps^++g_\eps^-}{2}\right) \tilde\psi M dv
		& =
		\frac 52 u\cdot\nabla_x\theta - \frac 52 \kappa \Delta_x\theta .
	\end{aligned}
\end{equation*}

\bigskip

On the whole, provided nonlinear terms remain stable in the limiting process, we obtain the following asymptotic system~:
\begin{equation}\label{asymptotic hydro system}
	\begin{cases}
		\begin{aligned}
			\d_t u +
			u\cdot\nabla_x u - \mu\Delta_x u
			& = -\nabla_x p
			+\frac 12\left[\frac\alpha{\eps}\right] n E +\frac 12\left[\frac\beta{\delta}\right] j \wedge B, \\
			\d_t \theta
			+
			u\cdot\nabla_x\theta - \kappa \Delta_x\theta
			& = 0,
			% \frac 15\left[\frac{\alpha}{\delta}\right] j \cdot E ,
		\end{aligned}
	\end{cases}
\end{equation}
with the constraints from \eqref{constraint1 two species}
\begin{equation}\label{asymptotic hydro constraints}
	\DIV u =0, \qquad \rho+\theta =0.
\end{equation}
Unfortunately, as will be discussed later on in Section \ref{stability existence 2}, the rigorous weak stability of the nonlinear terms $n_\eps E_\eps \rightharpoonup nE$ and $j_\eps\wedge B_\eps \rightharpoonup j\wedge B$ remains unclear, in general. This will be, in fact, one of the main reasons for the breakdown of the weak compactness method in the most singular cases of hydrodynamic limits of the two species Vlasov-Maxwell-Boltzmann system \eqref{scaled VMB two species}, which will lead us to develop new relative entropy methods and consider dissipative solutions (see Section \ref{laure-diogo} and Chapter \ref{entropy method}).

There only remains now to formally establish the asymptotic system for the electrodynamic variables $(n,j,w)$ and the electromagnetic field $(E,B)$, which we do next.

\subsection{Macroscopic electrodynamic constraints and evolution}\label{ohm section}

The constraint equations for the electrodynamic variables will be obtained from the analysis of the difference of both components of \eqref{boltz-lin two species}~:
\begin{equation}\label{difference ohm}
	\begin{aligned}
		\frac\eps\delta\d_t \left(g_\eps^+-g_\eps^-\right) & + \frac 1\delta v \cdot \nabla_x \left(g_\eps^+-g_\eps^-\right) \\
		& + \left(\frac\alpha\delta E_\eps+\frac\beta\delta v\wedge B_\eps\right) \cdot \nabla_v \left(g_\eps^++g_\eps^-\right)
		- {\alpha\over \delta\eps}  E_\eps \cdot v \left(2+\eps\left(g_\eps^++g_\eps^-\right)\right) \\
		% & = -
		% \mathcal{L}\left(h_\eps^+-h_\eps^-\right)
		% -\delta^2\mathcal{L}\left(h_\eps^+-h_\eps^- , h_\eps^--h_\eps^+\right)
		% +
		% \delta\left[\cQ(g_\eps^+,g_\eps^+) - \cQ(g_\eps^-,g_\eps^-)\right]
		% + \delta^3\left[ \cQ\left(g_\eps^+,g_\eps^-\right) - \cQ\left(g_\eps^-,g_\eps^+\right) \right] \\
		& = -\frac 1{\delta^2}
		\mathcal{L}\left(h_\eps^+-h_\eps^-\right)
		-\mathcal{L}\left(h_\eps^+-h_\eps^- , h_\eps^--h_\eps^+\right) \\
		& +\frac{1-\delta^2}{2\delta}\cQ\left(g_\eps^++g_\eps^-,n_\eps\right)
		+\frac{1+\delta^2}{2\delta}\cQ\left(n_\eps,g_\eps^++g_\eps^-\right) \\
		& + \eps \frac{1-\delta^2}{2\delta^2}\cQ\left(g_\eps^++g_\eps^-,h_\eps^+-h_\eps^-\right)
		+\eps \frac{1+\delta^2}{2\delta^2}\cQ\left(h_\eps^+-h_\eps^-,g_\eps^++g_\eps^-\right).
	\end{aligned}
\end{equation}
However, the analysis in the case $\delta\sim 1$ will slightly differ from the case $\delta=o(1)$, with $\frac\delta\eps$ unbounded.

We begin with the case $\delta\sim 1$ of strong interspecies interactions. First, integrating \eqref{difference ohm} in $Mdv$ and letting $\eps\rightarrow 0$ easily yields the continuity equation
\begin{equation}\label{asymptotic continuity}
	\partial_t n+\frac 1{[\delta]}\Div j=0.
\end{equation}
Moreover, the above equation \eqref{difference ohm} contains no singular term in this situation. Therefore, letting $\eps\rightarrow 0$ yields, employing \eqref{infinitesimal maxwellian 1},
\begin{equation*}
	\begin{aligned}
		\frac 1{[\delta]} v \cdot \nabla_x n
		& - 2\left[\frac\beta\delta\right]\left( u\wedge B\right) \cdot v
		- 2\left[{\alpha\over \delta\eps}\right]  E \cdot v \\
		& = -\frac 1{[\delta]^2}
		\mathcal{L}\left(h^+-h^-\right)
		-\mathcal{L}\left(h^+-h^- , h^--h^+\right) \\
		& +[\delta]n\cL\left(u\cdot v +\theta\frac{|v|^2}{2}, -u\cdot v -\theta\frac{|v|^2}{2}\right).
	\end{aligned}
\end{equation*}
Further defining the linear operator
\begin{equation}\label{L frak def}
	\mathfrak{L}g = \mathcal{L}(g,-g),
\end{equation}
we have
\begin{equation}\label{difference ohm 4}
	\begin{aligned}
		\left(\frac 1{[\delta]} \nabla_x n
		- 2\left[\frac\beta\delta\right]\left( u\wedge B\right)
		- 2\left[{\alpha\over \delta\eps}\right]  E \right)\cdot v & \\
		= [\delta]n\mathfrak{L}\left(u\cdot v +\theta\frac{|v|^2}{2}\right)
		- \frac 1{[\delta]^2} \mathcal{L} & \left(h^+-h^-\right) -\mathfrak{L}\left(h^+-h^-\right).
	\end{aligned}
\end{equation}

Now, it can be shown that, in general, the linear operator $\frac 1{[\delta]^2}\mathcal{L}+\mathfrak{L}$ is self-adjoint and Fredholm of index zero on $L^2(Mdv)$ (or a variant of it depending on the cross-section, see Propositions \ref{hilbert-prop} and \ref{hilbert-prop 3}). Therefore, its range is exactly the orthogonal complement of its kernel, which is composed of all constant functions (see Proposition \ref{coercivity 3}). It follows that $\Phi(v)=v\in L^2(Mdv)$ and $\Psi(v)=\frac{|v|^2}{2}-\frac 32\in L^2(Mdv)$ belong to the range of $\frac 1{[\delta]^2}\mathcal{L}+\mathfrak{L}$ and, thus, that there are inverses $\tilde\Phi\in L^2(Mdv)$ and $\tilde\Psi\in L^2(Mdv)$ such that
\begin{equation}\label{phi-psi-def inverses two species}
	\Phi =\frac 1{[\delta]^2}\mathcal{L} \tilde \Phi + \mathfrak{L} \tilde \Phi \qquad\text{and}\qquad \Psi= \frac 1{[\delta]^2}\mathcal{L} \tilde \Psi + \mathfrak{L} \tilde \Psi,
\end{equation}
which can be uniquely determined by the fact that they are orthogonal to the kernel of $\frac 1{[\delta]^2}\mathcal{L}+\mathfrak{L}$ (i.e.\ to constant functions). Furthermore, it can be shown that $\tilde\Phi$ and $\tilde\Psi$ have similar symmetry properties as $\Phi$ and $\Psi$, thanks to the rotational invariance of $\mathcal{L}$ and $\mathfrak{L}$. More precisely, employing methods from \cite{desvillettes} (see also \cite[Section 2.2.3]{golse0}), one verifies that there exist two scalar valued functions $\alpha,\beta:[0,\infty)\to\mathbb{R}$ such that
\begin{equation*}
	\tilde\Phi(v) = \alpha\left(\left|v\right|\right)\Phi(v)
	\qquad\text{and}\qquad
	\tilde\Psi(v) = \beta\left(\left|v\right|\right)\Psi(v),
\end{equation*}
which implies that
\begin{equation}\label{sigma 2}
	\int_{\mathbb{R}^3}\Phi_i\tilde\Phi_j Mdv =\frac 12 \sigma\delta_{ij},
\end{equation}
where
\begin{equation}\label{sigma 3}
	\sigma =\frac 23 \int_{\mathbb{R}^3}\Phi\cdot \tilde\Phi Mdv
\end{equation}
defines the electrical conductivity $\sigma>0$. For completeness, we also define the energy conductivity $\lambda>0$ by
\begin{equation}\label{lambda 2}
	\lambda = \int_{\mathbb{R}^3}\Psi \tilde\Psi Mdv.
\end{equation}

Then, multiplying \eqref{difference ohm 4} by $\tilde\Phi$, integrating in $Mdv$, exploiting the self-adjointness of $\mathcal{L}$ and $\mathfrak{L}$ and the limiting representation \eqref{h limit 2} of $\Pi h^\pm$, yields Ohm's law
\begin{equation}\label{ohm}
		j - [\delta]nu = \sigma\left( - \frac 1{2[\delta]} \nabla_x n
		+ \left[{\alpha\over \delta\eps}\right]  E
		+ \left[\frac\beta\delta\right]u\wedge B \right).
\end{equation}
Similarly, multiplying \eqref{difference ohm 4} by $\tilde\Psi$, we obtain the energy equivalence relation
\begin{equation}\label{ohm energy equivalence}
	w = [\delta] n\theta.
\end{equation}

Finally, in the case $\delta\sim 1$, the whole asymptotic system \eqref{asymptotic hydro system}-\eqref{asymptotic hydro constraints}-\eqref{asymptotic continuity}-\eqref{ohm}-\eqref{ohm energy equivalence} will be fully determined when considering the coupling with the limiting Maxwell's equations from \eqref{Maxwell-eps two species}~:
\begin{equation*}% \label{asymptotic Maxwell-eps two species}
	\begin{cases}
		\begin{aligned}
			[\gamma]\d_t E - \ROT B &= - \left[\frac{\beta}{\delta}\right]j,
			\\
			[\gamma]\d_t B + \ROT E & = 0,
			\\
			\DIV E &=\left[\frac{\alpha}{\epsilon}\right] n,
			\\
			\DIV B &=0.
		\end{aligned}
	\end{cases}
\end{equation*}

Let us focus now on the case $\delta=o(1)$, which turns out to be more complicated than the case $\delta\sim 1$, for it contains yet another singular limit, as we are about to see. Indeed, the most singular term in \eqref{difference ohm} being $-\frac 1{\delta^2} \mathcal{L}\left(h_\eps^+-h_\eps^-\right)$, we begin by projecting \eqref{difference ohm} onto the collision invariants in order to eliminate this singular term. This yields (this system may also be deduced directly from \eqref{two fluid system eps} by considering the difference of the equations for cations and anions)
\begin{equation}\label{difference ohm 2}
	\begin{cases}
		\begin{aligned}
			\d_t n_\eps + \frac 1{\delta} \DIV j_\eps & = 0,\\
			\frac{\eps^2}{\delta^2}\d_t j_\eps + \frac1\delta \nabla_x \left( n_\eps+\frac\eps\delta w_\eps \right) & -\frac{2\alpha}{\delta\eps}E_\eps
			+\int_{\mathbb{R}^3}\mathcal{L}\left(h_\eps^+-h_\eps^-,h_\eps^--h_\eps^+\right)vMdv \\
			& =
			\left(\frac{2\alpha}{\delta} \rho_\eps E_\eps +\frac{2\beta}{\delta} u_\eps \wedge B_\eps\right) \\
			& -\frac\eps\delta \DIV \int_{\mathbb{R}^3} \frac1\eps\left(g_\eps^+ - \Pi g_\eps^+ - g_\eps^- +\Pi g_\eps^-\right) \phi M dv \\
			& +\delta\int_{\mathbb{R}^3}\left[\mathcal{Q}\left(g_\eps^+,g_\eps^-\right)
			-\mathcal{Q}\left(g_\eps^-,g_\eps^+\right) \right] vMdv , \\
			\frac {3\eps^2}{2\delta^2} \d_t w_\eps + \frac \eps{\delta^2} \DIV j_\eps & +
			\int_{\mathbb{R}^3}\mathcal{L}\left(h_\eps^+-h_\eps^-,h_\eps^--h_\eps^+\right)\frac{|v|^2}{2}Mdv \\
			& =
			\frac{2\alpha}{\delta} u_\eps \cdot E_\eps \\
			& - \frac{\eps}{\delta} \DIV \int_{\mathbb{R}^3}\frac1\eps\left(g_\eps^+ - \Pi g_\eps^+ - g_\eps^- +\Pi g_\eps^-\right) \psi M dv \\
			& +\delta\int_{\mathbb{R}^3}\left[\mathcal{Q}\left(g_\eps^+,g_\eps^-\right)
			-\mathcal{Q}\left(g_\eps^-,g_\eps^+\right) \right] \frac{|v|^2}{2}Mdv.
		\end{aligned}
	\end{cases}
\end{equation}
Then, since $\Pi h^\pm_\eps = \pm \frac 12\left(j_\eps\cdot v + w_\eps\left(\frac{|v|^2}{2}-\frac 32\right)\right)$, straightforward computations based on symmetry of integrands show that
\begin{equation*}
	\begin{aligned}
		\int_{\mathbb{R}^3}\mathcal{L}\left(\Pi\left(h_\eps^+-h_\eps^-\right),\Pi\left(h_\eps^--h_\eps^+\right)\right)vMdv
		&=  \frac 2\sigma j_\eps \\
		\int_{\mathbb{R}^3}\mathcal{L}\left(\Pi\left(h_\eps^+-h_\eps^-\right),\Pi\left(h_\eps^--h_\eps^+\right)\right)\frac{|v|^2}{2}Mdv
		&= \frac 1\lambda w_\eps \\
	\end{aligned}
\end{equation*}
where the electrical conductivity $\sigma>0$ and the energy conductivity $\lambda>0$ are constants defined by
\begin{equation}\label{sigma}
	\begin{aligned}
		\frac 1\sigma & =\frac 16 \int_{\mathbb{R}^3}v\cdot\mathcal{L}\left(v , - v \right)Mdv \\
		& =\frac 16\int_{\mathbb{R}^3\times\mathbb{R}^3\times\mathbb{S}^2}\left|v-v'\right|^2 b(v-v_*,\sigma) MM_* dvdv_*d\sigma \\
		& =\frac 1{6}\int_{\mathbb{R}^3\times\mathbb{R}^3}\left|v-v_*\right|^2 
		m(v-v_*) MM_* dvdv_*,
	\end{aligned}
\end{equation}
and
\begin{equation}\label{lambda}
	\begin{aligned}
		\frac 1\lambda & =\frac 14 \int_{\mathbb{R}^3}|v|^2\mathcal{L}\left(|v|^2,
		-|v|^2\right)Mdv \\
		& =\frac 14 \int_{\mathbb{R}^3\times\mathbb{R}^3\times\mathbb{S}^2}\left(|v|^2-|v'|^2\right)^2 b(v-v_*,\sigma) MM_* dvdv_*d\sigma \\
		& =\frac 1{4}\int_{\mathbb{R}^3\times\mathbb{R}^3}\left(|v|^2-|v_*|^2\right)^2 
		m(v-v_*) MM_* dvdv_*,
	\end{aligned}
\end{equation}
where the cross-section for momentum and energy transfer $m(v-v_*)$ is defined in Proposition \ref{cross section transfer}. It follows that the system \eqref{difference ohm 2} may be rewritten as
\begin{equation}\label{difference ohm 3}
	\begin{cases}
		\begin{aligned}
			\d_t n_\eps & + \frac 1{\delta} \DIV j_\eps = 0,\\
			\frac 2\sigma j_\eps & + \frac1\delta \nabla_x \left( n_\eps+\frac\eps\delta w_\eps \right)
			-\frac{2\alpha}{\delta\eps}E_\eps
			\\
			& =
			\left(\frac{2\alpha}{\delta} \rho_\eps E_\eps +\frac{2\beta}{\delta} u_\eps \wedge B_\eps\right)
			-\frac{\eps^2}{\delta^2}\d_t j_\eps \\
			& -\frac\eps\delta \DIV \int_{\mathbb{R}^3} \frac1\eps\left(g_\eps^+ - \Pi g_\eps^+ - g_\eps^- +\Pi g_\eps^-\right) \phi M dv \\
			& - \int_{\mathbb{R}^3}\mathcal{L}\left(h_\eps^+ - \Pi h_\eps^+ - h_\eps^- +\Pi h_\eps^-,h_\eps^- - \Pi h_\eps^- - h_\eps^+ +\Pi h_\eps^+\right)vMdv \\
			& +\delta\int_{\mathbb{R}^3}\left[\mathcal{Q}\left(g_\eps^+,g_\eps^-\right)
			-\mathcal{Q}\left(g_\eps^-,g_\eps^+\right) \right] vMdv , \\
			\frac 1\lambda w_\eps & + \frac{\eps}{\delta^2}\Div j_\eps \\
			& =
			\frac{2\alpha}{\delta} u_\eps \cdot E_\eps - \frac {3\eps^2}{2\delta^2} \d_t w_\eps \\
			& - \frac{\eps}{\delta} \DIV \int_{\mathbb{R}^3}\frac1\eps\left(g_\eps^+ - \Pi g_\eps^+ - g_\eps^- +\Pi g_\eps^-\right) \psi M dv \\
			& - \int_{\mathbb{R}^3}\mathcal{L}\left(h_\eps^+ - \Pi h_\eps^+ - h_\eps^- +\Pi h_\eps^-,h_\eps^- - \Pi h_\eps^- - h_\eps^+ +\Pi h_\eps^+\right)\frac{|v|^2}{2}Mdv \\
			& +\delta\int_{\mathbb{R}^3}\left[\mathcal{Q}\left(g_\eps^+,g_\eps^-\right)
			-\mathcal{Q}\left(g_\eps^-,g_\eps^+\right) \right] \frac{|v|^2}{2}Mdv.
		\end{aligned}
	\end{cases}
\end{equation}

In fact, the system \eqref{difference ohm 3} coupled with Maxwell's equations \eqref{Maxwell-eps two species} still contains a singular perturbation, which will be treated much like the singular perturbation of the one species case in Sections \ref{macro constraint one species} and \ref{evolution}. Thus, by virtue of \eqref{h equilibrium}, it is readily seen that the system \eqref{difference ohm 3} may be further simplified to
\begin{equation}\label{difference ohm 5}
	\begin{cases}
		\begin{aligned}
			\DIV j_\eps & = O(\delta),\\
			\nabla_x \left( n_\eps+\frac\eps\delta w_\eps \right)
			& = \frac{2\alpha}{\eps}E_\eps + O(\delta), \\
			P j_\eps
			& =
			\sigma P\left(\frac{\alpha}{\delta\eps}E_\eps + \frac{\beta}{\delta} u_\eps \wedge B_\eps\right) + o(1), \\
			w_\eps & = o(1).
		\end{aligned}
	\end{cases}
\end{equation}
We discuss now the limit $\eps\rightarrow 0$ of the coupled system \eqref{Maxwell-eps two species}-\eqref{difference ohm 5}.
\begin{enumerate}
	\item When $\gamma\sim 1$ (so that $\alpha=O(\delta\eps)$), letting $\eps\rightarrow 0$ in \eqref{Maxwell-eps two species}-\eqref{difference ohm 5}, we obtain
	\begin{equation*}% \label{asymptotic EM constraints}
		\begin{cases}
			\begin{aligned}
				[\gamma]\d_t E - \ROT B &= - \left[\frac{\beta}{\delta}\right]j, & \DIV E &=0,
				\\
				[\gamma]\d_t B + \ROT E & = 0, & \DIV B &=0,
				\\
				j
				& =
				\sigma\left(-\nabla_x \bar p +\left[\frac{\alpha}{\delta\eps}\right]E + \left[\frac{\beta}{\delta}\right] u \wedge B\right),
				& \DIV j & = 0,
				\\
				n & = 0, & w & = 0,
			\end{aligned}
		\end{cases}
	\end{equation*}
	where $\bar p$ is an electrodynamic pressure.
	
	\item When $\gamma=o(1)$ and $\alpha=O(\delta\eps)$, letting $\eps\rightarrow 0$ in \eqref{Maxwell-eps two species}-\eqref{difference ohm 5}, we obtain
	\begin{equation*}
		\begin{cases}
			\begin{aligned}
				E & = 0, && \\
				\ROT B &= \left[\frac{\beta}{\delta}\right]j, & \DIV B &=0,
				\\
				j
				& =
				\sigma\left(-\nabla_x \bar p + \left[\frac{\beta}{\delta}\right] u \wedge B\right),
				& \DIV j & = 0,
				\\
				n & = 0, & w & = 0,
			\end{aligned}
		\end{cases}
	\end{equation*}
	where $\bar p$ is an electrodynamic pressure.
	
	\item When $\gamma=o(1)$ and $\frac\alpha{\delta\eps}$ is unbounded, we need to further use Faraday's equation from \eqref{Maxwell-eps two species}, as in Section \ref{evolution}, to write that
	\begin{equation*}
		\frac{\beta}{\delta}\d_t A_\eps + \frac{\alpha}{\delta\eps} PE_\eps = 0,
	\end{equation*}
	where $B_\eps=\rot A_\eps$ and $\Div A_\eps =0$. Thus, letting $\eps\rightarrow 0$ in \eqref{Maxwell-eps two species}-\eqref{difference ohm 5}, we obtain
	\begin{equation*}
		\begin{cases}
			\begin{aligned}
				\DIV E &=\left[\frac\alpha\eps\right]n, & \ROT E & = 0,
				\\
				\ROT B &= \left[\frac{\beta}{\delta}\right]j, & \DIV B &=0,
				\\
				\ROT A & =B, & \Div A & = 0,
				\\
				j
				& =
				\sigma\left(-\nabla_x \bar p - \left[\frac{\beta}{\delta}\right]\partial_t A + \left[\frac{\beta}{\delta}\right] u \wedge B\right),
				& \DIV j & = 0,
				\\
				\nabla_x n & = 2\left[\frac\alpha\eps\right]E, & w & = 0,
			\end{aligned}
		\end{cases}
	\end{equation*}
	where $A_\eps\rightarrow A$ and $\bar p$ is an electrodynamic pressure. Note that $\Delta_x n=2\left[\frac\alpha\eps\right]^2 n$, so that necessarily $n=0$ and $E=0$. The above system can be rewritten more explicitly by defining the adjusted electric field $\tilde E = - \partial_t A$. It then holds that
	\begin{equation*}
		\begin{cases}
			\begin{aligned}
				\ROT B &= \left[\frac{\beta}{\delta}\right]j, & \DIV B &=0,
				\\
				\partial_t B + \ROT \tilde E & =0, & \Div \tilde E &= 0 ,
				\\
				j
				& =
				\sigma\left(-\nabla_x \bar p + \left[\frac{\beta}{\delta}\right]\tilde E + \left[\frac{\beta}{\delta}\right] u \wedge B\right),
				& \DIV j & = 0,
				\\
				n & = 0, & w & = 0.
			\end{aligned}
		\end{cases}
	\end{equation*}
	% Notice, finally, that if further $\alpha=o(\eps)$, then the above system is greatly simplified and becomes
	% \begin{equation*}
	% 	\begin{cases}
	% 		\begin{aligned}
	% 			\ROT B &= \left[\frac{\beta}{\delta}\right]j, & \DIV B &=0,
	% 			\\
	% 			\partial_t B + \ROT \tilde E & =0, & \DIV \tilde E &=0,
	% 			\\
	% 			j
	% 			& =
	% 			\sigma\left(-\nabla_x \bar p + \left[\frac{\beta}{\delta}\right]\tilde E + \left[\frac{\beta}{\delta}\right] u \wedge B\right),
	% 			& \DIV j & = 0,
	% 			\\
	% 			n & = 0, & w & = 0.
	% 		\end{aligned}
	% 	\end{cases}
	% \end{equation*}
	
\end{enumerate}

\subsection{Summary}

\begin{figure}[!ht]
	\centering
	\includegraphics[scale=0.54]{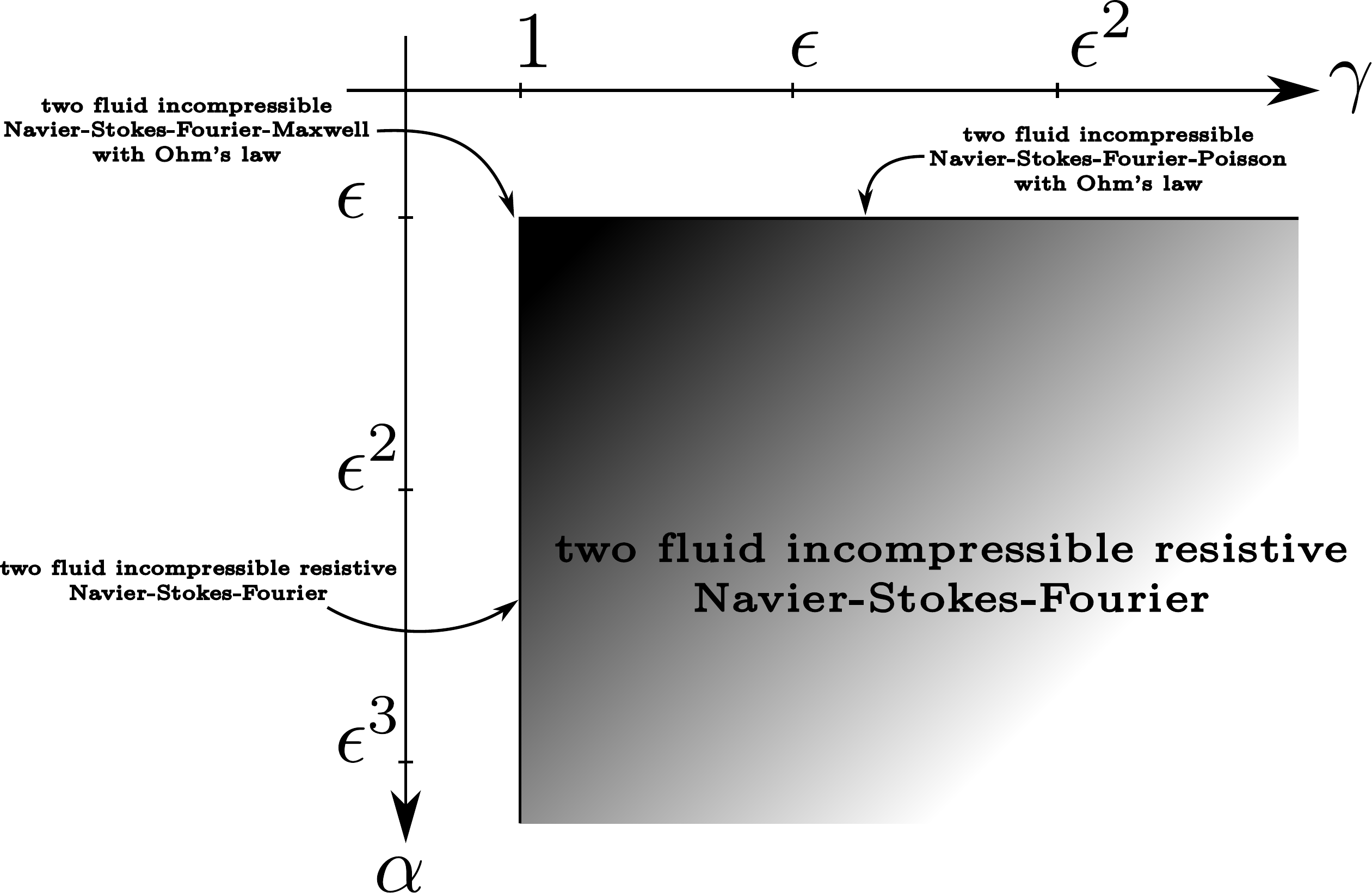}
	\caption{Asymptotic regimes of the two species Vlasov-Maxwell-Boltzmann system \eqref{scaled VMB two species} for strong interspecies interactions.}
	\label{figure 2}
\end{figure}

\begin{figure}[!ht]
	\centering
	\includegraphics[scale=0.54]{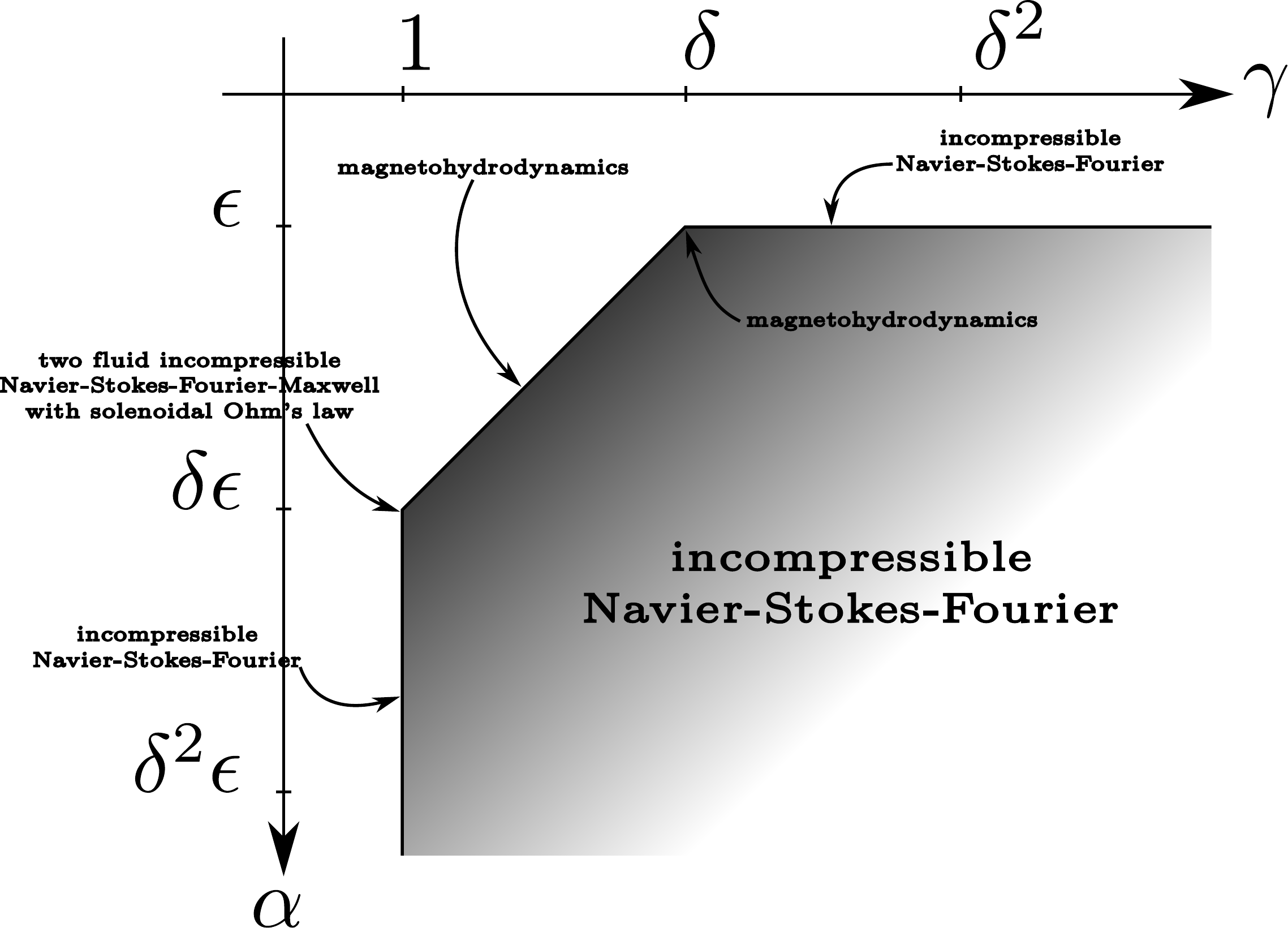}
	\caption{Asymptotic regimes of the two species Vlasov-Maxwell-Boltzmann system \eqref{scaled VMB two species} for weak interspecies interactions.}
	\label{figure 3}
\end{figure}

At last, we see that the asymptotics of the two species Vlasov-Maxwell-Boltzmann system \eqref{scaled VMB two species} can be depicted in terms of the limits of the following parameters~:
\begin{itemize}
	\item the strength of the electric induction $\alpha$,
	\item the strength of the magnetic induction $\beta=\frac{\alpha\gamma}{\eps}$,
	\item the ratio of the bulk velocity to the speed of light $\gamma$,
	\item the strength of the interspecies collisional interactions $\delta$.
\end{itemize}
The case of very weak interspecies collisions has already been discussed in Section \ref{very weak section} and is analogous to the one species case. Regarding the weak and strong interspecies collisions, Figures \ref{figure 2} and \ref{figure 3} summarize the different asymptotic regimes, on a logarithmic scale, of the two species Vlasov-Maxwell-Boltzmann system \eqref{scaled VMB two species}.

Thus, up to multiplicative constants, in the case of strong interactions $\delta=1$, we reach the following asymptotic systems of equations~:
\begin{enumerate}
	
	\item If $\alpha=o(\eps)$, we obtain the two fluid incompressible resistive Navier-Stokes-Fourier system~:
	\begin{equation*}
		\begin{cases}
			\begin{aligned}
				\d_t u +
				u\cdot\nabla_x u - \mu\Delta_x u
				& = -\nabla_x p , & \Div u & =0, && \\
				\d_t \theta
				+
				u\cdot\nabla_x\theta - \kappa \Delta_x\theta
				& = 0, & \rho+\theta & =0, && \\
				\d_t n
				+
				u\cdot\nabla_xn - \frac\sigma 2 \Delta_x n
				& = 0, & j-nu & =-\frac\sigma 2\nabla_x n, & w & = n\theta.
			\end{aligned}
		\end{cases}
	\end{equation*}
	This system satisfies the following formal energy conservation laws~:
	\begin{equation*}
		\begin{aligned}
			\frac 12\frac{d}{dt}\left\|u\right\|_{L^2_x}^2 + \mu \left\|\nabla_x u\right\|_{L^2_x}^2 & = 0, \\
			\frac 12\frac{d}{dt}\left\|\theta\right\|_{L^2_x}^2 + \kappa \left\|\nabla_x \theta\right\|_{L^2_x}^2 & = 0, \\
			\frac 12\frac{d}{dt}\left\|n\right\|_{L^2_x}^2 + \frac\sigma 2 \left\|\nabla_x n\right\|_{L^2_x}^2 & = 0.
		\end{aligned}
	\end{equation*}
	
	\item If $\alpha=\eps$ and $\gamma=1$, we obtain the two fluid incompressible Navier-Stokes-Fourier-Maxwell system with Ohm's law~:
	\begin{equation}\label{TFINSFMO}
		\begin{cases}
			\begin{aligned}
				\d_t u +
				u\cdot\nabla_x u - \mu\Delta_x u
				& = -\nabla_x p+
				\frac 12 \left(nE + j \wedge B\right) , & \Div u & = 0,\\
				\d_t \theta
				+
				u\cdot\nabla_x\theta - \kappa \Delta_x\theta
				& = 0, & \rho+\theta & = 0, \\
				\d_t E - \ROT B &= -  j, & \Div E & = n,
				\\
				\d_t B + \ROT E & = 0, & \Div B & = 0, \\
				j-nu & = \sigma\left(-\frac 12 \nabla_x n + E + u\wedge B\right), & w & =n\theta.
			\end{aligned}
		\end{cases}
	\end{equation}
	This system satisfies the following formal energy conservation laws (see Proposition \ref{energy estimate 2} for an explicit computation of the energy)~:
	\begin{equation*}
		\begin{aligned}
			\frac 14\frac{d}{dt}\left(2\left\|u\right\|_{L^2_x}^2 + \frac 12\left\|n\right\|_{L^2_x}^2 + \left\|E\right\|_{L^2_x}^2 + \left\|B\right\|_{L^2_x}^2 \right)
			+ \mu \left\|\nabla_x u\right\|_{L^2_x}^2 + \frac 1{2\sigma} \left\|j-nu\right\|_{L^2_x}^2 & = 0, \\
			\frac 12\frac{d}{dt}\left\|\theta\right\|_{L^2_x}^2 + \kappa \left\|\nabla_x \theta\right\|_{L^2_x}^2 & = 0.
		\end{aligned}
	\end{equation*}
	
	\item If $\alpha=\eps$ and $\gamma=o(1)$, we obtain the two fluid incompressible Navier-Stokes-Fourier-Poisson system with Ohm's law~:
	% \begin{equation*}
	% 	\begin{cases}
	% 		\begin{aligned}
	% 			\d_t u +
	% 			u\cdot\nabla_x u - \mu\Delta_x u
	% 			& = -\nabla_x p+
	% 			\frac 12 nE , & \Div u & = 0,\\
	% 			\d_t \theta
	% 			+
	% 			u\cdot\nabla_x\theta - \kappa \Delta_x\theta
	% 			& = 0, & \rho+\theta & = 0, \\
	% 			\d_t n
	% 			+
	% 			u\cdot\nabla_xn - \frac\sigma 2 \Delta_x n + \sigma n
	% 			& = 0, & \Div E & = n,
	% 			\\
	% 			\ROT E & = 0, & B & = 0, \\
	% 			j-nu & = \sigma\left(-\frac 12 \nabla_x n + E \right), & w & =n\theta.
	% 		\end{aligned}
	% 	\end{cases}
	% \end{equation*}
	\begin{equation*}
		\begin{cases}
			\begin{aligned}
				\d_t u +
				u\cdot\nabla_x u - \mu\Delta_x u
				& = -\nabla_x p+
				\frac 12 n\nabla_x\phi , & \Div u & = 0,\\
				\d_t \theta
				+
				u\cdot\nabla_x\theta - \kappa \Delta_x\theta
				& = 0, & \rho+\theta & = 0, \\
				\d_t n
				+
				u\cdot\nabla_xn - \frac\sigma 2 \Delta_x n + \sigma n
				& = 0, & \Delta_x \phi & = n,
				\\
				j-nu & = \sigma \nabla_x\left(\phi -\frac 12 n\right), & w & =n\theta.
			\end{aligned}
		\end{cases}
	\end{equation*}
	This system satisfies the following formal energy conservation laws~:
	\begin{equation*}
		\begin{aligned}
			\frac 14\frac{d}{dt}\left(2\left\|u\right\|_{L^2_x}^2 + \frac 12\left\|n\right\|_{L^2_x}^2 + \left\|\nabla_x\phi\right\|_{L^2_x}^2 \right)
			+ \mu \left\|\nabla_x u\right\|_{L^2_x}^2 + \frac 1{2\sigma} \left\|j-nu\right\|_{L^2_x}^2 & = 0, \\
			\frac 12\frac{d}{dt}\left\|n\right\|_{L^2_x}^2 + \frac\sigma 2 \left\|\nabla_x n\right\|_{L^2_x}^2
			+ \sigma \left\| n \right\|_{L^2_x}^2 & = 0,
			\\
			\frac 12\frac{d}{dt}\left\|\theta\right\|_{L^2_x}^2 + \kappa \left\|\nabla_x \theta\right\|_{L^2_x}^2 & = 0.
		\end{aligned}
	\end{equation*}
	
\end{enumerate}

Finally, up to multiplicative constants, in the case of weak interactions $\delta=o(1)$, we reach the following asymptotic systems of equations~:
\begin{enumerate}
	
	\item If $\beta=o(\delta)$, we obtain the incompressible Navier-Stokes-Fourier system~:
	\begin{equation*}
		\begin{cases}
			\begin{aligned}
				\d_t u +
				u\cdot\nabla_x u - \mu\Delta_x u
				& = -\nabla_x p , & \Div u & =0, && \\
				\d_t \theta
				+
				u\cdot\nabla_x\theta - \kappa \Delta_x\theta
				& = 0, & \rho+\theta & =0, && \\
				n & = 0, & j & = 0, & w & = 0.
			\end{aligned}
		\end{cases}
	\end{equation*}
	This system satisfies the following formal energy conservation laws~:
	\begin{equation*}
		\begin{aligned}
			\frac 12\frac{d}{dt}\left\|u\right\|_{L^2_x}^2 + \mu \left\|\nabla_x u\right\|_{L^2_x}^2 & = 0, \\
			\frac 12\frac{d}{dt}\left\|\theta\right\|_{L^2_x}^2 + \kappa \left\|\nabla_x \theta\right\|_{L^2_x}^2 & = 0.
		\end{aligned}
	\end{equation*}
	
	\item If $\alpha=\delta\eps$ and $\gamma=1$, we obtain the two fluid incompressible Navier-Stokes-Fourier-Maxwell system with solenoidal Ohm's law~:
	\begin{equation}\label{TFINSFMSO}
		\begin{cases}
			\begin{aligned}
				\d_t u +
				u\cdot\nabla_x u - \mu\Delta_x u
				& = -\nabla_x p+
				\frac 12 j \wedge B , & \Div u & = 0,\\
				\d_t \theta
				+
				u\cdot\nabla_x\theta - \kappa \Delta_x\theta
				& = 0, & \rho+\theta & = 0, \\
				\d_t E - \ROT B &= -  j, & \Div E & = 0,
				\\
				\d_t B + \ROT E & = 0, & \Div B & = 0, \\
				j & = \sigma\left(- \nabla_x \bar p + E + u\wedge B\right), & \Div j & =0,\\
				n&=0, & w&=0.
			\end{aligned}
		\end{cases}
	\end{equation}
	This system satisfies the following formal energy conservation laws (see Proposition \ref{energy estimate 2} for an explicit computation of the energy)~:
	\begin{equation*}
		\begin{aligned}
			\frac 14\frac{d}{dt}\left(2\left\|u\right\|_{L^2_x}^2 + \left\|E\right\|_{L^2_x}^2 + \left\|B\right\|_{L^2_x}^2 \right)
			+ \mu \left\|\nabla_x u\right\|_{L^2_x}^2 + \frac 1{2\sigma} \left\|j\right\|_{L^2_x}^2 & = 0, \\
			\frac 12\frac{d}{dt}\left\|\theta\right\|_{L^2_x}^2 + \kappa \left\|\nabla_x \theta\right\|_{L^2_x}^2 & = 0.
		\end{aligned}
	\end{equation*}
	
	\item If $\beta=\delta$ and $\gamma=o(1)$, we obtain the two fluid incompressible quasi-static Navier-Stokes-Fourier-Maxwell system with solenoidal Ohm's law~:
	\begin{equation*}
		\begin{cases}
			\begin{aligned}
				\d_t u +
				u\cdot\nabla_x u - \mu\Delta_x u
				& = -\nabla_x p+
				\frac 12 j \wedge B , & \Div u & = 0,\\
				\d_t \theta
				+
				u\cdot\nabla_x\theta - \kappa \Delta_x\theta
				& = 0, & \rho+\theta & = 0, \\
				\ROT B &= j, & \Div E & = 0,
				\\
				\d_t B + \ROT E & = 0, & \Div B & = 0, \\
				j & = \sigma\left(-\nabla_x \bar p + E + u\wedge B\right), &  \Div j & = 0 , \\
				n & = 0, & w & = 0.
			\end{aligned}
		\end{cases}
	\end{equation*}
	This system satisfies the following formal energy conservation laws~:
	\begin{equation*}
		\begin{aligned}
			\frac 14\frac{d}{dt}\left(2\left\|u\right\|_{L^2_x}^2 + \left\|B\right\|_{L^2_x}^2 \right)
			+ \mu \left\|\nabla_x u\right\|_{L^2_x}^2 + \frac 1{2\sigma} \left\|j\right\|_{L^2_x}^2 & = 0, \\
			\frac 12\frac{d}{dt}\left\|\theta\right\|_{L^2_x}^2 + \kappa \left\|\nabla_x \theta\right\|_{L^2_x}^2 & = 0.
		\end{aligned}
	\end{equation*}
	Here, the electric field is defined indirectly as a mere distribution, through Faraday's equation, by
	\begin{equation*}
		E=-\partial_t A,
	\end{equation*}
	where $B=\rot A$ and $\Div A=0$. Note that the above system can be rewritten as
	\begin{equation*}
		\begin{cases}
			\begin{aligned}
				\d_t u +
				u\cdot\nabla_x u - \mu\Delta_x u
				& = -\nabla_x p+
				\frac 12 \rot B \wedge B , & \Div u & = 0, \\
				\d_t \theta
				+
				u\cdot\nabla_x\theta - \kappa \Delta_x\theta
				& = 0, & \rho+\theta & = 0, \\
				\d_t B + u\cdot\nabla_x B-\frac1\sigma\Delta_x B & = B\cdot\nabla_x u, & \Div B & = 0,
			\end{aligned}
		\end{cases}
	\end{equation*}
	which is nothing but the well-known magnetohydrodynamic system. The rigorous derivation of this system starting from other macroscopic systems such as \eqref{two fluid} and \eqref{TFINSFMSO} as been investigated in \cite{arsenio6}.

\end{enumerate}

\subsection{The two species Vlasov-Poisson-Boltzmann system}

The two species Vlasov-Poisson-Boltzmann system describes the evolution of a gas of two species of charged particles (cations and anions) subject to an auto-induced electrostatic force. This system is obtained formally from the two species Vlasov-Maxwell-Boltzmann system by letting the speed of light tend to infinity while all other parameters remain fixed. Accordingly, setting $\gamma = 0$ in \eqref{scaled VMB two species} yields the scaled Vlasov-Poisson-Boltzmann system~:
\begin{equation}\label{scaledVMP two species}
	\begin{cases}
		\begin{aligned}
			\eps \d_t f_\eps^\pm + v \cdot \nabla_x f_\eps^\pm \pm \alpha \nabla_x\phi_\eps \cdot \nabla_v f_\eps^\pm &= \frac 1\eps Q(f_\eps^\pm,f_\eps^\pm) + \frac{\delta^2}{\eps} Q(f_\eps^\pm,f_\eps^\mp),
			\\
			f_\eps^\pm & =M\left(1+\epsilon g_\eps^\pm\right),
			\\
			\Delta_x\phi_\eps&=\frac{\alpha}{\epsilon}\int_{\mathbb{R}^3} \left(g_\eps^+-g_\eps^-\right) M dv,
		\end{aligned}
	\end{cases}
\end{equation}
Here, the plasma is subject to a self-induced electrostatic field $E_\eps$ determined by
\begin{equation*}
	\rot E_\eps = 0,\qquad \Div E_\eps=\frac{\alpha}{\epsilon}\int_{\mathbb{R}^3} \left(g_\eps^+-g_\eps^-\right) M dv,
\end{equation*}
hence
\begin{equation*}
	E_\eps=\nabla_x\phi_\eps.
\end{equation*}
The above system is supplemented with some initial data satisfying
\begin{equation*}
	\frac1{\eps^2} H\left(f_\eps^{+\mathrm{in}}\right)
	+ \frac1{\eps^2} H\left(f_\eps^{-\mathrm{in}}\right)
	+ \frac1{2}\int_{\mathbb{R}^3} |E_\eps^{\rm in}|^2  dx < \infty.
\end{equation*}
In particular, solutions of \eqref{scaledVMP two species} satisfy the corresponding scaled entropy inequality, where $t>0$,
\begin{equation*}
	\begin{aligned}
	 	\frac1{\eps^2} H\left(f_\eps^{+}\right)
		+ \frac1{\eps^2} H\left(f_\eps^{-}\right)
		& + \frac 1{2} \int_{\mathbb{R}^3} |E_\eps|^2 dx \\
		& +\frac{1}{\epsilon^4}\int_0^t\int_{\mathbb{R}^3}\left(D\left(f_\eps^+\right)+D\left(f_\eps^-\right)
		+ \delta^2 D\left(f_\eps^+,f_\eps^-\right)\right)(s) dx ds
		\\
		& \leq
		\frac1{\eps^2} H\left(f_\eps^{+\mathrm{in}}\right)
		+ \frac1{\eps^2} H\left(f_\eps^{-\mathrm{in}}\right)
		+ \frac1{2}\int_{\mathbb{R}^3} |E_\eps^{\rm in}|^2 dx,
	\end{aligned}
\end{equation*}

Thus, the formal asymptotic analysis of \eqref{scaledVMP two species} is contained in our analysis of the two species Vlasov-Maxwell-Boltzmann system \eqref{scaled VMB two species}. Specifically, setting $\gamma=\beta=0$ in the limiting systems first obtained in Section \ref{very weak section}, for very weak interspecies collisions, we see that the two species Vlasov-Poisson-Boltzmann system \eqref{scaledVMP two species} converges, when $\alpha=o(\eps)$ and $\delta=O(\eps)$, towards the two fluid incompressible Navier-Stokes-Fourier system in a Boussinesq regime, with $E=0$~:
\begin{equation*}
	\begin{cases}
		\begin{aligned}
			\d_t u^\pm +
			u^\pm\cdot\nabla_x u^\pm - \mu\Delta_x u^\pm
			\pm\left[\frac\delta\eps\right]^2\frac 1\sigma\left(u^+-u^-\right)
			& = -\nabla_x p^\pm , & \Div u^\pm & =0\\
			\d_t \theta^\pm
			+
			u^\pm\cdot\nabla_x\theta^\pm - \kappa \Delta_x\theta^\pm
			\pm\left[\frac\delta\eps\right]^2\frac 1\kappa\left(\theta^+-\theta^-\right)
			& = 0, & \rho^\pm+\theta^\pm & =0.
		\end{aligned}
	\end{cases}
\end{equation*}
While, when $\left[\frac\alpha\eps\right]\neq 0$ and $\delta=O(\eps)$, we find the convergence towards the two fluid incompressible Navier-Stokes-Fourier-Poisson system~:
\begin{equation*}
	\begin{cases}
		\begin{aligned}
			\d_t u^\pm +
			u^\pm\cdot\nabla_x u^\pm - \mu\Delta_x u^\pm
			\pm \left[\frac\delta\eps\right]^2\frac 1\sigma\left(u^+-u^-\right)
			& = -\nabla_x p^\pm +
			\rho^\pm \nabla_x\theta^\pm , \\ \Div u^\pm & = 0,\\
			\d_t \left(\frac32\theta^\pm-\rho^\pm\right)
			+
			u^\pm\cdot\nabla_x\left(\frac32\theta^\pm-\rho^\pm\right)
			- \frac 52 & \kappa \Delta_x\theta^\pm \\
			\pm \frac5{2}\left[\frac\delta\eps\right]^2\frac1\lambda \left(\theta^+-\theta^-\right)
			& = 0, \\
			\Delta_x\left(\rho^\pm+\theta^\pm\right) & =\pm\left[\frac{\alpha}{\eps}\right]^2\left(\rho^+-\rho^-\right),
		\end{aligned}
	\end{cases}
\end{equation*}
where the electrostatic field is determined by $\left[\frac\alpha\eps\right]E=\pm\nabla_x\left(\rho^\pm+\theta^\pm\right)$.

Regarding weak interspecies interactions, setting $\gamma=\beta=0$ in the corresponding limiting systems obtained in Sections \ref{evolution two species} and \ref{ohm section}, we see that the two species Vlasov-Poisson-Boltzmann system \eqref{scaledVMP two species} always converges, when $\delta=o(1)$ and $\frac\delta\eps$ is unbounded, towards the incompressible Navier-Stokes-Fourier system in a Boussinesq regime, with $E=0$~:
\begin{equation*}
	\begin{cases}
		\begin{aligned}
			\d_t u +
			u\cdot\nabla_x u - \mu\Delta_x u
			& = -\nabla_x p , & \Div u & =0, && \\
			\d_t \theta
			+
			u\cdot\nabla_x\theta - \kappa \Delta_x\theta
			& = 0, & \rho+\theta & =0, && \\
			n & = 0, & j & = 0, & w & = 0.
		\end{aligned}
	\end{cases}
\end{equation*}

Finally, in the case of strong interspecies interactions, setting $\gamma=\beta=0$ in corresponding the limiting systems obtained in Sections \ref{evolution two species} and \ref{ohm section}, we see that the two species Vlasov-Poisson-Boltzmann system \eqref{scaledVMP two species} converges, when $\alpha=o(\eps)$ and $\delta = 1$, towards the two fluid incompressible resistive Navier-Stokes-Fourier system in a Boussinesq regime, with $E=0$~:
\begin{equation*}
	\begin{cases}
		\begin{aligned}
			\d_t u +
			u\cdot\nabla_x u - \mu\Delta_x u
			& = -\nabla_x p , & \Div u & =0, && \\
			\d_t \theta
			+
			u\cdot\nabla_x\theta - \kappa \Delta_x\theta
			& = 0, & \rho+\theta & =0, && \\
			\d_t n
			+
			u\cdot\nabla_xn - \frac\sigma 2 \Delta_x n
			& = 0, & j-nu & =-\frac\sigma 2\nabla_x n, & w & = n\theta.
		\end{aligned}
	\end{cases}
\end{equation*}
While, when $\left[\frac\alpha\eps\right]\neq 0$ and $\delta=1$, we find the convergence towards the two fluid incompressible Navier-Stokes-Fourier-Poisson system with Ohm's law~:
\begin{equation*}
	\begin{cases}
		\begin{aligned}
			\d_t u +
			u\cdot\nabla_x u - \mu\Delta_x u
			& = -\nabla_x p+
			\frac 12\left[\frac\alpha\eps\right] n\nabla_x\phi , & \Div u & = 0,\\
			\d_t \theta
			+
			u\cdot\nabla_x\theta - \kappa \Delta_x\theta
			& = 0, & \rho+\theta & = 0, \\
			\d_t n
			+
			u\cdot\nabla_xn - \frac\sigma 2 \Delta_x n + \sigma \left[\frac\alpha\eps\right]^2 n
			& = 0, & \Delta_x \phi & = \left[\frac\alpha\eps\right]n,
			\\
			j-nu & = \sigma \nabla_x\left(\left[\frac\alpha\eps\right]\phi -\frac 12 n\right), & w & =n\theta,
		\end{aligned}
	\end{cases}
\end{equation*}
where the electrostatic field is determined by $E=\nabla_x\phi$.

In fact, the two species Vlasov-Poisson-Boltzmann system is inherently simpler than the two species Vlasov-Maxwell-Boltzmann system, because it couples the Vlasov-Boltzmann equations with a simple elliptic equation, namely Poisson's equation, while the two species Vlasov-Maxwell-Boltzmann system couples the Vlasov-Boltzmann equations with an hyperbolic system, namely Maxwell's system of equations. Thus, the rigorous mathematical analysis on the two species Vlasov-Maxwell-Boltzmann system, presented in the remainder of this work, will also apply to the two species Vlasov-Poisson-Boltzmann system and, therefore, analog results will hold.

%% file: stability.tex
\chapter{Weak stability of the limiting macroscopic systems}\label{weak stability}

In the previous chapter, we have formally derived numerous viscous incompressible systems for plasmas starting from Vlasov-Maxwell-Boltzmann systems and we intend to provide, in the remainder of our work, justifications of these derivations. Nevertheless, prior to any rigorous proof of hydrodynamic limit, it is crucial to understand the well-posedness of the asymptotic macroscopic models and to study their stability properties.

Describing the Cauchy problem of each single macroscopic system from Chapter \ref{formal-chap} would be unreasonable. Rather, we are now going to focus on the following three systems found therein~:
\begin{itemize}
	\item the incompressible quasi-static Navier-Stokes-Fourier-Maxwell-Poisson system \eqref{NSFMP},
	\item the two-fluid incompressible Navier-Stokes-Fourier-Maxwell system with Ohm's law \eqref{TFINSFMO},
	\item the two-fluid incompressible Navier-Stokes-Fourier-Maxwell system with solenoidal Ohm's law \eqref{TFINSFMSO},
\end{itemize}
and establish the existence of weak or dissipative solutions to their respective initial value problems. In fact, these three systems are among the most singular ones found in Chapter \ref{formal-chap}. Thus, we hope the reader will find it clear that the existence of appropriate weak or dissipative solutions to the remaining macroscopic systems from Chapter \ref{formal-chap} will then follow from straightforward adjustments of the existence theories presented here.

In the remaining Parts \ref{part 2}, \ref{part 3} and \ref{part 4} of our work, we will also focus on the three aforementioned systems and give complete justifications of their derivation from hydrodynamic limits of Vlasov-Maxwell-Boltzmann systems.

\section[The incompressible quasi-static Navier-Stokes-Fourier-\ldots]{The incompressible quasi-static Navier-Stokes-Fourier-Maxwell-Poisson system}\label{stability existence 1}

We are first concerned here with the incompressible quasi-static Navier-Stokes-Fourier-Maxwell-Poisson system \eqref{NSFMP}, which we rewrite, for mere convenience~:
\begin{equation}\label{NSFMP 3}
	\begin{cases}
		\begin{aligned}
			\d_t u
			+
			u\cdot\nabla_x u - \mu\Delta_x u
			& = -\nabla_x p + E
			+ \rho \nabla_x\theta + u \wedge B , \hspace{-20mm}&& \\
			&& \Div u & = 0,\\
			\d_t \left(\frac32\theta-\rho\right)
			+
			u\cdot\nabla_x\left(\frac32\theta-\rho\right)
			- \frac 52 \kappa \Delta_x\theta
			& = 0,
			& \Delta_x(\rho+\theta) & =\rho, \\
			\ROT B & = u, & \Div E & = \rho , \\
			\partial_t B + \rot E  & = 0, & \Div B & = 0.
		\end{aligned}
	\end{cases}
\end{equation}
Although it looks more complicated because it involves more terms, the system \eqref{NSFMP 3} has the same structure as the incompressible Navier-Stokes equations~: it is indeed a system of parabolic equations, in which the nonlinear advection terms are well-defined by the energy estimate.

The following formal proposition shows how to compute the energy.

\begin{prop}\label{energy estimate}
	Let $\left(\rho,u,\theta,B\right)$ be a smooth solution to the incompressible quasi-static Navier-Stokes-Fourier-Maxwell-Poisson system \eqref{NSFMP 3}.
	
	Then, the following global energy inequality holds~:
	\begin{equation}\label{energy}
		\begin{aligned}
			& \frac 12\left(\left\|\rho(t)\right\|_{L^2_x}^2 + \left\|u(t)\right\|_{L^2_x}^2 + \frac 32 \left\|\theta(t)\right\|_{L^2_x}^2
			+ \left\|\nabla_x\left(\rho+\theta\right)(t)\right\|_{L^2_x}^2 + \left\|B(t)\right\|_{L^2_x}^2 \right) \\
			& \hspace{13mm} + \int_0^t \mu \left\|\nabla_x u(s)\right\|_{L^2_x}^2 + \frac 52 \kappa \left\|\nabla_x \theta(s)\right\|_{L^2_x}^2 ds \\
			& \hspace{13mm} \leq \frac 12\left(\left\|\rho^\mathrm{in}\right\|_{L^2_x}^2 + \left\|u^\mathrm{in}\right\|_{L^2_x}^2 + \frac 32 \left\|\theta^\mathrm{in}\right\|_{L^2_x}^2
			+ \left\|\nabla_x\left(\rho^\mathrm{in}+\theta^\mathrm{in}\right)\right\|_{L^2_x}^2 + \left\|B^\mathrm{in}\right\|_{L^2_x}^2 \right),
		\end{aligned}
	\end{equation}
	where $\nabla_x(\rho+\theta)=P^\perp E\neq E$.
	% \begin{equation}
	% \begin{aligned}
	%  \frac12&\| u(t)\|^2_{L^2} +\frac14 \vvvert \theta(t) \vvvert^2 +\frac12 \| B(t)\|^2_{L^2}  \\
	%  &+\mu \int_0^t \| \nabla u (s)\|_{L^2}^2 ds +\frac52 \kappa \int_0^t \| \nabla \theta (s)\|_{L^2}^2 ds
	%  &\leq   \frac12\| u^0\|^2_{L^2}+\frac14 \vvvert \theta^0 \vvvert^2 +\frac12 \| B^0\|^2_{L^2} 
	% \end{aligned}
	% \end{equation}
	% where, by definition,
	% $$ \vvvert \theta \vvvert^2 = \int \theta ( I-\Delta)^{-1} (3I-5\Delta) \theta dx\,.$$
\end{prop}

\begin{proof}
Multiplying the equation expressing the conservation of momentum in \eqref{NSFMP 3} by $u$ and integrating with respect to space variables, we get
\begin{equation*}
	\begin{aligned}
		\frac12 {d\over dt} \| u\|_{L^2_x}^2 + \mu \left\|\nabla_x u \right\|_{L^2_x}^2 &= \int_{\mathbb{R}^3} u\cdot E + \rho  u\cdot \nabla_x\theta dx\\
		&= \int_{\mathbb{R}^3} \ROT B\cdot E + \rho  u\cdot \nabla_x\theta dx\\
		& = -\frac12 {d\over dt} \| B\|_{L^2_x}^2 +\int_{\mathbb{R}^3} \rho u\cdot \nabla_x \theta dx .
	\end{aligned}
\end{equation*}
Then, multiplying the equation expressing the conservation of energy by $\theta$, we get similarly
\begin{equation*}
	\begin{aligned}
		\frac34 {d\over dt} \| \theta\|_{L^2_x}^2 + \frac 52 \kappa \left\|\nabla_x \theta \right\|_{L^2_x}^2
		&= \int_{\mathbb{R}^3} \theta\partial_t\rho + \theta  u\cdot \nabla_x\rho dx\\
		& = -\frac12 {d\over dt} \| \rho\|_{L^2_x}^2
		+ \int_{\mathbb{R}^3} \left(\rho+\theta\right)\partial_t\rho - \rho  u\cdot \nabla_x\theta dx \\
		& = -\frac12 {d\over dt} \| \rho\|_{L^2_x}^2
		-\frac12 {d\over dt} \| \nabla_x\left(\rho+\theta\right) \|_{L^2_x}^2
		- \int_{\mathbb{R}^3} \rho  u\cdot \nabla_x\theta dx.
	\end{aligned}
\end{equation*}
% $$
% \begin{aligned}
% \frac14 {d\over dt} \int \theta (3\theta - 2 (I-\Delta)^{-1} \Delta \theta) dx&= \int \Big(- \theta u\cdot \nabla ( \frac32 \theta- \rho) -\frac52 \kappa|\nabla \theta|^2 \Big) dx\\
% &= -\frac52 \kappa \| \nabla \theta\|_{L^2}^2 +\int \theta u \cdot \nabla \rho dx  \,.
% \end{aligned}
% $$
Summing the above identities, we obtain the expected global conservation of energy.
\end{proof}

% Note that, using the Fourier transform and Plancherel's formula, we see easily that the norm $\vvvert \cdot \vvvert$ is equivalent to the $L^2$-norm.
% Furthermore, because of the constraint $\rho-\Delta \rho = \Delta \theta$ and elliptic regularity estimates, $\rho$ has the same regularity as $\theta$.

Using the a priori estimates provided by the energy inequality \eqref{energy} and reproducing the arguments of Leray \cite{leray}, we can easily establish the global existence of weak solutions. Indeed, combining first the bound on $\nabla_x(\rho+\theta)$ with the additional spatial regularity on $u$ and $\theta$, coming from the dissipation terms in the energy inequality \eqref{energy}, we infer that all three terms $\rho$, $u$ and $\theta$ enjoy some spatial regularity. More precisely, they are all uniformly bounded in $L^2_\mathrm{loc}\left(dt;H^1(dx)\right)$. Furthermore, recalling that $PE=-\partial_t A$, where $B=\ROT A$ with $\DIV A = 0$, some temporal regularity on $u+A$ and $\frac 32\theta-\rho$ is clearly inherited from the evolution equations, which allows us to establish, invoking a classical compactness result by Aubin and Lions \cite{aubin, lions6} (see also \cite{simon} for a sharp compactness criterion), that $u+A$ and $\frac 32\theta-\rho$ are strongly relatively compact in all variables in $L^2_\mathrm{loc}\left(dtdx\right)$. Finally, noticing that one may express
\begin{equation*}
	\begin{aligned}
		\rho & =\frac{2\Delta_x}{3-5\Delta_x}\left(\frac 32\theta-\rho\right),
		&
		u & = \frac{-\Delta_x}{1-\Delta_x}\left(u+A\right),
		\\
		\theta & =\frac{2-2\Delta_x}{3-5\Delta_x}\left(\frac 32\theta-\rho\right),
		&
		B & = \frac{\ROT}{1-\Delta_x}\left(u+A\right),
	\end{aligned}
\end{equation*}
using Poisson's equations
\begin{equation*}
	-\Delta_x A = u, \qquad \Delta_x\left(\rho+\theta\right)=\rho,
\end{equation*}
we easily find that all four observables $\rho$, $u$, $\theta$ and $B$ belong to a compact subset of $L^2_\mathrm{loc}\left(dtdx\right)$.

The above compactness properties, which also hold for the similar systems \eqref{IQSNSFM}, \eqref{INSFP}, \eqref{TFIQSNSFM},\eqref{TFIQSNSFMP} and \eqref{TFINSFP}, allow us to prove the weak stability of the nonlinear terms in \eqref{NSFMP 3} and, therefore, to take weak limits in any suitable approximation scheme to establish the existence of weak solutions. Analogous existence results hold for systems \eqref{IQSNSFM}, \eqref{INSFP}, \eqref{TFIQSNSFM},\eqref{TFIQSNSFMP} and \eqref{TFINSFP}, as well.

Henceforth, we will utilize the prefixes $\textit{w-}$ or $\textit{w$^*$-}$ to express that a given space is endowed with its weak or weak-$*$ topology, respectively.

\begin{thm}
	Let $\left(\rho^\mathrm{in},u^\mathrm{in},\theta^\mathrm{in},B^\mathrm{in}\right) \in L^2\left(\mathbb{R}^3,dx\right)$ be such that
	\begin{equation*}
		\Div u^\mathrm{in} = 0, \qquad \Div B^\mathrm{in} = 0, \qquad \ROT B^\mathrm{in} = u^\mathrm{in},\qquad \Delta_x \left(\rho^\mathrm{in} + \theta^\mathrm{in}\right) = \rho^\mathrm{in}.
	\end{equation*}

	Then, there exists (at least) one global weak solution (in the sense of Leray \cite{leray})
	\begin{equation*}
			\begin{aligned}
				\left(\rho, u, \theta, B\right) & \in C\left([0,\infty) ; \textit{w-}L^2\left(\mathbb{R}^3,dx\right)\right)
				\cap L^\infty \left( [0,\infty), dt ; L^2\left(\mathbb{R}^3, dx\right)\right), \\
				\left(u,\theta\right) & \in L^2\left([0,\infty), dt ; \dot H^1\left(\mathbb{R}^3, dx\right)\right),
			\end{aligned}
	\end{equation*}
	to the incompressible quasi-static Navier-Stokes-Fourier-Maxwell-Poisson system \eqref{NSFMP 3}. Furthermore, it satisfies the energy inequality \eqref{energy}.
\end{thm}

As usual for such weak solutions, uniqueness is not known to hold. To prove that the system \eqref{NSFMP 3} is well-posed in the sense of Hadamard, we would have to deal with a stronger notion of solution. Note however that, by modulating the energy inequality, we can establish some weak-strong uniqueness principle, meaning that if a somewhat regular solution to \eqref{NSFMP 3} is known to exist, then any weak solution with matching initial data coincides with the smooth one as long as it exists. We refer to the next Section \ref{stability existence 2} for details on how to modulate the energy and, thus, establish such weak-strong uniqueness principles.

\section[The two-fluid incompressible Navier-Stokes-Fourier-Maxwell\ldots]{The two-fluid incompressible Navier-Stokes-Fourier-Maxwell system with (solenoidal) Ohm's law}\label{stability existence 2}

We focus now on the two-fluid incompressible Navier-Stokes-Fourier-Maxwell system with Ohm's law \eqref{TFINSFMO}~:
\begin{equation}\label{TFINSFMO 3}
	\begin{cases}
		\begin{aligned}
			\d_t u +
			u\cdot\nabla_x u - \mu\Delta_x u
			& = -\nabla_x p+
			\frac 12 \left(nE + j \wedge B\right) , & \Div u & = 0,\\
			\d_t \theta
			+
			u\cdot\nabla_x\theta - \kappa \Delta_x\theta
			& = 0, & \rho+\theta & = 0, \\
			\d_t E - \ROT B &= -  j, & \Div E & = n,
			\\
			\d_t B + \ROT E & = 0, & \Div B & = 0, \\
			j-nu & = \sigma\left(-\frac 12 \nabla_x n + E + u\wedge B\right), & w & =n\theta,
		\end{aligned}
	\end{cases}
\end{equation}
and on the two-fluid incompressible Navier-Stokes-Fourier-Maxwell system with solenoidal Ohm's law \eqref{TFINSFMSO}~:
\begin{equation}\label{TFINSFMSO 3}
	\begin{cases}
		\begin{aligned}
			\d_t u +
			u\cdot\nabla_x u - \mu\Delta_x u
			& = -\nabla_x p+
			\frac 12 j \wedge B , & \Div u & = 0,\\
			\d_t \theta
			+
			u\cdot\nabla_x\theta - \kappa \Delta_x\theta
			& = 0, & \rho+\theta & = 0, \\
			\d_t E - \ROT B &= -  j, & \Div E & = 0,
			\\
			\d_t B + \ROT E & = 0, & \Div B & = 0, \\
			j & = \sigma\left(- \nabla_x \bar p + E + u\wedge B\right), & \Div j & =0,\\
			n&=0, & w&=0.
		\end{aligned}
	\end{cases}
\end{equation}
The above models \eqref{TFINSFMO 3} and \eqref{TFINSFMSO 3} are not stable under weak convergence in the energy space and, thus, share more similarities with the three-dimensional incompressible Euler equations, as we are about to discuss.

To this end, note first that the advection-diffusion equation on $\theta$ is not really coupled with the other equations on $(u,n,j,E,B)$ in \eqref{TFINSFMO 3} and \eqref{TFINSFMSO 3}, and that it is linear provided the velocity field $u$ is given. It is therefore sufficient to focus on the reduced systems of equations
\begin{equation}\label{TFINSMO}
	\begin{cases}
		\begin{aligned}
			\d_t u +
			u\cdot\nabla_x u - \mu\Delta_x u
			& = -\nabla_x p+
			\frac 12 \left(nE + j \wedge B\right) , & \Div u & = 0,\\
			\d_t E - \ROT B &= -  j, & \Div E & = n,
			\\
			\d_t B + \ROT E & = 0, & \Div B & = 0, \\
			j-nu & = \sigma\left(-\frac 12 \nabla_x n + E + u\wedge B\right), &&
		\end{aligned}
	\end{cases}
\end{equation}
and
\begin{equation}\label{TFINSMSO}
	\begin{cases}
		\begin{aligned}
			\d_t u +
			u\cdot\nabla_x u - \mu\Delta_x u
			& = -\nabla_x p+
			\frac 12 j \wedge B , & \Div u & = 0,\\
			\d_t E - \ROT B &= -  j, & \Div E & = 0,
			\\
			\d_t B + \ROT E & = 0, & \Div B & = 0, \\
			j & = \sigma\left(- \nabla_x \bar p + E + u\wedge B\right), & \Div j & =0.
		\end{aligned}
	\end{cases}
\end{equation}

\begin{rem}
	The system \eqref{TFINSMSO} can be viewed as an asymptotic regime of system \eqref{TFINSMO}. Indeed, at least formally, it is obtained, as $\delta\rightarrow 0$, from the system
	\begin{equation*}
		\begin{cases}
			\begin{aligned}
				\d_t u +
				u\cdot\nabla_x u - \mu\Delta_x u
				& = -\nabla_x p+
				\frac 12 \left(\delta nE + j \wedge B\right) , & \Div u & = 0,\\
				\d_t E - \ROT B &= -  j, & \Div E & = \delta n,
				\\
				\d_t B + \ROT E & = 0, & \Div B & = 0, \\
				j-\delta nu & = \sigma\left(-\frac 1{2\delta} \nabla_x n + E + u\wedge B\right), &&
			\end{aligned}
		\end{cases}
	\end{equation*}
	which is consistent with the formal derivations from Section \ref{formal two}.
\end{rem}

A natural framework to study these equations (coming from physics) should be the energy space, i.e.\ the functional space defined by the (formal) energy conservation. We indeed expect solutions in this space to be global.

The following formal proposition shows how to compute the energy of the two-fluid incompressible Navier-Stokes-Maxwell system with Ohm's law \eqref{TFINSMO}, or with solenoidal Ohm's law \eqref{TFINSMSO}.

\begin{prop}\label{energy estimate 2}
Let $(u,E,B)$ be a smooth solution to the two-fluid incompressible Navier-Stokes-Maxwell system with Ohm's law \eqref{TFINSMO}, or with solenoidal Ohm's law \eqref{TFINSMSO}.

Then the following global conservation of energy holds~:
\begin{equation*}
	\CE(t) +\int_0^t\CD(s)ds= \CE(0),\qquad\text{for all }t>0,
\end{equation*}
where the energy $\CE$ and the energy dissipation $\CD$ are given by, for the system \eqref{TFINSMO},
\begin{equation*}
	\begin{aligned}
		\CE(t) & =
		\frac 12 \left\|u(t)\right\|_{L^2_x}^2 + \frac 18\left\|n(t)\right\|_{L^2_x}^2
		+ \frac 14 \left\|E(t)\right\|_{L^2_x}^2 + \frac 14 \left\|B(t)\right\|_{L^2_x}^2, \\
		\CD(t) & =
		\mu \left\|\nabla_x u(t)\right\|_{L^2_x}^2 + \frac 1{2\sigma} \left\|\left(j-nu\right)(t)\right\|_{L^2_x}^2,
	\end{aligned}
\end{equation*}
or, for the system \eqref{TFINSMSO},
\begin{equation*}
	\begin{aligned}
		\CE(t) & =
		\frac 12 \left\|u(t)\right\|_{L^2_x}^2
		+ \frac 14 \left\|E(t)\right\|_{L^2_x}^2 + \frac 14 \left\|B(t)\right\|_{L^2_x}^2, \\
		\CD(t) & =
		\mu \left\|\nabla_x u(t)\right\|_{L^2_x}^2 + \frac 1{2\sigma} \left\|j(t)\right\|_{L^2_x}^2.
	\end{aligned}
\end{equation*}
\end{prop}

\begin{proof}
	We consider the system \eqref{TFINSMO} first. Multiplying the equation expressing the conservation of momentum in \eqref{TFINSMO} by $u$ and integrating with respect to space variables, we get
	\begin{equation*}
		\frac12 {d\over dt} \| u\|_{L^2_x}^2 + \mu \left\|\nabla_x u \right\|_{L^2_x}^2 =
		\frac 12 \int_{\mathbb{R}^3} \left(nE+j\wedge B\right)\cdot u dx ,
	\end{equation*}
	while multiplying Ohm's law in \eqref{TFINSMO} by $j-nu$ and integrating in space yields the identity
	\begin{equation*}
		\frac 1\sigma\left\|j-nu\right\|_{L^2_x}^2=
		\int_{\mathbb{R}^3} \left(\frac 12 n\DIV j + E\cdot j- \left(nE+j\wedge B\right)\cdot u\right) dx,
	\end{equation*}
	where we have employed the incompressibility of the velocity field. Hence, we obtain, further exploiting the continuity equation $\partial_t n+\DIV j=0$ (deduced by taking the divergence of Amp\`ere's equation and from Gauss' law), that
	\begin{equation*}
		\begin{aligned}
			\frac12 {d\over dt} \| u\|_{L^2_x}^2 + \mu \left\|\nabla_x u \right\|_{L^2_x}^2
			+ \frac 1{2\sigma}\left\|j-nu\right\|_{L^2_x}^2 & =
			\int_{\mathbb{R}^3} \frac 14 n\DIV j + \frac 12 E\cdot j dx \\
			& = -\frac 18\frac{d}{dt}\left\|n\right\|_{L^2_x}^2
			+ \int_{\mathbb{R}^3} \frac 12 E\cdot j dx.
		\end{aligned}
	\end{equation*}
	
	As for the system \eqref{TFINSMSO}, similar and, actually, simpler computations yield that
	\begin{equation*}
		\frac12 {d\over dt} \| u\|_{L^2_x}^2 + \mu \left\|\nabla_x u \right\|_{L^2_x}^2
		+ \frac 1{2\sigma}\left\|j\right\|_{L^2_x}^2 =
		\int_{\mathbb{R}^3} \frac 12 E\cdot j dx.
	\end{equation*}
	
	Next, for both systems \eqref{TFINSMO} and \eqref{TFINSMSO}, the conservation of the electromagnetic energy is given by Maxwell's equations
	\begin{equation*}
		\frac 12 \frac d{dt}\left(\left\|E\right\|_{L^2_x}^2+\left\|B\right\|_{L^2_x}^2\right)
		= -\int_{\mathbb{R}^3}E\cdot j dx.
	\end{equation*}
	Summing the above formal identities leads to the expected global conservation of energy.
\end{proof}

The uniform bounds resulting from the energy conservations in Proposition \ref{energy estimate 2} imply that all the terms in systems \eqref{TFINSMO} and \eqref{TFINSMSO} make sense, especially the nonlinear terms in the motion equations and in Ohm's laws.

Notice, however, that it is at first not clear that the Lorentz force $nE + j\wedge B$ in \eqref{TFINSMO} is a well-defined distribution, based on the natural a priori estimates provided by the energy and energy dissipation, because $j$ does not necessarily lie in $L^1_tL_x^2$. Nevertheless, it is possible to give it a rigorous sense by exploiting simple identities. A first approach consists in identifying the force term $nE+j\wedge B$ with the conservation law for the electromagnetic energy flux $E\wedge B$ (also called the Poynting vector, see \cite[Section 6.7]{jackson})
\begin{equation}\label{poynting}
	\partial_t\left(E\wedge B\right)+\frac 12\nabla_x\left(E^2+B^2\right)-\nabla_x\cdot\left(E\otimes E+B\otimes B\right)
	=
	-nE-j\wedge B
\end{equation}
derived directly from Maxwell's equations in \eqref{TFINSMO} (see the derivation of \eqref{maxwell poynting} and \eqref{maxwell poynting two species}), so that the force makes sense in some Sobolev space with negative regularity index. In fact, it will be much more appropriate to estimate the Lorentz force directly using Ohm's law from \eqref{TFINSMO} as follows
\begin{equation*}
	\begin{aligned}
		nE+j\wedge B & =
		(j-nu)\wedge B + n\left(E+u\wedge B\right) \\
		& =
		(j-nu)\wedge B + \frac 1\sigma n\left(j-nu\right)+\frac 14\nabla_x\left(n^2\right),
	\end{aligned}
\end{equation*}
so that the force is now understood as the sum of a locally integrable function and a pressure gradient. All other terms from \eqref{TFINSMO} and \eqref{TFINSMSO} are obviously well-defined.

\bigskip

Unfortunately, the uniform energy bounds do not guarantee the weak stability of the nonlinear terms $nE$ and $j\wedge B$ composing the Lorentz forces. This is a major obstacle to establishing the global existence of weak solutions in the spirit of Leray \cite{leray}, which are therefore not known to exist in general.

There are two evident strategies, which unfortunately turn out to be unsuccessful, that one would want to apply here in order to circumvent the lack of weak stability of the Lorentz forces in systems \eqref{TFINSMO} and \eqref{TFINSMSO}.

The first one consists in propagating strong compactness or regularity in Maxwell's equations, which are indeed the archetype of hyperbolic equations, meaning that singularities are propagated. In general, these singularities, or oscillations, may be created either by boundary data, by initial data or by the source terms, and they remain localized on the corresponding light cones. Here, we are not considering boundaries and the initial data can always be well-prepared. However, it remains unclear how to prevent the emergence of oscillations from the source term $-j$ in Maxwell's equations, which is determined by the nonlinear Ohm's laws in \eqref{TFINSMO} and \eqref{TFINSMSO}. Therefore, we do not expect to gain regularity (or even compactness) on the electromagnetic field $(E,B)$. So, this strategy fails in general.

It is to be noted, though, that this approach has been successfully applied by Masmoudi \cite{masmoudi} to a slightly different system coupling the incompressible Navier-Stokes equations with Maxwell's equations in the two-dimensional case. Since the equations studied therein are very similar to \eqref{TFINSMO} and \eqref{TFINSMSO}, we present Masmoudi's result below in Section \ref{nader} in order to emphasize the mathematical difficulties inherent to the coupling with Maxwell's equations through Ohm's law and its similarities with the two-dimensional Euler equations. Also, we believe that similar results on systems \eqref{TFINSMO} and \eqref{TFINSMSO} can be achieved.

The second strategy consists in utilizing the linear structure of Maxwell's equations with the specific quadratic structure of the Lorentz force to apply the theory of compensated compactness of Murat and Tartar \cite{murat, murat2, tartar} (see also \cite{struwe} for an introduction to the subject) and, thus, filter any undesired nonlinear resonances. This approach plainly fails and it seems that it can only potentially succeed by exploiting the full nonlinear structure of the whole systems \eqref{TFINSMO} and \eqref{TFINSMSO}. But we are not aware of such successful nonlinear treatment of resonances. We refer to \cite{arsenio6} for some more details about the failure of the method of compensated compactness in the electromagnetic setting.

\bigskip

Following the concise Section \ref{nader} below, where we present the main result from \cite{masmoudi} on the well-posedness of an incompressible Navier-Stokes-Maxwell system in two dimensions, we will discuss very briefly in Section \ref{slim} the well-posedness of the same system in three dimensions and for small initial data. Finally, in Section \ref{laure-diogo}, we will introduce the dissipative solutions of the systems \eqref{TFINSMO} and \eqref{TFINSMSO} and justify their global existence in any dimension, which will be particularly relevant to our work.

\subsection{Large global solutions in two dimensions}\label{nader}% {Global strong solutions in 2D}

In \cite{masmoudi}, Masmoudi studied the following incompressible Navier-Stokes-Maxwell system~:
\begin{equation}\label{INSM}
	\begin{cases}
		\begin{aligned}
			\d_t u +
			u\cdot\nabla_x u - \mu\Delta_x u
			& = -\nabla_x p+
			j \wedge B , & \Div u & = 0,\\
			\d_t E - \ROT B &= -  j, & j & = \sigma\left( E + u\wedge B\right),
			\\
			\d_t B + \ROT E & = 0, & \Div B & = 0,
		\end{aligned}
	\end{cases}
\end{equation}
which is somewhat related to the systems \eqref{TFINSMO} and \eqref{TFINSMSO}, and satisfies the formal energy conservation
\begin{equation}\label{INSM energy}
		\frac 12\frac{d}{dt}\left(\left\|u\right\|_{L^2_x}^2 + \left\|E\right\|_{L^2_x}^2 + \left\|B\right\|_{L^2_x}^2 \right)
		+ \mu \left\|\nabla_x u\right\|_{L^2_x}^2 + \frac 1{\sigma} \left\|j\right\|_{L^2_x}^2 = 0.
\end{equation}
Notice that, in this system, there is no constraint on $\Div E$ or $\Div j$. He restricted his analysis to the two-dimensional case, which is obtained by assuming that
\begin{equation*}
	u=\begin{pmatrix}u_1(x_1,x_2) \\ u_2(x_1,x_2) \\ 0\end{pmatrix},
	\quad
	E=\begin{pmatrix}E_1(x_1,x_2) \\ E_2(x_1,x_2) \\ 0\end{pmatrix}
	\quad\text{and}\quad
	B=\begin{pmatrix}0 \\ 0 \\ B_3(x_1,x_2)\end{pmatrix}.
\end{equation*}

\bigskip

In order to understand the propagation of singularities in Maxwell's system (in two or three dimensions), it is often convenient to express it using vector and scalar potentials in an equivalent form (see \cite[Sections 6.2 and 6.3]{jackson}). To this end, since the magnetic field $B$ is solenoidal, we may always write $B=\ROT A$, for some vector potential $A$. Moreover, taking into account Faraday's equation, we see that necessarily $E=-\nabla_x\varphi-\partial_t A$ for some scalar potential $\varphi$. As a matter of fact, the potentials $A$ and $\varphi$ are not uniquely determined. Indeed, the electromagnetic field is invariant under the so-called gauge transformation $(A,\varphi)\mapsto(A+\nabla_x\psi,\varphi-\partial_t\psi)$.

This gauge invariance allows us to impose a further condition on the potentials. Typically, one may impose the so-called Coulomb gauge $\DIV A=0$, which is simple and natural for stationary settings. Another classical example of gauge fixing includes the Lorenz (not to be confused with Lorentz) gauge $\DIV A = -\partial_t\varphi$, which usually yields an evolution for the potentials governed by decoupled wave equations.

Here, for the Maxwell system in \eqref{INSM}, we choose the slight variant of the Lorenz gauge
\begin{equation}\label{gauge}
	\DIV A = - \partial_t\varphi - \sigma \varphi,
\end{equation}
which yields the decoupled damped wave equation
\begin{equation}\label{damped wave}
	\partial_t^2 A + \sigma \partial_t A - \Delta_x A = \sigma u\wedge\left(\ROT A\right).
\end{equation}
Note that it is always possible to find $A$ and $\varphi$ satisfying \eqref{gauge}. Indeed, if \eqref{gauge} is not satisfied, one may always apply a gauge transformation with $\psi$ solving the damped wave equation $\partial_t^2\psi + \sigma \partial_t\psi - \Delta_x \psi = \DIV A + \partial_t\varphi + \sigma \varphi$ and produce new potentials for which \eqref{gauge} holds.

\bigskip

Now, if the velocity field $u$ is bounded in $L^1\left([0,T],dt;L^\infty(dx)\right)$, for some $T>0$, it is possible to show, through standard energy estimates, that the damped wave equation \eqref{damped wave}, which is linear in $A$, propagates the strong compactness of $\partial_t A$ and $\nabla_x A$ in $L^\infty\left([0,T],dt;L^2(dx)\right)$. This would obviously imply the propagation of strong compactness for the magnetic field $B$.

Unfortunately, in two dimensions of space, the $H^1$ estimate on the velocity field $u$ provided by the conservation of energy barely fails to yield, by Sobolev embedding, an $L^\infty$ bound on $u$. Masmoudi's idea was then to compensate this lack of critical embedding by placing the initial electromagnetic field in a better $H^s$ space, with $0<s<1$, and to propagate this initial regularity with Maxwell's equations at the same time that the parabolic regularity of the Stokes flow is employed to estimate the velocity field in a higher regularity space. This approach eventually allows to bound $u$ in $L^1_tL^\infty_x$ in terms of its $L^2_t\dot H^1_x$ norm with some logarithmic loss.

As a byproduct of these estimates, it is also possible to establish the exponential growth of the $H^s$ norms. In Masmoudi's own words~: ``One can compare this growth estimate with the double exponential growth estimate of the $H^s$ norms in the two-dimensional incompressible Euler system.''

Finally, it is interesting to note that Masmoudi's proof uses neither the divergence free condition of the magnetic field nor the decay property of the linear part coming from Maxwell's equations.

The following theorem contains the main well-posedness result from \cite{masmoudi}. Note that it gives the existence and uniqueness for initial data in a very large dense subspace of $L^2$, namely in $\cup_{0<s<1}H^s$, but it fails to guarantee the existence of a weak solution when the initial data lies merely in $L^2$.

\begin{thm}[\cite{masmoudi}]
	Take $0 < s < 1$,
	\begin{equation*}
		u^\mathrm{in}\in  L^2\left(\mathbb{R}^2\right)
		\qquad
		\text{and}
		\qquad
		E^\mathrm{in}, B^\mathrm{in} \in H^s\left(\mathbb{R}^2\right). 
	\end{equation*}
	
	Then, there exists a unique global solution $(u,E,B)$ of \eqref{INSM} such that for all $T > 0$,
	\begin{equation*}
		u\in C\left([0,T]; L^2\right) \cap L^2\left([0,T]; \dot H^1\right)
		\qquad
		\text{and}
		\qquad
		E,B \in C\left([0,T]; H^s\right).
	\end{equation*}
	Moreover,
	\begin{equation*}
		j \in L^2\left([0,T];L^2\right) \cap L^1\left([0,T]; H^s\right)
		\qquad
		\text{and}
		\qquad
		u\in L^1\left([0,T]; H^{s'}\right),
	\end{equation*}
	for each $ 1 < s' < \min \left(2s+1,2\right)$. In addition, the energy identity \eqref{INSM energy} holds and we have the following exponential growth estimate for all $t>0$~:
	\begin{equation*}
		\| u \|_{L^1\left([0,t]; H^{s'}\right)}+\left\| E(t)\right\|_{H^s}+\left\|B(t)\right\|_{H^s} \leq
		\left(1+\left\| E^\mathrm{in}\right\|_{H^s}+\left\|B^\mathrm{in}\right\|_{H^s}\right) e^{C^\mathrm{in}(1+t)},
	\end{equation*}
	where $C^\mathrm{in}=C\left(1+\left\| u^\mathrm{in}\right\|_{L^2}^2+\left\| E^\mathrm{in}\right\|_{L^2}^2
	+\left\|B^\mathrm{in}\right\|_{L^2}^2\right)$ for some constant $C$.
\end{thm}

\subsection{Small global solutions in three dimensions}\label{slim}% {Small solutions in 3D}

As we have seen, there are serious obstacles to the construction of global solutions of the system \eqref{INSM} for large initial data in the energy space. Nevertheless, it is in general possible to achieve the well-posedness of a system, globally in time, by showing its strong stability for small initial data in some space satisfying the same scaling invariance as the given system of equations.

This is precisely what Ibrahim and Keraani managed to obtain in \cite{ibrahim} for the three-dimensional incompressible Navier-Stokes-Maxwell system \eqref{INSM} using the strategy of Fujita and Kato \cite{fujita}, which is based on refined a priori estimates obtained by paradifferential calculus and some fixed point argument. Note that the results from \cite{ibrahim} do not imply the local existence of strong solutions for large data, which has been established in a separate work by Ibrahim and Yoneda in \cite{ibrahim2}.

These results have then been unified and extended to a more natural setting by Germain, Ibrahim and Masmoudi in \cite{germain1}.

We believe that the methods employed in \cite{germain1, ibrahim, ibrahim2} can potentially lead to similar results for the analogous incompressible Navier-Stokes-Maxwell systems \eqref{TFINSMO} and \eqref{TFINSMSO}. The main result in this three-dimensional setting is contained in the following theorem (we refer directly to \cite{germain1} for definitions of the functional spaces).

\begin{thm}[\cite{germain1, ibrahim, ibrahim2}]
	To any initial data
	\begin{equation*}
		u^\mathrm{in}, E^\mathrm{in}, B^\mathrm{in} \in \dot H^\frac{1}{2}\left(\mathbb{R}^3\right), 
	\end{equation*}
	there corresponds an existence time $T>0$ and a unique local solution of \eqref{INSM}
	\begin{equation*}
		u\in \tilde L^\infty\left((0,T); \dot H^\frac 12 \right) \cap L^2\left((0,T); \dot H^\frac{3}{2}\cap L^\infty\right)
		\quad
		\text{and}
		\quad
		E,B \in \tilde L^\infty \left((0,T); \dot H^\frac{1}{2}\right).
	\end{equation*}
	Furthermore, the solution is global (i.e.\ $T=\infty$) if the initial data is sufficiently small.
	% $E,B \in C \left([0,\infty); \dot H^\frac{1}{2}\right)$ and $E \in L^2 \left([0,\infty); \dot H^\frac{1}{2}\right)$.
\end{thm}

\subsection{Weak-strong stability and dissipative solutions}\label{laure-diogo}

On the one hand, As already explained, there is no known global well-posedness theory for the systems \eqref{TFINSMO}, \eqref{TFINSMSO} and \eqref{INSM} in the energy space, due to their lack of weak stability. On the other hand, in Sections \ref{nader} and \ref{slim}, we have briefly presented theorems on the existence and uniqueness of strong solutions to the system \eqref{INSM}.

Around such smooth solutions and in order to circumvent the lack of weak stability, we introduce now the dissipative solutions of these incompressible Navier-Stokes-Maxwell systems. Generally speaking, the concept of dissipative solutions is based on the weak-strong stability, when available, of a given system, i.e.\ the uniqueness of all weak solutions provided at least one strong solution exists. It seems that such weak-strong stability principles were first introduced by Dafermos \cite{dafermos} in the context of conservation laws.

Dissipative solutions are not new in fluid and gas dynamics. They are precisely employed to treat the instability of nonlinear terms in the energy space. Lions first defined them for the Boltzmann equation in \cite{lions2}. He then established their existence for the incompressible Euler system in \cite[Section 4.4]{lions7}, as an alternative to the very weak notion of measure-valued solutions introduced by DiPerna and Majda \cite{diperna4}, which have later been shown in \cite{brenier} to be actually stronger (at least not weaker, as each measure-valued solution is shown to be a dissipative solution, as well). It can more easily be shown that any weak solution of the incompressible Euler system is a dissipative solution (see \cite[Appendix B]{delellis3} for a proof). This, however, is not known to hold in general for renormalized solutions of the Boltzmann equation, i.e.\ renormalized solutions are not known, in general, to be dissipative solutions as defined by Lions in \cite{lions2}.

It is sometimes argued that dissipative solutions are too weak and that they do not express any physical reality, because they are not shown to be unique in general. Even so, they do enjoy certain definite qualities~:
\begin{itemize}
	\item they exist globally in time for large initial data in the energy space~;
	\item they coincide with the unique strong solution when the latter exists~;
	\item they allow energy dissipation phenomena to occur.
\end{itemize}
The last property above is especially significant in light of recent results on the energy dissipation in the incompressible Euler flow establishing, in particular, the existence of weak solutions with kinetic energy strictly decaying (or increasing, which is equivalent since the Euler flow is reversible) over time (see \cite{delellis3, shnirelman}). This energy dissipation cannot hold beyond a certain regularity threshold (see \cite{cheskidov, constantin} on Onsager's conjecture) and, therefore, it is crucial to consider rather low regularity weak solutions of the incompressible Euler system in order to understand energy dissipation and turbulent flow. In this context, we wish to mention the striking recent developments \cite{buckmaster2, buckmaster, buckmaster3, choffrut, delellis2, delellis} demonstrating the existence of energy-dissipating flows enjoying some H\"older regularity.

Dissipative solutions have found an important application in a wide range of asymptotic problems, for they are especially well adapted, through relative entropy methods (or modulated energy methods), to situations presenting a lack of compactness. In particular, they were employed by the second author in \cite{SR2, SR3} to establish the hydrodynamic convergence of renormalized solutions of the Boltzmann equation towards dissipative solutions of the incompressible Euler system (see also \cite{SR}). Another application by Brenier \cite{brenier2} concerns the convergence of the Vlasov-Poisson system towards the incompressible Euler equations in the quasi-neutral regime.

% \bigskip

\subsubsection{The incompressible Navier-Stokes-Maxwell system}

We explain now how the energy \eqref{INSM energy} can be modulated and establish a weak-strong stability principle, which will eventually lead to a suitable notion of dissipative solution for the incompressible Navier-Stokes-Maxwell system \eqref{INSM} in any dimension. We will then move on to apply the same strategy to the more complex systems \eqref{TFINSMO} and \eqref{TFINSMSO} by modulating the energies from Proposition \ref{energy estimate 2} and thus produce similar dissipative solutions.

\begin{prop}\label{modulated energy estimate 0}
	Let $(u,E,B)$ be a smooth solution to the incompressible Navier-Stokes-Maxwell system \eqref{INSM}. Further consider test functions $\left(\bar u, \bar j, \bar E, \bar B\right)\in C_c^\infty\left([0,\infty)\times\mathbb{R}^3\right)$ such that
	\begin{equation}\label{test constraints}
		\begin{cases}
			\begin{aligned}
				\d_t \bar E - \ROT \bar B &= -  \bar j, & \Div \bar u & = 0,
				\\
				\d_t \bar B + \ROT \bar E & = 0, & \Div \bar B & = 0.
			\end{aligned}
		\end{cases}
	\end{equation}
	We define the acceleration operator by
	\begin{equation*}
		\mathbf{A}\left(\bar u, \bar j, \bar E, \bar B\right)
		=
		\begin{pmatrix}
			-\d_t \bar u -
			P\left(\bar u\cdot\nabla_x \bar u\right) + \mu\Delta_x \bar u
			+ P \left(\bar j \wedge \bar B\right)
			\\
			- \frac1\sigma\bar j + \bar E+\bar u\wedge\bar B
		\end{pmatrix},
	\end{equation*}
	and the growth rate by
	\begin{equation*}
		\lambda(t)=
		\left(\frac 2\mu + \sigma\right)\left\|\bar u(t)\right\|_{L^\infty_x}^2
		+ \frac{2C_0^2}{\mu} \left\|\bar j(t)\right\|_{L^3_x}^2,
	\end{equation*}
	where $C_0>0$ denotes the operator norm of the Sobolev embedding $\dot H^1\left(\mathbb{R}^3\right) \hookrightarrow L^6\left(\mathbb{R}^3\right)$.
	
	Then, one has the stability inequality
	\begin{equation}\label{stability 0}
		\begin{aligned}
			\delta\CE(t) + & \frac 12 \int_0^t \delta\CD(s) e^{\int_s^t\lambda(\sigma)d\sigma}ds
			\\
			& \leq \delta\CE(0) e^{\int_0^t\lambda(s)ds}
			+\int_0^t
			\left[\int_{\mathbb{R}^3} \mathbf{A}\cdot
			\begin{pmatrix}
				u-\bar u \\ j - \bar j
			\end{pmatrix}
			dx\right](s)
			e^{\int_s^t\lambda(\sigma)d\sigma}ds,
		\end{aligned}
	\end{equation}
	where the modulated energy $\delta\CE$ and energy dissipation $\delta\CD$ are given by
	\begin{equation}\label{modulated defi}
		\begin{aligned}
			\delta\CE(t)
			& =
			\frac 12 \left\|\left(u-\bar u\right)(t)\right\|_{L^2_x}^2
			+ \frac 12 \left\|\left(E-\bar E\right)(t)\right\|_{L^2_x}^2
			+ \frac 12 \left\|\left(B-\bar B\right)(t)\right\|_{L^2_x}^2, \\
			\delta \CD(t)
			& =
			\mu \left\|\nabla_x (u-\bar u)(t)\right\|_{L^2_x}^2
			+ \frac 1{\sigma} \left\|\left(j-\bar j\right)(t)\right\|_{L^2_x}^2.
		\end{aligned}
	\end{equation}
\end{prop}

\begin{proof}
	We have already formally established in \eqref{INSM energy} the conservation of the energy for the system \eqref{INSM}. The very same computations applied to the test functions $\left(\bar u, \bar j, \bar E, \bar B\right)$ yield the identity
	\begin{equation*}
		\frac 12\frac{d}{dt}\left(\left\|\bar u\right\|_{L^2_x}^2 + \left\|\bar E\right\|_{L^2_x}^2 + \left\|\bar B\right\|_{L^2_x}^2 \right)
		+ \mu \left\|\nabla_x \bar u\right\|_{L^2_x}^2 + \frac 1{\sigma} \left\|\bar j\right\|_{L^2_x}^2 = - \int_{\mathbb{R}^3} \mathbf{A}\cdot
		\begin{pmatrix}
			\bar u \\ \bar j
		\end{pmatrix}
		dx.
	\end{equation*}
	
	Furthermore, another similar duality computation gives
	\begin{equation*}
		\begin{aligned}
			\frac{d}{dt}\int_{\mathbb{R}^3} & u\cdot\bar u + E\cdot \bar E+B\cdot\bar B dx
			+\int_{\mathbb{R}^3}2\mu \nabla_x u:\nabla_x\bar u
			+\frac 2\sigma j\cdot\bar j dx \\
			& =
			- \int_{\mathbb{R}^3}
			\bar u \otimes\left(u-\bar u\right):\nabla_x(u-\bar u) dx
			\\
			% +
			% \left(u-\bar u)^{\otimes 2}:\nabla_x\bar u
			& + \int_{\mathbb{R}^3}
			\left((j-\bar j)\wedge(B-\bar B)\right)\cdot\bar u
			+ \left((u-\bar u)\wedge(B-\bar B)\right)\cdot\bar j
			dx
			\\
			& - \int_{\mathbb{R}^3} \mathbf{A}\cdot
			\begin{pmatrix}
				u \\ j
			\end{pmatrix}
			dx.
		\end{aligned}
	\end{equation*}
	
	On the whole, combining the above identities with the formal energy conservation \eqref{INSM energy}, we find
	\begin{equation*}
		\begin{aligned}
			\frac 12 \frac{d}{dt} & \left(\left\|u-\bar u\right\|_{L^2_x}^2 + \left\|E-\bar E\right\|_{L^2_x}^2 + \left\|B-\bar B\right\|_{L^2_x}^2 \right)
			+ \mu \left\|\nabla_x (u-\bar u)\right\|_{L^2_x}^2 + \frac 1{\sigma} \left\|j-\bar j\right\|_{L^2_x}^2
			\\
			& =
			\int_{\mathbb{R}^3}
			\bar u \otimes\left(u-\bar u\right):\nabla_x(u-\bar u) dx
			\\
			% +
			% \left(u-\bar u)^{\otimes 2}:\nabla_x\bar u
			& + \int_{\mathbb{R}^3}
			\left(\bar j\wedge(B-\bar B)\right)\cdot(u-\bar u)
			- \left((j-\bar j)\wedge(B-\bar B)\right)\cdot\bar u
			dx
			\\
			& + \int_{\mathbb{R}^3} \mathbf{A}\cdot
			\begin{pmatrix}
				u-\bar u \\ j - \bar j
			\end{pmatrix}
			dx.
		\end{aligned}
	\end{equation*}
	The next step consists in estimating the terms in the right-hand side above that are nonlinear in $(u,j,E,B)$ and to absorb the resulting expressions with the modulated energy $\delta\CE(t)$ and the modulated energy dissipation $\delta\CD(t)$ by suitable uses of Young's inequality and Gr\"onwall's lemma. Thus, we obtain
	\begin{equation*}
		\begin{aligned}
			\frac{d}{dt}\delta\CE(t)
			& +\delta\CD(t)
			\\
			& \leq
			\left\|\bar u\right\|_{L^\infty_x} \left\| u-\bar u\right\|_{L^2_x}\left\|\nabla_x(u-\bar u)\right\|_{L^2_x}
			\\
			& +
			\left\|\bar j\right\|_{L^3_x}\left\|B-\bar B\right\|_{L^2_x}\left\|u-\bar u\right\|_{L^6_x}
			+ \left\|\bar u\right\|_{L^\infty_x}\left\|B-\bar B\right\|_{L^2_x}\left\|j-\bar j\right\|_{L^2_x}
			\\
			& + \int_{\mathbb{R}^3} \mathbf{A}\cdot
			\begin{pmatrix}
				u-\bar u \\ j - \bar j
			\end{pmatrix}
			dx
			\\
			& \leq
			\frac 1{\mu}\left\|\bar u\right\|_{L^\infty_x}^2 \left\| u-\bar u\right\|_{L^2_x}^2
			+\left(\frac\sigma 2\left\|\bar u\right\|_{L^\infty_x}^2
			+ \frac{C_0^2}{\mu} \left\|\bar j\right\|_{L^3_x}^2\right)
			\left\|B-\bar B\right\|_{L^2_x}^2
			\\
			& +\frac\mu 2\left\|\nabla_x(u-\bar u)\right\|_{L^2_x}^2
			+ \frac{1}{2\sigma}\left\|j-\bar j\right\|_{L^2_x}^2
			+ \int_{\mathbb{R}^3} \mathbf{A}\cdot
			\begin{pmatrix}
				u-\bar u \\ j - \bar j
			\end{pmatrix}
			dx.
		\end{aligned}
	\end{equation*}
	
	Hence,
	\begin{equation*}
		\frac{d}{dt}\delta\CE(t)
		+\frac 12\delta\CD(t)
		\leq
		\lambda(t)
		\delta\CE(t)
		+ \int_{\mathbb{R}^3} \mathbf{A}\cdot
		\begin{pmatrix}
			u-\bar u \\ j - \bar j
		\end{pmatrix}
		dx,
	\end{equation*}
	which concludes the proof of the proposition with a direct application of Gr\"onwall's lemma.
\end{proof}

Note that the test functions satisfying the linear constraints \eqref{test constraints} are easily constructed by considering scalar potentials $\bar \varphi\in C_c^\infty\left([0,\infty)\times\mathbb{R}^3\right)$ and vector potentials $\bar A\in C_c^\infty\left([0,\infty)\times\mathbb{R}^3\right)$ and then setting
\begin{equation}\label{test construction}
	\bar E = -\nabla_x\bar\varphi - \partial_t \bar A \qquad\text{and}\qquad \bar B = \ROT \bar A.
\end{equation}
One may prefer, for various reasons, to deal, in a completely equivalent manner, with test functions $\left(\bar u, \bar j, \bar E, \bar B\right)\in C_c^\infty\left([0,\infty)\times\mathbb{R}^3\right)$ satisfying the stationary constraints
\begin{equation*}
	\Div \bar u = 0,
	\qquad
	\Div \bar B = 0,
	\qquad
	\bar j = \sigma\left(\bar E + \bar u\wedge \bar B\right),
\end{equation*}
rather than the constraints \eqref{test constraints}. In this case, instead of \eqref{stability 0}, we obtain the stability inequality
\begin{equation*}
	\begin{aligned}
		\delta\CE(t) + & \frac 12 \int_0^t \delta\CD(s) e^{\int_s^t\lambda(\sigma)d\sigma}ds
		\\
		& \leq \delta\CE(0) e^{\int_0^t\lambda(s)ds}
		+\int_0^t
		\left[\int_{\mathbb{R}^3} \mathbf{A}\cdot
		\begin{pmatrix}
			u-\bar u \\ E - \bar E \\ B - \bar B
		\end{pmatrix}
		dx\right](s)
		e^{\int_s^t\lambda(\sigma)d\sigma}ds,
	\end{aligned}
\end{equation*}
where the acceleration operator is now defined by
\begin{equation*}
	\mathbf{A}\left(\bar u, \bar j, \bar E, \bar B\right)
	=
	\begin{pmatrix}
		-\d_t \bar u -
		P\left(\bar u\cdot\nabla_x \bar u\right) + \mu\Delta_x \bar u
		+ P \left(\bar j \wedge \bar B\right)
		\\
		-\partial_t \bar E + \ROT\bar B - \bar j
		\\
		-\partial_t \bar B - \ROT\bar E
	\end{pmatrix}.
\end{equation*}

The preceding proposition provides an important weak-strong stability property for the incompressible Navier-Stokes-Maxwell system \eqref{INSM}. Indeed, the stability inequality \eqref{stability 0} essentially implies that a solution $(\bar u,\bar j,\bar E,\bar B)$ of \eqref{INSM} such that $\bar u\in L^2_tL^\infty_x$ and $\bar j\in L^2_tL^3_x$, if it exists, is unique in the whole class of weak solutions in the energy space, for any given initial data.

\begin{rem}
	In order to impose minimal local integrability assumptions on the test function $\bar u$, it is tempting to employ the method of Lions and Masmoudi \cite{lions5} for estimating the nonlinear term $\left[\bar u\otimes\left(u-\bar u\right)\right]:\nabla_x\left(u-\bar u\right)$ by splitting $\bar u = \bar u \mathds{1}_{\left\{\left|\bar u\right|\leq K\right\}} + \bar u \mathds{1}_{\left\{\left|\bar u\right| > K\right\}}$, for some large $K>0$, which yields
\begin{equation*}
	\begin{aligned}
		& \left\|\left[\bar u\otimes\left(u-\bar u\right)\right]:\nabla_x\left(u-\bar u\right)\right\|_{L^1_x}
		\\
		& \hspace{10mm} \leq
		K \left\| u-\bar u\right\|_{L^2_x}\left\|\nabla_x(u-\bar u)\right\|_{L^2_x}
		+
		\left\|\bar u\mathds{1}_{\left\{\left|\bar u\right| > K\right\}} \right\|_{L^3_x}
		\left\| u-\bar u\right\|_{L^6_x}\left\|\nabla_x(u-\bar u)\right\|_{L^2_x}
		\\
		& \hspace{10mm} \leq
		\frac 1\nu\frac {K^2}4 \left\| u-\bar u\right\|_{L^2_x}^2
		+
		\nu \left\|\nabla_x(u-\bar u)\right\|_{L^2_x}^2
		+
		C_0 \left\|\bar u\mathds{1}_{\left\{\left|\bar u\right| > K\right\}} \right\|_{L^3_x}
		\left\|\nabla_x(u-\bar u)\right\|_{L^2_x}^2,
	\end{aligned}
\end{equation*}
for any $\nu>0$, where $C_0>0$ denotes the operator norm of the Sobolev embedding $\dot H^1_x \hookrightarrow L^6_x$. Then, choosing $K>0$ large enough so that $\left\|\bar u \mathds{1}_{\left\{\left|\bar u\right|> K\right\}} \right\|_{L^3_x}$ is arbitrarily small and setting $\nu>0$ small enough, it is readily seen that the last two terms above can be absorbed by the modulated entropy dissipation. Of course, the choice of the parameter $K$ is not uniform for all $\bar u\in L^3_x$. This approach definitely allows us to merely consider velocity fields $\bar u\in L^2_tL^\infty_x+L^\infty_tL^3_x$ when establishing weak-strong stability principles for the incompressible Navier-Stokes system (see \cite{lions5}). Here, however, considering the coupling of the fluid equations with Maxwell's system introduces other nonlinear terms in the estimates, which unfortunately require that $\bar u\in L^2_tL^\infty_x$ in order to be duly controlled.
\end{rem}

By analogy with Lions' dissipative solutions to the incompressible Euler system \cite[Section 4.4]{lions7}, we provide now a suitable notion of dissipative solution for the incompressible Navier-Stokes-Maxwell system \eqref{INSM}, based on Proposition \ref{modulated energy estimate 0}, and establish their existence next.

\begin{defi}
	We say that
	\begin{equation*}
		(u,E,B)\in L^\infty\left([0,\infty);L^2\left(\mathbb{R}^3\right)\right)
		\cap C\left([0,\infty);\textit{w-}L^2\left(\mathbb{R}^3\right)\right)
	\end{equation*}
	such that
	\begin{equation*}
		\DIV u = 0, \qquad \DIV B = 0,
	\end{equation*}
	is a \textbf{dissipative solution of the incompressible Navier-Stokes-Maxwell system \eqref{INSM}}, if it solves Maxwell's equations
	\begin{equation*}
		\begin{cases}
			\begin{aligned}
				\d_t E - \ROT B &= -  j, 
				\\
				\d_t B + \ROT E & = 0,
			\end{aligned}
		\end{cases}
	\end{equation*}
	with Ohm's law
	\begin{equation*}
		j = \sigma\left( E + u\wedge B\right),
	\end{equation*}
	in the sense of distributions, and if, for any test functions $\left(\bar u, \bar j, \bar E, \bar B\right)\in C_c^\infty\left([0,\infty)\times\mathbb{R}^3\right)$ satisfying the linear constraints \eqref{test constraints}, the stability inequality \eqref{stability 0} is verified.
\end{defi}

As previously mentioned, dissipative solutions define actual solutions in the sense that they coincide with the unique strong solution when the latter exists. The following theorem asserts their existence.

\begin{thm}\label{dissipative 1}
	For any initial data $\left(u^\mathrm{in},E^\mathrm{in},B^\mathrm{in}\right)\in L^2\left(\mathbb{R}^3\right)$ such that
	\begin{equation*}
		\DIV u^\mathrm{in} = 0, \qquad \DIV B^\mathrm{in} = 0,
	\end{equation*}
	there exists a dissipative solution to the incompressible Navier-Stokes-Maxwell system \eqref{INSM}.
\end{thm}

\begin{proof}
	Following Lions \cite{lions7}, we easily build the dissipative solutions by introducing viscous approximations of the system \eqref{INSM}. Thus, for each $\nu>0$, we consider weak solutions of the following system~:
	\begin{equation*}
		\begin{cases}
			\begin{aligned}
				\d_t u_\nu +
				u_\nu\cdot\nabla_x u_\nu - \mu\Delta_x u_\nu
				& = -\nabla_x p_\nu+
				j_\nu \wedge B_\nu , & \Div u_\nu & = 0,\\
				\d_t E_\nu - \ROT B_\nu &= -  j_\nu, & j_\nu & = \sigma\left( E_\nu + u_\nu\wedge B_\nu\right),
				\\
				\d_t B_\nu + \ROT E_\nu -\nu\Delta_x B_\nu & = 0, & \Div B_\nu & = 0,
			\end{aligned}
		\end{cases}
	\end{equation*}
	associated with the initial data $\left(u^\mathrm{in},E^\mathrm{in},B^\mathrm{in}\right)$ and satisfying the energy inequality, for all $t>0$,
	% \begin{equation*}
	% 	\begin{aligned}
	% 		\frac 12 \frac{d}{dt} & \left(\left\|u_\nu\right\|_{L^2_x}^2 + \left\|E_\nu\right\|_{L^2_x}^2 + \left\|B_\nu\right\|_{L^2_x}^2 \right)
	% 		\\
	% 		& + \left(\mu \left\|\nabla_x u_\nu\right\|_{L^2_x}^2
	% 		+ \frac 1{\sigma} \left\|j_\nu \right\|_{L^2_x}^2
	% 		+\nu \left\|\nabla_xE_\nu \right\|_{L^2_x}^2
	% 		+\nu \left\|\nabla_xB_\nu \right\|_{L^2_x}^2\right) \leq 0,
	% 	\end{aligned}
	% \end{equation*}
	% in the sense that, for every non-negative $\varphi\in C_c^\infty \left([0,\infty)\right)$,
	\begin{equation*}
		\begin{aligned}
			\frac 12 & \left(\left\|u_\nu\right\|_{L^2_x}^2 + \left\|E_\nu\right\|_{L^2_x}^2 + \left\|B_\nu\right\|_{L^2_x}^2 \right)(t)
			\\
			& \hspace{15mm} + \int_0^t \mu \left\|\nabla_x u_\nu(s)\right\|_{L^2_x}^2 + \frac 1{\sigma} \left\|j_\nu(s)\right\|_{L^2_x}^2
			+\nu \left\|\nabla_xB_\nu(s)\right\|_{L^2_x}^2 ds \\
			& \hspace{15mm} \leq
			\frac 12\left(\left\|u^\mathrm{in}\right\|_{L^2_x}^2
			+ \left\|E^\mathrm{in}\right\|_{L^2_x}^2 + \left\|B^\mathrm{in}\right\|_{L^2_x}^2 \right).
		\end{aligned}
	\end{equation*}
	% \begin{equation*}
	% 	\begin{aligned}
	% 		& - \int_0^\infty \varphi'(t) \frac 12 \left(\left\|u_\nu\right\|_{L^2_x}^2 + \left\|E_\nu\right\|_{L^2_x}^2 + \left\|B_\nu\right\|_{L^2_x}^2 \right)(t)dt
	% 		\\
	% 		& \hspace{10mm} + \int_0^\infty \varphi(t) \left(\mu \left\|\nabla_x u_\nu\right\|_{L^2_x}^2
	% 		+ \frac 1{\sigma} \left\|j_\nu \right\|_{L^2_x}^2
	% 		+\nu \left\|\nabla_xE_\nu \right\|_{L^2_x}^2
	% 		+\nu \left\|\nabla_xB_\nu \right\|_{L^2_x}^2\right)(t) dt \\
	% 		& \hspace{10mm} \leq \varphi(0)
	% 		\frac 12\left(\left\|u^\mathrm{in}\right\|_{L^2_x}^2
	% 		+ \left\|E^\mathrm{in}\right\|_{L^2_x}^2 + \left\|B^\mathrm{in}\right\|_{L^2_x}^2 \right).
	% 	\end{aligned}
	% \end{equation*}
	Such weak solutions are easily established following the method of Leray \cite{leray}, for the nonlinear term $j_\nu\wedge B_\nu$ is now stable with respect to weak convergence in the energy space defined by the above energy inequality, thanks to the dissipation on $B_\nu$.

	Then, repeating the computations of Proposition \ref{modulated energy estimate 0}, it is readily seen that
	\begin{equation*}
		\begin{aligned}
			\frac{d}{dt}\int_{\mathbb{R}^3} & u_\nu\cdot\bar u + E_\nu\cdot \bar E+B_\nu\cdot\bar B dx
			+\int_{\mathbb{R}^3}2\mu \nabla_x u_\nu:\nabla_x\bar u
			+\frac 2\sigma j_\nu\cdot\bar j dx \\
			& =
			- \int_{\mathbb{R}^3}
			\bar u \otimes\left(u_\nu-\bar u\right):\nabla_x(u_\nu-\bar u) dx
			\\
			% +
			% \left(u-\bar u)^{\otimes 2}:\nabla_x\bar u
			& + \int_{\mathbb{R}^3}
			\left((j_\nu-\bar j)\wedge(B_\nu-\bar B)\right)\cdot\bar u
			+ \left((u_\nu-\bar u)\wedge(B_\nu-\bar B)\right)\cdot\bar j
			dx
			\\
			& - \int_{\mathbb{R}^3} \mathbf{A}\cdot
			\begin{pmatrix}
				u_\nu \\ j_\nu
			\end{pmatrix}
			dx
			-\int_{\mathbb{R}^3}\nu \nabla_x B_\nu:\nabla_x\bar B dx.
		\end{aligned}
	\end{equation*}
	Hence, defining the modulated energy $\delta\CE_\nu(t)$ and modulated energy dissipation $\delta\CD_\nu(t)$ by simply replacing $(u,j,E,B)$ by $(u_\nu,j_\nu,E_\nu,B_\nu)$ in \eqref{modulated defi}, we infer that
	\begin{equation*}
		\begin{aligned}
			\delta\CE_\nu(t) & + \int_0^t \delta\CD_\nu(s)
			+ \nu \left\|\nabla_xB_\nu(s)\right\|_{L^2_x}^2 ds
			\\
			& \leq \delta\CE_\nu(0) +
			\int_0^t\int_{\mathbb{R}^3}
			\bar u \otimes\left(u_\nu-\bar u\right):\nabla_x(u_\nu-\bar u) dx ds
			\\
			% +
			% \left(u-\bar u)^{\otimes 2}:\nabla_x\bar u
			& + \int_0^t\int_{\mathbb{R}^3}
			\left(\bar j\wedge(B_\nu-\bar B)\right)\cdot(u_\nu-\bar u)
			- \left((j_\nu-\bar j)\wedge(B_\nu-\bar B)\right)\cdot\bar u
			dxds
			\\
			& + \int_0^t\int_{\mathbb{R}^3} \mathbf{A}\cdot
			\begin{pmatrix}
				u_\nu-\bar u \\ j_\nu - \bar j
			\end{pmatrix}
			+\nu \nabla_x B_\nu:\nabla_x\bar B dxds.
		\end{aligned}
	\end{equation*}
	% \begin{equation*}
	% 	\begin{aligned}
	% 		\frac{d}{dt} & \delta\CE_\nu(t) + \delta\CD_\nu(t)
	% 		+ \nu \left\|\nabla_xB_\nu\right\|_{L^2_x}^2
	% 		\\
	% 		& \leq
	% 		\int_{\mathbb{R}^3}
	% 		\bar u \otimes\left(u_\nu-\bar u\right):\nabla_x(u_\nu-\bar u) dx
	% 		\\
	% 		% +
	% 		% \left(u-\bar u)^{\otimes 2}:\nabla_x\bar u
	% 		& + \int_{\mathbb{R}^3}
	% 		\left(\bar j\wedge(B_\nu-\bar B)\right)\cdot(u_\nu-\bar u)
	% 		- \left((j_\nu-\bar j)\wedge(B_\nu-\bar B)\right)\cdot\bar u
	% 		dx
	% 		\\
	% 		& + \int_{\mathbb{R}^3} \mathbf{A}\cdot
	% 		\begin{pmatrix}
	% 			u_\nu-\bar u \\ j_\nu - \bar j
	% 		\end{pmatrix}
	% 		dx
	% 		+\int_{\mathbb{R}^3}\nu \nabla_x B_\nu:\nabla_x\bar B dx.
	% 	\end{aligned}
	% \end{equation*}
	
	Then, following the proof of Proposition \ref{modulated energy estimate 0}, we arrive at
	\begin{equation*}
		\begin{aligned}
			\delta\CE_\nu(t)
			& +\int_0^t \frac 12\delta\CD_\nu(s)
			+ \nu \left\|\nabla_xB_\nu(s)\right\|_{L^2_x}^2 ds
			\\
			& \leq \delta\CE_\nu(0) +
			\int_0^t \lambda(s)
			\delta\CE_\nu(s)
			+ \left[\int_{\mathbb{R}^3} \mathbf{A}\cdot
			\begin{pmatrix}
				u_\nu-\bar u \\ j_\nu - \bar j
			\end{pmatrix}
			+
			\nu \nabla_x B_\nu:\nabla_x\bar B dx\right]ds
			\\
			& \leq \delta\CE_\nu(0) +
			\int_0^t\lambda(s)
			\delta\CE_\nu(s)
			+ \left[\int_{\mathbb{R}^3} \mathbf{A}\cdot
			\begin{pmatrix}
				u_\nu-\bar u \\ j_\nu - \bar j
			\end{pmatrix}
			dx\right]ds\\
			& + \int_0^t \frac\nu 2 \left\|\nabla_x B_\nu\right\|_{L^2_x}^2 + \frac\nu 2\left\|\nabla_x\bar B\right\|_{L^2_x}^2 ds,
		\end{aligned}
	\end{equation*}
	% \begin{equation*}
	% 	\begin{aligned}
	% 		\frac{d}{dt}\delta\CE_\nu(t)
	% 		& +\frac 12\delta\CD_\nu(t)
	% 		+ \nu \left\|\nabla_xB_\nu\right\|_{L^2_x}^2
	% 		\\
	% 		& \leq
	% 		\lambda(t)
	% 		\delta\CE_\nu(t)
	% 		+ \int_{\mathbb{R}^3} \mathbf{A}\cdot
	% 		\begin{pmatrix}
	% 			u_\nu-\bar u \\ j_\nu - \bar j
	% 		\end{pmatrix}
	% 		+
	% 		\nu \nabla_x B_\nu:\nabla_x\bar B dx
	% 		\\
	% 		& \leq
	% 		\lambda(t)
	% 		\delta\CE_\nu(t)
	% 		+ \int_{\mathbb{R}^3} \mathbf{A}\cdot
	% 		\begin{pmatrix}
	% 			u_\nu-\bar u \\ j_\nu - \bar j
	% 		\end{pmatrix}
	% 		dx
	% 		+ \frac\nu 2 \left\|\nabla_x B_\nu\right\|_{L^2_x}^2 + \frac\nu 2\left\|\nabla_x\bar B\right\|_{L^2_x}^2,
	% 	\end{aligned}
	% \end{equation*}
	and an application of Gr\"onwall's lemma yields
	\begin{equation*}
		\begin{aligned}
			\delta\CE_\nu(t) + & \frac 12 \int_0^t \delta\CD_\nu(s) e^{\int_s^t\lambda(\sigma)d\sigma}ds
			\leq \delta\CE_\nu(0) e^{\int_0^t\lambda(s)ds}
			\\
			& +\int_0^t
			\left[\int_{\mathbb{R}^3} \mathbf{A}\cdot
			\begin{pmatrix}
				u_\nu-\bar u \\ j_\nu - \bar j
			\end{pmatrix}
			dx
			+\frac\nu 2\left\|\nabla_x\bar B\right\|_{L^2_x}^2\right](s)
			e^{\int_s^t\lambda(\sigma)d\sigma}ds.
		\end{aligned}
	\end{equation*}
	
	We may now pass to the limit in the above stability inequality. Thus, up to extraction of subsequences, we may assume that, as $\nu\rightarrow 0$,
	\begin{equation*}
		\begin{aligned}
			u_\nu & \stackrel{*}{\rightharpoonup} u & \text{in }& L^\infty_tL^2_x\cap L^2_t\dot H^1_x,  \\
			j_\nu & \rightharpoonup j & \text{in }& L^2_tL^2_x, \\
			E_\nu & \stackrel{*}{\rightharpoonup} E & \text{in }& L^\infty_t L^2_x, \\
			B_\nu & \stackrel{*}{\rightharpoonup} B & \text{in }& L^\infty_t L^2_x.
		\end{aligned}
	\end{equation*}
	Furthermore, noticing that $\partial_t u_\nu$, $\partial_t E_\nu$ and $\partial_t B_\nu$ are uniformly bounded, in $L^1_\mathrm{loc}$ in time and in some negative index Sobolev space in $x$, it is possible to show (see \cite[Appendix C]{lions7}) that $(u_\nu,E_\nu,B_\nu)$ converges to $(u,E,B)\in C\left([0,\infty);\textit{w-}L^2\left(\mathbb{R}^3\right)\right)$ weakly in $L^2_x$ uniformly locally in time. Then, by the weak lower semi-continuity of the norms, we obtain that, for every $t>0$,
	\begin{equation*}
		\delta\CE(t) + \frac 12 \int_0^t \delta\CD(s) e^{\int_s^t\lambda(\sigma)d\sigma}ds
		\leq
		\liminf_{\nu\rightarrow 0}
		\delta\CE_\nu(t) + \frac 12 \int_0^t \delta\CD_\nu(s) e^{\int_s^t\lambda(\sigma)d\sigma}ds.
	\end{equation*}
	Hence, the stability inequality \eqref{stability 0} holds.
	
	Finally, invoking a classical compactness result by Aubin and Lions \cite{aubin, lions6} (see also \cite{simon} for a sharp compactness criterion), we infer that the $u_\nu$'s converge towards $u$ strongly in $L^2_\mathrm{loc}\left(dtdx\right)$. Therefore, it is readily seen that Ohm's law is satisfied asymptotically, which concludes the proof of the theorem.
\end{proof}

\subsubsection{The two-fluid incompressible Navier-Stokes-Maxwell system with Ohm's law}

Following the strategy of Proposition \ref{modulated energy estimate 0}, the next result establishes a crucial weak-strong stability principle for the two-fluid incompressible Navier-Stokes-Maxwell system with Ohm's law \eqref{TFINSMO}.

\begin{prop}\label{modulated energy estimate 1}
	Let $(u,n,E,B)$ be a smooth solution to the two-fluid incompressible Navier-Stokes-Maxwell system with Ohm's law \eqref{TFINSMO}. Further consider test functions $\left(\bar u, \bar n, \bar j, \bar E, \bar B\right)\in C_c^\infty\left([0,\infty)\times\mathbb{R}^3\right)$ such that
	\begin{equation}\label{test constraints 1}
		\begin{cases}
			\begin{aligned}
				&& \Div \bar u & = 0, \\
				\d_t \bar E - \ROT \bar B &= -  \bar j, & \Div \bar E & = \bar n,
				\\
				\d_t \bar B + \ROT \bar E & = 0, & \Div \bar B & = 0.
			\end{aligned}
		\end{cases}
	\end{equation}
	% \begin{equation*}
	% 	\begin{gathered}
	% 		\DIV \bar u = 0 , \qquad \DIV \bar E = \bar n , \qquad \DIV \bar B = 0 , \\
	% 		\bar j-\bar n\bar u = \sigma\left(-\frac 12 \nabla_x \bar n + \bar E + \bar u\wedge \bar B\right)
	% 	\end{gathered}
	% \end{equation*}
	We define the acceleration operator by
	\begin{equation*}
		\mathbf{A}\left(\bar u, \bar n, \bar j, \bar E, \bar B\right)
		=
		\begin{pmatrix}
			-\d_t \bar u -
			P\left(\bar u\cdot\nabla_x \bar u\right) + \mu\Delta_x \bar u
			+ \frac 12 P \left(\bar n\bar E + \bar j \wedge \bar B\right)
			\\
			-\frac 1{2\sigma}\left(\bar j-\bar n\bar u\right) + \frac 12\left(-\frac 12 \nabla_x \bar n + \bar E + \bar u\wedge \bar B\right)
		\end{pmatrix},
	\end{equation*}
	and the growth rate by
	\begin{equation*}
		\lambda(t)=
		\left(\frac 3{\mu} + \frac 4{\sigma} + 2\sigma\right)
		\left\|\bar u(t)\right\|_{L^\infty_x}^2
		+\frac {3C_0^2}{\mu}\left(\left\|\left(\frac 12\nabla_x\bar n-\bar E\right)(t)\right\|_{L^3_x}^2
		+ \frac 12 \left\|\bar j(t)\right\|_{L^3_x}^2\right),
	\end{equation*}
	where $C_0>0$ denotes the operator norm of the Sobolev embedding $\dot H^1\left(\mathbb{R}^3\right) \hookrightarrow L^6\left(\mathbb{R}^3\right)$.

	Then, one has the stability inequality
	\begin{equation}\label{stability 1}
		\begin{aligned}
			\delta\CE(t) + & \frac 12 \int_0^t \delta\CD(s) e^{\int_s^t\lambda(\sigma)d\sigma}ds
			\\
			& \leq \delta\CE(0) e^{\int_0^t\lambda(s)ds}
			+\int_0^t
			\left[\int_{\mathbb{R}^3} \mathbf{A}\cdot
			\begin{pmatrix}
				u-\bar u \\ j - \bar j - n(u-\bar u)
			\end{pmatrix}
			dx\right](s)
			e^{\int_s^t\lambda(\sigma)d\sigma}ds,
		\end{aligned}
	\end{equation}
	where the modulated energy $\delta\CE$ and energy dissipation $\delta\CD$ are given by
	\begin{equation}\label{modulated defi 1}
		\begin{aligned}
			\delta\CE(t)
			& =
			\frac 12 \left\|\left(u-\bar u\right)(t)\right\|_{L^2_x}^2
			+ \frac 18 \left\|\left(n-\bar n\right)(t)\right\|_{L^2_x}^2 \\
			& + \frac 14 \left\|\left(E-\bar E\right)(t)\right\|_{L^2_x}^2
			+ \frac 14 \left\|\left(B-\bar B\right)(t)\right\|_{L^2_x}^2, \\
			\delta \CD(t)
			& =
			\mu \left\|\nabla_x (u-\bar u)(t)\right\|_{L^2_x}^2
			+ \frac 1{2\sigma} \left\|\left(j-nu-(\bar j-\bar n\bar u)\right)(t)\right\|_{L^2_x}^2.
		\end{aligned}
	\end{equation}
\end{prop}

\begin{proof}
	We have already formally established in Proposition \ref{energy estimate 2} the conservation of the energy for systems \eqref{TFINSMO}. The very same computations applied to the test functions $(\bar u, \bar n, \bar j, \bar E, \bar B)$ yield the identity
	\begin{equation}\label{energy test}
		\frac{d}{dt}\bar\CE(t)+\bar\CD(t)= - \int_{\mathbb{R}^3} \mathbf{A}\cdot
		\begin{pmatrix}
			\bar u \\ \bar j - \bar n \bar u
		\end{pmatrix}
		dx,
	\end{equation}
	where the energy $\bar\CE$ and energy dissipation $\bar\CD$ are obtained simply by replacing the unknowns by the test functions in the respective definitions of Proposition \ref{energy estimate 2}.
	
	Furthermore, other similar duality computations yield that
	\begin{equation*}
		\begin{aligned}
			\frac{d}{dt} & \int_{\mathbb{R}^3} \left( u\cdot\bar u +\frac 14 n\bar n
			+\frac 12 E\cdot \bar E+ \frac 12 B\cdot\bar B \right) dx
			+\int_{\mathbb{R}^3}2\mu \nabla_x u:\nabla_x\bar u dx \\
			& =
			- \int_{\mathbb{R}^3}
			\bar u \otimes (u-\bar u):\nabla_x(u-\bar u) dx
			\\
			& +\frac 12
			\int_{\mathbb{R}^3}
			\left(nE+j\wedge B\right)\cdot \bar u
			+\left(\bar n\bar E+\bar j\wedge \bar B\right)\cdot u dx
			\\
			& -\frac 12\int_{\mathbb{R}^3}
			j\cdot\left(\bar E-\frac 12\nabla_x \bar n\right)
			+
			\bar j\cdot\left(E-\frac 12\nabla_x n\right)
			dx
			- \int_{\mathbb{R}^3} \mathbf{A}\cdot
			\begin{pmatrix}
				u \\ 0
			\end{pmatrix}
			dx,
		\end{aligned}
	\end{equation*}
	and
	\begin{equation*}
		\begin{aligned}
			\frac 1{\sigma} & (j-nu)\cdot(\bar j-\bar n\bar u)
			\\
			& =
			\frac 1{2\sigma}(j-nu+(n-\bar n)\bar u)\cdot(\bar j-\bar n\bar u)
			+\frac 1{2\sigma}
			(n-\bar n)\bar u\cdot
			\left(j-nu - (\bar j-\bar n\bar u)\right)
			\\
			& +\frac 1{2\sigma}(j-nu)\cdot(\bar j-n\bar u)
			\\
			& =
			\frac 1{2}(j-nu+(n-\bar n)\bar u)\cdot\left(-\frac 12\nabla_x\bar n + \bar E + \bar u\wedge\bar B\right)
			- \mathbf{A}\cdot
			\begin{pmatrix}
				0 \\ j-nu+(n-\bar n)\bar u
			\end{pmatrix}
			\\
			& +\frac 1{2\sigma}
			(n-\bar n)\bar u\cdot
			\left(j-nu - (\bar j-\bar n\bar u)\right)
			+\frac 1{2}\left(-\frac 12\nabla_x n + E + u\wedge B\right)\cdot(\bar j-n\bar u),
		\end{aligned}
	\end{equation*}
	whence, considering the sum of the preceding relations,
	\begin{equation*}
		\begin{aligned}
			\frac{d}{dt} & \int_{\mathbb{R}^3} \left( u\cdot\bar u +\frac 14 n\bar n
			+\frac 12 E\cdot \bar E+ \frac 12 B\cdot\bar B \right) dx \\
			& +\int_{\mathbb{R}^3}2\mu \nabla_x u:\nabla_x\bar u
			+\frac 1{\sigma} (j-nu)\cdot(\bar j-\bar n\bar u) dx \\
			& =
			- \int_{\mathbb{R}^3}
			\bar u \otimes (u-\bar u):\nabla_x(u-\bar u) dx
			+\frac 12\int_{\mathbb{R}^3} (n-\bar n)(u-\bar u)\cdot\left(\frac 12\nabla_x\bar n-\bar E\right)dx
			\\
			& +\frac 12 \int_{\mathbb{R}^3}
			\left((j-nu-(\bar j-\bar n\bar u))\wedge(B-\bar B)\right)\cdot\bar u
			+ \left((u-\bar u)\wedge(B-\bar B)\right)\cdot\bar j
			dx
			\\
			& +\frac 1{2\sigma} \int_{\mathbb{R}^3}
			(n-\bar n)(j-nu-(\bar j-\bar n\bar u))\cdot\bar udx
			- \int_{\mathbb{R}^3} \mathbf{A}\cdot
			\begin{pmatrix}
				u \\ j-nu+(n-\bar n)\bar u
			\end{pmatrix}
			dx.
		\end{aligned}
	\end{equation*}

	On the whole, combining the above identities with the energy decay imposed by the formal energy conservations from Proposition \ref{energy estimate 2}, we find the following modulated energy inequality~:
	\begin{equation*}
		\begin{aligned}
			\frac{d}{dt} & \delta\CE(t) + \delta\CD(t) \\
			& \leq
			\int_{\mathbb{R}^3}
			\bar u \otimes(u-\bar u):\nabla_x(u-\bar u) dx
			-\frac 12\int_{\mathbb{R}^3} (n-\bar n)(u-\bar u)\cdot\left(\frac 12\nabla_x\bar n-\bar E\right)dx
			\\
			& -\frac 12 \int_{\mathbb{R}^3}
			\left((j-nu-(\bar j-\bar n\bar u))\wedge(B-\bar B)\right)\cdot\bar u
			+ \left((u-\bar u)\wedge(B-\bar B)\right)\cdot\bar j
			dx
			\\
			& -\frac 1{2\sigma} \int_{\mathbb{R}^3}
			(n-\bar n)(j-nu-(\bar j-\bar n\bar u))\cdot\bar udx
			+
			\int_{\mathbb{R}^3} \mathbf{A}\cdot
			\begin{pmatrix}
				u-\bar u \\ j-\bar j -n(u- \bar u)
			\end{pmatrix}
			dx.
		\end{aligned}
	\end{equation*}
	The next step consists in estimating the terms in the right-hand side above that are nonlinear in $(u,n,j,E,B)$ and to absorb the resulting expressions with the modulated energy $\delta\CE(t)$ and the modulated energy dissipation $\delta\CD(t)$ by suitable uses of Young's inequality and Gr\"onwall's lemma. Thus, we obtain
	\begin{equation*}
		\begin{aligned}
			\frac{d}{dt} & \delta\CE(t) + \delta\CD(t) \\
			& \leq
			\left\|\bar u\right\|_{L^\infty_x} \left\| u-\bar u\right\|_{L^2_x}\left\|\nabla_x(u-\bar u)\right\|_{L^2_x}
			+\frac 12\left\|\frac 12\nabla_x\bar n-\bar E\right\|_{L^3_x}\left\|n-\bar n\right\|_{L^2_x}\left\|u-\bar u\right\|_{L^6_x}
			\\
			& +\frac 1{2\sigma}\left\|\bar u\right\|_{L^\infty_x}
			\left(\left\|n-\bar n\right\|_{L^2_x}+\sigma\left\|B-\bar B\right\|_{L^2_x}\right)
			\left\|j-nu-(\bar j-\bar n\bar u)\right\|_{L^2_x}
			\\
			& + \frac 12 \left\|\bar j\right\|_{L^3_x}\left\|B-\bar B\right\|_{L^2_x}\left\|u-\bar u\right\|_{L^6_x}
			+
			\int_{\mathbb{R}^3} \mathbf{A}\cdot
			\begin{pmatrix}
				u-\bar u \\ j-\bar j -n(u- \bar u)
			\end{pmatrix}
			dx
			\\
			& \leq
			\frac 3{2\mu}\left\|\bar u\right\|_{L^\infty_x}^2 \left\| u-\bar u\right\|_{L^2_x}^2
			+\left(\frac {3C_0^2}{8\mu}\left\|\frac 12\nabla_x\bar n-\bar E\right\|_{L^3_x}^2
			+ \frac 1{2\sigma}\left\|\bar u\right\|_{L^\infty_x}^2\right)
			\left\|n-\bar n\right\|_{L^2_x}^2
			\\
			& +\left(\frac \sigma{2}\left\|\bar u\right\|_{L^\infty_x}^2 + \frac{3C_0^2}{8\mu} \left\|\bar j\right\|_{L^3_x}^2\right)
			\left\|B-\bar B\right\|_{L^2_x}^2
			\\
			& + \frac\mu 2\left\|\nabla_x(u-\bar u)\right\|_{L^2_x}^2
			+ \frac 1{4\sigma}\left\|j-nu-(\bar j-\bar n\bar u)\right\|_{L^2_x}^2
			\\
			& +
			\int_{\mathbb{R}^3} \mathbf{A}\cdot
			\begin{pmatrix}
				u-\bar u \\ j-\bar j -n(u- \bar u)
			\end{pmatrix}
			dx.
		\end{aligned}
	\end{equation*}
	Hence,
	\begin{equation*}
		\frac{d}{dt}\delta\CE(t)
		+\frac 12\delta\CD(t)
		\leq
		\lambda(t)
		\delta\CE(t)
		+ \int_{\mathbb{R}^3} \mathbf{A}\cdot
		\begin{pmatrix}
			u-\bar u \\ j - \bar j -n(u- \bar u)
		\end{pmatrix}
		dx,
	\end{equation*}
	which concludes the proof of the proposition with a direct application of Gr\"onwall's lemma.
\end{proof}

Again, note that the test functions satisfying the linear constraints \eqref{test constraints 1} are easily constructed employing the relations \eqref{test construction}. Now, one may prefer to deal, in a completely equivalent manner, with test functions $\left(\bar u, \bar n, \bar j, \bar E, \bar B\right)\in C_c^\infty\left([0,\infty)\times\mathbb{R}^3\right)$ satisfying the stationary constraints
\begin{equation*}
	\Div \bar u = 0,
	\qquad
	\Div \bar E = \bar n,
	\qquad
	\Div \bar B = 0,
	\qquad
	\bar j - \bar n\bar u = \sigma\left( -\frac 12\nabla_x \bar n + \bar E + \bar u\wedge \bar B\right),
\end{equation*}
rather than the constraints \eqref{test constraints 1}. In this case, instead of \eqref{stability 1}, we obtain the stability inequality
\begin{equation*}
	\begin{aligned}
		& \delta\CE(t) + \frac 12 \int_0^t \delta\CD(s) e^{\int_s^t\lambda(\sigma)d\sigma}ds
		\\
		& \leq \delta\CE(0) e^{\int_0^t\lambda(s)ds}
		+\int_0^t
		\left[\int_{\mathbb{R}^3} \mathbf{A}\cdot
		\begin{pmatrix}
			u-\bar u \\ E -\frac 12\nabla_x n - \left(\bar E-\frac 12\nabla_x\bar n\right) \\ B - \bar B
		\end{pmatrix}
		dx\right](s)
		e^{\int_s^t\lambda(\sigma)d\sigma}ds,
	\end{aligned}
\end{equation*}
where the acceleration operator is now defined by
\begin{equation*}
	\mathbf{A}\left(\bar u, \bar n, \bar j, \bar E, \bar B\right)
	=
	\begin{pmatrix}
		-\d_t \bar u -
		P\left(\bar u\cdot\nabla_x \bar u\right) + \mu\Delta_x \bar u
		+ \frac 12 P \left(\bar n\bar E + \bar j \wedge \bar B\right)
		\\
		\frac 12\left(-\partial_t \bar E + \ROT\bar B - \bar j\right)
		\\
		\frac 12\left(-\partial_t \bar B - \ROT\bar E\right)
	\end{pmatrix}.
\end{equation*}

The preceding proposition provides an important weak-strong stability property for the two-fluid incompressible Navier-Stokes-Maxwell system with Ohm's law \eqref{TFINSMO}. Indeed, the stability inequality \eqref{stability 1} essentially implies that a solution $(\bar u, \bar n, \bar j,\bar E,\bar B)$ of \eqref{TFINSMO} such that $\bar u\in L^2_tL^\infty_x$, $\bar j\in L^2_tL^3_x$ and $\frac 12\nabla_x\bar n-\bar E\in L^2_tL^3_x$, if it exists, is unique in the whole class of weak solutions in the energy space, for any given initial data.

\bigskip

By analogy with Lions' dissipative solutions to the incompressible Euler system \cite[Section 4.4]{lions7}, we provide now a suitable notion of dissipative solution for the two-fluid incompressible Navier-Stokes-Maxwell system with Ohm's law \eqref{TFINSMO}, based on Proposition \ref{modulated energy estimate 1}, and establish their existence next.

\begin{defi}
	We say that
	\begin{equation*}
		(u,n,E,B)\in L^\infty\left([0,\infty);L^2\left(\mathbb{R}^3\right)\right)
		\cap C\left([0,\infty);\textit{w-}L^2\left(\mathbb{R}^3\right)\right)
	\end{equation*}
	such that
	\begin{equation*}
		\DIV u = 0, \qquad \DIV E = n, \qquad \DIV B = 0,
	\end{equation*}
	is a \textbf{dissipative solution of the two-fluid incompressible Navier-Stokes-Maxwell system with Ohm's law \eqref{TFINSMO}}, if it solves Maxwell's equations
	\begin{equation*}
		\begin{cases}
			\begin{aligned}
				\d_t E - \ROT B &= -  j, 
				\\
				\d_t B + \ROT E & = 0,
			\end{aligned}
		\end{cases}
	\end{equation*}
	with Ohm's law
	\begin{equation*}
		j - nu = \sigma\left( -\frac 12\nabla_x n + E + u\wedge B\right),
	\end{equation*}
	in the sense of distributions, and if, for any test functions $\left(\bar u, \bar n, \bar j, \bar E, \bar B\right)\in C_c^\infty\left([0,\infty)\times\mathbb{R}^3\right)$ satisfying the linear constraints \eqref{test constraints 1}, the stability inequality \eqref{stability 1} is verified.
\end{defi}

As previously mentioned, dissipative solutions define actual solutions in the sense that they coincide with the unique strong solution when the latter exists. The following theorem asserts their existence.

\begin{thm}\label{dissipative 2}
	For any initial data $\left(u^\mathrm{in}, n^\mathrm{in} ,E^\mathrm{in},B^\mathrm{in}\right)\in L^2\left(\mathbb{R}^3\right)$ such that
	\begin{equation*}
		\DIV u^\mathrm{in} = 0, \qquad \DIV E^\mathrm{in} = n^\mathrm{in}, \qquad \DIV B^\mathrm{in} = 0,
	\end{equation*}
	there exists a dissipative solution to the two-fluid incompressible Navier-Stokes-Maxwell system with Ohm's law \eqref{TFINSMO}.
\end{thm}

\begin{proof}
	As in the proof of Theorem \ref{dissipative 1}, it is possible, here, to justify the existence of dissipative solutions by introducing viscous approximations of the system \eqref{TFINSMO}. Thus, for each $\nu>0$, we consider weak solutions of the following system~:
	\begin{equation}\label{viscous approximation}
		\begin{cases}
			\begin{aligned}
				\d_t u_\nu +
				u_\nu\cdot\nabla_x u_\nu - \mu\Delta_x u_\nu
				& = -\nabla_x p_\nu+
				\frac 12 \left(n_\nu E_\nu + j_\nu \wedge B_\nu\right) , & \Div u_\nu & = 0,\\
				\d_t E_\nu - \ROT B_\nu-\nu\Delta_x E_\nu &= -  j_\nu, & \Div E_\nu & = n_\nu,
				\\
				\d_t B_\nu + \ROT E_\nu-\nu\Delta_x B_\nu & = 0, & \Div B_\nu & = 0, \\
				j_\nu-n_\nu u_\nu & = \sigma\left(-\frac 12 \nabla_x n_\nu + E_\nu + u_\nu\wedge B_\nu\right), &&
			\end{aligned}
		\end{cases}
	\end{equation}
	associated with the initial data $\left(u^\mathrm{in},n^\mathrm{in},E^\mathrm{in},B^\mathrm{in}\right)$ and satisfying the energy inequality, for all $t>0$,
	% \begin{equation*}
	% 	\begin{aligned}
	% 		& \frac{d}{dt}\left(\frac 12\left\|u_\nu\right\|_{L^2_x}^2+\frac 18 \left\|n_\nu\right\|_{L^2_x}^2 + \frac 14\left\|E_\nu\right\|_{L^2_x}^2 + \frac 14 \left\|B_\nu\right\|_{L^2_x}^2 \right)
	% 		\\
	% 		& \hspace{30mm}+ \mu \left\|\nabla_x u_\nu\right\|_{L^2_x}^2 + \frac 1{2\sigma} \left\|j_\nu-n_\nu u_\nu\right\|_{L^2_x}^2
	% 		\\
	% 		& \hspace{30mm} +\frac\nu 4 \left\|\nabla_xn_\nu\right\|_{L^2_x}^2
	% 		+\frac\nu 2 \left\|\nabla_xE_\nu\right\|_{L^2_x}^2
	% 		+\frac\nu 2 \left\|\nabla_xB_\nu\right\|_{L^2_x}^2 \leq 0,
	% 	\end{aligned}
	% \end{equation*}
	% in the sense that, for every non-negative $\varphi\in C_c^\infty \left([0,\infty)\right)$,
	\begin{equation*}
		\begin{aligned}
			& \left(\frac 12\left\|u_\nu\right\|_{L^2_x}^2+\frac 18 \left\|n_\nu\right\|_{L^2_x}^2 + \frac 14\left\|E_\nu\right\|_{L^2_x}^2 + \frac 14 \left\|B_\nu\right\|_{L^2_x}^2 \right)(t)
			\\
			& \hspace{12mm} + \int_0^t \bigg( \mu \left\|\nabla_x u_\nu\right\|_{L^2_x}^2 + \frac 1{2\sigma} \left\|j_\nu-n_\nu u_\nu\right\|_{L^2_x}^2
			\\
			& \hspace{12mm}
			+\frac\nu 4 \left\|\nabla_xn_\nu\right\|_{L^2_x}^2
			+\frac\nu 2 \left\|\nabla_xE_\nu\right\|_{L^2_x}^2
			+\frac\nu 2 \left\|\nabla_xB_\nu\right\|_{L^2_x}^2 \bigg)(s) ds \\
			& \hspace{12mm} \leq
			\frac 12\left\|u^\mathrm{in}\right\|_{L^2_x}^2+\frac 18 \left\|n^\mathrm{in}\right\|_{L^2_x}^2 + \frac 14\left\|E^\mathrm{in}\right\|_{L^2_x}^2 + \frac 14 \left\|B^\mathrm{in}\right\|_{L^2_x}^2.
		\end{aligned}
	\end{equation*}
	Such weak solutions are easily established following the method of Leray \cite{leray}, for the nonlinear terms $n_\nu E_\nu$ and $j_\nu\wedge B_\nu$ are now stable with respect to weak convergence in the energy space defined by the above energy inequality, thanks to the dissipation on $n_\nu$, $E_\nu$ and $B_\nu$.

	Then, repeating the computations of Proposition \ref{modulated energy estimate 1}, it is readily seen that
	\begin{equation*}
		\begin{aligned}
			& \frac{d}{dt} \int_{\mathbb{R}^3} \left( u_\nu \cdot\bar u +\frac 14 n_\nu \bar n
			+\frac 12 E_\nu \cdot \bar E+ \frac 12 B_\nu \cdot\bar B \right) dx \\
			& +\int_{\mathbb{R}^3}2\mu \nabla_x u_\nu :\nabla_x\bar u
			+\frac 1{\sigma} (j_\nu -n_\nu u_\nu )\cdot(\bar j-\bar n\bar u) dx \\
			& =
			- \int_{\mathbb{R}^3}
			\bar u \otimes (u_\nu -\bar u):\nabla_x(u_\nu -\bar u) dx
			+\frac 12\int_{\mathbb{R}^3} (n_\nu -\bar n)(u_\nu -\bar u)\cdot\left(\frac 12\nabla_x\bar n-\bar E\right)dx
			\\
			& +\frac 12 \int_{\mathbb{R}^3}
			\left((j_\nu -n_\nu u_\nu -(\bar j-\bar n\bar u))\wedge(B_\nu -\bar B)\right)\cdot\bar u
			+ \left((u_\nu -\bar u)\wedge(B_\nu -\bar B)\right)\cdot\bar j
			dx
			\\
			& +\frac 1{2\sigma} \int_{\mathbb{R}^3}
			(n_\nu -\bar n)(j_\nu -n_\nu u_\nu -(\bar j-\bar n\bar u))\cdot\bar udx
			\\
			& - \int_{\mathbb{R}^3} \mathbf{A}\cdot
			\begin{pmatrix}
				u_\nu  \\ j_\nu -n_\nu u_\nu +(n_\nu -\bar n)\bar u
			\end{pmatrix}
			dx
			\\
			& -\int_{\mathbb{R}^3}\frac\nu 4 \nabla_x n_\nu:\nabla_x\bar n
			+
			\frac\nu 2 \nabla_x E_\nu:\nabla_x\bar E
			+
			\frac \nu 2 \nabla_x B_\nu:\nabla_x\bar B dx.
		\end{aligned}
	\end{equation*}
	Hence, defining the modulated energy $\delta\CE_\nu(t)$ and modulated energy dissipation $\delta\CD_\nu(t)$ by simply replacing $(u,n,j,E,B)$ by $(u_\nu,n_\nu,j_\nu,E_\nu,B_\nu)$ in \eqref{modulated defi 1}, we infer that
	\begin{equation*}
		\begin{aligned}
			\delta\CE_\nu(t) & + \int_0^t \delta\CD_\nu(s)
			+\frac\nu 2
			\left(\frac 1 2 \left\|\nabla_xn_\nu\right\|_{L^2_x}^2
			+ \left\|\nabla_xE_\nu\right\|_{L^2_x}^2
			+ \left\|\nabla_xB_\nu\right\|_{L^2_x}^2\right)(s)
			ds
			\\
			& \leq \delta\CE_\nu(0) +
			\int_0^t\int_{\mathbb{R}^3}
			\bar u \otimes (u_\nu -\bar u):\nabla_x(u_\nu -\bar u) dxds
			\\
			&-\frac 12\int_0^t\int_{\mathbb{R}^3} (n_\nu -\bar n)(u_\nu -\bar u)\cdot\left(\frac 12\nabla_x\bar n-\bar E\right)dxds
			\\
			& -\frac 12 \int_0^t \int_{\mathbb{R}^3}
			\left((j_\nu -n_\nu u_\nu -(\bar j-\bar n\bar u))\wedge(B_\nu -\bar B)\right)\cdot\bar u
			dxds
			\\
			&
			-\frac 12 \int_0^t\int_{\mathbb{R}^3}
			\left((u_\nu -\bar u)\wedge(B_\nu -\bar B)\right)\cdot\bar j
			dxds
			\\
			& -\frac 1{2\sigma} \int_0^t \int_{\mathbb{R}^3}
			(n_\nu -\bar n)(j_\nu -n_\nu u_\nu -(\bar j-\bar n\bar u))\cdot\bar udxds
			\\
			& + \int_0^t \int_{\mathbb{R}^3} \mathbf{A}\cdot
			\begin{pmatrix}
				u_\nu -\bar u  \\ j_\nu-\bar j -n_\nu (u_\nu-\bar u)
			\end{pmatrix}
			dxds
			\\
			& +\int_0^t\int_{\mathbb{R}^3}\frac\nu 4 \nabla_x n_\nu:\nabla_x\bar n
			+
			\frac\nu 2 \nabla_x E_\nu:\nabla_x\bar E
			+
			\frac \nu 2 \nabla_x B_\nu:\nabla_x\bar B dxds.
		\end{aligned}
	\end{equation*}
	
	Then, following the proof of Proposition \ref{modulated energy estimate 1}, we arrive at
	\begin{equation*}
		\begin{aligned}
			\delta\CE_\nu(t)
			& +\int_0^t \frac 12\delta\CD_\nu(s)
			+\frac \nu 2\left(\frac 1 2 \left\|\nabla_xn_\nu\right\|_{L^2_x}^2
			+ \left\|\nabla_xE_\nu\right\|_{L^2_x}^2
			+ \left\|\nabla_xB_\nu\right\|_{L^2_x}^2 \right)(s)ds
			\\
			& \leq \delta\CE_\nu(0) +
			\int_0^t
			\lambda(s)
			\delta\CE_\nu(s)
			+ \left[ \int_{\mathbb{R}^3} \mathbf{A}\cdot
			\begin{pmatrix}
				u_\nu-\bar u \\ j_\nu-\bar j -n_\nu (u_\nu-\bar u)
			\end{pmatrix}
			dx \right] ds
			\\
			& +\frac\nu 2\int_0^t\int_{\mathbb{R}^3}\frac 12 \nabla_x n_\nu:\nabla_x\bar n
			+
			\nabla_x E_\nu:\nabla_x\bar E
			+
			\nabla_x B_\nu:\nabla_x\bar B dxds
			\\
			& \leq \delta\CE_\nu(0) +
			\int_0^t
			\lambda(s)
			\delta\CE_\nu(s)
			+ \left[ \int_{\mathbb{R}^3} \mathbf{A}\cdot
			\begin{pmatrix}
				u_\nu-\bar u \\ j_\nu-\bar j -n_\nu (u_\nu-\bar u)
			\end{pmatrix}
			dx \right] ds
			\\
			& +\frac\nu 4 \int_0^t \frac 12 \left\|\nabla_xn_\nu\right\|_{L^2_x}^2
			+ \left\|\nabla_xE_\nu\right\|_{L^2_x}^2
			+ \left\|\nabla_xB_\nu\right\|_{L^2_x}^2 ds
			\\
			& +\frac\nu 4 \int_0^t \frac 12 \left\|\nabla_x\bar n\right\|_{L^2_x}^2
			+ \left\|\nabla_x\bar E\right\|_{L^2_x}^2
			+ \left\|\nabla_x\bar B\right\|_{L^2_x}^2 ds,
		\end{aligned}
	\end{equation*}
	and an application of Gr\"onwall's lemma yields
	\begin{equation*}
		\begin{aligned}
			\delta\CE_\nu(t) & + \frac 12 \int_0^t \delta\CD_\nu(s) e^{\int_s^t\lambda(\sigma)d\sigma}ds
			\leq \delta\CE_\nu(0) e^{\int_0^t\lambda(s)ds}
			\\
			& +\int_0^t
			\left[\int_{\mathbb{R}^3} \mathbf{A}\cdot
			\begin{pmatrix}
				u_\nu-\bar u \\ j_\nu-\bar j -n_\nu (u_\nu-\bar u)
			\end{pmatrix}
			dx\right](s)e^{\int_s^t\lambda(\sigma)d\sigma}ds
			\\
			& + \frac\nu 4 \int_0^t \left[\frac 12 \left\|\nabla_x\bar n\right\|_{L^2_x}^2
			+ \left\|\nabla_x\bar E\right\|_{L^2_x}^2
			+ \left\|\nabla_x\bar B\right\|_{L^2_x}^2
			\right](s)
			e^{\int_s^t\lambda(\sigma)d\sigma}ds.
		\end{aligned}
	\end{equation*}
	
	We may now pass to the limit in the above stability inequality. Thus, up to extraction of subsequences, we may assume that, as $\nu\rightarrow 0$,
	\begin{equation*}
		\begin{aligned}
			u_\nu & \stackrel{*}{\rightharpoonup} u & \text{in }& L^\infty_tL^2_x\cap L^2_t\dot H^1_x,  \\
			n_\nu & \stackrel{*}{\rightharpoonup} n & \text{in }& L^\infty_t L^2_x, \\
			E_\nu & \stackrel{*}{\rightharpoonup} E & \text{in }& L^\infty_t L^2_x, \\
			B_\nu & \stackrel{*}{\rightharpoonup} B & \text{in }& L^\infty_t L^2_x.
		\end{aligned}
	\end{equation*}
	Furthermore, noticing that $\partial_t u_\nu$, $\partial_t n_\nu$, $\partial_t E_\nu$ and $\partial_t B_\nu$ are uniformly bounded, in $L^1_\mathrm{loc}$ in time and in some negative index Sobolev space in $x$, it is possible to show (see \cite[Appendix C]{lions7}) that $(u_\nu, n_\nu,E_\nu,B_\nu)$ converges to $(u,n,E,B)\in C\left([0,\infty);\textit{w-}L^2\left(\mathbb{R}^3\right)\right)$ weakly in $L^2_x$ uniformly locally in time. Moreover, invoking a classical compactness result by Aubin and Lions \cite{aubin, lions6} (see also \cite{simon} for a sharp compactness criterion), we infer that the $u_\nu$'s converge towards $u$ strongly in $L^2_\mathrm{loc}\left(dtdx\right)$. In particular, it follows that, up to extraction,
	\begin{equation*}
		\begin{aligned}
			n_\nu u_\nu & \stackrel{*}{\rightharpoonup} nu & \text{in }& L^\infty_tL^1_x\cap L^2_tL^6_x,  \\
			j_\nu-n_\nu u_\nu & \rightharpoonup j-nu & \text{in }& L^2_tL^2_x.
		\end{aligned}
	\end{equation*}

	Then, by the weak lower semi-continuity of the norms, we obtain that, for every $t>0$,
	\begin{equation*}
		\delta\CE(t) + \frac 12 \int_0^t \delta\CD(s) e^{\int_s^t\lambda(\sigma)d\sigma}ds
		\leq
		\liminf_{\nu\rightarrow 0}
		\delta\CE_\nu(t) + \frac 12 \int_0^t \delta\CD_\nu(s) e^{\int_s^t\lambda(\sigma)d\sigma}ds.
	\end{equation*}
	Hence, the stability inequality \eqref{stability 1} holds. Finally, it is readily seen that Ohm's law is satisfied asymptotically, which concludes the proof of the theorem.
\end{proof}

We present now an alternative kind of stability inequality for the two-fluid incompressible Navier-Stokes-Maxwell system with Ohm's law \eqref{TFINSMO}, whose understanding will be crucial for the relative entropy method --~developed later on in Chapter \ref{entropy method}~-- in the hydrodynamic limit of the two species Vlasov-Maxwell-Boltzmann system \eqref{scaled VMB two species}. It is based on the identity \eqref{poynting} linking the Lorentz force with the Poynting vector $E\wedge B$, which will allow us to stabilize the modulated nonlinear terms solely with the modulated energy $\delta\mathcal{E}$ (i.e.\ without absorbing nonlinear terms with the modulated dissipation $\delta\mathcal{D}$~; note the different coefficient in front of $\delta\mathcal{D}$ in the stability inequalities \eqref{stability 1} and \eqref{stability 3}, below).

\begin{prop}\label{stability poynting}
	Let $(u,n,E,B)$ be a smooth solution to the two-fluid incompressible Navier-Stokes-Maxwell system with Ohm's law \eqref{TFINSMO}. Further consider test functions $\left(\bar u, \bar n, \bar j, \bar E, \bar B\right)\in C_c^\infty\left([0,\infty)\times\mathbb{R}^3\right)$ such that $\left\|\bar u\right\|_{L^\infty_{t,x}}<1$ and
	\begin{equation*}
		\begin{cases}
			\begin{aligned}
				&& \Div \bar u & = 0, \\
				\d_t \bar E - \ROT \bar B &= -  \bar j, & \Div \bar E & = \bar n,
				\\
				\d_t \bar B + \ROT \bar E & = 0, & \Div \bar B & = 0.
			\end{aligned}
		\end{cases}
	\end{equation*}
	We define the acceleration operator by
	\begin{equation*}
		\mathbf{A}\left(\bar u, \bar n, \bar j, \bar E, \bar B\right)
		=
		\begin{pmatrix}
			-\d_t \bar u -
			P\left(\bar u\cdot\nabla_x \bar u\right) + \mu\Delta_x \bar u
			+ \frac 12 P \left(\bar n\bar E + \bar j \wedge \bar B\right)
			\\
			-\frac 1{2\sigma}\left(\bar j-\bar n\bar u\right) + \frac 12\left(-\frac 12 \nabla_x \bar n + \bar E + \bar u\wedge \bar B\right)
		\end{pmatrix},
	\end{equation*}
	and the growth rate by
	\begin{equation*}
		\lambda(t) =
		\frac{2\left\|\nabla_{t,x}\bar u(t)\right\|_{L^\infty_x}}{1-\left\|\bar u(t)\right\|_{L^\infty_x}}
		+ \frac {\sqrt 2 \left\|\left(\bar j-\bar n\bar u\right)(t)\right\|_{L^\infty_x}}{2\left(1-\left\|\bar u(t)\right\|_{L^\infty_x}\right)}
		+\left\|\left(\frac 12\nabla_x\bar n-\bar E-\bar u\wedge\bar B\right)(t)\right\|_{L^\infty_x}.
	\end{equation*}
	
	Then, one has the stability inequality
	\begin{equation}\label{stability 3}
		\begin{aligned}
			\delta\CE(t) + & \int_0^t \delta\CD(s) e^{\int_s^t\lambda(\sigma)d\sigma}ds
			\\
			& \leq \delta\CE(0) e^{\int_0^t\lambda(s)ds}
			+\int_0^t
			\left[\int_{\mathbb{R}^3} \mathbf{A}\cdot
			\begin{pmatrix}
				u-\bar u \\ j-nu-(\bar j-\bar n\bar u)
			\end{pmatrix}
			dx\right](s)
			e^{\int_s^t\lambda(\sigma)d\sigma}ds,
		\end{aligned}
	\end{equation}
	where the modulated energy $\delta\CE$ and energy dissipation $\delta\CD$ are given by
	\begin{equation}\label{modulated defi 2}
		\begin{aligned}
			\delta\CE(t)
			& =
			\frac 12 \left\|\left(u-\bar u\right)(t)\right\|_{L^2_x}^2
			+ \frac 18 \left\|\left(n-\bar n\right)(t)\right\|_{L^2_x}^2 \\
			& + \frac 14 \left\|\left(E-\bar E\right)(t)\right\|_{L^2_x}^2
			+ \frac 14 \left\|\left(B-\bar B\right)(t)\right\|_{L^2_x}^2 \\
			& - \frac 12\int_{\mathbb{R}^3} \left(\left(E-\bar E\right)(t)\wedge \left(B-\bar B\right)(t)\right)\cdot\bar u(t) dx, \\
			\delta \CD(t)
			& =
			\mu \left\|\nabla_x (u-\bar u)(t)\right\|_{L^2_x}^2
			+ \frac 1{2\sigma} \left\|\left(j-nu-(\bar j-\bar n\bar u)\right)(t)\right\|_{L^2_x}^2.
		\end{aligned}
	\end{equation}
\end{prop}

\begin{proof}
	Following the proof of Proposition \ref{modulated energy estimate 1}, using that $\Div u = \Div \bar u = 0$, we consider first the identity
	\begin{equation*}
		\begin{aligned}
			\frac{d}{dt} & \int_{\mathbb{R}^3}\left( u\cdot\bar u +\frac 14 n\bar n
			+\frac 12 E\cdot \bar E+ \frac 12 B\cdot\bar B \right) dx
			+\int_{\mathbb{R}^3}2\mu \nabla_x u:\nabla_x\bar u dx \\
			& =
			\int_{\mathbb{R}^3}
			(u-\bar u) \otimes (u-\bar u):\nabla_x\bar u dx
			\\
			& +\frac 12
			\int_{\mathbb{R}^3}
			\left(nE+j\wedge B\right)\cdot \bar u
			+\left(\bar n\bar E+\bar j\wedge \bar B\right)\cdot u dx
			\\
			& -\frac 12\int_{\mathbb{R}^3}
			j\cdot\left(\bar E-\frac 12\nabla_x \bar n\right)
			+
			\bar j\cdot\left(E-\frac 12\nabla_x n\right)
			dx
			- \int_{\mathbb{R}^3} \mathbf{A}\cdot
			\begin{pmatrix}
				u \\ 0
			\end{pmatrix}
			dx
			\\
			& =
			\int_{\mathbb{R}^3}
			(u-\bar u) \otimes (u-\bar u):\nabla_x\bar u dx
			\\
			& +\frac 12
			\int_{\mathbb{R}^3}
			\left(\left(n-\bar n\right)\left(E-\bar E\right)+\left(j-\bar j\right)\wedge \left(B-\bar B\right)\right)\cdot \bar u
			dx
			\\
			& +\frac 12\int_{\mathbb{R}^3}
			\left(\bar j-\bar n\bar u\right)\cdot\left(\left(u-\bar u\right)\wedge\left(B-\bar B\right)\right)
			dx
			\\
			& +\frac 12\int_{\mathbb{R}^3}
			\left(n-\bar n\right)\left(u-\bar u\right)\cdot\left(\frac 12\nabla_x\bar n-\bar E - \bar u\wedge\bar B\right)
			dx
			\\
			& +\frac 12\int_{\mathbb{R}^3}
			\left(j-nu\right)\cdot\left(\frac 12\nabla_x\bar n - \bar E - \bar u\wedge\bar B\right)
			+
			\left(\bar j-\bar n\bar u\right)\cdot\left(\frac 12\nabla_x n-E-u\wedge B\right)
			dx
			\\
			& - \int_{\mathbb{R}^3} \mathbf{A}\cdot
			\begin{pmatrix}
				u \\ 0
			\end{pmatrix}
			dx.
		\end{aligned}
	\end{equation*}
	Further using Ohm's laws, we find
	\begin{equation}\label{computation}
		\begin{aligned}
			\frac{d}{dt} \int_{\mathbb{R}^3} & \left(u\cdot\bar u +\frac 14 n\bar n
			+\frac 12 E\cdot \bar E+ \frac 12 B\cdot\bar B\right) dx
			\\
			& +\int_{\mathbb{R}^3}2\mu \nabla_x u:\nabla_x\bar u
			+\frac 1{\sigma}
			\left(j-nu\right)\cdot\left(\bar j-\bar n\bar u\right)
			dx
			\\
			& =
			\int_{\mathbb{R}^3}
			(u-\bar u) \otimes (u-\bar u):\nabla_x\bar u dx
			\\
			& +\frac 12
			\int_{\mathbb{R}^3}
			\left(\left(n-\bar n\right)\left(E-\bar E\right)+\left(j-\bar j\right)\wedge \left(B-\bar B\right)\right)\cdot \bar u
			dx
			\\
			& +\frac 12\int_{\mathbb{R}^3}
			\left(\bar j-\bar n\bar u\right)\cdot\left(\left(u-\bar u\right)\wedge\left(B-\bar B\right)\right)
			dx
			\\
			& +\frac 12\int_{\mathbb{R}^3}
			\left(n-\bar n\right)\left(u-\bar u\right)\cdot\left(\frac 12\nabla_x\bar n-\bar E - \bar u\wedge\bar B\right)
			dx
			\\
			& - \int_{\mathbb{R}^3} \mathbf{A}\cdot
			\begin{pmatrix}
				u \\ j-nu
			\end{pmatrix}
			dx.
		\end{aligned}
	\end{equation}

	Then, expressing the modulated Lorentz force with a modulated Poynting vector as
	\begin{equation*}
		\begin{aligned}
			\partial_t\left(\left(E-\bar E\right)\wedge \left(B-\bar B\right)\right)
			&
			+\frac 12\nabla_x\left(\left|E-\bar E\right|^2+\left|B-\bar B\right|^2\right)
			\\
			& -\nabla_x\cdot\left(\left(E-\bar E\right)\otimes \left(E-\bar E\right)+\left(B-\bar B\right)\otimes \left(B-\bar B\right)\right)
			\\
			& =
			-\left(n-\bar n\right)\left(E-\bar E\right)-\left(j-\bar j\right)\wedge \left(B-\bar B\right),
			% \\
			% \partial_t\left(\bar E\wedge \bar B\right)+\frac 12\nabla_x\left(\bar E^2+\bar B^2\right)-\nabla_x\cdot\left(\bar E\otimes \bar E+\bar B\otimes \bar B\right)
			% & =
			% -\bar n\bar E-\bar j\wedge \bar B,
			% \\
			% \partial_t\left(E\wedge \bar B+\bar E\wedge B\right)+\nabla_x\left(E\cdot \bar E+B\cdot \bar B\right)
			% \hspace{20mm}
			% &
			% \\
			% -\nabla_x\cdot\left(E\otimes \bar E+\bar E\otimes E+B\otimes \bar B+\bar B\otimes B\right)
			% & =
			% -\left(n\bar E+j\wedge \bar B\right)
			% \\
			% & -\left(\bar n E+\bar j\wedge B\right),
		\end{aligned}
	\end{equation*}
	we arrive at the relation
	\begin{equation*}
		\begin{aligned}
			\frac{d}{dt} \int_{\mathbb{R}^3} & \left( u\cdot\bar u +\frac 14 n\bar n
			+\frac 12 E\cdot \bar E+ \frac 12 B\cdot\bar B
			+\frac 12 \left(\left(E-\bar E\right)\wedge \left(B-\bar B\right)\right)\cdot\bar u \right) dx
			\\
			& - \int_{\mathbb{R}^3}
			\frac 12 \left(\left(E-\bar E\right)\wedge \left(B-\bar B\right)\right)\cdot \partial_t \bar u
			dx
			\\
			& +\int_{\mathbb{R}^3}2\mu \nabla_x u:\nabla_x\bar u
			+\frac 1{\sigma}
			\left(j-nu\right)\cdot\left(\bar j-\bar n\bar u\right)
			dx
			\\
			& =
			\int_{\mathbb{R}^3}
			(u-\bar u) \otimes (u-\bar u):\nabla_x\bar u dx
			\\
			& -\frac 12
			\int_{\mathbb{R}^3}
			\left(\left(E-\bar E\right)\otimes \left(E-\bar E\right)+\left(B-\bar B\right)\otimes \left(B-\bar B\right)\right)
			: \nabla_x\bar u
			dx
			\\
			& +\frac 12\int_{\mathbb{R}^3}
			\left(\bar j-\bar n\bar u\right)\cdot\left(\left(u-\bar u\right)\wedge\left(B-\bar B\right)\right)
			dx
			\\
			& +\frac 12\int_{\mathbb{R}^3}
			\left(n-\bar n\right)\left(u-\bar u\right)\cdot\left(\frac 12\nabla_x\bar n-\bar E - \bar u\wedge\bar B\right)
			dx
			\\
			& - \int_{\mathbb{R}^3} \mathbf{A}\cdot
			\begin{pmatrix}
				u \\ j-nu
			\end{pmatrix}
			dx.
		\end{aligned}
	\end{equation*}

	On the whole, combining the preceding identity with the energy conservation law for test functions \eqref{energy test} and the energy decay imposed by the formal energy conservations from Proposition \ref{energy estimate 2}, we find the following modulated energy inequality~:
	\begin{equation*}
		\begin{aligned}
			\frac{d}{dt} & \delta\CE(t) + \delta\CD(t)
			\\
			& \leq
			-\int_{\mathbb{R}^3}
			(u-\bar u) \otimes (u-\bar u):\nabla_x\bar u dx
			- \int_{\mathbb{R}^3}
			\frac 12 \left(\left(E-\bar E\right)\wedge \left(B-\bar B\right)\right)\cdot \partial_t \bar u
			dx
			\\
			& +\frac 12
			\int_{\mathbb{R}^3}
			\left(\left(E-\bar E\right)\otimes \left(E-\bar E\right)+\left(B-\bar B\right)\otimes \left(B-\bar B\right)\right)
			: \nabla_x\bar u
			dx
			\\
			& -\frac 12\int_{\mathbb{R}^3}
			\left(\bar j-\bar n\bar u\right)\cdot\left(\left(u-\bar u\right)\wedge\left(B-\bar B\right)\right)
			dx
			\\
			& -\frac 12\int_{\mathbb{R}^3}
			\left(n-\bar n\right)\left(u-\bar u\right)\cdot\left(\frac 12\nabla_x\bar n-\bar E - \bar u\wedge\bar B\right)
			dx
			\\
			& + \int_{\mathbb{R}^3} \mathbf{A}\cdot
			\begin{pmatrix}
				u-\bar u \\ j-nu-(\bar j-\bar n\bar u)
			\end{pmatrix}
			dx.
		\end{aligned}
	\end{equation*}
	The next step consists in estimating the terms in the right-hand side above that are nonlinear in $(u,n,j,E,B)$ and to absorb the resulting expressions with the modulated energy $\delta\CE(t)$ by suitable uses of Young's inequality and Gr\"onwall's lemma. Thus, we obtain
	\begin{equation*}
		\begin{aligned}
			& \frac{d}{dt} \delta\CE(t) + \delta\CD(t) \\
			& \leq
			\left\|\nabla_{t,x} \bar u\right\|_{L^\infty_x}\left( \left\| u-\bar u\right\|_{L^2_x}^2
			+ \frac 12\left\|E-\bar E\right\|_{L^2_x}^2+\frac 12\left\|B-\bar B\right\|_{L^2_x}^2\right)
			\\
			& + \frac 12 \left(\left\|\bar j-\bar n\bar u\right\|_{L^\infty_x}\left\|B-\bar B\right\|_{L^2_x}
			+\left\|\frac 12\nabla_x\bar n-\bar E-\bar u\wedge \bar B\right\|_{L^\infty_x}\left\|n-\bar n\right\|_{L^2_x}\right)\left\|u-\bar u\right\|_{L^2_x}
			\\
			&
			+
			\int_{\mathbb{R}^3} \mathbf{A}\cdot
			\begin{pmatrix}
				u-\bar u \\ j-nu-(\bar j-\bar n\bar u)
			\end{pmatrix}
			dx
			\\
			& \leq
			\left\|\nabla_{t,x} \bar u\right\|_{L^\infty_x}\left( \left\| u-\bar u\right\|_{L^2_x}^2
			+ \frac 12\left\|E-\bar E\right\|_{L^2_x}^2+\frac 12\left\|B-\bar B\right\|_{L^2_x}^2\right)
			\\
			& + \frac{\sqrt 2}{4}\left\|\bar j-\bar n\bar u\right\|_{L^\infty_x}
			\left(
			\left\|u-\bar u\right\|_{L^2_x}^2
			+
			\frac 12\left\|B-\bar B\right\|_{L^2_x}^2\right)
			\\
			& +
			\frac 12\left\|\frac 12\nabla_x\bar n-\bar E-\bar u\wedge \bar B\right\|_{L^\infty_x}
			\left(
			\left\|u-\bar u\right\|_{L^2_x}^2
			+
			\frac 14\left\|n-\bar n\right\|_{L^2_x}^2\right)
			\\
			&
			+
			\int_{\mathbb{R}^3} \mathbf{A}\cdot
			\begin{pmatrix}
				u-\bar u \\ j-nu-(\bar j-\bar n\bar u)
			\end{pmatrix}
			dx.
% 			& \leq
% 			\frac 3{2\mu}\left\|\bar u\right\|_{L^\infty_x}^2 \left\| u-\bar u\right\|_{L^2_x}^2
% 			+ \frac 12\left\|\nabla_x \bar u\right\|_{L^\infty_x}
% 			\left(\left\|E-\bar E\right\|_{L^2_x}^2+\left\|B-\bar B\right\|_{L^2_x}^2\right)
% 			\\
% 			& +\frac{3C_0^2}{8\mu}\left( \left\|\bar j-\bar n\bar u\right\|_{L^3_x}^2
% 			\left\|B-\bar B\right\|_{L^2_x}^2
% 			+\left\|\frac 12\nabla_x\bar n-\bar E-\bar u\wedge \bar B\right\|_{L^3_x}^2
% 			\left\|n-\bar n\right\|_{L^2_x}^2\right)
% 			\\
% 			& + \frac\mu 2\left\|\nabla_x(u-\bar u)\right\|_{L^2_x}^2
% 			+
% 			\int_{\mathbb{R}^3} \mathbf{A}\cdot
% 			\begin{pmatrix}
% 				u-\bar u \\ j-nu-(\bar j-\bar n\bar u)
% 			\end{pmatrix}
% 			dx.
		\end{aligned}
	\end{equation*}
	Hence, further noticing that
	\begin{equation*}
		\begin{aligned}
			\left(1-\left\|\bar u\right\|_{L^\infty_{x}}\right)
			& \left(
			\frac 12\left\|u-\bar u\right\|_{L^2_x}^2
			+ \frac 18 \left\|n-\bar n\right\|_{L^2_x}^2
			+\frac 14
			\left\|E-\bar E\right\|_{L^2_x}^2+\frac 14\left\|B-\bar B\right\|_{L^2_x}^2\right)
			\\
			& \leq
			\frac12 \left\|u-\bar u\right\|_{L^2_x}^2
			+ \frac 18 \left\|n-\bar n\right\|_{L^2_x}^2
			+
			\frac 14\left\|E-\bar E\right\|_{L^2_x}^2+\frac 14\left\|B-\bar B\right\|_{L^2_x}^2
			\\
			& -\frac 12\left\|\bar u\right\|_{L^\infty_{x}}\int_{\mathbb{R}^3}\left|E-\bar E\right|\left|B-\bar B\right|dx
			\\
			& \leq
			\frac 12 \left\|u-\bar u\right\|_{L^2_x}^2
			+ \frac 18 \left\|n-\bar n\right\|_{L^2_x}^2
			+
			\frac 14\left\|E-\bar E\right\|_{L^2_x}^2+\frac 14\left\|B-\bar B\right\|_{L^2_x}^2
			\\
			& -\frac 12\int_{\mathbb{R}^3}\left(\left(E-\bar E\right)\wedge\left(B-\bar B\right)\right)\cdot\bar udx
			\\
			& = \delta\mathcal{E}(t),
		\end{aligned}
	\end{equation*}
	we find, since $\left\|\bar u\right\|_{L^\infty_{t,x}}<1$, that
	\begin{equation*}
		\frac{d}{dt}\delta\CE(t)
		+\delta\CD(t)
		\leq
		\lambda(t)
		\delta\CE(t)
		+ \int_{\mathbb{R}^3} \mathbf{A}\cdot
		\begin{pmatrix}
			u-\bar u \\ j-nu-(\bar j-\bar n\bar u)
		\end{pmatrix}
		dx,
	\end{equation*}
	which concludes the proof of the proposition with a direct application of Gr\"onwall's lemma.
\end{proof}

\begin{rem}
	Notice that the preceding method of modulation of the Poynting vector is not applicable to the incompressible Navier-Stokes-Maxwell system \eqref{INSM}, for the divergence of the electric field $E$ is not determined therein, i.e.\ Gauss' law $\Div E=n$ cannot be used to provide a bound on $\Div E$.
\end{rem}

The preceding proposition provides another weak-strong stability property for the two-fluid incompressible Navier-Stokes-Maxwell system with Ohm's law \eqref{TFINSMO}. Indeed, the stability inequality \eqref{stability 3} essentially implies that a solution $(\bar u, \bar n, \bar j,\bar E,\bar B)$ of \eqref{TFINSMO} such that $\bar u\in L^\infty_{t,x}$, $\nabla_{t,x}\bar u\in L^1_tL^\infty_x$, $\bar j-\bar n\bar u\in L^1_tL^\infty_x$ and $\left\|\bar u\right\|_{L^\infty_{t,x}}<1$, if it exists, is unique in the whole class of weak solutions in the energy space, for any given initial data.

We do not know whether the condition $\left\|\bar u\right\|_{L^\infty_{t,x}}<1$ in Proposition \ref{stability poynting} is actually necessary or merely a technical limitation. Nevertheless, this result shows that such a condition has a stabilizing effect on the two-fluid incompressible Navier-Stokes-Maxwell system with Ohm's law \eqref{TFINSMO}. Furthermore, this restriction is physically relevant since it imposes that the modulus of the bulk velocity $\bar u$ remains everywhere and at all times below the speed of light. More precisely, keeping track of the relevant physical constants in the formal derivations of Chapter \ref{formal-chap}, we see that the system \eqref{TFINSMO} can be recast as
\begin{equation*}
	\begin{cases}
		\begin{aligned}
			\d_t u +
			u\cdot\nabla_x u - \mu\Delta_x u
			& = -\nabla_x p+
			\frac 12 \left(c nE + j \wedge B\right) , & \Div u & = 0,\\
			\frac 1c\d_t E - \ROT B &= -  j, & \Div E & = c n,
			\\
			\frac 1c\d_t B + \ROT E & = 0, & \Div B & = 0, \\
			j-nu & = \sigma\left(-\frac 12 \nabla_x n + c E + u\wedge B\right), &&
		\end{aligned}
	\end{cases}
\end{equation*}
where the constant $c>0$ denotes the speed of light. Notice that the formal energy conservation law satisfied by this system is independent of $c>0$ and is thus given by Proposition \ref{energy estimate 2}. Moreover, expressing the Lorentz force with the Poynting vector as in \eqref{poynting} yields now
\begin{equation*}
	\frac 1c\partial_t\left(E\wedge B\right)+\frac 12\nabla_x\left(E^2+B^2\right)-\nabla_x\cdot\left(E\otimes E+B\otimes B\right)
	=
	-cnE-j\wedge B.
\end{equation*}
Therefore, applying the proof of Proposition \ref{stability poynting} to the preceding system, we arrive at a stability inequality valid under the restriction that the bulk velocity remains bounded by the speed of light $\left\|\bar u\right\|_{L^\infty_{t,x}}<c$, which is natural.

\bigskip

Following the previous developments, it is also possible to use now the stability inequality \eqref{stability 3} from Proposition \ref{stability poynting} to define another notion of dissipative solutions for the two-fluid incompressible Navier-Stokes-Maxwell system with Ohm's law \eqref{TFINSMO}, whose existence is then established by reproducing the arguments from Theorem \ref{dissipative 2}. Indeed, applying the computations from the proof of Proposition \ref{stability poynting} to the viscous approximation \eqref{viscous approximation} only produces new dissipative terms which are easily controlled in the limit $\nu\to 0$ (note that the condition $\left\|\bar u\right\|_{L^\infty_{t,x}}<1$ has to be used in order to absorb the dissipative terms produced by expressing the Lorentz force with the Poynting vector through the viscous Maxwell system from \eqref{viscous approximation}). Thus, the only remaining argument from the proof of Theorem \ref{dissipative 2} that needs special care in order to conclude the existence of dissipative solutions resides in the weak lower semi-continuity of the modulated energy $\delta\mathcal{E}(t)$ defined by \eqref{modulated defi 2}, which we establish now.

To this end, let us consider
\begin{equation*}
	\begin{aligned}
		E_\nu & {\rightharpoonup} E & \text{in }& L^2_x, \\
		B_\nu & {\rightharpoonup} B & \text{in }& L^2_x,
	\end{aligned}
\end{equation*}
as $\nu \to 0$. It is enough to show that
\begin{equation}\label{convex wedge}
	\begin{aligned}
		\frac 12 \left\|E\right\|_{L^2_x}^2
		+ \frac 12 \left\|B\right\|_{L^2_x}^2
		& - \int_{\mathbb{R}^3} \left(E\wedge B\right)\cdot\bar u dx
		\\
		& \leq
		\liminf_{\nu\to 0}\left(\frac 12 \left\|E_\nu\right\|_{L^2_x}^2
		+ \frac 12 \left\|B_\nu\right\|_{L^2_x}^2
		- \int_{\mathbb{R}^3} \left(E_\nu\wedge B_\nu\right)\cdot\bar u dx\right),
	\end{aligned}
\end{equation}
provided $\left\|\bar u\right\|_{L^\infty_x}<1$, which will follow from a convexity argument. Indeed, defining the bilinear form $\mathcal{B}:\mathbb{R}^6\times\mathbb{R}^6\to\mathbb{R}$ by
\begin{equation*}
	\mathcal{B}\left(
	\begin{pmatrix}
		E \\B
	\end{pmatrix},
	\begin{pmatrix}
		E_\nu\\B_\nu
	\end{pmatrix}
	\right)
	=E\cdot E_\nu+B\cdot B_\nu-\left(E\wedge B_\nu\right)\cdot \bar u
	-\left(E_\nu\wedge B\right)\cdot \bar u,
\end{equation*}
it is readily seen that $\mathcal{B}$ is symmetric and positive definite~:
\begin{equation*}
	\mathcal{B}\left(
	\begin{pmatrix}
		E\\B
	\end{pmatrix},
	\begin{pmatrix}
		E\\B
	\end{pmatrix}
	\right)
	=\left|E\right|^2+\left|B\right|^2-2\left(E\wedge B\right)\cdot \bar u\geq
	\left(1-|\bar u|\right)\left(\left|E\right|^2+\left|B\right|^2\right)\geq 0.
\end{equation*}
In particular, it follows that
\begin{equation*}
	\mathcal{B}\left(
	\begin{pmatrix}
		E\\B
	\end{pmatrix},
	\begin{pmatrix}
		E_\nu\\B_\nu
	\end{pmatrix}
	\right)
	\leq
	\frac 12
	\mathcal{B}\left(
	\begin{pmatrix}
		E\\B
	\end{pmatrix},
	\begin{pmatrix}
		E\\B
	\end{pmatrix}
	\right)
	+\frac 12
	\mathcal{B}\left(
	\begin{pmatrix}
		E_\nu\\B_\nu
	\end{pmatrix},
	\begin{pmatrix}
		E_\nu\\B_\nu
	\end{pmatrix}
	\right).
\end{equation*}
Hence, we deduce
\begin{equation*}
	\begin{aligned}
		\int_{\mathbb{R}^3}\mathcal{B}\left(
		\begin{pmatrix}
			E\\B
		\end{pmatrix},
		\begin{pmatrix}
			E\\B
		\end{pmatrix}
		\right)dx
		& =
		\lim_{\nu\to 0}\int_{\mathbb{R}^3}\mathcal{B}\left(
		\begin{pmatrix}
			E\\B
		\end{pmatrix},
		\begin{pmatrix}
			E_\nu\\B_\nu
		\end{pmatrix}
		\right)dx
		\\
		& \leq
		\frac 12
		\int_{\mathbb{R}^3}\mathcal{B}\left(
		\begin{pmatrix}
			E\\B
		\end{pmatrix},
		\begin{pmatrix}
			E\\B
		\end{pmatrix}
		\right)dx
		\\
		& +
		\lim_{\nu\to 0}\frac 12\int_{\mathbb{R}^3}\mathcal{B}\left(
		\begin{pmatrix}
			E_\nu\\B_\nu
		\end{pmatrix},
		\begin{pmatrix}
			E_\nu\\B_\nu
		\end{pmatrix}
		\right)dx,
	\end{aligned}
\end{equation*}
which establishes \eqref{convex wedge}.

\subsubsection{The two-fluid incompressible Navier-Stokes-Maxwell system with solenoidal Ohm's law}

Following the strategy of Propositions \ref{modulated energy estimate 0} and \ref{modulated energy estimate 1}, the next result establishes a crucial weak-strong stability principle for the two-fluid incompressible Navier-Stokes-Maxwell system with solenoidal Ohm's law \eqref{TFINSMSO}.

\begin{prop}\label{modulated energy estimate 2}
	Let $(u,E,B)$ be a smooth solution to the two-fluid incompressible Navier-Stokes-Maxwell system with solenoidal Ohm's law \eqref{TFINSMSO}. Further consider test functions $\left(\bar u, \bar j, \bar E, \bar B\right)\in C_c^\infty\left([0,\infty)\times\mathbb{R}^3\right)$ such that
	\begin{equation}\label{test constraints 2}
		\begin{cases}
			\begin{aligned}
				\Div \bar j & = 0, & \Div \bar u & = 0, \\
				\d_t \bar E - \ROT \bar B &= -  \bar j, & \Div \bar E & = 0,
				\\
				\d_t \bar B + \ROT \bar E & = 0, & \Div \bar B & = 0.
			\end{aligned}
		\end{cases}
	\end{equation}
	% \begin{equation*}
	% 	\begin{gathered}
	% 		\DIV \bar u = 0 , \qquad \DIV \bar E = 0 , \qquad \DIV \bar B = 0 , \\
	% 		\bar j = \sigma\left(- \nabla_x \bar p + \bar E + \bar u\wedge \bar B\right) , \qquad \Div \bar j = 0.
	% 	\end{gathered}
	% \end{equation*}
	We define the acceleration operator by
	\begin{equation*}
		\mathbf{A}\left(\bar u, \bar j, \bar E, \bar B\right)
		=
		\begin{pmatrix}
			-\d_t \bar u -
			P\left(\bar u\cdot\nabla_x \bar u\right) + \mu\Delta_x \bar u
			+ \frac 12 P \left(\bar j \wedge \bar B\right)
			\\
			- \frac 1{2\sigma}\bar j + \frac 12 P\left(\bar E + \bar u\wedge \bar B\right)
		\end{pmatrix},
	\end{equation*}
	and the growth rate by
	\begin{equation*}
		\lambda(t)=
		\left(\frac 2\mu + 4\sigma\right)\left\|\bar u(t)\right\|_{L^\infty_x}^2
		+ \frac{2C_0^2}{\mu} \left\|\bar j(t)\right\|_{L^3_x}^2,
	\end{equation*}
	where $C_0>0$ denotes the operator norm of the Sobolev embedding $\dot H^1\left(\mathbb{R}^3\right) \hookrightarrow L^6\left(\mathbb{R}^3\right)$.

	Then, one has the stability inequality
	\begin{equation}\label{stability 2}
		\begin{aligned}
			\delta\CE(t) + & \frac 12 \int_0^t \delta\CD(s) e^{\int_s^t\lambda(\sigma)d\sigma}ds
			\\
			& \leq \delta\CE(0) e^{\int_0^t\lambda(s)ds}
			+\int_0^t
			\left[\int_{\mathbb{R}^3} \mathbf{A}\cdot
			\begin{pmatrix}
				u-\bar u \\ j - \bar j
			\end{pmatrix}
			dx\right](s)
			e^{\int_s^t\lambda(\sigma)d\sigma}ds,
		\end{aligned}
	\end{equation}
	where the modulated energy $\delta\CE$ and energy dissipation $\delta\CD$ are given by
	\begin{equation*}
		\begin{aligned}
			\delta\CE(t)
			& =
			\frac 12 \left\|\left(u-\bar u\right)(t)\right\|_{L^2_x}^2
			+ \frac 14 \left\|\left(E-\bar E\right)(t)\right\|_{L^2_x}^2
			+ \frac 14 \left\|\left(B-\bar B\right)(t)\right\|_{L^2_x}^2, \\
			\delta \CD(t)
			& =
			\mu \left\|\nabla_x (u-\bar u)(t)\right\|_{L^2_x}^2
			+ \frac 1{2\sigma} \left\|\left(j-\bar j\right)(t)\right\|_{L^2_x}^2.
		\end{aligned}
	\end{equation*}
\end{prop}

\begin{proof}
	We have already formally established in Proposition \ref{energy estimate 2} the conservation of the energy for system \eqref{TFINSMSO}. The very same computations applied to the test functions $(\bar u, \bar j, \bar E, \bar B)$ yield the identity
	\begin{equation}\label{energy test 2}
		\frac{d}{dt}\bar\CE(t)+\bar\CD(t)= - \int_{\mathbb{R}^3} \mathbf{A}\cdot
		\begin{pmatrix}
			\bar u \\ \bar j
		\end{pmatrix}
		dx,
	\end{equation}
	where the energy $\bar\CE$ and energy dissipation $\bar\CD$ are obtained simply by replacing the unknowns by the test functions in the respective definitions of Proposition \ref{energy estimate 2}.
	
	Furthermore, another similar duality computation yields that
	\begin{equation}\label{bad modulation}
		\begin{aligned}
			\frac{d}{dt} & \int_{\mathbb{R}^3} \left(u\cdot\bar u
			+\frac 12 E\cdot \bar E+ \frac 12 B\cdot\bar B\right) dx
			+\int_{\mathbb{R}^3}2\mu \nabla_x u:\nabla_x\bar u
			+\frac 1{\sigma} j\cdot\bar j dx \\
			& =
			- \int_{\mathbb{R}^3}
			\bar u \otimes (u-\bar u):\nabla_x(u-\bar u) dx
			\\
			& +\frac 12 \int_{\mathbb{R}^3}
			\left((j -\bar j)\wedge(B-\bar B)\right)\cdot\bar u
			+ \left((u-\bar u)\wedge(B-\bar B)\right)\cdot\bar j
			dx
			\\
			&
			- \int_{\mathbb{R}^3} \mathbf{A}\cdot
			\begin{pmatrix}
				u \\ j
			\end{pmatrix}
			dx.
		\end{aligned}
	\end{equation}

	On the whole, combining the above identities with the energy decay imposed by the formal energy conservations from Proposition \ref{energy estimate 2}, we find the following modulated energy inequality~:
	\begin{equation*}
		\begin{aligned}
			\frac{d}{dt} & \delta\CE(t) + \delta\CD(t) \\
			& \leq
			\int_{\mathbb{R}^3}
			\bar u \otimes(u-\bar u):\nabla_x(u-\bar u) dx
			\\
			& -\frac 12 \int_{\mathbb{R}^3}
			\left((j - \bar j)\wedge(B-\bar B)\right)\cdot\bar u
			+ \left((u-\bar u)\wedge(B-\bar B)\right)\cdot\bar j
			dx
			\\
			&
			+
			\int_{\mathbb{R}^3} \mathbf{A}\cdot
			\begin{pmatrix}
				u-\bar u \\ j-\bar j
			\end{pmatrix}
			dx.
		\end{aligned}
	\end{equation*}
	The next step consists in estimating the terms in the right-hand side above that are nonlinear in $(u,j,E,B)$ and to absorb the resulting expressions with the modulated energy $\delta\CE(t)$ and the modulated energy dissipation $\delta\CD(t)$ by suitable uses of Young's inequality and Gr\"onwall's lemma. Thus, we obtain
	\begin{equation*}
		\begin{aligned}
			\frac{d}{dt} & \delta\CE(t) + \delta\CD(t) \\
			& \leq
			\left\|\bar u\right\|_{L^\infty_x} \left\| u-\bar u\right\|_{L^2_x}\left\|\nabla_x(u-\bar u)\right\|_{L^2_x}
			\\
			& +
			\frac 12 \left\|\bar j\right\|_{L^3_x}\left\|B-\bar B\right\|_{L^2_x}\left\|u-\bar u\right\|_{L^6_x}
			+ \frac 12 \left\|\bar u\right\|_{L^\infty_x}\left\|B-\bar B\right\|_{L^2_x}\left\|j-\bar j\right\|_{L^2_x}
			\\
			& + \int_{\mathbb{R}^3} \mathbf{A}\cdot
			\begin{pmatrix}
				u-\bar u \\ j - \bar j
			\end{pmatrix}
			dx
			\\
			& \leq
			\frac 1{\mu}\left\|\bar u\right\|_{L^\infty_x}^2 \left\| u-\bar u\right\|_{L^2_x}^2
			+\left(\sigma\left\|\bar u\right\|_{L^\infty_x}^2
			+ \frac{C_0^2}{2\mu} \left\|\bar j\right\|_{L^3_x}^2\right)
			\left\|B-\bar B\right\|_{L^2_x}^2
			\\
			& +\frac\mu 2\left\|\nabla_x(u-\bar u)\right\|_{L^2_x}^2
			+ \frac{1}{4\sigma}\left\|j-\bar j\right\|_{L^2_x}^2
			+ \int_{\mathbb{R}^3} \mathbf{A}\cdot
			\begin{pmatrix}
				u-\bar u \\ j - \bar j
			\end{pmatrix}
			dx.
		\end{aligned}
	\end{equation*}
	Hence,
	\begin{equation*}
		\frac{d}{dt}\delta\CE(t)
		+\frac 12\delta\CD(t)
		\leq
		\lambda(t)
		\delta\CE(t)
		+ \int_{\mathbb{R}^3} \mathbf{A}\cdot
		\begin{pmatrix}
			u-\bar u \\ j - \bar j
		\end{pmatrix}
		dx,
	\end{equation*}
	which concludes the proof of the proposition with a direct application of Gr\"onwall's lemma.
\end{proof}

Note that the test functions satisfying the linear constraints \eqref{test constraints 2} are easily constructed by considering vector potentials $\bar A\in C_c^\infty\left([0,\infty)\times\mathbb{R}^3\right)$ and then setting
\begin{equation*}
	\bar E = - \partial_t \ROT \bar A \qquad\text{and}\qquad \bar B = \ROT \ROT \bar A.
\end{equation*}
Now, one may prefer to deal, in a completely equivalent manner, with test functions $\left(\bar u, \bar E, \bar B\right)\in C_c^\infty\left([0,\infty)\times\mathbb{R}^3\right)$ and $\bar j \in C^\infty\left([0,\infty)\times\mathbb{R}^3\right)$ (here, we cannot impose that $\bar j$ be compactly supported) satisfying the stationary constraints
\begin{equation*}
	\Div \bar u = 0,
	\qquad
	\Div \bar E = 0,
	\qquad
	\Div \bar B = 0,
	\qquad
	\bar j = \sigma P\left( \bar E + \bar u\wedge \bar B\right),
\end{equation*}
rather than the constraints \eqref{test constraints 2}. In this case, instead of \eqref{stability 2}, we obtain the stability inequality
\begin{equation*}
	\begin{aligned}
		\delta\CE(t) + & \frac 12 \int_0^t \delta\CD(s) e^{\int_s^t\lambda(\sigma)d\sigma}ds
		\\
		& \leq \delta\CE(0) e^{\int_0^t\lambda(s)ds}
		+\int_0^t
		\left[\int_{\mathbb{R}^3} \mathbf{A}\cdot
		\begin{pmatrix}
			u-\bar u \\ E - \bar E \\ B - \bar B
		\end{pmatrix}
		dx\right](s)
		e^{\int_s^t\lambda(\sigma)d\sigma}ds,
	\end{aligned}
\end{equation*}
where the acceleration operator is now defined by
\begin{equation*}
	\mathbf{A}\left(\bar u, \bar j, \bar E, \bar B\right)
	=
	\begin{pmatrix}
		-\d_t \bar u -
		P\left(\bar u\cdot\nabla_x \bar u\right) + \mu\Delta_x \bar u
		+ \frac 12 P \left( \bar j \wedge \bar B\right)
		\\
		\frac 12\left(-\partial_t \bar E + \ROT\bar B - \bar j\right)
		\\
		\frac 12\left(-\partial_t \bar B - \ROT\bar E\right)
	\end{pmatrix}.
\end{equation*}

The preceding proposition provides an important weak-strong stability property for the two-fluid incompressible Navier-Stokes-Maxwell system with solenoidal Ohm's law \eqref{TFINSMSO}. Indeed, the stability inequality \eqref{stability 2} essentially implies that a solution $(\bar u, \bar j,\bar E,\bar B)$ of \eqref{TFINSMSO} such that $\bar u\in L^2_tL^\infty_x$ and $\bar j\in L^2_tL^3_x$, if it exists, is unique in the whole class of weak solutions in the energy space, for any given initial data.

\bigskip

By analogy with Lions' dissipative solutions to the incompressible Euler system \cite[Section 4.4]{lions7}, we provide now a suitable notion of dissipative solution for the two-fluid incompressible Navier-Stokes-Maxwell system with solenoidal Ohm's law \eqref{TFINSMSO}, based on Proposition \ref{modulated energy estimate 2}, and establish their existence next.

\begin{defi}
	We say that
	\begin{equation*}
		(u,E,B)\in L^\infty\left([0,\infty);L^2\left(\mathbb{R}^3\right)\right)
		\cap C\left([0,\infty);\textit{w-}L^2\left(\mathbb{R}^3\right)\right)
	\end{equation*}
	such that
	\begin{equation*}
		\DIV u = 0, \qquad \DIV E = 0, \qquad \DIV B = 0,
	\end{equation*}
	is a \textbf{dissipative solution of the two-fluid incompressible Navier-Stokes-Maxwell system with solenoidal Ohm's law \eqref{TFINSMSO}}, if it solves Maxwell's equations
	\begin{equation*}
		\begin{cases}
			\begin{aligned}
				\d_t E - \ROT B &= -  j, 
				\\
				\d_t B + \ROT E & = 0,
			\end{aligned}
		\end{cases}
	\end{equation*}
	with solenoidal Ohm's law
	\begin{equation*}
		j = \sigma\left(-\nabla_x\bar p + E + u\wedge B\right),
	\end{equation*}
	in the sense of distributions, and if, for any test functions $\left(\bar u, \bar j, \bar E, \bar B\right)\in C_c^\infty\left([0,\infty)\times\mathbb{R}^3\right)$ satisfying the linear constraints \eqref{test constraints 2}, the stability inequality \eqref{stability 2} is verified.
\end{defi}

As previously mentioned, dissipative solutions define actual solutions in the sense that they coincide with the unique strong solution when the latter exists. The following theorem asserts their existence.

\begin{thm}\label{dissipative 3}
	For any initial data $\left(u^\mathrm{in}, E^\mathrm{in},B^\mathrm{in}\right)\in L^2\left(\mathbb{R}^3\right)$ such that
	\begin{equation*}
		\DIV u^\mathrm{in} = 0, \qquad \DIV E^\mathrm{in} = 0, \qquad \DIV B^\mathrm{in} = 0,
	\end{equation*}
	there exists a dissipative solution to the two-fluid incompressible Navier-Stokes-Maxwell system with solenoidal Ohm's law \eqref{TFINSMSO}.
\end{thm}

\begin{proof}
	As in the proof of Theorems \ref{dissipative 1} and \ref{dissipative 2}, it is possible, here, to justify the existence of dissipative solutions by introducing viscous approximations of the system \eqref{TFINSMSO}. However, it will be much more judicious to recover the system \eqref{TFINSMSO} as an asymptotic regime of the two fluid incompressible Navier-Stokes-Maxwell system \eqref{two fluid} for very weak interspecies interactions, which we recast here, for all $\nu>0$, as
	\begin{equation}\label{two fluid 2}
		\begin{cases}
			\begin{aligned}
				\partial_t u_\nu^+ +
				u_\nu^+\cdot\nabla_x  u_\nu^+ - \mu\Delta_x  u_\nu^+ & && \\
				+ \frac 1{\sigma \nu^2}\left(u_\nu^+-u_\nu^-\right)
				& = -\nabla_x  p_\nu^+
				+ \frac1\nu\left(E_\nu + u_\nu^+ \wedge B_\nu \right) , & \DIV u_\nu^+ & = 0,\\
				\partial_t u_\nu^- +
				u_\nu^-\cdot\nabla_x  u_\nu^- - \mu\Delta_x  u_\nu^- & && \\
				- \frac 1{\sigma \nu^2}\left(u_\nu^+-u_\nu^-\right)
				& = -\nabla_x  p_\nu^-
				- \frac1\nu\left(E_\nu + u_\nu^- \wedge B_\nu \right) , & \DIV u_\nu^- & = 0,\\
				\partial_t E_\nu - \ROT B_\nu &= -  \frac 1{\nu}\left(u_\nu^+-u_\nu^-\right), & \DIV E_\nu & = 0,
				\\
				\partial_t B_\nu + \ROT E_\nu & = 0, & \DIV B_\nu & = 0,
			\end{aligned}
		\end{cases}
	\end{equation}
	associated with an initial data $\left(u^{\pm\mathrm{in}}_\nu,E^\mathrm{in},B^\mathrm{in}\right)$ satisfying
	\begin{equation*}
		u^\mathrm{in} = \frac{u_\nu^{+\mathrm{in}}+u_\nu^{-\mathrm{in}}}{2}.
	\end{equation*}
	The above two fluid system satisfies the energy inequality, for all $t>0$,
	\begin{equation*}
		\begin{aligned}
			\frac 12 & \left(\left\|u^+_\nu\right\|_{L^2_x}^2 + \left\|u^-_\nu\right\|_{L^2_x}^2
			+ \left\|E_\nu\right\|_{L^2_x}^2 + \left\|B_\nu\right\|_{L^2_x}^2 \right)(t) \\
			& \hspace{10mm} + \int_0^t \mu \left( \left\|\nabla_x u^+_\nu(s)\right\|_{L^2_x}^2 + \left\|\nabla_x u^-_\nu(s)\right\|_{L^2_x}^2 \right)
			+ \frac 1\sigma \left\|\frac{u^+_\nu(s)-u^-_\nu(s)}\nu\right\|_{L^2_x}^2 ds
			\\
			& \hspace{10mm} \leq \frac 12
			\left(\left\|u^{+\mathrm{in}}_\nu\right\|_{L^2_x}^2 + \left\|u^{-\mathrm{in}}_\nu\right\|_{L^2_x}^2
			+ \left\|E^\mathrm{in}\right\|_{L^2_x}^2 + \left\|B^\mathrm{in}\right\|_{L^2_x}^2 \right).
		\end{aligned}
	\end{equation*}
	Further defining the variables
	\begin{equation*}
		u_\nu=\frac{u_\nu^++u_\nu^-}{2}\qquad\text{and}\qquad j_\nu=\frac{u^+_\nu-u^-_\nu}{\nu},
	\end{equation*}
	the system \eqref{two fluid 2} can be rewritten as
	\begin{equation}\label{two fluid 3}
		\begin{cases}
			\begin{aligned}
				\partial_t u_\nu +
				u_\nu \cdot\nabla_x  u_\nu +\frac{\nu^2}{4} j_\nu\cdot\nabla_x  j_\nu - \mu\Delta_x  u_\nu \hspace{-15mm} & && \\
				& = -\nabla_x  p_\nu
				+ \frac12 j_\nu \wedge B_\nu , & \DIV u_\nu & = 0,\\
				\frac{\nu^2}{2}\left(\partial_t j_\nu +
				u_\nu\cdot\nabla_x  j_\nu + j_\nu\cdot\nabla_x  u_\nu - \mu \Delta_x  j_\nu \right) \hspace{-15mm}	 & && \\
				+ \frac 1{\sigma}j_\nu
				& = -\nabla_x  \bar p_\nu
				+ E_\nu + u_\nu \wedge B_\nu , & \DIV j_\nu & = 0,\\
				\partial_t E_\nu - \ROT B_\nu &= -  j_\nu, & \DIV E_\nu & = 0,
				\\
				\partial_t B_\nu + \ROT E_\nu & = 0, & \DIV B_\nu & = 0,
			\end{aligned}
		\end{cases}
	\end{equation}
	and the corresponding energy inequality becomes, for all $t>0$,
	\begin{equation*}
		\begin{aligned}
			& \left(\frac 12\left\|u_\nu\right\|_{L^2_x}^2 + \frac{\nu^2}{8}\left\|j_\nu\right\|_{L^2_x}^2
			+ \frac 14\left\|E_\nu\right\|_{L^2_x}^2 + \frac 14\left\|B_\nu\right\|_{L^2_x}^2 \right)(t) \\
			& \hspace{10mm} + \int_0^t \mu \left( \left\|\nabla_x u_\nu(s)\right\|_{L^2_x}^2 +\frac{\nu^2}{4} \left\|\nabla_x j_\nu(s)\right\|_{L^2_x}^2 \right)
			+ \frac 1{2\sigma} \left\|j_\nu(s)\right\|_{L^2_x}^2 ds
			\\
			& \hspace{10mm} \leq \frac 12 \left\|u^{\mathrm{in}}\right\|_{L^2_x}^2
			+ \frac{\nu^2}{8}\left\|j^{\mathrm{in}}_\nu\right\|_{L^2_x}^2
			+ \frac 14 \left\|E^\mathrm{in}\right\|_{L^2_x}^2 + \frac 14 \left\|B^\mathrm{in}\right\|_{L^2_x}^2 ,
		\end{aligned}
	\end{equation*}
	where $j^{\mathrm{in}}_\nu=\frac{u^{+\mathrm{in}}_\nu-u^{-\mathrm{in}}_\nu}{\nu}$.

	Weak solutions of the above systems \eqref{two fluid 2} and \eqref{two fluid 3} are easily established following the method of Leray \cite{leray}, for the nonlinear terms $u^\pm_\nu\wedge B_\nu$ (or, equivalently, $u_\nu\wedge B_\nu$ and $j_\nu\wedge B_\nu$) are stable with respect to weak convergence in the energy space defined by the above energy inequalities.

	Now, for any test functions $\left(\bar u, \bar j, \bar E, \bar B\right)\in C_c^\infty\left([0,\infty)\times\mathbb{R}^3\right)$ satisfying the linear constraints \eqref{test constraints 2}, we define the approximate acceleration operator by
	\begin{equation*}
		\mathbf{A}_\nu \left(\bar u, \bar j, \bar E, \bar B\right)
		=
		\mathbf{A} \left(\bar u, \bar j, \bar E, \bar B\right)
		-\frac{\nu^2}{4}
		\begin{pmatrix}
			P\left( \bar j\cdot\nabla_x  \bar j\right)
			\\
			\partial_t \bar j +
			P\left( \bar u \cdot\nabla_x \bar j + \bar j \cdot\nabla_x \bar u \right) - \mu \Delta_x \bar j
		\end{pmatrix}.
	\end{equation*}
	Then, a straightforward energy estimate yields that
	\begin{equation*}
		\begin{aligned}
			\frac{d}{dt}& \left(\frac 12\left\|\bar u\right\|_{L^2_x}^2 + \frac{\nu^2}{8}\left\|\bar j\right\|_{L^2_x}^2
			+ \frac 14\left\|\bar E\right\|_{L^2_x}^2 + \frac 14\left\|\bar B\right\|_{L^2_x}^2 \right) \\
			& \hspace{10mm} + \mu \left( \left\|\nabla_x \bar u\right\|_{L^2_x}^2 +\frac{\nu^2}{4} \left\|\nabla_x \bar j \right\|_{L^2_x}^2 \right)
			+ \frac 1{2\sigma} \left\|\bar j\right\|_{L^2_x}^2 = - \int_{\mathbb{R}^3} \mathbf{A}_\nu\cdot
			\begin{pmatrix}
				\bar u \\ \bar j
			\end{pmatrix}
			dx.
		\end{aligned}
	\end{equation*}
	Moreover, another similar duality computation gives that
	\begin{equation*}
		\begin{aligned}
			\frac{d}{dt}\int_{\mathbb{R}^3} & \left(u_\nu\cdot\bar u + \frac{\nu^2}{4}j_\nu\cdot\bar j
			+ \frac 12 E_\nu\cdot \bar E+ \frac 12 B_\nu\cdot\bar B\right) dx
			\\
			& + \int_{\mathbb{R}^3}2\mu \nabla_x u_\nu:\nabla_x\bar u
			+\frac{\mu\nu^2}{2} \nabla_x j_\nu:\nabla_x\bar j
			+\frac 1{\sigma} j_\nu\cdot\bar j dx \\
			& =
			- \int_{\mathbb{R}^3}
			\bar u \otimes\left(u_\nu-\bar u\right):\nabla_x(u_\nu-\bar u)
			+
			\frac{\nu^2}{4}
			\bar u \otimes\left(j_\nu-\bar j\right):\nabla_x(j_\nu-\bar j)
			dx
			\\
			& -\frac{\nu^2}{4}\int_{\mathbb{R}^3}
			\bar j \otimes\left(j_\nu-\bar j\right):\nabla_x(u_\nu-\bar u)
			+
			\bar j \otimes\left(u_\nu-\bar u\right):\nabla_x(j_\nu-\bar j)
			dx
			\\
			& + \frac 12 \int_{\mathbb{R}^3}
			\left((j_\nu-\bar j)\wedge(B_\nu-\bar B)\right)\cdot\bar u
			+ \left((u_\nu-\bar u)\wedge(B_\nu-\bar B)\right)\cdot\bar j
			dx
			\\
			& - \int_{\mathbb{R}^3} \mathbf{A}_\nu\cdot
			\begin{pmatrix}
				u_\nu \\ j_\nu
			\end{pmatrix}
			dx.
		\end{aligned}
	\end{equation*}
	Hence, defining the modulated energy $\delta\CE_\nu(t)$ and modulated energy dissipation $\delta\CD_\nu(t)$ by
	\begin{equation*}
		\begin{aligned}
			\delta\CE_\nu(t)
			& =
			\frac 12 \left\|\left(u_\nu-\bar u\right)(t)\right\|_{L^2_x}^2
			+ \frac {\nu^2}{8} \left\|\left(j_\nu-\bar j\right)(t)\right\|_{L^2_x}^2 \\
			& + \frac 14 \left\|\left(E_\nu-\bar E\right)(t)\right\|_{L^2_x}^2
			+ \frac 14 \left\|\left(B_\nu-\bar B\right)(t)\right\|_{L^2_x}^2, \\
			\delta \CD_\nu(t)
			& =
			\mu \left\|\nabla_x (u_\nu-\bar u)(t)\right\|_{L^2_x}^2
			+ \frac{\mu\nu^2}{4} \left\|\nabla_x (j_\nu-\bar j)(t)\right\|_{L^2_x}^2
			+ \frac 1{2\sigma} \left\|\left(j_\nu-\bar j\right)(t)\right\|_{L^2_x}^2,
		\end{aligned}
	\end{equation*}
	we find that
	\begin{equation*}
		\begin{aligned}
			& \delta\CE_\nu(t) + \int_0^t \delta\CD_\nu(s) ds \\
			& \leq \delta\CE_\nu(0) +
			\int_0^t\int_{\mathbb{R}^3}
			\bar u \otimes\left(u_\nu-\bar u\right):\nabla_x(u_\nu-\bar u)
			+
			\frac{\nu^2}{4}
			\bar u \otimes\left(j_\nu-\bar j\right):\nabla_x(j_\nu-\bar j)
			dxds
			\\
			& + \frac{\nu^2}{4}\int_0^t\int_{\mathbb{R}^3}
			\bar j \otimes\left(j_\nu-\bar j\right):\nabla_x(u_\nu-\bar u)
			+
			\bar j \otimes\left(u_\nu-\bar u\right):\nabla_x(j_\nu-\bar j)
			dxds
			\\
			& - \frac 12 \int_0^t \int_{\mathbb{R}^3}
			\left((j_\nu-\bar j)\wedge(B_\nu-\bar B)\right)\cdot\bar u
			+ \left((u_\nu-\bar u)\wedge(B_\nu-\bar B)\right)\cdot\bar j
			dxds
			\\
			& + \int_0^t \int_{\mathbb{R}^3} \mathbf{A}_\nu\cdot
			\begin{pmatrix}
				u_\nu -\bar u \\ j_\nu -\bar j
			\end{pmatrix}
			dxds.
		\end{aligned}
	\end{equation*}
	The next step consists in estimating the terms in the right-hand side above that are nonlinear in $(u_\nu,j_\nu,E_\nu,B_\nu)$ and to absorb the resulting expressions with the modulated energy $\delta\CE_\nu(t)$ and the modulated energy dissipation $\delta\CD_\nu(t)$ by suitable uses of Young's inequality and Gr\"onwall's lemma. Thus, we obtain
	\begin{equation*}
		\begin{aligned}
			& \delta\CE_\nu(t) + \int_0^t \delta\CD_\nu(s) ds \\
			& \leq \delta\CE_\nu(0) +
			\int_0^t \left\|\bar u\right\|_{L^\infty_x}
			\left\| u_\nu-\bar u\right\|_{L^2_x}\left\|\nabla_x(u_\nu-\bar u)\right\|_{L^2_x} ds
			\\
			& +\frac{\nu^2}{4}
			\int_0^t \left\|\bar u\right\|_{L^\infty_x}
			\left\| j_\nu-\bar j\right\|_{L^2_x}\left\|\nabla_x(j_\nu-\bar j)\right\|_{L^2_x} ds
			\\
			& + \frac{\nu^2}{4}
			\int_0^t
			\left\|\bar j\right\|_{L^\infty_x}
			\left\| j_\nu-\bar j\right\|_{L^2_x}\left\|\nabla_x(u_\nu-\bar u)\right\|_{L^2_x}
			+
			\left\|\bar j\right\|_{L^3_x}
			\left\| u_\nu-\bar u\right\|_{L^6_x}\left\|\nabla_x(j_\nu-\bar j)\right\|_{L^2_x}
			ds
			\\
			& +
			\frac 12 \int_0^t \left\|\bar j\right\|_{L^3_x}\left\|B_\nu-\bar B\right\|_{L^2_x}\left\|u_\nu-\bar u\right\|_{L^6_x}
			+ \left\|\bar u\right\|_{L^\infty_x}\left\|B_\nu-\bar B\right\|_{L^2_x}\left\|j_\nu-\bar j\right\|_{L^2_x}
			ds
			\\
			& + \int_0^t \int_{\mathbb{R}^3} \mathbf{A}_\nu\cdot
			\begin{pmatrix}
				u_\nu-\bar u \\ j_\nu - \bar j
			\end{pmatrix}
			dxds
			\\
			& \leq \delta\CE_\nu(0) +
			\int_0^t \frac 1{\mu}\left\|\bar u\right\|_{L^\infty_x}^2 \left\| u_\nu-\bar u\right\|_{L^2_x}^2
			+\left(\sigma\left\|\bar u\right\|_{L^\infty_x}^2
			+ \frac{C_0^2}{2\mu} \left\|\bar j\right\|_{L^3_x}^2\right)
			\left\|B_\nu-\bar B\right\|_{L^2_x}^2 ds
			\\
			& + \int_0^t \left(\frac\mu 2 + \frac{\mu\nu^2}{16} + \frac{\nu^2C_0^2}{4\mu}\left\|\bar j\right\|_{L^3_x}^2\right)\left\|\nabla_x(u_\nu-\bar u)\right\|_{L^2_x}^2
			+ \frac{\mu\nu^2}{8}\left\|\nabla_x(j_\nu-\bar j)\right\|_{L^2_x}^2 ds
			\\
			& + \int_0^t \left(\frac{1}{4\sigma}+\frac{\nu^2}{4\mu}\left\|\bar u\right\|_{L^\infty_x}^2
			+ \frac{\nu^2}{4\mu}
			\left\|\bar j\right\|_{L^\infty_x}^2\right)\left\|j_\nu-\bar j\right\|_{L^2_x}^2
			+ \left[\int_{\mathbb{R}^3} \mathbf{A}_\nu\cdot
			\begin{pmatrix}
				u_\nu-\bar u \\ j_\nu - \bar j
			\end{pmatrix}
			dx \right] ds
			\\
			& \leq \delta\CE_\nu(0) +
			\int_0^t
			\lambda(s)\delta\CE_\nu(s)
			+ \left[\int_{\mathbb{R}^3} \mathbf{A}_\nu\cdot
			\begin{pmatrix}
				u_\nu-\bar u \\ j_\nu - \bar j
			\end{pmatrix}
			dx\right]ds
			\\
			& +
			\int_0^t
			\left(\frac 12 + \nu^2 \beta(s)\right)\delta\CD_\nu(s) ds,
		\end{aligned}
	\end{equation*}
	where
	\begin{equation*}
		\beta(t) = \frac{1}{16} + \frac{C_0^2}{4\mu^2}\left\|\bar j\right\|_{L^3_x}^2
		+\frac{\sigma}{2\mu}
		\left(
		\left\|\bar u\right\|_{L^\infty_x}^2
		+
		\left\|\bar j\right\|_{L^\infty_x}^2
		\right).
	\end{equation*}
	Hence,
	\begin{equation*}
		\begin{aligned}
			\delta\CE_\nu(t) & + \int_0^t \frac 12 \delta\CD_\nu(s) ds \\
			& \leq \delta\CE_\nu(0) +
			\int_0^t \lambda(s)\delta\CE_\nu(s)
			+ \left[\int_{\mathbb{R}^3} \mathbf{A}_\nu\cdot
			\begin{pmatrix}
				u_\nu-\bar u \\ j_\nu - \bar j
			\end{pmatrix}
			dx\right] + \nu^2 \beta(s) \delta\CD_\nu(s)ds,
		\end{aligned}
	\end{equation*}
	and an application of Gr\"onwall's lemma yields
	\begin{equation*}
		\begin{aligned}
			\delta\CE_\nu(t) & + \frac 12 \int_0^t \delta\CD_\nu(s) e^{\int_s^t\lambda(\sigma)d\sigma}ds
			\leq \delta\CE_\nu(0) e^{\int_0^t\lambda(s)ds}
			\\
			& +\int_0^t
			\left[\left[\int_{\mathbb{R}^3} \mathbf{A}_\nu\cdot
			\begin{pmatrix}
				u_\nu-\bar u \\ j_\nu-\bar j
			\end{pmatrix}
			dx\right](s)
			+ \nu^2 \beta(s) \delta\CD_\nu(s)
			\right]
			e^{\int_s^t\lambda(\sigma)d\sigma}ds.
		\end{aligned}
	\end{equation*}

	We may now pass to the limit in the above stability inequality. Thus, up to extraction of subsequences, we may assume that, as $\nu\rightarrow 0$,
	\begin{equation*}
		\begin{aligned}
			u_\nu & \stackrel{*}{\rightharpoonup} u & \text{in }& L^\infty_tL^2_x\cap L^2_t\dot H^1_x,  \\
			j_\nu & \rightharpoonup j & \text{in }& L^2_t L^2_x, \\
			E_\nu & \stackrel{*}{\rightharpoonup} E & \text{in }& L^\infty_t L^2_x, \\
			B_\nu & \stackrel{*}{\rightharpoonup} B & \text{in }& L^\infty_t L^2_x.
		\end{aligned}
	\end{equation*}
	Furthermore, noticing that $\partial_t u_\nu$, $\partial_t E_\nu$ and $\partial_t B_\nu$ are uniformly bounded, in $L^1_\mathrm{loc}$ in time and in some negative index Sobolev space in $x$, it is possible to show (see \cite[Appendix C]{lions7}) that $(u_\nu,E_\nu,B_\nu)$ converges to $(u,E,B)\in C\left([0,\infty);\textit{w-}L^2\left(\mathbb{R}^3\right)\right)$ weakly in $L^2_x$ uniformly locally in time. Then, by the weak lower semi-continuity of the norms, we obtain that, for every $t>0$,
	\begin{equation*}
		\delta\CE(t) + \frac 12 \int_0^t \delta\CD(s) e^{\int_s^t\lambda(\sigma)d\sigma}ds
		\leq
		\liminf_{\nu\rightarrow 0}
		\delta\CE_\nu(t) + \frac 12 \int_0^t \delta\CD_\nu(s) e^{\int_s^t\lambda(\sigma)d\sigma}ds.
	\end{equation*}
	Hence, further assuming that $\delta\CE_\nu(0)\rightarrow \delta\CE(0)$, as $\nu\rightarrow 0$, the stability inequality \eqref{stability 2} holds. Notice that the convergence of the initial data is satisfied whenever $\left\|u^{+\mathrm{in}}_\nu-u^{-\mathrm{in}}_\nu\right\|_{L^2_x}\rightarrow 0$, as $\nu\rightarrow 0$.
	
	Finally, invoking a classical compactness result by Aubin and Lions \cite{aubin, lions6} (see also \cite{simon} for a sharp compactness criterion), we infer that the $u_\nu$'s converge towards $u$ strongly in $L^2_\mathrm{loc}\left(dtdx\right)$. Therefore, passing to the limit in the evolution equation for $j_\nu$ in \eqref{two fluid 3}, it is readily seen that Ohm's law is satisfied asymptotically, which concludes the proof of the theorem.
\end{proof}

As before, we present now an alternative kind of stability inequality for the two-fluid incompressible Navier-Stokes-Maxwell system with solenoidal Ohm's law \eqref{TFINSMSO}. It is a mere adaptation of Proposition \ref{stability poynting} to the present case, which relies on the interpretation of the Lorentz force with the Poynting vector. We recall that this method allows us to stabilize the modulated nonlinear terms solely with the modulated energy $\delta\mathcal{E}$.

\begin{prop}\label{stability poynting 2}
	Let $(u,E,B)$ be a smooth solution to the two-fluid incompressible Navier-Stokes-Maxwell system with solenoidal Ohm's law \eqref{TFINSMSO}. Further consider test functions $\left(\bar u, \bar j, \bar E, \bar B\right)\in C_c^\infty\left([0,\infty)\times\mathbb{R}^3\right)$ such that $\left\|\bar u\right\|_{L^\infty_{t,x}}<1$ and
	\begin{equation*}
		\begin{cases}
			\begin{aligned}
				\Div \bar j & = 0, & \Div \bar u & = 0, \\
				\d_t \bar E - \ROT \bar B &= -  \bar j, & \Div \bar E & = 0,
				\\
				\d_t \bar B + \ROT \bar E & = 0, & \Div \bar B & = 0.
			\end{aligned}
		\end{cases}
	\end{equation*}
	We define the acceleration operator by
	\begin{equation*}
		\mathbf{A}\left(\bar u, \bar j, \bar E, \bar B\right)
		=
		\begin{pmatrix}
			-\d_t \bar u -
			P\left(\bar u\cdot\nabla_x \bar u\right) + \mu\Delta_x \bar u
			+ \frac 12 P \left(\bar j \wedge \bar B\right)
			\\
			- \frac 1{2\sigma}\bar j + \frac 12 P\left(\bar E + \bar u\wedge \bar B\right)
		\end{pmatrix},
	\end{equation*}
	and the growth rate by
	\begin{equation*}
		\lambda(t) =
		\frac{2\left\|\nabla_{t,x}\bar u(t)\right\|_{L^\infty_x}}{1-\left\|\bar u(t)\right\|_{L^\infty_x}}
		+ \frac {\sqrt 2 \left\|\bar j(t)\right\|_{L^\infty_x}}{2\left(1-\left\|\bar u(t)\right\|_{L^\infty_x}\right)}.
	\end{equation*}
	
	Then, one has the stability inequality
	\begin{equation}\label{stability 4}
		\begin{aligned}
			\delta\CE(t) + & \int_0^t \delta\CD(s) e^{\int_s^t\lambda(\sigma)d\sigma}ds
			\\
			& \leq \delta\CE(0) e^{\int_0^t\lambda(s)ds}
			+\int_0^t
			\left[\int_{\mathbb{R}^3} \mathbf{A}\cdot
			\begin{pmatrix}
				u-\bar u \\ j-\bar j
			\end{pmatrix}
			dx\right](s)
			e^{\int_s^t\lambda(\sigma)d\sigma}ds,
		\end{aligned}
	\end{equation}
	where the modulated energy $\delta\CE$ and energy dissipation $\delta\CD$ are given by
	\begin{equation}\label{modulated defi 3}
		\begin{aligned}
			\delta\CE(t)
			& =
			\frac 12 \left\|\left(u-\bar u\right)(t)\right\|_{L^2_x}^2
			+ \frac 14 \left\|\left(E-\bar E\right)(t)\right\|_{L^2_x}^2
			+ \frac 14 \left\|\left(B-\bar B\right)(t)\right\|_{L^2_x}^2 \\
			& - \frac 12\int_{\mathbb{R}^3} \left(\left(E-\bar E\right)(t)\wedge \left(B-\bar B\right)(t)\right)\cdot\bar u(t) dx, \\
			\delta \CD(t)
			& =
			\mu \left\|\nabla_x (u-\bar u)(t)\right\|_{L^2_x}^2
			+ \frac 1{2\sigma} \left\|\left(j-\bar j\right)(t)\right\|_{L^2_x}^2.
		\end{aligned}
	\end{equation}
\end{prop}

\begin{proof}
	Following the proof of Proposition \ref{modulated energy estimate 2}, we consider first the identity
	\begin{equation*}
		\begin{aligned}
			\frac{d}{dt} & \int_{\mathbb{R}^3} \left(u\cdot\bar u
			+\frac 12 E\cdot \bar E+ \frac 12 B\cdot\bar B\right) dx
			+\int_{\mathbb{R}^3}2\mu \nabla_x u:\nabla_x\bar u
			+\frac 1{\sigma} j\cdot\bar j dx \\
			& =
			\int_{\mathbb{R}^3}
			(u-\bar u) \otimes (u-\bar u):\nabla_x\bar u dx
			\\
			& +\frac 12 \int_{\mathbb{R}^3}
			\left((j -\bar j)\wedge(B-\bar B)\right)\cdot\bar u
			+ \left((u-\bar u)\wedge(B-\bar B)\right)\cdot\bar j
			dx
			\\
			&
			- \int_{\mathbb{R}^3} \mathbf{A}\cdot
			\begin{pmatrix}
				u \\ j
			\end{pmatrix}
			dx.
		\end{aligned}
	\end{equation*}
	Note that this relation can be recovered by formally discarding all terms involving the charge density $n$ in \eqref{computation}.

	Then, expressing the modulated Lorentz force with a modulated Poynting vector as
	\begin{equation*}
		\begin{aligned}
			\partial_t\left(\left(E-\bar E\right)\wedge \left(B-\bar B\right)\right)
			&
			+\frac 12\nabla_x\left(\left|E-\bar E\right|^2+\left|B-\bar B\right|^2\right)
			\\
			& -\nabla_x\cdot\left(\left(E-\bar E\right)\otimes \left(E-\bar E\right)+\left(B-\bar B\right)\otimes \left(B-\bar B\right)\right)
			\\
			& =
			-\left(j-\bar j\right)\wedge \left(B-\bar B\right),
		\end{aligned}
	\end{equation*}
	we arrive at the relation
	\begin{equation*}
		\begin{aligned}
			\frac{d}{dt} \int_{\mathbb{R}^3} & \left(u\cdot\bar u
			+\frac 12 E\cdot \bar E+ \frac 12 B\cdot\bar B
			+\frac 12 \left(\left(E-\bar E\right)\wedge \left(B-\bar B\right)\right)\cdot\bar u\right) dx
			\\
			& - \int_{\mathbb{R}^3}
			\frac 12 \left(\left(E-\bar E\right)\wedge \left(B-\bar B\right)\right)\cdot \partial_t \bar u
			dx
			+\int_{\mathbb{R}^3}2\mu \nabla_x u:\nabla_x\bar u
			+\frac 1{\sigma}
			j\cdot\bar j
			dx
			\\
			& =
			\int_{\mathbb{R}^3}
			(u-\bar u) \otimes (u-\bar u):\nabla_x\bar u dx
			\\
			& -\frac 12
			\int_{\mathbb{R}^3}
			\left(\left(E-\bar E\right)\otimes \left(E-\bar E\right)+\left(B-\bar B\right)\otimes \left(B-\bar B\right)\right)
			: \nabla_x\bar u
			dx
			\\
			& +\frac 12\int_{\mathbb{R}^3}
			\bar j\cdot\left(\left(u-\bar u\right)\wedge\left(B-\bar B\right)\right)
			dx
			- \int_{\mathbb{R}^3} \mathbf{A}\cdot
			\begin{pmatrix}
				u \\ j
			\end{pmatrix}
			dx.
		\end{aligned}
	\end{equation*}

	On the whole, combining the preceding identity with the energy conservation law for test functions \eqref{energy test 2} and the energy decay imposed by the formal energy conservations from Proposition \ref{energy estimate 2}, we find the following modulated energy inequality~:
	\begin{equation*}
		\begin{aligned}
			\frac{d}{dt} & \delta\CE(t) + \delta\CD(t)
			\\
			& \leq
			-\int_{\mathbb{R}^3}
			(u-\bar u) \otimes (u-\bar u):\nabla_x\bar u dx
			- \int_{\mathbb{R}^3}
			\frac 12 \left(\left(E-\bar E\right)\wedge \left(B-\bar B\right)\right)\cdot \partial_t \bar u
			dx
			\\
			& +\frac 12
			\int_{\mathbb{R}^3}
			\left(\left(E-\bar E\right)\otimes \left(E-\bar E\right)+\left(B-\bar B\right)\otimes \left(B-\bar B\right)\right)
			: \nabla_x\bar u
			dx
			\\
			& -\frac 12\int_{\mathbb{R}^3}
			\bar j\cdot\left(\left(u-\bar u\right)\wedge\left(B-\bar B\right)\right)
			dx
			+ \int_{\mathbb{R}^3} \mathbf{A}\cdot
			\begin{pmatrix}
				u-\bar u \\ j-\bar j
			\end{pmatrix}
			dx.
		\end{aligned}
	\end{equation*}
	The next step consists in estimating the terms in the right-hand side above that are nonlinear in $(u,j,E,B)$ and to absorb the resulting expressions with the modulated energy $\delta\CE(t)$ by suitable uses of Young's inequality and Gr\"onwall's lemma. Thus, we obtain
	\begin{equation*}
		\begin{aligned}
			& \frac{d}{dt} \delta\CE(t) + \delta\CD(t) \\
			& \leq
			\left\|\nabla_{t,x} \bar u\right\|_{L^\infty_x}\left( \left\| u-\bar u\right\|_{L^2_x}^2
			+ \frac 12\left\|E-\bar E\right\|_{L^2_x}^2+\frac 12\left\|B-\bar B\right\|_{L^2_x}^2\right)
			\\
			& + \frac 12 \left\|\bar j\right\|_{L^\infty_x}\left\|B-\bar B\right\|_{L^2_x}
			\left\|u-\bar u\right\|_{L^2_x}
			+
			\int_{\mathbb{R}^3} \mathbf{A}\cdot
			\begin{pmatrix}
				u-\bar u \\ j-\bar j
			\end{pmatrix}
			dx
			\\
			& \leq
			\left\|\nabla_{t,x} \bar u\right\|_{L^\infty_x}\left( \left\| u-\bar u\right\|_{L^2_x}^2
			+ \frac 12\left\|E-\bar E\right\|_{L^2_x}^2+\frac 12\left\|B-\bar B\right\|_{L^2_x}^2\right)
			\\
			& + \frac{\sqrt 2}{4}\left\|\bar j\right\|_{L^\infty_x}
			\left(
			\left\|u-\bar u\right\|_{L^2_x}^2
			+
			\frac 12\left\|B-\bar B\right\|_{L^2_x}^2\right)
			+
			\int_{\mathbb{R}^3} \mathbf{A}\cdot
			\begin{pmatrix}
				u-\bar u \\ j-\bar j
			\end{pmatrix}
			dx.
		\end{aligned}
	\end{equation*}
	Hence, further noticing that
	\begin{equation*}
		\begin{aligned}
			\left(1-\left\|\bar u\right\|_{L^\infty_{x}}\right)
			& \left(
			\frac 12\left\|u-\bar u\right\|_{L^2_x}^2
			+\frac 14
			\left\|E-\bar E\right\|_{L^2_x}^2+\frac 14\left\|B-\bar B\right\|_{L^2_x}^2\right)
			\\
			& \leq
			\frac12 \left\|u-\bar u\right\|_{L^2_x}^2
			+
			\frac 14\left\|E-\bar E\right\|_{L^2_x}^2+\frac 14\left\|B-\bar B\right\|_{L^2_x}^2
			\\
			& -\frac 12\left\|\bar u\right\|_{L^\infty_{x}}\int_{\mathbb{R}^3}\left|E-\bar E\right|\left|B-\bar B\right|dx
			\\
			& \leq
			\frac 12 \left\|u-\bar u\right\|_{L^2_x}^2
			+
			\frac 14\left\|E-\bar E\right\|_{L^2_x}^2+\frac 14\left\|B-\bar B\right\|_{L^2_x}^2
			\\
			& -\frac 12\int_{\mathbb{R}^3}\left(\left(E-\bar E\right)\wedge\left(B-\bar B\right)\right)\cdot\bar udx
			\\
			& = \delta\mathcal{E}(t),
		\end{aligned}
	\end{equation*}
	we find, since $\left\|\bar u\right\|_{L^\infty_{t,x}}<1$, that
	\begin{equation*}
		\frac{d}{dt}\delta\CE(t)
		+\delta\CD(t)
		\leq
		\lambda(t)
		\delta\CE(t)
		+ \int_{\mathbb{R}^3} \mathbf{A}\cdot
		\begin{pmatrix}
			u-\bar u \\ j-\bar j
		\end{pmatrix}
		dx,
	\end{equation*}
	which concludes the proof of the proposition with a direct application of Gr\"onwall's lemma.
\end{proof}

The preceding proposition provides another weak-strong stability property for the two-fluid incompressible Navier-Stokes-Maxwell system with solenoidal Ohm's law \eqref{TFINSMSO}. Indeed, the stability inequality \eqref{stability 4} essentially implies that a solution $(\bar u, \bar j,\bar E,\bar B)$ of \eqref{TFINSMSO} such that $\bar u\in L^\infty_{t,x}$, $\nabla_{t,x}\bar u\in L^1_tL^\infty_x$, $\bar j\in L^1_tL^\infty_x$ and $\left\|\bar u\right\|_{L^\infty_{t,x}}<1$, if it exists, is unique in the whole class of weak solutions in the energy space, for any given initial data.

As in Proposition \ref{stability poynting}, the condition $\left\|\bar u\right\|_{L^\infty_{t,x}}<1$ in Proposition \ref{stability poynting 2} is physically relevant, for it imposes that the modulus of the bulk velocity $\bar u$ remains everywhere and at all times below the speed of light. More precisely, keeping track of the relevant physical constants in the formal derivations of Chapter \ref{formal-chap}, we see that the system \eqref{TFINSMSO} can be recast as
\begin{equation*}
	\begin{cases}
		\begin{aligned}
			\d_t u +
			u\cdot\nabla_x u - \mu\Delta_x u
			& = -\nabla_x p+
			\frac 12 j \wedge B , & \Div u & = 0,\\
			\frac 1c\d_t E - \ROT B &= -  j, & \Div E & = 0,
			\\
			\frac 1c\d_t B + \ROT E & = 0, & \Div B & = 0, \\
			j & = \sigma\left(- \nabla_x \bar p + c E + u\wedge B\right), & \Div j & = 0,
		\end{aligned}
	\end{cases}
\end{equation*}
where the constant $c>0$ denotes the speed of light. Then, applying the proof of Proposition \ref{stability poynting 2} to the preceding system, we arrive at a stability inequality valid under the restriction that the bulk velocity remains bounded by the speed of light $\left\|\bar u\right\|_{L^\infty_{t,x}}<c$, which is natural.

\bigskip

Following the previous developments, it is also possible to use now the stability inequality \eqref{stability 4} from Proposition \ref{stability poynting 2} to define another notion of dissipative solutions for the two-fluid incompressible Navier-Stokes-Maxwell system with solenoidal Ohm's law \eqref{TFINSMSO}, whose existence is then established by reproducing the arguments from Theorem \ref{dissipative 3}. The only argument from the proof of Theorem \ref{dissipative 3} that needs special care in order to conclude the existence of dissipative solutions resides in the weak lower semi-continuity of the modulated energy $\delta\mathcal{E}(t)$ defined by \eqref{modulated defi 3}, which we have already established in \eqref{convex wedge}.

%% file: convergence0.tex
\chapter{Two typical regimes}\label{convergence results}%{Two typical convergence results}

We will now focus on two specific regimes which are critical, in the sense that all the formal asymptotics mentioned in Chapter \ref{formal-chap} can be rigorously obtained by similar or even simpler arguments.

The first scaling we will investigate here is the one leading from the one species Vlasov-Maxwell-Boltzmann equations \eqref{scaledVMB} to the incompressible quasi-static Navier-Stokes-Fourier-Maxwell-Poisson system \eqref{NSFMP}. More precisely, we will set $\alpha=\eps$, $\beta=\eps$ and $\gamma=\eps$ in \eqref{scaledVMB}. As discussed in Section \ref{stability existence 1}, the resulting limiting model is then very similar to the incompressible Navier-Stokes equations and, thus, the usual methods of hydrodynamic limits will apply. We shall focus specifically on the influence of the electromagnetic field, which induces numerous technical complications.

The second regime we will study is more singular since the magnetic forcing is much stronger. Specifically, we will consider the scaling leading from the two species Vlasov-Maxwell-Boltzmann equations \eqref{scaled VMB two species} to the two-fluid incompressible Navier-Stokes-Fourier-Maxwell system with Ohm's law \eqref{TFINSFMO} in the case of strong interspecies collisions, or to the two-fluid incompressible Navier-Stokes-Fourier-Maxwell system with solenoidal Ohm's law \eqref{TFINSFMSO} in the case of weak interspecies collisions. More precisely, we will set $\alpha=\delta\eps$, $\beta=\delta$ and $\gamma=1$ in \eqref{scaled VMB two species}, with $\frac\delta\eps$ unbounded. Actually, as discussed in Section \ref{stability existence 2}, the corresponding limiting models \eqref{TFINSFMO} and \eqref{TFINSFMSO} are not stable under weak convergence in the energy space and, thus, share more similarities with the three-dimensional incompressible Euler equations. So will our proofs of hydrodynamic convergence in this setting.

\bigskip

All along this second part on rigorous hydrodynamic convergence proofs, we will consider renormalized solutions, whose definition we recall below in Section \ref{renorm sol cond}, of the Vlasov-Maxwell-Boltzmann systems for any number species. In fact, their existence is not established, which is precisely the reason why the convergence results presented here are deemed conditional, and remains a challenging open problem of outstanding difficulty.

Loosely speaking, the specific complexity of the Vlasov-Maxwell-Boltzmann system originates in the nonlinear coupling of the Vlasov-Boltzmann equation with a hyperbolic system, namely Maxwell's equations. This essential difficulty remains ubiquitous in our analysis of its hydrodynamic limits and is passed on to the most singular asymptotic models present in our work, such as the systems \eqref{TFINSFMO} and \eqref{TFINSFMSO}, whose well-posedness is not fully understood (see Section \ref{stability existence 2}) and contains very challenging open questions, as well.

% ==========================
% = Renormalized solutions =
% ==========================

\section{Renormalized solutions}\label{renorm sol cond}

We are now going to recall the notion of renormalized solutions for the Vlasov-Maxwell-Boltzmann systems \eqref{scaledVMB}
\begin{equation}\label{scaledVMB bis}
	\begin{cases}
		\begin{aligned}
			\d_t f + v \cdot \nabla_x f + \left( E + v \wedge B \right) \cdot \nabla_v f &= Q(f,f), \\
			\d_t E - \ROT B &= - \int_{\mathbb{R}^3} fv dv, \\
			\d_t B + \ROT E& = 0, \\
			\DIV E &=\int_{\mathbb{R}^3} fdv -1, \\
			\DIV B &=0,
		\end{aligned}
	\end{cases}
\end{equation}
and \eqref{scaled VMB two species}
\begin{equation}\label{scaled VMB two species bis}
	\begin{cases}
		\begin{aligned}
			\d_t f^\pm + v \cdot \nabla_x f^\pm \pm \left( E + v \wedge B \right) \cdot \nabla_v f^\pm &= Q(f^\pm,f^\pm) + Q(f^\pm,f^\mp), \\
			\d_t E - \ROT B &= -  \int_{\mathbb{R}^3} \left(f^+-f^-\right)v dv, \\
			\d_t B + \ROT E& = 0, \\
			\DIV E &=\int_{\mathbb{R}^3} \left(f^+-f^-\right)dv, \\
			\DIV B &=0,
		\end{aligned}
	\end{cases}
\end{equation}
where we have discarded the free parameters.

\subsection{The Vlasov-Boltzmann equation}\label{vlasov}

Let us focus first on the simpler Vlasov-Boltzmann equation~:
\begin{equation}\label{boltzmann}
	\partial_t f + v\cdot\nabla_x f + F\cdot\nabla_v f = Q\left(f,f\right),
\end{equation}
with a given force field $F(t,x,v)$ satisfying, at least,
\begin{equation*}
		F,\nabla_v\cdot F\in L^1_\mathrm{loc}\left(dtdx;L^1\left(M^\alpha dv\right)\right)
		\qquad\text{for all }\alpha>0.
\end{equation*}
The above conditions on the force field are minimal requirements so that it is possible to define renormalized solutions of \eqref{boltzmann} (see definition below). We will, however, further restrict the range of applicability of force fields~:
\begin{itemize}
	\item we assume that $\nabla_v\cdot F =0$, so that the local conservation of mass is verified~;
	\item we assume that $F\cdot v =0$, so that the global Maxwellian $M(v)$ is an equilibrium state of \eqref{boltzmann}.
\end{itemize}

Renormalized solutions of \eqref{boltzmann} are known to exist since the late eighties, thanks to DiPerna and Lions \cite{diperna} (at least for the Boltzmann equation, i.e.\ for the case $F=0$). We are going to briefly describe their derivation, their limitations and emphasize the main mathematical difficulties preventing their construction for the above Vlasov-Maxwell-Boltzmann systems.

\bigskip

Throughout this work, we are interested in the fluctuations of a density $f(t,x,v)$ around a global normalized Maxwellian $M(v)$, we will therefore conveniently employ the density $G(t,x,v)$ defined by $f=MG$. In this notation, the Vlasov-Boltzmann equation \eqref{boltzmann} reads
\begin{equation}\label{boltzmann2}
	\partial_t G + v\cdot\nabla_x G + F\cdot\nabla_v G = \mathcal{Q}\left(G,G\right),
\end{equation}
where we denote
\begin{equation*}
	\cQ(G,H) =\frac1M Q(MG,MH).
\end{equation*}

Thus, DiPerna and Lions formulated in \cite{diperna} the first theory yielding global solutions to the Boltzmann equation \eqref{boltzmann2}, with $F=0$, for large initial data $G(0,x,v)=G^{\mathrm{in}}(x,v)\geq0$. Their construction heavily relied on a new notion of solutions, namely the renormalized solutions.

Recall that we utilize the prefixes $\textit{w-}$ or $\textit{w$^*$-}$ to express that a given space is endowed with its weak or weak-$*$ topology, respectively.

\begin{defi}
	We say that a nonlinearity $\beta\in C^1\left([0,\infty);\mathbb{R}\right)$ is an admissible renormalization if it satisfies, for some $C>0$,
	\begin{equation*}
		\left|\beta'(z)\right|\leq \frac{C}{\left(1+z\right)^\frac{1}{2}} \qquad \text{for all }z\geq 0.
	\end{equation*}
	
	A density function $f(t,x,v)=MG(t,x,v)\geq0$, where $(t,x,v)\in [0,\infty)\times\mathbb{R}^3\times\mathbb{R}^3$, such that
	\begin{equation}\label{membership0}
		G\in C\left([0,\infty);\textit{w-}L^1_{\mathrm{loc}}\left(dxdv\right)\right)
		\cap L^\infty\left([0,\infty),dt;L^1_{\mathrm{loc}}\left(dx;L^1\left((1+|v|^2)M dv\right)\right)\right),
	\end{equation}
	% and
	% \begin{equation*}
	% 	\sup_{t\geq 0}H(f)=\sup_{t\geq 0}H(f|M)
	% 	=\sup_{t\geq 0} \int_{\mathbb{R}^3\times\mathbb{R}^3}\left(G\log G-G+1\right)M dxdv < \infty,
	% \end{equation*}
	is a \textbf{renormalized solution of the Vlasov-Boltzmann equation \eqref{boltzmann2}} if it solves
	\begin{equation}\label{boltzmann renorm}
		\partial_t\beta(G)+v\cdot\nabla_x\beta(G)
		+ F\cdot\nabla_v\beta(G)
		= \beta'(G)\mathcal{Q}(G,G)
	\end{equation}
	in the sense of distributions for any admissible renormalization, and satisfies the entropy inequality, for all $t>0$,
	\begin{equation*}
		H(f(t)) + \int_0^t \int_{\mathbb{R}^3} D(f(s)) dx ds\leq H(f^{\mathrm{in}})<\infty,
	\end{equation*}
	where $f^\mathrm{in}=MG^\mathrm{in}$ is the initial value of $f=MG$ and the relative entropy $H(f)=H(f|M)$ is defined in \eqref{def H}, while the entropy dissipation $D(f)$ is defined in \eqref{def D}.
\end{defi}

Note that the renormalized collision operator $\beta'(G)\mathcal{Q}(G,G)$ is well-defined in $L^1_\mathrm{loc}\left(dtdx;L^1\left(M^\alpha dv\right)\right)$, with $\alpha>0$, for any admissible renormalization, any function in \eqref{membership0} and any integrable cross-section $b(z,\sigma)\in L^1_{\mathrm{loc}}\left(\mathbb{R}^3\times\mathbb{S}^2\right)$ satisfying the so-called DiPerna-Lions assumption
\begin{equation}\label{diperna lions assumption}
	\lim_{|v|\rightarrow\infty}\frac{1}{|v|^2}\int_{K\times\mathbb{S}^2}b(v-v_*,\sigma) dv_*d\sigma = 0,
\end{equation}
for any compact subset $K\subset\mathbb{R}^3$.

Indeed, it is possible to show directly from \eqref{diperna lions assumption} that (see \cite{arsenio3}, for instance, for more details), for any $\alpha>0$,
\begin{equation*}
	\lim_{|v|\rightarrow\infty}\frac{1}{|v|^2}\int_{\mathbb{R}^3\times\mathbb{S}^2}b(v-v_*,\sigma) M_*^\alpha dv_*d\sigma = 0.
\end{equation*}
Therefore, considering first non-negative renormalizations satisfying $0 \leq \beta'(z) \leq \frac{C}{1+z}$, the renormalized loss part $\beta'(G)\mathcal{Q}^-(G,G)$ is easily estimated as
\begin{equation*}
	\begin{aligned}
		\int_{\mathbb{R}^3} & \beta'(G) \mathcal{Q}^-(G,G) M^\alpha dv \\
		& = \int_{\mathbb{R}^3} G_* (1+|v_*|^2) M_* \left[\frac{1}{1+|v_*|^2} \int_{\mathbb{R}^3\times\mathbb{S}^2} \beta'(G) G b(v-v_*,\sigma)
		M^\alpha dv d\sigma \right] dv_* \\
		& \leq C \left\|G\right\|_{L^1\left((1+|v|^2)Mdv\right)},
	\end{aligned}
\end{equation*}
while the renormalized gain term $\beta'(G)\mathcal{Q}^+(G,G)$ is well-defined in $L^1_\mathrm{loc}\left(dtdx;L^1\left(M^\alpha dv\right)\right)$ by the renormalized Vlasov-Boltzmann equation \eqref{boltzmann renorm} because it is the only unestimated expression remaining and it is non-negative. These controls are easily extended to signed renormalizations satisfying $\left|\beta'(z)\right| \leq \frac{C}{1+z}$, for the Vlasov-Boltzmann equation \eqref{boltzmann renorm} is linear with respect to renormalizations so that we may decompose $\beta'(z)$ with respect to its positive and negative parts.

Alternatively and as was originally performed in \cite{diperna}, we could also use the elementary inequality \eqref{weaker young r}, setting $z=\frac{G'G_*'}{GG_*}-1$ and $y=\log K$, with $K>1$, which implies that
\begin{equation*}
	\begin{aligned}
		\int_{\mathbb{R}^3}\beta'(G)&\mathcal{Q}^+ (G,G)Mdv\leq K\int_{\mathbb{R}^3}\beta'(G)\mathcal{Q}^-(G,G)Mdv\\
		&+\frac{1}{\log K}\int_{\mathbb{R}^3\times\mathbb{R}^3\times\mathbb{S}^2}\beta'(G)\left(G'G_*'-GG_*\right)\log\left(\frac{G'G_*'}{GG_*}\right)MM_*b dvdv_*d\sigma,
	\end{aligned}
\end{equation*}
to claim that the gain part belongs to $L^1_\mathrm{loc}\left(dtdx;L^1\left(M dv\right)\right)$, since it is natural to control the entropy dissipation term above.

Finally, it is possible to extend the definition of the renormalized collision operator $\beta'(G)\mathcal{Q}(G,G)$ to all admissible renormalizations by decomposing the renormalized collision integrand as
\begin{equation*}
	\begin{aligned}
		\beta'(G)\left(G'G_*'-GG_*\right) & = \beta'(G)\left(\sqrt{G'G_*'}-\sqrt{GG_*}\right)^2 \\
		& + 2\beta'(G)\sqrt{GG_*}\left(\sqrt{G'G_*'}-\sqrt{GG_*}\right),
	\end{aligned}
\end{equation*}
and noticing that
\begin{equation*}
	\begin{aligned}
		\int _{\mathbb{R}^3 \times \mathbb{R}^3 \times \mathbb{S}^2}
		& \left(\sqrt{f'f_*'}-\sqrt{ff_*}\right)^2
		b(v-v_*,\sigma)
		dvdv_* d\sigma \\
		& \leq
		\frac14 \int _{\mathbb{R}^3 \times \mathbb{R}^3 \times \mathbb{S}^2}
		\left(f'f'_* - ff_*\right)  \log \left({f'f'_*\over ff_*}\right) b(v-v_*,\sigma)
		dvdv_* d\sigma = D(f),
	\end{aligned}
\end{equation*}
which follows from the elementary inequality \eqref{sqrt dissip}.

Thus, by a solution $G$ of the renormalized equation \eqref{boltzmann renorm}, we naturally mean that $G$ should satisfy, for every $\alpha>0$ and any non-negative test functions $\rho(t,x)\in C_c^\infty\left([0,\infty)\times\mathbb{R}^3\right)$ and $\varphi(v)\in W^{1,\infty}\left(dv\right)$, that
\begin{equation*}
	\begin{gathered}
		\begin{aligned}
			-\int_{\mathbb{R}^3\times\mathbb{R}^3} & \beta \left(G^\mathrm{in}\right) \rho(0,x)\varphi(v) M^\alpha dxdv\\
			&-\int_{[0,\infty)\times\mathbb{R}^3\times\mathbb{R}^3} \beta\left(G\right)
			\left(\partial_t+v\cdot\nabla_x + F\cdot\nabla_v\right)\left[\rho(t,x)\varphi(v)
			M^\alpha\right] dtdxdv \\
			% &-\int_{[0,\infty)\times\mathbb{R}^3\times\mathbb{R}^3}
			% 	\left[\left(\nabla_v\cdot F\right)\beta\left(G\right)
			% 	+F\cdot v\beta'(G)G\right]
			% 	\rho(t,x)\varphi(v)
			% 	M^\alpha dtdxdv \\
			& =
			\int_{[0,\infty)\times\mathbb{R}^3\times\mathbb{R}^3} \beta'\left(G\right)\mathcal{Q}\left(G,G\right)\rho(t,x)\varphi(v) M^\alpha dtdxdv.
		\end{aligned}
	\end{gathered}
\end{equation*}

The following theorem is a modern formulation of the existence result found in \cite{diperna}. The existence of renormalized solutions for Vlasov-Boltzmann systems where the force field derives from a self-induced potential, such as the Vlasov-Poisson-Boltzmann system, has been established in \cite{lions3}, while the study of renormalized solutions close to Maxwellian equilibrium has been performed in \cite{lions8}.

\begin{thm}[\cite{diperna, diperna3}]\label{classical renorm sol}
	Let $b(z,\sigma)$ be a locally integrable collision kernel satisfying the DiPerna-Lions assumption \eqref{diperna lions assumption} and $F(t,x,v)\in L^1_\mathrm{loc}(dtdxdv)$ a given force field such that
	\begin{equation}\label{force field hypotheses}
		\nabla_v\cdot F=0,\qquad F\cdot v=0\qquad\text{and}\qquad F\in
		L^1_{\mathrm{loc}}\left(dt;W^{1,1}_{\mathrm{loc}}\left(dxdv\right)\right).
	\end{equation}
	
	Then, for any initial condition $f^\mathrm{in}=MG^{\mathrm{in}}\in L^1_\mathrm{loc}\left(dx;L^1\left((1+|v|^2) dv\right)\right)$ such that $f^\mathrm{in}=MG^\mathrm{in}\geq 0$ and
	\begin{equation*}
		H(f^\mathrm{in})=H(f^\mathrm{in}|M)=\int_{\mathbb{R}^3\times\mathbb{R}^3}\left(G^{\mathrm{in}}\log G^{\mathrm{in}}-G^{\mathrm{in}}+1\right)M dxdv < \infty,
	\end{equation*}
	there exists a renormalized solution $f(t,x,v)=MG(t,x,v)$ to the Vlasov-Boltzmann equation \eqref{boltzmann2}. Moreover, it satisfies the local conservation of mass
	\begin{equation*}
		\partial_t\int_{\mathbb{R}^3}fdv + \nabla_x\cdot\int_{\mathbb{R}^3}fvdv = 0,
	\end{equation*}
	and the global entropy inequality, for any $t\geq 0$,
	\begin{equation}\label{entropy inequality boltzmann}
		H(f(t)) + \int_0^t \int_{\mathbb{R}^3} D(f(s)) dx ds\leq H(f^{\mathrm{in}}).
	\end{equation}
\end{thm}

The proof of the above theorem follows the usual steps found in the analysis of weak solutions of partial differential equations, that is to say, solving an approximate truncated equation, establishing uniform a priori estimates and the weak compactness of the approximate solutions, and finally passing to the limit (by showing the weak stability of nonlinear terms) and, thus, recovering the original equation. It is often the case that these steps reduce to the study of the crucial weak stability of solutions. Thus, for the Vlasov-Boltzmann equation \eqref{boltzmann2}, the above theorem naturally follows from the weak stability of renormalized solutions, or, in other words, from the weak stability of weak solutions of the renormalized equation \eqref{boltzmann renorm} satisfying the uniform bounds provided by the entropy inequality \eqref{entropy inequality boltzmann}.

DiPerna and Lions showed the weak stability of the Boltzmann equation, i.e.\ when $F=0$, in \cite{diperna} and refined their result in \cite{diperna3} by establishing the entropy inequality \eqref{entropy inequality boltzmann}. Note that, since we are assuming $\nabla_v\cdot F=0$ and $F\cdot v=0$, the entropy inequality \eqref{entropy inequality boltzmann} easily follows from formal estimates on the Vlasov-Boltzmann equation \eqref{boltzmann2}, even when $F\neq 0$. Later, Lions improved the method of proof in \cite{lions, lions2, lions3}. We briefly explain now Lions' strategy, which relies on velocity averaging lemmas, heavy renormalization techniques and, most importantly, on the compactifying (even regularizing, in some cases) effect of the gain term $Q^+\left(f,f\right)$ of the collision operator.

To this end, let us consider a sequence $\left\{f_k\right\}_{k\in\mathbb{N}}$ of actual renormalized solutions to \eqref{boltzmann}, with initial data $\left\{f_k^\mathrm{in}\right\}_{k\in\mathbb{N}}$, which converges weakly (at least in $L^1_\mathrm{loc}$, say) as $k\rightarrow\infty$ to $f^\mathrm{in}$. We further assume that the initial data satisfies the following strong entropic convergence
\begin{equation*}
	\lim_{k\rightarrow\infty}
	H\left(f_k^{\mathrm{in}}\right)
	=
	H\left(f^{\mathrm{in}}\right),
\end{equation*}
so that the entropy inequality is uniformly satisfied
\begin{equation}\label{entropy inequality k Young 0}
		H\left(f_k(t)\right)
		+ \int_0^t\int_{\mathbb{R}^3} D\left(f_k(s)\right) dxds
		\leq H\left(f_k^{\mathrm{in}}\right).
\end{equation}

Notice that a uniform bound on the entropies $H\left(f_k(t)\right)$ yields, with a direct application of the elementary Young inequality \eqref{Young h}, a uniform bound on $f_k(t,x,v)$ in $L^\infty_{\mathrm{loc}}\left(dt;L^1_{\mathrm{loc}}\left(dx;L^1\left((1+v^2)dv\right)\right)\right)$. Moreover, it is possible to show, with a slightly more refined application of the Young inequality \eqref{Young h} with the Dunford-Pettis compactness criterion (see \cite{royden} and Section \ref{control rel entropy bound} for details), that the $f_k$'s are in fact weakly relatively compact in $L^1_\mathrm{loc}\left(dtdx;L^1\left(dv\right)\right)$. Therefore, up to extraction, we may assume that the sequence $\left\{f_k\right\}_{k\in\mathbb{N}}$ converges weakly, as $k\to \infty$, to some $f$ in $L^1_\mathrm{loc}\left(dtdx;L^1\left(dv\right)\right)$.

Similarly, uniform bounds on the nonlinear terms
\begin{equation*}
	\frac{Q^\pm\left(f_k,f_k\right)}{1+\delta f_k*_v\left[\int_{\mathbb{S}^2}b(\cdot,\sigma)d\sigma\right]}
	\qquad\text{and}\qquad
	\beta'\left(G_k\right)Q^\pm\left(f_k,f_k\right),
\end{equation*}
where $f_k=MG_k$, for any $\delta>0$ and any admissible nonlinearity $\beta(z)\in C^1\left([0,\infty);\mathbb{R}\right)$, are easily obtained from \eqref{entropy inequality k Young 0} through the standard use of the elementary inequality \eqref{weaker young r}, setting $z=\frac{f_k'f_{k*}'}{f_kf_{k*}}-1$ and $y=\log K$, with $K>1$, which yields
\begin{equation*}
	\begin{aligned}
		f_k'f_{k*}' \leq Kf_kf_{k*}
		+ \frac{1}{\log K}
		{(f_k'f_{k*}'-f_kf_{k*})\log\left(\frac{f_k'f_{k*}'}{f_kf_{k*}}\right)}.
	\end{aligned}
\end{equation*}
The above functional inequality further implies the weak compactness of the above nonlinear terms, thanks to the Dunford-Pettis compactness criterion (see \cite{royden}).

At this point, using the convexity methods from \cite{diperna3}, one can already establish the limiting entropy inequality \eqref{entropy inequality boltzmann}, passing to the limit in \eqref{entropy inequality k Young 0}.

Next, since each $f_k=MG_k$ is a weak solution of the renormalized equation \eqref{boltzmann renorm}, it is also possible to show, with a standard use of velocity averaging lemmas (one can use the results from \cite{diperna5}, for instance, treating $F\cdot\nabla_v\beta(G_k)=\nabla_v\cdot\left(F\beta(G_k)\right)$ as a source term), that, for any admissible nonlinearity $\beta(z)\in C^1\left([0,\infty);\mathbb{R}\right)$ and any cutoff $\varphi(v)\in C^\infty_c\left(\mathbb{R}^3\right)$,
\begin{equation}\label{velocity averages 0}
		\int_{\mathbb{R}^3} \beta(G_k)(t,x,v)\varphi(v)dv\quad
		\text{is relatively compact in }L^1_\mathrm{loc}\left(dtdx\right),
\end{equation}
and, up to further extraction of subsequences,
\begin{equation}\label{compactness loss 0}
	\begin{gathered}
		f_k*_v\left[\int_{\mathbb{S}^2}b(\cdot,\sigma)d\sigma\right]
		\rightarrow
		f*_v\left[\int_{\mathbb{S}^2}b(\cdot,\sigma)d\sigma\right]\\
		\text{as }k\rightarrow\infty \text{ in }L^1_\mathrm{loc}\left(dtdxdv\right) \text{ and almost everywhere,}
	\end{gathered}
\end{equation}
where $b(z,\sigma)$ may in fact be replaced by any collision kernel satisfying the DiPerna-Lions assumption \eqref{diperna lions assumption}. In particular, it follows that
\begin{equation*}
	\begin{gathered}
		\frac{Q^\pm\left(f_k,f_k\right)}{1+\delta f_k*_v\left[\int_{\mathbb{S}^2}b(\cdot,\sigma)d\sigma\right]}
		\varphi(t,x)
		\rightarrow
		\frac{Q^\pm\left(f,f\right)}{1+\delta f*_v\left[\int_{\mathbb{S}^2}b(\cdot,\sigma)d\sigma\right]}
		\varphi(t,x)\\
		\text{as }k\rightarrow\infty \text{ in }\textit{w-}L^1\left(dtdxdv\right),
	\end{gathered}
\end{equation*}
for any $\varphi(t,x)\in C_c^\infty\left([0,\infty)\times\mathbb{R}^3\right)$.

Lions showed in \cite{lions}, using Fourier integral operators, that the weak convergence of $f_k$ towards $f$ in $L^1_{\mathrm{loc}}(dtdxdv)$, the strong relative compactness of the velocity averages \eqref{velocity averages 0}-\eqref{compactness loss 0} and the uniform bounds from the entropy inequality \eqref{entropy inequality k Young 0} are sufficient to imply that, up to extraction of a subsequence, for every $\delta>0$,
\begin{equation*}
	\begin{gathered}
		\frac{Q^+\left(f_k,f_k\right)}{1+\delta f_k*_v\left[\int_{\mathbb{S}^2}b(\cdot,\sigma)d\sigma\right]}
		\varphi(t,x)
		\rightarrow
		\frac{Q^+\left(f,f\right)}{1+\delta f*_v\left[\int_{\mathbb{S}^2}b(\cdot,\sigma)d\sigma\right]}
		\varphi(t,x)
		\\
		\text{as }k\rightarrow\infty \text{ in }L^1\left(dtdxdv\right) \text{ and almost everywhere,}
	\end{gathered}
\end{equation*}
for any $\varphi(t,x)\in C_c^\infty\left([0,\infty)\times\mathbb{R}^3\right)$. Therefore, it holds in particular that
\begin{equation}\label{compactness gain 0}
		Q^+\left(f_k,f_k\right)
		\rightarrow
		Q^+\left(f,f\right)
		\text{ almost everywhere.}
\end{equation}

Following \cite{lions3}, we fix now the specific renormalization $\beta_\delta(z)=\frac{z}{1+\delta z}$, for any $0<\delta<1$, and we assume, without loss of generality, up to extraction of subsequences, that, as $k\rightarrow \infty$,
\begin{equation*}
	\begin{aligned}
		\beta_\delta(G_k) =\frac{G_k}{1+\delta G_k}  & \rightharpoonup \beta_\delta \leq \beta_\delta(G), \\
		\beta_\delta'(G_k) =\frac{1}{\left(1+\delta G_k\right)^2} & \rightharpoonup h_\delta \geq \beta'_\delta(G), \\
		\beta_\delta'(G_k)G_k =\frac{G_k}{\left(1+\delta G_k\right)^2}% \hspace{10mm} & \\
		=\beta_\delta(G_k)\left(1-\delta\beta_\delta(G_k)\right)
		& \rightharpoonup g_\delta
		\leq \beta_\delta\left(1-\delta\beta_\delta\right) ,
	\end{aligned}
\end{equation*}
in $\textit{w$^*$-}L^\infty_\mathrm{loc}(dtdxdv)$. Therefore, passing to the limit in \eqref{boltzmann renorm}, we obtain, in view of the strong convergences \eqref{compactness loss 0} and \eqref{compactness gain 0},
\begin{equation}\label{boltzmann renorm limit}
	\partial_t\beta_\delta+v\cdot\nabla_x\beta_\delta
	+ F\cdot\nabla_v\beta_\delta
	=
	h_\delta \mathcal{Q}^+(G,G)
	-g_\delta f*_v\left[\int_{\mathbb{S}^2}b(\cdot,\sigma)d\sigma\right],
\end{equation}
where the last term $h_\delta \mathcal{Q}^+(G,G)$ is well-defined in $L^1_\mathrm{loc}\left(dtdxdv\right)$ by its mere non-negativeness.

Note that, for any $\lambda>0$, choosing $K>0$ large enough so that, by equi-integrability of the $G_k$'s,
\begin{equation*}
	\sup_{k\in\mathbb{N}}\left\|G_k\mathds{1}_{\left\{G_k\geq K\right\}}\right\|_{L^1_\mathrm{loc}(dtdxdv)}\leq\lambda,
\end{equation*}
we find
\begin{equation*}
	\begin{aligned}
		\left\|G - \beta_\delta\right\|_{L^1_\mathrm{loc}(dtdxdv)}
		&\leq
		\liminf_{k\rightarrow\infty}\left\|G_k - \beta_\delta(G_k)\right\|_{L^1_\mathrm{loc}(dtdxdv)} \\
		&\leq
		\frac{\delta K}{1+\delta K} \sup_{k\in\mathbb{N}}\left\|G_k\right\|_{L^1_\mathrm{loc}(dtdxdv)}
		+\lambda,
	\end{aligned}
\end{equation*}
and
\begin{equation*}
	\begin{aligned}
		\left\|G - g_\delta\right\|_{L^1_\mathrm{loc}(dtdxdv)}
		&\leq
		\liminf_{k\rightarrow\infty}\left\|G_k - \beta_\delta'(G_k)G_k\right\|_{L^1_\mathrm{loc}(dtdxdv)} \\
		&\leq
		\frac{\delta^2 K^2+2\delta K}{\left(1+\delta K\right)^2} \sup_{k\in\mathbb{N}}\left\|G_k\right\|_{L^1_\mathrm{loc}(dtdxdv)}
		+\lambda.
	\end{aligned}
\end{equation*}
Hence, by the arbitrariness of $\lambda>0$,
\begin{equation*}
	\begin{aligned}
		\lim_{\delta\rightarrow 0}\left\|G - \beta_\delta\right\|_{L^1_\mathrm{loc}(dtdxdv)} & = 0, \\
		\lim_{\delta\rightarrow 0}\left\|G - g_\delta\right\|_{L^1_\mathrm{loc}(dtdxdv)} & = 0.
	\end{aligned}
\end{equation*}
Similarly, it is readily seen that, for any $1\leq p<\infty$,
\begin{equation*}
	\lim_{\delta\rightarrow 0}\left\|1 - h_\delta\right\|_{L^p_\mathrm{loc}(dtdxdv)} = 0.
\end{equation*}
Finally, notice that $\beta_\delta$, $g_\delta$ and $h_\delta$ are all increasing as $\delta$ vanishes. Hence, as $\delta\rightarrow 0$, both $\beta_\delta$ and $g_\delta$ converge towards $G$ almost everywhere, while $h_\delta$ converges toward a constant almost everywhere.

Now comes a fundamental idea of Lions from \cite{lions2, lions3}, which will be of particular interest to us and which has numerous qualitative consequences on renormalized solutions. This key idea consists in renormalizing equation \eqref{boltzmann renorm limit} over again according to the following simple yet crucial lemma from \cite{diperna6}.

\begin{lem}[{\cite[Theorem II.1, p. 516]{diperna6}}]\label{renorm transport}
	Let $f(t,x)\in L^\infty\left([0,T];L^p_\mathrm{loc}\left(\mathbb{R}^n\right)\right)$, with $1< p\leq\infty$, $T>0$ and $n\in\mathbb{N}$, be a solution of the linear transport equation
	\begin{equation}\label{basic transport}
		\partial_t f +b\cdot\nabla_x f +cf=h,
	\end{equation}
	where
	\begin{equation*}
		\begin{aligned}
			b & \in L^\gamma\left([0,T];W^{1,\alpha}_\mathrm{loc}\left(\mathbb{R}^n\right)\right), \\
			c & \in L^\gamma\left([0,T];L^{\alpha}_\mathrm{loc}\left(\mathbb{R}^n\right)\right), \\
			h & \in L^\gamma\left([0,T];L^{\beta}_\mathrm{loc}\left(\mathbb{R}^n\right)\right),
		\end{aligned}
	\end{equation*}
	for some $p'\leq \alpha < \infty$, $\frac 1p+\frac 1{p'}=1$, $1\leq\gamma<\infty$ and $1\leq\beta<\infty$ such that $\frac 1\beta=\frac 1\alpha + \frac 1p$.
	
	Then, for any $\chi_\delta(x)=\frac 1{\delta^n}\chi\left(\frac x\delta\right)$, with $\chi\in C_c^\infty\left(\mathbb{R}^n\right)$, $\chi\geq 0$, $\int_{\mathbb{R}^n}\chi(x)dx=1$ and $\delta>0$, the mollification $f_\delta=f * \chi_\delta$ satisfies
	\begin{equation*}
		\partial_t f_\delta +b\cdot\nabla_x f_\delta +cf_\delta=h+r_\delta,
	\end{equation*}
	where the remainder $r_\delta$ vanishes in $L^\gamma\left([0,T];L^{\beta}_\mathrm{loc}\left(\mathbb{R}^n\right)\right)$, as $\delta\rightarrow 0$.
	
	In particular, it follows that, for any renormalization $\beta\in C^1(\mathbb{R})$ such that $\beta'$ is bounded on $\mathbb{R}$,
	\begin{equation*}
		\partial_t \beta(f) +b\cdot\nabla_x \beta(f) +cf\beta'(f)=h\beta'(f).
	\end{equation*}
\end{lem}

The above lemma has fundamental consequences in transport theory and in the theory of ordinary differential equations. Indeed, as established by DiPerna and Lions in \cite{diperna6}, it can be shown that, loosely speaking, as soon as Lemma \ref{renorm transport} applies, weak solutions of \eqref{basic transport} are, in fact, renormalized solutions, unique and time continuous in the strong topology, and that the transport equation \eqref{basic transport} propagates strong compactness. In turn, the properties of the transport equation have important consequences on ordinary differential equations, and the existence and uniqueness of a Lagrangian flow was also established in \cite{diperna6} under very weak assumptions on the corresponding Eulerian flow, which should typically be in $L^1\left([0,T];W^{1,1}_\mathrm{loc}\left(\mathbb{R}^n\right)\right)$.

\bigskip

Thus, in view of the regularity hypothesis \eqref{force field hypotheses} on the force field
\begin{equation*}
	F\in L^1_{\mathrm{loc}}\left(dt;W^{1,1}_{\mathrm{loc}}\left(dxdv\right)\right),
\end{equation*}
applying Lemma \ref{renorm transport} to the transport equation \eqref{boltzmann renorm limit} (transport by the vector field $(v,F(t,x,v))\in \mathbb{R}^6$) yields that $\beta_\delta$ is a renormalized solutions of \eqref{boltzmann renorm limit}, that is to say, for any admissible renormalization $\beta$,
\begin{equation}\label{rerenorm}
	\begin{aligned}
		\partial_t\beta\left(\beta_\delta\right)+v\cdot\nabla_x\beta\left(\beta_\delta\right)
		& + F\cdot\nabla_v\beta\left(\beta_\delta\right) \\
		& =
		\beta'\left(\beta_\delta\right)h_\delta \mathcal{Q}^+(G,G)
		- \beta'\left(\beta_\delta\right)g_\delta f*_v\left[\int_{\mathbb{S}^2}b(\cdot,\sigma)d\sigma\right].
	\end{aligned}
\end{equation}

Finally, we let $\delta\rightarrow 0$ in the above renormalized equation. To this end, notice that $\beta'\left(\beta_\delta\right)g_\delta$ is bounded uniformly by a constant pointwise and converges almost everywhere to $\beta'(G)G$. Therefore, the last term above converges towards the expected renormalized loss term while it remains uniformly locally integrable. Moreover, we see that, integrating \eqref{rerenorm} locally in all variables, the gain term $\beta'\left(\beta_\delta\right)h_\delta \mathcal{Q}^+(G,G)$ remains uniformly locally integrable, so that it converges towards the expected renormalized gain term. On the whole, since the left-hand side of \eqref{rerenorm} is easily handled by the strong convergence of $\beta_\delta$ towards $G$, we conclude, letting $\delta\rightarrow 0$ in \eqref{rerenorm}, that $G$ solves \eqref{boltzmann renorm} in the sense of distributions, which completes the justification of Theorem \ref{classical renorm sol} according to \cite{diperna, diperna3, lions, lions2, lions3}.

It is to be emphasized that Theorem \ref{classical renorm sol} can be easily generalized to a system of Vlasov-Boltzmann equations for two species of particles.

\subsection{Coupling the Boltzmann equation with Maxwell's equations}

Thus, we see that the validity of Theorem \ref{classical renorm sol} rests crucially upon Lemma \ref{renorm transport} and, so, that the regularity hypothesis on the force field $F\in
L^1_{\mathrm{loc}}\left(dt;W^{1,1}_{\mathrm{loc}}\left(dxdv\right)\right)$ cannot be weakened, at least not with this method of proof. This is precisely the unique obstacle which prevents the construction of renormalized solutions for the Vlasov-Maxwell-Boltzmann systems \eqref{scaledVMB bis} and \eqref{scaled VMB two species bis}, whose force fields are not regular. 

As far as the existence theory of global solutions is concerned, notice that the nonlinear coupling of a kinetic equation with Maxwell's equations through the influence of a Lorentz force is not always a problem. In particular, it is possible to show the weak stability of the Vlasov-Maxwell system (without collisions) for densities in $L^\infty_tL^2_{x,v}$ and, therefore, to establish the existence of (non-renormalized) weak solutions for this system (see \cite{diperna2}). Indeed, neglecting the collision operators in \eqref{scaledVMB bis} and \eqref{scaled VMB two species bis}, the only remaining nonlinear terms are
\begin{equation*}
	\left( E + v \wedge B \right) \cdot \nabla_v f
	\qquad\text{and}\qquad
	\left( E + v \wedge B \right) \cdot \nabla_v f^\pm.
\end{equation*}
Since the densities $f$ in \eqref{scaledVMB bis} and $f^\pm$ in \eqref{scaled VMB two species bis} do enjoy some strong compactness (even some kind of regularity) in time and space by virtue of velocity averaging lemmas (see \cite{diperna5}, for instance), while the Lorentz force $E+v\wedge B$ is smooth in velocity (obviously, $E$ and $B$ do not depend on $v$), it is clear that the above nonlinear electromagnetic forcing terms are weakly stable as long as no renormalization is required. In conclusion, problematic difficulties arise when entering the realm of collisional kinetic theory, where renormalizing becomes a necessity. Nevertheless, it is to be noted that the existence of renormalized solutions for such collisionless Vlasov-Maxwell systems remains unknown, as well.

In contrast with the Vlasov-Maxwell-Boltzmann systems, the Vlasov-Poisson-Boltzmann systems \eqref{scaledVMP} and \eqref{scaledVMP two species} do enjoy the existence of renormalized solutions (see \cite{lions3}). Indeed, thanks to Poisson's equation, the force fields therein have enough regularity to apply Lemma \ref{renorm transport} and the strategy of proof of Theorem \ref{classical renorm sol} applies.

Of course, since then, there have been generalizations of Lemma \ref{renorm transport} and incidentally of the results from \cite{diperna6}, most notably by Ambrosio \cite{ambrosio}, where the local Sobolev regularity of the vector field was relaxed to a local $BV$ regularity, and by Le Bris and Lions \cite{lebris}, where a specific structure of the vector field, which unfortunately doesn't match the structure of \eqref{boltzmann renorm limit}, was used in order to impose a mere partial $W^{1,1}_\mathrm{loc}$ regularity on it. In any case, it is apparent, much like in the Cauchy-Lipschitz theorem on ordinary differential equations, that a minimum of a control on one full derivative of the vector field is necessary to crank the proof of Theorem \ref{classical renorm sol}, which is far from reach in the case of Vlasov-Maxwell-Boltzmann systems where $E,B\in L^\infty\left(dt;L^2(dx)\right)$, at best. This viewpoint is also corroborated by the counterexamples presented at the end of \cite{diperna6}.

Thus, it seems that any result confirming the existence of renormalized solutions for Vlasov-Maxwell-Boltzmann systems will have to exploit the very specific structure of the electromagnetic interaction within the plasma.

Surprisingly the situation is much better when the microscopic interactions described by the collision operator have infinite range so that the collisional cross-section has a singularity at grazing collisions~: the entropy dissipation indeed controls some derivative with respect to $v$ in this case. Using the hypoellipticity of the kinetic transport operator, we can then transfer part of this regularity onto the $x$ variable. Following the strategy by Alexandre and Villani \cite{alexandre}, and renormalizing the Vlasov-Boltzmann equation by concave functions, we thus get some global renormalized solutions involving a defect measure (which is formally $0$ because of the conservation of mass). This construction has been sketched in \cite{arsenio-sr-CRAS}. It will be detailed and used to obtain fully rigorous convergence results in Part \ref{part 3}.

An alternative approach based on Young measures, as introduced by Lions in \cite{lions3} will be the focus of our work in Part \ref{part 4}. We will see that, even though the notion of solution is very poor, the asymptotic analysis is robust and leads to similar convergence results.

Note that Parts \ref{part 3} and \ref{part 4} will be more technical as we will have to deal with very weak solutions. However the strategy of proof as well as the main arguments will be similar to the ones presented here, this is why we start with conditional results.

\subsection{The setting of our conditional study}%\label{conditional setting}

We provide now a precise definition of renormalized solutions for the Vlasov-Maxwell-Boltzmann systems \eqref{scaledVMB bis} and \eqref{scaled VMB two species bis}, even though their existence remains uncertain.

\begin{defi}
	We say that a density function $f(t,x,v)=MG(t,x,v)\geq0$ and electromagnetic vector fields $E(t,x)$ and $B(t,x)$, where $(t,x,v)\in [0,\infty)\times\mathbb{R}^3\times\mathbb{R}^3$, such that
	\begin{equation*}% \label{membership}
		\begin{aligned}
			G & \in C\left([0,\infty);\textit{w-}L^1_{\mathrm{loc}}\left(dxdv\right)\right)
			\cap L^\infty\left([0,\infty),dt;L^1_{\mathrm{loc}}\left(dx;L^1\left((1+|v|^2)M dv\right)\right)\right),\\
			E,B & \in
			C\left([0,\infty);\textit{w-}L^2\left(dx\right)\right)
			\cap L^\infty\left([0,\infty),dt;L^2\left(dx\right)\right),
		\end{aligned}
	\end{equation*}
	are a \textbf{renormalized solution of the one species Vlasov-Maxwell-Boltzmann system \eqref{scaledVMB bis}} if they solve
	\begin{equation*}
		\begin{cases}
		\begin{aligned}
			\d_t \beta\left(G\right) + v \cdot \nabla_x \beta\left(G\right)
			+ \left( E + v \wedge B \right) \cdot \nabla_v \beta\left(G\right)
			& - E \cdot v \beta'\left(G\right)G\\
			&= \beta'\left(G\right)\mathcal{Q}(G,G), \\
			\d_t E - \ROT B &= - \int_{\mathbb{R}^3} MGv dv,\\
			\d_t B + \ROT E& = 0,\\
			\DIV E &=\int_{\mathbb{R}^3} MG dv -1,\\
			\DIV B &=0,
		\end{aligned}
		\end{cases}
	\end{equation*}
	in the sense of distributions for any admissible renormalization, and satisfies the entropy inequality, for all $t>0$,
	\begin{equation}\label{renormalized solution 1 entropy}
		\begin{aligned}
			H\left(f \right)
			+ \frac 1{2} \int_{\mathbb{R}^3} \left(|E|^2+ |B|^2\right) dx
			& + \int_0^t\int_{\mathbb{R}^3}D(f)(s) dx ds
			\\
			& \leq H\left(f^{\mathrm{in}}\right)
			+ \frac1{2}\int_{\mathbb{R}^3} \left(|E^{\rm in}|^2+ |B^{\rm in}|^2\right) dx<\infty,
		\end{aligned}
	\end{equation}
	where $f^\mathrm{in}=MG^\mathrm{in}$ is the initial value of $f=MG$ and the relative entropy $H(f)=H(f|M)$ is defined in \eqref{def H}, while the entropy dissipation $D(f)$ is defined in \eqref{def D}.
\end{defi}
	
\begin{defi}
	We say that density functions $G^+(t,x,v)\geq0$ and $G^-(t,x,v)\geq0$, and electromagnetic vector fields $E(t,x)$ and $B(t,x)$, where $(t,x,v)\in [0,\infty)\times\mathbb{R}^3\times\mathbb{R}^3$, such that
	\begin{equation*}% \label{membership}
		\begin{aligned}
			G^\pm & \in C\left([0,\infty);\textit{w-}L^1_{\mathrm{loc}}\left(dxdv\right)\right)
			\cap L^\infty\left([0,\infty),dt;L^1_{\mathrm{loc}}\left(dx;L^1\left((1+|v|^2)M dv\right)\right)\right),\\
			E,B & \in
			C\left([0,\infty);\textit{w-}L^2\left(dx\right)\right)
			\cap L^\infty\left([0,\infty),dt;L^2\left(dx\right)\right),
		\end{aligned}
	\end{equation*}
	are a \textbf{renormalized solution of the two species Vlasov-Maxwell-Boltzmann system \eqref{scaled VMB two species bis}} if they solve
	\begin{equation*}
		\begin{cases}
		\begin{aligned}
			\d_t \beta\left(G^\pm\right) + v \cdot \nabla_x \beta\left(G^\pm\right)
			& \pm \left( E + v \wedge B \right) \cdot \nabla_v \beta\left(G^\pm\right)
			\mp E \cdot v \beta'\left(G^\pm\right)G^\pm\\
			&= \beta'\left(G^\pm\right)\mathcal{Q}(G^\pm,G^\pm) + \beta'\left(G^\pm\right)\mathcal{Q}(G^\pm,G^\mp), \\
			\d_t E - \ROT B &= - \int_{\mathbb{R}^3} M\left(G^+-G^-\right)v dv,\\
			\d_t B + \ROT E& = 0,\\
			\DIV E &=\int_{\mathbb{R}^3} M\left(G^+-G^-\right) dv,\\
			\DIV B &=0,
		\end{aligned}
		\end{cases}
	\end{equation*}
	in the sense of distributions for any admissible renormalization, and satisfy the entropy inequality, for all $t>0$,
	\begin{equation}\label{renormalized solution 2 entropy}
		\begin{aligned}
			H\left(f^{+}\right)
			+ H\left(f^{-}\right)
			& + \frac 1{2} \int_{\mathbb{R}^3} \left(|E|^2+ |B|^2\right) dx \\
			& + \int_0^t\int_{\mathbb{R}^3}\left(D\left(f^+\right)+D\left(f^-\right)
			+ D\left(f^+,f^-\right)\right)(s) dx ds
			\\
			& \leq
			H\left(f^{+\mathrm{in}}\right)
			+ H\left(f^{-\mathrm{in}}\right)
			+ \frac1{2}\int_{\mathbb{R}^3} \left(|E^{\rm in}|^2+ |B^{\rm in}|^2\right) dx<\infty,
		\end{aligned}
	\end{equation}
	where $f^{\pm\mathrm{in}}=MG^{\pm\mathrm{in}}$ is the initial value of $f^\pm=MG^\pm$ and the relative entropies $H\left(f^\pm\right)=H\left(f^\pm|M\right)$ are defined in \eqref{def H}, while the entropy dissipations $D\left(f^\pm\right)$ and $D\left(f^+,f^-\right)$ are defined in \eqref{def D} and \eqref{def D mixed}.
\end{defi}

\subsection{Macroscopic conservation laws}\label{macroscopic defects}

As already explained in Section \ref{macro properties}, the one species Vlasov-Maxwell-Boltzmann system \eqref{scaledVMB bis} formally satisfies the macroscopic conservation laws
\begin{equation*}
	\partial_t\int_{\mathbb{R}^3}f
	\begin{pmatrix} 1 \\ v \\ \frac{|v|^2}{2} \end{pmatrix}
	dv
	+
	\nabla_x\cdot\int_{\mathbb{R}^3}f
	\begin{pmatrix} v \\ v\otimes v \\ \frac{|v|^2}{2}v \end{pmatrix}
	dv
	=
	\int_{\mathbb{R}^3}f
	\begin{pmatrix} 0 \\ E+v\wedge B \\ E\cdot v \end{pmatrix}
	dv,
\end{equation*}
while the two species Vlasov-Maxwell-Boltzmann system \eqref{scaled VMB two species bis} formally satisfies the macroscopic conservation laws
\begin{equation*}
	\partial_t\int_{\mathbb{R}^3}f^\pm
	dv
	+
	\nabla_x\cdot\int_{\mathbb{R}^3}f^\pm v
	dv
	=0,
\end{equation*}
and
\begin{equation*}
	\begin{aligned}
		\partial_t\int_{\mathbb{R}^3}\left(f^++f^-\right)
		\begin{pmatrix} v \\ \frac{|v|^2}{2} \end{pmatrix}
		dv +
		\nabla_x & \cdot \int_{\mathbb{R}^3}\left(f^++f^-\right)
		\begin{pmatrix} v\otimes v \\ \frac{|v|^2}{2}v \end{pmatrix}
		dv \\
		& =	
		\int_{\mathbb{R}^3}\left(f^+-f^-\right)
		\begin{pmatrix} E+v\wedge B \\ E\cdot v \end{pmatrix}
		dv.
	\end{aligned}
\end{equation*}
However, it is at first unclear whether such formal laws are actually rigorously satisfied by the renormalized solutions defined in the previous section. It is therefore necessary to justify their validity.

To this end, we suppose now that such renormalized solutions $\left(f,E,B\right)$ and $\left(f^\pm,E,B\right)$ have been previously obtained through an approximation procedure as detailed in Section \ref{vlasov}. More precisely, we assume that there are sequences $\left\{\left(f_k,E_k,B_k\right)\right\}_{k\in\mathbb{N}}$ and $\left\{\left(f^\pm_k,E_k,B_k\right)\right\}_{k\in\mathbb{N}}$ of smooth solutions to \eqref{scaledVMB bis} and \eqref{scaled VMB two species bis} (or appropriate approximations of these systems), for some uniformly bounded initial data, satisfying all macroscopic conservation laws and respectively converging in some suitable weak sense towards $\left(f,E,B\right)$ and $\left(f^\pm,E,B\right)$. Therefore, in virtue of the uniform bounds provided by the entropy inequalities \eqref{renormalized solution 1 entropy} and \eqref{renormalized solution 2 entropy}, it is readily seen that the terms $f_k$, $f_kv$, $f_k^\pm$ and $f^\pm_k v$ are all respectively converging to $f$, $fv$, $f^\pm$ and $f^\pm v$ weakly in $L^1_\mathrm{loc}\left(dtdx;L^1\left(dv\right)\right)$.

It follows that the conservations of mass
\begin{equation*}
	\begin{aligned}
		\partial_t\int_{\mathbb{R}^3}f
		dv
		+
		\nabla_x\cdot\int_{\mathbb{R}^3}f v
		dv
		& = 0,
		\\
		\partial_t\int_{\mathbb{R}^3}f^\pm
		dv
		+
		\nabla_x\cdot\int_{\mathbb{R}^3}f^\pm v
		dv
		& = 0,
	\end{aligned}
\end{equation*}
are easily established for renormalized solutions. This is in general the case when dealing with collisional kinetic equations.

However, the same is unfortunately not true for the conservations of momentum and energy. Indeed, these laws involve higher moments of $f$ and $f^\pm$, which may be singular due to high velocities build-up, as well as products of electromagnetic fields with particle densities, which may not even make sense if not renormalized.

In order to account for large velocities, we introduce now, following \cite{lions4}, since the terms $f_k|v|^2$, $f|v|^2$, $f_k^\pm|v|^2$ and $f^\pm|v|^2$ are uniformly bounded in $L^1_\mathrm{loc}\left(dtdx;L^1\left(dv\right)\right)$, the Radon measures $m_{ij}\in \mathcal{M}_\mathrm{loc}\left([0,\infty)\times\mathbb{R}^3\right)$, $i,j=1,2,3$, defined as the following defects in the limit $k\to\infty$ (up to extraction of subsequences)~:
\begin{equation*}
	\int_{\mathbb{R}^3}f_kv_iv_jdv
	\stackrel{*}{\rightharpoonup}
	\int_{\mathbb{R}^3}fv_iv_jdv + m_{ij}
	\quad \text{in }
	\mathcal{M}_\mathrm{loc}\left([0,\infty)\times\mathbb{R}^3\right),
\end{equation*}
in the one species case, and
\begin{equation*}
	\int_{\mathbb{R}^3}\left(f_k^++f_k^-\right)v_iv_jdv
	\stackrel{*}{\rightharpoonup}
	\int_{\mathbb{R}^3}\left(f^++f^-\right)v_iv_jdv + m_{ij}
	\quad \text{in }
	\mathcal{M}_\mathrm{loc}\left([0,\infty)\times\mathbb{R}^3\right),
\end{equation*}
in the two species case. Note that the measures $m_{ij}$ are also defined by the limits, valid for any $R>0$,
\begin{equation*}
	\begin{aligned}
		\int_{\mathbb{R}^3}f_kv_iv_j\mathds{1}_{\left\{|v|\geq R\right\}}dv
		& \stackrel{*}{\rightharpoonup}
		\int_{\mathbb{R}^3}fv_iv_j\mathds{1}_{\left\{|v|\geq R\right\}}dv + m_{ij},
		\\
		\int_{\mathbb{R}^3}\left(f_k^++f_k^-\right)v_iv_j\mathds{1}_{\left\{|v|\geq R\right\}}dv
		& \stackrel{*}{\rightharpoonup}
		\int_{\mathbb{R}^3}\left(f^++f^-\right)v_iv_j\mathds{1}_{\left\{|v|\geq R\right\}}dv + m_{ij},
	\end{aligned}
\end{equation*}
in $\mathcal{M}_\mathrm{loc}\left([0,\infty)\times\mathbb{R}^3\right)$. In particular, it follows that the matrix measure $m=\left(m_{ij}\right)_{1\leq i,j\leq 3}$ is symmetric and positive definite in the sense that, for any $\varphi\in C_c\left([0,\infty)\times\mathbb{R}^3;\mathbb{R}^3\right)$,
\begin{equation*}
	\int_{[0,\infty)\times\mathbb{R}^3} \varphi^t (dm) \varphi=\sum_{i,j=1}^3\int_{[0,\infty)\times\mathbb{R}^3} \varphi_i\varphi_j dm_{ij}\geq 0,
\end{equation*}
whence, for any $\varphi,\psi\in C_c\left([0,\infty)\times\mathbb{R}^3;\mathbb{R}^3\right)$,
\begin{equation*}
	\left|\int_{[0,\infty)\times\mathbb{R}^3} \varphi^t (dm) \psi\right|
	\leq
	\frac 12
	\int_{[0,\infty)\times\mathbb{R}^3} \varphi^t (dm) \varphi
	+
	\frac 12
	\int_{[0,\infty)\times\mathbb{R}^3} \psi^t (dm) \psi.
\end{equation*}
Further substituting $\varphi$ and $\psi$ in the preceding inequality by $\lambda^\frac 12\varphi$ and $\lambda^{-\frac 12}\psi$, with $\lambda>0$, respectively, and then optimizing the resulting inequality in $\lambda>0$ yields that
\begin{equation}\label{measureCS1}
	\left|\int_{[0,\infty)\times\mathbb{R}^3} \varphi^t (dm) \psi\right|
	\leq
	\left(
	\int_{[0,\infty)\times\mathbb{R}^3} \varphi^t (dm) \varphi
	\right)^\frac 12
	\left(
	\int_{[0,\infty)\times\mathbb{R}^3} \psi^t (dm) \psi
	\right)^\frac 12,
\end{equation}
for any $\varphi,\psi\in C_c\left([0,\infty)\times\mathbb{R}^3;\mathbb{R}^3\right)$.

The matrix measure $m$ will be used to characterize the flux terms in the conservation of momentum and the density terms in the conservation of energy. However, the flux terms in the conservation of energy contain higher order moments which cannot be handled and we will therefore simply leave them out of the analysis by only considering the global conservation of energy.

As for the forcing terms involving the electromagnetic fields, they do not even make sense with the sole use of the a priori estimates provided by the entropy inequalities \eqref{renormalized solution 1 entropy} and \eqref{renormalized solution 2 entropy}. It is therefore necessary to use now the conservation laws of energy \eqref{maxwell energy}-\eqref{maxwell energy two species} and for the Poynting vector \eqref{maxwell poynting}-\eqref{maxwell poynting two species} in Maxwell's equations to recast these forcing terms with quadratic expressions involving the electromagnetic fields only.

Thus, using the identities \eqref{maxwell energy}, \eqref{maxwell energy two species}, \eqref{maxwell poynting} and \eqref{maxwell poynting two species}, the local conservations of momentum may be rewritten as
\begin{equation}\label{k conservation 1}
	\begin{aligned}
		\partial_t\left(\int_{\mathbb{R}^3}f_k
		v
		dv+E_k\wedge B_k\right)
		+
		\nabla_x\cdot\left(\int_{\mathbb{R}^3}f_k
		v\otimes v
		dv
		-E_k\otimes E_k-B_k\otimes B_k\right)
		&
		\\
		+\nabla_x\left(\frac{|E_k|^2+|B_k|^2}{2}\right)
		& =
		E_k,
	\end{aligned}
\end{equation}
in the one species case, and as
\begin{equation}\label{k conservation 2}
	\begin{aligned}
		\partial_t & \left(\int_{\mathbb{R}^3}\left(f_k^++f_k^-\right)
		v
		dv +E_k\wedge B_k \right)
		\\
		& +
		\nabla_x \cdot \left(
		\int_{\mathbb{R}^3}\left(f_k^++f_k^-\right)
		v\otimes v
		dv
		-E_k\otimes E_k-B_k\otimes B_k
		\right)
		\\
		& +\nabla_x\left(\frac{|E_k|^2+|B_k|^2}{2}\right)
		= 0,
	\end{aligned}
\end{equation}
in the two species case, whereas the global conservations of energy may be expessed as
\begin{equation}\label{k conservation 3}
	\frac{d}{dt}\int_{\mathbb{R}^3}\left(\int_{\mathbb{R}^3}f_k
	\frac{|v|^2}{2}
	dv+\frac{|E_k|^2+|B_k|^2}{2}\right)dx
	=0,
\end{equation}
in the one species case, and as
\begin{equation}\label{k conservation 4}
	\frac{d}{dt}\int_{\mathbb{R}^3}\left(\int_{\mathbb{R}^3}\left(f_k^++f_k^-\right)
	\frac{|v|^2}{2}
	dv+\frac{|E_k|^2+|B_k|^2}{2}\right)dx
	=0,
\end{equation}
in the two species case.

Passing to the limit $k\to \infty$ therefore requires the introduction of yet another set of Radon measures $a_{ij}\in \mathcal{M}_\mathrm{loc}\left([0,\infty)\times\mathbb{R}^3\right)$, $i,j=1,2,3,4,5,6$, where the matrix measure $a=\left(a_{ij}\right)_{1\leq i,j\leq 6}$ is defined as the following defect~:
\begin{equation*}
	\begin{pmatrix}
		E_k\\B_k
	\end{pmatrix}
	\otimes
	\begin{pmatrix}
		E_k\\B_k
	\end{pmatrix}
	\stackrel{*}{\rightharpoonup}
	\begin{pmatrix}
		E\\B
	\end{pmatrix}
	\otimes
	\begin{pmatrix}
		E\\B
	\end{pmatrix} + a
	\quad \text{in }
	\mathcal{M}_\mathrm{loc}\left([0,\infty)\times\mathbb{R}^3\right).
\end{equation*}
Note that the matrix measure $n$ is also defined by the limit
\begin{equation*}
	\begin{pmatrix}
		E_k-E\\B_k-B
	\end{pmatrix}
	\otimes
	\begin{pmatrix}
		E_k-E\\B_k-B
	\end{pmatrix}
	\stackrel{*}{\rightharpoonup}
	a
	\quad \text{in }
	\mathcal{M}_\mathrm{loc}\left([0,\infty)\times\mathbb{R}^3\right).
\end{equation*}
It then follows that, as before, the matrix measure $a$ is symmetric and positive definite in the sense that, for any $\varphi\in C_c\left([0,\infty)\times\mathbb{R}^3;\mathbb{R}^6\right)$,
\begin{equation*}
	\int_{[0,\infty)\times\mathbb{R}^3} \varphi^t (da) \varphi=\sum_{i,j=1}^6\int_{[0,\infty)\times\mathbb{R}^3} \varphi_i\varphi_j da_{ij}\geq 0,
\end{equation*}
whence (see the analogous inequality \eqref{measureCS1})
\begin{equation}\label{measureCS2}
	\left|\int_{[0,\infty)\times\mathbb{R}^3} \varphi^t (da) \psi\right|
	\leq
	\left(
	\int_{[0,\infty)\times\mathbb{R}^3} \varphi^t (da) \varphi
	\right)^\frac 12
	\left(
	\int_{[0,\infty)\times\mathbb{R}^3} \psi^t (da) \psi
	\right)^\frac 12,
\end{equation}
for any $\varphi,\psi\in C_c\left([0,\infty)\times\mathbb{R}^3;\mathbb{R}^6\right)$.

For mere convenience of notation, we further introduce the matrix measures $e=\left(a_{ij}\right)_{1\leq i,j\leq 3}$ and $b=\left(a_{(i+3)(j+3)}\right)_{1\leq i,j\leq 3}$.

Now, letting $k\to\infty$ in the conservation laws \eqref{k conservation 1}, \eqref{k conservation 2}, \eqref{k conservation 3} and \eqref{k conservation 4} respectively yields the local conservations of momentum
\begin{equation*}
	\begin{aligned}
		\partial_t & \left(\int_{\mathbb{R}^3}f
		v
		dv+E\wedge B +
		\begin{pmatrix}
			a_{26}-a_{35}\\a_{34}-a_{16}\\a_{15}-a_{24}
		\end{pmatrix}
		\right)
		\\
		& +
		\nabla_x\cdot\left(\int_{\mathbb{R}^3}f
		v\otimes v
		dv + m
		-E\otimes E-e-B\otimes B-b\right)
		\\
		& +\nabla_x\left(\frac{|E|^2+|B|^2+\operatorname{Tr}a}{2}\right)
		=
		E,
	\end{aligned}
\end{equation*}
in the one species case, and
\begin{equation*}
	\begin{aligned}
		\partial_t & \left(\int_{\mathbb{R}^3}\left(f^++f^-\right)
		v
		dv +E\wedge B
		+
		\begin{pmatrix}
			a_{26}-a_{35}\\a_{34}-a_{16}\\a_{15}-a_{24}
		\end{pmatrix}
		\right)
		\\
		& +
		\nabla_x \cdot \left(
		\int_{\mathbb{R}^3}\left(f^++f^-\right)
		v\otimes v
		dv + m
		-E\otimes E - e - B\otimes B - b
		\right)
		\\
		& +\nabla_x\left(\frac{|E|^2+|B|^2+\operatorname{Tr}a}{2}\right)
		= 0,
	\end{aligned}
\end{equation*}
in the two species case, as well as the global energy decay
\begin{equation*}
	\begin{aligned}
		\int_{\mathbb{R}^3} & \left(\int_{\mathbb{R}^3}f
		\frac{|v|^2}{2}
		dv+\frac{\operatorname{Tr}m}{2}+\frac{|E|^2+|B|^2+\operatorname{Tr}a}{2}\right)dx
		\\
		& \leq
		\int_{\mathbb{R}^3}\left(\int_{\mathbb{R}^3}f^\mathrm{in}
		\frac{|v|^2}{2}
		dv+\frac{|E^\mathrm{in}|^2+|B^\mathrm{in}|^2}{2}\right)dx,
	\end{aligned}
\end{equation*}
in the one species case, and
\begin{equation*}
	\begin{aligned}
		\int_{\mathbb{R}^3} & \left(\int_{\mathbb{R}^3}\left(f^++f^-\right)
		\frac{|v|^2}{2}
		dv +\frac{\operatorname{Tr}m}{2}+\frac{|E|^2+|B|^2+\operatorname{Tr}a}{2}\right)dx
		\\
		& \leq
		\int_{\mathbb{R}^3}\left(\int_{\mathbb{R}^3}\left(f^{+\mathrm{in}}+f^{-\mathrm{in}}\right)
		\frac{|v|^2}{2}
		dv+\frac{|E^\mathrm{in}|^2+|B^\mathrm{in}|^2}{2}\right)dx,
	\end{aligned}
\end{equation*}
in the two species case.

Note that the above global energy decay containing the defect measures may be incorporated into the entropy inequalities \eqref{renormalized solution 1 entropy} and \eqref{renormalized solution 2 entropy}, so that renormalized solutions of the Vlasov-Maxwell-Boltzmann systems \eqref{scaledVMB bis} and \eqref{scaled VMB two species bis} may be assumed to respectively satisfy the entropy inequalities
\begin{equation}\label{renormalized solution 1 entropy defect}
	\begin{aligned}
		H\left(f \right)
		+ \frac 1{2} \int_{\mathbb{R}^3} \left(|E|^2+ |B|^2\right) dx
		& + \frac 12 \int_{\mathbb{R}^3}\operatorname{Tr}\left(m+e+b\right)dx
		+ \int_0^t\int_{\mathbb{R}^3}D(f)(s) dx ds
		\\
		& \leq H\left(f^{\mathrm{in}}\right)
		+ \frac1{2}\int_{\mathbb{R}^3} \left(|E^{\rm in}|^2+ |B^{\rm in}|^2\right) dx<\infty,
	\end{aligned}
\end{equation}
in the one species case, and
\begin{equation}\label{renormalized solution 2 entropy defect}
	\begin{aligned}
		H\left(f^{+}\right)
		+ H\left(f^{-}\right)
		& + \frac 1{2} \int_{\mathbb{R}^3} \left(|E|^2+ |B|^2\right) dx
		+ \frac 12 \int_{\mathbb{R}^3}\operatorname{Tr}\left(m+e+b\right)dx
		\\
		& + \int_0^t\int_{\mathbb{R}^3}\left(D\left(f^+\right)+D\left(f^-\right)
		+ D\left(f^+,f^-\right)\right)(s) dx ds
		\\
		& \leq
		H\left(f^{+\mathrm{in}}\right)
		+ H\left(f^{-\mathrm{in}}\right)
		+ \frac1{2}\int_{\mathbb{R}^3} \left(|E^{\rm in}|^2+ |B^{\rm in}|^2\right) dx<\infty,
	\end{aligned}
\end{equation}
in the two species case.

\bigskip

The preceding characterization of defects in macroscopic conservation laws will not be of further use in our study of the hydrodynamic limits of the one species Vlasov-Maxwell-Boltzmann system \eqref{scaledVMB bis}. It will, however, be of crucial utility in the renormalized relative entropy method developed later on in Chapter \ref{entropy method} in relation with hydrodynamic limits of the two species Vlasov-Maxwell-Boltzmann system \eqref{scaled VMB two species bis}.

\bigskip

Notice, finally, that the symmetry and the positive definiteness of the matrix measures $m$ and $a$ imply that the bounds on the non-negative measures $\operatorname{Tr}\left(m+e+b\right)$, provided a priori by the entropy inequalities \eqref{renormalized solution 1 entropy defect} and \eqref{renormalized solution 2 entropy defect}, are sufficient to control all components of $m$ and $a$. Indeed, the inequalities \eqref{measureCS1} and \eqref{measureCS2} provide all necessary estimates of $m$ and $a$ in terms of $\operatorname{Tr}\left(m+e+b\right)$.

In particular, it is readily seen that, for any $u \in C_c\left(\mathbb{R}^3;\mathbb{R}^3\right)$,
\begin{equation*}
	\left|\int_{\mathbb{R}^3}
	u\cdot
	\begin{pmatrix}
		a_{26}-a_{35}\\a_{34}-a_{16}\\a_{15}-a_{24}
	\end{pmatrix}
	dx\right|
	\leq 6
	\left\|u\right\|_{L^\infty}
	\left(\int_{\mathbb{R}^3}\left(\operatorname{Tr}e\right) dx\right)^\frac 12
	\left(\int_{\mathbb{R}^3}\left(\operatorname{Tr}b\right) dx\right)^\frac 12.
\end{equation*}
Still, since the above defect stems from a vector product, it is possible to improve the constant in the preceding inequality. We record such improvement in the following result, for later use.

\begin{lem}\label{vector defect}
	For any $u \in C_c\left(\mathbb{R}^3;\mathbb{R}^3\right)$, it holds that
	\begin{equation*}
		\left|\int_{\mathbb{R}^3}
		u\cdot
		\begin{pmatrix}
			a_{26}-a_{35}\\a_{34}-a_{16}\\a_{15}-a_{24}
		\end{pmatrix}
		dx\right|
		\leq
		\left\|u\right\|_{L^\infty}
		\left(\int_{\mathbb{R}^3}\left(\operatorname{Tr}e\right) dx\right)^\frac 12
		\left(\int_{\mathbb{R}^3}\left(\operatorname{Tr}b\right) dx\right)^\frac 12.
	\end{equation*}
\end{lem}

\begin{proof}
	Using \eqref{measureCS2}, we first obtain
	\begin{equation*}
		\begin{aligned}
			\left|\int_{\mathbb{R}^3}
			u\cdot
			\begin{pmatrix}
				a_{26}-a_{35}\\a_{34}-a_{16}\\a_{15}-a_{24}
			\end{pmatrix}
			dx\right|
			& \leq
			\left|\int_{\mathbb{R}^3}
			\left(a_{15}u_3-a_{16}u_2\right)
			dx\right|
			+
			\left|\int_{\mathbb{R}^3}
			\left(a_{26}u_1-a_{24}u_3\right)
			dx\right|
			\\
			& +
			\left|\int_{\mathbb{R}^3}
			\left(a_{34}u_2-a_{35}u_1\right)
			dx\right|
			\\
			& \leq
			\left(\int_{\mathbb{R}^3}
			a_{11}
			dx\right)^\frac 12
			\left(\int_{\mathbb{R}^3}
			\left(a_{55}u_3^2-2a_{56}u_3u_2+a_{66}u_{2}^2\right)
			dx\right)^\frac 12
			\\
			& +
			\left(\int_{\mathbb{R}^3}
			a_{22}
			dx\right)^\frac 12
			\left(\int_{\mathbb{R}^3}
			\left(a_{66}u_1^2-2a_{64}u_1u_3+a_{44}u_{3}^2\right)
			dx\right)^\frac 12
			\\
			& +
			\left(\int_{\mathbb{R}^3}
			a_{33}
			dx\right)^\frac 12
			\left(\int_{\mathbb{R}^3}
			\left(a_{44}u_2^2-2a_{45}u_2u_1+a_{55}u_{1}^2\right)
			dx\right)^\frac 12.
		\end{aligned}
	\end{equation*}
	It follows that, for any $\alpha>0$,
	\begin{equation*}
		\begin{aligned}
			\left|\int_{\mathbb{R}^3}
			u\cdot
			\begin{pmatrix}
				a_{26}-a_{35}\\a_{34}-a_{16}\\a_{15}-a_{24}
			\end{pmatrix}
			dx\right|
			& \leq \frac\alpha 2
			\int_{\mathbb{R}^3}
			\left(\operatorname{Tr}e\right)
			dx
			\\
			& +
			\frac{1}{2\alpha}\int_{\mathbb{R}^3}
			\left(a_{55}u_3^2-2a_{56}u_3u_2+a_{66}u_{2}^2\right)
			dx
			\\
			& +\frac 1{2\alpha}
			\int_{\mathbb{R}^3}
			\left(a_{66}u_1^2-2a_{64}u_1u_3+a_{44}u_{3}^2\right)
			dx
			\\
			& +\frac 1{2\alpha}
			\int_{\mathbb{R}^3}
			\left(a_{44}u_2^2-2a_{45}u_2u_1+a_{55}u_{1}^2\right)
			dx
			\\
			& = \frac\alpha 2
			\int_{\mathbb{R}^3}
			\left(\operatorname{Tr}e\right)
			dx
			+
			\frac 1{2\alpha}
			\int_{\mathbb{R}^3}
			\left(\operatorname{Tr}b\right)|u|^2
			dx
			\\
			& -
			\frac{1}{2\alpha}\int_{\mathbb{R}^3}
			\left(u^tbu\right)
			dx
			\\
			& \leq \frac\alpha 2
			\int_{\mathbb{R}^3}
			\left(\operatorname{Tr}e\right)
			dx
			+
			\frac {\|u\|_{L^\infty}^2}{2\alpha}
			\int_{\mathbb{R}^3}
			\left(\operatorname{Tr}b\right)
			dx,
		\end{aligned}
	\end{equation*}
	which, upon optimizing in $\alpha>0$, concludes the justification of the lemma.
\end{proof}

% =================
% = Main result I =
% =================

\section[The incompressible quasi-static Navier-Stokes-Fourier-\ldots]{The incompressible quasi-static Navier-Stokes-Fourier-Maxwell-Poisson system}
\label{main moments method}
% \section{The electrostatic Navier-Stokes-Maxwell system}

Henceforth, in this second part of our work on conditional results, unless otherwise stated, we will focus, for the mere sake of technical simplicity, on some Maxwellian cross-section, say $b\equiv 1$. All other mathematically and physically pertinent cross-sections (deriving from hard, soft, short-range and long-range interaction potentials) will be discussed and treated in full generality in the remaining parts of our work on unconditional results.

\bigskip

Following Section \ref{formal one}, we first consider a plasma constituted of a gas of cations (positively charged ions), with a uniform background of heavy anions (negatively charged ions) assumed to be at statistical equilibrium. Elementary interactions are taken into account by both a mean field term (corresponding to long-range interactions) and a local collision term (associated to short-range interactions) involving possibly different mean free paths. Thus, the charged particles evolve under the coupled effect of the Lorentz force due to the self-induced electromagnetic field, and of the collisions with other particles, according to the following scaled Vlasov-Maxwell-Boltzmann system~:
\begin{equation}\label{VMB1}
	\begin{cases}
		\begin{aligned}
			\epsilon\d_t f_\eps + v \cdot \nabla_x f_\eps + \eps \left( E_\eps + v \wedge B_\eps \right) \cdot \nabla_v f_\eps &=
			\frac{1}{\epsilon}Q(f_\eps,f_\eps),
			\\
			f_\eps & =M\left(1+\epsilon g_\eps\right),
			\\
			\eps\d_t E_\eps - \ROT B_\eps &= - \int_{\mathbb{R}^3} g_\eps v M dv,
			\\
			\eps \d_t B_\eps + \ROT E_\eps & = 0,
			\\
			\DIV E_\eps &=\int_{\mathbb{R}^3} g_\eps M dv,
			\\
			\DIV B_\eps &=0.
		\end{aligned}
	\end{cases}
\end{equation}
In this scaling, the entropy inequality states that
\begin{equation}\label{entropy1}
	\begin{aligned}
		\frac1{\eps^2} H\left(f_\eps \right)
		+ \frac 1{2} \int_{\mathbb{R}^3} \left(|E_\eps|^2+ |B_\eps|^2\right) dx
		& +\frac{1}{\epsilon^4}\int_0^t\int_{\mathbb{R}^3}D(f_\eps)(s) dx ds
		\\
		& \leq  \frac1{\eps^2}H\left(f_\eps^{\mathrm{in}}\right)
		+ \frac1{2}\int_{\mathbb{R}^3} \left(|E_\eps^{\rm in}|^2+ |B_\eps^{\rm in}|^2\right) dx,
	\end{aligned}
\end{equation}
where $H\left(f_\eps\right)=H\left(f_\eps|M\right)$. In particular, it yields uniform bounds on $E_\eps$, $B_\eps$ and $g_\eps$.

Since we are interested in the limiting fluctuation $g_\eps\rightharpoonup g$, it is then natural to rewrite  the  kinetic equation in terms of the fluctuations $g_\eps$,
\begin{equation}\label{VMB-fluct1}
	\eps\d_t g_\eps +v\cdot \nabla_x g_\eps
	+\eps( E_\eps+ v\wedge B_\eps) \cdot \nabla_v g_\eps- E_\eps \cdot v\left(1+\epsilon g_\epsilon\right)
	=-\frac1\eps \cL g_\eps +\cQ(g_\eps,g_\eps).
\end{equation}

According to the formal analysis from Section \ref{formal one}, we then expect the limiting macroscopic observables to solve the incompressible quasi-static Navier-Stokes-Fourier-Maxwell-Poisson system~:
\begin{equation}\label{NSFMP 2}
	\begin{cases}
		\begin{aligned}
			\d_t u
			+
			u\cdot\nabla_x u - \mu\Delta_x u
			& = -\nabla_x p + E
			+ \rho \nabla_x\theta + u \wedge B , \hspace{-20mm}&& \\
			&& \Div u & = 0,\\
			\d_t \left(\frac32\theta-\rho\right)
			+
			u\cdot\nabla_x\left(\frac32\theta-\rho\right)
			- \frac 52 \kappa \Delta_x\theta
			& = 0,
			& \Delta_x(\rho+\theta) & =\rho, \\
			\ROT B & = u, & \Div E & = \rho , \\
			\partial_t B + \rot E  & = 0, & \Div B & = 0.
		\end{aligned}
	\end{cases}
\end{equation}
As discussed in Section \ref{stability existence 1}, this system is similar to the usual Navier-Stokes equations, since it is weakly stable in the class of functions of finite energy.

Because of this crucial weak stability property, the study of hydrodynamic limits follows closely what has been previously done for the incompressible Navier-Stokes-Fourier limit of the Boltzmann equation (see \cite{arsenio3, SR} and the references therein for a survey of related results). In particular, we will be able to prove a convergence result which~:
\begin{itemize}
	\item holds globally in time~;
	\item does not require any assumption on the initial velocity profile~;
	\item does not assume any constraint on the initial thermodynamic fields.
\end{itemize}
We would also be able to take into account boundary conditions, and describe their limiting form, but this point will not be dealt with here. We refer to \cite{masmoudi2, SR} for a complete treatment of boundary conditions in the viscous hydrodynamic limits of the Boltzmann equation, based on the renormalized solutions on bounded domains constructed by Mischler in \cite{mischler, mischler2}.

As we will see, if we assume that the Vlasov-Maxwell-Boltzmann system \eqref{VMB1} has renormalized solutions (which, again, is not known), the main challenge here lies in understanding the influence of the electromagnetic force both on hypoelliptic processes of the kinetic transport equation and on fast time oscillations.

Our goal, here, is to establish the convergence of scaled families of renormalized solutions to the one species Vlasov-Maxwell-Boltzmann system \eqref{VMB1} towards solutions of the incompressible quasi-static Navier-Stokes-Fourier-Maxwell-Poisson system \eqref{NSFMP 2}, without any restriction on their size, regularity or well-preparedness of the initial data.

Following the program proposed by Bardos, Golse and Levermore in \cite{BGL2} (which relies essentially on weak compactness arguments), we can prove the following theorem. Recall that, in this second part, we are only considering the Maxwellian cross-section $b\equiv 1$.

\begin{thm}\label{NS-WEAKCV}
	Let $\left(f_{\eps}^\mathrm{in}, E_\eps^\mathrm{in}, B_\eps^\mathrm{in}\right)$ be a family of initial data such that
	\begin{equation}\label{init-fluctuation}
		\frac1{\eps^2}H\left(f_\eps^{\mathrm{in}}\right)
		+ \frac1{2}\int_{\mathbb{R}^3} \left(|E_\eps^{\rm in}|^2+ |B_\eps^{\rm in}|^2\right) dx \leq C^\mathrm{in},
	\end{equation}
	for some $C^\mathrm{in}>0$, and
	\begin{equation}\label{initial gauss}
		\Div E_\eps^\mathrm{in}=\int_{\mathbb{R}^3}g_\eps^\mathrm{in} M dv,\qquad \Div B_\eps^\mathrm{in} = 0,
	\end{equation}
	where $f_\eps^\mathrm{in}=M\left(1+\eps g_\eps^\mathrm{in}\right)$. For any $\eps>0$, we assume the existence of a renormalized solution $\left(f_\eps, E_\eps, B_\eps\right)$ to the scaled one species Vlasov-Maxwell-Boltzmann system \eqref{VMB1} (for the Maxwellian cross-section $b\equiv 1$) with initial data $\left(f_{\eps}^\mathrm{in}, E_\eps^\mathrm{in}, B_\eps^\mathrm{in}\right)$. We define the macroscopic fluctuations of density $\rho_\eps$, bulk velocity $u_\eps$ and temperature $\theta_\eps$ by
	\begin{equation*}
		\begin{aligned}
			\rho_\eps & =\int_{\mathbb{R}^3}g_\eps Mdv,\\
			u_\eps & =\int_{\mathbb{R}^3}g_\eps v Mdv,\\
			\theta_\eps & =\int_{\mathbb{R}^3}g_\eps\left(\frac{|v|^2}{3}-1\right) Mdv,
		\end{aligned}
	\end{equation*}
	and denote their respective initial value by $\rho_\eps^\mathrm{in}$, $u_\eps^\mathrm{in}$ and $\theta_\eps^\mathrm{in}$.
	
	% \begin{equation}
	% \frac1{\eps^2} H(f_\eps^0|M)  + \frac12 \int (|E_\eps^0|^2+|B_\eps^0|^2)dx\leq C_0\,.
	% \end{equation}
	% and satisfying further the weak convergences
	% $$\frac1\eps  P\int f_{\eps, in} vdv \rightharpoonup u_{in},\quad \frac1\eps \int  (f_{\eps, in}-M)(\frac15|v|^2 -1)dv \rightharpoonup \theta_{in},$$
	% in $L^1_{loc}(\Omega)$, 
	% where $P$ denotes the Leray projection onto divergence free vector fields.

	Then, the family $\left(\rho_\eps,u_\eps,\theta_\eps, B_\eps\right)$ is weakly relatively compact in $L^1_\mathrm{loc}(dtdx)$ (while the family of initial data $\left(\rho_\eps^\mathrm{in},u_\eps^\mathrm{in},\theta_\eps^\mathrm{in}, B_\eps^\mathrm{in}\right)$ is weakly relatively compact in $L^1_\mathrm{loc}(dx)$) and any of its limit points $\left(\rho,u,\theta, B\right)$ is a weak solution of the incompressible quasi-static Navier-Stokes-Fourier-Maxwell-Poisson system \eqref{NSFMP 2} with initial data
	\begin{equation*}
		\begin{aligned}
			\rho^\mathrm{in} & = \frac{\Delta_x}{3-5\Delta_x}\left(3\theta^\mathrm{in}_0-2\rho^\mathrm{in}_0\right),
			&
			u^\mathrm{in} & = \frac{\ROT}{1-\Delta_x}\left(\ROT u^\mathrm{in}_0+B^\mathrm{in}_0\right),
			\\
			\theta^\mathrm{in} & = \frac{1-\Delta_x}{3-5\Delta_x}\left(3\theta^\mathrm{in}_0-2\rho^\mathrm{in}_0\right),
			&
			B^\mathrm{in} & = \frac{1}{1-\Delta_x}\left(\ROT u^\mathrm{in}_0+B^\mathrm{in}_0\right),
		\end{aligned}
	\end{equation*}
	where $\left(\rho_0^\mathrm{in},u_0^\mathrm{in},\theta_0^\mathrm{in}, B_0^\mathrm{in}\right)\in L^2(dx)$ is the weak limit of $\left(\rho_\eps^\mathrm{in},u_\eps^\mathrm{in},\theta_\eps^\mathrm{in}, B_\eps^\mathrm{in}\right)$.
\end{thm}

The proof of Theorem \ref{NS-WEAKCV} is built over the course of the coming chapters and is per se the subject of Chapter \ref{grad}.

\bigskip

Note that, strictly speaking, the weak solution we obtain in the limit (and which depends in general on the subsequence under consideration) is not necessarily a Leray solution of the system \eqref{NSFMP 2}, since it is does not satisfy the energy inequality \eqref{energy}, but only a bound. However, it is possible to obtain asymptotically a Leray solution by strengthening the initial well-preparedness of the data. More precisely, one would have to impose that the initial data converges entropically (as introduced in \cite{BGL2}) in the sense that
\begin{equation*}
	\begin{aligned}
		\frac1{\eps^2}H\left(f_\eps^{\mathrm{in}}\right)
		& + \frac1{2}\int_{\mathbb{R}^3} \left(|E_\eps^{\rm in}|^2+ |B_\eps^{\rm in}|^2\right) dx \\
		& \rightarrow
		\frac1{2}\int_{\mathbb{R}^3} \left({\rho^\mathrm{in}}^2+\left|u^\mathrm{in}\right|^2+\frac 32{\theta^\mathrm{in}}^2+\left|\nabla_x\left(\rho^\mathrm{in}+\theta^\mathrm{in}\right)\right|^2+ |B^{\rm in}|^2\right) dx.
	\end{aligned}
\end{equation*}

In fact, the above entropic convergence has rather strong implications on the initial data. Indeed, further denoting by $g^\mathrm{in}_0$ and $E_0^\mathrm{in}$ the weak limits of $g_\eps^\mathrm{in}$ and $E_\eps^\mathrm{in}$, respectively, standard convexity arguments on weak convergence (see Lemma \ref{L1-lem} below, or \cite[Proposition 3.1]{BGL2}) yield that
\begin{equation*}
	\begin{aligned}
		\frac1{2}\int_{\mathbb{R}^3} \left({\rho_0^\mathrm{in}}^2+\left|u_0^\mathrm{in}\right|^2+\frac 32{\theta_0^\mathrm{in}}^2\right) dx &
		=\frac 12\int_{\mathbb{R}^3\times\mathbb{R}^3}\left(\Pi g^\mathrm{in}_0\right)^2Mdvdx
		\\
		&\leq \frac 12 \int_{\mathbb{R}^3\times\mathbb{R}^3}\left(g^\mathrm{in}_0\right)^2Mdvdx
		\leq
		\liminf_{\eps\rightarrow 0}
		\frac1{\eps^2}H\left(f_\eps^{\mathrm{in}}\right),
	\end{aligned}
\end{equation*}
and
\begin{equation*}
	\frac1{2}\int_{\mathbb{R}^3} \left(\left|E_0^\mathrm{in}\right|^2+ \left|B^{\rm in}_0\right|^2\right) dx
	\leq
	\liminf_{\eps\rightarrow 0}
	\frac1{2}\int_{\mathbb{R}^3} \left(\left|E_\eps^{\rm in}\right|^2+ \left|B_\eps^{\rm in}\right|^2\right) dx.
\end{equation*}
Moreover, one easily verifies that
\begin{equation*}
	\begin{pmatrix}
		\rho_0^\mathrm{in} \\ \sqrt{\frac 32}\theta_0^\mathrm{in} \\ E_0^\mathrm{in}
	\end{pmatrix}
	\mapsto
	\begin{pmatrix}
		\frac{3}{3-5\Delta_x}\left(\rho_0^\mathrm{in}-\Div E_0^\mathrm{in}\right)
		+
		\frac{\Delta_x}{3-5\Delta_x}\left(3\theta_0^\mathrm{in}-2\rho_0^\mathrm{in}\right)
		\\
		\frac{2}{3-5\Delta_x}\sqrt{\frac 32}\left(\rho_0^\mathrm{in}-\Div E_0^\mathrm{in}\right)
		+
		\frac{1-\Delta_x}{3-5\Delta_x}\sqrt{\frac 32}\left(3\theta_0^\mathrm{in}-2\rho_0^\mathrm{in}\right)
		\\
		\frac{5}{3-5\Delta_x}\nabla_x\left(\rho_0^\mathrm{in}-\Div E_0^\mathrm{in}\right)
		+
		\frac{1}{3-5\Delta_x}\nabla_x\left(3\theta_0^\mathrm{in}-2\rho_0^\mathrm{in}\right)
	\end{pmatrix}
\end{equation*}
defines the orthogonal projection onto the subspace of $L^2(dx)$ defined by the constraint
\begin{equation*}
	\nabla_x\left(\rho_0^\mathrm{in}+\theta_0^\mathrm{in}\right) = E_0^\mathrm{in},
\end{equation*}
while
\begin{equation*}
	\begin{pmatrix}
		u_0^\mathrm{in} \\ B_0^\mathrm{in}
	\end{pmatrix}
	\mapsto
	\begin{pmatrix}
		\frac{\ROT}{1-\Delta_x}\left(\ROT u_0^\mathrm{in}+ P B_0^\mathrm{in}\right)
		\\
		\frac{1}{1-\Delta_x}\left( \ROT u_0^\mathrm{in}+ P B_0^\mathrm{in}\right)
	\end{pmatrix}
\end{equation*}
corresponds to the orthogonal projection onto the subspace of $L^2(dx)$ defined by the constraints
\begin{equation*}
	\ROT B_0^\mathrm{in} = u_0^\mathrm{in}
	\qquad
	\text{and}
	\qquad
	\Div B_0^\mathrm{in} = 0.
\end{equation*}
Therefore, it follows that
\begin{equation*}
	\begin{aligned}
		\frac1{2}\int_{\mathbb{R}^3} & \left({\rho^\mathrm{in}}^2+\left|u^\mathrm{in}\right|^2+\frac 32{\theta^\mathrm{in}}^2+\left|\nabla_x\left(\rho^\mathrm{in}+\theta^\mathrm{in}\right)\right|^2+ \left|B^{\rm in}\right|^2\right) dx &
		\\
		& \hspace{20mm} \leq \frac1{2}\int_{\mathbb{R}^3} \left({\rho_0^\mathrm{in}}^2+\left|u_0^\mathrm{in}\right|^2+\frac 32{\theta_0^\mathrm{in}}^2+\left|E_0^\mathrm{in}\right|^2+ \left|B^{\rm in}_0\right|^2\right) dx,
	\end{aligned}
\end{equation*}
with equality if and only if $\rho^\mathrm{in}=\rho^\mathrm{in}_0$, $u^\mathrm{in}=u^\mathrm{in}_0$, $\theta^\mathrm{in}=\theta^\mathrm{in}_0$, $\nabla_x\left(\rho^\mathrm{in}+\theta^\mathrm{in}\right)=E^\mathrm{in}_0$ and $B^\mathrm{in}=B^\mathrm{in}_0$, which, when combined with the above entropic convergence of the initial data and according to Proposition 4.11 from \cite{BGL2}, implies that $\left(\rho_\eps^\mathrm{in},u_\eps^\mathrm{in},\theta_\eps^\mathrm{in}, E_\eps^\mathrm{in}, B_\eps^\mathrm{in}\right)$ converges strongly to $\left(\rho^\mathrm{in},u^\mathrm{in},\theta^\mathrm{in}, \nabla_x\left(\rho^\mathrm{in}+\theta^\mathrm{in}\right), B^\mathrm{in}\right)$, where
\begin{equation}\label{constraints quasi static}
	\Div u^\mathrm{in} = 0, \qquad \Div B^\mathrm{in} = 0, \qquad \ROT B^\mathrm{in} = u^\mathrm{in},\qquad \Delta_x \left(\rho^\mathrm{in} + \theta^\mathrm{in}\right) = \rho^\mathrm{in},
\end{equation}
and that $g_\eps^\mathrm{in}$ converges strongly to
\begin{equation*}
	\rho^{\mathrm{in}}+u^\mathrm{in}\cdot v + \theta^\mathrm{in}\left(\frac{|v|^2}{2}-\frac 32\right),
	\qquad \text{in } L^1_\mathrm{loc}\left(dx;L^1\left(\left(1+|v|^2\right)Mdv\right)\right).
\end{equation*}

The strong convergence of $g_\eps^\mathrm{in}$ towards an infinitesimal Maxwellian implies the vanishing of the initial relaxation layer, while the strong convergence of the initial macroscopic observables towards initial data satisfying the constraints \eqref{constraints quasi static} implies that there are no acoustic-electromagnetic waves. In fact, the weak convergence result in Theorem \ref{NS-WEAKCV} could be strengthened into a strong convergence result, for well-prepared initial data and provided that the limiting system has a unique solution satisfying the energy equality (see \cite[Theorem 7.4]{BGL2} on the strong Navier-Stokes limit).

It is to be emphasized that the generalized relative entropy method, which is developed later on in Chapter \ref{entropy method} and is used to prove Theorems \ref{CV-OMHD} and \ref{CV-OMHDSTRONG} below, is also applicable to the asymptotic regime studied in Theorem \ref{NS-WEAKCV}. This method would provide some strong convergence result even for ill-prepared initial data provided we can build an approximate solution which is smooth and accounts for the corrections due to the initial layer and the acoustic-electromagnetic waves.

% The counterpart of these statements for our asymptotic problem will be~:
% \begin{itemize}
% 	\item the weak convergence of  observables in the fast relaxation limit (which uses crucially the weak stability of the limiting system)~;
% 	\item the absence of a rate of convergence for general initial data (else this would lead to some uniqueness criterion for the limiting system)~;
% 	\item the strong convergence whenever the limiting solution is at least Lipschitz (which is the natural extension of the weak-strong uniqueness principle).
% \end{itemize}

% ==================
% = Main result II =
% ==================

\section[The two-fluid incompressible Navier-Stokes-Fourier-Maxwell\ldots]{The two-fluid incompressible Navier-Stokes-Fourier-Maxwell system with (solenoidal) Ohm's law}\label{main relative entropy}

According to Section \ref{formal two}, we consider now a plasma constituted of two species of oppositely charged particles with approximately equal mass, namely cations (positively charged ions) and anions (negatively charged ions). Elementary interactions are taken into account by both mean field terms (corresponding to long-range interactions) and some local collision terms (associated to short-range interactions) involving possibly different mean free paths. Thus, the charged particles evolve under the coupled effect of the Lorentz force due to the self-induced electromagnetic field, and of the collisions with other particles, according to the following scaled two species Vlasov-Maxwell-Boltzmann system~:
\begin{equation}\label{VMB2}
	\begin{cases}
		\begin{aligned}
			\eps \d_t f_\eps^\pm + v \cdot \nabla_x f_\eps^\pm
			\pm \delta \left( \eps E_\eps + v \wedge B_\eps \right) \cdot \nabla_v f_\eps^\pm &= \frac 1\eps Q(f_\eps^\pm,f_\eps^\pm) + \frac{\delta^2}{\eps} Q(f_\eps^\pm,f_\eps^\mp),
			\\
			f_\eps^\pm & =M\left(1+\epsilon g_\eps^\pm\right),
			\\
			\d_t E_\eps - \ROT B_\eps &= - \frac{\delta}{\eps} \int_{\mathbb{R}^3} \left(g_\eps^+-g_\eps^-\right)v M dv,
			\\
			\d_t B_\eps + \ROT E_\eps& = 0,
			\\
			\DIV E_\eps &=\delta\int_{\mathbb{R}^3} \left(g_\eps^+-g_\eps^-\right) M dv,
			\\
			\DIV B_\eps &=0,
		\end{aligned}
	\end{cases}
\end{equation}
where $\frac \delta \eps$ is asymptotically unbounded. In this scaling, the entropy inequality states that
\begin{equation}\label{entropy2}
	\begin{aligned}
		\frac1{\eps^2} H\left(f_\eps^{+}\right)
		+ \frac1{\eps^2} H\left(f_\eps^{-}\right)
		&
		+ \frac 1{2\eps^2}\int_{\mathbb{R}^3}\operatorname{Tr}m_\eps dx
		+ \frac 1{2} \int_{\mathbb{R}^3} \left(|E_\eps|^2+ |B_\eps|^2 + \operatorname{Tr}a_\eps\right) dx \\
		& +\frac{1}{\epsilon^4}\int_0^t\int_{\mathbb{R}^3}\left(D\left(f_\eps^+\right)+D\left(f_\eps^-\right)
		+ \delta^2 D\left(f_\eps^+,f_\eps^-\right)\right)(s) dx ds
		\\
		& \leq
		\frac1{\eps^2} H\left(f_\eps^{+\mathrm{in}}\right)
		+ \frac1{\eps^2} H\left(f_\eps^{-\mathrm{in}}\right)
		+ \frac1{2}\int_{\mathbb{R}^3} \left(|E_\eps^{\rm in}|^2+ |B_\eps^{\rm in}|^2\right) dx,
	\end{aligned}
\end{equation}
where $H\left(f_\eps^{\pm}\right)=H\left(f_\eps^{\pm}|M\right)$ and the symmetric positive definite matrix measures $m_\eps$ and $a_\eps$ are the defects introduced in Section \ref{macroscopic defects} stemming from the terms $\int_{\mathbb{R}^3}\left(f_\eps^++f_\eps^-\right)v\otimes v dv$ and
	$\begin{pmatrix}
		E_\eps\\B_\eps
	\end{pmatrix}
	\otimes
	\begin{pmatrix}
		E_\eps\\B_\eps
	\end{pmatrix}$, respectively. In particular, it yields uniform bounds on $E_\eps$, $B_\eps$ and $g_\eps^\pm$.

Since we are interested in the limiting fluctuation $g_\eps^\pm\rightharpoonup g^\pm$, it is then natural to rewrite the kinetic equations in terms of the fluctuations $g_\eps^\pm$,
\begin{equation}\label{VMB-fluct}
	\begin{aligned}
		\eps\d_t \begin{pmatrix} g_\eps ^+\\ g_\eps^- \end{pmatrix} + v & \cdot \nabla_x \begin{pmatrix} g_\eps ^+\\ g_\eps^- \end{pmatrix}
		+ \delta (\eps E_\eps+ v\wedge B_\eps) \cdot \nabla_v \begin{pmatrix} g_\eps ^+\\ - g_\eps^- \end{pmatrix}
		- \delta  E_\eps \cdot v \begin{pmatrix} 1+\eps g_\eps ^+\\ -1-\eps g_\eps^- \end{pmatrix} \\
		% & = \frac 1\eps
		% \begin{pmatrix} - \cL g_\eps^+ +\eps\cQ(g_\eps^+,g_\eps^+) \\ - \cL g_\eps^- +\eps\cQ(g_\eps^-,g_\eps^-) \end{pmatrix}
		% +
		% \frac{\delta^2}{\eps}
		% \begin{pmatrix} - \cL \left(g_\eps^+,g_\eps^-\right) + \eps \cQ\left(g_\eps^+,g_\eps^-\right) \\ - \cL \left(g_\eps^-,g_\eps^+\right) + \eps \cQ\left(g_\eps^-,g_\eps^+\right) \end{pmatrix}\\
		& = -\frac 1\eps
		\begin{pmatrix} \cL g_\eps^+ + \delta^2\cL \left(g_\eps^+,g_\eps^-\right) \\
		\cL g_\eps^- + \delta^2 \cL \left(g_\eps^-,g_\eps^+\right) \end{pmatrix}
		+
		\begin{pmatrix} \cQ(g_\eps^+,g_\eps^+) + \delta^2 \cQ\left(g_\eps^+,g_\eps^-\right) \\ \cQ(g_\eps^-,g_\eps^-) + \delta^2 \cQ\left(g_\eps^-,g_\eps^+\right) \end{pmatrix}.
	\end{aligned}
\end{equation}

According to the formal analysis from Section \ref{formal two}, we then expect the limiting macroscopic observables to solve, in the case of strong interspecies collisions $\delta=1$, the two-fluid incompressible Navier-Stokes-Fourier-Maxwell system with Ohm's law~:
\begin{equation}\label{TFINSFMO 2}
	\begin{cases}
		\begin{aligned}
			\d_t u +
			u\cdot\nabla_x u - \mu\Delta_x u
			& = -\nabla_x p+
			\frac 12 \left(nE + j \wedge B\right) , & \Div u & = 0,\\
			\d_t \theta
			+
			u\cdot\nabla_x\theta - \kappa \Delta_x\theta
			& = 0, & \rho+\theta & = 0, \\
			\d_t E - \ROT B &= -  j, & \Div E & = n,
			\\
			\d_t B + \ROT E & = 0, & \Div B & = 0, \\
			j-nu & = \sigma\left(-\frac 12 \nabla_x n + E + u\wedge B\right), & w & =n\theta,
		\end{aligned}
	\end{cases}
\end{equation}
and, in the case of weak interspecies collisions $\delta=o(1)$, with $\frac\delta\eps$ unbounded, the two-fluid incompressible Navier-Stokes-Fourier-Maxwell system with solenoidal Ohm's law~:
\begin{equation}\label{TFINSFMSO 2}
	\begin{cases}
		\begin{aligned}
			\d_t u +
			u\cdot\nabla_x u - \mu\Delta_x u
			& = -\nabla_x p+
			\frac 12 j \wedge B , & \Div u & = 0,\\
			\d_t \theta
			+
			u\cdot\nabla_x\theta - \kappa \Delta_x\theta
			& = 0, & \rho+\theta & = 0, \\
			\d_t E - \ROT B &= -  j, & \Div E & = 0,
			\\
			\d_t B + \ROT E & = 0, & \Div B & = 0, \\
			j & = \sigma\left(- \nabla_x \bar p + E + u\wedge B\right), & \Div j & =0,\\
			n&=0, & w&=0.
		\end{aligned}
	\end{cases}
\end{equation}
As previously emphasized in Section \ref{stability existence 2}, the above limiting models \eqref{TFINSFMO 2} and \eqref{TFINSFMSO 2} are not stable under weak convergence in the energy space and, thus, share more similarities with the three-dimensional incompressible Euler equations.

Our goal, here, is to establish the convergence of scaled families of renormalized solutions to the two species Vlasov-Maxwell-Boltzmann system \eqref{VMB2} towards dissipative solutions of the two-fluid incompressible Navier-Stokes-Fourier-Maxwell system with Ohm's law \eqref{TFINSFMO 2}, in the case of strong interspecies interactions $\delta=1$, or with solenoidal Ohm's law \eqref{TFINSFMSO 2}, in the case of weak interspecies interactions $\delta=o(1)$, with $\frac\delta\eps$ unbounded, without any restriction on their size or regularity. We will, however, impose some well-preparedness of the initial data.

Improving on the program by the second author completed in \cite{SR2, SR3} (which relies essentially on modulated entropy arguments), we can prove the following theorems. Recall that, in this second part, we are only considering the Maxwellian cross-section $b\equiv 1$.

\subsection{Weak interactions}

In the case of weak interspecies interactions $\delta=o(1)$, with $\frac\delta\eps$ unbounded, we have the following result.

\begin{thm}\label{CV-OMHD}
	Let $\left(f_{\eps}^{\pm\mathrm{in}}, E_{\eps}^\mathrm{in}, B_{\eps}^\mathrm{in}\right)$ be a family of initial data such that
	\begin{equation}\label{init-fluctuation 2}
		\frac1{\eps^2} H\left(f_\eps^{+\mathrm{in}}\right)
		+ \frac1{\eps^2} H\left(f_\eps^{-\mathrm{in}}\right)
		+ \frac 1{2} \int_{\mathbb{R}^3} \left(|E_\eps^\mathrm{in}|^2+ |B_\eps^\mathrm{in}|^2\right) dx \leq C^\mathrm{in},
	\end{equation}
	for some $C^\mathrm{in}>0$, and
	\begin{equation}\label{initial gauss 2}
		\Div E_\eps^\mathrm{in}=\delta\int_{\mathbb{R}^3}\left(g_\eps^{+\mathrm{in}}-g_\eps^{-\mathrm{in}}\right) M dv,\qquad \Div B_\eps^\mathrm{in} = 0,
	\end{equation}
	where $f_\eps^{\pm\mathrm{in}}=M\left(1+\eps g_\eps^{\pm\mathrm{in}}\right)$. We further assume that the initial data is well-prepared in the sense that, as $\eps\rightarrow 0$,
	\begin{equation*}
		g_\eps^{\pm\mathrm{in}} \rightharpoonup g^{\mathrm{in}}
		= \rho^{\mathrm{in}}+u^\mathrm{in}\cdot v + \theta^\mathrm{in}\left(\frac{|v|^2}{2}-\frac 32\right)
		\qquad \text{in } \textit{w-}L^1_\mathrm{loc}\left(Mdxdv\right),
	\end{equation*}
	where $\rho^{\mathrm{in}}, u^{\mathrm{in}}, \theta^{\mathrm{in}} \in L^2\left(dx\right)$ satisfy the incompressibility and Boussinesq constraints
	\begin{equation*}
		\DIV u^\mathrm{in}=0,\qquad \rho^\mathrm{in}+\theta^{\mathrm{in}}=0,
	\end{equation*}
	and that the following strong convergences hold, as $\eps\rightarrow 0$,
	\begin{equation}\label{well-prepared init data}
		\begin{aligned}
			\frac1{\eps^2} H\left(f_\eps^{+\mathrm{in}}\right)
			+ \frac1{\eps^2} H\left(f_\eps^{-\mathrm{in}}\right)
			& \rightarrow \int_{\mathbb{R}^3\times\mathbb{R}^3} \left(g^{\mathrm{in}}\right)^2
			Mdvdx
			=\int_{\mathbb{R}^3} \left|u^\mathrm{in}\right|^2 +\frac 52\left(\theta^\mathrm{in}\right)^2dx, \\
			E_\eps^\mathrm{in} & \rightarrow E^\mathrm{in},
			\qquad \text{in } L^2\left(dx\right), \\
			B_\eps^\mathrm{in} & \rightarrow B^\mathrm{in},
			\qquad \text{in } L^2\left(dx\right),
		\end{aligned}
	\end{equation}
	for some $E^\mathrm{in},B^\mathrm{in}\in L^2(dx)$. In particular, in view of \eqref{initial gauss 2}, it necessarily holds that
	\begin{equation*}
		\Div E^\mathrm{in} = 0, \qquad \Div B^\mathrm{in}=0.
	\end{equation*}
	For any $\eps>0$, we assume the existence of a renormalized solution $\left(f_\eps^\pm, E_\eps, B_\eps\right)$ to the scaled two species Vlasov-Maxwell-Boltzmann system \eqref{VMB2} (for the Maxwellian cross-section $b\equiv 1$), where $\delta=o(1)$ and $\frac\delta\eps$ is asymptotically unbounded, with initial data $\left(f_{\eps}^{\pm\mathrm{in}}, E_{\eps}^\mathrm{in}, B_{\eps}^\mathrm{in}\right)$. We define the macroscopic fluctuations of density $\rho_\eps^\pm$, bulk velocity $u_\eps^\pm$ and temperature $\theta_\eps^\pm$ by
	\begin{equation*}
		\begin{aligned}
			\rho_\eps^\pm & =\int_{\mathbb{R}^3}g_\eps^\pm Mdv,\\
			u_\eps^\pm & =\int_{\mathbb{R}^3}g_\eps^\pm v Mdv,\\
			\theta_\eps^\pm & =\int_{\mathbb{R}^3}g_\eps^\pm\left(\frac{|v|^2}{3}-1\right) Mdv.
		\end{aligned}
	\end{equation*}
	We finally define the hydrodynamic variables
	\begin{equation*}
		\rho_\eps = \frac{\rho^+_\eps+\rho^-_\eps}{2}, \qquad
		u_\eps = \frac{u^+_\eps + u^-_\eps}{2}, \qquad
		\theta_\eps = \frac{\theta^+_\eps+\theta_\eps^-}{2},
	\end{equation*}
	and electrodynamic variables
	\begin{equation*}
		n_\eps = \rho^+_\eps - \rho^-_\eps, \qquad
		j_\eps = \frac\delta\eps\left(u^+_\eps - u^-_\eps\right), \qquad
		w_\eps = \frac\delta\eps\left(\theta^+_\eps -\theta^-_\eps\right).
	\end{equation*}

	Then, the family $\left(\rho_\eps,u_\eps,\theta_\eps, n_\eps, j_\eps, w_\eps, E_\eps, B_\eps\right)$ is weakly relatively compact in $L^1_\mathrm{loc}(dtdx)$ and any of its limit points $\left(\rho,u,\theta, n, j, w, E, B\right)$ is a dissipative solution of the two-fluid incompressible Navier-Stokes-Fourier-Maxwell system with solenoidal Ohm's law \eqref{TFINSFMSO 2} with initial data $\left(u^\mathrm{in},\theta^\mathrm{in}, E^\mathrm{in}, B^\mathrm{in}\right)$~--~that is, it verifies the energy inequality corresponding to \eqref{TFINSFMSO 2}, it enjoys the weak temporal continuity $\left(u,\theta,E,B\right)\in C\left([0,\infty);\textit{w-}L^2\left(\mathbb{R}^3\right)\right)$, it solves the system
	\begin{equation*}
		\begin{cases}
			\begin{aligned}
				\Div u & = 0, & \rho+\theta & =0,\\
				\d_t E - \ROT B &= -  j, & \Div E & = 0,
				\\
				\d_t B + \ROT E & = 0, & \Div B & = 0, \\
				j & = \sigma\left(- \nabla_x \bar p + E + u\wedge B\right), & \Div j & =0,\\
				n&=0, & w&=0,
			\end{aligned}
		\end{cases}
	\end{equation*}
	in the sense of distributions, and it satisfies the stability inequality
	\begin{equation*}
		\begin{aligned}
			& \delta\mathcal{E}(t) + \frac 12 \int_0^t \delta\mathcal{D}(s) e^{\int_s^t\lambda(\sigma)d\sigma}ds
			\\
			& \leq \delta\mathcal{E}(0) e^{\int_0^t\lambda(s)ds}
			+\int_0^t
			\int_{\mathbb{R}^3} \mathbf{A}\cdot
			\begin{pmatrix}
				u-\bar u \\ \frac 52\left(\theta-\bar\theta\right) \\
				j-\bar j
				\\ E-\bar E +\bar u \wedge \left(B-\bar B\right) \\ B-\bar B +\left(E-\bar E\right)\wedge\bar u
			\end{pmatrix}(s)
			dx
			e^{\int_s^t\lambda(\sigma)d\sigma}ds,
		\end{aligned}
	\end{equation*}
	for any test functions $\left(\bar u,\bar\theta,\bar j,\bar E,\bar B\right)\in C^\infty_c\left([0,\infty)\times\mathbb{R}^3\right)$ with
	\begin{equation*}
		\Div\bar u = \Div \bar j = \Div \bar E = \Div \bar B = 0
		\quad\text{and}\quad
		\left\|\bar u\right\|_{L^\infty(dtdx)}<1,
	\end{equation*}
	where the modulated energy and modulated energy dissipation are respectively given by
	\begin{equation*}
		\begin{aligned}
			\delta\mathcal{E}(t)
			& =
			\left\|u-\bar u\right\|_{L^2(dx)}^2+\frac 52\left\|\theta-\bar \theta\right\|_{L^2(dx)}^2
			+\frac 12\left\|E-\bar E\right\|_{L^2(dx)}^2
			+\frac 12\left\|B-\bar B\right\|_{L^2(dx)}^2
			\\
			& - \int_{\mathbb{R}^3}
			\left(\left(E_\eps-\bar E\right)\wedge\left(B_\eps-\bar B\right)
			\right)\cdot\bar u
			dx,
			\\
			\delta\mathcal{D}(t)
			& =
			2\mu
			\left\|\nabla_x \left(u-\bar u\right)\right\|_{L^2_x}^2
			+ 5\kappa
			\left\|\nabla_x\left(\theta-\bar\theta\right)\right\|_{L^2_x}^2
			+ \frac 1\sigma
			\left\|j-\bar j\right\|_{L^2_x}^2,
		\end{aligned}
	\end{equation*}
	the acceleration operator is defined by
	\begin{equation*}
		\mathbf{A} \left( \bar u, \bar \theta, \bar j, \bar E, \bar B\right)
		=
		\begin{pmatrix}
			-2\left(\d_t \bar u +
			P\left(\bar u\cdot\nabla_x \bar u\right) - \mu\Delta_x \bar u\right)
			+ P \left(\bar j \wedge \bar B\right)
			\\
			-2\left(\partial_t\bar\theta + \bar u \cdot\nabla_x\bar \theta - \kappa\Delta_x\bar \theta\right)
			\\
			- \frac 1{\sigma}\bar j + P\left(\bar E + \bar u\wedge \bar B\right)
			\\
			-\left(\partial_t\bar E - \rot\bar B + \bar j\right)
			\\
			-\left(\partial_t\bar B + \rot\bar E\right)
		\end{pmatrix},
	\end{equation*}
	and the growth rate is given by
	\begin{equation*}
		\begin{aligned}
			\lambda(t) =
			C\Bigg(\frac{\left\|\bar u(t)\right\|_{W^{1,\infty}\left(dx\right)}+\left\|\partial_t\bar u(t)\right\|_{L^\infty(dx)}+\left\|\bar j(t)\right\|_{L^\infty(dx)}}
			{1-\left\|\bar u(t)\right\|_{L^\infty(dx)}} &
			\\
			+\left\|\bar\theta(t)\right\|_{W^{1,\infty}(dx)} & +\left\|\bar\theta(t)\right\|_{W^{1,\infty}(dx)}^2\Bigg),
		\end{aligned}
	\end{equation*}
	with a constant $C>0$ independent of test functions.
	
	In particular, this dissipative solution coincides with the unique smooth solution with velocity field bounded pointwise by the speed of light (i.e.\ $\left\|u\right\|_{L^\infty(dtdx)}<1$) as long as the latter exists.
\end{thm}

The proof of Theorem \ref{CV-OMHD} is built over the course of the coming chapters and is per se the subject of Section \ref{proof of theorem weak}.

\subsection{Strong interactions}

In the case of strong interspecies interactions $\delta=1$, we have the following result.

\begin{thm}\label{CV-OMHDSTRONG}
	Let $\left(f_{\eps}^{\pm\mathrm{in}}, E_{\eps}^\mathrm{in}, B_{\eps}^\mathrm{in}\right)$ be a family of initial data such that
	\begin{equation}\label{init-fluctuation 2 strong}
		\frac1{\eps^2} H\left(f_\eps^{+\mathrm{in}}\right)
		+ \frac1{\eps^2} H\left(f_\eps^{-\mathrm{in}}\right)
		+ \frac 1{2} \int_{\mathbb{R}^3} \left(|E_\eps^\mathrm{in}|^2+ |B_\eps^\mathrm{in}|^2\right) dx \leq C^\mathrm{in},
	\end{equation}
	for some $C^\mathrm{in}>0$, and
	\begin{equation}\label{initial gauss 2 strong}
		\Div E_\eps^\mathrm{in}=\int_{\mathbb{R}^3}\left(g_\eps^{+\mathrm{in}}-g_\eps^{-\mathrm{in}}\right) M dv,\qquad \Div B_\eps^\mathrm{in} = 0,
	\end{equation}
	where $f_\eps^{\pm\mathrm{in}}=M\left(1+\eps g_\eps^{\pm\mathrm{in}}\right)$. We further assume that the initial data is well-prepared in the sense that, as $\eps\rightarrow 0$,
	\begin{equation*}
		g_\eps^{\pm\mathrm{in}} \rightharpoonup g^{\pm\mathrm{in}}
		= \rho^{\pm\mathrm{in}}+u^\mathrm{in}\cdot v + \theta^\mathrm{in}\left(\frac{|v|^2}{2}-\frac 32\right)
		\qquad \text{in } \textit{w-}L^1_\mathrm{loc}\left(Mdxdv\right),
	\end{equation*}
	where $\rho^{\pm\mathrm{in}}, u^{\mathrm{in}}, \theta^{\mathrm{in}} \in L^2\left(dx\right)$ satisfy the incompressibility and Boussinesq constraints, denoting $\rho^{\mathrm{in}}=\frac{\rho^{+\mathrm{in}}+\rho^{-\mathrm{in}}}{2}$,
	\begin{equation*}
		\DIV u^\mathrm{in}=0,\qquad \rho^\mathrm{in}+\theta^{\mathrm{in}}=0,
	\end{equation*}
	and that the following strong convergences hold, as $\eps\rightarrow 0$,
	\begin{equation}\label{well-prepared init data strong}
		\begin{aligned}
			\frac1{\eps^2} H\left(f_\eps^{+\mathrm{in}}\right)
			+ \frac1{\eps^2} H\left(f_\eps^{-\mathrm{in}}\right)
			& \rightarrow \int_{\mathbb{R}^3\times\mathbb{R}^3} \frac 12 \left(g^{+\mathrm{in}}\right)^2 + \frac 12 \left(g^{-\mathrm{in}}\right)^2
			Mdvdx
			\\
			& =\int_{\mathbb{R}^3} \frac 14 \left(n^\mathrm{in}\right)^2 + \left|u^\mathrm{in}\right|^2 +\frac 52\left(\theta^\mathrm{in}\right)^2dx, \\
			E_\eps^\mathrm{in} & \rightarrow E^\mathrm{in},
			\qquad \text{in } L^2\left(dx\right), \\
			B_\eps^\mathrm{in} & \rightarrow B^\mathrm{in},
			\qquad \text{in } L^2\left(dx\right),
		\end{aligned}
	\end{equation}
	for some $E^\mathrm{in},B^\mathrm{in}\in L^2(dx)$. In particular, in view of \eqref{initial gauss 2 strong}, it necessarily holds that, denoting $n^{\mathrm{in}}=\rho^{+\mathrm{in}}-\rho^{-\mathrm{in}}$,
	\begin{equation*}
		\Div E^\mathrm{in} = n^\mathrm{in}, \qquad \Div B^\mathrm{in}=0.
	\end{equation*}
	For any $\eps>0$, we assume the existence of a renormalized solution $\left(f_\eps^\pm, E_\eps, B_\eps\right)$ to the scaled two species Vlasov-Maxwell-Boltzmann system \eqref{VMB2} (for the Maxwellian cross-section $b\equiv 1$), where $\delta=1$, with initial data $\left(f_{\eps}^{\pm\mathrm{in}}, E_{\eps}^\mathrm{in}, B_{\eps}^\mathrm{in}\right)$. We define the macroscopic fluctuations of density $\rho_\eps^\pm$, bulk velocity $u_\eps^\pm$ and temperature $\theta_\eps^\pm$ by
	\begin{equation*}
		\begin{aligned}
			\rho_\eps^\pm & =\int_{\mathbb{R}^3}g_\eps^\pm Mdv,\\
			u_\eps^\pm & =\int_{\mathbb{R}^3}g_\eps^\pm v Mdv,\\
			\theta_\eps^\pm & =\int_{\mathbb{R}^3}g_\eps^\pm\left(\frac{|v|^2}{3}-1\right) Mdv.
		\end{aligned}
	\end{equation*}
	We finally define the hydrodynamic variables
	\begin{equation*}
		\rho_\eps = \frac{\rho^+_\eps+\rho^-_\eps}{2}, \qquad
		u_\eps = \frac{u^+_\eps + u^-_\eps}{2}, \qquad
		\theta_\eps = \frac{\theta^+_\eps+\theta_\eps^-}{2},
	\end{equation*}
	and electrodynamic variables
	\begin{equation*}
		n_\eps = \rho^+_\eps - \rho^-_\eps, \qquad
		j_\eps = \frac 1\eps\left(u^+_\eps - u^-_\eps\right), \qquad
		w_\eps = \frac 1\eps\left(\theta^+_\eps -\theta^-_\eps\right).
	\end{equation*}

	Then, the family $\left(\rho_\eps,u_\eps,\theta_\eps, n_\eps, j_\eps, E_\eps, B_\eps\right)$ (note that we have excluded the variable $w_\eps$) is relatively compact in the sense that for every sequence in this family there exists a subsequence such that
	\begin{equation*}
		\left(\rho_\eps,u_\eps,\theta_\eps, n_\eps, r_\eps j_\eps, E_\eps, B_\eps\right)
		\rightharpoonup
		\left(\rho,u,\theta, n, j, E, B\right) \qquad \text{in }\textit{w-}L^1_\mathrm{loc}(dtdx),
	\end{equation*}
	where $r_\eps(t,x)$ is a sequence of mesurable scalar functions converging almost everywhere towards the constant function $1$.
	
	Moreover, up to further extraction of subsequences, one also has the convergence
	\begin{equation*}
		r_\eps\frac 1\eps\left(g_\eps^+-g_\eps^--n_\eps\right)
		\rightharpoonup h
		\qquad \text{in }\textit{w-}L^1_\mathrm{loc}\left(dtdx;L^1\left(\left(1+|v|\right)Mdv\right)\right).
	\end{equation*}
	We set
	\begin{equation*}
		w=\int_{\mathbb{R}^3}h\left(\frac{|v|^2}{3}-1\right) Mdv.
	\end{equation*}
	Note that $w$ is not necessarily a limit point of $w_\eps$.

	Any such limit point $\left(\rho,u,\theta, n, j, w, E, B\right)$ is a dissipative solution of the two-fluid incompressible Navier-Stokes-Fourier-Maxwell system with Ohm's law \eqref{TFINSFMO 2} with initial data $\left(u^\mathrm{in},\theta^\mathrm{in},n^\mathrm{in}, E^\mathrm{in}, B^\mathrm{in}\right)$~--~that is, it verifies the energy inequality corresponding to \eqref{TFINSFMO 2}, it enjoys the weak temporal continuity $\left(u,\theta,n,E,B\right)\in C\left([0,\infty);\textit{w-}L^2\left(\mathbb{R}^3\right)\right)$, it solves the system
	\begin{equation*}
		\begin{cases}
			\begin{aligned}
				\Div u & = 0,
				& \rho+\theta & = 0, \\
				\Div E & = n, & \Div B & = 0, \\
				\d_t B + \ROT E & = 0, &
				j-nu & = \sigma\left(-\frac 12 \nabla_x n + E + u\wedge B\right), \\
				&& w & = n\theta,
			\end{aligned}
		\end{cases}
	\end{equation*}
	in the sense of distributions (note that we have voluntarily left Amp\`ere's equation out of the above system), and it satisfies the stability inequality
	\begin{equation*}
		\begin{aligned}
			\delta\mathcal{E}(t) & + \frac 12 \int_0^t \delta\mathcal{D}(s) e^{\int_s^t\lambda(\sigma)d\sigma}ds
			\\
			& \leq \delta\mathcal{E}(0) e^{\int_0^t\lambda(s)ds}
			\\
			& +\int_0^t
			\int_{\mathbb{R}^3} \mathbf{A}\cdot
			\begin{pmatrix}
				u-\bar u \\ \frac 52\left(\theta-\bar\theta\right) \\
				j - nu -\left(\bar j-\bar n\bar u\right)
				\\ E-\bar E +\bar u \wedge \left(B-\bar B\right) - \frac 12 \nabla_x\left(n-\bar n\right) \\ B-\bar B +\left(E-\bar E\right)\wedge\bar u
			\end{pmatrix}(s)
			dx
			e^{\int_s^t\lambda(\sigma)d\sigma}ds,
		\end{aligned}
	\end{equation*}
	for any test functions $\left(\bar u,\bar\theta,\bar n,\bar j,\bar E,\bar B\right)\in C^\infty_c\left([0,\infty)\times\mathbb{R}^3\right)$ with
	\begin{equation*}
		\Div\bar u = \Div \bar B = 0,\quad \Div \bar E= \bar n
		\quad\text{and}\quad
		\left\|\bar u\right\|_{L^\infty(dtdx)}<1,
	\end{equation*}
	where the modulated energy and modulated energy dissipation are respectively given by
	\begin{equation*}
		\begin{aligned}
			\delta\mathcal{E}(t)
			& = \frac 14 \left\|n-\bar n\right\|_{L^2(dx)}^2 +
			\left\|u-\bar u\right\|_{L^2(dx)}^2+\frac 52\left\|\theta-\bar \theta\right\|_{L^2(dx)}^2
			\\
			& +\frac 12\left\|E-\bar E\right\|_{L^2(dx)}^2
			+\frac 12\left\|B-\bar B\right\|_{L^2(dx)}^2
			- \int_{\mathbb{R}^3}
			\left(\left(E_\eps-\bar E\right)\wedge\left(B_\eps-\bar B\right)
			\right)\cdot\bar u
			dx,
			\\
			\delta\mathcal{D}(t)
			& =
			2\mu
			\left\|\nabla_x \left(u-\bar u\right)\right\|_{L^2_x}^2
			+ 5\kappa
			\left\|\nabla_x\left(\theta-\bar\theta\right)\right\|_{L^2_x}^2
			+ \frac 1\sigma
			\left\|\left(j-nu\right)-\left(\bar j-\bar n\bar u\right)\right\|_{L^2_x}^2,
		\end{aligned}
	\end{equation*}
	the acceleration operator is defined by
	\begin{equation*}
		\mathbf{A} \left( \bar u, \bar \theta, \bar j, \bar E, \bar B\right)
		=
		\begin{pmatrix}
			-2\left(\d_t \bar u +
			P\left(\bar u\cdot\nabla_x \bar u\right) - \mu\Delta_x \bar u\right)
			+ P \left(\bar n\bar E + \bar j \wedge \bar B\right)
			\\
			-2\left(\partial_t\bar\theta + \bar u \cdot\nabla_x\bar \theta - \kappa\Delta_x\bar \theta\right)
			\\
			- \frac 1{\sigma}\left(\bar j-\bar n\bar u\right) - \frac 12 \nabla_x\bar n + \bar E + \bar u\wedge \bar B
			\\
			-\left(\partial_t\bar E - \rot\bar B + \bar j\right)
			\\
			-\left(\partial_t\bar B + \rot\bar E\right)
		\end{pmatrix},
	\end{equation*}
	and the growth rate is given by
	\begin{equation*}
		\begin{aligned}
			& \lambda(t) =
			\\
			& C\Bigg(\frac{\left\|\bar u(t)\right\|_{W^{1,\infty}\left(dx\right)}+\left\|\partial_t\bar u(t)\right\|_{L^\infty(dx)}
			+\left\|\bar\theta(t)\right\|_{W^{1,\infty}(dx)}
			+\left\|\left(\bar j-\bar n\bar u\right)(t)\right\|_{L^\infty(dx)}}
			{1-\left\|\bar u(t)\right\|_{L^\infty(dx)}}
			\\
			& \hspace{35mm} +\left\|\bar\theta(t)\right\|_{W^{1,\infty}(dx)}^2
			+\left\|\left(\frac 12 \nabla_x\bar n - \bar E - \bar u\wedge\bar B\right)(t)\right\|_{L^\infty(dx)} \Bigg),
		\end{aligned}
	\end{equation*}
	with a constant $C>0$ independent of test functions.
	
	In particular, this dissipative solution coincides with the unique smooth solution with velocity field bounded pointwise by the speed of light (i.e.\ $\left\|u\right\|_{L^\infty(dtdx)}<1$) as long as the latter exists.
\end{thm}

The proof of Theorem \ref{CV-OMHDSTRONG} is built over the course of the coming chapters and is per se the subject of Section \ref{proof of theorem strong}.

\bigskip

In both Theorems \ref{CV-OMHD} and \ref{CV-OMHDSTRONG}, we focus on the case of well-prepared initial data. That is to say, we assume that the initial distribution has a velocity profile close to local thermodynamic equilibrium
\begin{equation*}
	g^{\pm \mathrm{in}} = \rho^{\pm\mathrm{in}}+u^\mathrm{in}\cdot v + \theta^\mathrm{in}\left(\frac{|v|^2}{2}-\frac 32\right),
\end{equation*}
(with $g^{+\mathrm{in}}=g^{-\mathrm{in}}$ in the case of weak interactions) in order that there is no relaxation layer, and that the  asymptotic initial thermodynamic fields  satisfy the incompressibility and Boussinesq constraints
\begin{equation*}
	\DIV u^\mathrm{in}=0,\qquad \rho^\mathrm{in}+\theta^{\mathrm{in}}=0,
\end{equation*}
which ensures that there are no acoustic waves.

The case of ill-prepared initial data could be handled by constructing an accurate approximate solution as in \cite{SR3}. The corresponding result should be even better in the present viscous incompressible regime because we can control conservation defects and fluxes without any additional integrability assumptions on renormalized solutions to \eqref{VMB2}. Note, however, that such a result would still be conditional as the existence of renormalized solutions to \eqref{VMB2} has to be assumed.

Relaxing the regularity assumption on the asymptotic solution would require new ideas~: the stability in the energy and entropy methods is indeed controlled by higher integrability or regularity norms of the limiting fields. As discussed in Section \ref{stability existence 2}, the two-fluid incompressible Navier-Stokes-Fourier-Maxwell systems with (solenoidal) Ohm's law \eqref{TFINSFMO 2} and \eqref{TFINSFMSO 2} are not known to have weak solutions, so that we do not expect to extend our convergence results for distributional solutions with low regularity.

\section{Outline of proofs}

We expect the Vlasov-Maxwell-Boltzmann systems \eqref{VMB1} and \eqref{VMB2} to exhibit very different qualitative behaviors in the three asymptotic scalings we consider~: one species, two species with weak interactions, and two species with strong interactions.
However, estimates coming directly from the entropy inequalities \eqref{entropy1} and \eqref{entropy2} and leading to weak compactness results are similar in all regimes, so we will gather them in Chapter \ref{weak bounds}. We will also obtain the thermodynamic equilibria coming from relaxation estimates in Chapter \ref{weak bounds}.

Then, Chapter \ref{constraints proof} will be devoted to the derivation of constraints which are stable under weak convergence and can be handled with the weak bounds from Chapter \ref{weak bounds}. These constraints include, for instance, some lower order macroscopic constraints (such as the Boussinesq and incompressibility relations) for one species or two species with weak interactions. We will also establish the limiting energy inequalities for one species and two species with weak interactions and discuss the limiting form of Maxwell's system. This chapter does not handle the constraints pertaining to two species with strong interactions, which will require the more advanced techniques of the following chapters.

A major difference between regimes appears in Chapter \ref{hypoellipticity} regarding spatial regularity. The basic idea is to use the hypoellipticity of the free transport operator --~as studied in \cite{arsenio}~-- to transfer regularity from the $v$ variable to the $x$ variable. But, because of the singularity in the Lorentz force, source terms in the kinetic equations are of different sizes so that different renormalizations of the kinetic equations will have to be considered for the three different regimes. Roughly speaking, we will be able to establish some strong compactness and equi-integrability on the fluctuations in the less singular regime with only one species, and only some weaker analog on some truncated fluctuations for two species (with a truncation depending on the asymptotic parameter $\delta$).

Another important difference comes from the nonlinear constraints (which occur only in the cases of two species). We will see in Chapter \ref{high constraints proof} that, for weak interactions, Ohm's law is obtained as a higher order singular perturbation. Its derivation will use the renormalized form of the kinetic equations as well as the partial equi-integrability established in Chapter \ref{hypoellipticity}. In the case of strong interactions, even though Ohm's law appears at leading order, proving its stability is much more intricate as it will require an additional macroscopic renormalization. Chapter \ref{high constraints proof} will also contain the derivation of the energy inequality for strong interactions, which will require the use of the strong compactness bounds from Chapter \ref{hypoellipticity}, as well.

The last pieces of information we will need to get the consistency of the hydrodynamic limits are the approximate conservation laws, for which we have to go even further in the asymptotic expansions (see Chapter \ref{conservation0-chap}). In the case of one species, this will require to use some suitable renormalization as well as the equi-integrability established in Chapter \ref{hypoellipticity}. In the case of two species, such strong equi-integrability properties are no longer available and we will have to rely on weaker bounds. Thus, in this case, we will only obtain a conditional result, in the sense that the conservation defects and remainders will be controlled by some modulated entropy. We will therefore need, later on (in Chapter \ref{entropy method}), some loop argument based on Gr\"onwall's lemma to prove both the consistency and the convergence in these regimes.

\bigskip

In view of these differences, the convergence proofs will follow different strategies.

\bigskip

In the case of one species, we will use a weak compactness method which relies on some precise study of acoustic and electromagnetic waves (described in Chapter \ref{oscillations}) and compensated compactness.

The core of the proof of convergence for one species will then be the content of Chapter \ref{grad}.

For two species with both weak and strong interactions, we will finally introduce in Chapter \ref{entropy method} a novel renormalized relative entropy method, which will allow to get some stability without any a priori spatial regularity.

%% file: weak0.tex
% ===================
% = First estimates =
% ===================

\chapter{Weak compactness and relaxation estimates}\label{weak bounds}

% \section{Similarities and main differences}
% \section{Fluctuation lemmas}

In this chapter, we establish and recall, from previous works on the hydrodynamic limit of the Boltzmann equation, essential weak compactness estimates on the fluctuations based on the uniform bounds provided by the scaled relative entropy inequalities \eqref{entropy1}, in the case of Theorem \ref{NS-WEAKCV}, and \eqref{entropy2}, in the case of Theorems \ref{CV-OMHD} and \ref{CV-OMHDSTRONG}.

The results presented here are somewhat preliminary to the core of the proofs of Theorems \ref{NS-WEAKCV}, \ref{CV-OMHD} and \ref{CV-OMHDSTRONG}. Thus, they include the first rigorous steps in the proofs of our main theorems and are sometimes straightforward adaptations of lemmas from previous works on the hydrodynamic limit of the Boltzmann equation, while some are new or non-trivial adaptations. In particular, the estimates for two species of particles presented below are all novel.

\bigskip

First, recall that we are considering, in Theorem \ref{NS-WEAKCV}, a sequence of renormalized solutions $\left(f_{\eps},E_\eps,B_\eps\right)$ of the scaled one species Vlasov-Maxwell-Boltzmann system \eqref{VMB1} with initial data $\left(f^{\mathrm{in}}_{\eps},E_\eps^\mathrm{in},B_\eps^\mathrm{in}\right)$ satisfying the uniform bound \eqref{init-fluctuation}, while in Theorems \ref{CV-OMHD} and \ref{CV-OMHDSTRONG}, we consider a sequence of renormalized solutions $\left(f_{\eps}^\pm,E_\eps,B_\eps\right)$ of the scaled two species Vlasov-Maxwell-Boltzmann system \eqref{VMB2} with initial data $\left(f^{\pm\mathrm{in}}_{\eps},E_\eps^\mathrm{in},B_\eps^\mathrm{in}\right)$ satisfying the uniform bound \eqref{init-fluctuation 2}.

We will conveniently employ the notations for fluctuations
\begin{equation*}
	\begin{aligned}
		f_\eps & = M G_\eps = M\left(1+\eps g_\eps\right), &
		f_\eps^\mathrm{in} & = M G_\eps^\mathrm{in} = M\left(1+\eps g_\eps^\mathrm{in}\right), \\
		f_\eps^\pm & = M G_\eps^\pm = M\left(1+\eps g_\eps^\pm\right), &
		f_\eps^{\pm\mathrm{in}} & = M G_\eps^{\pm\mathrm{in}} = M\left(1+\eps g_\eps^{\pm\mathrm{in}}\right),
	\end{aligned}
\end{equation*}
and for scaled collision integrands
\begin{equation*}
	\begin{aligned}
		q_\eps & =\frac{1}{\eps^2}\left(G_{\eps}'G_{\eps *}'-G_{\eps}G_{\eps *}\right), \\
		q_\eps^+ & =\frac{1}{\eps^2}\left({G_{\eps}^+}'{G_{\eps *}^+}'-G_{\eps}^+G_{\eps *}^+\right), \\
		q_\eps^- & =\frac{1}{\eps^2}\left({G_{\eps}^-}'{G_{\eps *}^-}'-G_{\eps}^-G_{\eps *}^-\right), \\
		q_\eps^{+,-} & =\frac{\delta}{\eps^2}\left({G_{\eps}^+}'{G_{\eps *}^-}'-G_{\eps}^+G_{\eps *}^-\right), \\
		q_\eps^{-,+} & =\frac{\delta}{\eps^2}\left({G_{\eps}^-}'{G_{\eps *}^+}'-G_{\eps}^-G_{\eps *}^+\right).
	\end{aligned}
\end{equation*}

% ==========================
% = Relative entropy bound =
% ==========================

\section{Controls from the relative entropy bound}\label{control rel entropy bound}

Following \cite{BGL2}, we introduce the non-negative convex function
\begin{equation*}
		h(z)=(1+z)\log(1+z)-z,
\end{equation*}
defined over $(-1,\infty)$. We may then recast the entropy inequalities \eqref{entropy1} and \eqref{entropy2} utilizing this notation to get the \textbf{relative entropy bounds}, for all $t\geq 0$,
\begin{equation}\label{entropy bound}
	\frac1{\eps^2} H\left(f_\eps\right)
	=
	\frac1{\eps^2} H\left(f_\eps|M\right)
	=
	\int_{\mathbb{R}^3\times\mathbb{R}^3}\frac 1{\eps^2}h\left(\eps g_\eps\right)Mdxdv
	\leq C^\mathrm{in},
\end{equation}
and
\begin{equation}\label{entropy bound 2}
	\frac1{\eps^2} H\left(f_\eps^\pm\right)
	=
	\frac1{\eps^2} H\left(f_\eps^\pm|M\right)
	=
	\int_{\mathbb{R}^3\times\mathbb{R}^3}\frac 1{\eps^2}h\left(\eps g_\eps^\pm\right)Mdxdv
	\leq C^\mathrm{in}.
\end{equation}

The relative entropy bounds are expected to control the size of the fluctuations $g_\eps$ and $g_\eps^\pm$ since
\begin{equation*}
	h(z) \sim\frac12 z^2, \quad \text{near } z = 0.
\end{equation*}
However, this behavior only holds asymptotically, as $z\rightarrow 0$, and thus, in order to exploit the relative entropy bounds, we will have to rely crucially on Young's inequality \eqref{Young h} for $h(z)$, presented in Appendix \ref{inequalities}.

The following lemma is a mere reformulation of Proposition 3.1 from \cite{BGL2}. It is a consequence solely of the fact that the fluctuations satisfy the entropy bounds \eqref{entropy bound} and \eqref{entropy bound 2}.

\begin{lem}\label{L1-lem} 
	Let $f_\eps(t,x,v)$ be a family of measurable, almost everywhere non-negative distribution functions such that, for all $t\geq 0$,
	\begin{equation*}
		\frac1{\eps^2} H\left(f_\eps\right)(t)
		=
		\int_{\mathbb{R}^3\times\mathbb{R}^3}\frac 1{\eps^2}h\left(\eps g_\eps\right)(t)Mdxdv
		\leq C^\mathrm{in}.
	\end{equation*}
	
	Then, as $\eps\rightarrow 0$~:
	\begin{enumerate}
		\item any subsequence of fluctuations $g_\eps$ is uniformly bounded in
		\begin{equation*}
			L^\infty\left(dt;L^1_\mathrm{loc}\left(dx;L^1\left(\left(1+\left|v\right|^2\right)Mdv\right)\right)\right),
		\end{equation*}
		and weakly relatively compact in
		\begin{equation*}
			L^1_\mathrm{loc}\left(dtdx ; L^1\left(\left(1+\left|v\right|^2\right)M dv\right)\right).
		\end{equation*}
		\item if $g$ is a weak limit point in $L^1_\mathrm{loc}\left(dtdx ; L^1\left(\left(1+\left|v\right|^2\right)M dv\right)\right)$ of the family of fluctuations $g_\eps$, then $g$ belongs to $L^\infty\left(dt;L^2\left(Mdxdv\right)\right)$ and satisfies, for almost every $t\geq 0$,
		\begin{equation*}
			\frac{1}{2}\int_{\mathbb{R}^3\times\mathbb{R}^3}g(t)^2 Mdxdv
			% \leq
			% \liminf_{\eps\rightarrow 0}\frac{1}{\eps^2}H(f_\eps)(t)
			\leq C^\mathrm{in}.
		\end{equation*}
	\end{enumerate}
\end{lem}

\begin{proof}
	For the sake of completeness, we recall the main ideas from the proof of Proposition 3.1 in \cite{BGL2}, which is based on an application of inequality \eqref{Young h}. Thus, setting $y=\frac{1}{4}\left(1+|v|^2\right)$, $z= g_\eps$, $\beta=4$ and $\alpha\geq 4\eps$ in \eqref{Young h} yields, almost everywhere in $(t,x,v)$,
	\begin{equation*}
		\left(1+|v|^2\right)\left|g_\eps\right|\leq
		\frac{\alpha}{\eps^2}h\left(\eps g_\eps\right)
		+
		\frac{16e^{\frac{1}{4}}}{\alpha} e^{\frac{|v|^2}{4}}.
	\end{equation*}
	This is then integrated in all variables on suitable sets to demonstrate, with the entropy bound, the equi-integrability and tightness of the sequences, and thus their weak compactness.

	We set $\alpha=1$ first. Then, for each measurable set $E\subset\mathbb{R}^3$ of finite measure, it holds that, for every $0<\eps\leq\frac 14$, i.e.\ for all but a finite number of $\eps$'s,
	\begin{equation*}
		\int_{E\times\mathbb{R}^3} (1+|v|^2)\left|g_{\eps}(t)\right|Mdxdv
		\leq C^\mathrm{in} + 16e^\frac{1}{4}|E| \int_{\mathbb{R}^3} e^{\frac{|v|^2}{4}} Mdv.
	\end{equation*}
	Hence, the family $\left(1+|v|^2\right)g_{\eps}$ is uniformly bounded in $L^\infty\left(dt ; L^1_\mathrm{loc}\left(dx;L^1\left(Mdv\right)\right)\right)$.

	Similarly, for arbitrary $\alpha\geq 4\eps$ and for any measurable set $E\subset [0,T]\times\mathbb{R}^3\times\mathbb{R}^3$, where $T>0$, one has that
	\begin{equation*}
		\int_{E} (1+|v|^2)\left|g_{\eps}\right| M dtdxdv
		\leq \alpha T C^\mathrm{in} + \frac{16e^\frac{1}{4}}{\alpha}|E|.
	\end{equation*}
	This shows, by the arbitrariness of $\alpha\geq 4 \eps$, that the family $(1+|v|^2)g_{\eps}M$ is uniformly integrable on $[0,T]\times\mathbb{R}^3\times\mathbb{R}^3$.
	
	% Therefore, setting $\alpha=4\left|E\right|^\frac 12$, we deduce that, for every $0<\eps\leq {\left|E\right|^\frac 12}$, i.e.\ for all but a finite number of $\eps$'s,
	% \begin{equation*}
	% 	\int_{E} (1+|v|^2)\left|g_{\eps}\right| M dtdxdv
	% 	\leq 4\left(T C^\mathrm{in} + e^\frac{1}{4}\right)|E|^\frac{1}{2}.
	% \end{equation*}

	Finally, for arbitrary $\alpha\geq 4\eps$, any time $T>0$, any compact set $K\subset \mathbb{R}^3$ and any large radius $R>0$, we find that
	\begin{equation*}
		\int_{[0,T]\times K \times \left\{|v|\geq R\right\}}\left(1+|v|^2\right)\left|g_\eps\right| M dtdxdv \leq
		\alpha T C^\mathrm{in}
		+
		\frac{16e^{\frac{1}{4}}}{\alpha}T|K| \int_{\left\{|v|\geq R\right\}}e^{\frac{|v|^2}{4}}M dv,
	\end{equation*}
	which, by the arbitrariness of $\alpha\geq 4\eps$, clearly implies the tightness in velocity of the family $(1+|v|^2)g_{\eps}M$.

	On the whole, by virtue of the Dunford-Pettis criterion \cite{royden}, we infer the weak relative compactness of the family $(1+|v|^2)g_{\eps}M$ in $L^1_\mathrm{loc}\left(dtdx;L^1\left(dv\right)\right)$, which concludes the demonstration of the first assertion of the lemma.

	The second assertion will follow from a convexity analysis of the relative entropy functional $H(f_\eps)$. Indeed, by convexity of $h(z)$, it holds that
	\begin{equation*}
		\frac 1 {\eps^2}h(\eps g) + \frac 1\eps h'(\eps g)\left(g_\eps - g \right) \leq \frac 1{\eps^2}h(\eps g_\eps).
	\end{equation*}
	Hence, for any large $\lambda>0$, any times $0\leq t_1 < t_2$ and any compact set $K\subset \mathbb{R}^3$, by the non-negativity of $h(z)$, we have that, for $\eps<\frac 1\lambda$ (so that $\eps g >-1$),
	\begin{equation*}
		\begin{aligned}
			\int_{t_1}^{t_2}\int_{K\times\mathbb{R}^3}
			\left( \frac 1 {\eps^2}h(\eps g) + \frac 1\eps h'(\eps g)\left(g_\eps - g \right) \right)
			& \mathds{1}_{\left\{|g|\leq \lambda\right\}} M
			dxdvdt
			\\
			\leq
			& \int_{t_1}^{t_2}\int_{\mathbb{R}^3\times\mathbb{R}^3}
			\frac 1{\eps^2}h(\eps g_\eps) M
			dxdvdt.
		\end{aligned}
	\end{equation*}
	Furthermore, notice that $h(z)=\frac 12 z^2 - \int_0^z\frac{1}{(1+y)^2}\frac{(y-z)^2}{2}dy$ and $h'(z)=z-\int_0^z\frac{1}{(1+y)^2}(y-z)dy$, from which we easily deduce the strong convergences
	\begin{equation*}
		\frac{1}{\eps^2}h(\eps g)\mathds{1}_{\left\{|g|\leq \lambda\right\}}\rightarrow
		\frac 12 g^2\mathds{1}_{\left\{|g|\leq \lambda\right\}}
		\text{ and }
		\frac{1}{\eps}h'(\eps g)\mathds{1}_{\left\{|g|\leq \lambda\right\}}\rightarrow
		g\mathds{1}_{\left\{|g|\leq \lambda\right\}}
		\text{ in }L^\infty\left(dtdxdv\right).
	\end{equation*}
	Therefore, taking weak limits in the above convexity inequality yields
	\begin{equation*}
		\begin{aligned}
			\int_{t_1}^{t_2}\int_{K\times\mathbb{R}^3}
			\frac 1 {2}g^2
			\mathds{1}_{\left\{|g|\leq \lambda\right\}} M
			dxdvdt
			& \leq
			\liminf_{\eps\rightarrow 0}\int_{t_1}^{t_2}\int_{\mathbb{R}^3\times\mathbb{R}^3}
			\frac 1{\eps^2}h(\eps g_\eps) M
			dxdvdt \\
			& \leq C^\mathrm{in}(t_2-t_1),
		\end{aligned}
	\end{equation*}
	which, by monotonicity of the integrands, gives
	\begin{equation*}
		\int_{t_1}^{t_2}\int_{\mathbb{R}^3\times\mathbb{R}^3}
		\frac 1 {2}g^2 M
		dxdvdt\leq C^\mathrm{in}(t_2-t_1).
	\end{equation*}
	Finally, the proof of the lemma is concluded by the arbitrariness of $t_1$ and $t_2$.
\end{proof}

The second assertion of the preceding lemma shows that, in the vanishing $\eps$ limit, the limiting fluctuation belongs to $L^2(Mdxdv)$ uniformly in $t$. Hence, the weighted $L^1$-bound implied by the first assertion of Lemma \ref{L1-lem} is certainly not optimal.

Thus, in order to refine our understanding of the limit $\eps\rightarrow 0$, we consider the following renormalized fluctuations
\begin{equation*}
	\begin{aligned}
		\sqrt{G_\eps} & = 1+\frac\eps 2 \hat g_\eps, &
		\sqrt{G_\eps^\mathrm{in}} & = 1+\frac\eps 2 \hat g_\eps^\mathrm{in}, \\
		\sqrt{G_\eps^\pm} & = 1+\frac\eps 2 \hat g_\eps^\pm, &
		\sqrt{G_\eps^{\pm\mathrm{in}}} & = 1+\frac\eps 2 \hat g_\eps^{\pm\mathrm{in}},
	\end{aligned}
\end{equation*}
or, equivalently,
\begin{equation}\label{hatg}
	\begin{aligned}
		\hat g_\eps & = \frac2{\eps}\left(\sqrt{G_\eps}-1\right),
		&
		\hat g_\eps^\mathrm{in} & = \frac2{\eps}\left(\sqrt{G_\eps^\mathrm{in}}-1\right),
		\\
		\hat g_\eps^\pm & = \frac2{\eps}\left(\sqrt{G_\eps^\pm}-1\right),
		&
		\hat g_\eps^{\mathrm{in}\pm} & = \frac2{\eps}\left(\sqrt{G_\eps^{\mathrm{in}\pm}}-1\right).
	\end{aligned}
\end{equation}
Such square root renormalizations have already been used in previous works on hydrodynamic limits. The advantages of these renormalized fluctuations over the original ones become apparent in the coming lemma, which is, essentially, a modern reformulation of Corollary 3.2 from \cite{BGL2}.

\begin{lem}\label{L2-lem} 
	Let $f_\eps(t,x,v)$ be a family of measurable, almost everywhere non-negative distribution functions such that, for all $t\geq 0$,
	\begin{equation*}
		\frac1{\eps^2} H\left(f_\eps\right)(t)
		=
		\int_{\mathbb{R}^3\times\mathbb{R}^3}\frac 1{\eps^2}h\left(\eps g_\eps\right)(t)Mdxdv
		\leq C^\mathrm{in}.
	\end{equation*}
	
	Then, as $\eps\rightarrow 0$, any subsequence of renormalized fluctuations $\hat g_\eps$ is uniformly bounded in $L^\infty\left(dt;L^2\left(Mdxdv\right)\right)$.
\end{lem}

\begin{proof}
	The elementary inequality \eqref{sqrt entropy} implies that, for all $t\geq 0$,
	\begin{equation}\label{g-L2}
		\begin{aligned}
			\int_{\mathbb{R}^3\times\mathbb{R}^3}{\hat g_\eps}^2(t)Mdxdv
			& =
			\int_{\mathbb{R}^3\times\mathbb{R}^3}\frac 4{\eps^2}\left(\sqrt{1+\eps g_\eps}-1\right)^2(t)Mdxdv \\
			& \leq
			\frac 4{\eps^2}H(f_\eps)(t) \leq 4C^\mathrm{in},
		\end{aligned}
	\end{equation}
	which is the announced result.
\end{proof}

The simple Lemma \ref{L2-lem} provides important information on any subsequence of fluctuations $g_\eps$. Indeed, a very natural application of this refined a priori estimate follows from decomposing the fluctuations as
\begin{equation}\label{fluct-decomposition}
	g_\eps = \hat g_\eps + \frac \eps 4 \hat g_\eps^2
	\in
	L^\infty\left(dt;L^2\left(Mdxdv\right)\right)
	+ \eps
	L^\infty\left(dt;L^1\left(Mdxdv\right)\right).
\end{equation}
Therefore, we see from Lemma \ref{L2-lem}, that the fluctuations $g_\eps$ are uniformly bounded in $L^\infty\left(dt;L^2\left(Mdxdv\right)\right)$, up to a remainder of order $\eps$ in $L^\infty\left(dt;L^1\left(Mdxdv\right)\right)$. In particular, according to Lemma \ref{L1-lem}, if $g$ is a weak limit point in $L^1_\mathrm{loc}\left(dtdx ; L^1\left(\left(1+\left|v\right|^2\right)M dv\right)\right)$ of a converging subsequence of fluctuations $g_\eps$, then $\hat g_\eps$ also converges towards $g$ in the weak-$*$ topology of $L^\infty\left(dt;L^2\left(Mdxdv\right)\right)$.

As we will see later on, it will be crucial to establish sharper properties of tightness and equi-integrability on the sequence of integrable functions $\hat g_\eps^2$. These refinements will follow from the joint control of the fluctuations by the entropy and the entropy dissipation bounds.

% =============================
% = Entropy dissipation bound =
% =============================

\section{Controls from the entropy dissipation bound}

Following \cite{BGL2}, again, we introduce the non-negative convex function
\begin{equation*}
	r(z)=z\log(1+z),
\end{equation*}
defined over $(-1,\infty)$. We may then recast the entropy inequalities \eqref{entropy1} and \eqref{entropy2} utilizing this notation to get the \textbf{entropy dissipation bounds} (here, exceptionally, we consider any cross-section $b\geq 0$, for more generality), for all $t\geq 0$,
\begin{equation*}
	\begin{aligned}
		\frac{1}{\epsilon^4}\int_0^t\int_{\mathbb{R}^3} & D(f_\eps)(s) dx ds
		\\
		& =
		\int_0^t \int_{\mathbb{R}^3}\int_{\mathbb{R}^3\times\mathbb{R}^3\times\mathbb{S}^2}
		\frac 1{4\eps^4}r\left(\frac{\eps^2 q_\eps}{G_\eps G_{\eps *}}\right)f_\eps f_{\eps *}bdvdv_* d\sigma dx ds
		\leq C^\mathrm{in},
	\end{aligned}
\end{equation*}
and
\begin{equation*}
	\begin{aligned}
		\frac{1}{\epsilon^4}\int_0^t\int_{\mathbb{R}^3} & \left(D\left(f_\eps^+\right)+D\left(f_\eps^-\right)
		+ \delta^2 D\left(f_\eps^+,f_\eps^-\right)\right)(s) dx ds
		\\
		& =
		\int_0^t \int_{\mathbb{R}^3}\int_{\mathbb{R}^3\times\mathbb{R}^3\times\mathbb{S}^2}
		\frac 1{4\eps^4}r\left(\frac{\eps^2 q_\eps^+}{G_\eps^+ G_{\eps *}^+}\right)f_\eps^+ f_{\eps *}^+ bdvdv_* d\sigma dx ds
		\\
		& +
		\int_0^t \int_{\mathbb{R}^3}\int_{\mathbb{R}^3\times\mathbb{R}^3\times\mathbb{S}^2}
		\frac 1{4\eps^4}r\left(\frac{\eps^2 q_\eps^-}{G_\eps^- G_{\eps *}^-}\right)f_\eps^- f_{\eps *}^- bdvdv_* d\sigma dx ds
		\\
		& +
		\int_0^t \int_{\mathbb{R}^3}\int_{\mathbb{R}^3\times\mathbb{R}^3\times\mathbb{S}^2}
		\frac {\delta^2}{2\eps^4}r\left(\frac{\eps^2 q_\eps^{+,-}}{\delta G_\eps^+ G_{\eps *}^-}\right)f_\eps^+f_{\eps *}^- bdvdv_* d\sigma dx ds
		\leq C^\mathrm{in}.
	\end{aligned}
\end{equation*}

The entropy dissipation bounds are expected to control the size of the collision integrands $q_\eps$, $q_\eps^+$, $q_\eps^-$, $q_\eps^\pm$ and $q_\eps^\mp$ since
\begin{equation*}
	r(z) \sim z^2, \quad \text{near } z = 0.
\end{equation*}
However, this behavior only holds asymptotically, as $z\rightarrow 0$, and thus, in order to exploit the entropy dissipation bounds, we will have to rely crucially on Young's inequality \eqref{Young r} for $r(z)$ and on inequality \eqref{sqrt dissip}, presented in Appendix \ref{inequalities}. Furthermore, when coupled with a coercivity estimate for some suitable non-singular linearized collision operator, the entropy dissipation bounds will actually provide some control on the relaxation to equilibrium of the fluctuations $g_\eps$ and $g_\eps^\pm$ (see Section \ref{relaxation} below).

In order to refine our understanding of the limit $\eps\rightarrow 0$, we consider the following renormalized collision integrands
\begin{equation}\label{hatq-def}
	\begin{aligned}
		\hat q_\eps & =\frac{2}{\eps^2}\left(\sqrt{G_{\eps}'G_{\eps *}'}-\sqrt{G_{\eps}G_{\eps *}}\right), \\
		\hat q_\eps^+ & =\frac{2}{\eps^2}\left(\sqrt{{G_{\eps}^{+\prime}}{G_{\eps *}^{+\prime}}}-\sqrt{G_{\eps}^+G_{\eps *}^+}\right), \\
		\hat q_\eps^- & =\frac{2}{\eps^2}\left(\sqrt{{G_{\eps}^{-\prime}}{G_{\eps *}^{-\prime}}}-\sqrt{G_{\eps}^-G_{\eps *}^-}\right), \\
		\hat q_\eps^{+,-} & =\frac{2\delta}{\eps^2}\left(\sqrt{{G_{\eps}^{+\prime}}{G_{\eps *}^{-\prime}}}-\sqrt{G_{\eps}^+G_{\eps *}^-}\right), \\
		\hat q_\eps^{-,+} & =\frac{2\delta}{\eps^2}\left(\sqrt{{G_{\eps}^{-\prime}}{G_{\eps *}^{+\prime}}}-\sqrt{G_{\eps}^-G_{\eps *}^+}\right).
	\end{aligned}
\end{equation}
The advantages of these renormalized collision integrands over the original ones become apparent in the coming lemma.

\begin{lem}\label{L2-qlem}
	Let $f_\eps^+(t,x,v)$ and $f_\eps^-(t,x,v)$ be two families of measurable, almost everywhere non-negative distribution functions such that, for all $t\geq 0$,
	\begin{equation*}
		\begin{aligned}
			\frac{\delta^2}{\epsilon^4}\int_0^t\int_{\mathbb{R}^3} &
			D\left(f_\eps^+,f_\eps^-\right)(s) dx ds
			\\
			& =
			\int_0^t \int_{\mathbb{R}^3}\int_{\mathbb{R}^3\times\mathbb{R}^3\times\mathbb{S}^2}
			\frac {\delta^2}{2\eps^4}r\left(\frac{\eps^2 q_\eps^{+,-}}{\delta G_\eps^+ G_{\eps *}^-}\right)f_\eps^+f_{\eps *}^- bdvdv_* d\sigma dx ds
			\leq C^\mathrm{in}.
		\end{aligned}
	\end{equation*}
	
	Then, as $\eps\rightarrow 0$, any subsequence of renormalized collision integrands $\hat q_\eps^{+,-}$ is uniformly bounded in $L^2\left(bMM_*dtdxdvdv_*d\sigma\right)$.
\end{lem}

\begin{proof}
	The elementary inequality \eqref{sqrt dissip} implies that, for all $t\geq 0$,
	\begin{equation}\label{q-est}
		\begin{aligned}
			\frac 12 \int_0^t \int_{\mathbb{R}^3} & \int_{\mathbb{R}^3\times\mathbb{R}^3\times\mathbb{S}^2}
			\left(\hat{q}_\eps^{+,-}\right)^2
			bMM_*dvdv_* d\sigma dx ds
			\\
			& = 2 \int_0^t \int_{\mathbb{R}^3}\int_{\mathbb{R}^3\times\mathbb{R}^3\times\mathbb{S}^2}
			\frac {\delta^2}{\eps^4} \left(\sqrt{1+\frac{\eps^2 q_\eps^{+,-}}{\delta G_\eps^+ G_{\eps *}^-}}-1\right)^2f_\eps^+f_{\eps *}^- bdvdv_* d\sigma dx ds
			\\
			& \leq \frac{\delta^2}{\epsilon^4}\int_0^t\int_{\mathbb{R}^3}
			D\left(f_\eps^+,f_\eps^-\right)(s) dx ds
			\leq C^\mathrm{in},
		\end{aligned}
	\end{equation}
	which is the announced result.
\end{proof}

The simple Lemma \ref{L2-qlem} provides important information on any subsequence of collision integrands $q^{+,-}_\eps$. Indeed, a very natural application of this refined a priori estimate follows from decomposing the collision integrands as
\begin{equation}\label{integrands-decomposition}
	q^{+,-}_\eps = \sqrt{G_\eps^+ G_{\eps *}^-} \hat q_\eps^{+,-} + \frac{\eps^2}{4\delta}\left(\hat q_\eps^{+,-}\right)^2.
\end{equation}
Therefore, at least in the simpler case of the Maxwellian cross-section $b\equiv 1$, we see from Lemma \ref{L2-qlem} that, for any admissible renormalization $\beta(z)$, the renormalized collision integrands $\beta'\left(G_\eps^+\right)q^{+,-}_\eps$ are uniformly bounded in $L^1_\mathrm{loc}\left(dtdx;L^1\left(MM_*dvdv_*d\sigma\right)\right)$, provided the natural entropy and entropy dissipation bounds are satisfied. Moreover, employing Lemma \ref{L2-lem} and Egorov's theorem, it is possible to show that
\begin{equation*}
	\beta'\left(G_\eps^+\right)\sqrt{G_\eps^+G_{\eps*}^-}
	\rightarrow \beta'(1)
	\qquad \text{in }
	L^2_\mathrm{loc}\left(dtdx;L^2\left(MM_*dvdv_*d\sigma\right)\right).
\end{equation*}
In particular, if $q^{+,-}$ is a weak limit point in $L^2\left(MM_*dtdxdvdv_*d\sigma\right)$ of a converging subsequence of collision integrands $\hat q_\eps^{+,-}$, then $\beta'\left(G_\eps^+\right)q^{+,-}_\eps$ converges --~up to extraction of jointly converging subsequences~-- towards $\beta'(1)q^{+,-}$ in the weak topology of $L^1_\mathrm{loc}\left(dtdx;L^1\left(MM_*dvdv_*d\sigma\right)\right)$.

% ==============
% = relaxation =
% ==============

\section{Relaxation towards thermodynamic equilibrium}\label{relaxation}

In this section, we establish the relaxation of fluctuations towards thermodynamic equilibrium as a consequence of the relative entropy and the entropy dissipation bounds. As we consider fluctuations around a global equilibrium, the linearized collision operator $\cL$, defined in \eqref{def L and Q} and \eqref{def L and Q two species}, is expected to play here a fundamental role.

We recall that, for the sake of simplicity, we restrict our attention, in this part of our work, to the case of Maxwellian molecules, that is to constant collision cross-sections, say $b \equiv 1$. Nevertheless, up to additional technical difficulties, the results in this section will be extended to general cross-sections in the remaining parts of the present work.

\bigskip

The spectral analysis of the linearized collision operator uses crucially the following decomposition based on a clever change of variables sometimes called ``Carleman's collision parametrization'', although it goes back to Hilbert \cite{hilbert} (see equation (17) in \cite{hilbert} and the computations therein). We refer to \cite{SR} for more details and to \cite{levermore2} for a modern and general treatment of the linearized Boltzmann operator.

\begin{prop}[Hilbert's decomposition of $\cL$]\label{hilbert-prop}
	The linearized collision operator $\cL$ defined by
	\begin{equation*}% \label{L-def}
		\cL g =-\frac1M \left(Q(Mg,M)+Q(M,Mg)\right) = \int_{\mathbb{R}^3\times\mathbb{S}^2} \left(g+g_* - g'-g'_*\right)M_* dv_*d\sigma,
	\end{equation*}
	can be decomposed as
	\begin{equation*}
		\cL g =g-\cK g,
	\end{equation*}
	where $\cK$ is a compact integral operator on $L^2(Mdv)$.
\end{prop}

As an immediate consequence of the preceding proposition, the operator $\cL$ satisfies the Fredholm alternative, as well as some coercivity estimate, which will be used to control the relaxation process. We refer to \cite{levermore2} or \cite{SR} for details and justifications of the following proposition, or to the proof of the more general Proposition \ref{coercivity 2} below.

\begin{prop}[Coercivity of $\cL$]\label{coercivity}
	The linear collision operator $\cL$ is a non-negative self-adjoint operator on $L^2(Mdv)$ with nullspace
	\begin{equation*}
		\Ker(\cL)=\Span\left\{1,v_1,v_2,v_3,|v|^2\right\}.
	\end{equation*}
	Moreover, the following coercivity estimate holds~: there exists $C>0$ such that, for each $g\in \Ker(\cL)^\perp\subset L^2(Mdv)$,
	\begin{equation*}
		\|g\|_{L^2(M dv)}^2
		\leq C
		\int_{\mathbb{R}^3} g\cL g(v)M(v)dv.
	\end{equation*}
	In particular, for any $g\in\Ker(\cL)^\perp\subset L^2(Mdv)$,
	\begin{equation*}
		\|g\|_{L^2(M dv)}
		\leq C
		\left\|\cL g\right\|_{L^2(Mdv)}.
	\end{equation*}
\end{prop}

We will also need the generalization of the preceding propositions to the linearized collision operator for two species of particles $\L$. In fact, employing the results from \cite{levermore2}, we easily obtain the following Hilbert's decomposition for $\L$.

\begin{prop}[Hilbert's decomposition of $\L$]\label{hilbert-prop 2}
	The linearized collision operator for two species $\L$ defined by
	\begin{equation*}
		\L
		\begin{pmatrix} g \\ h \end{pmatrix}
		=
		\begin{pmatrix}
			\cL g + \cL \left(g,h\right) \\
			\cL h + \cL \left(h,g\right)
		\end{pmatrix},
		% =
		% \begin{pmatrix}
		% 	\cL g + \int_{\mathbb{R}^3\times\mathbb{S}^2} \left(g+h_*-g'-h'_*\right) b M_*dv_*d\sigma\\
		% 	\cL h + \int_{\mathbb{R}^3\times\mathbb{S}^2} \left(h+g_*-h'-g_*'\right)  b M_*dv_*d\sigma
		% \end{pmatrix}
	\end{equation*}
	where
	\begin{equation*}
		\cL \left(g,h\right) =-\frac1M \left(Q(Mg,M)+Q(M,Mh)\right) = \int_{\mathbb{R}^3\times\mathbb{S}^2} \left(g+h_* - g'-h'_*\right)M_* dv_*d\sigma,
	\end{equation*}
	can be decomposed as
	\begin{equation*}
		\L
		\begin{pmatrix} g \\ h \end{pmatrix}
		=2
		\begin{pmatrix} g \\ h \end{pmatrix}
		-\K
		\begin{pmatrix} g \\ h \end{pmatrix},
	\end{equation*}
	where $\K$ is a compact integral operator on $L^2(Mdv)$.
\end{prop}

% Note that the definition of $\L$ above slightly differs from the definitions \eqref{vector L 2} and \eqref{vector L} of the linearized operators. However, for the sake of simplicity, our notation will not distinguish these definitions and we will, from now on, stick to the definition of $\L$ from Proposition \ref{hilbert-prop 2}.

As an immediate consequence of the preceding proposition, the operator $\L$ satisfies the Fredholm alternative, as well as some coercivity estimate, which will be used to control the relaxation process for two species of particles. For the sake of completeness, we provide here a brief justification of the following proposition.

\begin{prop}[Coercivity of $\L$]\label{coercivity 2}
	The linear collision operator $\L$ is a non-negative self-adjoint operator on $L^2(Mdv)$ with nullspace
	\begin{equation*}
		\Ker(\L)=\Span
		\left\{
		\begin{pmatrix}
			1\\0
		\end{pmatrix},
		\begin{pmatrix}
			0\\1
		\end{pmatrix},
		\begin{pmatrix}
			v_1\\v_1
		\end{pmatrix},
		\begin{pmatrix}
			v_2\\v_2
		\end{pmatrix},
		\begin{pmatrix}
			v_3\\v_3
		\end{pmatrix},
		\begin{pmatrix}
			|v|^2\\|v|^2
		\end{pmatrix}
		\right\}.
	\end{equation*}
	Moreover, the following coercivity estimate holds~: there exists $C>0$ such that, for each $\begin{pmatrix} g \\ h \end{pmatrix} \in \Ker(\L)^\perp\subset L^2(Mdv)$,
	\begin{equation*}
		\left\|\begin{pmatrix} g \\ h \end{pmatrix}\right\|_{L^2(M dv)}^2
		\leq C
		\int_{\mathbb{R}^3} \begin{pmatrix} g \\ h \end{pmatrix}\cdot
		\L \begin{pmatrix} g \\ h \end{pmatrix}(v)
		M(v)dv.
	\end{equation*}
	In particular, for any $\begin{pmatrix} g \\ h \end{pmatrix} \in \Ker(\L)^\perp\subset L^2(Mdv)$,
	\begin{equation*}
		\left\|\begin{pmatrix} g \\ h \end{pmatrix}\right\|_{L^2(M dv)}
		\leq C
		\left\|\L\begin{pmatrix} g \\ h \end{pmatrix}\right\|_{L^2(M dv)}.
	\end{equation*}
\end{prop}

\begin{proof}
	The non-negativity and the self-adjointness of $\L$ easily follow from a standard use of the collision symmetries by showing that
	\begin{equation}\label{vector symmetries}
		\begin{aligned}
			\int_{\mathbb{R}^3} & \begin{pmatrix} g \\ h \end{pmatrix}\cdot
			\L \begin{pmatrix} \bar g \\ \bar h \end{pmatrix}(v)
			M(v)dv \\
			& =
			\frac{1}{4}
			\int_{\mathbb{R}^3\times\mathbb{R}^3\times\mathbb{S}^2}
			\left(g+g_*-g'-g_*'\right)\left(\bar g+\bar g_*-\bar g'-\bar g_*'\right) MM_* dvdv_*d\sigma \\
			& +
			\frac{1}{4}
			\int_{\mathbb{R}^3\times\mathbb{R}^3\times\mathbb{S}^2}
			\left(h+h_*-h'-h_*'\right)\left(\bar h+\bar h_*-\bar h'-\bar h_*'\right) MM_* dvdv_*d\sigma \\
			& +
			\frac{1}{2}
			\int_{\mathbb{R}^3\times\mathbb{R}^3\times\mathbb{S}^2}
			\left(g+h_*-g'-h_*'\right)\left(\bar g+\bar h_*-\bar g'-\bar h_*'\right) MM_* dvdv_*d\sigma.
		\end{aligned}
	\end{equation}

	Next, consider $\begin{pmatrix} g \\ h \end{pmatrix} \in \Ker(\L)$. We deduce from \eqref{vector symmetries} that, necessarily, $g=\Pi g$, $h=\Pi h$ and
	\begin{equation*}
		\int_{\mathbb{R}^3\times\mathbb{R}^3\times\mathbb{S}^2}
		\left(\Pi \left(g-h\right)-\left(\Pi \left(g-h\right)\right)'\right)^2  MM_* dvdv_*d\sigma =0.
	\end{equation*}
	A simple and direct computation shows then that $h$ and $g$ have the same bulk velocity and temperature, which completes the characterization of the kernel of $\L$. In particular, the orthogonal projection onto the kernel of $\L$ in $L^2(Mdv)$ is explicitly given by
	\begin{equation*}
		\mathbb{P}
		\begin{pmatrix}
			g \\ h
		\end{pmatrix}
		=
		\begin{pmatrix}
			\frac 12\int_{\mathbb{R}^3}(g-h)Mdv + \Pi\frac{g+h}{2} \\ -\frac 12\int_{\mathbb{R}^3}(g-h)Mdv + \Pi\frac{g + h}{2}
		\end{pmatrix}.
	\end{equation*}

	Finally, since $\L$ is positive definite, self-adjoint and satisfies Hilbert's decomposition from Proposition \ref{hilbert-prop 2}, we easily obtain, by the spectral theorem for compact self-adjoint operators, writing $\begin{pmatrix} g \\ h \end{pmatrix}$ in the Hilbert basis of eigenvectors of $\L$, that
	\begin{equation*}
		\left\| \begin{pmatrix} g\\h \end{pmatrix} -\P \begin{pmatrix} g\\h \end{pmatrix}\right\|_{L^2( Mdv)}^2
		\leq C
		\int_{\mathbb{R}^3} \begin{pmatrix} g\\h \end{pmatrix} \cdot\L \begin{pmatrix} g\\h \end{pmatrix} Mdv,
	\end{equation*}
	which concludes the justification of the proposition.
\end{proof}

Finally, we also extend the preceding propositions to the linearized collision operator $\mathfrak{L}$. In fact, employing the results from \cite{levermore2}, we easily obtain the following Hilbert's decomposition for $\mathfrak{L}$.

\begin{prop}[Hilbert's decomposition of $\mathfrak{L}$]\label{hilbert-prop 3}
	The linearized collision operator $\mathfrak{L}$ defined by
	\begin{equation*}
		\begin{aligned}
			\mathfrak{L} g = \mathcal{L}(g,-g) & =-\frac1M \left(Q(Mg,M)-Q(M,Mg)\right)
			\\
			& = \int_{\mathbb{R}^3\times\mathbb{S}^2} \left(g-g_* - g'+g'_*\right)M_* dv_*d\sigma,
		\end{aligned}
	\end{equation*}
	can be decomposed as
	\begin{equation*}
		\mathfrak{L} g =g-\mathfrak{K} g,
	\end{equation*}
	where $\mathfrak{K}$ is a compact integral operator on $L^2(Mdv)$.
\end{prop}

Note that the definition of $\mathfrak{L}$ above coincides with \eqref{L frak def}.

As an immediate consequence of the preceding proposition, the operator $\mathfrak{L}$ satisfies the Fredholm alternative, as well as some coercivity estimate, which will be used to control the relaxation process for two species of particles. For the sake of completeness, we provide here a brief justification of the following proposition.

\begin{prop}[Coercivity of $\mathfrak{L}$]\label{coercivity 3}
	The linear collision operator $\mathfrak{L}$ is a non-negative self-adjoint operator on $L^2(Mdv)$ with nullspace
	\begin{equation*}
		\Ker(\mathfrak{L})=\Span\left\{1\right\}.
	\end{equation*}
	Moreover, the following coercivity estimate holds~: there exists $C>0$ such that, for each $g\in \Ker(\mathfrak{L})^\perp\subset L^2(Mdv)$,
	\begin{equation*}
		\|g\|_{L^2(M dv)}^2
		\leq C
		\int_{\mathbb{R}^3} g\mathfrak{L} g(v)M(v)dv.
	\end{equation*}
	In particular, for any $g\in\Ker(\mathfrak{L})^\perp\subset L^2(Mdv)$,
	\begin{equation*}
		\|g\|_{L^2(M dv)}
		\leq C
		\left\|\mathfrak{L} g\right\|_{L^2(Mdv)}.
	\end{equation*}
\end{prop}

\begin{proof}
	The non-negativity and the self-adjointness of $\mathfrak{L}$ easily follow from a standard use of the collision symmetries by showing that
	\begin{equation}\label{frak symmetries}
		\begin{aligned}
			\int_{\mathbb{R}^3} & g
			\mathfrak{L} h(v)
			M(v)dv \\
			& =
			\frac{1}{4}
			\int_{\mathbb{R}^3\times\mathbb{R}^3\times\mathbb{S}^2}
			\left(g-g_*-g'+g_*'\right)\left(h-h_*-h'+h_*'\right) MM_* dvdv_*d\sigma.
		\end{aligned}
	\end{equation}

	Next, consider $g \in \Ker(\mathfrak{L})$. We deduce from \eqref{frak symmetries} that, necessarily,
	\begin{equation*}
		g-g_*=g'-g_*',
	\end{equation*}
	almost everywhere. Hence, since the change of variable $\sigma\mapsto -\sigma$ merely exchanges $v'$ and $v_*'$, we find, averaging over $\sigma\in\mathbb{S}^2$,
	\begin{equation*}
		g-g_*=\frac{1}{\left|\mathbb{S}^2\right|}\int_{\mathbb{S}^2}\left(g'-g_*'\right)d\sigma = 0,
	\end{equation*}
	for every $v,v_*\in\mathbb{R}^3$. It follows that $g$ is a constant function.

	Finally, since $\mathfrak{L}$ is positive definite, self-adjoint and satisfies Hilbert's decomposition from Proposition \ref{hilbert-prop 3}, we easily obtain, by the spectral theorem for compact self-adjoint operators, writing $g$ in the Hilbert basis of eigenvectors of $\mathfrak{L}$, that
	\begin{equation*}
		\left\| g -\int_{\mathbb{R}^3}g_*M_*dv_*\right\|_{L^2( Mdv)}^2
		\leq C
		\int_{\mathbb{R}^3} g \mathfrak{L} g Mdv,
	\end{equation*}
	which concludes the justification of the proposition.
\end{proof}

It is to be emphasized that, since we are only considering here the case of Maxwellian molecules $b\equiv 1$, the linearized operator $\mathfrak{L}$ can be explicitly rewritten, using that $\int_{\mathbb{S}^2}\left(g'-g_*'\right)d\sigma = 0$, as
\begin{equation*}
	\mathfrak{L}g=\left|\mathbb{S}^2\right|\left(g-\int_{\mathbb{R}^3}g_*M_*dv_*\right),
\end{equation*}
which renders the proofs of Propositions \ref{hilbert-prop 3} and \ref{coercivity 3} trivial. However, we chose to provide more robust justifications of both propositions, which work in more general settings of hard and soft potentials, as well.

\subsection{Infinitesimal Maxwellians}

Using the usual relative entropy and entropy dissipation bounds together with the coercivity of the linearized collision operator, we easily get that each species of particles reaches almost instantaneously the local thermodynamic equilibrium in the fast relaxation limit. More precisely, we have the following lemma.

\begin{lem}\label{relaxation-control}
	Let $f_\eps(t,x,v)$ be a family of measurable, almost everywhere non-negative distribution functions such that, for all $t\geq 0$,
	\begin{equation*}
		\frac1{\eps^2} H\left(f_\eps\right)(t)
		+
		\frac{1}{\epsilon^4}\int_0^t\int_{\mathbb{R}^3}
		D\left(f_\eps\right)(s) dx ds
		\leq C^\mathrm{in}.
	\end{equation*}
	
	Then, as $\eps\rightarrow 0$, any subsequence of renormalized fluctuations $\hat g_\eps$ satisfies the relaxation estimate
	\begin{equation}\label{relaxation-est}
		\left\| \hat g_\eps -\Pi \hat g_\eps \right\|_{L^2( Mdv)}
		\leq
		O(\eps)
		\left\| \hat g_\eps \right\|_{L^2(Mdv)} ^2
		+O\left(\eps\right)_{L^2\left(dtdx\right)},
	\end{equation}
	where $\Pi$ denotes the orthogonal projection on $\Ker \cL$ in $L^2(Mdv)$.
\end{lem}

\begin{proof}
	We start from the elementary decomposition
	\begin{equation}\label{quadratique}
		\cL \hat g_\eps =
		\frac\eps 2  \cQ\left(\hat g_\eps ,\hat g_\eps \right)-\frac 2{\eps}  \cQ \left(\sqrt{G_\eps},\sqrt{G_\eps}\right),
	\end{equation}
	and we estimate each term in the right-hand side separately.
	
	First, since $b\equiv 1$, it is readily seen that the quadratic collision operator is continuous on $L^2(Mdv)$~:
	\begin{equation*}
		\begin{aligned}
			\left\| \cQ(\hat g_\eps ,\hat g_\eps )\right\|_{L^2( Mdv)}
			& =
			\left\| \int_{\mathbb{R}^3\times\mathbb{S}^2} \left(\hat g_\eps ' \hat g_{\eps*}'-\hat g_\eps \hat g_{\eps *}\right) M_* dv_*d\sigma\right\|_{L^2( Mdv)}
			\\
			& \leq
			\left|\mathbb{S}^2\right|^\frac{1}{2}\left\| \left(\int_{\mathbb{R}^3\times\mathbb{S}^2} \left(\hat g_\eps ' \hat g_{\eps*}'-\hat g_\eps \hat g_{\eps *}\right)^2
			M_* dv_*d\sigma\right)^\frac{1}{2}\right\|_{L^2( Mdv)}
			\\
			& \leq 2\left|\mathbb{S}^2\right| \left\| \hat g_\eps \right\|_{L^2(Mdv)} ^2,
		\end{aligned}
	\end{equation*}
	which, when combined with the bound \eqref{g-L2} from the proof of Lemma \ref{L2-lem}, yields
	\begin{equation*}
		\left\| \cQ(\hat g_\eps ,\hat g_\eps )\right\|_{L^1\left(dx;L^2(Mdv)\right)}
		\leq 8\left|\mathbb{S}^2\right| C^\mathrm{in}.
	\end{equation*}
	
	Furthermore, employing the uniform $L^2$-estimate \eqref{q-est} from the proof of Lemma \ref{L2-qlem} on the renormalized collision integrands $\hat q_\eps$ and the Cauchy-Schwarz inequality, we deduce that
	\begin{equation*}
		\begin{aligned}
		\int_0^t\int_{\mathbb{R}^3}\int_{\mathbb{R}^3}
		& \left( \frac1{\eps^2} \cQ \left(\sqrt{G_\eps},\sqrt{G_\eps}\right)\right)^2
		Mdvdxds
		\\
		& =\frac 14
		\int_0^t\int_{\mathbb{R}^3}\int_{\mathbb{R}^3}
		\left( \int_{\mathbb{R}^3\times\mathbb{S}^2} \hat q_\eps M_* dv_*d\sigma \right)^2
		Mdvdxds
		\\
		& \leq \frac{\left|\mathbb{S}^2\right|}{4}
		\int_0^t\int_{\mathbb{R}^3}
		\int_{\mathbb{R}^3\times\mathbb{R}^3\times\mathbb{S}^2} \left(\hat q_\eps\right)^2 M M_* dvdv_*d\sigma dxds
		\\
		& \leq {\left|\mathbb{S}^2\right|\over \eps^4}\int_0^t \int_{\mathbb{R}^3} D(f_\eps )(s) dxds
		\leq \left|\mathbb{S}^2\right| C^\mathrm{in}.
		\end{aligned}
	\end{equation*}

	Therefore, combining \eqref{quadratique} with the coercivity estimate from Proposition \ref{coercivity} leads to
	\begin{equation*}
		\begin{aligned}
			\left\| \hat g_\eps -\Pi \hat g_\eps \right\|_{L^2(Mdv)}
			& \leq C\left\| \mathcal{L}\hat g_\eps \right\|_{L^2(Mdv)}
			\\
			& \leq C\eps\left(
			\left\| \cQ(\hat g_\eps ,\hat g_\eps )\right\|_{L^2( Mdv)}
			+
			\left\|\frac{1}{\eps^2} \cQ\left(\sqrt{G_\eps},\sqrt{G_\eps}\right)\right\|_{L^2( Mdv)}
			\right)
			\\
			& \leq C\eps\left(
			\left\| \hat g_\eps \right\|_{L^2(Mdv)} ^2
			+
			\left\|\frac{1}{\eps^2} \cQ\left(\sqrt{G_\eps},\sqrt{G_\eps}\right)\right\|_{L^2( Mdv)}
			\right)
			\\
			& = C\eps
			\left\| \hat g_\eps \right\|_{L^2(Mdv)} ^2
			+O\left(\eps\right)_{L^2\left(dtdx\right)},
		\end{aligned}
	\end{equation*}
	which concludes the proof of the lemma.
\end{proof}

\subsection{Bulk velocity and temperature}\label{bulk}

When considering the two species Vlasov-Maxwell-Boltzmann system \eqref{VMB2}, we have an additional relaxation estimate on bulk velocities and temperatures coming from the mixed entropy dissipation.

\begin{lem}\label{relaxation2-control}
	Let $f_\eps^+(t,x,v)$ and $f_\eps^-(t,x,v)$ be two families of measurable, almost everywhere non-negative distribution functions such that, for all $t\geq 0$,
	\begin{equation*}
		\begin{aligned}
			\frac1{\eps^2} H\left(f_\eps^{+}\right)
			& + \frac1{\eps^2} H\left(f_\eps^{-}\right)
			\\
			& +\frac{1}{\epsilon^4}\int_0^t\int_{\mathbb{R}^3}\left(D\left(f_\eps^+\right)+D\left(f_\eps^-\right)
			+ \delta^2 D\left(f_\eps^+,f_\eps^-\right)\right)(s) dx ds
			\leq C^\mathrm{in}.
		\end{aligned}
	\end{equation*}
	
	Then, as $\eps\rightarrow 0$, any subsequence of renormalized fluctuations $\hat g_\eps^\pm$ satisfies the relaxation estimate
	\begin{equation}\label{relaxation2-control 1}
		\left\|
		\begin{pmatrix} \hat g_\eps^+ \\ \hat g_\eps^- \end{pmatrix}
		-\P
		\begin{pmatrix} \hat g_\eps^+ \\ \hat g_\eps^- \end{pmatrix}
		\right\|_{L^2(Mdv)}
		\leq
		O(\eps)
		\left\|
		\begin{pmatrix} \hat g_\eps^+ \\ \hat g_\eps^- \end{pmatrix}
		\right\|_{L^2(Mdv)}^2
		+O\left(\frac\eps\delta\right)_{L^2_{\mathrm{loc}}\left(dt;L^2\left(dx\right)\right)},
	\end{equation}
	where $\P$ denotes the orthogonal projection on $\Ker \L$ in $L^2(Mdv)$.
	
	In particular, further considering the densities $\hat\rho_\eps^\pm$, bulk velocities $\hat u_\eps^\pm$ and temperatures $\hat \theta_\eps^\pm$ respectively associated with the renormalized fluctuations $\hat g_\eps^\pm$, it holds that
	\begin{equation}\label{relaxation estimate}
		\hat h_\eps = \frac{\delta}{\eps}\left[\left(\hat g_\eps^+-\hat g_\eps^-\right) - \hat n_\eps\right]
		\quad\text{is uniformly bounded in }L^1_{\mathrm{loc}}\left(dtdx;L^2(Mdv)\right),
	\end{equation}
	where $\hat n_\eps = \hat\rho_\eps^+ - \hat\rho_\eps^-$, and
	\begin{equation*}
		\begin{gathered}
			\hat j_\eps=\frac\delta\eps\left(\hat u_\eps^+-\hat u_\eps^-\right)
			\qquad\text{and}\qquad
			\hat w_\eps=\frac\delta\eps\left(\hat \theta_\eps^+-\hat\theta_\eps^-\right)
			\\
			\text{are uniformly bounded in }L^1_{\mathrm{loc}}\left(dtdx\right).
		\end{gathered}
	\end{equation*}
	
	Finally, one also has the refined relaxation estimate
	\begin{equation}\label{relaxation2-control 2}
		\begin{aligned}
			& \left\|\hat h_\eps-\frac\delta 2\hat n_\eps\left(\hat g_\eps^\pm - \hat \rho_\eps^\pm\right)\right\|_{L^2(Mdv)}
			\\
			& \hspace{20mm} \leq
			O(\delta) \left\|\hat g_\eps^+ - \hat g_\eps^- -\hat n_\eps\right\|_{L^2(Mdv)}
			\left\| \hat g_\eps^\pm \right\|_{L^2(Mdv)}
			+O\left(1\right)_{L^2_{\mathrm{loc}}\left(dt;L^2\left(dx\right)\right)}.
		\end{aligned}
	\end{equation}
\end{lem}

\begin{proof}
	First, a direct application of Lemma \ref{relaxation-control} yields
	\begin{equation}\label{pre quadratique2}
		\left\| \hat g_\eps^\pm -\Pi \hat g_\eps^\pm \right\|_{L^2( Mdv)}
		\leq
		O(\eps)
		\left\| \hat g_\eps^\pm \right\|_{L^2(Mdv)} ^2
		+O\left(\eps\right)_{L^2_{\mathrm{loc}}\left(dt;L^2\left(dx\right)\right)}.
	\end{equation}
	
	Next, we apply similar arguments from the proof of Lemma \ref{relaxation-control} to the mixed entropy dissipation $D\left(f_\eps^+,f_\eps^-\right)$. Thus, according to the definitions of $\cL (g,h)$ and $\cQ(g,h)$ in \eqref{def L and Q two species}, we start from the elementary decomposition
	\begin{equation}\label{quadratique2}
		\cL \left(\hat g_\eps^\pm,\hat g_\eps^\mp\right) =
		\frac\eps 2  \cQ\left(\hat g_\eps^\pm ,\hat g_\eps^\mp \right)-\frac 2{\eps}  \cQ \left(\sqrt{G_\eps^\pm},\sqrt{G_\eps^\mp}\right),
	\end{equation}
	and we estimate each term in the right-hand side separately.
	
	Since $b\equiv 1$, it is readily seen that the quadratic collision operator is continuous on $L^2(Mdv)$~:
	\begin{equation}\label{Q continuous}
		\begin{aligned}
			\left\| \cQ\left(\hat g_\eps^\pm ,\hat g_\eps^\mp \right)\right\|_{L^2( Mdv)}
			& =
			\left\| \int_{\mathbb{R}^3\times\mathbb{S}^2} \left(\hat g_\eps^{\pm\prime} \hat g_{\eps*}^{\mp\prime}-\hat g_\eps^\pm \hat g_{\eps *}^\mp\right) M_* dv_*d\sigma\right\|_{L^2( Mdv)}
			\\
			& \leq
			\left|\mathbb{S}^2\right|^\frac{1}{2}\left\| \left(\int_{\mathbb{R}^3\times\mathbb{S}^2} \left(\hat g_\eps^{\pm\prime} \hat g_{\eps*}^{\mp\prime}-\hat g_\eps^\pm \hat g_{\eps *}^\mp\right)^2
			M_* dv_*d\sigma\right)^\frac{1}{2}\right\|_{L^2( Mdv)}
			\\
			& \leq 2\left|\mathbb{S}^2\right| \left\| \hat g_\eps^+ \right\|_{L^2(Mdv)} \left\| \hat g_\eps^- \right\|_{L^2(Mdv)},
		\end{aligned}
	\end{equation}
	which, when combined with the bound \eqref{g-L2} from the proof of Lemma \ref{L2-lem}, yields
	\begin{equation*}
		\left\| \cQ\left(\hat g_\eps^\pm ,\hat g_\eps^\mp \right)\right\|_{L^1\left(dx;L^2(Mdv)\right)}
		\leq 8\left|\mathbb{S}^2\right| C^\mathrm{in}.
	\end{equation*}

	Furthermore, employing the uniform $L^2$-estimate \eqref{q-est} from the proof of Lemma \ref{L2-qlem} on the renormalized collision integrands $\hat q_\eps^\pm$ and the Cauchy-Schwarz inequality, we deduce that
	\begin{equation}\label{q plus}
		\begin{aligned}
		\int_0^t\int_{\mathbb{R}^3}\int_{\mathbb{R}^3}
		& \left( \frac\delta{\eps^2} \cQ \left(\sqrt{G_\eps^+},\sqrt{G_\eps^-}\right)\right)^2
		Mdvdxds
		\\
		& =\frac 14
		\int_0^t\int_{\mathbb{R}^3}\int_{\mathbb{R}^3}
		\left( \int_{\mathbb{R}^3\times\mathbb{S}^2} \hat q_\eps^{+,-} M_* dv_*d\sigma \right)^2
		Mdvdxds
		\\
		& \leq \frac{\left|\mathbb{S}^2\right|}{4}
		\int_0^t\int_{\mathbb{R}^3}
		\int_{\mathbb{R}^3\times\mathbb{R}^3\times\mathbb{S}^2} \left(\hat q_\eps^{+,-}\right)^2 M M_* dvdv_*d\sigma dxds
		\\
		& \leq {\left|\mathbb{S}^2\right|\delta^2\over 2\eps^4}\int_0^t \int_{\mathbb{R}^3} D\left(f_\eps^+,f_\eps^- \right)(s) dxds
		\leq \frac{\left|\mathbb{S}^2\right|}{2} C^\mathrm{in}.
		\end{aligned}
	\end{equation}
	Notice that the same estimate holds on the renormalized collision integrands $\hat q_\eps^{-,+}$, which yields
	\begin{equation}\label{q minus}
		\begin{aligned}
		\int_0^t\int_{\mathbb{R}^3}\int_{\mathbb{R}^3}
		& \left( \frac\delta{\eps^2} \cQ \left(\sqrt{G_\eps^-},\sqrt{G_\eps^+}\right)\right)^2
		Mdvdxds
		\\
		% & =\frac 14
		% \int_0^t\int_{\mathbb{R}^3}\int_{\mathbb{R}^3}
		% \left( \int_{\mathbb{R}^3\times\mathbb{S}^2} \hat q_\eps^{-,+} M_* dv_*d\sigma \right)^2
		% Mdvdxds
		% \\
		% & \leq \frac{\left|\mathbb{S}^2\right|}{4}
		% \int_0^t\int_{\mathbb{R}^3}
		% \int_{\mathbb{R}^3\times\mathbb{R}^3\times\mathbb{S}^2} \left(\hat q_\eps^{-,+}\right)^2 M M_* dvdv_*d\sigma dxds
		% \\
		& \leq {\left|\mathbb{S}^2\right|\delta^2\over 2\eps^4}\int_0^t \int_{\mathbb{R}^3} D\left(f_\eps^+,f_\eps^- \right)(s) dxds
		\leq \frac{\left|\mathbb{S}^2\right|}{2} C^\mathrm{in}.
		\end{aligned}
	\end{equation}

	Therefore, combining \eqref{pre quadratique2} and \eqref{quadratique2} with the coercivity estimate from Proposition \ref{coercivity 2} leads to
	\begin{equation*}
		\begin{aligned}
			& \left\|
			\begin{pmatrix} \hat g_\eps^+ \\ \hat g_\eps^- \end{pmatrix}
			-\P
			\begin{pmatrix} \hat g_\eps^+ \\ \hat g_\eps^- \end{pmatrix}
			\right\|_{L^2(Mdv)}
			\\
			& \leq C\left\|\L
			\begin{pmatrix} \hat g_\eps^+ \\ \hat g_\eps^- \end{pmatrix}
			\right\|_{L^2(Mdv)}
			\leq C
			\left\|
			\begin{pmatrix} \mathcal{L}\hat g_\eps^+ \\ \mathcal{L}\hat g_\eps^- \end{pmatrix}
			\right\|_{L^2(Mdv)}
			+ C
			\left\|
			\begin{pmatrix} \mathcal{L}\left(\hat g_\eps^+,\hat g_\eps^-\right) \\ \mathcal{L}\left(\hat g_\eps^-,\hat g_\eps^+\right) \end{pmatrix}
			\right\|_{L^2(Mdv)}
			\\
			& \leq C
			\left\|
			\begin{pmatrix} \hat g_\eps^+ - \Pi\hat g_\eps^+ \\ \hat g_\eps^- - \Pi\hat g_\eps^- \end{pmatrix}
			\right\|_{L^2(Mdv)}
			+ C\eps
			\left\|
			\begin{pmatrix} \mathcal{Q}\left(\hat g_\eps^+,\hat g_\eps^-\right) \\ \mathcal{Q}\left(\hat g_\eps^-,\hat g_\eps^+\right) \end{pmatrix}
			\right\|_{L^2(Mdv)}
			\\
			& +C\frac\eps\delta
			\left\|\frac{\delta}{\eps^2}
			\begin{pmatrix} \mathcal{Q}\left(\sqrt{G_\eps^+},\sqrt{G_\eps^-}\right)
			\\
			\mathcal{Q}\left(\sqrt{G_\eps^-},\sqrt{G_\eps^+}\right) \end{pmatrix}
			\right\|_{L^2(Mdv)}
			\\
			& \leq O(\eps)
			\left\|
			\begin{pmatrix} \hat g_\eps^+ \\ \hat g_\eps^- \end{pmatrix}
			\right\|_{L^2(Mdv)}^2 + O(\eps)_{L^2_\mathrm{loc}\left(dt;L^2(dx)\right)}
			+ O\left(\frac\eps\delta\right)_{L^2_\mathrm{loc}\left(dt;L^2(dx)\right)},
		\end{aligned}
	\end{equation*}
	which concludes the proof of the relaxation estimate \eqref{relaxation2-control 1}.
	
	Then, in order to deduce the control of $\hat h_\eps$, $\hat j_\eps$ and $\hat w_\eps$, it suffices to notice that
	\begin{equation*}
		\begin{pmatrix}
			1\\-1
		\end{pmatrix}\cdot
		\left[\begin{pmatrix} \hat g_\eps^+ \\ \hat g_\eps^- \end{pmatrix}
		-\P
		\begin{pmatrix} \hat g_\eps^+ \\ \hat g_\eps^- \end{pmatrix}\right]
		= \left(\hat g_\eps^+-\hat g_\eps^-\right)-\left(\hat\rho_\eps^+ - \hat\rho_\eps^-\right),
	\end{equation*}
	and
	\begin{equation*}
		\begin{aligned}
			\Pi\left[
			\begin{pmatrix} \hat g_\eps^+ \\ \hat g_\eps^- \end{pmatrix}
			-\P
			\begin{pmatrix} \hat g_\eps^+ \\ \hat g_\eps^- \end{pmatrix}
			\right]
			& =
			\begin{pmatrix} \Pi \frac{\hat g_\eps^+ - \hat g_\eps^- }{2}
			- \frac 12\int_{\mathbb{R}^3}\left(\hat g_\eps^+-\hat g_\eps^-\right)Mdv
			\\
			\Pi \frac{\hat g_\eps^- - \hat g_\eps^+ }{2}
			- \frac 12\int_{\mathbb{R}^3}\left(\hat g_\eps^--\hat g_\eps^+\right)Mdv \end{pmatrix}
			\\
			& =
			\begin{pmatrix}
				\frac{\hat u_\eps^+-\hat u_\eps^-}{2}\cdot v +\frac{\hat\theta_\eps^+ - \hat\theta_\eps^-}{2}\left(\frac{|v|^2}{2}-\frac 32\right)
				\\
				\frac{\hat u_\eps^--\hat u_\eps^+}{2}\cdot v +\frac{\hat\theta_\eps^- - \hat\theta_\eps^+}{2}\left(\frac{|v|^2}{2}-\frac 32\right)
			\end{pmatrix},
		\end{aligned}
	\end{equation*}
	whence
	\begin{equation*}
		\left\|\Pi\left[
		\begin{pmatrix} \hat g_\eps^+ \\ \hat g_\eps^- \end{pmatrix}
		-\P
		\begin{pmatrix} \hat g_\eps^+ \\ \hat g_\eps^- \end{pmatrix}
		\right]\right\|_{L^2(Mdv)}^2
		=
		\frac 12\left(\hat u_\eps^+-\hat u_\eps^-\right)^2 +
		\frac 34\left(\hat\theta_\eps^+ - \hat\theta_\eps^-\right)^2.
	\end{equation*}

	There only remains to establish the more precise relaxation estimate \eqref{relaxation2-control 2} on $\hat h_\eps$, which is achieved by employing the coercivity of the operator $\mathfrak{L}$. To this end, we use the identities \eqref{quadratique2} to decompose
	\begin{equation*}
		\begin{aligned}
			\mathfrak{L}\hat h_\eps & =
			\frac\delta\eps\mathfrak{L}\left(\hat g_\eps^+-\hat g_\eps^-\right)
			\\
			& = \frac\delta 2  \left[\cQ\left(\hat g_\eps^+ ,\hat g_\eps^- \right) - \cQ\left(\hat g_\eps^- ,\hat g_\eps^+\right)\right]
			-\frac {2\delta}{\eps^2} \left[\cQ \left(\sqrt{G_\eps^+},\sqrt{G_\eps^-}\right)
			- \cQ \left(\sqrt{G_\eps^-},\sqrt{G_\eps^+}\right)\right]
			\\
			& =
			\frac\delta 2  \left[\cQ\left(\hat g_\eps^+-\hat g_\eps^- ,\hat g_\eps^\pm \right)
			- \cQ\left(\hat g_\eps^\pm ,\hat g_\eps^+-\hat g_\eps^-\right)\right]
			-
			\int_{\mathbb{R}^3\times\mathbb{S}^2} \left(\hat q_\eps^{+,-}-\hat q_\eps^{-,+}\right) M_* dv_*d\sigma.
		\end{aligned}
	\end{equation*}
	It follows that
	\begin{equation*}
		\begin{aligned}
			\mathfrak{L}\left(\hat h_\eps - \frac\delta 2 \hat n_\eps \hat g_\eps^\pm\right)
			& =
			\frac\eps 2  \left[\cQ\left(\hat h_\eps ,\hat g_\eps^\pm \right)
			- \cQ\left(\hat g_\eps^\pm ,\hat h_\eps\right)\right]
			\\
			& -
			\int_{\mathbb{R}^3\times\mathbb{S}^2} \left(\hat q_\eps^{+,-}-\hat q_\eps^{-,+}\right) M_* dv_*d\sigma.
		\end{aligned}
	\end{equation*}
	
	Therefore, repeating the estimates \eqref{Q continuous}, \eqref{q plus} and \eqref{q minus}, we find that
	\begin{equation*}
		\begin{aligned}
			\left\|\mathfrak{L}\left(\hat h_\eps - \frac\delta 2 \hat n_\eps \hat g_\eps^\pm\right)\right\|_{L^2(Mdv)}
			& \leq
			O(\eps) \left\|\hat h_\eps \right\|_{L^2(Mdv)}
			\left\| \hat g_\eps^\pm \right\|_{L^2(Mdv)}
			+ O\left(1\right)_{L^2_{\mathrm{loc}}\left(dt;L^2\left(dx\right)\right)}
			\\
			& =
			O(\delta) \left\|\hat g_\eps^+ - \hat g_\eps^- - \hat n_\eps\right\|_{L^2(Mdv)}
			\left\| \hat g_\eps^\pm \right\|_{L^2(Mdv)} \\
			& + O\left(1\right)_{L^2_{\mathrm{loc}}\left(dt;L^2\left(dx\right)\right)}.
		\end{aligned}
	\end{equation*}
	Finally, employing the coercivity of $\mathfrak{L}$ from Proposition \ref{coercivity 3}, we easily conclude that \eqref{relaxation2-control 2} holds, which concludes the proof of the lemma.
\end{proof}

\begin{rem}
	Under the hypotheses of the preceding lemma, it is possible to obtain a very explicit identity providing some improved information on the relaxation of $\hat h_\eps$. To this end, we decompose
	\begin{equation*}
		\begin{aligned}
			\int_{\mathbb{R}^3\times\mathbb{S}^2} & \left(\hat q_\eps^{+,-}-\hat q_\eps^{-,+}\right) M_* dv_*d\sigma
			\\
			& =
			\frac{2\delta}{\eps^2}\mathcal{Q}\left(\sqrt{G^+_\eps},\sqrt{G^-_\eps}\right)
			-
			\frac{2\delta}{\eps^2}\mathcal{Q}\left(\sqrt{G^-_\eps},\sqrt{G^+_\eps}\right)
			\\
			& =
			\frac{2\delta}{\eps^2}\mathcal{Q}\left(\sqrt{G^+_\eps}-\sqrt{G^-_\eps}
			-\frac\eps 2  \frac{\hat n_\eps}{1+\frac\eps 2\hat\rho_\eps^\pm}\sqrt{G^\pm_\eps},\sqrt{G^\pm_\eps}\right)
			\\
			& -
			\frac{2\delta}{\eps^2}\mathcal{Q}\left(\sqrt{G^\pm_\eps},\sqrt{G^+_\eps}-\sqrt{G^-_\eps}
			-\frac\eps 2  \frac{\hat n_\eps}{1+\frac\eps 2\hat\rho_\eps^\pm}\sqrt{G^\pm_\eps}\right)
			\\
			& =
			\mathcal{Q}\left(
			\hat h_\eps - \frac\delta 2 \hat n_\eps \frac{\hat g_\eps^\pm-\hat \rho_\eps^\pm}{1+\frac \eps 2\hat\rho_\eps^\pm}
			,\sqrt{G^\pm_\eps}\right) -
			\mathcal{Q}\left(\sqrt{G^\pm_\eps},
			\hat h_\eps - \frac\delta 2 \hat n_\eps \frac{\hat g_\eps^\pm-\hat \rho_\eps^\pm}{1+\frac \eps 2\hat\rho_\eps^\pm}
			\right).
		\end{aligned}
	\end{equation*}
	Then, since we are only considering here the Maxwellian cross-section $b\equiv 1$, notice that the gain terms cancel each other out
	\begin{equation*}
		\mathcal{Q}^+\left(
		\hat h_\eps - \frac\delta 2 \hat n_\eps \frac{\hat g_\eps^\pm-\hat \rho_\eps^\pm}{1+\frac \eps 2\hat\rho_\eps^\pm}
		,\sqrt{G^\pm_\eps}\right) -
		\mathcal{Q}^+\left(\sqrt{G^\pm_\eps},
		\hat h_\eps - \frac\delta 2 \hat n_\eps \frac{\hat g_\eps^\pm-\hat \rho_\eps^\pm}{1+\frac \eps 2\hat\rho_\eps^\pm}
		\right)=0,
	\end{equation*}
	for the change of variable $\sigma\mapsto -\sigma$ merely exchanges $v'$ and $v_*'$, and that one of the two loss terms vanishes
	\begin{equation*}
		\mathcal{Q}^-\left(\sqrt{G^\pm_\eps},
		\hat h_\eps - \frac\delta 2 \hat n_\eps \frac{\hat g_\eps^\pm-\hat \rho_\eps^\pm}{1+\frac \eps 2\hat\rho_\eps^\pm}
		\right)=0.
	\end{equation*}
	Thus, on the whole, we are left with the identity
	\begin{equation*}
		\begin{aligned}
			\left|\mathbb{S}^2\right|\left(\left(1+\frac \eps 2\hat\rho_\eps^\pm\right)
			\hat h_\eps - \frac\delta 2 \hat n_\eps \left(\hat g_\eps^\pm-\hat \rho_\eps^\pm\right)
			\right)
			& =
			\mathcal{Q}^-\left(
			\hat h_\eps - \frac\delta 2 \hat n_\eps \frac{\hat g_\eps^\pm-\hat \rho_\eps^\pm}{1+\frac \eps 2\hat\rho_\eps^\pm}
			,\sqrt{G^\pm_\eps}\right)
			\\
			& =
			- \int_{\mathbb{R}^3\times\mathbb{S}^2} \left(\hat q_\eps^{+,-}-\hat q_\eps^{-,+}\right) M_* dv_*d\sigma,
		\end{aligned}
	\end{equation*}
	which yields the control, in view of Lemma \ref{L2-qlem},
	\begin{equation*}
		\left(1+\frac \eps 2\hat\rho_\eps^\pm\right)
		\hat h_\eps - \frac\delta 2 \hat n_\eps \left(\hat g_\eps^\pm-\hat \rho_\eps^\pm\right)
		= O(1)_{L^2\left(Mdtdxdv\right)}.
	\end{equation*}
	In particular, further integrating against $vMdv$ and $\left(\frac{|v|^2}{3}-1\right)Mdv$, we obtain that
	\begin{equation*}
		\begin{aligned}
			\left(1+\frac \eps 2\hat\rho_\eps^\pm\right)
			\hat j_\eps - \frac\delta 2 \hat n_\eps \hat u_\eps^\pm
			& = O(1)_{L^2\left(dtdx\right)},
			\\
			\left(1+\frac \eps 2\hat\rho_\eps^\pm\right)
			\hat w_\eps - \frac\delta 2 \hat n_\eps \hat \theta_\eps^\pm
			& = O(1)_{L^2\left(dtdx\right)}.
		\end{aligned}
	\end{equation*}
	These estimates are slightly more precise than \eqref{relaxation2-control 2}. However, their significance is unclear.
\end{rem}

% ===================
% = High velocities =
% ===================

\section{Improved integrability in velocity}\label{velocity integrability}

Another important consequence of the control on the relaxation is to provide further integrability on the renormalized fluctuations $\hat g_\eps$ and $\hat g_\eps^\pm$ with respect to the $v$ variable at infinity. The following result, which was first established as such in \cite{SR} (see Lemma 3.2.5 therein) and \cite{golse5} (see Proposition 3.2 therein), improves Lemma \ref{L2-lem} and is a direct consequence of Lemma \ref{relaxation-control}. It constitutes a significant simplification with respect to earlier works on hydrodynamic limits of the Boltzmann equation, which required convoluted estimates to establish some improved integrability in the $v$ variable.

\begin{lem}\label{v2-int}
	Let $f_\eps(t,x,v)$ be a family of measurable, almost everywhere non-negative distribution functions such that, for all $t\geq 0$,
	\begin{equation*}
		\frac1{\eps^2} H\left(f_\eps\right)(t)
		+
		\frac{1}{\epsilon^4}\int_0^t\int_{\mathbb{R}^3}
		D\left(f_\eps\right)(s) dx ds
		\leq C^\mathrm{in}.
	\end{equation*}

	Then, as $\eps\rightarrow 0$, any subsequence of renormalized fluctuations $\hat g_\eps$ is uniformly bounded in $L^2_{\mathrm{loc}}\left(dtdx;L^2\left(\left(1+|v|^2\right)Mdv\right)\right)$.
	
	Furthermore, the family $|\hat g_\eps|^2$ is equi-integrable in $v$ (or uniformly integrable in $v$) in the sense that, for any $\eta>0$ and every compact subset $K\subset [0,\infty)\times\mathbb{R}^3\times\mathbb{R}^3$, there exists $\gamma>0$ such that, if $A\subset K$ is a measurable set satisfying
	\begin{equation*}
		\sup_{(t,x)\in [0,\infty)\times\mathbb{R}^3}
		\int_{\mathbb{R}^3} \mathds{1}_{A}(t,x,v) dv < \gamma ,
	\end{equation*}
	then
	\begin{equation*}
		\sup_{\eps>0}\int_A |\hat g_\eps|^2 dtdxdv < \eta.
	\end{equation*}
	
	We also have that, for any $\lambda>0$ and any $1\leq p<2$, the families $\mathds{1}_{\left\{\lambda \eps |\hat g_\eps| \leq 1 \right\}} |\hat g_\eps|^2$ and $\frac{|\hat g_\eps|^2}{1+\lambda\sqrt{G_\eps}}$ are uniformly bounded in $L^1_{\mathrm{loc}}\left(dtdx;L^p\left(Mdv\right)\right)$.
\end{lem}

\begin{proof}
	The crucial idea behind these results rests upon decomposing $\hat g_\eps$ according to
	\begin{equation}\label{Pi decomp}
		\hat g_\eps=\left(\hat g_\eps -\Pi \hat g_\eps\right)+\Pi \hat g_\eps,
	\end{equation}
	and then using the control on the relaxation provided by Lemma \ref{relaxation-control}~:
	\begin{equation*}
		\hat g_\eps - \Pi \hat g_\eps = O(\eps)_{L^1_\mathrm{loc}\left(dtdx ; L^2(Mdv)\right)}.
	\end{equation*}

	\noindent$\bullet$ First, we establish the uniform control on the high speed tails of $|\hat g_\eps|^2$, i.e.\ the uniform weighted integrability estimate in $L^2_{\mathrm{loc}}\left(dtdx;L^2\left(\left(1+|v|^2\right)Mdv\right)\right)$. To this end, we start from the decomposition
	\begin{equation}\label{new-decomposition}
		 (1+|v|)^2\left|\hat g_\eps\right|^2
		 =\hat g_\eps
		\left[(1+|v|)^2 \Pi \hat g_\eps\right]
		+ \left[(1+|v|)\hat g_\eps\right]
		\left[(1+|v|) \left(\hat g_\eps -\Pi \hat g_\eps \right)\right].
	\end{equation}
	
	Next, recalling from Lemma \ref{L2-lem} that
	\begin{equation*}
		\hat g_\eps
		\text{ is uniformly bounded in }
		L^\infty\left(dt;L^2\left(Mdxdv\right)\right),
	\end{equation*}
	we see, by definition of the hydrodynamic projection $\Pi$, for any $1\leq p<\infty$, that
	\begin{equation}\label{est0}
		(1+|v|)^2\Pi\hat g_\eps
		\text{ is uniformly bounded in }
		L^\infty\left(dt;L^2\left(dx;L^p\left(Mdv\right)\right)\right),
	\end{equation}
	whence, for any $1\leq r <2$,
	\begin{equation}\label{est1}
		\hat g_\eps \left[(1+|v|)^2 \Pi \hat g_\eps\right]
		\text{ is uniformly bounded in }
		L^\infty\left(dt;L^1\left(dx;L^r\left(Mdv\right)\right)\right),
	\end{equation}
	which takes care of the first term in right-hand side of \eqref{new-decomposition}.
	
	In order to estimate the second term in the right-hand side of \eqref{new-decomposition}, we first apply Young's inequality \eqref{Young h} with $z=g_\eps$, $y={(1+|v|)^2\over\gamma}$, $\alpha=4\gamma$, $\beta = \frac {4\gamma}\eps$ and $\gamma>0$, to get, employing the elementary inequality \eqref{sqrt abs},
	\begin{equation}\label{v2 trick}
		\begin{aligned}
		(1+|v|)^{2} \hat g_\eps ^2 
		&\leq \left(1+|v|\right)^2\frac{4} {\eps^2} \left|{G_\eps} -1\right|
		= \frac{4\gamma} {\eps} \left|g_\eps\right| {(1+|v|)^2\over\gamma}
		\\
		& \leq {4\gamma \over \eps^2} h\left(\eps g_\eps\right)
		+ \frac{4\gamma}{\eps^2}e^{\frac{(1+|v|)^2}{\gamma}}.
		\end{aligned}
	\end{equation}
	% \begin{equation}\label{v2 trick}
	% 	\begin{aligned}
	% 	(1+|v|)^{2} \hat g_\eps ^2
	% 	&\leq \left(1+|v|\right)^2\frac{4} {\eps^2} \left|{G_\eps} -1\right|
	% 	= \frac{16} {\eps} \left|g_\eps\right| {(1+|v|)^2\over4}
	% 	\\
	% 	& \leq {16 \over \eps^2} h\left(\eps g_\eps\right)
	% 	+ \frac{16}{\eps^2}e^{\frac{(1+|v|)^2}{4}}.
	% 	\end{aligned}
	% \end{equation}
	Therefore, setting $\gamma=4$ in \eqref{v2 trick}, we obtain that
	\begin{equation*}
		\begin{aligned}
			|(1+|v| & ) \hat g_\eps
			(1+|v|) \left(\hat g_\eps -\Pi \hat g_\eps \right)|
			\\
			& \leq
			{4 \over \eps} \sqrt{h\left(\eps g_\eps\right)}(1+|v|) \left|\hat g_\eps -\Pi \hat g_\eps \right|
			+
			\frac{4}{\eps}e^{\frac{(1+|v|)^2}{8}}(1+|v|) \left|\hat g_\eps -\Pi \hat g_\eps \right|
			\\
			& \leq
			{16 \over \eps^2} {h\left(\eps g_\eps\right)}+\frac 14 (1+|v|)^2 \left|\hat g_\eps -\Pi \hat g_\eps \right|^2
			+
			\frac{4}{\eps}e^{\frac{(1+|v|)^2}{8}}(1+|v|) \left|\hat g_\eps -\Pi \hat g_\eps \right|
			\\
			& \leq
			{16 \over \eps^2} {h\left(\eps g_\eps\right)}+\frac 12 (1+|v|)^2 \left(\left|\hat g_\eps\right|^2 +\left|\Pi \hat g_\eps\right|^2 \right)
			+
			\frac{4}{\eps}e^{\frac{(1+|v|)^2}{8}}(1+|v|) \left|\hat g_\eps -\Pi \hat g_\eps \right|,
		\end{aligned}
	\end{equation*}
	which, by virtue of the uniform entropy bound, the uniform estimate \eqref{est0} and the relaxation estimate \eqref{relaxation-est} from Lemma \ref{relaxation-control}, yields
	\begin{equation}\label{est2}
		|(1+|v|) \hat g_\eps
		(1+|v|) \left(\hat g_\eps -\Pi \hat g_\eps \right)|
		\leq O(1)_{L^1_{\mathrm{loc}}\left(dtdx;L^1\left(Mdv\right)\right)}
		+\frac 12 (1+|v|)^2 \left|\hat g_\eps\right|^2.
	\end{equation}

	% In other words, setting
	% \begin{equation*}
	% 	\begin{aligned}
	% 		\left(1+|v|\right)^2\left(\hat g_\eps^{(1)}\right)^2 & =\left[\left(1+|v|\right)^2\hat g_\eps^2 - \frac{16}{\eps^2}e^{\frac{(1+|v|)^2}{4}}\right]
	% 		\mathds{1}_{\left\{\frac{16}{\eps^2}e^{\frac{(1+|v|)^2}{4}} < \left(1+|v|\right)^2\hat g_\eps^2\right\}},
	% 		\\
	% 		\left(1+|v|\right)^2\left(\hat g_\eps^{(2)}\right)^2 & =\left(1+|v|\right)^2\hat g_\eps^2
	% 		\mathds{1}_{\left\{\left(1+|v|\right)^2\hat g_\eps^2 \leq \frac{16}{\eps^2}e^{\frac{(1+|v|)^2}{4}} \right\}}
	% 		\\
	% 		& +
	% 		\frac{16}{\eps^2}e^{\frac{(1+|v|)^2}{4}}
	% 		\mathds{1}_{\left\{\frac{16}{\eps^2}e^{\frac{(1+|v|)^2}{4}} < \left(1+|v|\right)^2\hat g_\eps^2\right\}},
	% 	\end{aligned}
	% \end{equation*}
	% we have established a decomposition $\hat g_\eps^2 =\left(\hat g_\eps^{(1)}\right)^2 +\left(\hat g_\eps^{(2)}\right)^2$ such that
	% \begin{equation*}
	% 	\begin{aligned}
	% 		\left\| (1+|v|) \hat g_\eps^{(1)}\right\|^2_{L^2\left(Mdxdv\right)} & \leq \frac{16}{\eps^2} H(f_\eps), \\
	% 		\text{and}\qquad (1+|v|) \left|g_\eps^{(2)}\right| & \leq \frac{4}{\eps}e^{\frac{(1+|v|)^2}{8}}.
	% 	\end{aligned}
	% \end{equation*}

	On the whole, incorporating \eqref{est1} and \eqref{est2} into the decomposition \eqref{new-decomposition}, we deduce that
	\begin{equation*}
		(1+|v|)^2 |\hat g_\eps|^2
		\leq O(1)_{L^1_{\mathrm{loc}}\left(dtdx;L^1\left(Mdv\right)\right)}
		+\frac 12 (1+|v|)^2 \left|\hat g_\eps\right|^2.
	\end{equation*}
	Hence, we conclude
	\begin{equation*}
		(1+|v|)^2 |\hat g_\eps|^2
		= O(1)_{L^1_{\mathrm{loc}}\left(dtdx;L^1\left(Mdv\right)\right)},
	\end{equation*}
	which is the expected result.

	\noindent$\bullet$ We establish now the uniform integrability statement of the lemma. To this end, we start from the decomposition, for any large $\lambda>e$,
	\begin{equation}\label{old-decomposition}
		\left|\hat g_\eps\right|^2 =
		\mathds{1}_{\left\{G_\eps >\lambda\right\}} \left|\hat g_\eps\right|^2
		+\mathds{1}_{\left\{G_\eps \leq \lambda\right\}} \hat g_\eps \Pi \hat g_\eps
		+\mathds{1}_{\left\{G_\eps \leq \lambda\right\}}\hat g_\eps\left(\hat g_\eps-\Pi \hat g_\eps\right).
	\end{equation}
	
	Then, we use the relative entropy bound, a pointwise estimate of $h(z)$, for $z>\lambda$, and the elementary inequality \eqref{sqrt abs} to control the large tails of $G_\eps$ as follows
	\begin{equation}\label{large values}
		\begin{aligned}
			\mathds{1}_{\left\{G_\eps >\lambda\right\}} \left|\hat g_\eps\right|^2
			& \leq\frac 4{\eps^2}
			\mathds{1}_{\left\{G_\eps >\lambda\right\}} \left|G_\eps -1\right|
			\leq\frac 8{\eps^2}
			\mathds{1}_{\left\{G_\eps >\lambda\right\}} G_\eps
			\\
			& \leq \frac{8}{\eps^2\left(\log\lambda -1\right)}\mathds{1}_{\left\{G_\eps >\lambda\right\}}
			G_\eps\left(\log G_\eps - 1\right)
			\leq \frac{8}{\eps^2\left(\log\lambda -1\right)}
			h\left(\eps g_\eps\right)
			\\
			& =O\left({1\over \log \lambda}\right)_{L^\infty\left(dt;L^1\left(Mdxdv\right)\right)},
		\end{aligned}
	\end{equation}
	which takes care of the first term in the right-hand side of \eqref{old-decomposition}, while the second term is handled by estimate \eqref{est1}. As for the remaining term in \eqref{old-decomposition}, we deduce from the relaxation estimate \eqref{relaxation-est} in Lemma \ref{relaxation-control} and from the pointwise estimate $\left|\mathds{1}_{\left\{G_\eps \leq \lambda\right\}}\hat g_\eps\right|\leq\frac{2}{\eps}\left(1+\sqrt{\lambda}\right)$, which follows straightforwardly from \eqref{sqrt abs}, that
	\begin{equation*}
		\mathds{1}_{\left\{G_\eps \leq \lambda\right\}} \hat g_\eps \left(\hat g_\eps-\Pi \hat g_\eps\right)
		=O\left(\sqrt{\lambda}\right)_{L^1_{\mathrm{loc}}\left(dtdx;L^2\left(Mdv\right)\right)}.
	\end{equation*}
	
	Thus, on the whole, we have established from the decomposition \eqref{old-decomposition} that, for any arbitrarily large $\lambda$ and each $1\leq r <2$,
	\begin{equation*}
		\left|\hat g_\eps\right|^2=
		O\left({1\over \log \lambda}\right)_{L^\infty\left(dt;L^1\left(Mdxdv\right)\right)}
		+
		O\left(\sqrt{\lambda}\right)_{L^1_{\mathrm{loc}}\left(dtdx;L^r\left(Mdv\right)\right)},
	\end{equation*}
	which clearly implies that $\left|\hat g_\eps\right|^2$ is locally uniformly integrable in $v$.

	\noindent$\bullet$ The final statement is easily obtained by combining decomposition \eqref{Pi decomp} with the bounds
	\begin{equation*}
		\begin{aligned}
			% \left|\hat g_\eps\right| \mathds{1}_{\left\{\lambda\eps |\hat g_\eps| \leq 1 \right\}}
			% \leq
			\frac{\hat g_\eps}{1+\lambda\eps\left|\hat g_\eps\right|}
			& = O(1)_{L^\infty(dt;L^2(Mdxdv))},
			\\
			% \left|\hat g_\eps\right| \mathds{1}_{\left\{\lambda\eps |\hat g_\eps| \leq 1 \right\}}
			% \leq
			\frac{\hat g_\eps}{1+\lambda\eps\left|\hat g_\eps\right|}
			& = O\left( \frac 1{\lambda\eps}\right)_{L^\infty(dtdxdv)},
		\end{aligned}
	\end{equation*}
	for any $\lambda>0$, and noticing that $1+\lambda\eps\left|\hat g_\eps\right|\leq \max\left\{2,1+2\lambda\right\}\left(1+\lambda\sqrt{G_\eps}\right)$. We indeed find, for any $1\leq p <2$, that
	\begin{equation*}
		\begin{aligned}
			\left\|\frac{\hat g_\eps^2}{1+\lambda\eps\left|\hat g_\eps\right|}\right\|_{L^p(Mdv)}
			& = \left\|\left(\hat g_\eps -\Pi \hat g_\eps\right)\frac{\hat g_\eps}{1+\lambda\eps\left|\hat g_\eps\right|}
			+ \Pi \hat g_\eps \frac{\hat g_\eps}{1+\lambda\eps\left|\hat g_\eps\right|}\right\|_{L^p(Mdv)}
			\\
			& \leq \frac 1{\lambda\eps} \left\|\hat g_\eps -\Pi \hat g_\eps\right\|_{L^p(Mdv)}
			+ \left\| \hat g_\eps \Pi \hat g_\eps \right\|_{L^p(Mdv)}
			\\
			& \leq \frac C{\lambda\eps} \left\|\hat g_\eps -\Pi \hat g_\eps\right\|_{L^2(Mdv)}
			+ C\left\| \hat g_\eps \right\|_{L^2(Mdv)}^2
			\\
			& \leq C \left\| \hat g_\eps \right\|_{L^2(Mdv)}^2 + \left\|\hat q_\eps \right\|_{L^2\left(MM_*dvdv_*d\sigma\right)}
			\\
			& = O(1)_{L^\infty\left(dt;L^1\left( dx \right)\right)} + O(1)_{L^2(dtdx)},
		\end{aligned}
	\end{equation*}
	% \begin{equation*}
	% 	\begin{aligned}
	% 		\hat g_\eps^2 \mathds{1}_{\left\{\eps |\hat g_\eps| \leq \lambda \right\}}
	% 		& = \left(\hat g_\eps -\Pi \hat g_\eps\right)\hat g_\eps \mathds{1}_{\left\{\eps |\hat g_\eps| \leq \lambda \right\}}
	% 		+ \Pi \hat g_\eps \hat g_\eps \mathds{1}_{\left\{\eps |\hat g_\eps| \leq \lambda \right\}}
	% 		\\
	% 		& = O(\lambda)_{L^1_\mathrm{loc}\left(dtdx ; L^2(Mdv)\right)} + O(1)_{L^\infty\left(dt;L^1\left( dx; L^p(Mdv)\right)\right)},
	% 	\end{aligned}
	% \end{equation*}
	which concludes the proof of the lemma.
\end{proof}

In the two species case, the preceding lemma has simple but important consequences on the integrability of the difference of fluctuations, which is the content of the next lemmas.

\begin{lem}\label{bound hjw}
	Let $f_\eps^+(t,x,v)$ and $f_\eps^-(t,x,v)$ be two families of measurable, almost everywhere non-negative distribution functions such that, for all $t\geq 0$,
	\begin{equation*}
		\begin{aligned}
			\frac1{\eps^2} H\left(f_\eps^{+}\right)
			& + \frac1{\eps^2} H\left(f_\eps^{-}\right)
			\\
			& +\frac{1}{\epsilon^4}\int_0^t\int_{\mathbb{R}^3}\left(D\left(f_\eps^+\right)+D\left(f_\eps^-\right)
			+ \delta^2 D\left(f_\eps^+,f_\eps^-\right)\right)(s) dx ds
			\leq C^\mathrm{in}.
		\end{aligned}
	\end{equation*}
	
	Then, as $\eps\rightarrow 0$, considering the densities $\rho_\eps^\pm$, bulk velocities $u_\eps^\pm$ and temperatures $\theta_\eps^\pm$ respectively associated with any subsequence of fluctuations $g_\eps^\pm$, it holds that
	\begin{equation}\label{def h}
		\begin{gathered}
			h_\eps = \frac{\delta}{\eps}\left[\left(g_\eps^+-g_\eps^-\right) - n_\eps\right]
			\\
			\text{is uniformly bounded in }L^1_{\mathrm{loc}}\left(dtdx;L^1\left(\left(1+|v|^2\right)Mdv\right)\right),
		\end{gathered}
	\end{equation}
	where $n_\eps = \rho_\eps^+ - \rho_\eps^-$, and
	\begin{equation*}
		\begin{gathered}
			j_\eps=\frac\delta\eps\left(u_\eps^+-u_\eps^-\right)
			\qquad\text{and}\qquad
			w_\eps=\frac\delta\eps\left(\theta_\eps^+-\theta_\eps^-\right)
			\\
			\text{are uniformly bounded in }L^1_{\mathrm{loc}}\left(dtdx\right).
		\end{gathered}
	\end{equation*}
\end{lem}

\begin{proof}
	According to the decomposition \eqref{fluct-decomposition}, it is readily seen that
	\begin{equation*}
		h_\eps = \hat h_\eps + \frac{\delta}{4}\left[\left|\hat g^+_\eps\right|^2-\left|\hat g^-_\eps\right|^2
		+ \int_{\mathbb{R}^3}\left(\left|\hat g^+_\eps\right|^2-\left|\hat g^-_\eps\right|^2\right)Mdv\right],
	\end{equation*}
	whence, by virtue of Lemma \ref{v2-int},
	\begin{equation*}
		h_\eps = \hat h_\eps + O(\delta)_{L^1_\mathrm{loc}\left(dtdx ; L^1\left((1+|v|^2)Mdv\right)\right)},
	\end{equation*}
	which establishes the uniform bound on $h_\eps$, thanks to the uniform bound on $\hat h_\eps$ from Lemma \ref{relaxation2-control}.
	
	Finally, integrating the above decomposition against $vMdv$ and $\left(\frac{|v|^2}{3}-1\right)Mdv$ clearly yields
	\begin{equation*}
		j_\eps = \hat j_\eps + O(\delta)_{L^1_\mathrm{loc}\left(dtdx\right)}
		\qquad\text{and}\qquad
		w_\eps = \hat w_\eps + O(\delta)_{L^1_\mathrm{loc}\left(dtdx\right)},
	\end{equation*}
	which, employing the uniform bounds on $\hat j_\eps$ and $\hat w_\eps$ from Lemma \ref{relaxation2-control}, concludes the justification of the lemma.
\end{proof}

\begin{lem}\label{weak compactness h}
	Let $f_\eps^+(t,x,v)$ and $f_\eps^-(t,x,v)$ be two families of measurable, almost everywhere non-negative distribution functions such that, for all $t\geq 0$,
	\begin{equation*}
		\begin{aligned}
			\frac1{\eps^2} H\left(f_\eps^{+}\right)
			& + \frac1{\eps^2} H\left(f_\eps^{-}\right)
			\\
			& +\frac{1}{\epsilon^4}\int_0^t\int_{\mathbb{R}^3}\left(D\left(f_\eps^+\right)+D\left(f_\eps^-\right)
			+ \delta^2 D\left(f_\eps^+,f_\eps^-\right)\right)(s) dx ds
			\leq C^\mathrm{in}.
		\end{aligned}
	\end{equation*}
	
	Then, as $\eps\rightarrow 0$, in the case of weak interspecies interactions, i.e.\ when $\delta = o(1)$ and $\frac\delta\eps$ is unbounded, any subsequences of fluctuations $g_\eps^\pm$ and renormalized fluctuations $\hat g_\eps^\pm$ satisfy that
	\begin{equation*}
		\begin{gathered}
			h_\eps = \frac{\delta}{\eps}\left[\left(g_\eps^+-g_\eps^-\right) - n_\eps\right]
			\qquad\text{and}\qquad
			\hat h_\eps = \frac{\delta}{\eps}\left[\left(\hat g_\eps^+-\hat g_\eps^-\right) - \hat n_\eps\right]
			\\
			\text{are weakly relatively compact in }
			L^1_{\mathrm{loc}}\left(dtdx;L^1\left(\left(1+|v|^2\right)Mdv\right)\right),
		\end{gathered}
	\end{equation*}
	and
	\begin{equation*}
		\begin{gathered}
			\begin{aligned}
				j_\eps & =\frac\delta\eps\left(u_\eps^+-u_\eps^-\right),
				&
				w_\eps & =\frac\delta\eps\left(\theta_\eps^+-\theta_\eps^-\right),
				\\
				\hat j_\eps & =\frac\delta\eps\left(\hat u_\eps^+-\hat u_\eps^-\right),
				&
				\hat w_\eps & =\frac\delta\eps\left(\hat \theta_\eps^+-\hat\theta_\eps^-\right),
			\end{aligned}
			\\
			\text{are weakly relatively compact in }L^1_{\mathrm{loc}}\left(dtdx\right).
		\end{gathered}
	\end{equation*}
	
	Moreover, as $\eps\rightarrow 0$, in the case of strong interspecies interactions, i.e.\ when $\delta = 1$, any subsequences of fluctuations $g_\eps^\pm$ and renormalized fluctuations $\hat g_\eps^\pm$ satisfy that
	% \begin{equation*}
	% 	\begin{gathered}
	% 		\frac{\hat h_\eps}{1+\left\| \hat g_\eps^+-\hat g_\eps^- \right\|_{L^2(Mdv)}}
	% 		=
	% 		\frac{\left(\hat g_\eps^+-\hat g_\eps^-\right) - \hat n_\eps}
	% 		{\eps\left(1+\left\| \hat g_\eps^+-\hat g_\eps^- \right\|_{L^2(Mdv)}\right)}
	% 		\\
	% 		\text{is uniformly bounded in }L^2_{\mathrm{loc}}\left(dt;L^2\left(Mdxdv\right)\right),
	% 	\end{gathered}
	% \end{equation*}
	\begin{equation*}
		\frac{\hat h_\eps}{1+\left\| \hat g_\eps^+-\hat g_\eps^- \right\|_{L^2(Mdv)}}
		\text{ is uniformly bounded in }L^2_{\mathrm{loc}}\left(dt;L^2\left(Mdxdv\right)\right),
	\end{equation*}
	and
	% \begin{equation*}
	% 	\begin{gathered}
	% 		\frac{\hat j_\eps}{1+\left\| \hat g_\eps^+-\hat g_\eps^- \right\|_{L^2(Mdv)}}
	% 		=
	% 		\frac{ \hat u_\eps^+-\hat u_\eps^- }
	% 		{\eps\left(1+\left\| \hat g_\eps^+-\hat g_\eps^- \right\|_{L^2(Mdv)}\right)},
	% 		\\
	% 		\frac{\hat w_\eps}{1+\left\| \hat g_\eps^+-\hat g_\eps^- \right\|_{L^2(Mdv)}}
	% 		=
	% 		\frac{\hat \theta_\eps^+-\hat \theta_\eps^- }
	% 		{\eps\left(1+\left\| \hat g_\eps^+-\hat g_\eps^- \right\|_{L^2(Mdv)}\right)}
	% 		\\
	% 		\text{are uniformly bounded in }L^2_{\mathrm{loc}}\left(dt;L^2\left(dx\right)\right),
	% 	\end{gathered}
	% \end{equation*}
	\begin{equation*}
		\begin{gathered}
			\frac{\hat j_\eps}{1+\left\| \hat g_\eps^+-\hat g_\eps^- \right\|_{L^2(Mdv)}}
			\qquad\text{and}\qquad
			\frac{\hat w_\eps}{1+\left\| \hat g_\eps^+-\hat g_\eps^- \right\|_{L^2(Mdv)}}
			\\
			\text{are uniformly bounded in }L^2_{\mathrm{loc}}\left(dt;L^2\left(dx\right)\right),
		\end{gathered}
	\end{equation*}
	while
	% \begin{equation*}
	% 	\begin{gathered}
	% 		\frac{h_\eps}{1+\left\| \hat g_\eps^+- \hat g_\eps^- \right\|_{L^2(Mdv)}}
	% 		=
	% 		\frac{\left(g_\eps^+-g_\eps^-\right) - n_\eps}
	% 		{\eps\left(1+\left\| \hat g_\eps^+- \hat g_\eps^- \right\|_{L^2(Mdv)}\right)}
	% 		\\
	% 		\text{is weakly relatively compact in }L^1_{\mathrm{loc}}\left(dtdx;L^1\left(\left(1+|v|\right)Mdv\right)\right),
	% 	\end{gathered}
	% \end{equation*}
	\begin{equation*}
		\begin{gathered}
			\frac{h_\eps}{1+\left\| \hat g_\eps^+- \hat g_\eps^- \right\|_{L^2(Mdv)}}
			\text{ is}
			\\
			\text{uniformly bounded in }L^2_{\mathrm{loc}}\left(dtdx;L^1((1+|v|)Mdv)\right)
			\\
			\text{and weakly relatively compact in }L^1_{\mathrm{loc}}\left(dtdx;L^1\left(\left(1+|v|\right)Mdv\right)\right),
		\end{gathered}
	\end{equation*}
	and
	% \begin{equation*}
	% 	\begin{gathered}
	% 		\frac{j_\eps}{1+\left\| \hat g_\eps^+- \hat g_\eps^- \right\|_{L^2(Mdv)}}
	% 		=
	% 		\frac{ u_\eps^+- u_\eps^- }
	% 		{\eps\left(1+\left\| \hat g_\eps^+- \hat g_\eps^- \right\|_{L^2(Mdv)}\right)},
	% 		\\
	% 		\text{is uniformly bounded in }L^2_{\mathrm{loc}}\left(dtdx\right).
	% 	\end{gathered}
	% \end{equation*}
	\begin{equation*}
		\frac{j_\eps}{1+\left\| \hat g_\eps^+- \hat g_\eps^- \right\|_{L^2(Mdv)}}
		\text{ is uniformly bounded in }L^2_{\mathrm{loc}}\left(dtdx\right).
	\end{equation*}
\end{lem}

\begin{proof}
	We handle the case $\delta=o(1)$ first.
	
	We have already established the uniform boundedness of $\hat h_\eps$, $\hat j_\eps$ and $\hat w_\eps$ in Lemma \ref{relaxation2-control}, while the uniform boundedness of $h_\eps$, $j_\eps$ and $w_\eps$ comes from Lemma \ref{bound hjw}. Moreover, the tightness in $v$ of $\hat h_\eps$ is easily deduced from the bound \eqref{relaxation estimate} from Lemma \ref{relaxation2-control}. Therefore, according to the Dunford-Pettis compactness criterion (see \cite{royden}), it suffices to show that $h_\eps$, $j_\eps$, $w_\eps$, $\hat h_\eps$, $\hat j_\eps$ and $\hat w_\eps$ are uniformly integrable in all variables and that $h_\eps$ is tight in $v$.

	We deal with $\hat h_\eps$, $\hat j_\eps$ and $\hat w_\eps$ first. To this end, simply notice Lemma \ref{relaxation2-control} provides the control
	\begin{equation*}
			\left|\hat j_\eps \right| +
			\left|\hat w_\eps \right|
			\leq
			C\left\|\hat h_\eps\right\|_{L^2(Mdv)}
			\leq
			O(\delta)
			\left\|
			\begin{pmatrix} \hat g_\eps^+ \\ \hat g_\eps^- \end{pmatrix}
			\right\|_{L^2(Mdv)}^2
			+O\left(1\right)_{L^2_{\mathrm{loc}}\left(dt;L^2\left(dx\right)\right)},
	\end{equation*}
	whence, by Lemma \ref{L2-lem},
	\begin{equation*}
			\left|\hat j_\eps \right| +
			\left|\hat w_\eps \right|
			\leq
			C\left\|\hat h_\eps\right\|_{L^2(Mdv)}
			\leq
			O(\delta)_{L^\infty\left(dt;L^1(dx)\right)}
			+O\left(1\right)_{L^2_{\mathrm{loc}}\left(dt;L^2\left(dx\right)\right)},
	\end{equation*}
	which establishes the equi-integrability of $\hat j_\eps$, $\hat w_\eps$ and $\left\|\hat h_\eps\right\|_{L^2(Mdv)}$ in $t$ and $x$. Furthermore, since $\hat h_\eps$ is clearly equi-integrable in $v$ thanks to the bound \eqref{relaxation estimate} from Lemma \ref{relaxation2-control}, a direct application of Lemma 5.2 from \cite{golse3} yields that $\hat h_\eps$ is equi-integrable in all variables.

	Next, we deduce the relative weak compactness of $h_\eps$, $j_\eps$ and $w_\eps$ from the relative weak compactness of $\hat h_\eps$, $\hat j_\eps$ and $\hat w_\eps$ employing the decomposition \eqref{fluct-decomposition}, which clearly yields
	\begin{equation*}
		h_\eps = \hat h_\eps + \frac{\delta}{4}\left[\left|\hat g^+_\eps\right|^2-\left|\hat g^-_\eps\right|^2
		- \int_{\mathbb{R}^3}\left(\left|\hat g^+_\eps\right|^2-\left|\hat g^-_\eps\right|^2\right)Mdv\right].
	\end{equation*}
	Therefore, since $\delta = o(1)$, it is readily seen, by virtue of the uniform integrability of $\hat j_\eps$, $\hat w_\eps$ and $\hat h_\eps$ and the uniform boundedness of $\hat g_\eps^\pm$ in $L^2_\mathrm{loc}\left(dtdx;L^2\left(\left(1+|v|^2\right)Mdv\right)\right)$ from Lemma \ref{v2-int}, that $j_\eps$, $w_\eps$ and $h_\eps$ are uniformly integrable in all variables, as well, and that $h_\eps$ is tight in $v$ in $L^1_\mathrm{loc}\left(dtdx;L^1\left(\left(1+|v|^2\right)Mdv\right)\right)$.
	
	\bigskip
	
	We turn now to the case $\delta=1$.
	
	It is readily seen that the estimate \eqref{relaxation2-control 2} from Lemma \ref{relaxation2-control} provides the refined control
	\begin{equation*}
		\left\|\hat h_\eps \right\|_{L^2(Mdv)}
		\leq
		C \left\|\hat g_\eps^+ - \hat g_\eps^-\right\|_{L^2(Mdv)}
		\left\| \hat g_\eps^\pm \right\|_{L^2(Mdv)}
		+O\left(1\right)_{L^2_{\mathrm{loc}}\left(dt;L^2\left(dx\right)\right)}.
	\end{equation*}
	In particular, in view of the boundedness of $\hat g_\eps^\pm$ in $L^\infty\left(dt;L^2\left(Mdxdv\right)\right)$ from Lemma \ref{L2-lem}, it follows that
	\begin{equation*}
			\frac{\hat h_\eps}{1+\left\| \hat g_\eps^+-\hat g_\eps^- \right\|_{L^2(Mdv)}}
			\qquad
			\text{is uniformly bounded in }L^2_{\mathrm{loc}}\left(dt;L^2\left(Mdxdv\right)\right),
	\end{equation*}
	and, incidentally, that
	\begin{equation*}
			\frac{\left|\hat j_\eps \right| +
			\left|\hat w_\eps \right|}{1+\left\| \hat g_\eps^+-\hat g_\eps^- \right\|_{L^2(Mdv)}}
			\qquad
			\text{is uniformly bounded in }L^2_{\mathrm{loc}}\left(dt;L^2\left(dx\right)\right).
	\end{equation*}

	Next, we deduce the relative weak compactness of $h_\eps$ and the uniform bound on $j_\eps$ from the uniform bound on $\hat h_\eps$ employing the decomposition \eqref{fluct-decomposition}, which clearly yields
	\begin{equation}\label{decomp h}
		h_\eps = \hat h_\eps + \frac{1}{4}\left[\left(\hat g^+_\eps-\hat g^-_\eps\right)\left(\hat g^+_\eps+\hat g^-_\eps\right)
		- \int_{\mathbb{R}^3}\left(\hat g^+_\eps-\hat g^-_\eps\right)\left(\hat g^+_\eps+\hat g^-_\eps\right)Mdv\right].
	\end{equation}
	Therefore, it is readily seen that
	\begin{equation*}
		\begin{aligned}
			\left\|h_\eps\right\|_{L^1\left((1+|v|)Mdv\right)} & \leq C \left\|\hat h_\eps\right\|_{L^2\left(Mdv\right)}
			\\
			& + C \left\|\hat g^+_\eps-\hat g^-_\eps\right\|_{L^2\left(Mdv\right)}
			\left\|\hat g^+_\eps+\hat g^-_\eps\right\|_{L^2\left((1+|v|^2)Mdv\right)},
		\end{aligned}
	\end{equation*}
	whence
	\begin{equation*}
			\frac{h_\eps}{1+\left\| \hat g_\eps^+-\hat g_\eps^- \right\|_{L^2(Mdv)}}
			\qquad
			\text{is uniformly bounded in }L^2_{\mathrm{loc}}\left(dtdx;L^1((1+|v|)Mdv)\right),
	\end{equation*}
	and, incidentally,
	\begin{equation*}
			\frac{ j_\eps}{1+\left\| \hat g_\eps^+-\hat g_\eps^- \right\|_{L^2(Mdv)}}
			\qquad
			\text{is uniformly bounded in }L^2_{\mathrm{loc}}\left(dtdx\right).
	\end{equation*}

		Now, it is readily seen from the uniform bound on $\hat h_\eps$ in $L^1_{\mathrm{loc}}\left(dtdx;L^2(Mdv)\right)$ established in Lemma \ref{relaxation2-control} and from the equi-integrability in $v$ of the families $\left|\hat g_\eps^\pm\right|^2$ established in Lemma \ref{v2-int} that the decomposition \eqref{decomp h} yields that $h_\eps$ is equi-integrable in $v$, as well. Consequently, a direct application of Lemma 5.2 from \cite{golse3} yields that the family $\frac{h_\eps}{1+\left\| \hat g_\eps^+-\hat g_\eps^- \right\|_{L^2(Mdv)}}$ is equi-integrable in all variables, which, according to the Dunford-Pettis compactness criterion (see \cite{royden}), implies its weak relative compactness in $L^1_{\mathrm{loc}}\left(dtdx;L^1(Mdv)\right)$. Finally, further using the uniform bound \eqref{def h} from Lemma \ref{bound hjw}, we deduce that the family $\frac{h_\eps}{1+\left\| \hat g_\eps^+-\hat g_\eps^- \right\|_{L^2(Mdv)}}$ is weakly relatively compact in $L^1_{\mathrm{loc}}\left(dtdx;L^1\left(\left(1+|v|\right)Mdv\right)\right)$, which concludes the proof of the lemma.
\end{proof}

%% file: constraint0.tex
\chapter[Lower order linear constraint equations and energy\ldots]{Lower order linear constraint equations and energy inequalities}\label{constraints proof}

In the preceding chapter, we have established uniform estimates and controls on the fluctuations and collision integrands by analyzing the relative entropy and entropy dissipation bounds. At this stage, we have now all the necessary tools to derive the asymptotic lower order linear constraint equations and energy inequalities from Theorems \ref{NS-WEAKCV} and \ref{CV-OMHD}. This first part of the rigorous convergence proofs is therefore very similar for both theorems.

The derivation of higher order and nonlinear constraint equations --~in particular, constraints pertaining to Theorem \ref{CV-OMHDSTRONG}~-- is performed in Chapter \ref{high constraints proof} and will require more advanced methods and refined properties on the fluctuations. More precisely, strong compactness and nonlinear weak compactness properties of the fluctuations, established later on in Chapter \ref{hypoellipticity}, will allow us to obtain the remaining constraint equations such as Ohm's law.

% ===========================
% = macroscopic constraints =
% ===========================

\section{Macroscopic constraint equations for one species}\label{macro constraint}

The macroscopic constraint equations are obtained by integrating the limiting kinetic equation against the collision invariants. In the simplest case of a one species plasma, taking limits in the kinetic equation is straightforward once we introduce the suitable renormalization.

\begin{prop}\label{weak-comp}
	Let $\left(f_\eps, E_\eps, B_\eps\right)$ be the sequence of renormalized solutions to the scaled one species Vlasov-Maxwell-Boltzmann system \eqref{VMB1} considered in Theorem \ref{NS-WEAKCV}. In accordance with Lemmas \ref{L1-lem}, \ref{L2-lem} and \ref{L2-qlem}, denote by
	\begin{equation*}
		\begin{gathered}
			g\in L^\infty\left(dt;L^2\left(Mdxdv\right)\right),
			\qquad
			q \in L^2\left(MM_*dtdxdvdv_*d\sigma\right),
			\\
			\text{and}\qquad
			E,B\in L^\infty\left(dt;L^2\left(dx\right)\right),
		\end{gathered}
	\end{equation*}
	any joint limit points of the families $\hat g_\eps$ and $\hat q_\eps$ defined by \eqref{hatg} and \eqref{hatq-def}, $E_\eps$ and $B_\eps$, respectively.
	
	Then, one has
	\begin{equation}\label{q phi psi}
		\int_{\mathbb{R}^3\times\mathbb{S}^2} q M_*dv_*d\sigma
		= \phi:\nabla_x u +\psi\cdot \nabla_x \theta,
	 \end{equation}
	where $u$ and $\theta $ are, respectively, the bulk velocity and temperature associated with the limiting fluctuation $g$, and $\phi$ and $\psi$ are the kinetic fluxes defined by \eqref{phi-psi-def}. Furthermore, $\rho$, $u$, $\theta$ and $E$ satisfy the following constraints
	\begin{equation}\label{constraints1}
		\Div u = 0, \qquad \nabla_x\left(\rho+\theta\right)-E = 0,
	\end{equation}
	where $\rho$ is the density associated with the limiting fluctuation $g$.
\end{prop}

\begin{proof}
	We start from some square root renormalization of the scaled Vlasov-Boltzmann equation \eqref{VMB1}. More precisely, we choose the admissible renormalization
	\begin{equation*}
		\beta(z)=\frac{\sqrt{z+\eps^a}-1}{\eps},
	\end{equation*}
	for some given $1<a<4$, which yields, using the decomposition of collision integrands \eqref{integrands-decomposition},
	\begin{equation}\label{renormalized}
		\begin{aligned}
			( \eps \d_t + & v \cdot \nabla_x + \eps \left(E_\eps+ v\wedge B_\eps\right)\cdot \nabla_v )
			{\sqrt{ G_\eps+\eps^a }-1 \over \eps}
			- E_\eps \cdot v {G_\eps \over 2\sqrt{ G_\eps+\eps^a }} \\
			& =
			{\sqrt{G_\eps} \over 2\sqrt{ G_\eps+\eps^a }}
			\int_{\mathbb{R}^3\times\mathbb{S}^2} \sqrt{ G_{\eps *} } \hat q_\eps M_*dv_*d\sigma
			+
			{\eps^2 \over 8\sqrt{G_\eps+\eps^a  }}
			\int_{\mathbb{R}^3\times\mathbb{S}^2} \hat q_\eps^2 M_*dv_*d\sigma
			\\
			& \eqdefa Q_\eps^1+Q_\eps^2 .
		\end{aligned}
	\end{equation}

	Then, thanks to Lemma \ref{L2-lem}, it holds that
	\begin{equation*}
		\sqrt{G_{\eps *}} = 1 + O(\eps)_{L^2_\mathrm{loc}\left(dt ; L^2\left(M_* dxdv_*\right)\right)},
	\end{equation*}
	whence, employing the uniform bound $\hat q_\eps\in L^2\left(MM_*dtdxdvdv_*d\sigma\right)$ from Lemma \ref{L2-qlem},
	\begin{equation}\label{Q 1}
		\begin{aligned}
			Q_\eps^1 & =
			{\sqrt{G_\eps} \over 2\sqrt{ G_\eps+\eps^a }}
			\int_{\mathbb{R}^3\times\mathbb{S}^2} \hat q_\eps M_*dv_*d\sigma
			+ O(\eps)_{L^1_\mathrm{loc}\left(dtdx ; L^2\left(Mdv\right)\right)},
			\\
			Q_\eps^2 & = O\left(\eps^{2-\frac a2}\right)_{L^1\left(Mdtdxdv\right)}.
		\end{aligned}
	\end{equation}

	Next, since, decomposing according to the tails of $G_\eps$,
	\begin{equation*}
		\begin{aligned}
			\left|\frac{\sqrt{G_\eps+\eps^a}-1}{\eps}-\frac{\sqrt{G_\eps}-1}{\eps}\right|
			& =
			\frac{\eps^{a-1}}{\sqrt{\eps^a+G_\eps}+\sqrt{G_\eps}}
			\\
			& \leq
			\frac{\eps^{a-1}}{\sqrt{\eps^a+G_\eps}+\sqrt{G_\eps}}
			\mathds{1}_{\left\{G_\eps>\frac 12\right\}}
			+
			\eps^{\frac a2-1}\mathds{1}_{\left\{G_\eps\leq\frac 12\right\}}
			\\
			& \leq
			O\left(\eps^{a-1}\right)_{L^\infty(dtdxdv)}
			+\eps^{\frac a2}\left(2+\sqrt{2}\right) \left|\frac{\sqrt{G_\eps}-1}{\eps}\right|,
		\end{aligned}
	\end{equation*}
	one proves, by virtue of Lemma \ref{L2-lem} or Lemma \ref{v2-int}, that
	\begin{equation}\label{error-alpha}
		\begin{aligned}
		2{\sqrt{G_\eps+\eps^a }-1\over \eps}-\hat g_\eps
		& = O\left(\eps^{\frac a2 -1}\right)_{L^\infty(dtdxdv)},
		\\
		2{\sqrt{G_\eps+\eps^a }-1\over \eps}-\hat g_\eps
		& = O\left(\eps^{a-1}\right)_{L^\infty(dtdxdv)} +O\left(\eps^{\frac a2}\right)_{L^\infty\left(dt;L^2\left(M dxdv\right)\right)},
		\\
		2{\sqrt{G_\eps+\eps^a }-1\over \eps}-\hat g_\eps
		& = O\left(\eps^{a-1}\right)_{L^\infty(dtdxdv)} +O\left(\eps^{\frac a2}\right)_{L^2_\mathrm{loc}\left(dtdx;L^2\left(\left(1+|v|^2\right)Mdv\right)\right)},
		\end{aligned}
	\end{equation}
	and
	\begin{equation*}% \label{error-alpha 2}
		\begin{aligned}
			E_\eps \cdot v & {G_\eps \over \sqrt{G_\eps+\eps^a }} \\
			& =  E_\eps \cdot v  \left(1+\eps {\sqrt{G_\eps+\eps^a }-1\over \eps}\right)
			- E_\eps \cdot v {\eps^a \over \sqrt{G_\eps+\eps^a }} \\
			& = E_\eps \cdot v +O\left( \eps+\eps^{\frac a2+1} \right)_{L^\infty\left(dt;L^1\left(M dxdv\right)\right)}
			+O\left(\eps^{a}+\eps^{\frac a2}\right)_{L^\infty\left(dt;L^2\left(M dxdv\right)\right)}.
		\end{aligned}
	\end{equation*}

	In particular, employing \eqref{error-alpha} to deduce that
	\begin{equation*}
		\begin{aligned}
			\left|{\sqrt{G_\eps} \over \sqrt{ G_\eps+\eps^a }} - 1\right|
			& =
			{\eps \over \sqrt{ G_\eps+\eps^a }}
			\left|{\sqrt{ G_\eps+\eps^a } - 1 \over \eps} - \frac 12 \hat g_\eps\right|
			\\
			& \leq
			\eps^{1-\frac a2}
			\left|{\sqrt{ G_\eps+\eps^a } - 1 \over \eps} - \frac 12 \hat g_\eps\right|
			\\
			& \leq
			O\left(\eps^{\frac a2}\right)_{L^\infty(dtdxdv)} +O\left(\eps\right)_{L^2_\mathrm{loc}\left(dtdx;L^2\left(\left(1+|v|^2\right)Mdv\right)\right)},
		\end{aligned}
	\end{equation*}
	we obtain the following refinement of \eqref{Q 1}
	\begin{equation*}% \label{Q 2}
		\begin{aligned}
			Q_\eps^1 & =
			\frac 12
			\int_{\mathbb{R}^3\times\mathbb{S}^2} \hat q_\eps M_*dv_*d\sigma
			+ O\left(\eps^{\frac a2}\right)_{L^2(Mdtdxdv)}
			+ O(\eps)_{L^1_\mathrm{loc}\left(dtdx ; L^1\left(\left(1+|v|\right)Mdv\right)\right)},
			\\
			Q_\eps^2 & = O\left(\eps^{2-\frac a2}\right)_{L^1\left(Mdtdxdv\right)}.
		\end{aligned}
	\end{equation*}

	Therefore, taking weak limits in \eqref{renormalized} leads to
	\begin{equation}\label{kinetic equation 1}
		v\cdot \nabla_x g-E\cdot v =\int_{\mathbb{R}^3\times\mathbb{S}^2} q M_* dv_* d\sigma,
	\end{equation}
	which, together with the fact, according to Lemma \ref{relaxation-control}, that $g$ is an infinitesimal Maxwellian, provides that
	\begin{equation*}
		\begin{aligned}
			\int_{\mathbb{R}^3\times\mathbb{S}^2} q M_*dv_*d\sigma
			& = \Div\left((\rho+\theta)v+\frac{|v|^2}{3}u+\phi u+\theta\psi\right)
			- E \cdot v \\
			& = \left( \phi:\nabla_x u +\psi\cdot \nabla_x \theta \right)
			+ \left(\nabla_x(\rho+\theta)-E\right)\cdot v
			+\frac13 \left(\Div u\right)|v|^2.
		 \end{aligned}
	 \end{equation*}
	Then, remarking that $q$ inherits the collisional symmetries of $q_\eps$ and $\hat q_\eps$, we get
	\begin{equation*}
		\int_{\mathbb{R}^3\times\mathbb{R}^3\times\mathbb{S}^2}
		q
		\begin{pmatrix}
			1 \\ v \\ \frac{|v|^2}{2}
		\end{pmatrix}
		MM_* dvdv_*d\sigma = 0,
	\end{equation*}
	so that, since $\phi(v)$ and $\psi(v)$ are orthogonal to the collisional invariants, the constraints \eqref{constraints1} hold.
	
	The proof of the proposition is complete.
\end{proof}

\section[Macroscopic constraint equations for two species\ldots]{Macroscopic constraint equations for two species, weak interactions}\label{macro constraint 2}

In the case of a two species plasma, the renormalization process is more complicated because there are two different distributions. Nevertheless, for weak interspecies interactions, i.e.\ $\delta=o(1)$ and $\frac\delta\eps$ unbounded, we have a result quite similar to the preceding proposition. As for strong interspecies interactions, i.e.\ $\delta=1$, even the lowest order constraints will require the dealing with nonlinear terms and, therefore, will be handled with more advanced techniques in Chapter \ref{high constraints proof} (see Propositions \ref{high weak-comp2} and \ref{high weak-comp3}).

\begin{prop}\label{weak-comp2}
	Let $\left(f_\eps^\pm, E_\eps, B_\eps\right)$ be the sequence of renormalized solutions to the scaled two species Vlasov-Maxwell-Boltzmann system \eqref{VMB2} considered in Theorem \ref{CV-OMHD} for weak interspecies interactions, i.e.\ $\delta=o(1)$ and $\frac\delta\eps$ unbounded. In accordance with Lemmas \ref{L1-lem}, \ref{L2-lem} and \ref{L2-qlem}, denote by
		\begin{equation*}
			\begin{gathered}
				g^\pm\in L^\infty\left(dt;L^2\left(Mdxdv\right)\right),
				\qquad
				q^\pm, q^{\pm,\mp} \in L^2\left(MM_*dtdxdvdv_*d\sigma\right),
				\\
				\text{and}\qquad
				E,B\in L^\infty\left(dt;L^2\left(dx\right)\right),
			\end{gathered}
		\end{equation*}
	any joint limit points of the families $\hat g_\eps^\pm$, $\hat q_\eps^\pm$ and $\hat q_\eps^{\pm,\mp}$ defined by \eqref{hatg} and \eqref{hatq-def}, $E_\eps$ and $B_\eps$, respectively.
	
	Then, one has
	\begin{equation}\label{q phi psi 2}
		\int_{\mathbb{R}^3\times\mathbb{S}^2} q^\pm M_*dv_*d\sigma
		= \phi:\nabla_x u +\psi\cdot \nabla_x \theta ,
	\end{equation}
	where $u$ and $\theta $ are, respectively, the bulk velocity and temperature associated with the limiting fluctuations $g^\pm$, and $\phi$ and $\psi$ are the kinetic fluxes defined by \eqref{phi-psi-def}. Furthermore, $\rho^\pm$, $u$ and $\theta$ satisfy the following constraints
	\begin{equation}\label{constraints2-alpha}
		\Div u = 0, \qquad \nabla_x\left(\rho^\pm+\theta\right) = 0,
	\end{equation}
	where $\rho^\pm$ are the densities respectively associated with the limiting fluctuations $g^\pm$. In particular, the strong Boussinesq relation $\rho^\pm+\theta=0$ holds and, moreover, since $\nabla_x\left(\rho^+-\rho^-\right)=0$, it also holds that $\rho^+=\rho^-$.
\end{prop}

\begin{proof}
	We start from some square root renormalization of the scaled Vlasov-Boltzmann equation \eqref{VMB2}. More precisely, we choose, as previously in the proof of Proposition \ref{weak-comp}, the admissible renormalization
	\begin{equation*}
		\beta(z)=\frac{\sqrt{z+\eps^a}-1}{\eps},
	\end{equation*}
	for some given $1<a<4$, which yields, recalling the definitions \eqref{hatq-def} of renormalized collision integrands and using the decomposition of collision integrands \eqref{integrands-decomposition},
	\begin{equation}\label{renormalized2}
		\begin{aligned}
			( \eps \d_t + v \cdot \nabla_x \pm \delta \left(\eps E_\eps+ v\wedge B_\eps\right) & \cdot \nabla_v )
			{\sqrt{ G_\eps^\pm +\eps^a }-1 \over \eps}
			\mp \delta E_\eps \cdot v {G_\eps^\pm \over 2\sqrt{ G_\eps^\pm +\eps^a }} \\
			& =
			{\sqrt{G_\eps^\pm} \over 2\sqrt{ G_\eps^\pm+\eps^a }}
			\int_{\mathbb{R}^3\times\mathbb{S}^2} \sqrt{ G_{\eps *}^\pm } \hat q_\eps^\pm M_*dv_*d\sigma
			\\
			& +
			{\eps^2 \over 8\sqrt{G_\eps^\pm+\eps^a  }}
			\int_{\mathbb{R}^3\times\mathbb{S}^2} \left(\hat q_\eps^\pm\right)^2 M_*dv_*d\sigma
			\\
			& +
			{\delta\sqrt{G_\eps^\pm} \over 2\sqrt{ G_\eps^\pm+\eps^a }}
			\int_{\mathbb{R}^3\times\mathbb{S}^2} \sqrt{ G_{\eps *}^\mp } \hat q_\eps^{\pm,\mp} M_*dv_*d\sigma
			\\
			& +
			{\eps^2 \over 8\sqrt{G_\eps^\pm+\eps^a  }}
			\int_{\mathbb{R}^3\times\mathbb{S}^2} \left(\hat q_\eps^{\pm,\mp}\right)^2 M_*dv_*d\sigma.
		\end{aligned}
	\end{equation}
	% \begin{equation}
	% \begin{array}l
	% \ds(\eps \d_t +v\cdot \nabla_x + \eps^{\frac{1-q}2}(\eps E_\eps+  v\wedge B_\eps)\cdot \nabla_v ){\sqrt{G^+_\eps+\eps^a }-1\over \eps} - \eps^{\frac{1-q}2} E_\eps \cdot v {G^+_\eps \over 2\sqrt{G^+_\eps+\eps^a }}\\
	% \ds= \frac1{2\eps^2} {1\over \sqrt{G_\eps^++\eps^a  }} \iint \left( \hat q_\eps ^+\right)^2  d\sigma dv_*+\frac1{\eps^2} {\sqrt{G^+_\eps}\over \sqrt{G^+_\eps+\eps^a  }} \iint \hat q_\eps^+\sqrt{ G^+_{\eps *}}  d\sigma M_* dv_*\\
	% \ds +  \frac1{2\eps^{1+q}} {1\over \sqrt{G_\eps^++\eps^a  }} \iint \left( \hat q_\eps ^{+,-}\right)^2  d\sigma dv_*+\frac1{\eps^{1+q}} {\sqrt{G^+_\eps}\over \sqrt{G^+_\eps+\eps^a  }} \iint \hat q_\eps^{+,-} \sqrt{ G^+_{\eps *}}  d\sigma M_* dv_*\\
	% \ds(\eps \d_t +v\cdot \nabla_x - \eps^{\frac{1-q}2} (\eps E_\eps+  v\wedge B_\eps)\cdot \nabla_v ){\sqrt{G^-_\eps+\eps^a }-1\over \eps} + \eps^{\frac{1-q}2} E_\eps \cdot v {G^-_\eps \over 2\sqrt{G^-_\eps+\eps^a }}\\
	% \ds= \frac1{2\eps^2} {1\over \sqrt{G_\eps^-+\eps^a  }} \iint \left( \hat q_\eps ^-\right)^2  d\sigma dv_*+\frac1{\eps^2} {\sqrt{G^-_\eps}\over \sqrt{G^-_\eps+\eps^a  }} \iint \hat q_\eps^-\sqrt{ G^-_{\eps *}}  d\sigma M_* dv_*\\
	% \ds +  \frac1{2\eps^{1+q}} {1\over \sqrt{G_\eps^-+\eps^a  }} \iint \left( \hat q_\eps ^{-,=}\right)^2  d\sigma dv_*+\frac1{\eps^{1+q}} {\sqrt{G^-_\eps}\over \sqrt{G^-_\eps+\eps^a  }} \iint \hat q_\eps^{-,+} \sqrt{ G^-_{\eps *}}  d\sigma M_* dv_*
	% \end{array}
	% \end{equation}
	
	 The proof follows then the exact same lines as the proof of Proposition \ref{weak-comp}. In particular, we obtain without any additional difficulty the fact that the right-hand side of \eqref{renormalized2} converges weakly to
	\begin{equation*}
		% \begin{aligned}
			\frac 12 \int_{\mathbb{R}^3\times\mathbb{S}^2} q^\pm M_*dv_*d\sigma,
			% && \text{if }\delta=o(1), \\
			% & \frac 12 \int_{\mathbb{R}^3\times\mathbb{S}^2} \left(q^\pm + q^{\pm,\mp}\right) M_*dv_*d\sigma,
			% && \text{if }\delta=1,
		% \end{aligned}
	\end{equation*}
	while the renormalized densities satisfy, following \eqref{error-alpha},% and \eqref{error-alpha 2},
	\begin{equation*}% \label{error-alpha 3}
		\begin{aligned}
			2{\sqrt{G_\eps^\pm+\eps^a }-1\over \eps}-\hat g_\eps^\pm
			& = O\left(\eps^{\frac a2 -1}\right)_{L^\infty(dtdxdv)},
			\\
			2{\sqrt{G_\eps^\pm+\eps^a }-1\over \eps}-\hat g_\eps^\pm
			& = O\left(\eps^{a-1}\right)_{L^\infty(dtdxdv)} +O\left(\eps^{\frac a2}\right)_{L^\infty\left(dt;L^2\left(M dxdv\right)\right)},
			\\
			2{\sqrt{G_\eps^\pm+\eps^a }-1\over \eps}-\hat g_\eps^\pm
			& = O\left(\eps^{a-1}\right)_{L^\infty(dtdxdv)} +O\left(\eps^{\frac a2}\right)_{L^2_\mathrm{loc}\left(dtdx;L^2\left(\left(1+|v|^2\right)Mdv\right)\right)}.
		\end{aligned}
	\end{equation*}
	% and
	% \begin{equation*}% \label{error-alpha 4}
	% 	\begin{aligned}
	% 		E_\eps \cdot v & {G_\eps^\pm \over \sqrt{G_\eps^\pm+\eps^a }} \\
	% 		& =  E_\eps \cdot v  \left(1+\eps {\sqrt{G_\eps^\pm+\eps^a }-1\over \eps}\right)
	% 		- E_\eps \cdot v {\eps^a \over \sqrt{G_\eps^\pm+\eps^a }} \\
	% 		& = E_\eps \cdot v +O\left( \eps+\eps^{\frac a2+1} \right)_{L^\infty\left(dt;L^1\left(M dxdv\right)\right)}
	% 		+O\left(\eps^{a}+\eps^{\frac a2}\right)_{L^\infty\left(dt;L^2\left(M dxdv\right)\right)}.
	% 	\end{aligned}
	% \end{equation*}

	Next, using the uniform $L^\infty\left(dt;L^2(dx)\right)$ bounds on $E_\eps$ and $B_\eps$, as well as the $L^2_\mathrm{loc}\left(dtdx; L^2(Mdv)\right)$ bounds on the renormalized fluctuations ${\sqrt{G^\pm_\eps+\eps^a }-1\over \eps}$ and ${G^\pm_\eps \over 2\sqrt{G^\pm_\eps+\eps^a }}$, we easily obtain that all the terms coming from the Lorentz force in \eqref{renormalized2} vanish in the weak limit~:
		\begin{equation*}
			\pm \delta \left(\eps E_\eps+ v\wedge B_\eps\right)\cdot \nabla_v
			{\sqrt{ G_\eps^\pm +\eps^a }-1 \over \eps}
			\mp \delta E_\eps \cdot v {G_\eps^\pm \over 2\sqrt{ G_\eps^\pm +\eps^a }}
			\rightarrow 0.
		\end{equation*}
		Therefore, taking weak limits in \eqref{renormalized2} leads to
		\begin{equation}\label{kinetic equation 2}
			v\cdot \nabla_x g^\pm = \int_{\mathbb{R}^3\times\mathbb{S}^2} q^\pm M_* dv_* d\sigma,
		\end{equation}
		which, together with the fact, according to Lemmas \ref{relaxation-control} and \ref{relaxation2-control}, that $g^+$ and $g^-$ are infinitesimal Maxwellians, which differ only by their densities $\rho^+$ and $\rho^-$, provides that
		\begin{equation*}
			\begin{aligned}
				\int_{\mathbb{R}^3\times\mathbb{S}^2} q^\pm M_*dv_*d\sigma
				& = \Div\left((\rho^\pm+\theta)v+\frac{|v|^2}{3}u+\phi u+\theta\psi\right)
				\\
				& = \left( \phi:\nabla_x u +\psi\cdot \nabla_x \theta \right)
				+ \nabla_x(\rho^\pm+\theta)\cdot v
				+\frac13 \left(\Div u\right)|v|^2.
			 \end{aligned}
		 \end{equation*}
		Then, remarking that $q^\pm$ inherits the collisional symmetries of $q_\eps^\pm$ and $\hat q_\eps^\pm$, we get
		\begin{equation*}
			\int_{\mathbb{R}^3\times\mathbb{R}^3\times\mathbb{S}^2}
			q^\pm
			\begin{pmatrix}
				1 \\ v \\ \frac{|v|^2}{2}
			\end{pmatrix}
			MM_* dvdv_*d\sigma = 0,
		\end{equation*}
		so that, since $\phi(v)$ and $\psi(v)$ are orthogonal to the collisional invariants, the constraints \eqref{constraints2-alpha} hold.

	The proof of the proposition is complete.
\end{proof}

% \begin{rem}
% 	Note that, in the case $\delta=1$ of the preceding proof, we are not able to identify separately the limits of the terms coming from the Lorentz force. However, these terms are always bounded (even without the equi-integrability assumption \eqref{equiintegrability-ch4}). In particular, we have uniform bounds on the transport operator
% 	\begin{equation*}
% 		\left(\eps \d_t +v\cdot \nabla_x \pm (\eps E_\eps+  v\wedge B_\eps)\cdot \nabla_v \right){\sqrt{G^\pm_\eps+\eps^a }-1\over \eps}.
% 	\end{equation*}
% \end{rem}

\begin{prop}\label{weak-comp3}
	Let $\left(f_\eps^\pm, E_\eps, B_\eps\right)$ be the sequence of renormalized solutions to the scaled two species Vlasov-Maxwell-Boltzmann system \eqref{VMB2} considered in Theorem \ref{CV-OMHD} for weak interspecies interactions, i.e.\ $\delta=o(1)$ and $\frac\delta\eps$ unbounded. In accordance with Lemmas \ref{L1-lem}, \ref{L2-lem}, \ref{L2-qlem}, \ref{bound hjw} and \ref{weak compactness h} denote by
	\begin{equation*}
		\begin{aligned}
			g^\pm & \in L^\infty\left(dt;L^2\left(Mdxdv\right)\right),\\
			q^{\pm,\mp} & \in L^2\left(MM_*dtdxdvdv_*d\sigma\right),\\
			h & \in L^1_\mathrm{loc}\left(dtdx;L^1\left((1+|v|^2)Mdv\right)\right),
		\end{aligned}
	\end{equation*}
	any joint limit points of the families $\hat g_\eps^\pm$, $\hat q_\eps^{\pm,\mp}$ and $h_\eps$ defined by \eqref{hatg}, \eqref{hatq-def} and \eqref{def h}, respectively.
	
	Then, one has $h= j\cdot v + w\left(\frac{|v|^2}{2}-\frac 32\right)$ and
	\begin{equation}\label{mixed q phi psi}
		\pm 2 \int_{\mathbb{R}^3\times\mathbb{S}^2} q^{\pm,\mp} M_* dv_*d\sigma
		=-\mathfrak{L}(h)=
		-j\cdot \mathfrak{L} \left(v\right)
		-w\mathfrak{L} \left(\frac{|v|^2}{2}\right),
	\end{equation}
	where $j$ and $w$ are, respectively, the bulk velocity and temperature associated with the limiting fluctuation $h$, i.e.\ $j$ is the electric current and $w$ is the internal electric energy.
\end{prop}

\begin{proof}
	We start from the decomposition
	\begin{equation}\label{h decomposition}
		h_\eps = \hat h_\eps + \frac{\delta}{4}\left[\left|\hat g^+_\eps\right|^2-\left|\hat g^-_\eps\right|^2
		- \int_{\mathbb{R}^3}\left(\left|\hat g^+_\eps\right|^2-\left|\hat g^-_\eps\right|^2\right)Mdv\right],
	\end{equation}
	which follows from the decomposition \eqref{fluct-decomposition} of fluctuations. In particular, integrating \eqref{h decomposition} against $vMdv$ and $\left(\frac{|v|^2}{3}-1\right)Mdv$ yields
	\begin{equation*}
		\begin{aligned}
			j_\eps & = \hat j_\eps + \frac{\delta}{4}\int_{\mathbb{R}^3}\left(\left|\hat g^+_\eps\right|^2-\left|\hat g^-_\eps\right|^2 \right)vMdv,
			\\
			w_\eps & = \hat w_\eps + \frac{\delta}{4}\int_{\mathbb{R}^3}\left(\left|\hat g^+_\eps\right|^2-\left|\hat g^-_\eps\right|^2 \right)\left(\frac{|v|^2}{3}-1\right)Mdv.
		\end{aligned}
	\end{equation*}

	According to Lemma \ref{weak compactness h}, we consider now weakly convergent subsequences
	\begin{equation*}
		h_\eps\rightharpoonup h,\qquad
		\hat h_\eps\rightharpoonup \hat h,
	\end{equation*}
	in $L^1_\mathrm{loc}\left(dtdx;L^1\left(\left(1+|v|^2\right)Mdv\right)\right)$, and
	\begin{equation*}
		j_\eps\rightharpoonup j,\qquad
		w_\eps\rightharpoonup w, \qquad
		\hat j_\eps\rightharpoonup \hat j,\qquad
		\hat w_\eps\rightharpoonup \hat w,
	\end{equation*}
	in $L^1_\mathrm{loc}\left(dtdx\right)$. Clearly, since $\delta=o(1)$, we easily obtain, in view of the uniform $L^2_{\mathrm{loc}}\left(dtdx;L^2\left(\left(1+|v|^2\right)Mdv\right)\right)$ bound on $\hat g_\eps^\pm$ provided by Lemma \ref{v2-int}, passing to the limit in \eqref{h decomposition}, that
	\begin{equation*}
		h=\hat h, \qquad j=\hat j \qquad\text{and}\qquad w=\hat w.
	\end{equation*}
	Furthermore, using the relaxation estimate \eqref{relaxation-est} from Lemma \ref{relaxation-control}, it holds that
	\begin{equation*}
		\hat h_\eps-\Pi\hat h_\eps=
		\frac\delta\eps\left(\hat g_\eps^+-\Pi\hat g_\eps^+ - \hat g_\eps^-+\Pi\hat g_\eps^+\right)
		=O(\delta)_{L^1_\mathrm{loc}\left(dtdx;L^2\left(Mdv\right)\right)},
	\end{equation*}
	whence $\hat h=\Pi\hat h$, for $\delta$ vanishes asymptotically, and, therefore,
	\begin{equation}\label{h maxwellian}
		h=\hat h = j\cdot v + w\left(\frac{|v|^2}{2}-\frac 32\right).
	\end{equation}

	Next, it is readily seen that the elementary decompositions
	\begin{equation*}
		\begin{aligned}
			\cL \left(\hat g_\eps^\pm,\hat g_\eps^\mp\right) & =
			\frac\eps 2  \cQ\left(\hat g_\eps^\pm ,\hat g_\eps^\mp \right)-\frac 2{\eps}  \cQ \left(\sqrt{G_\eps^\pm},\sqrt{G_\eps^\mp}\right),
			\\
			\mathcal{L}\left(\hat g_\eps^\pm,\hat g_\eps^\mp\right) & =
			\pm \frac \eps{2\delta} \mathfrak{L}\left(\hat h_\eps\right)
			+ \frac 12 \mathcal{L}\left(\hat g_\eps^+ + \hat g_\eps^-\right),
		\end{aligned}
	\end{equation*}
	yield that
	\begin{equation}\label{h decomp 2}
		\begin{aligned}
			\mathfrak{L}\left(\hat h_\eps\right)
			& =
			\mp \frac\delta\eps\mathcal{L}\left(\hat g_\eps^+ + \hat g_\eps^-\right)
			\pm \delta \cQ\left(\hat g_\eps^\pm ,\hat g_\eps^\mp \right)
			\mp\frac {4\delta}{\eps^2}  \cQ \left(\sqrt{G_\eps^\pm},\sqrt{G_\eps^\mp}\right) \\
			& =
			\mp \delta\left(\mathcal{L}\left(\frac{\hat g_\eps^+-\Pi \hat g_\eps^+}{\eps}\right)
			+\mathcal{L}\left(\frac{\hat g_\eps^--\Pi \hat g_\eps^-}{\eps}\right)
			- \cQ\left(\hat g_\eps^\pm ,\hat g_\eps^\mp \right)\right)
			\\
			& \mp 2 \int_{\mathbb{R}^3\times\mathbb{S}^2} \hat q_\eps^{\pm,\mp} M_* dv_*d\sigma.
		\end{aligned}
	\end{equation}
	% \begin{equation}
	% 	\begin{aligned}
	% 		% \frac\eps 2 \cQ\left(\hat g_\eps^+-\hat g_\eps^- ,\hat h_\eps \right)
	% 		% 		+ \frac\eps 2 \cQ\left(\hat h_\eps ,\hat n_\eps \right)
	% 		\frac\delta 2 \cQ\left(\hat g_\eps^+-\hat g_\eps^- ,\hat g_\eps^+-\hat g_\eps^- \right)
	% 		% \\
	% 		% 			& = \frac\delta 2  \left[\cQ\left(\hat g_\eps^+ ,\hat g_\eps^+ \right) + \cQ\left(\hat g_\eps^- ,\hat g_\eps^- \right)
	% 		% 			- \cQ\left(\hat g_\eps^+ ,\hat g_\eps^-\right) - \cQ\left(\hat g_\eps^- ,\hat g_\eps^+\right)\right]
	% 		& =\frac {2\delta}{\eps^2} \left[\cQ \left(\sqrt{G_\eps^+},\sqrt{G_\eps^+}\right)
	% 		+ \cQ \left(\sqrt{G_\eps^-},\sqrt{G_\eps^-}\right)\right]
	% 		\\
	% 		& - \frac {2\delta}{\eps^2} \left[ \cQ \left(\sqrt{G_\eps^+},\sqrt{G_\eps^-}\right)
	% 		+ \cQ \left(\sqrt{G_\eps^-},\sqrt{G_\eps^+}\right)\right]
	% 		\\
	% 		& =
	% 		\int_{\mathbb{R}^3\times\mathbb{S}^2} \left(\delta\hat q_\eps^{+}+\delta\hat q_\eps^{-}-\hat q_\eps^{+,-}-\hat q_\eps^{-,+}\right) M_* dv_*d\sigma.
	% 	\end{aligned}
	% \end{equation}
	Therefore, passing to the limit $\eps\rightarrow 0$ in \eqref{h decomp 2}, we find, in view of the control \eqref{relaxation-est} from Lemma \ref{relaxation-control} and since the linear and quadratic collision operators are continuous on $L^2(Mdv)$ (see \eqref{Q continuous}), that
	\begin{equation*}
		\mathfrak{L} \left(h\right)
		= \mp 2
		\int_{\mathbb{R}^3\times\mathbb{S}^2} q^{\pm,\mp} M_* dv_*d\sigma.
	\end{equation*}
	Further employing the infinitesimal Maxwellian expression of $h$ from \eqref{h maxwellian}, we arrive at
	\begin{equation*}
		j\cdot\mathfrak{L} \left(v\right)
		+\frac 12 w\mathfrak{L} \left(|v|^2\right)
		= \mp 2
		\int_{\mathbb{R}^3\times\mathbb{S}^2} q^{\pm,\mp} M_* dv_*d\sigma,
	\end{equation*}
	which concludes the proof of the proposition.
\end{proof}

% =====================
% = energy inequality =
% =====================

\section{Energy inequalities}

In view of the results from Sections \ref{macro constraint} and \ref{macro constraint 2}, we are now able to establish the limiting energy inequalities for one species and for two species in the case of weak interactions only. The limiting energy inequality for strong interactions will require the results from Section \ref{high constraints 1} and, thus, will be treated later on in Section \ref{energy inequality singular}.

\begin{prop}\label{energy ineq 1}
	Let $\left(f_\eps, E_\eps, B_\eps\right)$ be the sequence of renormalized solutions to the scaled one species Vlasov-Maxwell-Boltzmann system \eqref{VMB1} considered in Theorem \ref{NS-WEAKCV}. In accordance with Lemmas \ref{L1-lem}, \ref{L2-lem} and \ref{L2-qlem}, denote by
		\begin{equation*}
			\begin{gathered}
				g\in L^\infty\left(dt;L^2\left(Mdxdv\right)\right),
				\qquad
				q \in L^2\left(MM_*dtdxdvdv_*d\sigma\right),
				\\
				\text{and}\qquad
				E,B\in L^\infty\left(dt;L^2\left(dx\right)\right),
			\end{gathered}
		\end{equation*}
	any joint limit points of the families $\hat g_\eps$ and $\hat q_\eps$ defined by \eqref{hatg} and \eqref{hatq-def}, $E_\eps$ and $B_\eps$, respectively.
	
	Then, one has the energy inequality, for almost every $t\geq 0$,
	\begin{equation*}
		\begin{aligned}
			& \frac 12\left(\left\|\rho\right\|_{L^2_x}^2 + \left\|u\right\|_{L^2_x}^2
			+ \frac 32\left\|\theta\right\|_{L^2_x}^2 + \left\|E\right\|_{L^2_x}^2
			+ \left\|B\right\|_{L^2_x}^2 \right)(t)
			\\
			& \hspace{30mm} +
			\int_0^t \left(\mu
			\left\|\nabla_x u\right\|_{L^2_x}^2
			+ \frac 52\kappa
			\left\|\nabla_x\theta\right\|_{L^2_x}^2\right)(s) ds
			\leq C^\mathrm{in},
		\end{aligned}
	\end{equation*}
	where $\rho$, $u$ and $\theta $ are, respectively, the density, bulk velocity and temperature associated with the limiting fluctuation $g$, and the viscosity $\mu>0$ and thermal conductivity $\kappa>0$ are defined by \eqref{mu kappa}.
\end{prop}

\begin{proof}
	First, by the estimate \eqref{q-est} from Lemma \ref{L2-qlem} and the weak sequential lower semi-continuity of convex functionals, we find that, for all $t\geq 0$,
	\begin{equation*}
		\begin{aligned}
			\frac 14 \int_0^t \int_{\mathbb{R}^3} \int_{\mathbb{R}^3\times\mathbb{R}^3\times\mathbb{S}^2} &
			{q}^2 MM_*dvdv_* d\sigma dx ds
			\\
			& \leq \liminf_{\eps\rightarrow 0}
			\frac 14 \int_0^t \int_{\mathbb{R}^3} \int_{\mathbb{R}^3\times\mathbb{R}^3\times\mathbb{S}^2}
			\hat{q}_\eps^2 MM_*dvdv_* d\sigma dx ds
			\\
			& \leq \liminf_{\eps\rightarrow 0} \frac{1}{\epsilon^4}\int_0^t\int_{\mathbb{R}^3}
			D\left(f_\eps\right)(s) dx ds,
		\end{aligned}
	\end{equation*}
	which, when combined with Lemma \ref{L1-lem}, yields, passing to the limit in the entropy inequality \eqref{entropy1}, for almost every $t\geq 0$,
	\begin{equation*}
		\begin{aligned}
		\frac 12\int_{\mathbb{R}^3\times\mathbb{R}^3} {g}^2(t)Mdxdv
		& + \frac 1{2} \int_{\mathbb{R}^3} \left(|E|^2+ |B|^2\right)(t) dx
		\\
		& +
		\frac 14 \int_0^t \int_{\mathbb{R}^3} \int_{\mathbb{R}^3\times\mathbb{R}^3\times\mathbb{S}^2}
		{q}^2 MM_*dvdv_* d\sigma dx ds
		\leq C^\mathrm{in}.
		\end{aligned}
	\end{equation*}
	Since, according to Lemma \ref{relaxation-control}, the limiting fluctuation $g=\rho+u\cdot v + \theta \left(\frac{|v|^2}{2}-\frac 32\right)$ is an infinitesimal Maxwellian, we easily compute that
	\begin{equation*}
		\int_{\mathbb{R}^3}{g}^2Mdv
		= \rho^2 + |u|^2 + \frac 32 \theta^2,
	\end{equation*}
	which implies
	\begin{equation}\label{entropy1 limit}
		\begin{aligned}
			& \frac 12\left(\left\|\rho\right\|_{L^2_x}^2 + \left\|u\right\|_{L^2_x}^2
			+ \frac 32\left\|\theta\right\|_{L^2_x}^2 + \left\|E\right\|_{L^2_x}^2
			+ \left\|B\right\|_{L^2_x}^2 \right)
			\\
			& \hspace{30mm} +
			\frac 14 \int_0^t \int_{\mathbb{R}^3} \int_{\mathbb{R}^3\times\mathbb{R}^3\times\mathbb{S}^2}
			{q}^2 MM_*dvdv_* d\sigma dx ds
			\leq C^\mathrm{in}.
		\end{aligned}
	\end{equation}

	There only remains to evaluate the contribution of the entropy dissipation in \eqref{entropy1 limit}, which will result from a direct application of the following Bessel inequality, established in \cite[Lemma 4.7]{BGL2}~:
	\begin{equation}\label{bessel 2}
		\begin{aligned}
			\frac 2{\mu}
			\left|\int_{\mathbb{R}^3\times\mathbb{R}^3\times\mathbb{S}^2} q
			\tilde\phi MM_*dvdv_*d\sigma\right|^2
			& + \frac{8}{5\kappa}
			\left|\int_{\mathbb{R}^3\times\mathbb{R}^3\times\mathbb{S}^2}q
			\tilde\psi MM_*dvdv_*d\sigma\right|^2
			\\
			& \leq
			\int_{\mathbb{R}^3\times\mathbb{R}^3\times\mathbb{S}^2}q^2 MM_*dvdv_*d\sigma,
		\end{aligned}
	\end{equation}
	where $\tilde \phi$ and $\tilde \psi$ are defined by \eqref{phi-psi-def inverses}.

	For the sake of completeness and for later reference, we provide a short justification of \eqref{bessel 2} below. But prior to this, let us conclude the proof of the present proposition. To this end, we employ the identity \eqref{q phi psi} from Proposition \ref{weak-comp} in combination with the relations \eqref{delta identities}, which we reproduce here for the mere convenience of the reader~:
	\begin{equation*}
		\begin{aligned}
			\int_{\mathbb{R}^3}\phi_{ij} \tilde\phi_{kl} Mdv &
			=\mu \left(\delta_{ik}\delta_{jl}+\delta_{il}\delta_{jk} - \frac 23 \delta_{ij}\delta_{kl} \right),
			\\
			\int_{\mathbb{R}^3}\psi_{i} \tilde\psi_{j} Mdv & = \frac 52 \kappa \delta_{ij},
		\end{aligned}
	\end{equation*}
	to deduce from the inequality \eqref{bessel 2} that
	\begin{equation*}
			2\mu
			\left|\nabla_x u + \nabla_x^t u - \frac 23 (\Div u)\operatorname{Id}\right|^2
			+ 10\kappa
			\left|\nabla_x\theta\right|^2
			\leq
			\int_{\mathbb{R}^3\times\mathbb{R}^3\times\mathbb{S}^2}q^2 MM_*dvdv_*d\sigma,
	\end{equation*}
	whence, thanks to the solenoidal constraint on $u$ established in \eqref{constraints1},
	\begin{equation*}
		\begin{aligned}
			\int_0^t & \left(\mu
			\left\|\nabla_x u\right\|_{L^2_x}^2
			+ \frac 52\kappa
			\left\|\nabla_x\theta\right\|_{L^2_x}^2\right)(s) ds \\
			& \hspace{30mm} \leq \frac 14
			\int_0^t\int_{\mathbb{R}^3}\int_{\mathbb{R}^3\times\mathbb{R}^3\times\mathbb{S}^2}q^2 MM_*dvdv_*d\sigma dxds.
		\end{aligned}
	\end{equation*}
	Combining this with \eqref{entropy1 limit} concludes the proof of the proposition.

	Now, as announced above, we give a short proof of \eqref{bessel 2}. To this end, following \cite[Lemma 4.7]{BGL2}, we recall that, for any traceless symmetric matrix $A\in\mathbb{R}^{3\times 3}$ and any vector $a\in\mathbb{R}^3$, one computes straightforwardly, employing the identities \eqref{delta identities} (reproduced above, for convenience) and the collisional symmetries, that
	\begin{equation*}
		\begin{aligned}
			\frac{1}{16}\int_{\mathbb{R}^3\times\mathbb{R}^3\times\mathbb{S}^2}
			& \left(A:\left(\tilde\phi+\tilde\phi_*-\tilde\phi'-\tilde\phi_*'\right)
			+
			a\cdot \left(\tilde\psi+\tilde\psi_*-\tilde\psi'-\tilde\psi_*'\right)\right)^2 MM_*dvdv_*d\sigma
			\\
			& =\frac 14 \left(A\otimes A\right):\int_{\mathbb{R}^3}\left(\phi\otimes\tilde\phi\right) Mdv
			+ \frac 14 \left(a\otimes a\right):\int_{\mathbb{R}^3}\left(\psi\otimes\tilde\psi\right) Mdv
			\\
			& =\frac 1 2\mu A:A + \frac{5}{8}\kappa a\cdot a.
		\end{aligned}
	\end{equation*}
	Therefore, defining, for any $q_0\in L^2\left(MM_*dvdv_*d\sigma\right)$, the projection
	\begin{equation*}
			\bar q_0 =
			A_0 :\frac 14 \left(\tilde\phi+\tilde\phi_*-\tilde\phi'-\tilde\phi_*'\right)
			+
			a_0 \cdot \frac 14 \left(\tilde\psi+\tilde\psi_*-\tilde\psi'-\tilde\psi_*'\right),
	\end{equation*}
	where
	\begin{equation*}
		\begin{aligned}
			A_0 & =\frac 1{2\mu}\int_{\mathbb{R}^3\times\mathbb{R}^3\times\mathbb{S}^2}q_0
			\left(\tilde\phi+\tilde\phi_*-\tilde\phi'-\tilde\phi_*'\right) MM_*dvdv_*d\sigma,
			\\
			a_0 & =\frac 2{5\kappa} \int_{\mathbb{R}^3\times\mathbb{R}^3\times\mathbb{S}^2}q_0
			\left(\tilde\psi+\tilde\psi_*-\tilde\psi'-\tilde\psi_*'\right) MM_*dvdv_*d\sigma,
		\end{aligned}
	\end{equation*}
	we find that
	\begin{equation*}
		\int_{\mathbb{R}^3\times\mathbb{R}^3\times\mathbb{S}^2}q_0\bar q_0 MM_*dvdv_*d\sigma
		=\frac 1 2\mu A_0:A_0 + \frac{5}{8}\kappa a_0\cdot a_0=
		\int_{\mathbb{R}^3\times\mathbb{R}^3\times\mathbb{S}^2}\bar q_0^2 MM_*dvdv_*d\sigma.
	\end{equation*}
	Hence the Bessel inequality
	\begin{equation}\label{bessel}
		\begin{aligned}
			\frac 1 2\mu A_0:A_0 + \frac{5}{8}\kappa a_0\cdot a_0
			& =
			\int_{\mathbb{R}^3\times\mathbb{R}^3\times\mathbb{S}^2}\bar q_0^2 MM_*dvdv_*d\sigma
			\\
			& \leq
			\int_{\mathbb{R}^3\times\mathbb{R}^3\times\mathbb{S}^2}q_0^2 MM_*dvdv_*d\sigma.
		\end{aligned}
	\end{equation}

	Therefore, setting $q_0=q$ in \eqref{bessel}, we find, exploiting the collisional symmetries of $q$, that
	\begin{equation*}
		\begin{aligned}
			\frac 2{\mu}
			\int_{\mathbb{R}^3\times\mathbb{R}^3\times\mathbb{S}^2} & q
			\tilde\phi MM_*dvdv_*d\sigma
			: \int_{\mathbb{R}^3\times\mathbb{R}^3\times\mathbb{S}^2}q
			\tilde\phi MM_*dvdv_*d\sigma
			\\
			& + \frac{8}{5\kappa}
			\int_{\mathbb{R}^3\times\mathbb{R}^3\times\mathbb{S}^2}q
			\tilde\psi MM_*dvdv_*d\sigma \cdot
			\int_{\mathbb{R}^3\times\mathbb{R}^3\times\mathbb{S}^2}q
			\tilde\psi MM_*dvdv_*d\sigma
			\\
			& \leq
			\int_{\mathbb{R}^3\times\mathbb{R}^3\times\mathbb{S}^2}q^2 MM_*dvdv_*d\sigma,
		\end{aligned}
	\end{equation*}
	which concludes the justification of \eqref{bessel 2}.
\end{proof}

\begin{prop}\label{energy inequality weak interactions}
	Let $\left(f_\eps^\pm, E_\eps, B_\eps\right)$ be the sequence of renormalized solutions to the scaled two species Vlasov-Maxwell-Boltzmann system \eqref{VMB2} considered in Theorem \ref{CV-OMHD} for weak interspecies interactions, i.e.\ $\delta=o(1)$ and $\frac\delta\eps$ unbounded. In accordance with Lemmas \ref{L1-lem}, \ref{L2-lem}, \ref{L2-qlem} and \ref{weak compactness h}, denote by
	\begin{equation*}
		\begin{gathered}
			g^\pm\in L^\infty\left(dt;L^2\left(Mdxdv\right)\right),
			\qquad
			h \in L^1_{\mathrm{loc}}\left(dtdx;L^1\left(\left(1+|v|^2\right)Mdv\right)\right),
			\\
			q^\pm, q^{\pm,\mp} \in L^2\left(MM_*dtdxdvdv_*d\sigma\right)
			\qquad\text{and}\qquad
			E,B\in L^\infty\left(dt;L^2\left(dx\right)\right),
		\end{gathered}
	\end{equation*}
	any joint limit points of the families $\hat g_\eps^\pm$,  $h_\eps$, $\hat q_\eps^{\pm}$ and $\hat q_\eps^{\pm,\mp}$ defined by \eqref{hatg}, \eqref{def h} and \eqref{hatq-def}, $E_\eps$ and $B_\eps$, respectively.
	
	Then, one has the energy inequality, for almost every $t\geq 0$,
	\begin{equation*}
		\begin{aligned}
			& \frac 12\left( 2 \left\|u\right\|_{L^2_x}^2
			+ 5 \left\|\theta\right\|_{L^2_x}^2 + \left\|E\right\|_{L^2_x}^2
			+ \left\|B\right\|_{L^2_x}^2 \right)(t)
			\\
			& \hspace{10mm} +
			\int_0^t \left(2\mu
			\left\|\nabla_x u\right\|_{L^2_x}^2
			+ 5 \kappa
			\left\|\nabla_x\theta\right\|_{L^2_x}^2
			+\frac 1{\sigma}
			\left\|j\right\|_{L^2_x}^2
			+ \frac 1{2\lambda}
			\left\|w\right\|_{L^2_x}^2\right)(s) ds \leq C^\mathrm{in},
		\end{aligned}
	\end{equation*}
	where $\rho$, $u$ and $\theta $ are, respectively, the density, bulk velocity and temperature associated with the limiting fluctuation $g$, while $j$ and $w$ are, respectively, the electric current and the internal electric energy associated with the limiting fluctuation $h$, and, finally, the viscosity $\mu>0$, thermal conductivity $\kappa>0$, electric conductivity $\sigma>0$ and energy conductivity $\lambda>0$ are respectively defined by \eqref{mu kappa 2}, \eqref{sigma} and \eqref{lambda}.
\end{prop}

\begin{proof}
	First, by the estimate \eqref{q-est} from Lemma \ref{L2-qlem} and the weak sequential lower semi-continuity of convex functionals, we find that, for all $t\geq 0$,
	\begin{equation*}
		\begin{aligned}
			\frac 14 \int_0^t \int_{\mathbb{R}^3} \int_{\mathbb{R}^3\times\mathbb{R}^3\times\mathbb{S}^2} &
			\left(q^\pm\right)^2 MM_*dvdv_* d\sigma dx ds
			\\
			& \leq \liminf_{\eps\rightarrow 0}
			\frac 14 \int_0^t \int_{\mathbb{R}^3} \int_{\mathbb{R}^3\times\mathbb{R}^3\times\mathbb{S}^2}
			\left(\hat{q}_\eps^\pm\right)^2 MM_*dvdv_* d\sigma dx ds
			\\
			& \leq \liminf_{\eps\rightarrow 0} \frac{1}{\epsilon^4}\int_0^t\int_{\mathbb{R}^3}
			D\left(f_\eps^\pm\right)(s) dx ds,
		\end{aligned}
	\end{equation*}
	and
	\begin{equation*}
		\begin{aligned}
			\frac 12 \int_0^t \int_{\mathbb{R}^3} \int_{\mathbb{R}^3\times\mathbb{R}^3\times\mathbb{S}^2} &
			\left(q^{\pm,\mp}\right)^2 MM_*dvdv_* d\sigma dx ds
			\\
			& \leq \liminf_{\eps\rightarrow 0}
			\frac 12 \int_0^t \int_{\mathbb{R}^3} \int_{\mathbb{R}^3\times\mathbb{R}^3\times\mathbb{S}^2}
			\left(\hat{q}_\eps^{\pm,\mp}\right)^2 MM_*dvdv_* d\sigma dx ds
			\\
			& \leq \liminf_{\eps\rightarrow 0} \frac{\delta^2}{\epsilon^4}\int_0^t\int_{\mathbb{R}^3}
			D\left(f_\eps^+,f_\eps^-\right)(s) dx ds,
		\end{aligned}
	\end{equation*}
	which, when combined with Lemma \ref{L1-lem}, yields, passing to the limit in the entropy inequality \eqref{entropy2}, for almost every $t\geq 0$,
	\begin{equation*}
		\begin{aligned}
		\frac 12 & \int_{\mathbb{R}^3\times\mathbb{R}^3} \left(\left(g^{+}\right)^2 + \left(g^{-}\right)^2\right)(t)Mdxdv
		+ \frac 1{2} \int_{\mathbb{R}^3} \left(|E|^2+ |B|^2\right)(t) dx
		\\
		& +
		\frac 14 \int_0^t \int_{\mathbb{R}^3} \int_{\mathbb{R}^3\times\mathbb{R}^3\times\mathbb{S}^2}
		\left(\left(q^+\right)^2 + \left(q^-\right)^2 + \left(q^{+,-}\right)^2+ \left(q^{-,+}\right)^2\right) MM_*dvdv_* d\sigma dx ds
		\\
		& \leq C^\mathrm{in}.
		\end{aligned}
	\end{equation*}
	Since, according to Lemmas \ref{relaxation-control} and \ref{relaxation2-control} and Proposition \ref{weak-comp2}, the limiting fluctuations $g^\pm=\rho+u\cdot v + \theta \left(\frac{|v|^2}{2}-\frac 32\right)$ are infinitesimal Maxwellians which coincide, we easily compute that, in view of the strong Boussinesq relation $\rho+\theta=0$ following from \eqref{constraints2-alpha},
	\begin{equation*}
		\frac 12 \int_{\mathbb{R}^3}\left(\left(g^{+}\right)^2 + \left(g^{-}\right)^2\right)Mdv
		= \rho^2 + |u|^2 + \frac 32 \theta^2
		= |u|^2 + \frac 52 \theta^2,
	\end{equation*}
	which implies
	\begin{equation}\label{entropy2 limit}
		\begin{aligned}
			& \left(\left\|u\right\|_{L^2_x}^2
			+ \frac 52\left\|\theta\right\|_{L^2_x}^2 + \frac 12 \left\|E\right\|_{L^2_x}^2
			+ \frac 12 \left\|B\right\|_{L^2_x}^2 \right)
			\\
			& +
			\frac 14 \int_0^t \int_{\mathbb{R}^3} \int_{\mathbb{R}^3\times\mathbb{R}^3\times\mathbb{S}^2}
			\left(\left(q^+\right)^2 + \left(q^-\right)^2 + \left(q^{+,-}\right)^2+ \left(q^{-,+}\right)^2\right) MM_*dvdv_* d\sigma dx ds
			\\
			& \leq C^\mathrm{in}.
		\end{aligned}
	\end{equation}

	There only remains to evaluate the contribution of the entropy dissipation in \eqref{entropy2 limit}. To this end, applying the method of proof of Proposition \ref{energy ineq 1}, based on the Bessel inequality \eqref{bessel 2}, with the constraints \eqref{q phi psi 2} and \eqref{constraints2-alpha} from Proposition \ref{weak-comp2}, note that it holds
	\begin{equation}\label{entropy2 1 limit}
		\begin{aligned}
			\int_0^t & \left(\mu
			\left\|\nabla_x u\right\|_{L^2_x}^2
			+ \frac 52\kappa
			\left\|\nabla_x\theta\right\|_{L^2_x}^2\right)(s) ds \\
			& \hspace{30mm} \leq \frac 14
			\int_0^t\int_{\mathbb{R}^3}\int_{\mathbb{R}^3\times\mathbb{R}^3\times\mathbb{S}^2}\left(q^\pm\right)^2 MM_*dvdv_*d\sigma dxds.
		\end{aligned}
	\end{equation}
	Next, the contributions of the mixed entropy dissipations $q^{\pm,\mp}$ will be evaluated through a direct application of the following Bessel inequality~:
	\begin{equation}\label{bessel 4}
		\begin{aligned}
			2\sigma
			\left|\int_{\mathbb{R}^3\times\mathbb{R}^3\times\mathbb{S}^2} q^{\pm,\mp}
			v MM_*dvdv_*d\sigma\right|^2
			& + \lambda
			\left|
			\int_{\mathbb{R}^3\times\mathbb{R}^3\times\mathbb{S}^2} q^{\pm,\mp}
			|v|^2 MM_*dvdv_*d\sigma\right|^2
			\\
			& \leq
			\int_{\mathbb{R}^3\times\mathbb{R}^3\times\mathbb{S}^2}\left(q^{\pm,\mp}\right)^2 MM_*dvdv_*d\sigma.
		\end{aligned}
	\end{equation}

	For the sake of completeness, we provide a short justification of \eqref{bessel 4} below. But prior to this, let us conclude the proof of the present proposition. To this end, we employ the identity \eqref{mixed q phi psi} from Proposition \ref{weak-comp3} in combination with the relations \eqref{sigma} and \eqref{lambda} to deduce from the inequality \eqref{bessel 4} that
	\begin{equation}\label{entropy2 1 limit 2}
		\begin{aligned}
			\frac 2\sigma
			\left| j \right|^2
			+ \frac 1\lambda
			\left|w\right|^2
			\leq
			\int_{\mathbb{R}^3\times\mathbb{R}^3\times\mathbb{S}^2}\left(q^{\pm,\mp}\right)^2 MM_*dvdv_*d\sigma.
		\end{aligned}
	\end{equation}
	Combining this with \eqref{entropy2 limit} and \eqref{entropy2 1 limit} concludes the proof of the proposition.

	Now, as announced, we give a short proof of \eqref{bessel 4}. To this end, for any vector $A\in\mathbb{R}^{3}$ and any scalar $a\in\mathbb{R}$, one computes straightforwardly, employing Proposition \ref{cross section transfer} and the collisional symmetries, that
	\begin{equation*}
		\begin{aligned}
			\int_{\mathbb{R}^3\times\mathbb{R}^3\times\mathbb{S}^2}
			& \left(A\cdot\left(v - v_* -v' + v_*'\right)
			+
			a \left(|v|^2 - |v_*|^2-|v'|^2 + |v_*'|^2\right)\right)^2 MM_*dvdv_*d\sigma
			\\
			& =
			4 \int_{\mathbb{R}^3\times\mathbb{R}^3\times\mathbb{S}^2}
			\left(A\cdot\left(v-v'\right)
			+
			a \left(|v|^2-|v'|^2\right)\right)^2 MM_*dvdv_*d\sigma
			\\
			& = \frac 8\sigma |A|^2 + \frac {16}\lambda a^2.
		\end{aligned}
	\end{equation*}
	Therefore, defining, for any $q_0\in L^2\left(MM_*dvdv_*d\sigma\right)$, the projection
	\begin{equation*}
			\bar q_0 =
			A_0 \cdot \left(v - v_* -v' + v_*'\right)
			+
			a_0 \left(|v|^2 - |v_*|^2-|v'|^2 + |v_*'|^2\right),
	\end{equation*}
	where
	\begin{equation*}
		\begin{aligned}
			A_0 & =\frac\sigma 8 \int_{\mathbb{R}^3\times\mathbb{R}^3\times\mathbb{S}^2}q_0
			\left(v - v_* -v' + v_*'\right) MM_*dvdv_*d\sigma,
			\\
			a_0 & =\frac\lambda {16} \int_{\mathbb{R}^3\times\mathbb{R}^3\times\mathbb{S}^2}q_0
			\left(|v|^2 - |v_*|^2-|v'|^2 + |v_*'|^2\right)
			MM_*dvdv_*d\sigma,
		\end{aligned}
	\end{equation*}
	we find that
	\begin{equation*}
		\int_{\mathbb{R}^3\times\mathbb{R}^3\times\mathbb{S}^2}q_0\bar q_0 MM_*dvdv_*d\sigma
		=\frac 8\sigma \left|A_0\right|^2+\frac {16}\lambda a_0^2
		=
		\int_{\mathbb{R}^3\times\mathbb{R}^3\times\mathbb{S}^2}\bar q_0^2 MM_*dvdv_*d\sigma.
	\end{equation*}
	Hence the Bessel inequality
	\begin{equation}\label{bessel 3}
		\begin{aligned}
			\frac 8\sigma \left|A_0\right|^2+\frac {16}\lambda a_0^2
			& =
			\int_{\mathbb{R}^3\times\mathbb{R}^3\times\mathbb{S}^2}\bar q_0^2 MM_*dvdv_*d\sigma
			\\
			& \leq
			\int_{\mathbb{R}^3\times\mathbb{R}^3\times\mathbb{S}^2}q_0^2 MM_*dvdv_*d\sigma.
		\end{aligned}
	\end{equation}

	Therefore, setting $q_0=q^{\pm,\mp}$ in \eqref{bessel 3}, we find, exploiting the collisional symmetries of $q^{\pm,\mp}$, that
	\begin{equation*}
		\begin{aligned}
			2\sigma
			\left|\int_{\mathbb{R}^3\times\mathbb{R}^3\times\mathbb{S}^2} q^{\pm,\mp}
			v MM_*dvdv_*d\sigma\right|^2
			& + \lambda
			\left(
			\int_{\mathbb{R}^3\times\mathbb{R}^3\times\mathbb{S}^2} q^{\pm,\mp}
			|v|^2 MM_*dvdv_*d\sigma\right)^2
			\\
			& \leq
			\int_{\mathbb{R}^3\times\mathbb{R}^3\times\mathbb{S}^2}\left(q^{\pm,\mp}\right)^2 MM_*dvdv_*d\sigma,
		\end{aligned}
	\end{equation*}
	which concludes the justification of \eqref{bessel 4}.
\end{proof}

\section{The limiting Maxwell's equations}\label{limit maxwell}

Using the uniform $L^\infty\left(dt;L^2\left(dx\right)\right)$ bounds on the electromagnetic fields $E_\eps$ and $B_\eps$, and the controls from Chapter \ref{weak bounds} on the fluctuations, we can also take limits in the full Maxwell system for one species and for two species in the case of weak interactions only. Because of the scaling of the light speed, we obtain different kinds of limiting systems in the two regimes to be considered, but there is no particular difficulty here, for everything remains linear.

As for the case of two species with strong interactions, we will not be able to pass to the limit in Maxwell's equations. Indeed, Amp\`ere's equation is nonlinear in this setting, which is a major obstacle to the weak stability of the system. More comments on this issue are provided below.

\bigskip

In the regime leading to the incompressible quasi-static Navier-Stokes-Fourier-Maxwell-Poisson system, considered in Section \ref{main moments method}, we start from
\begin{equation*}
	\begin{cases}
		\begin{aligned}
			\eps\d_t E_\eps - \ROT B_\eps &= - \int_{\mathbb{R}^3} g_\eps v M dv,
			\\
			\eps \d_t B_\eps + \ROT E_\eps & = 0,
			\\
			\DIV E_\eps &=\int_{\mathbb{R}^3} g_\eps M dv,
			\\
			\DIV B_\eps &=0.
		\end{aligned}
	\end{cases}
\end{equation*}
Then, the weak compactness of the fluctuations from Lemma \ref{L1-lem}, inherited from the scaled entropy inequality \eqref{entropy1}, allows us to consider converging subsequences
\begin{equation*}
	\begin{aligned}
		E_\eps & \stackrel{*}{\rightharpoonup} E & & \text{in } L^\infty\left(dt;L^2\left(dx\right)\right),
		\\
		B_\eps & \stackrel{*}{\rightharpoonup} B & & \text{in } L^\infty\left(dt;L^2\left(dx\right)\right),
		\\
		g_\eps & \rightharpoonup g & & \text{in } L^1_\mathrm{loc}\left(dtdx ; L^1\left(\left(1+\left|v\right|^2\right)M dv\right)\right),
	\end{aligned}
\end{equation*}
which easily leads to
\begin{equation*}
	\ROT E = 0, \qquad \Div E = \rho, \qquad  \ROT B = u, \qquad \Div B = 0,
\end{equation*}
where $\rho$ and $u$ respectively denote the density and bulk velocity associated to the limiting fluctuation $g$. Formally, this limit amounts to discarding the terms involving time derivatives in Maxwell's equations, which accounts for the terminology of ``quasi-static approximation'' since temporal variations are neglected.

\bigskip

Next, in the regime leading to the two-fluid incompressible Navier-Stokes-Fourier-Maxwell system with (solenoidal) Ohm's law, considered in Section \ref{main relative entropy}, we start from
\begin{equation*}
	\begin{cases}
		\begin{aligned}
			\d_t E_\eps - \ROT B_\eps &= - \frac{\delta}{\eps} \int_{\mathbb{R}^3} \left(g_\eps^+-g_\eps^-\right)v M dv
			=-\int_{\mathbb{R}^3}h_\eps vMdv,
			\\
			\d_t B_\eps + \ROT E_\eps& = 0,
			\\
			\DIV E_\eps &=\delta\int_{\mathbb{R}^3} \left(g_\eps^+-g_\eps^-\right) M dv,
			\\
			\DIV B_\eps &=0.
		\end{aligned}
	\end{cases}
\end{equation*}

We consider first the simpler case of weak interspecies collisions, i.e.\ the case $\delta=o(1)$ and $\frac\delta\eps$ unbounded. The weak compactness of the fluctuations from Lemmas \ref{L1-lem} and \ref{weak compactness h}, inherited from the scaled entropy inequality \eqref{entropy2}, allows us to consider converging subsequences
\begin{equation*}
	\begin{aligned}
		E_\eps & \stackrel{*}{\rightharpoonup} E & & \text{in } L^\infty\left(dt;L^2\left(dx\right)\right),
		\\
		B_\eps & \stackrel{*}{\rightharpoonup} B & & \text{in } L^\infty\left(dt;L^2\left(dx\right)\right),
		\\
		g_\eps^\pm & \rightharpoonup g^\pm & & \text{in } L^1_\mathrm{loc}\left(dtdx ; L^1\left(\left(1+\left|v\right|^2\right)M dv\right)\right),
		\\
		h_\eps & \rightharpoonup h & & \text{in } L^1_\mathrm{loc}\left(dtdx ; L^1\left(\left(1+\left|v\right|^2\right)M dv\right)\right),
	\end{aligned}
\end{equation*}
which easily leads to
\begin{equation*}
	\begin{cases}
		\begin{aligned}
			\d_t E - \ROT B &= -j,
			\\
			\d_t B + \ROT E& = 0,
			\\
			\DIV E &=0,
			\\
			\DIV B &=0,
		\end{aligned}
	\end{cases}
\end{equation*}
where $j$ denotes the electric current, that is the bulk velocity associated to the limiting fluctuation $h$.

Now, we see that in the case of strong interspecies collisions, i.e.\ $\delta=1$, Lemma \ref{weak compactness h} provides no longer enough compactness on $h_\eps$ to take weak limits in Amp\`ere's equation. Indeed, in view of Lemma \ref{bound hjw}, it holds, at best, that the $h_\eps$'s are uniformly bounded in $L^1_{\mathrm{loc}}\left(dtdx;L^1\left(\left(1+|v|^2\right)Mdv\right)\right)$, but nothing prevents the fluctuations $h_\eps$ from concentrating on small sets and, therefore, to converge towards a singular measure. Thus, in this asymptotic regime, Amp\`ere's equation will not be satisfied in the sense of distributions but only in a dissipative sense, which will be encoded in the inequality defining the dissipative solutions obtained in Chapter \ref{entropy method} through a generalized relative entropy method.

A closer inspection of Amp\`ere's equation in the limiting system \eqref{TFINSFMO 2} (or \eqref{TFINSFMO}) provides some insight on its lack of weak stability in the hydrodynamic limit. Indeed, even though Maxwell's system in \eqref{TFINSFMO 2} is linear in the variables $(E,B,n,j)$, the energy inequality associated with \eqref{TFINSFMO 2} suggests that the right mathematical variables are rather $(E,B,j-nu,n,u)$, which renders Amp\`ere's equation nonlinear.

Nevertheless, the rest of Maxwell's system remains linear and we can easily pass to the limit in Faraday's equation and Gauss' laws. Indeed, the weak compactness of the fluctuations from Lemma \ref{L1-lem}, inherited from the scaled entropy inequality \eqref{entropy2}, allows us to consider converging subsequences
\begin{equation*}
	\begin{aligned}
		E_\eps & \stackrel{*}{\rightharpoonup} E & & \text{in } L^\infty\left(dt;L^2\left(dx\right)\right),
		\\
		B_\eps & \stackrel{*}{\rightharpoonup} B & & \text{in } L^\infty\left(dt;L^2\left(dx\right)\right),
		\\
		g_\eps^\pm & \rightharpoonup g^\pm & & \text{in } L^1_\mathrm{loc}\left(dtdx ; L^1\left(\left(1+\left|v\right|^2\right)M dv\right)\right),
	\end{aligned}
\end{equation*}
which easily leads to
\begin{equation*}
	\begin{cases}
		\begin{aligned}
			\d_t B + \ROT E& = 0,
			\\
			\DIV E &=n,
			\\
			\DIV B &=0,
		\end{aligned}
	\end{cases}
\end{equation*}
where $n=\rho^+-\rho^-$ is the electric charge associated with the limiting fluctuations $g^\pm$, i.e.\ $\rho^\pm$ are the macroscopic densities of $g^\pm$.

%% file: strong0.tex
\chapter{Strong compactness and hypoellipticity}\label{hypoellipticity}

In Chapter \ref{weak bounds}, we have established uniform bounds and relaxation estimates on the fluctuations and collision integrands as consequences of the scaled relative entropy inequalities \eqref{entropy1} and \eqref{entropy2}. This is sufficient to handle linear terms. Thus, in Chapter \ref{constraints proof}, we exploited these uniform estimates to derive limiting constraint equations and energy inequalities.

In order to go any further in the rigorous derivation of the hydrodynamic limits under study, we need now to obtain precise strong compactness estimates on the fluctuations through a refined understanding of the Vlasov-Boltzmann equations from \eqref{VMB1} and \eqref{VMB2}. More precisely, in the present chapter, we are going to introduce {\bf mathematical tools used to study the dependence in $x$ and $v$} of the families of fluctuations and, then, deduce important strong compactness properties of these fluctuations.

The first and simplest step, performed in Section \ref{velocity compactness} below, consists in understanding the {\bf dependence of fluctuations with respect to the velocity variable}, which is essentially controlled by the relaxation mechanism. Since these estimates in $v$ are based only on results from functional analysis and on the relative entropy and entropy dissipation bounds, they will hold similarly in both regimes \eqref{VMB1} and \eqref{VMB2}.

This first step is novel and differs considerably from previous works on hydrodynamic limits of Boltzmann equations with cutoff assumptions in that it shows strong compactness of the fluctuations in velocity, whereas former results only employed weak bounds in $v$, such as the equi-integrability in $v$ from Lemma \ref{v2-int}. This strong compactness is crucial in order to carry out the next stage of the proof in Section \ref{hypoelliptic}. Note that strong velocity compactness has also been used in \cite{arsenio3} to treat hydrodynamic limits of the Boltzmann equation without any cutoff assumptions. The approach therein heavily relied on the smoothing effect in $v$ peculiar to long-range interactions, though. In fact, the methods developed here can also be used to improve the results from \cite{arsenio3} (see Part \ref{part 3}).

The second, more convoluted step, performed in Section \ref{hypoelliptic}, uses then the {\bf hypoellipticity in kinetic transport equations} studied in \cite{arsenio} to transfer strong compactness from the velocity variable $v$ to the space variable $x$. Some non-trivial technical care will be required in order to extend the results from \cite{arsenio}, which mainly concern the stationary kinetic transport equation, to the non-stationary transport equation with a vanishing time derivative.

Note that this second step also differs substantially from previous works on the subject, for these traditionally relied on classical velocity averaging lemmas to show some strong space compactness of the moments of the fluctuations (not the fluctuations themselves).

It is to be emphasized that the compactness properties for the two species regime \eqref{VMB2} obtained in Section \ref{compactness two species}, below, are substantially weaker than those corresponding to the one species regime \eqref{VMB1} and derived in Section \ref{compactness one species}. Essentially, the two species regime considered here being quite singular, the corresponding fluctuations cannot be shown to enjoy as much equi-integrability as in the one species regime, which will lead to significant difficulties in the remainder of our proofs.

The results from the present chapter constitute a crucial and difficult step in the rigorous proofs of hydrodynamic convergence. They will allow us to obtain higher order nonlinear constraint and evolution equations in the coming chapters. Finally, note that the results obtained here are only concerned with the compactness properties of fluctuations in $x$ and $v$, but not in $t$. In fact, there may be oscillations in time and the temporal behavior of fluctuations will be analyzed later on in Chapter \ref{oscillations}.

% ========================
% = vvvvvvvvvvvvvvvvvvvv =
% ========================

\section{Compactness with respect to $v$}\label{velocity compactness}

We have already shown in Section \ref{velocity integrability} how the relaxation process towards statistical equilibrium provides improved integrability in $v$ on the fluctuations. We show now how it further yields dissipative properties in the velocity variable. Loosely speaking, such a dissipation mechanism stems from the fact that the entropy dissipation controls the distance from the solutions to the set of statistical equilibria, which are in general smooth distributions in velocity.

As we consider fluctuations around a global equilibrium, the linearized collision operator will play a fundamental role, just as in Section \ref{relaxation} on the relaxation. In order to get strong compactness results, we will further need to control the correctors coming from the nonlinear part of the collision operator. To this end, we recall now the important regularizing effects of the gain term of the Boltzmann collision operator. This property will be crucial in our proof of compactness.

The results presented in Section \ref{section compactness gain 1}, below, concern general cross-sections satisfying some integrability assumptions. The properties from Section \ref{section compactness gain 2}, however, only concern the Maxwellian collision kernel $b\equiv 1$. The corresponding results for general cross-sections will be discussed in the remaining parts of our work.

\subsection{Compactness of the gain term}\label{section compactness gain 1}

In \cite{lions}, Lions exhibited the compactifying and regularizing effects of the gain term of the Boltzmann collision operator. The essential result contained therein establishes the regularity of the gain term for a smooth and truncated collision kernel. The precise result from \cite{lions} which is of interest to us is recalled in the following theorem. Variants and refinements of this result were obtained in \cite{bouchut, wennberg}. In particular, a simple argument based on the Fourier transform, due to Bouchut and Desvillettes in \cite{bouchut}, also provides a convenient compactness result.

\begin{thm}[\cite{lions}]\label{reg gain term}
	Let $b(z,\sigma)=b\left(|z|,\frac{z}{|z|}\cdot\sigma\right)\in C_c^\infty\left(\left(0,\infty\right)\times\left(0,\pi\right)\right)$ be a smooth compactly supported collision kernel.
	
	Then, there exists a finite $C>0$ such that
	\begin{equation*}
		\left\|Q^+\left(f,g\right)\right\|_{H^1\left(\mathbb{R}^3\right)}
		\leq
		C
		\left\|f\right\|_{L^2\left(\mathbb{R}^3\right)}
		\left\|g\right\|_{L^1\left(\mathbb{R}^3\right)},
	\end{equation*}
	for any $f\in L^2\left(\mathbb{R}^3\right)$ and $g\in L^1\left(\mathbb{R}^3\right)$, and
	\begin{equation*}
		\left\|Q^+\left(f,g\right)\right\|_{H^1\left(\mathbb{R}^3\right)}
		\leq
		C
		\left\|f\right\|_{L^1\left(\mathbb{R}^3\right)}
		\left\|g\right\|_{L^2\left(\mathbb{R}^3\right)},
	\end{equation*}
	for any $f\in L^1\left(\mathbb{R}^3\right)$ and $g\in L^2\left(\mathbb{R}^3\right)$.
\end{thm}

% \begin{prop}{\bf Regularity of $Q^+$}\label{Q+-prop}\marginpar{\small\color{blue}According to \cite{bouchut}, it seems that this version is wrong, unless we truncate the cross section.}
% Denote by $Q^+$ the gain part of the collision operator (with $b\equiv 1$)
%  Then, for any $f \in L^2((1+|v|)dv)$,
% $$\| Q^+(f,f)\|_{H^1(\R^3)} \leq C \| f \|_{L^2((1+|v|)dv)}^2\,.$$
% In particular, with the notation
% $$ \cQ^+(g,g) =\frac1M Q^+(Mg,Mg)\,,$$
% we have for any $g \in L^2(Mdv)$,
% $$\| M\cQ^+ (g,g)\|_{H^1(\R^3)}\leq C \| g \| _{L^2(Mdv)}\,.$$
% \end{prop}
% 
% Note that this result extends actually the study by Lions \cite{lions} to much more general situations excluding however singular 
% cross-sections. We will explain in the second part how to deal with the non cutoff case.

Note that, in the statement of the above theorem, we have carefully avoided the endpoints on the domain of definition of the collision kernel in order to restrict the compact support of $b(z,\sigma)$. More precisely, Lions' result only considers smooth kernels whose support is contained in $\left\{\lambda<|z|<\frac 1\lambda, \left|\frac{z}{|z|}\cdot\sigma\right| < 1-\lambda \right\}$, for some small $\lambda>0$. This hypothesis is definitely not optimal, but at least some truncation is clearly required in order to obtain the optimal gain of regularity for $Q^+(f,g)$.

For more general collision kernels, it is still possible to obtain some compactness of the gain operator by standard approximation procedures based on convolution inequalities for the gain term $Q^+(f,g)$. See for instance \cite{arsenio4} for such general convolution inequalities. Here, we will merely use an elementary version of these inequalities which we presently recall for convenience.

Thus, let $1\leq s\leq p,q\leq r'\leq\infty$ be such that
\begin{equation*}
	1+\frac{1}{s}=\frac{1}{p}+\frac{1}{q}+\frac{1}{r}
\end{equation*}
and consider $f\in L^p\left(\mathbb{R}^3\right)$, $g\in L^q\left(\mathbb{R}^3\right)$, $\varphi\in L^{s'}\left(\mathbb{R}^3\right)$ and a general collision kernel $b(z,\sigma)\in L^r\left(\mathbb{R}^3;L^1\left(\mathbb{S}^2\right)\right)$. Then, employing the collision symmetries with H\"older's and Young's inequalities and using the change of variables $v_*\mapsto V=v-v_*$, we find
\begin{equation*}
	\begin{aligned}
		\bigg|\int_{\mathbb{R}^3}Q^+ & (f,g)(v)\varphi(v) dv\bigg|
		\\
		& \leq
		\int_{\mathbb{R}^3\times\mathbb{R}^3\times\mathbb{S}^2} \left| fg_*\varphi' \right| b(v-v_*,\sigma) dvdv_*d\sigma
		\\
		&=
		\int_{\mathbb{R}^3\times\mathbb{R}^3\times\mathbb{S}^2} \left| f(v)g(v-V)\varphi\left(v-\frac{V}{2}+\frac{|V|}{2}\sigma\right) \right| b(V,\sigma) dVdvd\sigma
		\\
		&\leq
		\left\|\varphi\right\|_{L^{s'}}
		\int_{\mathbb{R}^3} \left\| f(v)g(v-V) \right\|_{L^s_v} \int_{\mathbb{S}^2}b(V,\sigma) d\sigma dV
		\\
		&\leq
		\left\|\varphi\right\|_{L^{s'}}
		\left\|\left\| f(v)g(v-V) \right\|_{L^s_v}\right\|_{L^{r'}_V}
		\left\|\int_{\mathbb{S}^2}b(V,\sigma) d\sigma\right\|_{L^r_V}
		\\
		&=
		\left\|\varphi\right\|_{L^{s'}}
		\left\| \int_{\mathbb{R}^3} \left|f(v)\right|^s\left|g(v-V)\right|^s dv \right\|_{L^{\frac{r'}{s}}_V}^\frac{1}{s}
		\left\|\int_{\mathbb{S}^2}b(V,\sigma) d\sigma\right\|_{L^r_V}
		\\
		&\leq
		\left\|\varphi\right\|_{L^{s'}}
		\left\|f\right\|_{L^{p}}\left\|g\right\|_{L^{q}}
		\left\|\int_{\mathbb{S}^2}b(V,\sigma) d\sigma\right\|_{L^r_V}.
	\end{aligned}
\end{equation*}
Notice that the exact same reasoning can be applied to the loss operator $Q^-(f,g)$, so that, considering the supremum over all $\varphi\in L^{s'}\left(\mathbb{R}^3\right)$, we arrive at the following estimate
\begin{equation}\label{convolution}
		\left\| Q^\pm(f,g)\right\|_{L^s}
		\leq
		\left\|f\right\|_{L^{p}}\left\|g\right\|_{L^{q}}
		\left\|\int_{\mathbb{S}^2}b(z,\sigma) d\sigma\right\|_{L^r_z}.
\end{equation}

It turns out that it is possible to extend the above inequality to the full range of parameters $1\leq p,q,r,s\leq\infty$ for the gain term $Q^+(f,g)$ only, provided we have a better control on the angular collision kernel. This is consistent with the fact that $Q^+(f,g)$ behaves nicely and better than $Q^-(f,g)$. Such results can be found in \cite{arsenio4}.

Thus, combining the regularizing properties from Theorem \ref{reg gain term} with the convolution inequalities \eqref{convolution}, we obtain the following convenient proposition.

\begin{prop}\label{compactness gain term}
	Let $b(z,\sigma)$ be a cross-section such that
	\begin{equation*}
		\int_{\mathbb{S}^2}b(z,\sigma)d\sigma \in L^2\left(\mathbb{R}^3, M(z)^\alpha dz\right),
	\end{equation*}
	for some given $\alpha < \frac{1}{2}$.
	
	Then, the bilinear operator
	\begin{equation*}
		\begin{array}{ccc}
			L^2\left(\mathbb{R}^3,Mdv\right)\times L^2\left(\mathbb{R}^3,Mdv\right)
			&  \longrightarrow & L^2\left(\mathbb{R}^3, M^{1+2\alpha}dv\right) \\
			(f,g) & \longmapsto & \mathcal{Q}^+\left(f,g\right)
		\end{array}
	\end{equation*}
	is locally compact. That is to say, it maps bounded subsets of $L^2\left(Mdv\right)\times L^2\left(Mdv\right)$ into relatively compact subsets of $L_{\mathrm{loc}}^2\left(dv\right)$.
\end{prop}

\begin{proof}
	First, it is easy to check that
	\begin{equation*}
		\begin{aligned}
			\left(Mf\right)'\left(Mg\right)_*'M^{\alpha-\frac 12}
			&=\left(\sqrt{M}f\right)'\left(\sqrt{M}g\right)_*'M^{\alpha}M_*^\frac{1}{2}\\
			&\leq C \left(\sqrt{M}f\right)'\left(\sqrt{M}g\right)_*'\left(MM_*\right)^\alpha\\
			&\leq C
			\left(\sqrt{M}|f|\right)'\left(\sqrt{M}|g|\right)_*' \left(M(v-v_*)\right)^{\frac\alpha 2},
		\end{aligned}
	\end{equation*}
	for some $C>0$. Hence, in virtue of the convolution inequality \eqref{convolution}, we obtain
	\begin{equation}\label{compactess tightness gain operator}
		\begin{aligned}
			& \left\| \mathcal{Q}^+\left(f,g\right)M^{\alpha}\right\|_{L^2(Mdv)}\\
			&\hspace{7mm} \leq C
			\left\| \int_{\mathbb{R}^3\times\mathbb{S}^2}
			\left(\sqrt{M}|f|\right)'\left(\sqrt{M}|g|\right)_*' b(v-v_*,\sigma)M(v-v_*)^\frac{\alpha}{2} dv_* d\sigma
			\right\|_{L^2(dv)}\\
			&\hspace{7mm} \leq C
			\left\|f\right\|_{L^{2}\left(Mdv\right)}\left\|g\right\|_{L^{2}\left(Mdv\right)}
			\left\|\int_{\mathbb{S}^2}b(z,\sigma) d\sigma\right\|_{L^2\left(M^{\alpha}dz\right)},
		\end{aligned}
	\end{equation}
	which establishes the boundedness of the quadratic operator.
	
	Next, in order to show the local compactness of the operator, we consider any bounded sequences $\left\{f_n\right\}_{n\in\mathbb{N}}, \left\{g_n\right\}_{n\in\mathbb{N}}\subset L^2\left(Mdv\right)$.
	% Then, the tightness of $\left\{Q^+\left(Mf_n,Mg_n\right)\right\}_{n\in\mathbb{N}}$ in $L^2\left(dv\right)$ easily follows from \eqref{compactess tightness gain operator}, by considering any $0< \gamma\leq \frac 12 - \alpha$.
	Then, defining $\mathcal{Q}_\lambda^+\left(f_n,g_n\right)$, for any $\lambda>0$, by simply replacing $b(z,\sigma)$ by some smooth kernel $b_\lambda(z,\sigma)=b_\lambda\left(|z|,\frac{z}{|z|}\cdot\sigma\right)\in C_c^\infty\left(\left(0,\infty\right)\times\left(0,\pi\right)\right)$ such that $0\leq b_\lambda\leq b$ and
	\begin{equation*}
		\left\|\int_{\mathbb{S}^2}|b-b_\lambda|(z,\sigma) d\sigma\right\|_{L^2\left(M^{\alpha}dz\right)} < \lambda,
	\end{equation*}
	we deduce, thanks to \eqref{compactess tightness gain operator}, that the $\mathcal{Q}^+\left(f_n,g_n\right)$'s can be uniformly approximated by the $\mathcal{Q}_\lambda^+\left(f_n,g_n\right)$'s in $L^2\left(M^{1+2\alpha}dv\right)$. Since, by Theorem \ref{reg gain term}, the $Q_\lambda^+\left(Mf_n,Mg_n\right)$'s are relatively compact in $L_\mathrm{loc}^2(dv)$, we conclude that the original sequence $\left\{Q^+\left(Mf_n,Mg_n\right)\right\}_{n\in\mathbb{N}}$ is relatively compact in $L^2_\mathrm{loc}\left(dv\right)$, which concludes the justification of the proposition.
\end{proof}

\subsection{Relative entropy, entropy dissipation and strong compactness}\label{section compactness gain 2}

Combining the uniform controls from the relative entropy and entropy dissipation with the compactness of the gain term presented in the previous section, we establish now the following result, valid for the Maxwellian collision kernel $b\equiv 1$.

\begin{lem}\label{v-compactness}
	Let $f_\eps(t,x,v)$ be a family of measurable, almost everywhere non-negative distribution functions such that, for all $t\geq 0$,
	\begin{equation*}
		\frac1{\eps^2} H\left(f_\eps\right)(t)
		+
		\frac{1}{\epsilon^4}\int_0^t\int_{\mathbb{R}^3}
		D\left(f_\eps\right)(s) dx ds
		\leq C^\mathrm{in}.
	\end{equation*}
	
	Then, as $\eps\rightarrow 0$, any subsequence of renormalized fluctuations $\hat g_\eps$ is locally relatively compact in $v$ in $L^2\left(dtdxdv\right)$ in the sense that, for any $\eta>0$ and every compact subset $K\subset [0,\infty)\times\mathbb{R}^3\times\mathbb{R}^3$, there exists $\gamma>0$ such that, if $h\in\mathbb{R}^3$ satisfies $|h|<\gamma$, then
	\begin{equation*}
		\sup_{\eps>0} \left\|\hat g_\eps(t,x,v+h) - \hat g_\eps(t,x,v)\right\|_{L^2\left(K,dtdxdv\right)}
		<\eta.
	\end{equation*}
\end{lem}

\begin{proof}
	Loosely speaking, the present proof can be summarized in three main steps, each corresponding to a decomposition of ${\hat g_\eps}$. First, we will show how to control the very large values of ${\hat g_\eps}$, i.e.\ values larger than $\frac{1}{\eps}$, with the entropy bound. This is a rather standard and simple estimate. Second, we obtain the strong compactness in velocity of ${\hat g_\eps}$ from the entropy dissipation bound and the compactness of the gain term (see Proposition \ref{compactness gain term}) making sure that we remain away from vacuum, i.e.\ away from the values $\int \left( 1 + \frac\eps 2 \hat g_\eps \right) M dv \approx 0$, for our estimates degenerate in this case. Finally, we deduce the strong compactness near vacuum arguing that the vacuum state $ \hat g_\eps \equiv  -\frac{2}{\eps}$ is actually smooth since it is constant.

	\noindent {\bf Control of very large values.} Using the entropy inequality, we first introduce some microscopic truncation of large values. For any fixed small $0<\lambda<\frac 12$ and any cutoff $\chi(r)\in C_c^\infty\left(\mathbb{R}\right)$ such that $\mathds{1}_{\left\{|r|\leq 1 \right\}}\leq\chi(r)\leq\mathds{1}_{\left\{|r|\leq 2\right\}}$, we have, since $\eps^2\hat g_\eps^2\leq 4G_\eps$ as soon as $\eps\hat g_\eps\geq -1$,
	\begin{equation*}
		\begin{aligned}
			\left|{\hat g_\eps} \left ( 1- \chi\left( \lambda \eps \hat g_\eps \right)\right)\right|^2
			& \leq
			\hat g_\eps^2 \mathds{1}_{\left\{\eps\hat g_\eps\geq\frac 1\lambda>2\right\}}
			\\ & \leq {4\over \eps^2} G_\eps \mathds{1}_{\left\{ G_\eps > \frac{1}{4\lambda^2}\right\}}
			\\ & \leq {2\over \eps^2\left|\log 2\lambda\right|} G_\eps\log G_\eps \mathds{1}_{\left\{ G_\eps > \frac{1}{4\lambda^2}\right\}},
		\end{aligned}
	\end{equation*}
	so that, by the relative entropy bound,
	\begin{equation}\label{micro-truncation}
		{\hat g_\eps} \left ( 1- \chi\left( \lambda \eps \hat g_\eps \right)\right)
		=O\left({1\over |\log \lambda|^\frac{1}{2}}\right)_{L^\infty \left(dt ; L^2\left(Mdxdv\right)\right)},
	\end{equation}
	as $\lambda \to 0$, uniformly in $\eps$.

	\noindent {\bf Away from vacuum.} We use now Hilbert's decomposition (see Proposition \ref{hilbert-prop}) for the Maxwellian cross-section $b\equiv 1$~:
	\begin{equation*}
		\cL \hat g_\eps =\hat g_\eps-\cK \hat g_\eps,
	\end{equation*}
	where $\cK$ is a compact integral operator on $L^2(Mdv)$. Then, from the identity
	\begin{equation*}
		\cL \hat g_\eps =
		\frac\eps 2  \cQ\left(\hat g_\eps ,\hat g_\eps \right)-\frac 2{\eps}  \cQ \left(\sqrt{G_\eps},\sqrt{G_\eps}\right),
	\end{equation*}
	we deduce that
	\begin{equation*}
		\begin{aligned}
			\hat g_\eps &=  \cL \hat g_\eps + \cK \hat g_\eps \\
			&=  \frac\eps 2 \mathcal{Q}^+\left(\hat g_\eps,\hat g_\eps\right)
			-\frac\eps 2 \hat g_\eps \int_{\mathbb{R}^3\times\mathbb{S}^2} \hat g_{\eps *} M_*dv_*d\sigma
			-\eps \int_{\mathbb{R}^3\times\mathbb{S}^2} \hat q_\eps M_*dv_*d\sigma
			+ \cK \hat g_\eps,
		\end{aligned}
	\end{equation*}
	or, equivalently,
	\begin{equation}\label{comp-decomposition}
			\hat g_\eps \left(1
			+\frac\eps 2 \int_{\mathbb{R}^3} \hat g_{\eps *} M_*dv_* \right)
			=
			-\eps \int_{\mathbb{R}^3\times\mathbb{S}^2} \hat q_\eps M_*dv_*d\sigma
			+ \cK \hat g_\eps
			+ \frac\eps 2 \mathcal{Q}^+\left(\hat g_\eps,\hat g_\eps\right).
	\end{equation}
	We are now going to control each term in the right-hand side above separately.

	The first term is easily estimated employing the uniform $L^2$-estimate from Lemma \ref{L2-qlem}. It yields that
	\begin{equation}\label{comp-decomposition 1}
		\eps \int_{\mathbb{R}^3\times\mathbb{S}^2} \hat q_\eps M_*dv_*d\sigma
		=O\left(\eps\right)_{L^2\left(Mdtdxdv\right)}.
	\end{equation}
	The second term $\cK \hat g_\eps$ satisfies the bound
	\begin{equation*}
		\left\|\cK \hat g_\eps\right\|_{L^2(Mdv)}\leq C \left\|\hat g_\eps\right\|_{L^2(Mdv)},
	\end{equation*}
	so that, in view of Lemma \ref{L2-lem} and by the compactness of the operator $\cK$,
	\begin{equation}\label{comp-decomposition 2}
		\cK \hat g_\eps = O(1)_{L^\infty\left(dt;L^2\left(dx;\mathcal{C} L^2(Mdv) \right)\right)},
	\end{equation}
	where we have used the notation $\mathcal{C} L^2(Mdv)$ to indicate that it is relatively compact with respect to the velocity variable in $L^\infty\left(dt;L^2\left(Mdxdv\right)\right)$. Similarly, the third term $\frac\eps 2 \mathcal{Q}^+\left(\hat g_\eps,\hat g_\eps\right)$ satisfies, in virtue of Proposition \ref{compactness gain term}, the control
	\begin{equation}\label{comp-decomposition 3}
		\frac\eps 2 \mathcal{Q}^+\left(\hat g_\eps,\hat g_\eps\right) = O(\eps)_{L^\infty\left(dt;L^1\left(dx;\mathcal{C} L^2_\mathrm{loc}(dv) \right)\right)},
	\end{equation}
	where, again, we have used the notation $\mathcal{C} L_\mathrm{loc}^2(dv)$ to indicate that it is relatively compact with respect to the velocity variable in $L^\infty\left(dt;L^1\left(dx;L_\mathrm{loc}^2\left(dv\right)\right)\right)$.

	On the whole, incorporating the controls \eqref{comp-decomposition 1}, \eqref{comp-decomposition 2} and \eqref{comp-decomposition 3} into the decomposition \eqref{comp-decomposition}, we have established that
	\begin{equation}\label{comp-decomposition O}
		\begin{aligned}
			\hat g_\eps \left(1
			+\frac\eps 2 \int_{\mathbb{R}^3} \hat g_{\eps *} M_*dv_* \right)
			&=
			O\left(\eps\right)_{L^2\left(Mdtdxdv\right)}
			\\
			& +O(1)_{L^\infty\left(dt;L^2\left(dx;\mathcal{C} L^2(Mdv) \right)\right)}
			\\
			&+O(\eps)_{L^\infty\left(dt;L^1\left(dx;\mathcal{C} L_\mathrm{loc}^2(dv) \right)\right)}.
		\end{aligned}
	\end{equation}

	Next, since the left-hand side of the above decomposition degenerates close to vacuum, i.e.\ whenever the density $\frac\eps 2 \int_{\mathbb{R}^3} \hat g_{\eps *} M_*dv_*$ is close to $-1$, we introduce a macroscopic truncation
	\begin{equation*}
		\chi_{\eps,r}(t,x) =\mathds{1}_{ \left\{ 1+\frac\eps 2\int M(v)\hat g_\eps(t,x,v) dv \geq r \right \} },
	\end{equation*}
	for some small $r>0$, thus excluding the domain where this degeneracy is present. It then follows, dividing \eqref{comp-decomposition O} by $1+\frac\eps 2 \int_{\mathbb{R}^3} \hat g_{\eps *} M_*dv_*$, that
	\begin{equation}\label{comp-decomposition*}
		\begin{aligned}
			\chi_{\eps,r} \hat g_\eps
			&=
			O\left(\frac{\eps}{r}\right)_{L^2\left(Mdtdxdv\right)}
			\\
			& +O\left(\frac{1}{r}\right)_{L^\infty\left(dt;L^2\left(dx;\mathcal{C} L^2(Mdv) \right)\right)}
			\\
			&+O\left(\frac{\eps}{r}\right)_{L^\infty\left(dt;L^1\left(dx;\mathcal{C} L_\mathrm{loc}^2(dv) \right)\right)}.
		\end{aligned}
	\end{equation}

	Next, for any small $h\in\mathbb{R}^3$ and any compact subset $K\subset [0,\infty)\times\mathbb{R}^3\times\mathbb{R}^3$, it holds that
	\begin{equation*}
		\begin{aligned}
			\int_K & \chi_{\eps,r } \left| \chi( \lambda \eps \hat g_\eps )\hat g_\eps(t,x,v+h)
			-\chi( \lambda \eps \hat g_\eps )\hat g_\eps (t,x,v)\right|^2 dtdxdv \\
			& = 2
			\int_K \chi_{\eps,r } \chi( \lambda \eps \hat g_\eps )\hat g_\eps(t,x,v)
			\left[ \chi( \lambda \eps \hat g_\eps )\hat g_\eps(t,x,v)
			-\chi( \lambda \eps \hat g_\eps )\hat g_\eps (t,x,v+h)\right] dtdxdv \\
			&\leq C \int_K \chi_{\eps,r }
			\left|\chi( \lambda \eps \hat g_\eps )\hat g_\eps(t,x,v)
			\left(\hat g_\eps(t,x,v+h)-\hat g_\eps(t,x,v)\right)\right|
			dtdxdv.
		\end{aligned}
	\end{equation*}
	Therefore, since
	\begin{equation*}
		\begin{aligned}
			\left\|\chi( \lambda \eps \hat g_\eps )\hat g_\eps\right\|_{L^2_\mathrm{loc}(dtdxdv)} & \leq C,
			\\
			\left\|\chi( \lambda \eps \hat g_\eps )\hat g_\eps\right\|_{L^\infty(dtdxdv)} & \leq \frac{2}{\lambda\eps},
		\end{aligned}
	\end{equation*}
	we conclude from \eqref{comp-decomposition*} that, for any fixed $0<\lambda,r<1$,
	\begin{equation*}
			\limsup_{|h|\rightarrow 0}\limsup_{\eps\rightarrow 0}
			\int_K \chi_{\eps,r } \left| \chi( \lambda \eps \hat g_\eps )\hat g_\eps(t,x,v+h)
			-\chi( \lambda \eps \hat g_\eps )\hat g_\eps (t,x,v)\right|^2 dt dxdv
			= 0,
	\end{equation*}
	which is the expected relative compactness statement away from vacuum on $\chi_{\eps,r}\chi( \lambda \eps \hat g_\eps )\hat g_\eps$. Consequently, combining this result with the control \eqref{micro-truncation} on the very large values of $\hat g_\eps$ yields that, for any given small $r>0$,
	\begin{equation}\label{farfromvacuum}
		\begin{gathered}
			\chi_{\eps,r } \hat g_\eps \text{ is relatively compact}\\
			\text{with respect to the velocity variable in } L_\mathrm{loc}^2(dtdxdv).
		\end{gathered}
	\end{equation}

	\noindent {\bf Near vacuum.} It only remains then to get a compactness estimate near vacuum on $(1-\chi_{\eps,r }) \hat g_\eps$. To this end, we simply decompose, for any given small $\lambda>0$,
	\begin{equation}\label{decomp near vacuum}
		\chi( \lambda \eps \hat g_\eps )\hat g_\eps
		=
		\chi( \lambda \eps \hat g_\eps )\frac{2}{\eps}
		\left(1+\frac{\eps}{2}\hat g_\eps\right)
		+
		\left(1-\chi( \lambda \eps \hat g_\eps )\right)\frac{2}{\eps}
		-
		\frac{2}{\eps},
	\end{equation}
	and we control each term in the right-hand side above individually.

	Thus, noticing that, for any $0<r\leq\frac 12$,
	\begin{equation*}
		\frac{1}{\eps}(1-\chi_{\eps,r})=
		\frac{1}{\eps}\mathds{1}_{\left\{\frac{2}{\eps}<\frac{-\int M\hat g_\eps dv}{1-r}\right\}}
		\leq
		\frac{-\int M\hat g_\eps dv}{2(1-r)}
		\leq
		\left\|\hat g_\eps\right\|_{L^2\left(Mdv\right)}=O(1)_{L^\infty\left(dt;L^2(dx)\right)},
	\end{equation*}
	and that, on the support of $(1-\chi_{\eps,r })$,
	\begin{equation*}
		\begin{aligned}
			\int_{\mathbb{R}^3}
			\left|\chi( \lambda \eps \hat g_\eps )
			\left(1+\frac{\eps}{2}\hat g_\eps\right)
			\right|^2Mdv
			&\leq
			\int_{\mathbb{R}^3}M
			\left(1+\frac{\eps}{2}\hat g_\eps\right)dv
			\left(1+\frac{1}{\lambda}\right)\\
			&\leq r\left(1+\frac{1}{\lambda}\right)=O\left(\frac r\lambda\right)_{L^\infty\left(dtdx\right)},
		\end{aligned}
	\end{equation*}
	we obtain concerning the first term in the right-hand side of \eqref{decomp near vacuum} that
	\begin{equation}\label{decomp near vacuum 1}
		(1-\chi_{\eps,r })\chi( \lambda \eps \hat g_\eps )\frac{2}{\eps}
		\left(1+\frac{\eps}{2}\hat g_\eps\right)
		=O\left(\sqrt{\frac r\lambda}\right)_{L^\infty\left(dt;L^2\left(Mdxdv\right)\right)}.
	\end{equation}

	Then, the second term is easily handled through the estimate
	\begin{equation*}
		\frac{1}{\eps}\left(1-\chi( \lambda \eps \hat g_\eps )\right)
		\leq
		\frac{1}{\eps}\mathds{1}_{\left\{\left|\lambda \eps \hat g_\eps\right|> 1\right\}}
		\leq \lambda \left|\hat g_\eps\right|,
	\end{equation*}
	whereby
	\begin{equation}\label{decomp near vacuum 2}
		(1-\chi_{\eps,r })\left(1-\chi( \lambda \eps \hat g_\eps )\right)\frac{2}{\eps}
		=O\left(\lambda\right)_{L^\infty\left(dt;L^2\left(Mdxdv\right)\right)}.
	\end{equation}

	The remaining term in the right-hand side of \eqref{decomp near vacuum} is constant, in particular it is smooth, and so there is no need to further control it, so that, on the whole, incorporating \eqref{decomp near vacuum 1} and \eqref{decomp near vacuum 2} into \eqref{decomp near vacuum}, we find
	\begin{equation*}
		(1-\chi_{\eps,r })\chi( \lambda \eps \hat g_\eps )\hat g_\eps
		=
		O\left(\sqrt{\frac r\lambda}\right)_{L^\infty\left(dt;L^2\left(Mdxdv\right)\right)}
		+
		O\left(\lambda\right)_{L^\infty\left(dt;L^2\left(Mdxdv\right)\right)}
		-
		\frac{2}{\eps}(1-\chi_{\eps,r }).
	\end{equation*}

	We therefore conclude, for any small $h\in\mathbb{R}^3$ and any compact subset $K\subset [0,\infty)\times\mathbb{R}^3\times\mathbb{R}^3$, that
	\begin{equation}\label{nearvacuum}
		\begin{aligned}
			\limsup_{|h|\rightarrow 0}\limsup_{\eps\rightarrow 0}
			\int_K (1-\chi_{\eps,r }) \left| \chi( \lambda \eps \hat g_\eps )\hat g_\eps(t,x,v+h)
			-\chi( \lambda \eps \hat g_\eps )\hat g_\eps (t,x,v)\right|^2 dt dxdv & \\
			\leq
			O\left(\frac r\lambda\right)+
			O\left(\lambda^2\right). &
		\end{aligned}
	\end{equation}

	\noindent {\bf Conclusion of proof.} On the whole, combining the above estimate \eqref{nearvacuum} near vacuum with the control \eqref{micro-truncation} on the very large values of $\hat g_\eps$ and the compactness statement \eqref{farfromvacuum} away from vacuum, we finally arrive at the control, for any compact subset $K\subset [0,\infty)\times\mathbb{R}^3\times\mathbb{R}^3$,
	\begin{equation*}
		\begin{aligned}
			\limsup_{|h|\rightarrow 0}\limsup_{\eps\rightarrow 0}
			\int_K \left| \hat g_\eps(t,x,v+h)
			-\hat g_\eps (t,x,v)\right|^2 dt dxdv & \\
			\leq
			O\left(\frac r\lambda\right)+
			O\left(\lambda^2\right)&
			+O\left(\frac{1}{|\log\lambda|}\right),
		\end{aligned}
	\end{equation*}
	which, by the arbitrary smallness of $r>0$ and $\lambda>0$, clearly implies that $\hat g_\eps$ is locally relatively compact with respect to the velocity variable in $L^2(dtdxdv)$ and thus concludes the proof of the lemma.
\end{proof}

% ========================
% = xxxxxxxxxxxxxxxxxxxx =
% ========================

\section{Compactness with respect to $x$}\label{hypoelliptic}

In order to get a refined description of the dependence of the fluctuations $g_\eps$ and $g_\eps^\pm$ with respect to $x$, we will use the compactness properties of the free transport operator $v\cdot \nabla_x$. More precisely, there are two types of mechanisms at play here~:
\begin{itemize}
	\item the {\bf transfer of compactness}, which expresses the fact that the free transport mixes the spatial and velocity variables and is a consequence of hypoellipticity~;
	\item the {\bf averaging lemma}, which predicts some regularizing effect for the averages with respect to $v$, due to the fact that the symbol of the free transport is elliptic on a large microlocal subset.
\end{itemize}

Of course both mechanisms require that we have a good control on the advection terms $v\cdot \nabla_x \hat g_\eps$ and $v\cdot\nabla_x \hat g_\eps^\pm$ or some similar quantity (since the square root renormalization is singular at the origin, and thus is not admissible~; see proofs of Lemmas \ref{x-compactness1 0} and \ref{x-compactness2 0} below). In particular, we see at this point that the situation is quite different in the one species scaling leading to the incompressible quasi-static Navier-Stokes-Fourier-Maxwell-Poisson system, and in the two species scaling leading to the two-fluid incompressible Navier-Stokes-Fourier-Maxwell system with (solenoidal) Ohm's law.

In the first case, at the formal level, it is natural to expect from \eqref{VMB-fluct1} that, up to some suitable renormalization, the advection term $v\cdot \nabla_x g_\eps$ is uniformly bounded. However, for the multi-species model, we see that \eqref{VMB-fluct} provides, at least formally, some uniform control on the advection terms $v\cdot \nabla_x g_\eps^\pm$, and therefore some strong compactness on the hydrodynamic variables $\rho_\eps^\pm$, $u_\eps^\pm$ and $\theta_\eps^\pm$, but it does not provide any information on $v\cdot \nabla_x \frac\delta\eps\left(g_\eps^+-g_\eps^-\right)$ which controls the electrodynamic variables $j_\eps$ and $w_\eps$.

As for the limiting systems \eqref{TFINSFMO 2} and \eqref{TFINSFMSO 2}, we therefore do not expect the electromagnetic terms $j \wedge B$ from the Lorentz force to be weakly stable (unfortunately, compensated compactness methods also fail here~; see \cite{arsenio6} for some details on this issue). In this case, we will use, later on in Chapter \ref{entropy method}, some weak-strong stability principle instead of a priori estimates. In other words, the dependence with respect to $x$ is partially understood a posteriori, by comparison with the solutions to the limiting systems.

\subsection{Hypoellipticity and the transfer of compactness}\label{space compactness}

We first explain our global strategy, presenting the main abstract results we will use on the free transport operator.

As mentioned in the introduction of the present chapter, the key idea here is to transfer the compactness with respect to $v$ inherited from the structure of the collision operator (see Lemma \ref{v-compactness}) onto the spatial variable $x$. To this end, we need the following result.

\begin{thm}[\cite{arsenio}]\label{hypo-thm 0}
	Let the bounded family of functions
	\begin{equation*}
		\left\{\phi_\lambda(t,x,v)\right\}_{\lambda\in\Lambda}
		\subset L^p\left(\mathbb{R}_t\times\mathbb{R}^3_x\times\mathbb{R}^3_v\right),
	\end{equation*}
	for some $1<p<\infty$, be locally relatively compact in $v$ and such that
	\begin{equation*}
		\left(\d_t + v\cdot\nabla_x\right) \phi_\lambda =
		\left(1-\Delta_{t,x}\right)^\frac{\beta}{2}\left(1-\Delta_v\right)^\frac{\alpha}{2}S_\lambda,
	\end{equation*}
	for all $\lambda\in\Lambda$ and for some bounded family
	\begin{equation*}
		\left\{S_\lambda(t,x,v)\right\}_{\lambda\in\Lambda}
		\subset L^p\left(\mathbb{R}_t\times\mathbb{R}^3_x\times\mathbb{R}^3_v\right),
	\end{equation*}
	where $\alpha\geq 0$ and $0\leq \beta<1$.
	
	Then, $\left\{\phi_\lambda(t,x,v)\right\}_{\lambda\in\Lambda}$ is locally relatively compact in $L^p\left(\mathbb{R}_t\times\mathbb{R}^3_x\times\mathbb{R}^3_v\right)$ (in all variables).
\end{thm}

The above result was formulated in \cite{arsenio}. It may also be deduced from the methods of \cite{bouchut2} or from the use of standard averaging lemmas from \cite{diperna5} for instance. However, it is to be emphasized that the methods from \cite{arsenio} are more natural and direct.

It turns out that, for the sake of the rigorous derivation of hydrodynamic limits, it is crucial to understand what happens to Theorem \ref{hypo-thm 0} when $p=1$ (carefully note that this case is not covered by the above theorem). To be precise, Theorem \ref{hypo-thm 0} will be sufficient to control oscillations but not concentrations. The basic result in this direction is given by the $L^1$ mixing lemma obtained by Golse and the second author in \cite{golse3}, which allows to transfer equi-integrability from $v$ to $x$ when the source term of the kinetic transport equation is locally integrable.

The point here is that, because of the electromagnetic force, the source term involves derivatives with respect to $v$ and, therefore, is not locally integrable. An analogous situation has been dealt with by the first author in \cite{arsenio3} when considering non-cutoff collision operators, which behave as nonlinear fractional derivatives with respect to $v$. In this singular setting, we are then able to transfer strong compactness, but --~to the best of our knowledge~-- not mere weak compactness, as the results from \cite{golse3} do not apply. More precisely, we have the following statement.

\begin{thm}[\cite{arsenio}]\label{hypo-thm}
	Let the bounded family of non-negative functions
	\begin{equation*}
		\left\{\phi_\lambda(t,x,v)\right\}_{\lambda\in\Lambda}
		\subset L^1\left(\mathbb{R}_t\times\mathbb{R}^3_x; L^r\left(\mathbb{R}^3_v\right)\right),
	\end{equation*}
	for some $1<r<\infty$, be locally relatively compact in $v$ and such that
	\begin{equation*}
		\left(\d_t + v\cdot\nabla_x\right) \phi_\lambda =
		\left(1-\Delta_{t,x}\right)^\frac{\beta}{2}\left(1-\Delta_v\right)^\frac{\alpha}{2}S_\lambda,
	\end{equation*}
	for all $\lambda\in\Lambda$ and for some bounded family
	\begin{equation*}
		\left\{S_\lambda(t,x,v)\right\}_{\lambda\in\Lambda}
		\subset L^{1}\left(\mathbb{R}_t\times\mathbb{R}^3_x; L^r(\mathbb{R}^3_v)\right),
	\end{equation*}
	where $\alpha\geq 0$ and $0\leq \beta<1$.
	
	Then, $\left\{\phi_\lambda(t,x,v)\right\}_{\lambda\in\Lambda}$ is locally relatively compact in $L^{1}\left(\mathbb{R}_t\times\mathbb{R}^3_x\times\mathbb{R}^3_v\right)$ (in all variables).
\end{thm}

The crucial idea behind such hypoelliptic results is that the free transport operator is ``invariant'' by Fourier transform in $(x,v)$, so that frequencies are transported by the semi-group. The argument relies then on a good interpolation formula which expresses both the transport and the elliptic nature of the transport operator away from the characteristic manifold. Nevertheless, because $L^1$ is not a convenient space for Fourier analysis, the proof is quite complex and requires in particular the use of singular integral operators, as well as a characterization of equi-integrability in terms of compactness in weak Hardy spaces. We refer to \cite{arsenio} for a complete discussion of the subject.

\bigskip

Note that, in the problem we consider in this work, the time derivative of the kinetic equations has a factor $\eps$, so that we cannot expect to establish temporal strong compactness and the above theorems cannot be applied as such. However, it is possible to get strong compactness with respect to the fast time variable $\frac t\eps$, but this does not provide any information on the slow dynamics. Thus, we reformulate now the preceding theorems in the following lemmas in order that they be directly applicable to our problem.

\begin{lem}\label{hypo-thm 2}
	Let the bounded family of functions
	\begin{equation*}
		\left\{\phi_\eps(t,x,v)\right\}_{\eps>0}
		\subset L^p\left(\mathbb{R}_t\times\mathbb{R}^3_x\times\mathbb{R}^3_v\right),
	\end{equation*}
	for some $1<p<\infty$, be locally relatively compact in $v$ and such that
	\begin{equation*}
		\left(\eps\d_t + v\cdot\nabla_x\right) \phi_\eps =
		\left(1-\Delta_{x}\right)^\frac{\beta}{2}\left(1-\Delta_v\right)^\frac{\alpha}{2}S_\eps,
	\end{equation*}
	for all $\eps>0$ and for some bounded family
	\begin{equation*}
		\left\{S_\eps(t,x,v)\right\}_{\eps>0}
		\subset L^p\left(\mathbb{R}_t\times\mathbb{R}^3_x\times\mathbb{R}^3_v\right),
	\end{equation*}
	where $\alpha\geq 0$ and $0\leq \beta<1$.
	
	Then, $\left\{\phi_\eps(t,x,v)\right\}_{\eps>0}$ is locally relatively compact in $L^p\left(\mathbb{R}_t\times\mathbb{R}^3_x\times\mathbb{R}^3_v\right)$ in $x$ and $v$ (but not necessarily in $t$).
\end{lem}

\begin{proof}
	This result is directly deduced from Theorem \ref{hypo-thm 0} (and its proof).
	To this end, we define
	\begin{equation*}
		\begin{aligned}
			\tilde\phi_\eps(t,x,v) & = \eps^\frac{1}{p} \phi_\eps(\eps t, x, v), \\
			\tilde S_\eps(t,x,v) & = \eps^{\frac 1p} S_\eps(\eps t, x, v),
		\end{aligned}
	\end{equation*}
	so that
	\begin{equation*}
		\left(\d_t + v\cdot\nabla_x\right) \tilde \phi_\eps =
		\left(1-\Delta_{x}\right)^\frac{\beta}{2}\left(1-\Delta_v\right)^\frac{\alpha}{2} \tilde S_\eps,
	\end{equation*}
	and $\tilde\phi_\eps$ and $\tilde S_\eps$ are uniformly bounded in $L^p\left(\mathbb{R}_t\times\mathbb{R}^3_x\times\mathbb{R}^3_v\right)$.
	
	We apply now Theorem \ref{hypo-thm 0} to the above transport equation to deduce that $\left\{\tilde \phi_\eps\right\}_{\eps>0}$ is locally relatively compact in $L^p\left(\mathbb{R}_t\times\mathbb{R}^3_x\times\mathbb{R}^3_v\right)$. In fact, a closer inspection of the proof of this theorem in \cite{arsenio} reveals that, by possibly localizing without loss of generality the above functions in $v$ only, one has the following global estimate~:
	\begin{equation*}
		\lim_{\gamma\rightarrow 0} \sup_{\eps>0} \sup_{ \left|k\right|+\left|h\right|+\left|l\right|<\gamma }
		\left\|
		\tilde\phi_\eps(t+k,x+h,v+l) - \tilde\phi_\eps(t,x,v)
		\right\|_{L^p\left(\mathbb{R}_t\times\mathbb{R}^3_x\times\mathbb{R}^3_v\right)} = 0.
	\end{equation*}
	It follows that
	\begin{equation*}
		\lim_{\gamma\rightarrow 0} \sup_{\eps>0} \sup_{ \left|h\right|+\left|l\right|<\gamma }
		\left\|
		\phi_\eps(t,x+h,v+l) - \phi_\eps(t,x,v)
		\right\|_{L^p\left(\mathbb{R}_t\times\mathbb{R}^3_x\times\mathbb{R}^3_v\right)} = 0,
	\end{equation*}
	which concludes the proof of the lemma.
\end{proof}

\begin{lem}\label{HYPO-THM 3}
	Let the bounded family of non-negative functions
	\begin{equation*}
		\left\{\phi_\eps(t,x,v)\right\}_{\eps>0}
		\subset L^1\left(\mathbb{R}_t\times\mathbb{R}^3_x; L^r\left(\mathbb{R}^3_v\right)\right),
	\end{equation*}
	for some $1<r<\infty$, be locally relatively compact in $v$ and such that
	\begin{equation*}
		\left(\eps\d_t + v\cdot\nabla_x\right) \phi_\eps =
		\left(1-\Delta_{x}\right)^\frac{\beta}{2}\left(1-\Delta_v\right)^\frac{\alpha}{2}S_\eps,
	\end{equation*}
	for all $\eps>0$ and for some bounded family
	\begin{equation*}
		\left\{S_\eps(t,x,v)\right\}_{\eps>0}
		\subset L^{1}\left(\mathbb{R}_t\times\mathbb{R}^3_x; L^r(\mathbb{R}^3_v)\right),
	\end{equation*}
	where $\alpha\geq 0$ and $0\leq \beta<1$. We further assume that, for any compact set $K\subset\mathbb{R}^3\times\mathbb{R}^3$,
	\begin{equation*}
		\left\{\int_{K}\phi_\eps(t,x,v) dxdv\right\}_{\eps>0}\quad \text{is equi-integrable (in $t$).}
	\end{equation*}
	
	Then, $\left\{\phi_\eps(t,x,v)\right\}_{\eps>0}$ is equi-integrable (in all variables) and locally relatively compact in $L^{1}\left(\mathbb{R}_t\times\mathbb{R}^3_x\times\mathbb{R}^3_v\right)$ in $x$ and $v$ (but not necessarily in $t$).
	
	Moreover, if the $\phi_\eps$'s are signed (in the sense that the functions may assume both positive and negative values), the conclusion still holds true, i.e.\ $\left\{\phi_\eps(t,x,v)\right\}_{\eps>0}$ is locally relatively compact in $L^{1}\left(\mathbb{R}_t\times\mathbb{R}^3_x\times\mathbb{R}^3_v\right)$ in $x$ and $v$ (but not necessarily in $t$), provided $\left\{\phi_\eps(t,x,v)\right\}_{\eps>0}$ is equi-integrable (in all variables) a priori.
\end{lem}

\begin{proof}
	When the $\phi_\eps$'s are signed and a priori equi-integrable (in all variables), this result is deduced from Theorem \ref{hypo-thm} (and its proof) utilizing the strategy of proof of Lemma \ref{hypo-thm 2}, that is by dilation of the time variable.
	To this end, we define
	\begin{equation*}
		\begin{aligned}
			\tilde\phi_\eps(t,x,v) & = \eps \phi_\eps(\eps t, x, v), \\
			\tilde S_\eps(t,x,v) & = \eps S_\eps(\eps t, x, v),
		\end{aligned}
	\end{equation*}
	so that
	\begin{equation*}
		\left(\d_t + v\cdot\nabla_x\right) \tilde \phi_\eps =
		\left(1-\Delta_{x}\right)^\frac{\beta}{2}\left(1-\Delta_v\right)^\frac{\alpha}{2} \tilde S_\eps,
	\end{equation*}
	and $\tilde\phi_\eps$ and $\tilde S_\eps$ are uniformly bounded in $L^1\left(\mathbb{R}_t\times\mathbb{R}^3_x; L^r\left(\mathbb{R}^3_v\right)\right)$.
	
	We apply now Theorem \ref{hypo-thm} to the above transport equation to deduce that $\left\{\tilde \phi_\eps\right\}_{\eps>0}$ is locally relatively compact in $L^1\left(\mathbb{R}_t\times\mathbb{R}^3_x\times\mathbb{R}^3_v\right)$. In fact, a closer inspection of the proof of this theorem in \cite{arsenio} reveals that, by possibly localizing without loss of generality the above functions in $v$ only, one has the following global estimate~:
	\begin{equation*}
		\lim_{\gamma\rightarrow 0} \sup_{\eps>0} \sup_{ \left|k\right|+\left|h\right|+\left|l\right|<\gamma }
		\left\|
		\tilde\phi_\eps(t+k,x+h,v+l) - \tilde\phi_\eps(t,x,v)
		\right\|_{L^{1,\infty}\left(\mathbb{R}_t\times\mathbb{R}^3_x\times\mathbb{R}^3_v\right)} = 0,
	\end{equation*}
	where $L^{1,\infty}$ denotes the standard weak Lebesgue space (or Lorentz space). Note that $L^{1,\infty}$ has the same homogeneity as the Lebesgue space $L^1$.
	It follows that
	\begin{equation*}
		\lim_{\gamma\rightarrow 0} \sup_{\eps>0} \sup_{ \left|h\right|+\left|l\right|<\gamma }
		\left\|
		\phi_\eps(t,x+h,v+l) - \phi_\eps(t,x,v)
		\right\|_{L^{1,\infty}\left(\mathbb{R}_t\times\mathbb{R}^3_x\times\mathbb{R}^3_v\right)} = 0.
	\end{equation*}
	
	Next, for any compact set $K\subset \mathbb{R}_t\times\mathbb{R}^3_x\times\mathbb{R}^3_v$ and any large $R>1$, we have that
	\begin{equation*}
		\begin{aligned}
			\|
			\phi_\eps & (t,x+h,v+l) - \phi_\eps(t,x,v)
			\|_{L^{1}\left(K\right)} \\
			& = \int_0^\infty \left|\set{(t,x,v)\in K}{\left|\phi_\eps(t,x+h,v+l) - \phi_\eps(t,x,v)\right|>\lambda}\right| d\lambda \\
			& \leq \int_{\frac 1R}^R \left|\set{(t,x,v)\in K}{\left|\phi_\eps(t,x+h,v+l) - \phi_\eps(t,x,v)\right|>\lambda}\right| d\lambda
			+ \frac{|K|}{R} \\
			& + \left\| \left(\phi_\eps(t,x+h,v+l) - \phi_\eps(t,x,v)\right)
			\mathds{1}_{\left\{\left|\phi_\eps(t,x+h,v+l) - \phi_\eps(t,x,v)\right|>R\right\}} \right\|_{L^1\left(K\right)} \\
			& \leq 2\log R \left\|
			\phi_\eps(t,x+h,v+l) - \phi_\eps(t,x,v)
			\right\|_{L^{1,\infty}\left(K\right)}
			+ \frac{|K|}{R} \\
			& + \left\| \left(\phi_\eps(t,x+h,v+l) - \phi_\eps(t,x,v)\right)
			\mathds{1}_{\left\{\left|\phi_\eps(t,x+h,v+l) - \phi_\eps(t,x,v)\right|>R\right\}} \right\|_{L^1\left(K\right)}.
		\end{aligned}
	\end{equation*}
	Hence, we deduce, provided the $\phi_\eps$'s are equi-integrable in all variables and by the arbitrariness of $R>1$, that
	\begin{equation*}
		\lim_{\gamma\rightarrow 0} \sup_{\eps>0} \sup_{ \left|h\right|+\left|l\right|<\gamma }
		\left\|
		\phi_\eps(t,x+h,v+l) - \phi_\eps(t,x,v)
		\right\|_{L^{1}_\mathrm{loc}\left(\mathbb{R}_t\times\mathbb{R}^3_x\times\mathbb{R}^3_v\right)} = 0,
	\end{equation*}
	which concludes the proof of the lemma when the $\phi_\eps$'s are signed and a priori equi-integrable.

	Therefore, there only remains to establish the equi-integrability of $\left\{\phi_\eps(t,x,v)\right\}_{\eps>0}$ when it is not already known a priori and when each $\phi_\eps$ is non-negative. However, the preceding strategy based on time-dilations to deduce results from Theorem \ref{hypo-thm} cannot be repeated here, for the notion of equi-integrability does unfortunately not behave suitably under partial dilations. Instead, the proof of Theorem 2.4 from \cite{arsenio} has to be adapted to treat the present setting, which is rather involved. Therefore, in order to provide a self-contained justification based on \cite{arsenio} and for the sake of clarity, we have moved the remainder of the proof of the present lemma to Appendix \ref{endofproof}.
\end{proof}

\subsection{Compactness of fluctuations for one species}\label{compactness one species}

The next step consists in combining the velocity compactness result from Lemma \ref{v-compactness} with the hypoelliptic transfer of compactness contained in Lemma \ref{HYPO-THM 3} to infer the compactness in $x$ and $v$ (but not in $t$) of the fluctuations $\hat g_\eps$ and $\hat g_\eps^\pm$. To this end, we will need to consider the action of the transport operator $\left(\eps\partial_t+v\cdot\nabla_x\right)$ on the fluctuations $\hat g_\eps$, $\hat g_\eps^\pm$ (to control oscillations) and their square $\hat g_\eps^2$, $\hat g_\eps^{\pm 2}$ (to control concentrations), or truncated versions of these fluctuations.

Again, note that this strategy differs from the methods developed in previous works on hydrodynamic limits, since we do not use classical averaging lemma. Indeed, we prove below that the fluctuations themselves, and not only their moments with respect to $v$, are strongly compact in $x$ and $v$.

Let us first focus on the regime considered in Theorem \ref{NS-WEAKCV} (with one species) leading to the incompressible quasi-static Navier-Stokes-Fourier-Maxwell-Poisson system \eqref{NSFMP 2}. In this case, we have the following lemma.

\begin{lem}\label{x-compactness1 0}
	Let $\left(f_\eps, E_\eps, B_\eps\right)$ be the sequence of renormalized solutions to the scaled one species Vlasov-Maxwell-Boltzmann system \eqref{VMB1} considered in Theorem \ref{NS-WEAKCV}.
	
	Then, as $\eps\rightarrow 0$, any subsequence of renormalized fluctuations $\hat g_\eps$ is locally relatively compact in $(x,v)$ in $L^2\left(dtdxdv\right)$ in the sense that, for any $\eta>0$ and every compact subset $K\subset [0,\infty)\times\mathbb{R}^3\times\mathbb{R}^3$, there exists $\gamma>0$ such that, if $h,l\in\mathbb{R}^3$ satisfy $|h|+|l|<\gamma$, then
	\begin{equation*}
		\sup_{\eps>0} \left\|\hat g_\eps(t,x+h,v+l) - \hat g_\eps(t,x,v)\right\|_{L^2\left(K,dtdxdv\right)}
		<\eta.
	\end{equation*}
	 In particular, the family $|\hat g_\eps|^2$ is equi-integrable (in all variables $t$, $x$ and $v$).
\end{lem}

\begin{proof}
	The proof of this lemma proceeds with two main steps. The first one establishes the compactness of $\hat g_\eps$ in $x$ and $v$ in $L^1_\mathrm{loc}$, while the second one shows the equi-integrability of $\hat g_\eps^2$ in all variables. The combination of these two steps will eventually allow us to conclude the proof.

	\noindent {\bf An admissible renormalization $\beta_1(G_\eps)$.} We consider first the admissible square root renormalization
	\begin{equation*}
		\beta_1(z)=\frac{\sqrt{z+\eps^a}-1}{\eps},
	\end{equation*}
	for some given $1<a<4$. This renormalization is introduced to circumvent the fact that the natural renormalization $2\frac{\sqrt{z}-1}{\eps}$ corresponding to $\hat g_\eps$ is not admissible for Vlasov-Boltzmann equations, for it is singular at $z=0$, i.e.\ $\left(2\frac{\sqrt{z}-1}{\eps}\right)' \to \infty$ as $z\to 0$.

	In fact, we have already used a similar strategy in the proof of Proposition \ref{weak-comp}, where we showed, as a consequence of the entropy and entropy dissipation bounds, that (see \eqref{error-alpha})
	\begin{equation}\label{error-alpha 5}
		\begin{aligned}
		2\frac{\sqrt{G_\eps+\eps^a}-1}{\eps}-\hat g_\eps
		& = O\left(\eps^{\frac a2 -1}\right)_{L^\infty(dtdxdv)},
		\\
		2\frac{\sqrt{G_\eps+\eps^a}-1}{\eps}-\hat g_\eps
		& = O\left(\eps^{a-1}\right)_{L^\infty(dtdxdv)} +O\left(\eps^{\frac a2}\right)_{L^\infty\left(dt;L^2\left(M dxdv\right)\right)},
		\\
		2\frac{\sqrt{G_\eps+\eps^a}-1}{\eps}-\hat g_\eps
		& = O\left(\eps^{a-1}\right)_{L^\infty(dtdxdv)} +O\left(\eps^{\frac a2}\right)_{L^2_\mathrm{loc}\left(dtdx;L^2\left(\left(1+|v|^2\right)Mdv\right)\right)},
		\end{aligned}
	\end{equation}
	so that strong compactness properties of $\frac{\sqrt{G_\eps+\eps^a}-1}{\eps}$ in $L^p_{\mathrm{loc}}\left(dtdxdv\right)$, for any given $1\leq p\leq 2$, will entail similar properties on $\hat g_\eps$ in the same space and vice versa.

	\noindent {\bf Compactness of $\beta_1(G_\eps)$ in $v$.} In particular, since, in view of Lemma \ref{v-compactness}, the renormalized fluctuations $\hat g_\eps$ are locally relatively compact in $v$ in $L^2_\mathrm{loc}(dtdxdv)$, the same holds true for any subsequence of $\frac{\sqrt{G_\eps+\eps^a}-1}{\eps}$ in the sense that, for any $\eta>0$ and every compact subset $K\subset [0,\infty)\times\mathbb{R}^3\times\mathbb{R}^3$, there exists $\gamma>0$ such that, if $h\in\mathbb{R}^3$ satisfies $|h|<\gamma$, then
	\begin{equation}\label{beta 1 v comp}
		\sup_{\eps>0} \left\|\beta_1\left(G_\eps\right)(t,x,v+h) - \beta_1\left(G_\eps\right)(t,x,v)\right\|_{L^2\left(K,dtdxdv\right)}
		<\eta.
	\end{equation}

	\noindent {\bf Action of the transport operator on $\beta_1(G_\eps)$.} Thus, using $\beta_1(z)$ to renormalize the Vlasov-Boltzmann equation in \eqref{VMB1} and decomposing the collision integrands according to \eqref{integrands-decomposition}, we find that (see \eqref{renormalized})
	\begin{equation}\label{beta 1 transport 0}
		\begin{aligned}
			( \eps \d_t + & v \cdot \nabla_x + \eps \left(E_\eps+ v\wedge B_\eps\right)\cdot \nabla_v )
			{\sqrt{ G_\eps+\eps^a }-1 \over \eps}
			- E_\eps \cdot v {G_\eps \over 2\sqrt{ G_\eps+\eps^a }} \\
			& =
			{\sqrt{G_\eps} \over 2\sqrt{ G_\eps+\eps^a }}
			\int_{\mathbb{R}^3\times\mathbb{S}^2} \sqrt{ G_{\eps *} } \hat q_\eps M_*dv_*d\sigma
			+
			{\eps^2 \over 8\sqrt{G_\eps+\eps^a  }}
			\int_{\mathbb{R}^3\times\mathbb{S}^2} \hat q_\eps^2 M_*dv_*d\sigma.
		\end{aligned}
	\end{equation}
	It follows, employing the uniform bounds $\hat g_\eps\in L^\infty\left(dt;L^2\left(Mdxdv\right)\right)$ and $\hat q_\eps\in L^2\left(MM_*dtdxdvdv_*d\sigma\right)$ from Lemmas \ref{L2-lem} and \ref{L2-qlem}, respectively, that (see \eqref{Q 1})
	\begin{equation}\label{beta 1 transport}
		\begin{aligned}
			( \eps \d_t + & v \cdot \nabla_x) {\sqrt{ G_\eps+\eps^a }-1 \over \eps}
			\\
			& =
			{\sqrt{G_\eps} \over 2\sqrt{ G_\eps+\eps^a }}
			E_\eps \cdot v \left(1+\frac\eps 2\hat g_\eps\right)
			-
			\eps\nabla_v\cdot
			\left(E_\eps+ v\wedge B_\eps\right)
			{\sqrt{ G_\eps+\eps^a }-1 \over \eps}
			\\
			&
			+ {\sqrt{G_\eps} \over 2\sqrt{ G_\eps+\eps^a }}
			\int_{\mathbb{R}^3\times\mathbb{S}^2} \hat q_\eps M_*dv_*d\sigma
			+
			{\eps \sqrt{G_\eps} \over 4\sqrt{ G_\eps+\eps^a }}
			\int_{\mathbb{R}^3\times\mathbb{S}^2} \hat g_{\eps *} \hat q_\eps M_*dv_*d\sigma
			\\
			& +
			{\eps^2 \over 8\sqrt{G_\eps+\eps^a  }}
			\int_{\mathbb{R}^3\times\mathbb{S}^2} \hat q_\eps^2 M_*dv_*d\sigma
			\\
			& = O(1)_{L^1_{\mathrm{loc}}\left(dtdxdv\right)} + O(\eps)_{L^1_\mathrm{loc}\left(dtdx ; W^{-1,1}_{\mathrm{loc}}\left(dv\right)\right)}.
		\end{aligned}
	\end{equation}

	\noindent {\bf Compactness of $\beta_1(G_\eps)$ in $(x,v)$.} On the whole, we have established the compactness in velocity of $\beta_1(G_\eps)$ in \eqref{beta 1 v comp} and a bound on the transport operator acting on $\beta_1(G_\eps)$ in \eqref{beta 1 transport}. Therefore, a direct application of Lemma \ref{HYPO-THM 3} yields that $\beta_1(G_\eps)=\frac{\sqrt{G_\eps+\eps^a}-1}{\eps}$ is locally relatively compact in $(x,v)$ in $L^1_\mathrm{loc}\left(dtdxdv\right)$.

	Combining this result with \eqref{error-alpha 5}, implies that the renormalized fluctuations $\hat g_\eps$ are relatively compact in $(x,v)$ in $L_{\mathrm{loc}}^1\left(dtdxdv\right)$ as well, in the sense that, for any $\eta>0$ and every compact subset $K\subset [0,\infty)\times\mathbb{R}^3\times\mathbb{R}^3$, there exists $\gamma>0$ such that, if $h,l\in\mathbb{R}^3$ satisfy $|h|+|l|<\gamma$, then
	\begin{equation}\label{L1 xv compactness}
		\sup_{\eps>0} \left\|\hat g_\eps(t,x+h,v+l) - \hat g_\eps(t,x,v)\right\|_{L^1\left(K,dtdxdv\right)}
		<\eta.
	\end{equation}

	\noindent {\bf An admissible renormalization $\beta_2(G_\eps)$.} In order to improve this local strong compactness in $x$ and $v$ from $L^1_\mathrm{loc}$ to $L^2_{\mathrm{loc}}$, we only have to show now that $\hat g_\eps^2$ is locally equi-integrable in all variables, which will also be seen as a consequence of Lemma \ref{HYPO-THM 3}.

	To this end, we consider now the admissible renormalization
	\begin{equation*}
		\beta_2(z) = \left(\frac{\sqrt{z+\eps^a}-1}{\eps}\right)^2
		\gamma\left(\lambda\left(\sqrt{z+\eps^a}-1\right)\right),
	\end{equation*}
	where $\gamma(z)\in C^1\left(\mathbb{R}\right)$ is a cutoff satisfying $\mathds{1}_{[-1,1]}\leq\gamma(z)\leq\mathds{1}_{[-2,2]}$, for some given $1<a <4$ (in fact, we will further restrict the range of $a$ so that necessarily $a=2$ below) and any small enough $\lambda>0$. As before, this renormalization is introduced to circumvent the fact that the natural renormalization $\left(2\frac{\sqrt{z}-1}{\eps}\right)^2$ corresponding to $\hat g_\eps^2$ is not admissible for Vlasov-Boltzmann equations, for it is singular at $z=0$, i.e.\ $\left[\left(2\frac{\sqrt{z}-1}{\eps}\right)^2\right]' \to -\infty$ as $z\to 0$, and its growth at infinity is not admissible, i.e.\ $\left(2\frac{\sqrt{z}-1}{\eps}\right)^2\sim z$ as $z\to \infty$.

	Next, it is readily seen that \eqref{error-alpha 5} implies that
	\begin{equation}\label{error-alpha 6}
		\begin{aligned}
		\left(2\frac{\sqrt{G_\eps+\eps^a}-1}{\eps}\right)^2 - \hat g_\eps^2
		& =
		2\hat g_\eps\left(2\frac{\sqrt{G_\eps+\eps^a}-1}{\eps} - \hat g_\eps\right)
		+
		\left(2\frac{\sqrt{G_\eps+\eps^a}-1}{\eps} - \hat g_\eps\right)^2
		\\
		& = O\left(\eps^{a-1}+\eps^\frac a2\right)_{L^\infty\left(dt;L^1_\mathrm{loc}\left(dx;L^1\left(Mdv\right)\right)\right)},
		\end{aligned}
	\end{equation}
	so that the equi-integrability of $\beta_2(G_\eps)=\left(\frac{\sqrt{G_\eps+\eps^a}-1}{\eps}\right)^2\gamma\left(\lambda\left(\sqrt{G_\eps+\eps^a}-1\right)\right)$ will entail the equi-integrability of $\hat g_\eps^2\gamma \left(\lambda\left(\sqrt{G_\eps+\eps^a}-1\right)\right)$, which, when combined with the following control on the very large values of fluctuations (see \eqref{large values})~:
	\begin{equation}\label{large values 2}
		\begin{aligned}
			\left(1-\gamma\right)
			\left(\lambda\left(\sqrt{G_\eps+\eps^a}-1\right)\right)
			\left|\hat g_\eps\right|^2
			& \leq
			\mathds{1}_{\left\{\lambda\left(\sqrt{G_\eps+\eps^a}-1\right) > 1\right\}}
			\left|\hat g_\eps\right|^2
			% \\
			% & =
			% \mathds{1}_{\left\{{G_\eps} > \left(\frac{1+\lambda}{\lambda}\right)^2-\eps^a\right\}}
			% \left|\hat g_\eps\right|^2
			\\
			& =O\left({1\over \left|\log \lambda\right|}\right)_{L^\infty\left(dt;L^1\left(Mdxdv\right)\right)},
		\end{aligned}
	\end{equation}
	for any small enough $\lambda>0$, will eventually imply the equi-integrability of $\hat g_\eps^2$.

	\noindent {\bf Compactness of $\beta_2(G_\eps)$ in $v$.} Furthermore, note that the velocity compactness stated in \eqref{beta 1 v comp} implies a corresponding property for $\beta_2(G_\eps)$. Indeed, expressing $\beta_2(z)=\beta_1(z)^2\gamma\left(\eps\lambda \beta_1(z)\right)$ and noticing that, for any $z_1,z_2\in\mathbb{R}$,
	\begin{equation}\label{gamma}
		\begin{aligned}
			\left|z_1^2\gamma(\eps\lambda z_1)-z_2^2\gamma(\eps\lambda z_2)\right|
			& =
			\left|(z_1-z_2)z_1\gamma(\eps\lambda z_1)+z_2\left(z_1\gamma(\eps\lambda z_1)-z_2\gamma(\eps\lambda z_2)\right)\right|
			\\
			& \leq
			C|z_1-z_2|\left(|z_1|+|z_2|\right),
		\end{aligned}
	\end{equation}
	we deduce, for any $h\in\mathbb{R}^3$, that
	\begin{equation*}
		\begin{aligned}
			& \left\|\beta_2\left(G_\eps\right)(t,x,v+h) - \beta_2\left(G_\eps\right)(t,x,v)\right\|_{L^1_\mathrm{loc}\left(dtdxdv\right)}
			\\
			& \leq C
			\left\|\beta_1\left(G_\eps\right)\right\|_{L^2_\mathrm{loc}\left(dtdxdv\right)}
			\left\|\beta_1\left(G_\eps\right)(t,x,v+h) - \beta_1\left(G_\eps\right)(t,x,v)\right\|_{L^2_\mathrm{loc}\left(dtdxdv\right)}.
		\end{aligned}
	\end{equation*}
	It then follows from \eqref{beta 1 v comp} that any subsequence of $\beta_2(G_\eps)$ is locally relatively compact in $v$ in $L^1\left(dtdxdv\right)$ in the sense that, for any $\eta>0$ and every compact subset $K\subset [0,\infty)\times\mathbb{R}^3\times\mathbb{R}^3$, there exists $\gamma>0$ such that, if $h\in\mathbb{R}^3$ satisfies $|h|<\gamma$, then
	\begin{equation}\label{beta 2 v comp 0}
		\sup_{\eps>0} \left\|\beta_2(G_\eps)(t,x,v+h) - \beta_2(G_\eps)(t,x,v)\right\|_{L^1\left(K,dtdxdv\right)}
		<\eta.
	\end{equation}

	However, in order to use the above velocity compactness of $\beta_2(G_\eps)$ in Lemma \ref{HYPO-THM 3}, we still need to show that $\beta_2(G_\eps)$ enjoys an improved integrability with respect to the velocity variable, namely that $\beta_2(G_\eps)$ is locally bounded in $L^1\left(dtdx;L^r(dx)\right)$, for some $r>1$. To this end, we introduce the following decomposition
	\begin{equation*}
		\begin{aligned}
			\beta_2(G_\eps) & = \frac12 \Pi \hat g_\eps  {\sqrt{G_\eps+\eps^a }-1\over \eps}
			\gamma \left(\lambda \left( \sqrt{G_\eps+\eps^a }-1 \right)\right)
			\\
			& + \left({\sqrt{G_\eps+\eps^a }-1 \over \eps}
			- \frac 12 \hat g_\eps \right)
			{\sqrt{G_\eps+\eps^a }-1\over \eps}
			\gamma \left( \lambda \left( \sqrt{G_\eps+\eps^a }-1 \right)\right)
			\\
			& + \frac12( \hat g_\eps -\Pi \hat g_\eps)
			{\sqrt{G_\eps+\eps^a }-1\over \eps}
			\gamma \left( \lambda \left( \sqrt{G_\eps+\eps^a }-1 \right)\right).
		\end{aligned}
	\end{equation*}
	Since $\Pi \hat g_\eps$ belongs to $L^\infty\left(dt;L^2\left(dx;L^p\left(Mdv\right)\right)\right)$, for any $1\leq p<\infty$, we therefore get, for any $1\leq r <2$,
	\begin{equation*}
		\beta_2(G_\eps) \leq O(1)_{L^1_\mathrm{loc}\left(dtdx;L^r\left(Mdv\right)\right)}
		+\frac{C}{\lambda\eps}
		\left|{\sqrt{G_\eps+\eps^a }-1 \over \eps}
		- \frac 12 \hat g_\eps \right|
		+ \frac C{\lambda\eps}\left| \hat g_\eps -\Pi \hat g_\eps \right|.
	\end{equation*}
	Then, by \eqref{error-alpha 5} and by the relaxation estimate \eqref{relaxation-est}, we obtain, provided $2\leq a<4$,
	\begin{equation*}
		\begin{aligned}
			\beta_2(G_\eps) & = O(1)_{L^1_\mathrm{loc}\left(dtdx;L^r\left(Mdv\right)\right)}
			\\
			& + O\left(\frac {\eps^\frac{a-2}{2}}\lambda\right)_{L^1_\mathrm{loc}\left(dtdx;L^2\left(Mdv\right)\right)}
			+ O\left(\frac 1\lambda\right)_{L^1_\mathrm{loc}\left(dtdx;L^2\left(Mdv\right)\right)}
			\\
			& =
			O\left(\frac 1\lambda\right)_{L^1_\mathrm{loc}\left(dtdx;L^r(Mdv)\right)},
		\end{aligned}
	\end{equation*}
	for any $1\leq r<2$.

	Finally, combining the preceding estimate with the compactness estimate \eqref{beta 2 v comp 0}, we deduce, for any $1\leq r<2$, that $\beta_2(G_\eps)$ is locally relatively compact in $v$ in $L^1\left(dtdx;L^r(dv)\right)$ in the sense that, for any $\eta>0$ and every compact subset $K\subset [0,\infty)\times\mathbb{R}^3\times\mathbb{R}^3$, there exists $\gamma>0$ such that, if $h\in\mathbb{R}^3$ satisfies $|h|<\gamma$, then
	\begin{equation}\label{beta 2 v comp}
		\sup_{\eps>0} \left\|\left(\beta_2(G_\eps)(t,x,v+h) - \beta_2(G_\eps)(t,x,v)\right)\mathds{1}_{K}(t,x,v)\right\|_{L^1\left(dtdx;L^r(dv)\right)}
		<\eta.
	\end{equation}

	\noindent {\bf Action of the transport operator on $\beta_2(G_\eps)$.} Thus, using $\beta_2(z)$ to renormalize the Vlasov-Boltzmann equation in \eqref{VMB1}, decomposing the collision integrands according to \eqref{integrands-decomposition} and writing for convenience
	\begin{equation*}
		\Gamma(z)=z\left(2\gamma(\lambda\eps z)+\lambda\eps z\gamma'(\lambda\eps z)\right),
	\end{equation*}
	so that $\beta_2'(z)=\beta_1'(z)\Gamma\left(\beta_1(z)\right)$, we have now that (note that this renormalization procedure amounts to multiplying \eqref{beta 1 transport 0} by $\Gamma\left(\beta_1(G_\eps)\right)$)
	\begin{equation*}
		\begin{aligned}
			( \eps & \d_t + v \cdot \nabla_x + \eps \left(E_\eps+ v\wedge B_\eps\right)\cdot \nabla_v )
			\beta_2(G_\eps)
			- E_\eps \cdot v {G_\eps\Gamma\left(\beta_1(G_\eps)\right) \over 2\sqrt{ G_\eps+\eps^a }} \\
			& =
			{\sqrt{G_\eps} \Gamma\left(\beta_1(G_\eps)\right) \over 2\sqrt{ G_\eps+\eps^a }}
			\int_{\mathbb{R}^3\times\mathbb{S}^2} \sqrt{ G_{\eps *} } \hat q_\eps M_*dv_*d\sigma
			+
			{\eps^2 \Gamma\left(\beta_1(G_\eps)\right) \over 8\sqrt{G_\eps+\eps^a  }}
			\int_{\mathbb{R}^3\times\mathbb{S}^2} \hat q_\eps^2 M_*dv_*d\sigma.
		\end{aligned}
	\end{equation*}
	It follows, employing the uniform bounds $\hat g_\eps\in L^\infty\left(dt;L^2\left(Mdxdv\right)\right)$ and $\hat q_\eps\in L^2\left(MM_*dtdxdvdv_*d\sigma\right)$ from Lemmas \ref{L2-lem} and \ref{L2-qlem}, respectively, and the direct estimates
	\begin{equation*}
		\begin{aligned}
			\left|\Gamma\left(\beta_1(G_\eps)\right)\right| & \leq C \left|\beta_1\left(G_\eps\right)\right|
			=O\left(1\right)_{L^2_\mathrm{loc}(dtdxdv)},
			\\
			\Gamma\left(\beta_1(G_\eps)\right) & =O\left(\frac{1}{\lambda\eps}\right)_{L^\infty(dtdxdv)},
			\\
			\left|\beta_2(G_\eps)\right| & \leq \frac{C}{\lambda\eps} \left|\beta_1\left(G_\eps\right)\right|
			=O\left(\frac 1 {\lambda\eps}\right)_{L^2_\mathrm{loc}(dtdxdv)},
		\end{aligned}
	\end{equation*}
	that, provided $1<a\leq 2$,
	\begin{equation}\label{beta 2 transport}
		\begin{aligned}
			( \eps \d_t + & v \cdot \nabla_x) \beta_2\left(G_\eps\right)
			\\
			& =
			{\sqrt{G_\eps} \Gamma\left(\beta_1(G_\eps)\right) \over 2\sqrt{ G_\eps+\eps^a }}
			E_\eps \cdot v \left(1+\frac\eps 2\hat g_\eps\right)
			\\
			& -
			\eps\nabla_v\cdot
			\left[\left(E_\eps+ v\wedge B_\eps\right)
			\beta_2\left(G_\eps\right)\right]
			\\
			&
			+ {\sqrt{G_\eps} \Gamma\left(\beta_1(G_\eps)\right) \over 2\sqrt{ G_\eps+\eps^a }}
			\int_{\mathbb{R}^3\times\mathbb{S}^2} \hat q_\eps M_*dv_*d\sigma
			\\
			& +
			{\eps \sqrt{G_\eps} \Gamma\left(\beta_1(G_\eps)\right) \over 4\sqrt{ G_\eps+\eps^a }}
			\int_{\mathbb{R}^3\times\mathbb{S}^2} \hat g_{\eps *} \hat q_\eps M_*dv_*d\sigma
			+
			{\eps^2 \Gamma\left(\beta_1(G_\eps)\right) \over 8\sqrt{G_\eps+\eps^a  }}
			\int_{\mathbb{R}^3\times\mathbb{S}^2} \hat q_\eps^2 M_*dv_*d\sigma
			\\
			& = % O\left(\frac 1\lambda\right)_{L^1_{\mathrm{loc}}\left(dtdxdv\right)} +
			O\left(\frac 1\lambda\right)_{L^1_\mathrm{loc}\left(dtdx ; W^{-1,1}_{\mathrm{loc}}\left(dv\right)\right)}.
		\end{aligned}
	\end{equation}

	\noindent {\bf Equi-integrability of $\beta_2(G_\eps)$ in $(t,x,v)$.} On the whole, we have established the compactness in velocity of $\beta_2(G_\eps)$ in \eqref{beta 2 v comp} and a bound on the transport operator acting on $\beta_2(G_\eps)$ in \eqref{beta 2 transport}. Therefore, further noticing that $\beta_2(G_\eps)$ is non-negative and, recalling $\hat g_\eps^2\in L^\infty\left(dt;L^1\left(Mdxdv\right)\right)$ and the error estimate \eqref{error-alpha 6}, that the family $\int_K\beta_2(G_\eps)dxdv$ is equi-integrable, for any compact set $K\subset \mathbb{R}^3\times\mathbb{R}^3$, a direct application of Lemma \ref{HYPO-THM 3} yields that $\beta_2(G_\eps)=\left(\frac{\sqrt{G_\eps+\eps^a}-1}{\eps}\right)^2\gamma\left(\lambda\left(\sqrt{G_\eps+\eps^a}-1\right)\right)$ is equi-integrable in all variables $(t,x,v)$.

	Combining this result with \eqref{error-alpha 6} and \eqref{large values 2}, implies that the renormalized fluctuations $\hat g_\eps^2$ are equi-integrable in all variables $(t,x,v)$ as well.

	Finally, further combining the equi-integrability of $\hat g_\eps^2$ with the local strong compactness estimate \eqref{L1 xv compactness}, we deduce that the renormalized fluctuations $\hat g_\eps$ are relatively compact in $(x,v)$ in $L_{\mathrm{loc}}^2\left(dtdxdv\right)$, in the sense that, for any $\eta>0$ and every compact subset $K\subset [0,\infty)\times\mathbb{R}^3\times\mathbb{R}^3$, there exists $\gamma>0$ such that, if $h,l\in\mathbb{R}^3$ satisfy $|h|+|l|<\gamma$, then
	\begin{equation*}
		\sup_{\eps>0} \left\|\hat g_\eps(t,x+h,v+l) - \hat g_\eps(t,x,v)\right\|_{L^2\left(K,dtdxdv\right)}
		<\eta,
	\end{equation*}
	which concludes the proof of the lemma.
\end{proof}

Immediate consequences of the preceding strong compactness lemma are~:
\begin{itemize}
	\item the relative compactness in $(x,v)$ in $L^1\left(dtdxdv\right)$ of any subsequence of renormalized fluctuations $\hat g_\eps^2$ in the sense that, for any $\eta>0$ and every compact subset $K\subset [0,\infty)\times\mathbb{R}^3\times\mathbb{R}^3$, there exists $\gamma>0$ such that, if $h,l\in\mathbb{R}^3$ satisfy $|h|+|l|<\gamma$, then
	\begin{equation}\label{square strong comp}
		\sup_{\eps>0} \left\|\hat g_\eps^2(t,x+h,v+l) - \hat g_\eps^2(t,x,v)\right\|_{L^1\left(K,dtdxdv\right)}
		<\eta.
	\end{equation}
	
	\item the nonlinear weak compactness property, for any $p<2$,
	\begin{equation}\label{equiintegrability} % \label{equicontinuity-prop}
		\left(1+|v|^p\right)\hat g_\eps^2 \text{ is weakly relatively compact in }L^1_\mathrm{loc}\left(dtdx ; L^1\left(Mdv\right)\right),
	\end{equation}
	which follows from the Dunford-Pettis compactness criterion (see \cite{royden}) by deducing the equi-integrability of $\hat g_\eps^2$ from Lemma \ref{x-compactness1 0} and the tightness of $\left(1+|v|^p\right)\hat g_\eps^2$ from Lemma \ref{v2-int}.
	
	\item the strong spatial compactness of the moments $\int_{\mathbb{R}^3}\hat g_\eps \varphi (v)Mdv$ in $L^2_\mathrm{loc}(dtdx)$, for any $\varphi(v)\in L^2\left(\left(1+|v|^2\right)^{-1}Mdv\right)$, in particular
	\begin{equation}\label{xreg-hatmoments}
		\begin{aligned}
			\lim_{|h|\rightarrow 0}\sup_{\eps>0} \left\| \hat \rho_\eps(t,x+h) - \hat \rho_\eps(t,x)\right\|_{L^2_\mathrm{loc}(dtdx)} & = 0, \\
			\lim_{|h|\rightarrow 0}\sup_{\eps>0} \left\| \hat u_\eps(t,x+h) - \hat u_\eps(t,x)\right\|_{L^2_\mathrm{loc}(dtdx)} & = 0, \\
			\lim_{|h|\rightarrow 0}\sup_{\eps>0} \left\| \hat \theta_\eps(t,x+h) - \hat \theta_\eps(t,x)\right\|_{L^2_\mathrm{loc}(dtdx)} & = 0.
		\end{aligned}
	\end{equation}
\end{itemize}

The next lemma is also a direct consequence of the strong compactness properties from the preceding lemma and concerns a refinement of the relaxation estimate \eqref{relaxation-est} to $L^2_\mathrm{loc}(dtdxdv)$.

\begin{lem}\label{L2 relaxation}
	Let $\left(f_\eps, E_\eps, B_\eps\right)$ be the sequence of renormalized solutions to the scaled one species Vlasov-Maxwell-Boltzmann system \eqref{VMB1} considered in Theorem \ref{NS-WEAKCV}.
	
	Then, as $\eps\rightarrow 0$, any subsequence of renormalized fluctuations $\hat g_\eps$ satisfies the relaxation estimate
	\begin{equation*}
		\hat g_\eps - \Pi \hat g_\eps \rightarrow 0
		\qquad\text{in }L^2_\mathrm{loc}\left(dtdx ; L^2\left(Mdv\right)\right).
	\end{equation*}
\end{lem}

\begin{proof}
	On the one hand, we already know from Lemma \ref{relaxation-control} that
	\begin{equation}\label{L2 relaxation 1}
		\hat g_\eps - \Pi \hat g_\eps = O(\eps)_{L^1_\mathrm{loc}\left(dtdx ; L^2\left(Mdv\right)\right)}.
	\end{equation}
	
	On the other hand, the uniform integrability in all variables of $|\hat g_\eps|^2$ from Lemma \ref{x-compactness1 0} and the tightness in $v$ of $|\hat g_\eps|^2M$ implied by Lemma \ref{v2-int} shows that
	\begin{equation*}
		\int_{\mathbb{R}^3} \left|\Pi \hat g_\eps\right|^2 Mdv \leq C \int_{\mathbb{R}^3}\left|\hat g_\eps\right|^2 M dv
		\qquad
		\text{is uniformly integrable in $t$ and $x$.}
	\end{equation*}
	Therefore, we deduce that
	\begin{equation}\label{L2 relaxation 2}
		\left\|\hat g_\eps - \Pi\hat g_\eps\right\|_{L^2(Mdv)}^2
		\qquad
		\text{is uniformly integrable in $t$ and $x$.}
	\end{equation}

	Then, decomposing, for any large $\lambda>0$,
	\begin{equation*}
		\hat g_\eps - \Pi \hat g_\eps
		=
		\left(\hat g_\eps - \Pi \hat g_\eps\right) \mathds{1}_{\left\{\left\|\hat g_\eps - \Pi\hat g_\eps\right\|_{L^2(Mdv)} \leq \lambda \right\}}
		+
		\left(\hat g_\eps - \Pi \hat g_\eps\right) \mathds{1}_{\left\{\left\|\hat g_\eps - \Pi\hat g_\eps\right\|_{L^2(Mdv)} > \lambda \right\}},
	\end{equation*}
	we find that
	\begin{equation*}
		\begin{aligned}
			\left\|\hat g_\eps - \Pi \hat g_\eps\right\|_{L^2_\mathrm{loc}\left(dtdx ; L^2\left(Mdv\right)\right)}
			& \leq \sqrt\lambda
			\left\|\hat g_\eps - \Pi \hat g_\eps\right\|_{L^1_\mathrm{loc}\left(dtdx ; L^2\left(Mdv\right)\right)}^\frac{1}{2}
			\\
			& +
			\left\|\left(\hat g_\eps - \Pi \hat g_\eps\right)
			\mathds{1}_{\left\{\left\|\hat g_\eps - \Pi\hat g_\eps\right\|_{L^2(Mdv)} > \lambda \right\}}
			\right\|_{L^2_\mathrm{loc}\left(dtdx ; L^2\left(Mdv\right)\right)},
		\end{aligned}
	\end{equation*}
	whence, by virtue of \eqref{L2 relaxation 1},
	\begin{equation*}
		\begin{aligned}
		\limsup_{\eps\rightarrow 0} & \left\|\hat g_\eps - \Pi \hat g_\eps\right\|_{L^2_\mathrm{loc}\left(dtdx ; L^2\left(Mdv\right)\right)}
		\\
		& \leq
		\sup_{\eps>0}\left\|\left(\hat g_\eps - \Pi \hat g_\eps\right)
		\mathds{1}_{\left\{\left\|\hat g_\eps - \Pi\hat g_\eps\right\|_{L^2(Mdv)} > \lambda \right\}}
		\right\|_{L^2_\mathrm{loc}\left(dtdx ; L^2\left(Mdv\right)\right)}.
		\end{aligned}
	\end{equation*}
	Finally, thanks to the uniform integrability \eqref{L2 relaxation 2} and by the arbitrariness of $\lambda$, we infer
	\begin{equation*}
		\lim_{\eps\rightarrow 0} \left\|\hat g_\eps - \Pi \hat g_\eps\right\|_{L^2_\mathrm{loc}\left(dtdx ; L^2\left(Mdv\right)\right)} = 0,
	\end{equation*}
	which concludes the proof of the lemma.
\end{proof}

\subsection{Compactness of fluctuations for two species}\label{compactness two species}

We move on now to the study of strong compactness properties of the fluctuations considered in Theorems \ref{CV-OMHD} and \ref{CV-OMHDSTRONG} leading to the two-fluid incompressible Navier-Stokes-Fourier-Maxwell systems with Ohm's laws \eqref{TFINSFMO 2} and \eqref{TFINSFMSO 2}.

Unlike the estimates from Chapter \ref{weak bounds} infered from entropy and entropy dissipation bounds, here, we cannot deduce results for the two species case from results for the one species case. In fact, the regimes considered in Theorems \ref{CV-OMHD} and \ref{CV-OMHDSTRONG} are much more singular than the regime studied in Theorem \ref{NS-WEAKCV} and, as a result, the compactness properties asserted in Lemma \ref{x-compactness1 0} may not hold in the two species case.

It is to be emphasized that this lack of compactness is one the main drawbacks and difficulties preventing the improvement of Theorems \ref{CV-OMHD} and \ref{CV-OMHDSTRONG} to a weak compactness result similar to Theorem \ref{NS-WEAKCV}. Recall, however, that such an improvement is not to be expected so readily since the limiting systems \eqref{TFINSFMO 2} and \eqref{TFINSFMSO 2} are not stable under weak convergence in the energy space and, in particular, are not known to have global weak solutions (see corresponding discussion in Section \ref{stability existence 2}).

Thus, in the two species case, we only have the following weaker strong compactness result.

\begin{lem}\label{x-compactness2 0}
	Let $\left(f_\eps^\pm, E_\eps, B_\eps\right)$ be the sequence of renormalized solutions to the scaled two species Vlasov-Maxwell-Boltzmann system \eqref{VMB2} considered in Theorems \ref{CV-OMHD} and \ref{CV-OMHDSTRONG}.
	
	Then, as $\eps\rightarrow 0$, any subsequence of renormalized fluctuations $\hat g_\eps^\pm$ is locally relatively compact in $(x,v)$ in $L^p\left(dtdxdv\right)$, for any $1\leq p <2$, in the sense that, for any $\eta>0$ and every compact subset $K\subset [0,\infty)\times\mathbb{R}^3\times\mathbb{R}^3$, there exists $\gamma>0$ such that, if $h,l\in\mathbb{R}^3$ satisfy $|h|+|l|<\gamma$, then
	\begin{equation*}
		\sup_{\eps>0} \left\|\hat g_\eps^\pm(t,x+h,v+l) - \hat g_\eps^\pm(t,x,v)\right\|_{L^p\left(K,dtdxdv\right)}
		<\eta.
	\end{equation*}
	
	Furthermore, for any $\lambda>0$, the families $|\hat g_\eps^\pm|^2\mathds{1}_{\left\{\delta\lambda|\hat g_\eps^\pm|\leq 1 \right\}}$ are equi-integrable (in all variables $t$, $x$ and $v$).
\end{lem}

\begin{rem}
	We do not know whether the families $|\hat g_\eps^\pm|^2$ are equi-integrable (in all variables $t$, $x$ and $v$) or not.
\end{rem}

\begin{proof}
	The method of proof of this lemma is similar to the strategy used in the proof of Lemma \ref{x-compactness1 0}.

	\noindent {\bf An admissible renormalization $\beta_1(G_\eps^\pm)$.} As in the proof of Lemma \ref{x-compactness1 0}, we consider first the admissible square root renormalization
	\begin{equation*}
		\beta_1(z)=\frac{\sqrt{z+\eps^a}-1}{\eps},
	\end{equation*}
	for some given $1<a<4$.% Again, this renormalization is introduced to circumvent the fact that the natural renormalization $2\frac{\sqrt{z}-1}{\eps}$ corresponding to $\hat g_\eps^\pm$ is not admissible for Vlasov-Boltzmann equations, for it is singular at $z=0$, i.e.\ $\left(2\frac{\sqrt{z}-1}{\eps}\right)' \to \infty$ as $z\to 0$.

	Similarly to \eqref{error-alpha 5}, as a consequence of the entropy and entropy dissipation bounds, we have now that
	\begin{equation}\label{error-alpha 7}
		\begin{aligned}
		2\frac{\sqrt{G_\eps^\pm+\eps^a}-1}{\eps}-\hat g_\eps^\pm
		& = O\left(\eps^{\frac a2 -1}\right)_{L^\infty(dtdxdv)},
		\\
		2\frac{\sqrt{G_\eps^\pm+\eps^a}-1}{\eps}-\hat g_\eps^\pm
		& = O\left(\eps^{a-1}\right)_{L^\infty(dtdxdv)} +O\left(\eps^{\frac a2}\right)_{L^\infty\left(dt;L^2\left(M dxdv\right)\right)},
		\\
		2\frac{\sqrt{G_\eps^\pm+\eps^a}-1}{\eps}-\hat g_\eps^\pm
		& = O\left(\eps^{a-1}\right)_{L^\infty(dtdxdv)} +O\left(\eps^{\frac a2}\right)_{L^2_\mathrm{loc}\left(dtdx;L^2\left(\left(1+|v|^2\right)Mdv\right)\right)},
		\end{aligned}
	\end{equation}
	so that strong compactness properties of $\frac{\sqrt{G_\eps+\eps^a}-1}{\eps}$ in $L^p_{\mathrm{loc}}\left(dtdxdv\right)$, for any given $1\leq p\leq 2$, will entail similar properties on $\hat g_\eps^\pm$ in the same space and vice versa.

	\noindent {\bf Compactness of $\beta_1(G_\eps^\pm)$ in $v$.} In particular, since, in view of Lemma \ref{v-compactness}, the renormalized fluctuations $\hat g_\eps^\pm$ are locally relatively compact in $v$ in $L^2_\mathrm{loc}(dtdxdv)$, the same holds true for any subsequence of $\frac{\sqrt{G_\eps^\pm+\eps^a}-1}{\eps}$ in the sense that, for any $\eta>0$ and every compact subset $K\subset [0,\infty)\times\mathbb{R}^3\times\mathbb{R}^3$, there exists $\gamma>0$ such that, if $h\in\mathbb{R}^3$ satisfies $|h|<\gamma$, then
	\begin{equation}\label{beta 1 v comp 2}
		\sup_{\eps>0} \left\|\beta_1\left(G_\eps^\pm\right)(t,x,v+h) - \beta_1\left(G_\eps^\pm\right)(t,x,v)\right\|_{L^2\left(K,dtdxdv\right)}
		<\eta.
	\end{equation}

	\noindent {\bf Action of the transport operator on $\beta_1(G_\eps^\pm)$.} Thus, using $\beta_1(z)$ to renormalize the Vlasov-Boltzmann equations in \eqref{VMB2} and decomposing the collision integrands according to \eqref{integrands-decomposition}, we find that (see \eqref{renormalized2})
	\begin{equation*}
		\begin{aligned}
			( \eps \d_t + v \cdot \nabla_x \pm \delta \left(\eps E_\eps+ v\wedge B_\eps\right) & \cdot \nabla_v )
			{\sqrt{ G_\eps^\pm +\eps^a }-1 \over \eps}
			\mp \delta E_\eps \cdot v {G_\eps^\pm \over 2\sqrt{ G_\eps^\pm +\eps^a }} \\
			& =
			{\sqrt{G_\eps^\pm} \over 2\sqrt{ G_\eps^\pm+\eps^a }}
			\int_{\mathbb{R}^3\times\mathbb{S}^2} \sqrt{ G_{\eps *}^\pm } \hat q_\eps^\pm M_*dv_*d\sigma
			\\
			& +
			{\eps^2 \over 8\sqrt{G_\eps^\pm+\eps^a  }}
			\int_{\mathbb{R}^3\times\mathbb{S}^2} \left(\hat q_\eps^\pm\right)^2 M_*dv_*d\sigma
			\\
			& +
			{\delta\sqrt{G_\eps^\pm} \over 2\sqrt{ G_\eps^\pm+\eps^a }}
			\int_{\mathbb{R}^3\times\mathbb{S}^2} \sqrt{ G_{\eps *}^\mp } \hat q_\eps^{\pm,\mp} M_*dv_*d\sigma
			\\
			& +
			{\eps^2 \over 8\sqrt{G_\eps^\pm+\eps^a  }}
			\int_{\mathbb{R}^3\times\mathbb{S}^2} \left(\hat q_\eps^{\pm,\mp}\right)^2 M_*dv_*d\sigma.
		\end{aligned}
	\end{equation*}
	It follows, employing the uniform bounds $\hat g_\eps^\pm\in L^\infty\left(dt;L^2\left(Mdxdv\right)\right)$ and $\hat q_\eps^\pm,\hat q_\eps^{\pm,\mp}\in L^2\left(MM_*dtdxdvdv_*d\sigma\right)$ from Lemmas \ref{L2-lem} and \ref{L2-qlem}, respectively, that
	\begin{equation}\label{beta 1 transport 2}
		\begin{aligned}
			( \eps \d_t + & v \cdot \nabla_x )
			{\sqrt{ G_\eps^\pm +\eps^a }-1 \over \eps}
			\\
			& =
			\pm {\delta \sqrt{G_\eps^\pm} \over 2\sqrt{ G_\eps^\pm +\eps^a }} E_\eps \cdot v \left(1+\frac\eps 2\hat g_\eps\right)
			\mp \delta \nabla_v\cdot \left(\eps E_\eps+ v\wedge B_\eps\right)
			{\sqrt{ G_\eps^\pm +\eps^a }-1 \over \eps}
			\\
			& +
			{\sqrt{G_\eps^\pm} \over 2\sqrt{ G_\eps^\pm+\eps^a }}
			\int_{\mathbb{R}^3\times\mathbb{S}^2} \hat q_\eps^\pm M_*dv_*d\sigma
			+
			{\eps \sqrt{G_\eps^\pm} \over 4\sqrt{ G_\eps^\pm+\eps^a }}
			\int_{\mathbb{R}^3\times\mathbb{S}^2} \hat g_{\eps *}^\pm \hat q_\eps^\pm M_*dv_*d\sigma
			\\
			& +
			{\eps^2 \over 8\sqrt{G_\eps^\pm+\eps^a  }}
			\int_{\mathbb{R}^3\times\mathbb{S}^2} \left(\hat q_\eps^\pm\right)^2 M_*dv_*d\sigma
			\\
			& +
			{\delta\sqrt{G_\eps^\pm} \over 2\sqrt{ G_\eps^\pm+\eps^a }}
			\int_{\mathbb{R}^3\times\mathbb{S}^2} \hat q_\eps^{\pm,\mp} M_*dv_*d\sigma
			+
			{\delta\eps\sqrt{G_\eps^\pm} \over 4\sqrt{ G_\eps^\pm+\eps^a }}
			\int_{\mathbb{R}^3\times\mathbb{S}^2} \hat g_{\eps *}^\mp \hat q_\eps^{\pm,\mp} M_*dv_*d\sigma
			\\
			& +
			{\eps^2 \over 8\sqrt{G_\eps^\pm+\eps^a  }}
			\int_{\mathbb{R}^3\times\mathbb{S}^2} \left(\hat q_\eps^{\pm,\mp}\right)^2 M_*dv_*d\sigma
			\\
			& = O(1)_{L^1_{\mathrm{loc}}\left(dtdxdv\right)} + O(\delta)_{L^1_\mathrm{loc}\left(dtdx ; W^{-1,1}_{\mathrm{loc}}\left(dv\right)\right)}.
		\end{aligned}
	\end{equation}

	\noindent {\bf Compactness of $\beta_1(G_\eps^\pm)$ in $(x,v)$.} On the whole, we have established the compactness in velocity of $\beta_1(G_\eps^\pm)$ in \eqref{beta 1 v comp 2} and a bound on the transport operator acting on $\beta_1(G_\eps^\pm)$ in \eqref{beta 1 transport 2}. Therefore, a direct application of Lemma \ref{HYPO-THM 3} yields that $\beta_1(G_\eps^\pm)=\frac{\sqrt{G_\eps^\pm+\eps^a}-1}{\eps}$ is locally relatively compact in $(x,v)$ in $L^1_\mathrm{loc}\left(dtdxdv\right)$.

	Combining this result with \eqref{error-alpha 7}, implies that the renormalized fluctuations $\hat g_\eps^\pm$ are relatively compact in $(x,v)$ in $L_{\mathrm{loc}}^1\left(dtdxdv\right)$ as well. Moreover, since the families $\hat g_\eps^\pm$ are uniformly bounded in $L^\infty\left(dt;L^2\left(Mdxdv\right)\right)$, we easily deduce that the renormalized fluctuations $\hat g_\eps^\pm$ are relatively compact in $(x,v)$ in $L_{\mathrm{loc}}^p\left(dtdxdv\right)$, for any $1\leq p<2$, in the sense that, for any $\eta>0$ and every compact subset $K\subset [0,\infty)\times\mathbb{R}^3\times\mathbb{R}^3$, there exists $\gamma>0$ such that, if $h,l\in\mathbb{R}^3$ satisfy $|h|+|l|<\gamma$, then
	\begin{equation*}%\label{L1 xv compactness 2}
		\sup_{\eps>0} \left\|\hat g_\eps^\pm(t,x+h,v+l) - \hat g_\eps^\pm(t,x,v)\right\|_{L^p\left(K,dtdxdv\right)}
		<\eta,
	\end{equation*}
	which concludes the proof of the first part of the lemma.

	\noindent {\bf An admissible renormalization $\beta_2(G_\eps^\pm)$.} We proceed now to showing that, for any $\lambda>0$, the families $|\hat g_\eps^\pm|^2\mathds{1}_{\left\{\delta\lambda|\hat g_\eps^\pm|\leq 1 \right\}}$ are equi-integrable (in all variables $t$, $x$ and $v$), which will also be seen as a consequence of Lemma \ref{HYPO-THM 3}.

	To this end, we consider now the admissible renormalization
	\begin{equation*}
		\beta_2(z) = \left(\frac{\sqrt{z+\eps^a}-1}{\eps}\right)^2
		\gamma\left(\delta\lambda\left(\frac{\sqrt{z+\eps^a}-1}{\eps}\right)\right),
	\end{equation*}
	where $\gamma(z)\in C^1\left(\mathbb{R}\right)$ is a cutoff satisfying $\mathds{1}_{[-1,1]}\leq\gamma(z)\leq\mathds{1}_{[-2,2]}$, for some given $1<a <4$ (in fact, we will further restrict the range of $a$ so that necessarily $a=2$ below) and any small enough $\lambda>0$.% Again, this renormalization is introduced to circumvent the fact that the natural renormalization $\left(2\frac{\sqrt{z}-1}{\eps}\right)^2$ corresponding to $\hat g_\eps^{\pm 2}$ is not admissible for Vlasov-Boltzmann equations, for it is singular at $z=0$, i.e.\ $\left[\left(2\frac{\sqrt{z}-1}{\eps}\right)^2\right]' \to -\infty$ as $z\to 0$, and its growth at infinity is not admissible, i.e.\ $\left(2\frac{\sqrt{z}-1}{\eps}\right)^2\sim z$ as $z\to \infty$.

	Next, expressing $\beta_2(z)=\beta_1(z)^2\gamma\left(\delta\lambda \beta_1(z)\right)$ and using a slight variant of \eqref{gamma},
	% noticing that, for any $z_1,z_2\in\mathbb{R}$,
	% \begin{equation}\label{gamma}
	% 	\begin{aligned}
	% 		\left|z_1^2\gamma(\delta\lambda z_1)-z_2^2\gamma(\delta\lambda z_2)\right|
	% 		& =
	% 		\left|(z_1-z_2)z_1\gamma(\delta\lambda z_1)+z_2\left(z_1\gamma(\delta\lambda z_1)-z_2\gamma(\delta\lambda z_2)\right)\right|
	% 		\\
	% 		& \leq
	% 		C|z_1-z_2|\left(|z_1|+|z_2|\right),
	% 	\end{aligned}
	% \end{equation}
	it is readily seen that \eqref{error-alpha 7} implies that
	\begin{equation}\label{error-alpha 8}
		\begin{aligned}
		\left|4\beta_2\left(G_\eps^\pm\right) - \hat g_\eps^{\pm 2} \gamma\left(\frac{\delta\lambda}{2}\hat g_\eps^\pm\right)\right|
		& \leq C
		\left|2\beta_1\left(G_\eps^\pm\right)-\hat g_\eps^\pm\right|
		\left(\left|\beta_1\left(G_\eps^\pm\right)\right|+\left|\hat g_\eps^\pm\right|\right)
		\\
		& = O\left(\eps^{a-1}+\eps^\frac a2\right)_{L^\infty\left(dt;L^1_\mathrm{loc}\left(dx;L^1\left(Mdv\right)\right)\right)},
		\end{aligned}
	\end{equation}
	so that the equi-integrability of $\beta_2(G_\eps^\pm)=\left(\frac{\sqrt{G_\eps^\pm+\eps^a}-1}{\eps}\right)^2\gamma\left(\delta\lambda\left(\frac{\sqrt{G_\eps^\pm+\eps^a}-1}{\eps}\right)\right)$ will clearly entail the equi-integrability of $|\hat g_\eps^\pm|^2\mathds{1}_{\left\{\delta\lambda|\hat g_\eps^\pm|\leq 1 \right\}}$, for any $\lambda>0$.

	\noindent {\bf Compactness of $\beta_2(G_\eps^\pm)$ in $v$.} Furthermore, note that the velocity compactness stated in \eqref{beta 1 v comp 2} implies a corresponding property for $\beta_2(G_\eps^\pm)$. Indeed, using a slight modification of \eqref{gamma} again, we deduce, for any $h\in\mathbb{R}^3$, that
	\begin{equation*}
		\begin{aligned}
			& \left\|\beta_2\left(G_\eps^\pm\right)(t,x,v+h) - \beta_2\left(G_\eps^\pm\right)(t,x,v)\right\|_{L^1_\mathrm{loc}\left(dtdxdv\right)}
			\\
			& \leq C
			\left\|\beta_1\left(G_\eps^\pm\right)\right\|_{L^2_\mathrm{loc}\left(dtdxdv\right)}
			\left\|\beta_1\left(G_\eps^\pm\right)(t,x,v+h) - \beta_1\left(G_\eps^\pm\right)(t,x,v)\right\|_{L^2_\mathrm{loc}\left(dtdxdv\right)}.
		\end{aligned}
	\end{equation*}
	It then follows from \eqref{beta 1 v comp 2} that any subsequence of $\beta_2(G_\eps^\pm)$ is locally relatively compact in $v$ in $L^1\left(dtdxdv\right)$ in the sense that, for any $\eta>0$ and every compact subset $K\subset [0,\infty)\times\mathbb{R}^3\times\mathbb{R}^3$, there exists $\gamma>0$ such that, if $h\in\mathbb{R}^3$ satisfies $|h|<\gamma$, then
	\begin{equation}\label{beta 2 v comp 0 2}
		\sup_{\eps>0} \left\|\beta_2(G_\eps^\pm)(t,x,v+h) - \beta_2(G_\eps^\pm)(t,x,v)\right\|_{L^1\left(K,dtdxdv\right)}
		<\eta.
	\end{equation}

	However, in order to use the above velocity compactness of $\beta_2(G_\eps^\pm)$ in Lemma \ref{HYPO-THM 3}, we still need to show that $\beta_2(G_\eps^\pm)$ enjoys an improved integrability with respect to the velocity variable. To this end, we introduce the following decomposition
	\begin{equation*}
		\begin{aligned}
			\beta_2(G_\eps^\pm) & = \frac12 \Pi \hat g_\eps^\pm  {\sqrt{G_\eps^\pm+\eps^a }-1\over \eps}
			\gamma \left(\delta\lambda \left( \frac{\sqrt{G_\eps^\pm+\eps^a }-1}{\eps} \right)\right)
			\\
			& + \left({\sqrt{G_\eps^\pm+\eps^a }-1 \over \eps}
			- \frac 12 \hat g_\eps^\pm \right)
			{\sqrt{G_\eps^\pm+\eps^a }-1\over \eps}
			\gamma \left( \delta\lambda \left( \frac{\sqrt{G_\eps^\pm+\eps^a }-1}{\eps} \right)\right)
			\\
			& + \frac12( \hat g_\eps^\pm -\Pi \hat g_\eps^\pm)
			{\sqrt{G_\eps^\pm+\eps^a }-1\over \eps}
			\gamma \left( \delta\lambda \left( \frac{\sqrt{G_\eps^\pm+\eps^a }-1}{\eps} \right)\right).
		\end{aligned}
	\end{equation*}
	Since $\Pi \hat g_\eps^\pm$ belongs to $L^\infty\left(dt;L^2\left(dx;L^p\left(Mdv\right)\right)\right)$, for any $1\leq p<\infty$, we therefore get, for any $1\leq r <2$,
	\begin{equation*}
		\beta_2(G_\eps^\pm) \leq O(1)_{L^1_\mathrm{loc}\left(dtdx;L^r\left(Mdv\right)\right)}
		+\frac{C}{\lambda\delta}
		\left|{\sqrt{G_\eps^\pm+\eps^a }-1 \over \eps}
		- \frac 12 \hat g_\eps^\pm \right|
		+ \frac C{\lambda\delta}\left| \hat g_\eps^\pm -\Pi \hat g_\eps^\pm \right|.
	\end{equation*}
	Then, by \eqref{error-alpha 7} and by the relaxation estimate \eqref{relaxation-est}, we obtain, provided $2\leq a<4$,
	\begin{equation*}
		\begin{aligned}
			\beta_2(G_\eps^\pm) & = O(1)_{L^1_\mathrm{loc}\left(dtdx;L^r\left(Mdv\right)\right)}
			\\
			& + O\left(\frac {\eps^\frac{a}{2}}{\lambda\delta}\right)_{L^1_\mathrm{loc}\left(dtdx;L^2\left(Mdv\right)\right)}
			+ O\left(\frac \eps{\lambda\delta}\right)_{L^1_\mathrm{loc}\left(dtdx;L^2\left(Mdv\right)\right)}
			\\
			& =
			O\left(\frac 1\lambda\right)_{L^1_\mathrm{loc}\left(dtdx;L^r(Mdv)\right)},
		\end{aligned}
	\end{equation*}
	for any $1\leq r<2$.

	Finally, combining the preceding estimate with the compactness estimate \eqref{beta 2 v comp 0 2}, we deduce, for any $1\leq r<2$, that $\beta_2(G_\eps^\pm)$ is locally relatively compact in $v$ in $L^1\left(dtdx;L^r(dv)\right)$ in the sense that, for any $\eta>0$ and every compact subset $K\subset [0,\infty)\times\mathbb{R}^3\times\mathbb{R}^3$, there exists $\gamma>0$ such that, if $h\in\mathbb{R}^3$ satisfies $|h|<\gamma$, then
	\begin{equation}\label{beta 2 v comp 2}
		\sup_{\eps>0} \left\|\left(\beta_2(G_\eps^\pm)(t,x,v+h) - \beta_2(G_\eps^\pm)(t,x,v)\right)\mathds{1}_{K}(t,x,v)\right\|_{L^1\left(dtdx;L^r(dv)\right)}
		<\eta.
	\end{equation}

	\noindent {\bf Action of the transport operator on $\beta_2(G_\eps^\pm)$.} Thus, using $\beta_2(z)$ to renormalize the Vlasov-Boltzmann equations in \eqref{VMB2}, decomposing the collision integrands according to \eqref{integrands-decomposition} and writing for convenience
	\begin{equation*}
		\Gamma(z)=z\left(2\gamma(\lambda\delta z)+\lambda\delta z\gamma'(\lambda\delta z)\right),
	\end{equation*}
	so that $\beta_2'(z)=\beta_1'(z)\Gamma\left(\beta_1(z)\right)$, we have now that
	\begin{equation*}
		\begin{aligned}
			( \eps \d_t + v \cdot \nabla_x \pm \delta \left(\eps E_\eps+ v\wedge B_\eps\right) & \cdot \nabla_v )
			\beta_2(G_\eps^\pm)
			\mp \delta E_\eps \cdot v {G_\eps^\pm \Gamma\left(\beta_1(G_\eps^\pm)\right) \over 2\sqrt{ G_\eps^\pm +\eps^a }} \\
			& =
			{\sqrt{G_\eps^\pm} \Gamma\left(\beta_1(G_\eps^\pm)\right) \over 2\sqrt{ G_\eps^\pm+\eps^a }}
			\int_{\mathbb{R}^3\times\mathbb{S}^2} \sqrt{ G_{\eps *}^\pm } \hat q_\eps^\pm M_*dv_*d\sigma
			\\
			& +
			{\eps^2 \Gamma\left(\beta_1(G_\eps^\pm)\right) \over 8\sqrt{G_\eps^\pm+\eps^a  }}
			\int_{\mathbb{R}^3\times\mathbb{S}^2} \left(\hat q_\eps^\pm\right)^2 M_*dv_*d\sigma
			\\
			& +
			{\delta\sqrt{G_\eps^\pm} \Gamma\left(\beta_1(G_\eps^\pm)\right) \over 2\sqrt{ G_\eps^\pm+\eps^a }}
			\int_{\mathbb{R}^3\times\mathbb{S}^2} \sqrt{ G_{\eps *}^\mp } \hat q_\eps^{\pm,\mp} M_*dv_*d\sigma
			\\
			& +
			{\eps^2 \Gamma\left(\beta_1(G_\eps^\pm)\right) \over 8\sqrt{G_\eps^\pm+\eps^a  }}
			\int_{\mathbb{R}^3\times\mathbb{S}^2} \left(\hat q_\eps^{\pm,\mp}\right)^2 M_*dv_*d\sigma.
		\end{aligned}
	\end{equation*}
	It follows, employing the uniform bounds $\hat g_\eps^\pm \in L^\infty\left(dt;L^2\left(Mdxdv\right)\right)$ and $\hat q_\eps^{\pm},\hat q_\eps^{\pm,\mp}\in L^2\left(MM_*dtdxdvdv_*d\sigma\right)$ from Lemmas \ref{L2-lem} and \ref{L2-qlem}, respectively, and the direct estimates
	\begin{equation*}
		\begin{aligned}
			\left|\Gamma\left(\beta_1(G_\eps^\pm)\right)\right| & \leq C \left|\beta_1\left(G_\eps^\pm\right)\right|
			=O\left(1\right)_{L^2_\mathrm{loc}(dtdxdv)},
			\\
			\Gamma\left(\beta_1(G_\eps^\pm)\right) & =O\left(\frac{1}{\lambda\delta}\right)_{L^\infty(dtdxdv)},
			\\
			\left|\beta_2(G_\eps^\pm)\right| & \leq \frac{C}{\lambda\delta} \left|\beta_1\left(G_\eps^\pm\right)\right|
			=O\left(\frac 1 {\lambda\delta}\right)_{L^2_\mathrm{loc}(dtdxdv)},
		\end{aligned}
	\end{equation*}
	that, provided $1<a\leq 2$,
	\begin{equation}\label{beta 2 transport 2}
		\begin{aligned}
			( \eps \d_t + & v \cdot \nabla_x )
			\beta_2(G_\eps^\pm)
			\\
			& =
			\pm {\delta \sqrt{G_\eps^\pm} \Gamma\left(\beta_1(G_\eps^\pm)\right) \over 2\sqrt{ G_\eps^\pm +\eps^a }} E_\eps \cdot v \left(1+\frac\eps 2\hat g_\eps\right)
			\mp \delta \nabla_v\cdot \left[\left(\eps E_\eps+ v\wedge B_\eps\right)
			\beta_2(G_\eps^\pm)\right]
			\\
			& +
			{\sqrt{G_\eps^\pm}\Gamma\left(\beta_1(G_\eps^\pm)\right) \over 2\sqrt{ G_\eps^\pm+\eps^a }}
			\int_{\mathbb{R}^3\times\mathbb{S}^2} \hat q_\eps^\pm M_*dv_*d\sigma
			\\
			& +
			{\eps \sqrt{G_\eps^\pm}\Gamma\left(\beta_1(G_\eps^\pm)\right) \over 4\sqrt{ G_\eps^\pm+\eps^a }}
			\int_{\mathbb{R}^3\times\mathbb{S}^2} \hat g_{\eps *}^\pm \hat q_\eps^\pm M_*dv_*d\sigma
			\\
			& +
			{\eps^2\Gamma\left(\beta_1(G_\eps^\pm)\right) \over 8\sqrt{G_\eps^\pm+\eps^a  }}
			\int_{\mathbb{R}^3\times\mathbb{S}^2} \left(\hat q_\eps^\pm\right)^2 M_*dv_*d\sigma
			\\
			& +
			{\delta\sqrt{G_\eps^\pm}\Gamma\left(\beta_1(G_\eps^\pm)\right) \over 2\sqrt{ G_\eps^\pm+\eps^a }}
			\int_{\mathbb{R}^3\times\mathbb{S}^2} \hat q_\eps^{\pm,\mp} M_*dv_*d\sigma
			\\
			& +
			{\delta\eps\sqrt{G_\eps^\pm}\Gamma\left(\beta_1(G_\eps^\pm)\right) \over 4\sqrt{ G_\eps^\pm+\eps^a }}
			\int_{\mathbb{R}^3\times\mathbb{S}^2} \hat g_{\eps *}^\mp \hat q_\eps^{\pm,\mp} M_*dv_*d\sigma
			\\
			& +
			{\eps^2\Gamma\left(\beta_1(G_\eps^\pm)\right) \over 8\sqrt{G_\eps^\pm+\eps^a  }}
			\int_{\mathbb{R}^3\times\mathbb{S}^2} \left(\hat q_\eps^{\pm,\mp}\right)^2 M_*dv_*d\sigma
			\\
			& = % O\left(\frac 1\lambda\right)_{L^1_{\mathrm{loc}}\left(dtdxdv\right)} +
			O\left(\frac 1\lambda\right)_{L^1_\mathrm{loc}\left(dtdx ; W^{-1,1}_{\mathrm{loc}}\left(dv\right)\right)}.
		\end{aligned}
	\end{equation}
	Note here that the critical term in \eqref{beta 2 transport 2} preventing a better control on the transport acting on $\beta_2(G_\eps^\pm)$ and, thereby, on the concentrations of $\left|\hat g_\eps^{\pm}\right|^2$, is precisely
	\begin{equation*}
		\mp \delta \nabla_v\cdot \left(v\wedge B_\eps\right)
		\beta_2(G_\eps^\pm).
	\end{equation*}

	\noindent {\bf Equi-integrability of $\beta_2(G_\eps^\pm)$ in $(t,x,v)$.} On the whole, we have established the compactness in velocity of $\beta_2(G_\eps^\pm)$ in \eqref{beta 2 v comp 2} and a bound on the transport operator acting on $\beta_2(G_\eps^\pm)$ in \eqref{beta 2 transport 2}. Therefore, further noticing that $\beta_2(G_\eps^\pm)$ is non-negative and, recalling $\hat g_\eps^{\pm 2}\in L^\infty\left(dt;L^1\left(Mdxdv\right)\right)$ and the error estimate \eqref{error-alpha 8}, that the family $\int_K\beta_2(G_\eps^\pm)dxdv$ is equi-integrable, for any compact set $K\subset \mathbb{R}^3\times\mathbb{R}^3$, a direct application of Lemma \ref{HYPO-THM 3} yields that $\beta_2(G_\eps^\pm)=\left(\frac{\sqrt{G_\eps^\pm+\eps^a}-1}{\eps}\right)^2\gamma\left(\delta\lambda\left(\frac{\sqrt{G_\eps^\pm+\eps^a}-1}{\eps}\right)\right)$ is equi-integrable in all variables $(t,x,v)$.

	Finally, combining this result with \eqref{error-alpha 8}, implies that the renormalized fluctuations $|\hat g_\eps^\pm|^2\mathds{1}_{\left\{\delta\lambda|\hat g_\eps^\pm|\leq 1 \right\}}$, for any $\lambda>0$, are equi-integrable in all variables $(t,x,v)$ as well, which concludes the proof of the lemma.
\end{proof}

\begin{rem}
	Note that the estimates for one species \eqref{square strong comp}, \eqref{equiintegrability} and \eqref{xreg-hatmoments}, which were deduced directly from Lemma \ref{x-compactness1 0}, are no longer valid for two species in the settings of weak or strong interactions considered here.
	
	Notice, however, that the control \eqref{large values} on the very large values of fluctuations holds in all cases, for it is a mere consequence of the entropy bound only. Thus, on the one hand, estimate \eqref{large values} implies, for any $\lambda_1>0$ small enough, that
	\begin{equation*}
		\left|\hat g_\eps^\pm\right|^2 \mathds{1}_{\left\{\eps\lambda_1 \left|\hat g_\eps^\pm\right| >1\right\}}
		=O\left({1\over \left|\log \lambda_1\right|}\right)_{L^\infty\left(dt;L^1\left(Mdxdv\right)\right)},
	\end{equation*}
	while, on the other hand, we showed in Lemma \ref{x-compactness2 0} by controlling the action of the transport operator on the fluctuations that, for any $\lambda_2>0$,
	\begin{equation*}
		|\hat g_\eps^\pm|^2\mathds{1}_{\left\{\delta\lambda_2|\hat g_\eps^\pm|\leq 1 \right\}}
		\qquad
		\text{is equi-integrable (in all variables $t$, $x$ and $v$).}
	\end{equation*}
	
	In order to establish the equi-integrability of $|\hat g_\eps^\pm|^2$, there would therefore only remain to control the quantity
	\begin{equation*}
		|\hat g_\eps^\pm|^2\mathds{1}_{\left\{\frac{1}{\delta\lambda_2}< |\hat g_\eps^\pm|\leq \frac{1}{\eps\lambda_1} \right\}},
	\end{equation*}
	by showing that it is equi-integrable or uniformly small in $L^1_\mathrm{loc}(dtdxdv)$ as $\lambda_2\rightarrow 0$. But nothing seems to imply such a control. At least, we do not know how to prove it.
\end{rem}

The compactness results stated in Lemma \ref{x-compactness2 0} are valid in both regimes of weak and strong interspecies interactions. However, it is to be emphasized that the equi-integrability of $|\hat g_\eps^\pm|^2\mathds{1}_{\left\{\delta\lambda|\hat g_\eps^\pm|\leq 1 \right\}}$ contained therein definitely becomes a weaker compactness property as the parameter $\delta$ converges slower to zero. In other words, the more singular the regime, the weaker the compactness. In particular, in the extreme case of strong interspecies interactions considered in Theorem \ref{CV-OMHDSTRONG}, i.e.\ in the case $\delta=1$, the above equi-integrability statement is void, for $|\hat g_\eps^\pm|^2\mathds{1}_{\left\{\delta\lambda|\hat g_\eps^\pm|\leq 1 \right\}}$ is then uniformly bounded pointwise by $\lambda^{-2}$.

On the other hand, the more singular the regime, the stronger the bounds on the fluctuations provided by the entropy dissipation. This fact is epitomized by Lemmas \ref{bound hjw} and \ref{weak compactness h}, where it is apparent that the fluctuations $g_\eps^+-g_\eps^--n_\eps$ and $\hat g_\eps^+-\hat g_\eps^--\hat n_\eps$ vanish faster as the parameter $\delta$ converges slower to zero. From this perspective, the extreme case of strong interspecies interactions, i.e.\ the case $\delta=1$, enjoys better convergence properties than other less singular settings. In particular, employing the electrodynamic continuity equation from \eqref{VMB2}
\begin{equation*}
	\partial_t n_\eps + \frac 1\delta\nabla_x\cdot j_\eps = 0,
\end{equation*}
it is possible to deduce some strong compactness of $n_\eps$ in both $t$ and $x$ when $\delta=1$. This fails whenever $\delta=o(1)$. This crucial compactness will then allow us to consider the renormalized convergence of $h_\eps$ and $\hat h_\eps$, which is the content of the following lemma.

\begin{lem}\label{strong n}
	Let $\left(f_\eps^\pm, E_\eps, B_\eps\right)$ be the sequence of renormalized solutions to the scaled two species Vlasov-Maxwell-Boltzmann system \eqref{VMB2} considered in Theorem \ref{CV-OMHDSTRONG} for strong interspecies interactions, i.e.\ in the case $\delta = 1$.
	
	Then, as $\eps\rightarrow 0$, any subsequence of renormalized fluctuations $\hat g_\eps^\pm$ satisfies that $\hat n_\eps$ is relatively compact in $L^p_{\mathrm{loc}}\left(dtdx\right)$, for any $1\leq p<2$, and that $\hat g_\eps^+-\hat g_\eps^-$ is relatively compact in $L^p_{\mathrm{loc}}\left(dtdx;L^2(Mdv)\right)$, for any $1\leq p<2$. In particular, $\left\|\hat g_\eps^+-\hat g_\eps^-\right\|_{L^2(Mdv)} - \left|\hat n_\eps\right|\rightarrow 0$ in $L^p_{\mathrm{loc}}\left(dtdx\right)$, for any $1\leq p<2$.
	
	Furthermore, let $n\in L^\infty(dt;L^2(dx))$ be a limit point of $\hat n_\eps$ and, according to Lemma \ref{weak compactness h}, let $H\in L^1_{\mathrm{loc}}\left(dtdx;L^1\left(\left(1+|v|\right)Mdv\right)\right)$ and $\hat H\in L^2_{\mathrm{loc}}\left(dtdx;L^2(Mdv)\right)$ be limit points of $\frac{h_\eps}{1+\left\| \hat g_\eps^+- \hat g_\eps^- \right\|_{L^2(Mdv)}}$ and $\frac{\hat h_\eps}{1+\left\| \hat g_\eps^+- \hat g_\eps^- \right\|_{L^2(Mdv)}}$, respectively.
	
	Then, there exist $h\in L^1_{\mathrm{loc}}\left(dtdx;L^1\left(\left(1+|v|^2\right)Mdv\right)\right)$ and $\hat h\in L^1_{\mathrm{loc}}\left(dtdx;L^2\left(Mdv\right)\right)$ such that
	\begin{equation*}
		H=\frac{h}{1+|n|}\qquad\text{and}\qquad \hat H=\frac{\hat h}{1+|n|}.
	\end{equation*}
\end{lem}

\begin{proof}
	It is readily seen from Lemma \ref{x-compactness2 0}, that the family $\hat n_\eps$ is locally relatively compact in $x$ in $L_\mathrm{loc}^p(dtdx)$. Therefore, according to the decomposition \eqref{fluct-decomposition}, so is $n_\eps$ in $L_\mathrm{loc}^1(dtdx)$.
	
	Furthermore, taking the divergence of Amp\`ere's equation in \eqref{VMB2} yields the continuity equation
	\begin{equation*}
		\partial_t n_\eps + \DIV j_\eps = 0,
	\end{equation*}
	which, since $j_\eps$ is uniformly bounded in $L_\mathrm{loc}^1(dtdx)$ by Lemma \ref{bound hjw}, yields some temporal regularity on $n_\eps$ allowing us to establish, invoking a classical compactness result by Aubin and Lions \cite{aubin, lions6} (see also \cite{simon} for a sharp compactness criterion), that the family $n_\eps$ is strongly relatively compact in all variables in $L^1_\mathrm{loc}\left(dtdx\right)$. Employing the decomposition \eqref{fluct-decomposition}, again, we deduce that $\hat n_\eps$ is strongly relatively compact in all variables in $L^1_\mathrm{loc}\left(dtdx\right)$.
	
	Then, by virtue of the uniform bounds on $\hat g_\eps^\pm$ from Lemma \ref{L2-lem}, which clearly imply that $\hat n_\eps$ is uniformly bounded in $L^\infty(dt;L^2(dx))$, we conclude, by interpolation, that $\hat n_\eps$ is strongly relatively compact in all variables in $L^p_\mathrm{loc}\left(dtdx\right)$, for any $1\leq p<2$.
	
	Finally, we decompose
	\begin{equation*}
		\hat g_\eps^+-\hat g_\eps^- = \left(\hat g_\eps^+-\hat g_\eps^--\hat n_\eps\right) + \hat n_\eps
		= \eps \hat h_\eps + \hat n_\eps,
	\end{equation*}
	to deduce, using the bound \eqref{relaxation estimate} on $\hat h_\eps$ from Lemma \ref{relaxation2-control}, that $\hat g_\eps^+-\hat g_\eps^-$ is relatively strongly compact in $L^1_\mathrm{loc}\left(dtdx;L^2(Mdv)\right)$ and that
	\begin{equation*}
		\left\|\hat g_\eps^+-\hat g_\eps^-\right\|_{L^2(Mdv)} - \left|\hat n_\eps\right|\rightarrow 0
		\qquad\text{in } L^1_{\mathrm{loc}}\left(dtdx\right).
	\end{equation*}
	Then, again, by virtue of the uniform bounds on $\hat g_\eps^\pm$ from Lemma \ref{L2-lem}, which clearly imply that $\left\|\hat g_\eps^+-\hat g_\eps^-\right\|_{L^2(Mdv)}$ is uniformly bounded in $L^\infty(dt;L^2(dx))$, we conclude, by interpolation, that $\hat g_\eps^+-\hat g_\eps^-$ is strongly relatively compact in all variables in $L^p_\mathrm{loc}\left(dtdx;L^2(Mdv)\right)$ and that
	\begin{equation*}
		\left\|\hat g_\eps^+-\hat g_\eps^-\right\|_{L^2(Mdv)} - \left|\hat n_\eps\right|\rightarrow 0
		\qquad\text{in } L^p_{\mathrm{loc}}\left(dtdx\right),
	\end{equation*}
	for any $1\leq p<2$.

	There only remains to characterize the weak limits of $\frac{h_\eps}{1+\left\| \hat g_\eps^+- \hat g_\eps^- \right\|_{L^2(Mdv)}}$ and $\frac{\hat h_\eps}{1+\left\| \hat g_\eps^+- \hat g_\eps^- \right\|_{L^2(Mdv)}}$. To this end, we first assume, up to extraction of subsequences, that $\left\| \hat g_\eps^+- \hat g_\eps^- \right\|_{L^2(Mdv)}$ converges almost everywhere to $|n|$. Therefore, by the weak compactness results from Lemma \ref{weak compactness h} and the Product Limit Theorem (see \cite[Appendix B]{BGL2} and \cite[Appendix A]{SR}), we obtain that, for every $\lambda>0$,
	\begin{equation*}
		\begin{aligned}
			\frac{h_\eps}{1+\lambda \left\| \hat g_\eps^+- \hat g_\eps^- \right\|_{L^2(Mdv)}}
			& =
			\frac{1+\left\| \hat g_\eps^+- \hat g_\eps^- \right\|_{L^2(Mdv)}}{1+\lambda\left\| \hat g_\eps^+- \hat g_\eps^- \right\|_{L^2(Mdv)}}
			\frac{h_\eps}{1+\left\| \hat g_\eps^+- \hat g_\eps^- \right\|_{L^2(Mdv)}}
			\\
			&\rightharpoonup \frac{1+|n|}{1+\lambda|n|}H
			\qquad\text{in }\textit{w-}L^1_{\mathrm{loc}}\left(dtdx;\textit{w-}L^1\left(\left(1+|v|\right)Mdv\right)\right),
		\end{aligned}
	\end{equation*}
	and, similarly,
	\begin{equation*}
		\begin{aligned}
			\frac{\hat h_\eps}{1+\lambda \left\| \hat g_\eps^+- \hat g_\eps^- \right\|_{L^2(Mdv)}}
			& =
			\frac{1+\left\| \hat g_\eps^+- \hat g_\eps^- \right\|_{L^2(Mdv)}}{1+\lambda\left\| \hat g_\eps^+- \hat g_\eps^- \right\|_{L^2(Mdv)}}
			\frac{\hat h_\eps}{1+\left\| \hat g_\eps^+- \hat g_\eps^- \right\|_{L^2(Mdv)}}
			\\
			&\rightharpoonup \frac{1+|n|}{1+\lambda|n|}\hat H
			\qquad\text{in }\textit{w-}L^2_{\mathrm{loc}}\left(dtdx;\textit{w-}L^2(Mdv)\right).
		\end{aligned}
	\end{equation*}
	Therefore, for any $\varphi(v)\in C_c^\infty\left(\mathbb{R}^3\right)$ such that $\left|\varphi(v)\right|\leq (1+|v|^2)$, we find
	\begin{equation*}
		\left\|\frac{1+|n|}{1+\lambda|n|}H\varphi\right\|_{L^1_{\mathrm{loc}}\left(dtdx;L^1(Mdv)\right)}
		\leq \liminf_{\eps\rightarrow 0}
		\left\|h_\eps\right\|_{L^1_{\mathrm{loc}}\left(dtdx;L^1\left(\left(1+|v|^2\right)Mdv\right)\right)},
	\end{equation*}
	and
	\begin{equation*}
		\left\|\frac{1+|n|}{1+\lambda|n|}\hat H\right\|_{L^1_{\mathrm{loc}}\left(dtdx;L^2(Mdv)\right)}
		\leq \liminf_{\eps\rightarrow 0}
		\left\|\hat h_\eps\right\|_{L^1_{\mathrm{loc}}\left(dtdx;L^2\left(Mdv\right)\right)},
	\end{equation*}
	so that, in view of the bound on $h_\eps$ from Lemma \ref{bound hjw} and the bound on $\hat h_\eps$ from Lemma \ref{relaxation2-control} and by the arbitrariness of $\lambda$ and $\varphi$, we conclude
	\begin{equation*}
		\left(1+|n|\right)H \in L^1_{\mathrm{loc}}\left(dtdx;L^1\left(\left(1+|v|^2\right)Mdv\right)\right),
	\end{equation*}
	and
	\begin{equation*}
		\left(1+|n|\right)\hat H \in L^1_{\mathrm{loc}}\left(dtdx;L^2\left(Mdv\right)\right).
	\end{equation*}

	The justification of the lemma is complete.
\end{proof}

%% file: highconstraint0.tex
\chapter{Higher order and nonlinear constraint equations} \label{high constraints proof}

In Chapter \ref{constraints proof}, using weak compactness methods, we have derived lower order linear macroscopic constraint equations for one species and for two species in a weak interactions regime. For the one species case considered in Theorem \ref{NS-WEAKCV}, this is sufficient to obtain all constraint equations contained in the limiting system \eqref{NSFMP 2}. As for the two species case considered in Theorems \ref{CV-OMHD} and \ref{CV-OMHDSTRONG}, the corresponding limiting systems \eqref{TFINSFMSO 2} and \eqref{TFINSFMO 2}, respectively, contain higher order constraint equations (appearing as singular perturbations of the equations of motion) and nonlinear constraint equations, namely the (solenoidal) Ohm's law and the internal electric energy constraint, which cannot be deduced solely from the weak compactness bounds established in Chapter \ref{weak bounds}. We address now these singular limits employing the strong compactness bounds obtained in Chapter \ref{hypoellipticity}.

\section[Macroscopic constraint equations for two species\ldots]{Macroscopic constraint equations for two species, weak interactions}

As seen in Section \ref{macro constraint 2} (see \eqref{kinetic equation 2} in the proof of Proposition \ref{weak-comp2}), it is possible to derive limiting kinetic equations of the type
\begin{equation}\label{kinetic equation 4}
	v\cdot \nabla_x g^\pm = \int_{\mathbb{R}^3\times\mathbb{S}^2} q^\pm M_* dv_* d\sigma,
\end{equation}
from \eqref{VMB2} when $\delta=o(1)$.

What we intend to do next is to take advantage of the symmetries of the collision integrands $q^\pm$ and $q^{\pm,\mp}$ and of the strong compactness bounds from Chapter \ref{hypoellipticity} to go one order further and, thus, to derive a singular limit in the regime considered in Theorem \ref{CV-OMHD}. This singular limit is precisely the content of Proposition \ref{high weak-comp}, which will eventually yield the solenoidal Ohm's law and internal electric energy constraint from \eqref{TFINSFMSO 2} in Proposition \ref{solenoidalOhm} below. Of course, since we are considering renormalized fluctuations, we do not expect that the integrals in $v$ of the right-hand sides of the Vlasov-Boltzmann equations in \eqref{VMB2} against collision invariants are zero, but they should converge to zero as $\eps \to 0$ provided that we choose some appropriate renormalization which is sufficiently close to the identity. To estimate the ensuing conservation defects, we will also need to truncate large velocities.

Note that, even if conservation laws were known to hold for renormalized solutions of \eqref{VMB2}, we would have to introduce similar truncations of large tails and large velocities in order to control uniformly the flux and acceleration terms.

The main result in this section concerning the derivation of higher order nonlinear constraint equations in the regime considered in Theorem \ref{CV-OMHD} is contained in the following proposition.

\begin{prop}\label{high weak-comp}
	Let $\left(f_\eps^\pm, E_\eps, B_\eps\right)$ be the sequence of renormalized solutions to the scaled two species Vlasov-Maxwell-Boltzmann system \eqref{VMB2} considered in Theorem \ref{CV-OMHD} for weak interspecies interactions, i.e.\ $\delta=o(1)$ and $\frac\delta\eps$ unbounded. In accordance with Lemmas \ref{L1-lem}, \ref{L2-lem} and \ref{L2-qlem}, denote by
		\begin{equation*}
			\begin{gathered}
				g^\pm\in L^\infty\left(dt;L^2\left(Mdxdv\right)\right),
				\qquad
				q^{\pm,\mp} \in L^2\left(MM_*dtdxdvdv_*d\sigma\right),
				\\
				\text{and}\qquad
				E,B\in L^\infty\left(dt;L^2\left(dx\right)\right),
			\end{gathered}
		\end{equation*}
	any joint limit points of the families $\hat g_\eps^\pm$ and $\hat q_\eps^{\pm,\mp}$ defined by \eqref{hatg} and \eqref{hatq-def}, $E_\eps$ and $B_\eps$, respectively.
	
	Then, one has
	\begin{equation}\label{pre Ohm}
		\begin{aligned}
			\pm
			\int_{\mathbb{R}^3\times\mathbb{R}^3\times\mathbb{S}^2}
			q^{\pm,\mp}
			v MM_*dvdv_*d\sigma & = \nabla_x\bar p -\left(E + u\wedge B\right),
			% \\
			% P \left(E + u\wedge B
			% \pm
			% \int_{\mathbb{R}^3\times\mathbb{R}^3\times\mathbb{S}^2}
			% q^{\pm,\mp}
			% v MM_*dvdv_*d\sigma \right) & = 0,
			\\
			\int_{\mathbb{R}^3\times\mathbb{R}^3\times\mathbb{S}^2}
			q^{\pm,\mp}
			\left(\frac{|v|^2}{2}-\frac 52\right) MM_*dvdv_*d\sigma
			& = 0,
		\end{aligned}
	\end{equation}
	where $u$ is the bulk velocity associated with the limiting fluctuations $g^\pm$ and $\bar p\in L^1_\mathrm{loc}\left(dtdx\right)$ is a pressure.
\end{prop}

The proof of Proposition \ref{high weak-comp} is lengthy and contains several steps. Therefore, for the sake of clarity, it is deferred to Section \ref{proof prop ohm weak}, below.

As a direct consequence of the previous proposition, we derive in the next result the solenoidal Ohm's law and the internal electric energy constraint from \eqref{TFINSFMSO 2}.

\begin{prop}\label{solenoidalOhm}
	Let $\left(f_\eps^\pm, E_\eps, B_\eps\right)$ be the sequence of renormalized solutions to the scaled two species Vlasov-Maxwell-Boltzmann system \eqref{VMB2} considered in Theorem \ref{CV-OMHD} for weak interspecies interactions, i.e.\ $\delta=o(1)$ and $\frac\delta\eps$ unbounded. In accordance with Lemmas \ref{L1-lem}, \ref{L2-lem}, \ref{bound hjw} and \ref{weak compactness h} denote by
		\begin{equation*}
			\begin{gathered}
				g^\pm \in L^\infty\left(dt;L^2\left(Mdxdv\right)\right),
				\qquad
				h \in L^1_\mathrm{loc}\left(dtdx;L^1\left((1+|v|^2)Mdv\right)\right),
				\\
				\text{and}\qquad
				E,B\in L^\infty\left(dt;L^2\left(dx\right)\right),
			\end{gathered}
		\end{equation*}
	any joint limit points of the families $\hat g_\eps^\pm$ and $h_\eps$ defined by \eqref{hatg} and \eqref{def h}, $E_\eps$ and $B_\eps$, respectively.
	
	Then, one has
	\begin{equation*}
		j = \sigma\left(-\nabla_x \bar p + E + u\wedge B \right)
		\qquad\text{and}\qquad
		w = 0,
	\end{equation*}
	where $u$ is the bulk velocity associated with the limiting fluctuations $g^\pm$, $j$ and $w$ are, respectively, the electric current and the internal electric energy associated with the limiting fluctuation $h$, $\bar p\in L^1_\mathrm{loc}\left(dtdx\right)$ is a pressure and the electric conductivity $\sigma>0$ is defined by \eqref{sigma}.
\end{prop}

\begin{proof}
	By Proposition \ref{high weak-comp}, we have that
	\begin{equation*}
		\begin{aligned}
			\pm
			\int_{\mathbb{R}^3\times\mathbb{R}^3\times\mathbb{S}^2}
			q^{\pm,\mp}
			v MM_*dvdv_*d\sigma & = \nabla_x\bar p -\left(E + u\wedge B\right),
			% \\
			% P \left(E + u\wedge B
			% \pm
			% \int_{\mathbb{R}^3\times\mathbb{R}^3\times\mathbb{S}^2}
			% q^{\pm,\mp}
			% v MM_*dvdv_*d\sigma \right) & = 0,
			\\
			\int_{\mathbb{R}^3\times\mathbb{R}^3\times\mathbb{S}^2}
			q^{\pm,\mp}
			\left(\frac{|v|^2}{2}-\frac 52\right) MM_*dvdv_*d\sigma
			& = 0.
		\end{aligned}
	\end{equation*}
	Then, further incorporating identity \eqref{mixed q phi psi} from Proposition \ref{weak-comp3} into the above relations yields that
	\begin{equation*}
		\begin{aligned}
			E + u\wedge B
			-\frac 12 \int_{\mathbb{R}^3} j\cdot \mathfrak{L} \left(v\right) vMdv
			-\frac 12 \int_{\mathbb{R}^3} w\mathfrak{L} \left(\frac{|v|^2}{2}\right) vMdv
			& = \nabla_x \bar p,
			\\
			\int_{\mathbb{R}^3}
			j\cdot \mathfrak{L} \left(v\right)
			\left(\frac{|v|^2}{2}-\frac 52\right)Mdv
			+
			\int_{\mathbb{R}^3}
			w\mathfrak{L} \left(\frac{|v|^2}{2}\right)
			\left(\frac{|v|^2}{2}-\frac 52\right)
			Mdv
			& = 0.
		\end{aligned}
	\end{equation*}
	Finally, since, by symmetry considerations, $\int_{\mathbb{R}^3} \mathfrak{L} \left(v\right) \left(\frac{|v|^2}{2}-\frac 52\right)Mdv = \int_{\mathbb{R}^3} \mathfrak{L} \left(\frac{|v|^2}{2}\right) vMdv = 0$ and $\int_{\mathbb{R}^3} \mathfrak{L} \left(v_i\right) v_j Mdv = 0 $, if $i\neq j$, we compute that
	\begin{equation*}
		E + u\wedge B
		-\frac 1\sigma j
		= \nabla_x \bar p
		\qquad\text{and}\qquad
		\frac 1\lambda w
		= 0,
	\end{equation*}
	where $\sigma>0$ and $\lambda>0$ are defined in \eqref{sigma} and \eqref{lambda}, respectively, which concludes the proof of the proposition.
\end{proof}

\subsection{Proof of Proposition \ref{high weak-comp}}\label{proof prop ohm weak}

	Here, we analyze the equations \eqref{VMB-fluct}, which have to be renormalized, at a higher order. This becomes more complicated than the previous asymptotic analysis of \eqref{renormalized2} in the proof of Proposition \ref{weak-comp2}, because we do not have enough strong compactness to take limits in the nonlinear terms
	\begin{equation*}
		\pm \left(v\wedge B_\eps\right)\cdot \nabla_v
		{\sqrt{ G_\eps^\pm +\eps^a }-1 \over \eps},
	\end{equation*}
	therein. More precisely, we are not able to control the concentrations of $|\hat g_\eps^\pm|^2$ (see Lemma \ref{x-compactness2 0}).

	We will therefore consider a stronger renormalization of the equation for the fluctuations of density and exploit the equi-integrability from Lemma \ref{x-compactness2 0}, however weak it may be.

	\subsubsection{An admissible renormalization} We introduce the admissible renormalization $\Gamma_\lambda(z)$ defined by
	\begin{equation*}
		\Gamma_\lambda(z)-1=(z-1)\gamma\left(\lambda\delta\frac{z-1}{\eps}\right),
	\end{equation*}
	where $\frac\eps\delta\leq \lambda \leq 1$ is small and $\gamma \in C^1\left(\mathbb{R}\right)$ satisfies that
	\begin{equation*}
		\mathds{1}_{[-1,1]}(z)\leq \gamma(z) \leq \mathds{1}_{[-2,2]}(z), \qquad \text{for all }z\in\mathbb{R}.
	\end{equation*}
	Without distinguishing, for simplicity, the notation for cations and anions, we denote $\gamma_\eps^\lambda$ for $\gamma\left(\lambda\delta g_\eps^\pm\right)$ and $\hat\gamma_\eps^\lambda$ for $\Gamma_\lambda'\left(G_\eps^\pm\right)$. Thus, renormalizing the Vlasov-Boltzmann equation from \eqref{VMB2} with respect to $\Gamma_\lambda(z)$ yields
	\begin{equation*}
		\begin{aligned}
			\left( \frac\eps\delta \d_t + \frac 1\delta v \cdot \nabla_x \pm \left(\eps E_\eps+ v\wedge B_\eps\right)\cdot \nabla_v \right) &
			g_\eps^\pm\gamma_\eps^\lambda
			\mp E_\eps \cdot v G_\eps^\pm \hat\gamma_\eps^\lambda
			\\
			& = \frac 1{\delta\eps^2}\hat\gamma_\eps^\lambda\mathcal{Q}\left(G_\eps^\pm,G_\eps^\pm\right)
			+
			\frac \delta{\eps^2}\hat\gamma_\eps^\lambda\mathcal{Q}\left(G_\eps^\pm,G_\eps^\mp\right).
			% \\&
			% 		\frac 1\delta \hat\gamma_\eps^\lambda\sqrt{G_\eps^\pm}
			% 		\int_{\mathbb{R}^3\times\mathbb{S}^2} \sqrt{ G_{\eps *}^\pm } \hat q_\eps^\pm M_*dv_*d\sigma
			% 		+
			% 		\frac{\eps^2}{4\delta}\hat\gamma_\eps^\lambda
			% 		\int_{\mathbb{R}^3\times\mathbb{S}^2} \left(\hat q_\eps^\pm\right)^2 M_*dv_*d\sigma
			% 		\\
			% 		& +
			% 		\hat\gamma_\eps^\lambda\sqrt{G_\eps^\pm}
			% 		\int_{\mathbb{R}^3\times\mathbb{S}^2} \sqrt{ G_{\eps *}^\mp } \hat q_\eps^{\pm,\mp} M_*dv_*d\sigma
			% 		+
			% 		\frac{\eps^2}{4\delta}\hat\gamma_\eps^\lambda
			% 		\int_{\mathbb{R}^3\times\mathbb{S}^2} \left(\hat q_\eps^{\pm,\mp}\right)^2 M_*dv_*d\sigma.
		\end{aligned}
	\end{equation*}

	Notice here that there are two singular terms in the equations above, namely $\frac 1\delta v \cdot \nabla_x g_\eps^\pm \gamma_\eps^\lambda$ and the first term in the right-hand side (as shown below, the second term in the right-hand side is not singular). Therefore, in order to annihilate asymptotically these singular expressions, we integrate now the above equations against $\varphi(v)\chi\left(\frac{|v|^2}{K_\delta}\right)Mdv$, with $K_\delta=K\left|\log\delta\right|$, for some large $K>0$ to be fixed later on, for any collision invariant $\varphi(v)$ and some smooth compactly supported truncation $\chi\in C_c^\infty\left([0,\infty)\right)$ satisfying $\mathds{1}_{[0,1]}\leq\chi\leq\mathds{1}_{[0,2]}$, which leads to
	\begin{equation}\label{renormalized5}
		\begin{aligned}
			\frac\eps\delta\d_t \int_{\mathbb{R}^3} & g_\eps^\pm \gamma_\eps^\lambda \varphi \chi\left( {|v|^2\over K_\delta} \right) Mdv
			+ \frac1\delta \nabla_x  \cdot \int_{\mathbb{R}^3}
			g_\eps^\pm \gamma_\eps^\lambda\varphi \chi\left( {|v|^2\over K_\delta}\right)v Mdv\\
			& \mp \int_{\mathbb{R}^3}
			g_\eps^\pm \gamma_\eps^\lambda
			(\eps E_\eps+v\wedge B_\eps) \cdot \nabla_v \left( \varphi \chi\left( {|v|^2\over K_\delta}\right) M \right) dv \\
			& \mp
			E_\eps \cdot
			\int_{\mathbb{R}^3} \left(1+\eps g_\eps^\pm\right) \hat \gamma_\eps^\lambda \varphi\chi\left( {|v|^2\over K_\delta}\right) vM dv \\
			& =\frac 1{\delta\eps^2}
			\int_{\mathbb{R}^3}
			\hat\gamma_\eps^\lambda\mathcal{Q}\left(G_\eps^\pm,G_\eps^\pm\right)
			\varphi \chi\left( {|v|^2\over K_\delta} \right) Mdv
			\\
			& +\frac \delta{\eps^2}
			\int_{\mathbb{R}^3}
			\hat\gamma_\eps^\lambda\mathcal{Q}\left(G_\eps^\pm,G_\eps^\mp\right)
			\varphi \chi\left( {|v|^2\over K_\delta} \right) Mdv.
		\end{aligned}
	\end{equation}

	\subsubsection{Convergence of collision integrals} Let us focus on the right-hand side of \eqref{renormalized5} first. One has
	\begin{equation}\label{renormalized5 2}
		\begin{aligned}
			\frac \delta{\eps^2}
			\int_{\mathbb{R}^3}
			\hat\gamma_\eps^\lambda & \mathcal{Q}\left(G_\eps^\pm,G_\eps^\mp\right)
			\varphi \chi\left( {|v|^2\over K_\delta} \right) Mdv
			\\
			& =
			\int_{\mathbb{R}^3\times\mathbb{R}^3\times\mathbb{S}^2} \hat\gamma_\eps^\lambda\sqrt{G_\eps^\pm G_{\eps *}^\mp }
			\hat q_\eps^{\pm,\mp}
			\varphi\chi\left( {|v|^2\over K_\delta}\right)MM_*dvdv_*d\sigma
			\\
			& +
			\frac{\eps^2}{4\delta}
			\int_{\mathbb{R}^3\times\mathbb{R}^3\times\mathbb{S}^2} \hat\gamma_\eps^\lambda \left(\hat q_\eps^{\pm,\mp}\right)^2
			\varphi\chi\left( {|v|^2\over K_\delta}\right)MM_*dvdv_*d\sigma.
		\end{aligned}
	\end{equation}
	Since $\frac\eps\delta$ vanishes, $\varphi\chi\left( {|v|^2\over K_\delta}\right)$ is bounded pointwise by a constant multiple of $\left|\log \delta\right|$ and the collision integrands $\hat q_\eps^{\pm,\mp}$ are uniformly bounded in $L^2\left(MM_*dtdxdvdv_*d\sigma\right)$, according to Lemma \ref{L2-qlem}, we find that the second term from the right-hand side of \eqref{renormalized5 2} vanishes in $L^1\left(dtdx\right)$. Further utilizing that, thanks to Lemma \ref{L2-lem},
	\begin{equation*}%\label{renormalized5 1}
		\begin{aligned}
			\sqrt{G_{\eps}^\pm} & = 1 + O(\eps)_{L^2_\mathrm{loc}\left(dt ; L^2\left(M dxdv\right)\right)}, \\
			\sqrt{G_{\eps *}^\pm} & = 1 + O(\eps)_{L^2_\mathrm{loc}\left(dt ; L^2\left(M_* dxdv_*\right)\right)},
		\end{aligned}
	\end{equation*}
	and that $\hat\gamma_\eps^\lambda\sqrt{G_\eps^\pm}$ is uniformly bounded pointwise, it is readily seen that the weak limit of the first term of the right-hand side of \eqref{renormalized5 2} coincides with the weak limit of
	\begin{equation*}
		\int_{\mathbb{R}^3\times\mathbb{R}^3\times\mathbb{S}^2} \hat\gamma_\eps^\lambda
		\hat q_\eps^{\pm,\mp}
		\varphi\chi\left( {|v|^2\over K_\delta}\right)MM_*dvdv_*d\sigma,
	\end{equation*}
	which, since $\hat\gamma_\eps^\lambda\varphi\chi\left( {|v|^2\over K_\eps}\right)$ is dominated by $|\varphi|$ and converges almost everywhere to $\varphi$, is easily shown to converge weakly in $L^2(dtdx)$ towards
	\begin{equation*}
		\int_{\mathbb{R}^3\times\mathbb{R}^3\times\mathbb{S}^2}
		q^{\pm,\mp}
		\varphi MM_*dvdv_*d\sigma.
	\end{equation*}
	Thus, so far, we have established that the second term in the right-hand side of \eqref{renormalized5} satisfies
	\begin{equation}\label{renormalized5 second term}
		\begin{gathered}
			\frac \delta{\eps^2}
			\int_{\mathbb{R}^3}
			\hat\gamma_\eps^\lambda\mathcal{Q}\left(G_\eps^\pm,G_\eps^\mp\right)
			\varphi \chi\left( {|v|^2\over K_\delta} \right) Mdv
			\rightharpoonup
			\int_{\mathbb{R}^3\times\mathbb{R}^3\times\mathbb{S}^2}
			q^{\pm,\mp}
			\varphi MM_*dvdv_*d\sigma
			\\
			\text{in }
			L^1_\mathrm{loc}\left(dtdx\right).
		\end{gathered}
	\end{equation}

	The first term of the right-hand side of \eqref{renormalized5} is more singular and, therefore, harder to control. One has, in this case, taking advantage of collisional symmetries, that
	\begin{equation}\label{first term}
		\begin{aligned}
			\frac 1{\delta\eps^2}
			\int_{\mathbb{R}^3}
			\hat\gamma_\eps^\lambda & \mathcal{Q}\left(G_\eps^\pm,G_\eps^\pm\right)
			\varphi \chi\left( {|v|^2\over K_\delta} \right) Mdv
			\\
			& =
			\frac{\eps^2}{4\delta}
			\int_{\mathbb{R}^3\times\mathbb{R}^3\times\mathbb{S}^2} \hat\gamma_\eps^\lambda
			\varphi\chi\left( {|v|^2\over K_\delta}\right)
			\left(\hat q_\eps^\pm\right)^2 MM_*dvdv_*d\sigma
			\\
			& -
			\frac 1\delta
			\int_{\mathbb{R}^3\times\mathbb{R}^3\times\mathbb{S}^2}
			\hat\gamma_\eps^\lambda
			\varphi\left(1-\chi\left( {|v|^2\over K_\delta}\right)\right)
			\hat q_\eps^\pm \sqrt{G_\eps^\pm G_{\eps *}^\pm } MM_*dvdv_*d\sigma
			\\
			& +
			\frac 1\delta
			\int_{\mathbb{R}^3\times\mathbb{R}^3\times\mathbb{S}^2}
			\hat\gamma_\eps^\lambda\left(1-\hat\gamma_{\eps*}^\lambda\right)
			\varphi
			\hat q_\eps^\pm \sqrt{G_\eps^\pm G_{\eps *}^\pm } MM_*dvdv_*d\sigma
			\\
			& +
			\frac 1\delta
			\int_{\mathbb{R}^3\times\mathbb{R}^3\times\mathbb{S}^2}
			\hat\gamma_\eps^\lambda \hat\gamma_{\eps*}^\lambda
			\left(1-\hat\gamma_{\eps}^{\lambda\prime}\hat\gamma_{\eps*}^{\lambda\prime}\right)
			\varphi
			\hat q_\eps^\pm \sqrt{G_\eps^\pm G_{\eps *}^\pm } MM_*dvdv_*d\sigma
			\\
			& -\frac {\eps^2}{4\delta}
			\int_{\mathbb{R}^3\times\mathbb{R}^3\times\mathbb{S}^2}
			\hat\gamma_\eps^\lambda \hat\gamma_{\eps*}^\lambda
			\hat\gamma_{\eps}^{\lambda\prime}\hat\gamma_{\eps*}^{\lambda\prime}
			\varphi
			\left(\hat q_\eps^\pm\right)^2 MM_*dvdv_*d\sigma
			\\
			& \eqdefa D^1_\eps(\varphi)+D^2_\eps(\varphi)+D^3_\eps(\varphi)+D^4_\eps(\varphi)+D^5_\eps(\varphi),
		\end{aligned}
	\end{equation}
	where we have used that $\varphi$ is a collision invariant to symmetrize the last term.
	
	Now, we control each term $D^i_\eps(\varphi)$, $i=1,\ldots,5$, separately.

	\noindent$\bullet$ The vanishing of the first term $D^1_\eps(\varphi)$, for any function $\varphi(v)$ growing at most quadratically at infinity, easily follows, using Lemma \ref{L2-qlem}, from the estimate
	\begin{equation}\label{first term D1}
		\begin{aligned}
			\left\|D_\eps^1(\varphi)\right\|_{L^1(dtdx)} & \leq
			\frac{\eps^2}{4\delta}\left\| \hat q_\eps^\pm \right\|^2_{L^2\left(MM_* dtdxdvdv_* d\sigma\right)}
			\left\| \hat \gamma_\eps^\lambda \right\|_{L^\infty}
			\left\| \chi \left( {|v|^2\over K_\delta} \right)\varphi \right\|_{L^\infty} \\
			& \leq C \frac{\eps^2}{\delta} K_\delta = CK \frac{\eps^2}{\delta}|\log\delta|.
		\end{aligned}
	\end{equation}

	\noindent$\bullet$ The second term $D^2_\eps(\varphi)$ is controlled by the following estimate on the tails of Gaussian distributions~: for any $p\in\mathbb{R}$, as $R\rightarrow \infty$,
	\begin{equation}\label{gaussian-decay 0}
		\int_{\left\{|v|^2 >R\right\}} |v|^p M(v) dv
		\sim \sqrt{\frac2\pi} R^\frac{p+1}{2} e^{-\frac{R}{2}},
	\end{equation}
	in the sense that the quotient of both sides converges to $1$ as $R\rightarrow\infty$, which is easily established by applying Bernoulli-l'Hospital's rule.

	We have indeed
	\begin{equation*}
		\begin{aligned}
			& \left|D_\eps^2(\varphi)\right|
			\\
			& \leq
			\frac 1\delta  \left\| \hat q_\eps^\pm \right\|_{L^2\left(MM_*dvdv_*d\sigma\right)}
			\left\| \hat \gamma_\eps^\lambda \varphi \mathds{1}_{\left\{|v|^2 \geq K_\delta\right\}} \sqrt{G_\eps^\pm G_{\eps *}^\pm}\right\|_{L^2\left(MM_*dvdv_*d\sigma\right)}
			\\
			& \leq
			\frac C\delta \left\| \hat q_\eps^\pm \right\|_{L^2\left(MM_*dvdv_*d\sigma\right)}
			\left\| \hat \gamma_\eps^\lambda \sqrt{G_\eps^\pm} \right\|_{L^\infty}
			\left\| \sqrt{G_\eps^\pm} \right\|_{ L^2\left(Mdv\right)}
			\left\| \mathds{1}_{\left\{|v|^2 \geq K_\delta\right\}} \varphi \right\|_{L^2\left(Mdv\right)}.
		\end{aligned}
	\end{equation*}
	Thus, using the bound from Lemma \ref{L2-qlem}, the pointwise boundedness of $\Gamma_\lambda'(z)\sqrt{z}$ and the Gaussian decay estimate \eqref{gaussian-decay 0}, we get, for all $\varphi(v)$ growing at most quadratically at infinity,
	\begin{equation}\label{first term D2}
		D_\eps^2(\varphi)
		=O\left(\delta^{\frac{K}{4}-1}\left|\log\delta\right|^\frac{5}{4}\right)_{L^1_\mathrm{loc}\left(dtdx\right)},
	\end{equation}
	which tends to zero as soon as $K>4$.

	\noindent$\bullet$ The last term $D^5_\eps (\varphi)$ is mastered using the same tools. For high energies, i.e.\ when $|v|^2 \geq K|\log \delta|$, we obtain
	\begin{equation*}
		\begin{aligned}
			& D_\eps^{5>}(\varphi) \\
			& \eqdefa
			\frac{\eps^2}{4\delta}\int_{\mathbb{R}^3\times\mathbb{R}^3\times\mathbb{S}^2}
			\hat\gamma_\eps^\lambda \hat\gamma_{\eps*}^\lambda
			\hat\gamma_{\eps}^{\lambda\prime}\hat\gamma_{\eps*}^{\lambda\prime}
			\varphi
			\mathds{1}_{\left\{|v|^2\geq K_\delta\right\}}
			\left(\hat q_\eps^\pm\right)^2 MM_* dvdv_*d\sigma
			\\
			& \leq \frac 1\delta
			\left\| \hat \gamma_\eps^\lambda \sqrt{G_\eps^\pm} \right\|_{L^\infty}^2
			\left\| \hat \gamma_\eps^\lambda \right\|_{L^\infty}^2
			\left\| \hat q_\eps^\pm \right\|_{L^2\left(MM_*dvdv_*d\sigma\right)}
			\left\| \varphi \mathds{1}_{\left\{|v|^2 \geq K_\delta\right\}} \right\|_{L^2\left(MM_*dvdv_*d\sigma\right)},
		\end{aligned}
	\end{equation*}
	so that, using the estimate \eqref{gaussian-decay 0} on the tails of Gaussian distributions and the bound on $\hat q_\eps^\pm$ from Lemma \ref{L2-qlem},
	\begin{equation}\label{first term D5 1}
		D_\eps^{5>}(\varphi)=O\left(\delta^{\frac{K}{4}-1}\left|\log\delta\right|^\frac{5}{4}\right)_{L^2\left(dtdx\right)},
	\end{equation}
	which tends to zero as soon as $K>4$.
	
	For moderate energies, i.e.\ when $|v|^2 <  K|\log \delta|$, we easily find
	\begin{equation*}
		\begin{aligned}
			D_\eps^{5<}(\varphi)
			& \eqdefa
			\frac{\eps^2}{4\delta}\int_{\mathbb{R}^3\times\mathbb{R}^3\times\mathbb{S}^2}
			\hat\gamma_\eps^\lambda \hat\gamma_{\eps*}^\lambda
			\hat\gamma_{\eps}^{\lambda\prime}\hat\gamma_{\eps*}^{\lambda\prime}
			\varphi
			\mathds{1}_{\left\{|v|^2 < K_\delta\right\}}
			\left(\hat q_\eps^\pm\right)^2 MM_* dvdv_*d\sigma
			\\
			& \leq CK\frac{\eps^2}{\delta}\left|\log\delta\right|\left\|\hat q_\eps^\pm\right\|_{L^2\left(MM_*dvdv_*d\sigma\right)}^2,
		\end{aligned}
	\end{equation*}
	so that the entropy dissipation bound from Lemma \ref{L2-qlem} provides
	\begin{equation}\label{first term D5 2}
		D_\eps^{5<}(\varphi) = O\left(\frac{\eps^2}{\delta}\left|\log\delta\right|\right)_{L^1\left(dtdx\right)}.
	\end{equation}

	\noindent$\bullet$ The handling of $D_\eps^3(\varphi)$ accounts for the introduction of the small parameter $\lambda$. First, one has, by the Cauchy-Schwarz inequality,
	\begin{equation*}
		\begin{aligned}
			& \left|D_\eps^3(\varphi)\right|
			\\
			& \leq
			\frac 1\delta  \left\| \hat q_\eps^\pm \right\|_{L^2\left(MM_*dvdv_*d\sigma\right)}
			\left\| \hat \gamma_\eps^\lambda \left(1-\hat\gamma_{\eps *}^\lambda\right) \varphi
			\sqrt{G_\eps^\pm G_{\eps *}^\pm}\right\|_{L^2\left(MM_*dvdv_*d\sigma\right)}
			\\
			& \leq C
			\left\| \hat q_\eps^\pm \right\|_{L^2\left(MM_*dvdv_*d\sigma\right)}
			\left\| \hat \gamma_\eps^\lambda \sqrt{G_\eps^\pm} \right\|_{L^\infty}
			\left\| \frac{1}{\delta}\left(1-\hat\gamma_\eps^\lambda\right)\sqrt{G_\eps^\pm} \right\|_{ L^2\left(Mdv\right)}
			\left\| \varphi \right\|_{L^2\left(Mdv\right)}
			\\
			& \leq C
			\left\| \hat q_\eps^\pm \right\|_{L^2\left(MM_*dvdv_*d\sigma\right)}
			\left\| \hat \gamma_\eps^\lambda \sqrt{G_\eps^\pm} \right\|_{L^\infty}
			\left\| \frac{1}{\delta}\left(1-\hat\gamma_\eps^\lambda\right)\right\|_{ L^2\left(Mdv\right)}
			\left\| \varphi \right\|_{L^2\left(Mdv\right)}
			\\
			& +C \frac\eps\delta
			\left\| \hat q_\eps^\pm \right\|_{L^2\left(MM_*dvdv_*d\sigma\right)}
			\left\| \hat \gamma_\eps^\lambda \sqrt{G_\eps^\pm} \right\|_{L^\infty}
			\left\| \hat g_\eps^\pm \right\|_{ L^2\left(Mdv\right)}
			\left\| \varphi \right\|_{L^2\left(Mdv\right)}.
		\end{aligned}
	\end{equation*}
	Moreover, in view of the hypotheses on $\gamma(z)$, the support of $\Gamma_\lambda'(z)-1=\gamma\left(\lambda\delta\frac{z-1}{\eps}\right)-1+\lambda\delta \frac{z-1}{\eps}\gamma'\left(\lambda\delta\frac{z-1}{\eps}\right)$ is clearly restricted to $\lambda\delta\frac{|z-1|}{\eps} \in [1,\infty)$, so that, employing the decomposition \eqref{fluct-decomposition},
	\begin{equation}\label{gamma bound}
		\begin{aligned}
			\frac1\delta \left|1-\hat \gamma_\eps^\lambda\right|
			& =
			\frac1\delta \left|1-\hat \gamma_\eps^\lambda\right| \mathds{1}_{\left\{\left|\eps\hat g_\eps^\pm\right|\leq 1\right\}}
			+
			\frac1\delta \left|1-\hat \gamma_\eps^\lambda\right| \mathds{1}_{\left\{\left|\eps\hat g_\eps^\pm\right|> 1\right\}}
			\\
			& \leq
			\lambda \left|1-\hat \gamma_\eps^\lambda\right| \left|g_\eps^\pm\right| \mathds{1}_{\left\{\left|\eps\hat g_\eps^\pm\right|\leq 1\right\}}
			+
			\frac\eps\delta \left|1-\hat \gamma_\eps^\lambda\right| \left|\hat g_\eps^\pm\right| \mathds{1}_{\left\{\left|\eps\hat g_\eps^\pm\right|> 1\right\}}
			\\
			& \leq
			C\lambda \left|\hat g_\eps^\pm+\frac\eps 4\hat g_\eps^{\pm 2}\right| \mathds{1}_{\left\{\left|\eps\hat g_\eps^\pm\right|\leq 1\right\}}
			+
			C\frac\eps\delta \left|\hat g_\eps^\pm\right|
			\\
			& \leq
			C\left(\lambda+\frac\eps\delta\right) \left|\hat g_\eps^\pm\right|.
		\end{aligned}
	\end{equation}
	Therefore, thanks to the bound on $\hat q_\eps$ from Lemma \ref{L2-qlem}, we infer
	\begin{equation*}
		\left\|D_\eps^3(\varphi)\right\|_{L^1_\mathrm{loc}\left(dtdx\right)}
		\leq C
		\left(\lambda + \frac\eps\delta\right)\left\| \hat g_\eps ^\pm\right\|_{ L^2_\mathrm{loc}\left(dtdx;L^2\left(Mdv\right)\right)}.
	\end{equation*}
	Thus, we conclude that
	\begin{equation}\label{first term D3}
		\left\|D_\eps^3(\varphi)\right\|_{L^1_\mathrm{loc}\left(dtdx\right)}\leq C\lambda.
	\end{equation}

	\noindent$\bullet$ A similar argument provides the convergence of the remaining term $D_\eps^4(\varphi)$. Thus, one has by the Cauchy-Schwarz inequality, for any $2<p<\infty$,
	\begin{equation*}% \label{eq-use2}
		\begin{aligned}
			& \left|D_\eps^4(\varphi)\right|
			\\
			& \leq
			\frac 1\delta  \left\| \hat q_\eps^\pm \right\|_{L^2\left(MM_*dvdv_*d\sigma\right)}
			\left\| \hat \gamma_\eps^\lambda \hat\gamma_{\eps *}^\lambda\left(1-\hat\gamma_{\eps}^{\lambda\prime}\hat\gamma_{\eps*}^{\lambda\prime}\right) \varphi
			\sqrt{G^\pm_\eps G^\pm_{\eps *}}\right\|_{L^2\left(MM_*dvdv_*d\sigma\right)}
			\\
			& \leq C
			\left\| \hat q_\eps^\pm \right\|_{L^2\left(MM_*dvdv_*d\sigma\right)}
			\left\| \hat \gamma_\eps^\lambda \sqrt{G_\eps^\pm} \right\|_{L^\infty}^2
			\left\| \frac 1\delta\left(1-\hat\gamma_{\eps}^{\lambda\prime}\hat\gamma_{\eps*}^{\lambda\prime}\right) \varphi \right\|_{L^2\left(MM_*dvdv_*d\sigma\right)}
			\\
			& \leq C_p
			\left\| \hat q_\eps^\pm \right\|_{L^2\left(MM_*dvdv_*d\sigma\right)}
			\left\| \frac 1\delta \left(1-\hat\gamma^\lambda_{\eps}\right) \right\|_{L^p\left(Mdv\right)}.
		\end{aligned}
	\end{equation*}
	Therefore, thanks to the bound on $\hat q_\eps^\pm$ from Lemma \ref{L2-qlem}, we infer, for any $2<p<\infty$,
	\begin{equation*}
		\left\|D_\eps^4(\varphi)\right\|_{L^1_\mathrm{loc}\left(dtdx\right)}
		\leq C
		\left\| \frac 1\delta\left(1-\hat\gamma_\eps^\lambda\right) \right\|_{ L^2_\mathrm{loc}\left(dtdx;L^p\left(Mdv\right)\right)}.
	\end{equation*}
	
	Next, using estimate \eqref{relaxation-est} from Lemma \ref{relaxation-control} and the bound \eqref{gamma bound}, we find that
	\begin{equation*}
		\begin{aligned}
			\frac1{\delta^2} & \left|1-\hat \gamma_\eps^\lambda\right|^2
			\\
			& \leq C
			\left(\lambda +\frac\eps\delta\right) \left|\hat g_\eps^\pm\right| \frac1\delta \left|1-\hat \gamma_\eps^\lambda\right|
			\\
			& \leq C
			\left(\lambda +\frac\eps\delta\right) \left(\left|\Pi\hat g_\eps^\pm\right| +\left|\hat g_\eps^\pm-\Pi\hat g_\eps^\pm\right|\right) \frac1\delta \left|1-\hat \gamma_\eps^\lambda\right|
			\\
			& \leq C
			\left(\lambda +\frac\eps\delta\right)^2 \left|\Pi\hat g_\eps^\pm\right| \left|\hat g_\eps^\pm\right|
			+
			C
			\left(\lambda +\frac\eps\delta\right) \frac 1\delta \left|\hat g_\eps^\pm-\Pi\hat g_\eps^\pm\right|
			\\
			&\leq
			O\left(\left(\lambda +\frac\eps\delta\right)^2\right)_{L^1_\mathrm{loc}\left(dtdx;L^r\left(Mdv\right)\right)}
			+ O\left(\left(\lambda +\frac\eps\delta\right)\frac\eps\delta\right)_{L^1_\mathrm{loc}\left(dtdx;L^2\left(Mdv\right)\right)},
		\end{aligned}
	\end{equation*}
	for any $1\leq r<2$.
	
	Then, we end up with
	\begin{equation}\label{first term D4}
		\left\|D_\eps^4(\varphi)\right\|_{L^1_\mathrm{loc}\left(dtdx\right)}\leq C\lambda.
	\end{equation}

	Finally, incorporating \eqref{first term D1},\eqref{first term D2}, \eqref{first term D5 1}, \eqref{first term D5 2}, \eqref{first term D3} and \eqref{first term D4} into \eqref{first term}, we have shown that the first term in the right-hand side of \eqref{renormalized5} is uniformly bounded in $L^1_\mathrm{loc}(dtdx)$ and satisfies, for any $\lambda>0$,
	\begin{equation}\label{renormalized5 first term}
		\limsup_{\eps\rightarrow 0}
		\left\|
		\frac 1{\delta\eps^2}
		\int_{\mathbb{R}^3}
		\hat\gamma_\eps^\lambda \mathcal{Q}\left(G_\eps^\pm,G_\eps^\pm\right)
		\varphi \chi\left( {|v|^2\over K_\delta} \right) Mdv
		\right\|_{L^1_\mathrm{loc}\left(dtdx\right)}
		\leq C\lambda.
	\end{equation}

	In particular, combining \eqref{renormalized5 second term} and \eqref{renormalized5 first term}, it follows from the Banach-Alaoglu theorem, up to further extraction of subsequences, that the right-hand side of \eqref{renormalized5} converges in the weak-$*$ topology of Radon measures $\mathcal{M}_\mathrm{loc}\left([0,\infty)\times\mathbb{R}^3\right)$ towards
	\begin{equation}\label{D limit}
		\int_{\mathbb{R}^3\times\mathbb{R}^3\times\mathbb{S}^2}
		q^{\pm,\mp}
		\varphi MM_*dvdv_*d\sigma
		+
		Q^\lambda(\varphi),
	\end{equation}
	where the Radon measure $Q^\lambda(\varphi)\in\mathcal{M}_{\mathrm{loc}}\left([0,\infty)\times\mathbb{R}^3\right)$ satisfies the control
	\begin{equation}\label{D limit Radon}
		\left\|Q^\lambda(\varphi)\right\|_{\mathcal{M}_\mathrm{loc}\left([0,\infty)\times\mathbb{R}^3\right)}\leq C\lambda.
	\end{equation}

	\subsubsection{Decomposition of flux terms} Next, we treat the convergence in \eqref{renormalized5} of the flux terms
	\begin{equation*}
		\frac1\delta \nabla_x  \cdot \int_{\mathbb{R}^3}
		g_\eps^\pm \gamma_\eps^\lambda\varphi \chi\left( {|v|^2\over K_\delta}\right)v Mdv.
	\end{equation*}
	To this end, we use the decomposition \eqref{fluct-decomposition} to write
	\begin{equation*}
		g_\eps^\pm = \hat g_\eps^\pm +\frac\eps4 \left(\hat g_\eps^\pm\right)^2 = \Pi \hat g^\pm_\eps + \left(\hat g^\pm_\eps -\Pi \hat g^\pm_\eps\right)+\frac\eps4 \left(\hat g_\eps^\pm\right)^2,
	\end{equation*}
	where $\Pi$ is the orthogonal projection on $\Ker\cL$ in $L^2\left(Mdv\right)$, which yields the decomposition of flux terms
	\begin{equation}\label{F0}
		\begin{aligned}
			\frac1\delta \nabla_x \cdot \int_{\mathbb{R}^3}
			g_\eps^\pm \gamma_\eps^\lambda\varphi \chi\left( {|v|^2\over K_\delta}\right)v Mdv
			& =
			\frac\eps{4\delta} \nabla_x \cdot \int_{\mathbb{R}^3}
			\left(\hat g_\eps^\pm \right)^2\gamma_\eps^\lambda\varphi \chi\left( {|v|^2\over K_\delta}\right)v Mdv
			\\
			& + \frac1\delta \nabla_x  \cdot \int_{\mathbb{R}^3}
			\left(\hat g_\eps^\pm-\Pi \hat g^\pm_\eps\right) \gamma_\eps^\lambda\varphi \chi\left( {|v|^2\over K_\delta}\right)v Mdv
			\\
			& +
			\frac1\delta \nabla_x \cdot \int_{\mathbb{R}^3}
			\Pi \hat g_\eps^\pm \varphi \left( \chi\left( {|v|^2\over K_\delta}\right) \gamma_\eps^\lambda-1\right) v Mdv
			\\
			& + \frac1\delta \nabla_x  \cdot \int_{\mathbb{R}^3}
			\Pi \hat g_\eps^\pm\varphi v Mdv
			\\
			& \eqdefa F^1_\eps(\varphi)+F^2_\eps(\varphi)+F^3_\eps(\varphi)+F^4_\eps(\varphi).
		\end{aligned}
	\end{equation}

	Then, from the condition on the support of $\gamma$ and since $\left|\hat g_\eps^\pm\right|=\left|\frac{g_\eps^\pm}{1+\frac\eps 4 \hat g_\eps^\pm}\right|\leq 2\left|g_\eps^\pm\right|$, we deduce that, on the support of $\gamma_\eps^\lambda$,
	\begin{equation*}
		\left|\hat g_\eps^\pm\right|
		\leq 2
		\left|g_\eps^\pm\right| \leq \frac4{\lambda \delta} \leq \frac4\eps.
	\end{equation*}
	Therefore, by Lemma \ref{v2-int}, for any $1\leq p<2$, we have that $\left(\hat g_\eps ^\pm\right)^2 \gamma_\eps^\lambda$ is uniformly bounded in $L^1_\mathrm{loc}\left(dtdx ; L^p(Mdv)\right)$. Hence, we conclude that
	\begin{equation}\label{F1}
		\begin{aligned}
			\| F^1_\eps(\varphi)\| _{W^{-1,1}_\mathrm{loc} (dtdx)}
			\leq C\frac\eps{\delta}
			\| (\hat g_\eps ^\pm)^2  \gamma_\eps^\lambda \|_{L^1_\mathrm{loc}\left(dtdx ; L^p(Mdv)\right)}
			\| v\varphi\| _{L^{p'}(Mdv)} \leq  C \frac\eps{\delta} .
		\end{aligned}
	\end{equation}

	Moreover, by \eqref{relaxation-est}, we easily get
	\begin{equation}\label{F2}
		\begin{aligned}
			\| F^2_\eps(\varphi)\| _{W^{-1,1}_\mathrm{loc} (dtdx)}
			\leq \frac C{\delta} \| \hat g_\eps^\pm -\Pi \hat g_\eps^\pm \|_{L^1_\mathrm{loc}(dtdx ; L^2(Mdv))}
			\| \gamma\|_\infty
			\| v\varphi\| _{L^2(Mdv)}\leq  C \frac\eps{\delta},
		\end{aligned}
	\end{equation}
	which handles the second term.

	Further note that, by definition of the projection $\Pi$, we have, for any $2< p< \infty$,
	\begin{equation*}
		\| \Pi \hat g_\eps^\pm\|_{L^\infty\left(dt;L^2\left(dx;L^p\left(Mdv\right)\right)\right)} \leq C_p \| \hat g_\eps^\pm \|_{L^\infty\left(dt;L^2\left(Mdxdv\right)\right)},
	\end{equation*}
	whence
	\begin{equation*}
		\begin{aligned}
			\| F^3_\eps(\varphi) & \| _{W^{-1,1}_\mathrm{loc} (dtdx)}
			\\
			& \leq C_p \| \hat g_\eps^\pm \|_{L^\infty\left(dt;L^2\left(Mdxdv\right)\right)}
			\left\| {1-\gamma_\eps^\lambda \over \delta} \right\|_{L^2_\mathrm{loc}(dtdx ; L^2(Mdv))} \| v\varphi\|_{L^q(Mdv)}
			\\
			& + C  \frac1\delta \| \hat g_\eps^\pm \|_{L^\infty\left(dt;L^2(Mdxdv)\right)} \left\| v\varphi \left( 1- \chi\left( {|v|^2\over K_\delta}\right)\right) \right\| _{L^{2}(Mdv)},
		\end{aligned}
	\end{equation*}
	with $\frac 1q = \frac 12 - \frac 1p$. Then, using estimate \eqref{gamma bound} (with $\gamma_\eps^\lambda$ instead of $\hat \gamma_\eps^\lambda$) and the control of Gaussian tails \eqref{gaussian-decay 0} to respectively bound the first and second terms in the right-hand side above, we infer that
	\begin{equation*}
		\| F^3_\eps(\varphi)\| _{W^{-1,1}_\mathrm{loc} (dtdx)}\leq C\lambda + C \delta^{\frac K4 -1} |\log \delta|^{\frac74},
	\end{equation*}
	which is small provided that $K>4$. Therefore, up to further extraction of subsequences, we deduce that
	\begin{equation}\label{F3}
		F^3_\eps(\varphi)\rightharpoonup \nabla_x \cdot R^\lambda(\varphi)
		\quad
		\text{in the sense of distributions,}
	\end{equation}
	where the Radon measure $R^\lambda(\varphi)\in\mathcal{M}_{\mathrm{loc}}\left([0,\infty)\times\mathbb{R}^3\right)$ satisfies the control
	\begin{equation}\label{F limit Radon}
		\left\|R^\lambda(\varphi)\right\|_{\mathcal{M}_\mathrm{loc}\left([0,\infty)\times\mathbb{R}^3\right)}\leq C\lambda.
	\end{equation}

	The form of the last remaining flux term $F^4_\eps(\varphi)$ depends on the collision invariant $\varphi$~:
	\begin{itemize}
		\item If $\varphi (v) = v$, we get, using that $\phi(v)=v\otimes v - \frac{|v|^2}{3}\operatorname{Id}$ is orthogonal to the collision invariants,
		\begin{equation}\label{F4 1}
			\begin{aligned}
				F^4_\eps(\varphi) & = \frac1\delta \nabla_x  \cdot \int_{\mathbb{R}^3}
			\Pi \hat g_\eps^\pm v\otimes v Mdv
			\\
			& = \frac1{3\delta} \nabla_x \int_{\mathbb{R}^3}
			\Pi \hat g_\eps^\pm|v|^2 Mdv
			\\
			& = \frac 1\delta\nabla_x\left(\hat\rho_\eps^\pm + \hat \theta_\eps^\pm \right),
			\end{aligned}
		\end{equation}
		where $\hat \rho_\eps^\pm$ and $\hat\theta_\eps^\pm$ denote the densities and temperatures respectively associated with the renormalized fluctuations $\hat g_\eps^\pm$. Thus, this term takes the form of a gradient and will therefore vanish upon integrating it against divergence free vector fields as required by the theory of weak solutions of Leray.
		
		\item If $\varphi(v) = \frac{|v|^2}{2}-\frac 52$, we obtain
		\begin{equation}\label{F4 2}
			F^4_\eps(\varphi) = \frac1\delta \nabla_x \cdot \int_{\mathbb{R}^3}
			\Pi \hat g_\eps^\pm \left(\frac{|v|^2}{2}-\frac 52\right) v Mdv = 0,
		\end{equation}
		for $\psi(v)=\left(\frac{|v|^2}{2}-\frac 52\right)v$ is orthogonal to the collision invariants.
	\end{itemize}

	Thus, on the whole incorporating \eqref{F1}, \eqref{F2}, \eqref{F3}, \eqref{F4 1} and \eqref{F4 2} into \eqref{F0}, we conclude that the flux terms satisfy the following convergences in the sense of distributions~:
	\begin{equation}\label{F limit}
		\begin{aligned}
			P\left(\frac1\delta \nabla_x  \cdot \int_{\mathbb{R}^3}
			g_\eps^\pm \gamma_\eps^\lambda \chi\left( {|v|^2\over K_\delta}\right)v\otimes v Mdv\right)
			&\rightharpoonup
			P\left(\nabla_x\cdot R^\lambda(v)\right),
			\\
			\frac1\delta \nabla_x  \cdot \int_{\mathbb{R}^3}
			g_\eps^\pm \gamma_\eps^\lambda \chi\left( {|v|^2\over K_\delta}\right)\left(\frac{|v|^2}{2}-\frac 52\right)v Mdv
			&\rightharpoonup
			\nabla_x\cdot R^\lambda\left(\frac{|v|^2}{2}-\frac 52\right),
		\end{aligned}
	\end{equation}
	where $P$ denotes the Leray projector onto solenoidal vector fields.

	\subsubsection{Decomposition of acceleration terms} It only remains to deal with the terms involving the electromagnetic field in \eqref{renormalized5}, which we decompose as
	\begin{equation}\label{A0}
		\begin{aligned}
			& \int_{\mathbb{R}^3}
			g_\eps^\pm \gamma_\eps^\lambda
			(\eps E_\eps+v\wedge B_\eps) \cdot \nabla_v \left( \varphi \chi\left( {|v|^2\over K_\delta}\right) M \right) dv
			\\
			& +
			E_\eps \cdot
			\int_{\mathbb{R}^3} \left(1+\eps g_\eps^\pm\right) \hat \gamma_\eps^\lambda \varphi\chi\left( {|v|^2\over K_\delta}\right) vM dv
			\\
			& =
			\left[\eps E_\eps \cdot \int_{\mathbb{R}^3}
			g_\eps^\pm \gamma_\eps^\lambda
			\nabla_v \left( \varphi \chi\left( {|v|^2\over K_\delta}\right) M \right) dv
			+ \eps E_\eps \cdot
			\int_{\mathbb{R}^3} g_\eps^\pm \hat \gamma_\eps^\lambda \varphi\chi\left( {|v|^2\over K_\delta}\right) vM dv
			\right]
			\\
			& + \left[
			E_\eps \cdot \int_{\mathbb{R}^3} \varphi vM dv
			+
			E_\eps \cdot
			\int_{\mathbb{R}^3} \left( \chi\left( {|v|^2\over K_\delta}\right) \hat \gamma_\eps^\lambda-1\right) \varphi vM dv
			\right]
			\\
			& -\left[
			B_\eps \cdot \int_{\mathbb{R}^3}
			 g_\eps^\pm \gamma_\eps^\lambda
			v\wedge \nabla_v  \varphi M dv
			+ B_\eps \cdot \int_{\mathbb{R}^3}
			g_\eps^\pm \gamma_\eps^\lambda
			v\wedge (\nabla_v \varphi)
			\left( \chi\left( {|v|^2\over K_\delta}\right)-1 \right) Mdv
			\right]
			\\
			& \eqdefa
			A^1_\eps(\varphi)+A^2_\eps(\varphi)-A^3_\eps(\varphi).
		\end{aligned}
	\end{equation}

	From the condition on the support of $\gamma$ and since $\left|\hat g_\eps^\pm\right|=\left|\frac{g_\eps^\pm}{1+\frac\eps 4 \hat g_\eps^\pm}\right|\leq 2\left|g_\eps^\pm\right|$, we clearly have that
	\begin{equation}\label{A1}
		\begin{aligned}
			& \| A^1_\eps(\varphi) \|_{L^\infty\left(dt;L^1(dx)\right)}
			\\
			& \leq C \eps \| E_\eps\|_{L^\infty\left(dt;L^2(dx)\right)}
			\left\| g_\eps^\pm \mathds{1}_{\left\{\left|\lambda\delta g_\eps^\pm\right|\leq 2\right\}} \right\|_{L^\infty\left(dt;L^2(Mdxdv)\right)}
			\left\| \varphi \chi\left( {|v|^2\over K_\delta}\right)v\right\|_{H^1(Mdv)}
			\\
			& \leq C \eps
			\left\| \left(\hat g_\eps^\pm+\frac\eps 4\hat g_\eps^{\pm 2}\right) \mathds{1}_{\left\{\left|\lambda\delta \hat g_\eps^\pm\right|\leq 4\right\}} \right\|_{L^\infty\left(dt;L^2(Mdxdv)\right)}
			\leq C \eps
			\left\| \hat g_\eps^\pm \right\|_{L^\infty\left(dt;L^2(Mdxdv)\right)},
		\end{aligned}
	\end{equation}
	which handles the first acceleration term.

	Then, in order to deal with the second acceleration term, we estimate first, using \eqref{gamma bound} and the control of Gaussian tails \eqref{gaussian-decay 0}, that
	\begin{equation*}
		\begin{aligned}
			& \left\|\int_{\mathbb{R}^3} \left( \chi\left( {|v|^2\over K_\delta}\right) \hat \gamma_\eps^\lambda-1\right) \varphi vM dv\right\|_{L^2_\mathrm{loc}(dtdx)}
			\\
			& \leq
			\left\|\int_{\mathbb{R}^3} \left( \hat \gamma_\eps^\lambda-1\right) \varphi vM dv\right\|_{L^2_\mathrm{loc}(dtdx)}
			+
			\left\|\int_{\mathbb{R}^3} \left( \chi\left( {|v|^2\over K_\delta}\right)-1\right) \varphi vM dv\right\|_{L^2_\mathrm{loc}(dtdx)}
			\\
			& \leq C\lambda\delta
			\left\|\hat g_\eps^\pm\right\|_{L^2_\mathrm{loc}(dtdx;L^2(Mdv))}
			+ C\delta^\frac K2 |\log \delta|^2,
		\end{aligned}
	\end{equation*}
	whence, as $\eps\rightarrow 0$,
	\begin{equation}\label{A2}
		A_\eps^2(\varphi)\stackrel{*}{\rightharpoonup} E\cdot\int_{\mathbb{R}^3}\varphi v Mdv
		\qquad\text{in }L^\infty\left(dt;L^2(dx)\right).
	\end{equation}

	As for the remaining term $A_\eps^3(\varphi)$, note first, using \eqref{gaussian-decay 0} again, that
	\begin{equation*}
		\begin{aligned}
			& \left\|
			B_\eps \cdot \int_{\mathbb{R}^3}
			g_\eps^\pm \gamma_\eps^\lambda
			v\wedge (\nabla_v \varphi)
			\left( \chi\left( {|v|^2\over K_\delta}\right)-1 \right) Mdv
			\right\|_{L^\infty(dt;L^2(dx))}
			\\
			& \leq \frac{C}{\lambda\delta}\left\| B_\eps\right\|_{L^\infty(dt;L^2(dx))}
			\int_{\mathbb{R}^3}
			\left( \chi\left( {|v|^2\over K_\delta}\right)-1 \right)
			|v|^2
			Mdv
			\leq
			\frac C\lambda \delta^{\frac K2 - 1}|\log\delta|^\frac 32,
		\end{aligned}
	\end{equation*}
	which is small as soon as $K>2$. Moreover, in view of \eqref{fluct-decomposition}, we find
	\begin{equation*}
		\begin{aligned}
			& \left\|
			B_\eps \cdot \int_{\mathbb{R}^3}
			\left(g_\eps^\pm - \hat g_\eps^\pm\right) \gamma_\eps^\lambda
			v\wedge \nabla_v \varphi M dv
			\right\|_{L^1_\mathrm{loc}(dtdx)}
			\\
			& \hspace{30mm} =
			\left\|B_\eps \cdot \int_{\mathbb{R}^3}
			\frac \eps 4 \hat g_\eps^{\pm 2} \gamma_\eps^\lambda
			v\wedge \nabla_v \varphi M dv
			\right\|_{L^1_\mathrm{loc}(dtdx)}
			\\
			& \hspace{30mm} \leq C\frac \eps{\lambda\delta}
			\left\|B_\eps\right\|_{L^2_\mathrm{loc}(dtdx)}\left\|\hat g_\eps^{\pm}\right\|_{L^2_\mathrm{loc}\left(dtdx;L^2(Mdv)\right)},
		\end{aligned}
	\end{equation*}
	so that, on the whole, the weak limit of $A_\eps^3(\varphi)$ will coincide with the weak limit of
	\begin{equation*}
		B_\eps \cdot \int_{\mathbb{R}^3}
		\hat g_\eps^\pm \gamma_\eps^\lambda
		v\wedge \nabla_v \varphi M dv.
	\end{equation*}
	In order to take the weak limit of the preceding term, notice, in view of Lemma \ref{x-compactness2 0} through a straighforward application of the mean value theorem to the function $z\gamma\left(\lambda\delta\left(z+\frac\eps 4z^2\right)\right)$, that $\hat g_\eps^\pm \gamma_\eps^\lambda=\hat g_\eps^\pm \gamma\left(\lambda\delta\left(\hat g_\eps^\pm+\frac \eps 4\hat g_\eps^{\pm 2}\right)\right)$ is locally relatively compact in $(x,v)$ in $L^p(dtdxdv)$, for any $1\leq p<2$. In fact, Lemma \ref{x-compactness2 0} further implies that $\left|\hat g_\eps^\pm \gamma_\eps^\lambda\right|^2\leq \left|\hat g_\eps^\pm\right|^2\mathds{1}_{\left\{\lambda\delta\left|\hat g_\eps^\pm\right|\leq 4\right\}}$ is equi-integrable. Therefore, we conclude that $\hat g_\eps^\pm \gamma_\eps^\lambda$ is locally relatively compact in $(x,v)$ in $L^2(dtdxdv)$. In particular, for any fixed $\lambda>0$, it is possible to approximate $\hat g_\eps^\pm \gamma_\eps^\lambda$, uniformly in $\eps>0$, in $L^2_\mathrm{loc}(dtdxdv)$ by its regularized version $\left(\hat g_\eps^\pm \gamma_\eps^\lambda\right) *_{x,v}\chi_{a}$, where $a>0$ and $\chi_a(x,v)=\frac{1}{a^6}\chi\left(\frac{x}{a},\frac{v}{a}\right)$ is an approximate identity, with $\chi\in C_c^\infty\left(\mathbb{R}^3\times\mathbb{R}^3\right)$ such that $\int_{\mathbb{R}^3\times\mathbb{R}^3}\chi(x,v)dxdv=1$.

	We use now compensated compactness in the following form. From the Faraday equation in \eqref{VMB2}, we deduce that
	\begin{equation*}
		\d_t B_\eps \in {L^\infty\left(dt;H^{-1}(dx)\right)},
	\end{equation*}
	so that $B_\eps$ enjoys some strong compactness with respect to the time variable. We then deduce, up to extraction of subsequences, that
	\begin{equation*}
		B_\eps \left(\hat g_\eps^\pm \gamma_\eps^\lambda\right) *_{x,v}\chi_{a}
		\rightharpoonup
		B g^{\pm}*_{x,v}\chi_{a},
	\end{equation*}
	where $g^\pm$ is the weak limit of $\hat g_\eps^\pm \gamma_\eps^\lambda$, which coincides with the weak limits of $\hat g_\eps^\pm$ and $g_\eps^\pm$ (note that $\gamma_\eps^\lambda\rightarrow 1$ almost everywhere). Incidentally, by the uniformity of the approximation of $\hat g_\eps^\pm \gamma_\eps^\lambda$ by $\left(\hat g_\eps^\pm \gamma_\eps^\lambda\right) *_{x,v}\chi_{a}$ in $L_\mathrm{loc}^2(dtdxdv)$, we infer that
	\begin{equation*}
		B_\eps \hat g_\eps^\pm \gamma_\eps^\lambda
		\rightharpoonup
		B g^{\pm},
	\end{equation*}
	in $L^1_\mathrm{loc}(dtdxdv)$, whence
	\begin{equation}\label{A3}
		A_\eps^3(\varphi)\rightharpoonup
		B \cdot \int_{\mathbb{R}^3}
		g^\pm
		v\wedge \nabla_v \varphi M dv
		\qquad\text{in }L^1_\mathrm{loc}(dtdx).
	\end{equation}

	Thus, incorporating \eqref{A1}, \eqref{A2} and \eqref{A3} into \eqref{A0}, we finally conclude that
	\begin{equation}\label{A limit}
		\begin{aligned}
			\int_{\mathbb{R}^3}
			g_\eps^\pm \gamma_\eps^\lambda
			(\eps E_\eps+v\wedge B_\eps) \cdot \nabla_v & \left( \varphi \chi\left( {|v|^2\over K_\delta}\right) M \right) dv
			\\
			& +
			E_\eps \cdot
			\int_{\mathbb{R}^3} \left(1+\eps g_\eps^\pm\right) \hat \gamma_\eps^\lambda \varphi\chi\left( {|v|^2\over K_\delta}\right) vM dv
			\\
			& \rightharpoonup
			E\cdot\int_{\mathbb{R}^3}\varphi v Mdv-B \cdot \int_{\mathbb{R}^3}
			g^\pm
			v\wedge \nabla_v \varphi M dv,
		\end{aligned}
	\end{equation}
	in $L^1_\mathrm{loc}(dtdx)$.

	\subsubsection{Convergence} We are now in a position to pass to the limit in \eqref{renormalized5}. To this end, note first that, since $\frac\eps\delta\rightarrow 0$ and $g_\eps^\pm$ is uniformly bounded in $L^1_\mathrm{loc}\left(dtdx;L^1\left(\left(1+|v|^2\right)Mdv\right)\right)$ by Lemma \ref{L1-lem}, the density term
	\begin{equation*}
		\frac\eps\delta\d_t \int_{\mathbb{R}^3} g_\eps^\pm \gamma_\eps^\lambda \varphi \chi\left( {|v|^2\over K_\delta} \right) Mdv,
	\end{equation*}
	vanishes as $\eps\rightarrow 0$ and thus brings no contribution to the weak limit. Therefore, according to the weak limits of conservation defects \eqref{D limit}, flux terms \eqref{F limit} and acceleration terms \eqref{A limit}, we conclude, letting $\eps\rightarrow 0$ in \eqref{renormalized5} in the sense of distributions for the collision invariants $\varphi=v$ and $\varphi=\frac{|v|^2}{2}-\frac 52$, that
	\begin{equation*}
		\begin{aligned}
			P \big(\nabla_x\cdot R^\lambda & (v)
			\mp \left(E + u\wedge B\right) \big)
			\\
			& =P\left(
			\int_{\mathbb{R}^3\times\mathbb{R}^3\times\mathbb{S}^2}
			q^{\pm,\mp}
			v MM_*dvdv_*d\sigma
			+
			Q^\lambda(v)\right),
			\\
			\nabla_x\cdot R^\lambda & \left(\frac{|v|^2}{2}-\frac 52\right)
			\\
			& =
			\int_{\mathbb{R}^3\times\mathbb{R}^3\times\mathbb{S}^2}
			q^{\pm,\mp}
			\left(\frac{|v|^2}{2}-\frac 52\right) MM_*dvdv_*d\sigma
			+
			Q^\lambda\left(\frac{|v|^2}{2}-\frac 52\right),
		\end{aligned}
	\end{equation*}
	where $u$ denotes the bulk velocity associated with the limiting fluctuations $g^\pm$ (recall that, according to Lemmas \ref{relaxation-control} and \ref{relaxation2-control}, $g^+$ and $g^-$ are infinitesimal Maxwellians which differ only by their densities).

	Next, in view of the bounds \eqref{D limit Radon} and \eqref{F limit Radon} on the Radon measures $Q^\lambda$ and $R^\lambda$, respectively, we deduce, by the arbitrariness of $\lambda>0$, that
	\begin{equation*}
		\begin{aligned}
			P \left(E + u\wedge B
			\pm
			\int_{\mathbb{R}^3\times\mathbb{R}^3\times\mathbb{S}^2}
			q^{\pm,\mp}
			v MM_*dvdv_*d\sigma \right) & = 0,
			\\
			\int_{\mathbb{R}^3\times\mathbb{R}^3\times\mathbb{S}^2}
			q^{\pm,\mp}
			\left(\frac{|v|^2}{2}-\frac 52\right) MM_*dvdv_*d\sigma
			& = 0,
		\end{aligned}
	\end{equation*}
	which concludes the proof of Proposition \ref{high weak-comp}.\qed

\section[Macroscopic constraint equations for two species\ldots]{Macroscopic constraint equations for two species, strong interactions}\label{high constraints 1}

We move on now to the regime of strong interspecies interactions considered in Theorem \ref{CV-OMHDSTRONG}. In this setting, the derivation of even the simplest macroscopic constraint equations (such as the incompressibility and Boussinesq contraints) involves the handling of nonlinear terms. Indeed, the limiting kinetic equation \eqref{kinetic equation 4} for weak interspecies interactions (obtained in Chapter \ref{constraints proof} with weak compactness methods) corresponds now, for strong interspecies interactions, to the nonlinear equation
\begin{equation*}
	( v \cdot \nabla_x \pm \left(v\wedge B\right)\cdot \nabla_v )
	g^\pm
	\mp E \cdot v
	=
	\int_{\mathbb{R}^3\times\mathbb{S}^2} \left(q^\pm + q^{\pm,\mp}\right) M_*dv_*d\sigma,
\end{equation*}
which is obtained in the proof of Proposition \ref{high weak-comp2}, below (see \eqref{kinetic equation 3}), and requires the compactness properties established in Chapter \ref{hypoellipticity}.

The next result fully characterizes the limiting kinetic equations in the regime of strong interactions.

\begin{prop}\label{high weak-comp2}
	Let $\left(f_\eps^\pm, E_\eps, B_\eps\right)$ be the sequence of renormalized solutions to the scaled two species Vlasov-Maxwell-Boltzmann system \eqref{VMB2} considered in Theorem \ref{CV-OMHDSTRONG} for strong interspecies interactions, i.e.\ $\delta=1$. In accordance with Lemmas \ref{L1-lem}, \ref{L2-lem} and \ref{L2-qlem}, denote by
	\begin{equation*}
		g^\pm\in L^\infty\left(dt;L^2\left(Mdxdv\right)\right)
		\qquad\text{and}\qquad
		q^\pm,q^{\pm,\mp} \in L^2\left(MM_*dtdxdvdv_*d\sigma\right)
	\end{equation*}
	any joint limit points of the families $\hat g_\eps^\pm$, $\hat q_\eps^\pm$ and $\hat q_\eps^{\pm,\mp}$ defined by \eqref{hatg} and \eqref{hatq-def}, respectively, and by
	\begin{equation*}
		E,B\in L^\infty\left(dt;L^2\left(dx\right)\right)
	\end{equation*}
	any joint limit points of the families $E_\eps$ and $B_\eps$, respectively.
	
	Then, one has
	% , at this stage, only a conditional result. Further assume that
	% 	\begin{equation}\label{high equiintegrability-ch4}
	% 		(\hat g_\eps^+)^2 \text{ and } (\hat g_\eps^+)^2 \text{ are weakly relatively compact in } L^1_\mathrm{loc}\left(dtdx, L^1(Mdv)\right),
	% 	\end{equation}
	% 	then
	\begin{equation}\label{q phi psi 3}
		\frac 12 \int_{\mathbb{R}^3\times\mathbb{S}^2} \left(q^++q^-+q^{+,-}+q^{-,+}\right) M_*dv_*d\sigma
		= \phi:\nabla_x u +\psi\cdot \nabla_x \theta ,
	 \end{equation}
	and
	\begin{equation}\label{pre Ohm strong}
		\begin{aligned}
			\frac 12 \int_{\mathbb{R}^3\times\mathbb{S}^2} \left(q^+-q^-+q^{+,-}-q^{-,+}\right) & M_*dv_*d\sigma
			\\
			& = \left(\frac 12 \nabla_x(\rho^+ - \rho^-) - \left(E+u\wedge B\right)\right)\cdot v,
		 \end{aligned}
	 \end{equation}
	where $\rho^\pm$, $u$ and $\theta $ are, respectively, the densities, bulk velocity and temperature associated with the limiting fluctuations $g^\pm$, and $\phi$ and $\psi$ are the kinetic fluxes defined by \eqref{phi-psi-def}. Furthermore, $\rho^\pm$, $u$ and $\theta$ satisfy the following constraints
	\begin{equation}\label{constraints2-1}
		\Div u = 0, \qquad \nabla_x\left(\frac{\rho^++\rho^-}{2}+\theta\right) = 0.
	\end{equation}
	In particular, the strong Boussinesq relation $\frac{\rho^++\rho^-}{2}+\theta=0$ holds.
\end{prop}

\begin{proof}
		The case $\delta=1$ is more complicated because we do not have enough strong compactness to take limits in the nonlinear terms
		\begin{equation*}
			\pm \left(v\wedge B_\eps\right)\cdot \nabla_v
			{\sqrt{ G_\eps^\pm +\eps^a }-1 \over \eps}.
		\end{equation*}
		More precisely, we are not able to control the concentrations of $|\hat g_\eps^\pm|^2$ (see Lemma \ref{x-compactness2 0}).
		
		The idea is therefore to consider a stronger renormalization of the equation for the fluctuations of density. To this end, we introduce the admissible renormalization $\Gamma_\lambda(z)$ defined by
		\begin{equation*}
			\Gamma_\lambda(z)-1=(z-1)\gamma\left(\lambda\frac{z-1}{\eps}\right),
		\end{equation*}
		where $\lambda>0$ is small and $\gamma \in C^1\left(\mathbb{R}\right)$ satisfies that
		\begin{equation*}
			\mathds{1}_{[-1,1]}(z)\leq \gamma(z) \leq \mathds{1}_{[-2,2]}(z), \qquad \text{for all }z\in\mathbb{R}.
		\end{equation*}
		Without distinguishing, for simplicity, the notation for cations and anions, we denote $\gamma_\eps^\lambda$ for $\gamma\left(\lambda g_\eps^\pm\right)$ and $\hat\gamma_\eps^\lambda$ for $\Gamma_\lambda'\left(G_\eps^\pm\right)$. Thus, renormalizing the Vlasov-Boltzmann equation from \eqref{VMB2} with respect to $\Gamma_\lambda(z)$ yields
		\begin{equation}\label{renormalized3}
			\begin{aligned}
				( \eps \d_t & + v \cdot \nabla_x \pm \left(\eps E_\eps+ v\wedge B_\eps\right)\cdot \nabla_v )
				g_\eps^\pm\gamma_\eps^\lambda
				\mp E_\eps \cdot v G_\eps^\pm \hat\gamma_\eps^\lambda
				\\
				& =
				\hat\gamma_\eps^\lambda\sqrt{G_\eps^\pm}
				\int_{\mathbb{R}^3\times\mathbb{S}^2} \sqrt{ G_{\eps *}^\pm } \hat q_\eps^\pm M_*dv_*d\sigma
				+
				\frac{\eps^2}{4}\hat\gamma_\eps^\lambda
				\int_{\mathbb{R}^3\times\mathbb{S}^2} \left(\hat q_\eps^\pm\right)^2 M_*dv_*d\sigma
				\\
				& +
				\hat\gamma_\eps^\lambda\sqrt{G_\eps^\pm}
				\int_{\mathbb{R}^3\times\mathbb{S}^2} \sqrt{ G_{\eps *}^\mp } \hat q_\eps^{\pm,\mp} M_*dv_*d\sigma
				+
				\frac{\eps^2}{4}\hat\gamma_\eps^\lambda
				\int_{\mathbb{R}^3\times\mathbb{S}^2} \left(\hat q_\eps^{\pm,\mp}\right)^2 M_*dv_*d\sigma.
			\end{aligned}
		\end{equation}

		Next, employing a strategy similar to the proof of Proposition \ref{weak-comp}, in particular, since $\hat\gamma_\eps^\lambda\sqrt{G_\eps^\pm}$ is uniformly bounded pointwise, utilizing that, thanks to Lemma \ref{L2-lem},
		\begin{equation*}
			\begin{aligned}
				\sqrt{G_{\eps}^\pm} & = 1 + O(\eps)_{L^2_\mathrm{loc}\left(dt ; L^2\left(M dxdv\right)\right)}, \\
				\sqrt{G_{\eps *}^\pm} & = 1 + O(\eps)_{L^2_\mathrm{loc}\left(dt ; L^2\left(M_* dxdv_*\right)\right)},
			\end{aligned}
		\end{equation*}
		and that, thanks to Lemma \ref{L2-qlem}, the collision integrands $\hat q_\eps^{\pm}$ and $\hat q_\eps^{\pm,\mp}$ are uniformly bounded in $L^2\left(MM_*dtdxdvdv_*d\sigma\right)$, we see that the weak limit of the right-hand side of \eqref{renormalized3} coincides with the weak limit of
		\begin{equation*}
			\int_{\mathbb{R}^3\times\mathbb{S}^2} \left(\hat q_\eps^\pm + \hat q_\eps^{\pm,\mp}\right) M_*dv_*d\sigma
			+ Q_\eps^\lambda,
		\end{equation*}
		where we denote the remainder
		\begin{equation*}
			Q_\eps^\lambda
			=
			\left[\hat\gamma_\eps^\lambda - 1\right]
			\int_{\mathbb{R}^3\times\mathbb{S}^2} \left(\hat q_\eps^\pm + \hat q_\eps^{\pm,\mp}\right) M_*dv_*d\sigma.
		\end{equation*}
		Then, since $\left|\Gamma_\lambda'(z)-1\right|\leq C_r\left|\lambda\frac{z-1}{\eps}\right|^\frac{1}{r}$, for any given $1\leq r\leq\infty$ and for every $z\geq 0$, it holds that, employing the uniform bounds from Lemmas \ref{L1-lem} and \ref{L2-qlem}, for any $2\leq r\leq\infty$,
		\begin{equation*}
			\begin{aligned}
				\left|Q_\eps^\lambda\right|
				& \leq C\lambda^{\frac 1r}\left|g_\eps^\pm\right|^{\frac 1r}
				\int_{\mathbb{R}^3\times\mathbb{S}^2} \left(\hat q_\eps^\pm + \hat q_\eps^{\pm,\mp}\right) M_*dv_*d\sigma
				\\
				& \leq C\lambda^{\frac 1r}\left|g_\eps^\pm\right|^{\frac 1r}
				\left(\int_{\mathbb{R}^3\times\mathbb{S}^2} \left(\hat q_\eps^\pm + \hat q_\eps^{\pm,\mp}\right)^2 M_*dv_*d\sigma\right)^\frac{1}{2}
				\\
				& = O\left(\lambda^\frac{1}{r}\right)_{L^\frac{2r}{2+r}_\mathrm{loc}(dtdxdv)}.
			\end{aligned}
		\end{equation*}
		Moreover, when $r=2$, it is readily seen, in view of the weak relative compactness of $g_\eps^\pm$ in $L^1_\mathrm{loc}(dtdxdv)$ established in Lemma \ref{L1-lem} and employing the Dunford-Pettis compactness criterion (see \cite{royden}), that $Q_\eps^\lambda$ is weakly relatively compact in $L^1_\mathrm{loc}(dtdxdv)$, as well.
		
		Therefore, up to extraction of a further subsequence as $\eps\rightarrow 0$, we may assume that $Q_\eps^\lambda$ converges weakly in $L^\frac{2r}{2+r}_\mathrm{loc}(dtdxdv)$ to some $Q^\lambda\in L^\frac{2r}{2+r}_\mathrm{loc}(dtdxdv)$, for any $2\leq r\leq\infty$, whose magnitude is at most of order $\lambda^\frac{1}{r}$. On the whole, we have evaluated that the right-hand side of \eqref{renormalized3} converges weakly towards
		\begin{equation*}
			\begin{aligned}
				\int_{\mathbb{R}^3\times\mathbb{S}^2} & \left( q^\pm +  q^{\pm,\mp}\right) M_*dv_*d\sigma
				+ Q^\lambda
				\\
				& =
				\int_{\mathbb{R}^3\times\mathbb{S}^2} \left( q^\pm + q^{\pm,\mp}\right) M_*dv_*d\sigma
				+ O\left(\lambda^\frac{1}{r}\right)_{L^\frac{2r}{2+r}_\mathrm{loc}(dtdxdv)}.
			\end{aligned}
		\end{equation*}

		As for the left-hand side of \eqref{renormalized3}, we first have that
		\begin{equation}\label{renormalized3 1}
			E_\eps \cdot v G_\eps^\pm \hat \gamma_\eps^\lambda
			=E_\eps \cdot v
			+E_\eps \cdot v \left(\hat \gamma_\eps^\lambda - 1\right)
			+ \eps E_\eps \cdot v g_\eps^\pm \hat \gamma_\eps^\lambda.
		\end{equation}
		Therefore, since $g_\eps^\pm \hat \gamma_\eps^\lambda$ is uniformly bounded pointwise and $\left|\Gamma_\lambda'(z)-1\right|\leq C_r\left|\lambda\frac{z-1}{\eps}\right|^\frac{1}{r}$, for any given $1\leq r\leq\infty$ and for every $z\geq 0$, we find, for any $2\leq r\leq\infty$,
		\begin{equation*}
			E_\eps \cdot v G_\eps^\pm \hat \gamma_\eps^\lambda
			=E_\eps \cdot v
			+O\left(\lambda^\frac{1}{r}\right)_{L^\frac{2r}{2+r}_\mathrm{loc}(dtdxdv)}
			+O\left(\frac\eps\lambda\right)_{L^2_{\mathrm{loc}}\left(dtdxdv\right)},
		\end{equation*}
		so that the expression from \eqref{renormalized3 1} converges weakly, for any $2\leq r\leq\infty$, towards
		\begin{equation*}
			E\cdot v + O\left(\lambda^\frac{1}{r}\right)_{L^\frac{2r}{2+r}_\mathrm{loc}(dtdxdv)}.
		\end{equation*}

		Next, since $g_\eps^\pm$ is weakly relatively compact in $L^1_\mathrm{loc}(dtdxdv)$, by Lemma \ref{L1-lem}, it holds that $g_\eps^\pm \left(1-\gamma_\eps^\lambda\right)$ is uniformly small in $L^1_\mathrm{loc}(dtdxdv)$, when $\lambda>0$ is small. In particular, we deduce that the family $g_\eps^\pm \gamma_\eps^\lambda$, which is weakly relatively compact in $L^p_\mathrm{loc}\left(dtdxdv\right)$, for any $1\leq p<\infty$, converges weakly towards some $g^{\pm,\lambda}\in L^p_\mathrm{loc}\left(dtdxdv\right)$ such that
		\begin{equation*}
			g^{\pm,\lambda} = g^\pm+o(1)_{L^1_\mathrm{loc}(dtdxdv)},
		\end{equation*}
		as $\lambda\rightarrow 0$. In fact, we claim that this can be improved to
		\begin{equation}\label{endpoint claim}
			g^{\pm,\lambda} = g^\pm+o(1)_{L^2_\mathrm{loc}(dtdxdv)},
		\end{equation}
		as $\lambda\rightarrow 0$. Indeed, up to extraction of subsequences, denoting by $r^\pm\in L^2_\mathrm{loc}(dtdxdv)$ the weak limit of $\left|g_\eps^\pm\right|$ in $L^1_\mathrm{loc}(dtdxdv)$, which coincides, in view of \eqref{fluct-decomposition}, with the weak limit of $\left|\hat g_\eps^\pm\right|$ in $L^2_\mathrm{loc}(dtdxdv)$, it clearly holds that, for any non-negative $\varphi\in C_c^\infty\left([0,\infty)\times\mathbb{R}^3\times\mathbb{R}^3\right)$,
		\begin{equation*}
			\begin{aligned}
				\int_{[0,\infty)\times\mathbb{R}^3\times\mathbb{R}^3} &
				\left|g^{\pm} - g^{\pm,\lambda} \right|\varphi dtdxdv
				\\
				& \leq \liminf_{\eps\rightarrow 0}
				\int_{[0,\infty)\times\mathbb{R}^3\times\mathbb{R}^3}
				\left|g_\eps^\pm \left(1-\gamma_\eps^\lambda\right)\right|\varphi dtdxdv
				\\
				& \leq \liminf_{\eps\rightarrow 0}
				\int_{[0,\infty)\times\mathbb{R}^3\times\mathbb{R}^3}
				\left| g_\eps^\pm \right|\varphi dtdxdv
				\\
				& = \int_{[0,\infty)\times\mathbb{R}^3\times\mathbb{R}^3}
				r^\pm \varphi dtdxdv,
			\end{aligned}
		\end{equation*}
		whence $\left|g^{\pm} - g^{\pm,\lambda}\right|^2$ is dominated by the integrable function $\left(r^\pm\right)^2$. It then follows from a direct application of Lebesgue's dominated convergence theorem that \eqref{endpoint claim} holds.

		At last, we deal with the convergence of the problematic nonlinear term $B_\eps g_\eps^\pm \gamma_\eps^\lambda$ from \eqref{renormalized3}. To this end, note first, according to \eqref{fluct-decomposition} and Lemma \ref{x-compactness2 0}, for any given $\lambda>0$, that the family $g_\eps^\pm \gamma_\eps^\lambda$ is locally relatively compact in $(x,v)$ in $L^2(dtdxdv)$. In particular, for any fixed $\lambda>0$, it is possible to approximate $g_\eps^\pm \gamma_\eps^\lambda$, uniformly in $\eps>0$, in $L^2_\mathrm{loc}(dtdxdv)$ by its regularized version $\left(g_\eps^\pm \gamma_\eps^\lambda\right) *_{x,v}\chi_{a}$, where $a>0$ and $\chi_a(x,v)=\frac{1}{a^6}\chi\left(\frac{x}{a},\frac{v}{a}\right)$ is an approximate identity, with $\chi\in C_c^\infty\left(\mathbb{R}^3\times\mathbb{R}^3\right)$ such that $\int_{\mathbb{R}^3\times\mathbb{R}^3}\chi(x,v)dxdv=1$.

		We use now compensated compactness in the following form. From the Faraday equation in \eqref{VMB2}, we deduce that
		\begin{equation*}
			\d_t B_\eps \in {L^\infty\left(dt;H^{-1}(dx)\right)},
		\end{equation*}
		so that $B_\eps$ enjoys some strong compactness with respect to the time variable. We then deduce, up to extraction of subsequences, that
		\begin{equation*}
			B_\eps \left(g_\eps^\pm \gamma_\eps^\lambda\right) *_{x,v}\chi_{a}
			\rightharpoonup
			B g^{\pm,\lambda}*_{x,v}\chi_{a},
		\end{equation*}
		and, incidentally, by the uniformity of the approximation of $g_\eps^\pm \gamma_\eps^\lambda$ by $\left(g_\eps^\pm \gamma_\eps^\lambda\right) *_{x,v}\chi_{a}$ in $L_\mathrm{loc}^2(dtdxdv)$, that
		\begin{equation*}
			B_\eps g_\eps^\pm \gamma_\eps^\lambda
			\rightharpoonup
			B g^{\pm,\lambda},
		\end{equation*}
		in $L^1_\mathrm{loc}(dtdxdv)$.

		We may now take weak limits in \eqref{renormalized3} to infer that, for any given $2\leq r\leq \infty$,
		\begin{equation*}
			\begin{aligned}
				( v \cdot \nabla_x \pm & \left(v\wedge B\right)\cdot \nabla_v )
				g^{\pm,\lambda}
				\mp E \cdot v \\
				& =
				\int_{\mathbb{R}^3\times\mathbb{S}^2} \left( q^\pm + q^{\pm,\mp}\right) M_*dv_*d\sigma
				+ O\left(\lambda^\frac{1}{r}\right)_{L^\frac{2r}{2+r}_\mathrm{loc}(dtdxdv)}.
			\end{aligned}
		\end{equation*}
		Finally, in view of \eqref{endpoint claim}, letting $\lambda\rightarrow 0$, we arrive at
		\begin{equation}\label{kinetic equation 3}
			( v \cdot \nabla_x \pm \left(v\wedge B\right)\cdot \nabla_v )
			g^\pm
			\mp E \cdot v
			=
			\int_{\mathbb{R}^3\times\mathbb{S}^2} \left(q^\pm + q^{\pm,\mp}\right) M_*dv_*d\sigma,
		\end{equation}
		which, together with the fact, according to Lemmas \ref{relaxation-control} and \ref{relaxation2-control}, that $g^+$ and $g^-$ are infinitesimal Maxwellians, which differ only by their densities $\rho^+$ and $\rho^-$, provides that
		\begin{equation*}
			\begin{aligned}
				\int_{\mathbb{R}^3\times\mathbb{S}^2} & \left(q^\pm + q^{\pm,\mp}\right) M_*dv_*d\sigma
				\\
				& = \Div\left((\rho^\pm+\theta)v+\frac{|v|^2}{3}u+\phi u+\theta\psi\right)
				\mp \left(E+u\wedge B\right)\cdot v
				\\
				& = \left( \phi:\nabla_x u +\psi\cdot \nabla_x \theta \right)
				+ \left(\nabla_x(\rho^\pm+\theta)\mp \left(E+u\wedge B\right)\right)\cdot v
				+\frac13 \left(\Div u\right)|v|^2.
			 \end{aligned}
		 \end{equation*}
		Equivalently, we find that
		\begin{equation*}
			\begin{aligned}
				\frac 12\int_{\mathbb{R}^3\times\mathbb{S}^2} & \left(q^++q^-+q^{+,-}+q^{-,+}\right) M_*dv_*d\sigma
				\\
				& = \left( \phi:\nabla_x u +\psi\cdot \nabla_x \theta \right)
				+ \nabla_x\left(\frac{\rho^++\rho^-}{2}+\theta\right)\cdot v
				+\frac13 \left(\Div u\right)|v|^2,
			 \end{aligned}
		\end{equation*}
		and
		\begin{equation*}
			\begin{aligned}
				\frac 12 \int_{\mathbb{R}^3\times\mathbb{S}^2} \left(q^+-q^-+q^{+,-}-q^{-,+}\right) & M_*dv_*d\sigma
				\\
				& = \left(\frac 12 \nabla_x(\rho^+ - \rho^-) - \left(E+u\wedge B\right)\right)\cdot v.
			 \end{aligned}
		 \end{equation*}
		
		Then, remarking that $q^\pm$ and $q^{\pm,\mp}$ inherit the collisional symmetries of $q_\eps^\pm$, $q_\eps^{\pm,\mp}$, $\hat q_\eps^\pm$ and $\hat q_\eps^{\pm,\mp}$, we get
		\begin{equation*}
			\int_{\mathbb{R}^3\times\mathbb{R}^3\times\mathbb{S}^2}
			\left(q^++q^-+q^{+,-}+q^{-,+}\right)
			\begin{pmatrix}
				1 \\ v \\ \frac{|v|^2}{2}
			\end{pmatrix}
			MM_* dvdv_*d\sigma = 0,
		\end{equation*}
		so that, since $\phi(v)$ and $\psi(v)$ are orthogonal to the collisional invariants, the constraints \eqref{constraints2-1} hold.

	The proof of the proposition is now complete.
\end{proof}

The next proposition further characterizes the limiting collision integrands.

\begin{prop}\label{high weak-comp3}
	Let $\left(f_\eps^\pm, E_\eps, B_\eps\right)$ be the sequence of renormalized solutions to the scaled two species Vlasov-Maxwell-Boltzmann system \eqref{VMB2} considered in Theorem \ref{CV-OMHDSTRONG} for strong interspecies interactions, i.e.\ $\delta=1$. In accordance with Lemmas \ref{L1-lem}, \ref{L2-lem}, \ref{L2-qlem}, \ref{bound hjw}, \ref{weak compactness h} and \ref{strong n} denote by
	\begin{equation*}
		\begin{aligned}
			g^\pm & \in L^\infty\left(dt;L^2\left(Mdxdv\right)\right),\\
			q^{\pm,\mp} & \in L^2\left(MM_*dtdxdvdv_*d\sigma\right),\\
			h & \in L^1_\mathrm{loc}\left(dtdx;L^1\left((1+|v|^2)Mdv\right)\right),
		\end{aligned}
	\end{equation*}
	any joint limit points of the families $\hat g_\eps^\pm$, $\hat q_\eps^{\pm,\mp}$ and $h_\eps$ defined by \eqref{hatg}, \eqref{hatq-def} and \eqref{def h}, respectively.

	Then, one has
	\begin{equation}\label{mixed q phi psi 2}
		\int_{\mathbb{R}^3\times\mathbb{S}^2} \left(q^{+}-q^{-}\right) M_* dv_*d\sigma
		=
		-\mathcal{L} \left(h\right),
	\end{equation}
	\begin{equation}\label{mixed q phi psi 3}
		\int_{\mathbb{R}^3\times\mathbb{S}^2} \left(q^{+,-}-q^{-,+}\right) M_* dv_*d\sigma
		=
		n  u\cdot \mathfrak{L}(v)
		+
		n  \theta \mathfrak{L}\left(\frac{|v|^2}{2}\right)
		-\mathfrak{L} \left(h\right),
	\end{equation}
	and
	\begin{equation}\label{null dissipation}
		q^{+}+q^{-}-q^{+,-}-q^{-,+} =0,
	\end{equation}
	where $n=\rho^+-\rho^-$, $u$ and $\theta $ are, respectively, the charge density, bulk velocity and temperature associated with the limiting fluctuations $g^\pm$.
\end{prop}

\begin{proof}
	We start from the decomposition
	\begin{equation}\label{asymptotic sigma strong 2}
		\begin{aligned}
			h_\eps
			& = \hat h_\eps + \frac{1}{4}\left[\left|\hat g^+_\eps\right|^2-\left|\hat g^-_\eps\right|^2
			- \int_{\mathbb{R}^3}\left(\left|\hat g^+_\eps\right|^2-\left|\hat g^-_\eps\right|^2\right)Mdv\right]
			\\
			& = \hat h_\eps + \frac{1}{4}\left[\left(\hat g^+_\eps - \hat g^-_\eps - \hat n_\eps \right)\left(\hat g^+_\eps + \hat g^-_\eps\right)
			+\hat n_\eps \left(\hat g^+_\eps + \hat g^-_\eps\right)\right]
			\\
			& - \frac 14 \left[\int_{\mathbb{R}^3}\left(\hat g^+_\eps - \hat g^-_\eps - \hat n_\eps \right)\left(\hat g^+_\eps + \hat g^-_\eps\right)Mdv + \hat n_\eps\left(\hat \rho_\eps^++\hat \rho_\eps^-\right)\right],
		\end{aligned}
	\end{equation}
	which follows from the decomposition \eqref{fluct-decomposition} of fluctuations.
	% which, when integrated against $vMdv$ and $\left(\frac{|v|^2}{3}-1\right)Mdv$, yields
	% \begin{equation*}
	% 	\begin{aligned}
	% 		j_\eps & = \hat j_\eps + \frac{1}{4}\left[\int_{\mathbb{R}^3}\left(\hat g^+_\eps - \hat g^-_\eps - \hat n_\eps \right)\left(\hat g^+_\eps + \hat g^-_\eps\right)vMdv
	% 		+\hat n_\eps \left(\hat u^+_\eps + \hat u^-_\eps\right)\right],\\
	% 		w_\eps & = \hat w_\eps + \frac{1}{4}\left[\int_{\mathbb{R}^3}\left(\hat g^+_\eps - \hat g^-_\eps - \hat n_\eps \right)\left(\hat g^+_\eps + \hat g^-_\eps\right)\left(\frac{|v|^2}{3}-1\right)Mdv
	% 		+\hat n_\eps \left(\hat \theta^+_\eps + \hat \theta^-_\eps\right)\right].
	% 	\end{aligned}
	% \end{equation*}
	In order to apply the compactness results from Lemmas \ref{weak compactness h} and \ref{strong n}, we consider the following renormalization of the above decomposition~:
	\begin{equation}\label{renormalized decomposition}
		\begin{aligned}
			\frac{h_\eps}{R_\eps}
			= \frac{\hat h_\eps}{R_\eps}
			& + \frac{1}{4}\left[\eps\frac{\hat h_\eps}{R_\eps}
			\left(\hat g^+_\eps + \hat g^-_\eps\right)
			+\frac{\hat n_\eps}{R_\eps}
			\left(\hat g^+_\eps + \hat g^-_\eps\right)\right]
			\\
			& - \frac 14 \left[\eps\int_{\mathbb{R}^3}\frac{\hat h_\eps}{R_\eps}
			\left(\hat g^+_\eps + \hat g^-_\eps\right)Mdv + \frac{\hat n_\eps}{R_\eps}
			\left(\hat \rho_\eps^++\hat \rho_\eps^-\right)\right],
		\end{aligned}
	\end{equation}
	where we have written $R_\eps={1+\left\| \hat g_\eps^+-\hat g_\eps^- \right\|_{L^2(Mdv)}}$, for convenience.

	Then, according to Lemmas \ref{weak compactness h} and \ref{strong n}, we have now weakly convergent subsequences
	\begin{equation*}
		\frac{h_\eps}{R_\eps}\rightharpoonup \frac{h}{1+|n|}
		\qquad\text{in }
		\textit{w-}L^1_\mathrm{loc}\left(dtdx;\textit{w-}L^1\left(Mdv\right)\right),
	\end{equation*}
	where $h\in L^1_{\mathrm{loc}}\left(dtdx;L^1\left(\left(1+|v|^2\right)Mdv\right)\right)$, and
	\begin{equation*}
		\frac{\hat h_\eps}{R_\eps}\rightharpoonup \frac{\hat h}{1+|n|}
		\qquad\text{in }
		\textit{w-}L^2_\mathrm{loc}\left(dtdx;\textit{w-}L^2\left(Mdv\right)\right),
	\end{equation*}
	where $\hat h\in L^1_{\mathrm{loc}}\left(dtdx;L^2\left(Mdv\right)\right)$. Hence, taking weak limits in \eqref{renormalized decomposition}, we find, further utilizing the strong convergence $\frac{\hat n_\eps}{R_\eps}\rightarrow \frac{n}{1+|n|}$ in $L^2_{\mathrm{loc}}\left(dtdx\right)$ from Lemma \ref{strong n}, that
	\begin{equation}\label{h hat h}
		h = \hat h + \frac{1}{2} n \left(u\cdot v+\theta\left(\frac{|v|^2}{2}-\frac 32\right)\right).
	\end{equation}
	% and
	% \begin{equation*}
	% 	j = \hat j + \frac{1}{2} nu.
	% \end{equation*}

	Next, it is readily seen that the elementary decompositions
	\begin{equation*}
		\begin{aligned}
			\cL \left(\hat g_\eps^\pm\right) & =
			\frac\eps 2  \cQ\left(\hat g_\eps^\pm ,\hat g_\eps^\pm \right)-\frac 2{\eps}  \cQ \left(\sqrt{G_\eps^\pm},\sqrt{G_\eps^\pm}\right),
			\\
			\cL \left(\hat g_\eps^\pm,\hat g_\eps^\mp\right) & =
			\frac\eps 2  \cQ\left(\hat g_\eps^\pm ,\hat g_\eps^\mp \right)-\frac 2{\eps}  \cQ \left(\sqrt{G_\eps^\pm},\sqrt{G_\eps^\mp}\right),
		\end{aligned}
	\end{equation*}
	yield that
	\begin{equation}\label{asymptotic sigma strong 3}
		\begin{aligned}
			\mathcal{L} \left(\hat h_\eps\right) & =
			\frac 1 2  \left[\cQ\left(\hat g_\eps^+ ,\hat g_\eps^+ \right) - \cQ\left(\hat g_\eps^- ,\hat g_\eps^-\right)\right]
			\\
			& -\frac {2}{\eps^2} \left[\cQ \left(\sqrt{G_\eps^+},\sqrt{G_\eps^+}\right)
			- \cQ \left(\sqrt{G_\eps^-},\sqrt{G_\eps^-}\right)\right]
			\\
			% & =
			% 		\frac 1 2  \left[\cQ\left(\hat g_\eps^+ ,\hat g_\eps^+ \right) - \cQ\left(\hat g_\eps^- ,\hat g_\eps^-\right)\right]
			% 		-
			% 		\int_{\mathbb{R}^3\times\mathbb{S}^2} \left(\hat q_\eps^{+}-\hat q_\eps^{-}\right) M_* dv_*d\sigma
			% 		\\
			& =
			\frac 1 2  \left[\cQ\left(\hat g_\eps^+-\hat g_\eps^- - \hat n_\eps ,\hat g_\eps^+ \right) + \cQ\left(\hat g_\eps^- ,\hat g_\eps^+-\hat g_\eps^- - \hat n_\eps \right)\right]
			\\
			& +
			\frac 1 2 \hat n_\eps  \left[\cQ\left(1 ,\hat g_\eps^+ \right) + \cQ\left(\hat g_\eps^- ,1 \right)\right]
			-\int_{\mathbb{R}^3\times\mathbb{S}^2} \left(\hat q_\eps^{+}-\hat q_\eps^{-}\right) M_* dv_*d\sigma,
		\end{aligned}
	\end{equation}
	and
	\begin{equation}\label{asymptotic sigma strong 4}
		\begin{aligned}
			\mathfrak{L} \left(\hat h_\eps\right) & =
			\frac 1 2  \left[\cQ\left(\hat g_\eps^+ ,\hat g_\eps^- \right) - \cQ\left(\hat g_\eps^- ,\hat g_\eps^+\right)\right]
			\\
			& -\frac {2}{\eps^2} \left[\cQ \left(\sqrt{G_\eps^+},\sqrt{G_\eps^-}\right)
			- \cQ \left(\sqrt{G_\eps^-},\sqrt{G_\eps^+}\right)\right]
			\\
			% & =
			% 		\frac 1 2  \left[\cQ\left(\hat g_\eps^+ ,\hat g_\eps^- \right) - \cQ\left(\hat g_\eps^- ,\hat g_\eps^+\right)\right]
			% 		-
			% 		\int_{\mathbb{R}^3\times\mathbb{S}^2} \left(\hat q_\eps^{+,-}-\hat q_\eps^{-,+}\right) M_* dv_*d\sigma
			% 		\\
			& =
			\frac 1 2  \left[\cQ\left(\hat g_\eps^+-\hat g_\eps^- - \hat n_\eps ,\hat g_\eps^- \right) - \cQ\left(\hat g_\eps^- ,\hat g_\eps^+-\hat g_\eps^- - \hat n_\eps \right)\right]
			\\
			& +
			\frac 1 2 \hat n_\eps  \left[\cQ\left(1 ,\hat g_\eps^- \right) - \cQ\left(\hat g_\eps^- ,1 \right)\right]
			-\int_{\mathbb{R}^3\times\mathbb{S}^2} \left(\hat q_\eps^{+,-}-\hat q_\eps^{-,+}\right) M_* dv_*d\sigma.
		\end{aligned}
	\end{equation}
	We also have the simple decomposition
	\begin{equation*}% \label{h decomp 4}
		\begin{aligned}
			% \frac\eps 2 \cQ\left(\hat g_\eps^+-\hat g_\eps^- ,\hat h_\eps \right)
			% 		+ \frac\eps 2 \cQ\left(\hat h_\eps ,\hat n_\eps \right)
			\hat q_\eps^{+}+\hat q_\eps^{-}-\hat q_\eps^{+,-}-\hat q_\eps^{-,+} & =
			\frac 1 2 \left(\hat g_\eps^+-\hat g_\eps^-\right)'\left(\hat g_\eps^+-\hat g_\eps^- \right)_*'
			-\frac 1 2 \left(\hat g_\eps^+-\hat g_\eps^-\right) \left(\hat g_\eps^+-\hat g_\eps^- \right)_*
			\\
			& =
			\frac\eps 2 \left(\hat g_\eps^+-\hat g_\eps^-\right)' \hat h_{\eps *}'
			-\frac\eps 2 \left(\hat g_\eps^+-\hat g_\eps^-\right) \hat h_{\eps *}
			+ \frac\eps 2 \hat h_{\eps}' \hat n_{\eps}
			-\frac\eps 2 \hat h_{\eps} \hat n_{\eps}.
		\end{aligned}
	\end{equation*}
	
	% \begin{equation*}
	% 	\begin{aligned}
	% 		\frac 1 2 \cQ\left(\hat g_\eps^+-\hat g_\eps^- ,\hat g_\eps^+-\hat g_\eps^- -\hat n_\eps \right)
	% 		& + \frac 1 2 \cQ\left(\hat g_\eps^+-\hat g_\eps^- -\hat n_\eps ,\hat n_\eps \right)
	% 		\\
	% 		& =
	% 		\int_{\mathbb{R}^3\times\mathbb{S}^2} \left(\hat q_\eps^{+}+\hat q_\eps^{-}-\hat q_\eps^{+,-}-\hat q_\eps^{-,+}\right) M_* dv_*d\sigma,
	% 	\end{aligned}
	% \end{equation*}
	As previously, we renormalize the above identities into
	\begin{equation*}
		\begin{aligned}
			\mathcal{L} \left(\frac{\hat h_\eps}{R_\eps}\right) & =
			\frac \eps 2  \left[\cQ\left(\frac{\hat h_\eps}{R_\eps} ,\hat g_\eps^+ \right) + \cQ\left(\hat g_\eps^- ,\frac{\hat h_\eps}{R_\eps}\right)\right]
			\\
			& +
			\frac{\hat n_\eps}{2 R_\eps}  \left[\cQ\left(1 ,\hat g_\eps^+ \right) + \cQ\left(\hat g_\eps^- ,1 \right)\right]
			-\frac{1}{R_\eps}\int_{\mathbb{R}^3\times\mathbb{S}^2} \left(\hat q_\eps^{+}-\hat q_\eps^{-}\right) M_* dv_*d\sigma,
			\\
			\mathfrak{L} \left(\frac{\hat h_\eps}{R_\eps}\right) & =
			\frac \eps 2  \left[\cQ\left(\frac{\hat h_\eps}{R_\eps} ,\hat g_\eps^- \right) - \cQ\left(\hat g_\eps^- ,\frac{\hat h_\eps}{R_\eps}\right)\right]
			\\
			& +
			\frac{\hat n_\eps}{2R_\eps}  \left[\cQ\left(1 ,\hat g_\eps^- \right) - \cQ\left(\hat g_\eps^- ,1 \right)\right]
			-\frac{1}{R_\eps}\int_{\mathbb{R}^3\times\mathbb{S}^2} \left(\hat q_\eps^{+,-}-\hat q_\eps^{-,+}\right) M_* dv_*d\sigma,
		\end{aligned}
	\end{equation*}
	and
	\begin{equation*}
		\begin{aligned}
			\frac{1}{R_\eps}\left(\hat q_\eps^{+}+\hat q_\eps^{-}-\hat q_\eps^{+,-}-\hat q_\eps^{-,+}\right)
			& =
			\frac\eps 2 \left(\hat g_\eps^+-\hat g_\eps^-\right)' \left(\frac{\hat h_{\eps}}{R_\eps}\right)_*'
			-\frac\eps 2 \left(\hat g_\eps^+-\hat g_\eps^-\right) \left(\frac{\hat h_{\eps}}{R_\eps}\right)_*
			\\
			& + \frac\eps 2 \left(\frac{\hat h_{\eps}}{R_\eps}\right)' \hat n_{\eps}
			-\frac\eps 2 \frac{\hat h_{\eps}}{R_\eps} \hat n_{\eps}.
		\end{aligned}
	\end{equation*}
	% \begin{equation*}
	% 	\begin{aligned}
	% 		\frac\eps 2 \cQ\left(\hat g_\eps^+-\hat g_\eps^- ,\frac{\hat h_\eps}{R(t,x)} \right)
	% 		& + \frac\eps 2 \cQ\left(\frac{\hat h_\eps}{R(t,x)} ,\hat n_\eps \right)
	% 		\\
	% 		& =\frac{1}{R(t,x)}
	% 		\int_{\mathbb{R}^3\times\mathbb{S}^2} \left(\hat q_\eps^{+}+\hat q_\eps^{-}-\hat q_\eps^{+,-}-\hat q_\eps^{-,+}\right) M_* dv_*d\sigma.
	% 	\end{aligned}
	% \end{equation*}
	
	Finally, taking weak limits, we find, utilizing again the strong convergence from Lemma \ref{strong n}, that
	\begin{equation*}
		\begin{aligned}
			\mathcal{L} \left(\hat h\right)
			& =
			-\int_{\mathbb{R}^3\times\mathbb{S}^2} \left(q^{+}-q^{-}\right) M_* dv_*d\sigma,
			\\
			\mathfrak{L} \left(\hat h\right) & =
			\frac 1 2 n  u\cdot\left[\cQ\left(1 ,v \right) - \cQ\left(v ,1 \right)\right]
			+
			\frac 1 4 n  \theta \left[\cQ\left(1 ,|v|^2 \right) - \cQ\left(|v|^2 ,1 \right)\right]
			\\
			& -\int_{\mathbb{R}^3\times\mathbb{S}^2} \left(q^{+,-}-q^{-,+}\right) M_* dv_*d\sigma
			\\
			& =
			\frac 1 2 n  u\cdot \mathfrak{L}(v)
			+
			\frac 1 4 n  \theta \mathfrak{L}\left(|v|^2\right)
			-\int_{\mathbb{R}^3\times\mathbb{S}^2} \left(q^{+,-}-q^{-,+}\right) M_* dv_*d\sigma,
		\end{aligned}
	\end{equation*}
	and
	\begin{equation*}
		q^{+}+q^{-}-q^{+,-}-q^{-,+}=0,
	\end{equation*}
	which yields, in view of \eqref{h hat h}, that
	\begin{equation*}
			\mathcal{L} \left( h\right)
			=
			-\int_{\mathbb{R}^3\times\mathbb{S}^2} \left(q^{+}-q^{-}\right) M_* dv_*d\sigma,
	\end{equation*}
	and
	\begin{equation*}
		\begin{aligned}
			\mathfrak{L} \left(h\right)
			=
			n  u\cdot \mathfrak{L}(v)
			+
			\frac 1 2 n  \theta \mathfrak{L}\left(|v|^2\right)
			-\int_{\mathbb{R}^3\times\mathbb{S}^2} \left(q^{+,-}-q^{-,+}\right) M_* dv_*d\sigma,
		\end{aligned}
	\end{equation*}
	and concludes the proof of the proposition.
\end{proof}

As a direct consequence of the previous propositions, we derive in the next result Ohm's law and the internal electric energy constraint from \eqref{TFINSFMO 2}.

\begin{prop}\label{strongOhm}
	Let $\left(f_\eps^\pm, E_\eps, B_\eps\right)$ be the sequence of renormalized solutions to the scaled two species Vlasov-Maxwell-Boltzmann system \eqref{VMB2} considered in Theorem \ref{CV-OMHDSTRONG} for strong interspecies interactions, i.e.\ $\delta=1$. In accordance with Lemmas \ref{L1-lem}, \ref{L2-lem}, \ref{bound hjw}, \ref{weak compactness h} and \ref{strong n} denote by
		\begin{equation*}
			\begin{gathered}
				g^\pm \in L^\infty\left(dt;L^2\left(Mdxdv\right)\right),
				\qquad
				h \in L^1_\mathrm{loc}\left(dtdx;L^1\left((1+|v|^2)Mdv\right)\right),
				\\
				\text{and}\qquad
				E,B\in L^\infty\left(dt;L^2\left(dx\right)\right),
			\end{gathered}
		\end{equation*}
	any joint limit points of the families $\hat g_\eps^\pm$ and $h_\eps$ defined by \eqref{hatg} and \eqref{def h}, $E_\eps$ and $B_\eps$, respectively.
	
	Then, one has
	\begin{equation*}
		j-nu = \sigma\left(-\frac 12 \nabla_x \left(\rho^+-\rho^-\right) + E + u\wedge B \right)
		\qquad\text{and}\qquad
		w = n\theta,
	\end{equation*}
	where $\rho^\pm$, $u$ and $\theta $ are, respectively, the densities, bulk velocity and temperature associated with the limiting fluctuations $g^\pm$, $j$ and $w$ are, respectively, the electric current and the internal electric energy associated with the limiting fluctuation $h$ and the electric conductivity $\sigma>0$ is defined by \eqref{sigma 3}.
\end{prop}

\begin{proof}
	By Proposition \ref{high weak-comp2}, we have that
	\begin{equation*}
		\begin{aligned}
			\int_{\mathbb{R}^3\times\mathbb{R}^3\times\mathbb{S}^2} \left(q^+-q^-+q^{+,-}-q^{-,+}\right) & \tilde\Phi MM_*dvdv_*d\sigma
			\\
			& = \sigma\left(\frac 12 \nabla_x(\rho^+ - \rho^-) - \left(E+u\wedge B\right)\right),
			\\
			\int_{\mathbb{R}^3\times\mathbb{R}^3\times\mathbb{S}^2} \left(q^+-q^-+q^{+,-}-q^{-,+}\right) & \tilde\Psi MM_*dvdv_*d\sigma
			\\
			& = 0,
		 \end{aligned}
	\end{equation*}
	where we have used the identity \eqref{sigma 2} and $\tilde\Phi$ and $\tilde\Psi$ are defined by \eqref{phi-psi-def inverses two species}. Then, further incorporating identities \eqref{mixed q phi psi 2} and \eqref{mixed q phi psi 3} from Proposition \ref{high weak-comp3} into the above relations yields that
	\begin{equation*}
		\begin{aligned}
			\int_{\mathbb{R}^3} \left(n  u\cdot \mathfrak{L}(v) + n \theta \mathfrak{L}\left(\frac{|v|^2}{2}\right) - \left(\mathcal{L} + \mathfrak{L}\right)(h)\right) & \tilde\Phi Mdv
			\\
			= & \sigma\left(\frac 12 \nabla_x(\rho^+ - \rho^-) - \left(E+u\wedge B\right)\right),
			\\
			\int_{\mathbb{R}^3} \left(n  u\cdot \mathfrak{L}(v) + n \theta \mathfrak{L}\left(\frac{|v|^2}{2}\right) - \left(\mathcal{L} + \mathfrak{L}\right)(h)\right) & \tilde\Psi Mdv
			\\
			= & 0.
		 \end{aligned}
	\end{equation*}
	Finally, using \eqref{phi-psi-def inverses two species} and the self-adjointness of $\mathcal{L}+\mathfrak{L}$, we deduce that
	\begin{equation*}
		\begin{aligned}
			nu-j
			& = \int_{\mathbb{R}^3} \left(nu\cdot\Phi + n\theta\Psi - h\right) \Phi Mdv
			\\
			& = \int_{\mathbb{R}^3} \left(\mathcal{L} + \mathfrak{L}\right)\left(nu\cdot\Phi + n\theta\Psi - h\right) \tilde\Phi Mdv
			\\
			& = \int_{\mathbb{R}^3} \left(n  u\cdot \mathfrak{L}(v) + n \theta \mathfrak{L}\left(\frac{|v|^2}{2}\right) - \left(\mathcal{L} + \mathfrak{L}\right)(h)\right) \tilde\Phi Mdv
			\\
			& = \sigma\left(\frac 12 \nabla_x(\rho^+ - \rho^-) - \left(E+u\wedge B\right)\right),
			\\
			\frac 32\left(n\theta-w\right)
			& = \int_{\mathbb{R}^3} \left(nu\cdot\Phi + n\theta\Psi - h\right) \Psi Mdv
			\\
			& = \int_{\mathbb{R}^3} \left(\mathcal{L} + \mathfrak{L}\right)\left(nu\cdot\Phi + n\theta\Psi - h\right) \tilde\Psi Mdv
			\\
			& = \int_{\mathbb{R}^3} \left(n  u\cdot \mathfrak{L}(v) + n \theta \mathfrak{L}\left(\frac{|v|^2}{2}\right) - \left(\mathcal{L} + \mathfrak{L}\right)(h)\right) \tilde\Psi Mdv
			\\
			& = 0,
		 \end{aligned}
	\end{equation*}
	which concludes the proof of the proposition.
\end{proof}

\section{Energy inequalities}\label{energy inequality singular}

In view of the results from Section \ref{high constraints 1}, we are now able to establish the limiting energy inequality for two species in the case of strong interactions.

\begin{prop}\label{energy inequality strong interactions}
	Let $\left(f_\eps^\pm, E_\eps, B_\eps\right)$ be the sequence of renormalized solutions to the scaled two species Vlasov-Maxwell-Boltzmann system \eqref{VMB2} considered in Theorem \ref{CV-OMHDSTRONG} for strong interspecies interactions, i.e.\ $\delta=1$. In accordance with Lemmas \ref{L1-lem}, \ref{L2-lem}, \ref{L2-qlem} and \ref{strong n}, denote by
	\begin{equation*}
		\begin{gathered}
			g^\pm\in L^\infty\left(dt;L^2\left(Mdxdv\right)\right),
			\qquad
			h \in L^1_{\mathrm{loc}}\left(dtdx;L^1\left(\left(1+|v|^2\right)Mdv\right)\right)
			\\
			\text{and}\qquad
			q^\pm, q^{\pm,\mp} \in L^2\left(MM_*dtdxdvdv_*d\sigma\right)
		\end{gathered}
	\end{equation*}
	any joint limit points of the families $\hat g_\eps^\pm$,  $h_\eps$, $\hat q_\eps^{\pm}$ and $\hat q_\eps^{\pm,\mp}$ defined by \eqref{hatg}, \eqref{def h} and \eqref{hatq-def}, respectively, and by
	\begin{equation*}
		E,B\in L^\infty\left(dt;L^2\left(dx\right)\right)
	\end{equation*}
	any joint limit points of the families $E_\eps$ and $B_\eps$, respectively.
	
	Then, one has the energy inequality, for almost every $t\geq 0$,
	\begin{equation*}
		\begin{aligned}
			& \frac 12\left(\frac 12\left\|n\right\|_{L^2_x}^2+2\left\|u\right\|_{L^2_x}^2
			+ 5\left\|\theta\right\|_{L^2_x}^2 + \left\|E\right\|_{L^2_x}^2
			+ \left\|B\right\|_{L^2_x}^2 \right)(t)
			\\
			& +
			\int_0^t \left(2\mu
			\left\|\nabla_x u\right\|_{L^2_x}^2
			+ 5\kappa
			\left\|\nabla_x\theta\right\|_{L^2_x}^2 +
			\frac 1\sigma\left\|j-nu\right\|_{L^2_x}^2 +
			\frac{9}{8\lambda}\left\|w-n\theta\right\|_{L^2_x}^2\right)(s) ds
			\\
			& \leq C^\mathrm{in},
		\end{aligned}
	\end{equation*}
	where $\rho^\pm$, $u$ and $\theta $ are, respectively, the densities, bulk velocity and temperature associated with the limiting fluctuations $g^\pm$ and the charge density is given by $n=\rho^+-\rho^-$, while $j$ and $w$ are, respectively, the electric current and the internal electric energy associated with the limiting fluctuation $h$, and, finally, the viscosity $\mu>0$, thermal conductivity $\kappa>0$, electric conductivity $\sigma>0$ and energy conductivity $\lambda>0$ are respectively defined by \eqref{mu kappa 2}, \eqref{sigma 3} and \eqref{lambda 2}.
\end{prop}

\begin{proof}
	First, by the estimate \eqref{q-est} from Lemma \ref{L2-qlem} and the weak sequential lower semi-continuity of convex functionals, we find that, for all $t\geq 0$,
	\begin{equation*}
		\begin{aligned}
			\frac 14 \int_0^t \int_{\mathbb{R}^3} \int_{\mathbb{R}^3\times\mathbb{R}^3\times\mathbb{S}^2} &
			\left(q^\pm\right)^2 MM_*dvdv_* d\sigma dx ds
			\\
			& \leq \liminf_{\eps\rightarrow 0}
			\frac 14 \int_0^t \int_{\mathbb{R}^3} \int_{\mathbb{R}^3\times\mathbb{R}^3\times\mathbb{S}^2}
			\left(\hat{q}_\eps^\pm\right)^2 MM_*dvdv_* d\sigma dx ds
			\\
			& \leq \liminf_{\eps\rightarrow 0} \frac{1}{\epsilon^4}\int_0^t\int_{\mathbb{R}^3}
			D\left(f_\eps^\pm\right)(s) dx ds,
		\end{aligned}
	\end{equation*}
	and
	\begin{equation*}
		\begin{aligned}
			\frac 12 \int_0^t \int_{\mathbb{R}^3} \int_{\mathbb{R}^3\times\mathbb{R}^3\times\mathbb{S}^2} &
			\left(q^{\pm,\mp}\right)^2 MM_*dvdv_* d\sigma dx ds
			\\
			& \leq \liminf_{\eps\rightarrow 0}
			\frac 12 \int_0^t \int_{\mathbb{R}^3} \int_{\mathbb{R}^3\times\mathbb{R}^3\times\mathbb{S}^2}
			\left(\hat{q}_\eps^{\pm,\mp}\right)^2 MM_*dvdv_* d\sigma dx ds
			\\
			& \leq \liminf_{\eps\rightarrow 0} \frac{1}{\epsilon^4}\int_0^t\int_{\mathbb{R}^3}
			D\left(f_\eps^+,f_\eps^-\right)(s) dx ds,
		\end{aligned}
	\end{equation*}
	which, when combined with Lemma \ref{L1-lem}, yields, passing to the limit in the entropy inequality \eqref{entropy2}, for almost every $t\geq 0$,
	\begin{equation*}
		\begin{aligned}
		\frac 12 & \int_{\mathbb{R}^3\times\mathbb{R}^3} \left(\left(g^{+}\right)^2 + \left(g^{-}\right)^2\right)(t)Mdxdv
		+ \frac 1{2} \int_{\mathbb{R}^3} \left(|E|^2+ |B|^2\right)(t) dx
		\\
		& +
		\frac 14 \int_0^t \int_{\mathbb{R}^3} \int_{\mathbb{R}^3\times\mathbb{R}^3\times\mathbb{S}^2}
		\left(\left(q^+\right)^2 + \left(q^-\right)^2 + \left(q^{+,-}\right)^2+ \left(q^{-,+}\right)^2\right) MM_*dvdv_* d\sigma dx ds
		\\
		& \leq C^\mathrm{in}.
		\end{aligned}
	\end{equation*}
	Since, according to Lemmas \ref{relaxation-control} and \ref{relaxation2-control}, the limiting fluctuations $g^\pm=\rho^\pm+u\cdot v + \theta \left(\frac{|v|^2}{2}-\frac 32\right)$ are infinitesimal Maxwellians which differ only by their densities $\rho^+$ and $\rho^-$, we easily compute that, in view of the strong Boussinesq relation $\frac{\rho^++\rho^-}{2}+\theta=0$ following from \eqref{constraints2-1},
	\begin{equation*}
		\begin{aligned}
			\frac 12 \int_{\mathbb{R}^3}\left(\left(g^{+}\right)^2 + \left(g^{-}\right)^2\right)Mdv
			& = \frac{\left(\rho^+\right)^2+\left(\rho^-\right)^2}{2} + |u|^2 + \frac 32 \theta^2 \\
			& = \frac 14 n^2 + |u|^2 + \frac 52 \theta^2,
		\end{aligned}
	\end{equation*}
	where $n=\rho^+-\rho^-$, which implies
	\begin{equation}\label{entropy3 limit}
		\begin{aligned}
			& \left(\frac 14\left\|n\right\|_{L^2_x}^2+\left\|u\right\|_{L^2_x}^2
			+ \frac 52\left\|\theta\right\|_{L^2_x}^2 + \frac 12 \left\|E\right\|_{L^2_x}^2
			+ \frac 12 \left\|B\right\|_{L^2_x}^2 \right)
			\\
			& +
			\frac 1{16} \int_0^t \int_{\mathbb{R}^3} \int_{\mathbb{R}^3\times\mathbb{R}^3\times\mathbb{S}^2}
			\left(q^++q^{+,-}+q^-+q^{-,+}\right)^2 MM_*dvdv_* d\sigma dx ds
			% \\
			% 			& +
			% 			\frac 1{16} \int_0^t \int_{\mathbb{R}^3} \int_{\mathbb{R}^3\times\mathbb{R}^3\times\mathbb{S}^2}
			% 			\left(q^+-q^{+,-}+q^--q^{-,+}\right)^2 MM_*dvdv_* d\sigma dx ds
			\\
			& +
			\frac 1{8} \int_0^t \int_{\mathbb{R}^3} \int_{\mathbb{R}^3\times\mathbb{R}^3\times\mathbb{S}^2}
			\left(\left(q^+-q^-\right)^2
			+\left(q^{+,-}-q^{-,+}\right)^2\right) MM_*dvdv_* d\sigma dx ds
			\\
			& \leq C^\mathrm{in},
		\end{aligned}
	\end{equation}
	where we have used the identity \eqref{null dissipation}.

	There only remains to evaluate the contribution of the entropy dissipation in \eqref{entropy3 limit}. To this end, applying the method of proof of Proposition \ref{energy ineq 1}, based on the Bessel inequality \eqref{bessel 2} (where $\mu$ and $\kappa$ are now defined by \eqref{mu kappa 2} with $\delta=1$ instead of \eqref{mu kappa}, which introduces a factor $2$ in \eqref{bessel 2}), with the constraints \eqref{q phi psi 3} and \eqref{constraints2-1} from Proposition \ref{high weak-comp2}, note that it holds
	\begin{equation}\label{entropy3 1 limit}
		\begin{aligned}
			\int_0^t & \left(2\mu
			\left\|\nabla_x u\right\|_{L^2_x}^2
			+ 5\kappa
			\left\|\nabla_x\theta\right\|_{L^2_x}^2\right)(s) ds \\
			& \leq \frac 1{16} \int_0^t \int_{\mathbb{R}^3} \int_{\mathbb{R}^3\times\mathbb{R}^3\times\mathbb{S}^2}
			\left(q^++q^{+,-}+q^-+q^{-,+}\right)^2 MM_*dvdv_* d\sigma dx ds.
		\end{aligned}
	\end{equation}
	Next, the remaining contributions in the entropy dissipation will be evaluated through a direct application of the following Bessel inequality~:
	\begin{equation}\label{bessel 5}
		\begin{aligned}
			\frac 8{\sigma} &
			\left|\int_{\mathbb{R}^3\times\mathbb{R}^3\times\mathbb{S}^2}
			\left(q^++q^{+,-}-q^--q^{-,+}\right)
			\tilde\Phi MM_*dvdv_*d\sigma
			\right|^2
			\\
			& + \frac 4{\lambda}
			\left(\int_{\mathbb{R}^3\times\mathbb{R}^3\times\mathbb{S}^2}
			\left(q^++q^{+,-}-q^--q^{-,+}\right)
			\tilde\Psi MM_*dvdv_*d\sigma\right)^2
			\\
			& \leq
			\int_{\mathbb{R}^3\times\mathbb{R}^3\times\mathbb{S}^2}
			\left(
			\left(q^{+}-q^-\right)^2
			+
			\left(q^{+,-}-q^{-,+}\right)^2
			\right)
			MM_*dvdv_*d\sigma,
		\end{aligned}
	\end{equation}
	where $\tilde\Phi$ and $\tilde\Psi$ are defined by \eqref{phi-psi-def inverses two species}.

	For the sake of completeness, we provide a short justification of \eqref{bessel 5} below. But prior to this, let us conclude the proof of the present proposition. To this end, we employ the identities \eqref{mixed q phi psi 2} and \eqref{mixed q phi psi 3} from Proposition \ref{high weak-comp3} in combination with the relations \eqref{phi-psi-def inverses two species} and the self-adjointness of $\mathcal{L}+\mathfrak{L}$ to deduce from the inequality \eqref{bessel 5} that
	\begin{equation}\label{entropy3 1 limit 2}
		\begin{aligned}
			\frac8\sigma
			\left|j-nu
			\right|^2
			+&  \frac9{\lambda}
			\left(w-n\theta\right)^2
			\\
			& \leq
			\int_{\mathbb{R}^3\times\mathbb{R}^3\times\mathbb{S}^2}
			\left(
			\left(q^{+}-q^-\right)^2
			+
			\left(q^{+,-}-q^{-,+}\right)^2
			\right)
			MM_*dvdv_*d\sigma.
		\end{aligned}
	\end{equation}
	Combining this with \eqref{entropy3 limit} and \eqref{entropy3 1 limit} concludes the proof of the proposition.

	Now, as announced above, we give a short proof of \eqref{bessel 5}. To this end, for any vector $A\in\mathbb{R}^{3}$ and any scalar $a\in\mathbb{R}$, one computes straightforwardly, employing the identities \eqref{sigma 2} and \eqref{lambda 2}, and the collisional symmetries, that
	\begin{equation*}
		\begin{aligned}
			& \int_{\mathbb{R}^3\times\mathbb{R}^3\times\mathbb{S}^2}
			\left|
			\begin{pmatrix}
				A\cdot \left(\tilde\Phi+\tilde\Phi_*-\tilde\Phi'-\tilde\Phi_*'\right)+
				a\left(\tilde\Psi+\tilde\Psi_*-\tilde\Psi'-\tilde\Psi_*'\right)\\
				A\cdot \left(\tilde\Phi-\tilde\Phi_*-\tilde\Phi'+\tilde\Phi_*'\right)+
				a\left(\tilde\Psi-\tilde\Psi_*-\tilde\Psi'+\tilde\Psi_*'\right)
			\end{pmatrix}
			\right|^2 MM_*dvdv_*d\sigma
			\\
			& = \int_{\mathbb{R}^3\times\mathbb{R}^3\times\mathbb{S}^2}
			\bigg(
			\left|A\cdot\left(\tilde\Phi+\tilde\Phi_*-\tilde\Phi'-\tilde\Phi_*'\right)\right|^2
			+
			\left|A\cdot\left(\tilde\Phi-\tilde\Phi_*-\tilde\Phi'+\tilde\Phi_*'\right)\right|^2
			\\
			& +
			\left|a\left(\tilde\Psi+\tilde\Psi_*-\tilde\Psi'-\tilde\Psi_*'\right)\right|^2
			+
			\left|a\left(\tilde\Psi-\tilde\Psi_*-\tilde\Psi'+\tilde\Psi_*'\right)\right|^2
			\bigg)MM_*dvdv_*d\sigma
			\\
			& = 4\left(A\otimes A\right)
			:\int_{\mathbb{R}^3}\left(\Phi\otimes\tilde\Phi\right) Mdv
			+ 4a^2
			\int_{\mathbb{R}^3}\left(\Psi \tilde\Psi\right) Mdv
			\\
			& =2 \sigma A\cdot A + 4\lambda a^2.
		\end{aligned}
	\end{equation*}
	Therefore, defining, for any $q_0,q_1\in L^2\left(MM_*dvdv_*d\sigma\right)$, the projection
	\begin{equation*}
			\begin{pmatrix}
				\bar q_0 \\ \bar q_1
			\end{pmatrix}
			=
			\begin{pmatrix}
				A_0\cdot \left(\tilde\Phi+\tilde\Phi_*-\tilde\Phi'-\tilde\Phi_*'\right)+
				a_0 \left(\tilde\Psi+\tilde\Psi_*-\tilde\Psi'-\tilde\Psi_*'\right)\\
				A_0\cdot \left(\tilde\Phi-\tilde\Phi_*-\tilde\Phi'+\tilde\Phi_*'\right)+
				a_0 \left(\tilde\Psi-\tilde\Psi_*-\tilde\Psi'+\tilde\Psi_*'\right)
			\end{pmatrix},
	\end{equation*}
	where
	\begin{equation*}
		\begin{aligned}
			A_0 & =\frac1{2\sigma}\int_{\mathbb{R}^3\times\mathbb{R}^3\times\mathbb{S}^2}
			q_0
			\left(\tilde\Phi+\tilde\Phi_*-\tilde\Phi'-\tilde\Phi_*'\right) MM_*dvdv_*d\sigma
			\\
			& +\frac 1{2\sigma}\int_{\mathbb{R}^3\times\mathbb{R}^3\times\mathbb{S}^2}
			q_1
			\left(\tilde\Phi-\tilde\Phi_*-\tilde\Phi'+\tilde\Phi_*'\right) MM_*dvdv_*d\sigma,
			\\
			a_0 & =\frac1{4\lambda}\int_{\mathbb{R}^3\times\mathbb{R}^3\times\mathbb{S}^2}
			q_0
			\left(\tilde\Psi+\tilde\Psi_*-\tilde\Psi'-\tilde\Psi_*'\right) MM_*dvdv_*d\sigma
			\\
			& +\frac1{4\lambda}\int_{\mathbb{R}^3\times\mathbb{R}^3\times\mathbb{S}^2}
			q_1
			\left(\tilde\Psi-\tilde\Psi_*-\tilde\Psi'+\tilde\Psi_*'\right) MM_*dvdv_*d\sigma,
		\end{aligned}
	\end{equation*}
	we find that
	\begin{equation*}
		\begin{aligned}
			\int_{\mathbb{R}^3\times\mathbb{R}^3\times\mathbb{S}^2}
			\begin{pmatrix}
				q_0 \\ q_1
			\end{pmatrix}
			\cdot
			\begin{pmatrix}
				\bar q_0 \\ \bar q_1
			\end{pmatrix}
			MM_*dvdv_*d\sigma
			& =
			2\sigma A_0\cdot A_0 + 4 \lambda a_0^2
			\\
			& =
			\int_{\mathbb{R}^3\times\mathbb{R}^3\times\mathbb{S}^2}
			\left|
			\begin{pmatrix}
				\bar q_0 \\ \bar q_1
			\end{pmatrix}
			\right|^2
			MM_*dvdv_*d\sigma.
		\end{aligned}
	\end{equation*}
	Hence the Bessel inequality
	\begin{equation}\label{bessel 6}
		\begin{aligned}
			2\sigma A_0\cdot A_0 + 4 \lambda a_0^2
			& =
			\int_{\mathbb{R}^3\times\mathbb{R}^3\times\mathbb{S}^2}
			\left|
			\begin{pmatrix}
				\bar q_0 \\ \bar q_1
			\end{pmatrix}
			\right|^2
			MM_*dvdv_*d\sigma
			\\
			& \leq
			\int_{\mathbb{R}^3\times\mathbb{R}^3\times\mathbb{S}^2}
			\left|
			\begin{pmatrix}
				q_0 \\ q_1
			\end{pmatrix}
			\right|^2
			MM_*dvdv_*d\sigma.
		\end{aligned}
	\end{equation}

	Therefore, setting $q_0=q^+-q^-$ and $q_1=q^{+,-}-q^{-,+}$ in \eqref{bessel 6}, we find, exploiting the collisional symmetries of $q^\pm$ and $q^{\pm,\mp}$, that
	\begin{equation*}
		\begin{aligned}
			\frac 8{\sigma} &
			\left|\int_{\mathbb{R}^3\times\mathbb{R}^3\times\mathbb{S}^2}
			\left(q^++q^{+,-}-q^--q^{-,+}\right)
			\tilde\Phi MM_*dvdv_*d\sigma
			\right|^2
			\\
			& + \frac 4{\lambda}
			\left(\int_{\mathbb{R}^3\times\mathbb{R}^3\times\mathbb{S}^2}
			\left(q^++q^{+,-}-q^--q^{-,+}\right)
			\tilde\Psi MM_*dvdv_*d\sigma\right)^2
			\\
			& \leq
			\int_{\mathbb{R}^3\times\mathbb{R}^3\times\mathbb{S}^2}
			\left(
			\left(q^{+}-q^-\right)^2
			+
			\left(q^{+,-}-q^{-,+}\right)^2
			\right)
			MM_*dvdv_*d\sigma,
		\end{aligned}
	\end{equation*}
	which concludes the justification of \eqref{bessel 5}.
\end{proof}

%% file: conservation0.tex
% \chapter{Approximate moment equations and conservation defects}
% \chapter{Approximate macroscopic equations and conservation defects}
\chapter{Approximate macroscopic equations}\label{conservation0-chap}

The most difficult part of the asymptotic analysis consists in deriving the evolution equations for the bulk velocity and temperature insofar as they involve a singular limit and nonlinear advection terms. In particular, we expect the situation to be very different according to the asymptotic regime from Theorems \ref{NS-WEAKCV}, \ref{CV-OMHD} and \ref{CV-OMHDSTRONG} under consideration. Indeed, the corresponding limiting systems, \eqref{NSFMP 2} and \eqref{TFINSFMSO 2}, respectively, do not enjoy the same stability properties~: as explained in Chapter \ref{weak stability}, the incompressible quasi-static Navier-Stokes-Fourier-Maxwell-Poisson system \eqref{NSFMP 2} is weakly stable in the energy space, which is not the case for the two-fluid incompressible Navier-Stokes-Fourier-Maxwell system with solenoidal Ohm's law \eqref{TFINSFMSO 2}.

Before focusing on this question of stability, we will first investigate the {\bf consistency of the electro-magneto-hydrodynamic approximation}. For renormalized solutions (even though we cannot prove their existence, see Section \ref{renorm sol cond}), it is not known that conservation laws are satisfied.

We therefore have to prove that approximate conservation laws hold and control their conservation defects. However, the uniform bounds established in Chapter \ref{weak bounds} are not sufficient to do so~:
\begin{itemize}
	\item In the regime of Theorem \ref{NS-WEAKCV} leading to the incompressible quasi-static Navier-Stokes-Fourier-Maxwell-Poisson system \eqref{NSFMP 2}, we will also use the nonlinear weak compactness contained in Lemma \ref{x-compactness1 0}.
	\item In the more singular regimes of Theorems \ref{CV-OMHD} and \ref{CV-OMHDSTRONG} leading to the two-fluid incompressible Navier-Stokes-Fourier-Maxwell systems with (solenoidal) Ohm's laws \eqref{TFINSFMO 2} and \eqref{TFINSFMSO 2}, we are not able to establish such an a priori nonlinear control (compare Lemma \ref{x-compactness1 0} to the weaker nonlinear compactness statement of Lemma \ref{x-compactness2 0}). The idea is therefore to use a modulated energy (or relative entropy) argument, in the same spirit as the weak-strong stability results of Chapter \ref{weak stability}.
\end{itemize}

In order to simplify the presentation, we will first detail the decompositions and convergence proof in the one species case of Theorem \ref{NS-WEAKCV}, thus enlightening the points where the equi-integrability from Lemma \ref{x-compactness1 0} is required. We will then explain how to adapt these parts of the proof to the more singular regimes of Theorems \ref{CV-OMHD} and \ref{CV-OMHDSTRONG}.

% \section{Conservation of densities, total momentum and total energy for \eqref{VMB1}}
% \section{Approximate conservation of mass, momentum and energy for \eqref{VMB1}}
\section[Approximate conservation of mass, momentum and energy\ldots]{Approximate conservation of mass, momentum and energy for one species}\label{conservation defects 1 species}

In fact, in Chapter \ref{high constraints proof}, we have already treated a similar singular limit (of lower order, though) in the regime of weak interactions for two species, which led to the derivation of the solenoidal Ohm's law and internal electric energy constraint (see Proposition \ref{solenoidalOhm}). Here, in order to deduce the limiting evolution equations for one species, we are confronted with an even more singular limit and face similar difficulties, which we briefly recall now.

We have seen in Section \ref{macro constraint} that it is possible to derive limiting kinetic equations of the type
\begin{equation*}
	v\cdot \nabla_x g -E\cdot v = \int_{\mathbb{R}^3\times\mathbb{S}^2} qM_*dv_* d\sigma,
\end{equation*}
from \eqref{VMB1} (see \eqref{kinetic equation 1} in the proof of Proposition \ref{weak-comp}). Here, we intend to take advantage of the symmetries of the collision integrand $q$ to go one order further and, thus, to derive a singular limit. Of course, since we are considering renormalized fluctuations, we do not expect that the integrals in $v$ of the right-hand side of the Vlasov-Boltzmann equation in \eqref{VMB1} against collision invariants are zero, but they should converge to zero as $\eps \to 0$ provided that we choose some appropriate renormalization which is sufficiently close to the identity. To estimate the ensuing conservation defects, we will also need to truncate large velocities. The precise construction is detailed below and will be essentially the same, later on in Section \ref{conservation defects 2 species}, for approximate conservation laws of mass, momentum and energy associated with \eqref{VMB2}.

Note that, even if conservation laws were known to hold for renormalized solutions of \eqref{VMB1} and \eqref{VMB2}, we would have to introduce similar truncations of large tails and large velocities in order to control uniformly the flux and acceleration terms.

Thus, similarly to the proof of Proposition \ref{high weak-comp}, we start from the Vlasov-Boltzmann equation from \eqref{VMB1} renormalized with the admissible nonlinearity $\Gamma(z)$ defined by
\begin{equation*}
	\Gamma(z)-1 = (z-1)\gamma(z),
\end{equation*}
where $\gamma\in C^1\left([0,\infty);\mathbb{R}\right)$ satisfies the following assumptions, for some given $C>0$~:
\begin{equation}\label{gamma-def}
	\begin{aligned}
		\gamma(z) & \equiv 1, && \hbox{for all } z \in [0, 2],\\
		\gamma(z) & \rightarrow 0, && \hbox{as } z\rightarrow\infty,\\
		\left|\gamma'(z)\right| & \leq {C \over \left(1+ z\right)^\frac{3}{2}}, && \hbox{for all } z \in [0,\infty).
	\end{aligned}
\end{equation}
Note that necessarily $\left|\gamma(z)\right| \leq {2C \over \left(1+ z\right)^\frac{1}{2}}$.

With the notation $\gamma_\eps$ for $\gamma(G_\eps)$ and $\hat \gamma_\eps$ for $\Gamma'(G_\eps)$, the scaled Vlasov-Boltzmann equation in \eqref{VMB1} renormalized relatively to the Maxwellian $M$ with the admissible nonlinearity $\Gamma(z)$ reads
\begin{equation}\label{moment1}
	\begin{aligned}
		\d_t \left(g_\eps \gamma_\eps\right) +\frac1\eps v\cdot \nabla_x \left(g_\eps \gamma_\eps\right)
		+\left(E_\eps +v\wedge B_\eps\right)\cdot \nabla_v \left(g_\eps \gamma_\eps \right)
		- & \frac1\eps E_\eps \cdot v G_\eps \hat \gamma_\eps
		\\
		& = \frac1{\eps^3} \hat\gamma_\eps
		\cQ\left(G_\eps, G_\eps\right).
	\end{aligned}
\end{equation}

We also introduce a truncation of large velocities $\chi\left( {|v|^2 \over K_\eps} \right)$, with $K_\eps =K|\log \eps|$, for some large $K>0$ to be fixed later on, and $\chi\in C_c^\infty\left([0,\infty)\right)$ a smooth compactly supported function such that $\mathds{1}_{[0,1]}\leq \chi \leq \mathds{1}_{[0,2]}$.

Thus, multiplying each side of the above equation by $\varphi(v) \chi\left( {|v|^2 \over K_\eps} \right)$, where $\varphi$ is a collision invariant, and averaging with respect to $Mdv$ leads to the approximate conservation laws
% \begin{equation*}
% 	\begin{aligned}
% 		\d_t \int_{\mathbb{R}^3} g_\eps \gamma_\eps \varphi \chi\left( {|v|^2\over K_\eps} \right) Mdv
% 		& + \frac1\eps \nabla_x  \cdot \int_{\mathbb{R}^3}
% 		g_\eps \gamma_\eps\varphi \chi\left( {|v|^2\over K_\eps}\right)v Mdv\\
% 		& = \frac1\eps E_\eps \cdot
% 		\int_{\mathbb{R}^3} (1+\eps g_\eps) \hat \gamma_\eps \varphi\chi\left( {|v|^2\over K_\eps}\right) vM dv \\
% 		& + \int_{\mathbb{R}^3}
% 		g_\eps \gamma_\eps (E_\eps+v\wedge B_\eps) \cdot \nabla_v \left( \varphi \chi\left( {|v|^2\over K_\eps}\right) M \right) dv \\
% 		& +\frac1{\eps^3} \int_{\mathbb{R}^3}
% 		\hat \gamma_\eps \cQ(G_\eps, G_\eps) \varphi \chi\left( {|v|^2\over K_\eps}\right) M dv.
% 	\end{aligned}
% \end{equation*}
\begin{equation}\label{conservation FAD}
	\d_t \int_{\mathbb{R}^3} g_\eps \gamma_\eps \varphi \chi\left( {|v|^2\over K_\eps} \right) Mdv
	+ \nabla_x\cdot F_\eps(\varphi) = A_\eps(\varphi)+ D_\eps(\varphi),
\end{equation}
with the notations
\begin{equation}\label{F-def}
	F_\eps (\varphi)=\frac1\eps \int_{\mathbb{R}^3}
	g_\eps \gamma_\eps \varphi \chi\left( {|v|^2\over K_\eps}\right) v M dv,
\end{equation}
for the fluxes, 
\begin{equation}\label{A-def}
	\begin{aligned}
		A_\eps (\varphi)
		& = \frac1\eps E_\eps \cdot
		\int_{\mathbb{R}^3} (1+\eps g_\eps) \hat \gamma_\eps \varphi\chi\left( {|v|^2\over K_\eps}\right) vM dv \\
		& + \int_{\mathbb{R}^3}
		g_\eps \gamma_\eps (E_\eps+v\wedge B_\eps) \cdot \nabla_v \left( \varphi \chi\left( {|v|^2\over K_\eps}\right) M \right) dv
		\\
		& = \frac1\eps E_\eps \cdot
		\int_{\mathbb{R}^3} \hat \gamma_\eps \varphi\chi\left( {|v|^2\over K_\eps}\right) vM dv
		+ E_\eps \cdot
		\int_{\mathbb{R}^3} g_\eps \hat \gamma_\eps \varphi\chi\left( {|v|^2\over K_\eps}\right) vM dv
		\\
		& + E_\eps\cdot\int_{\mathbb{R}^3}
		g_\eps \gamma_\eps \varphi \left( \frac{2}{K_\eps}\chi'\left( {|v|^2\over K_\eps}\right) - \chi\left( {|v|^2\over K_\eps}\right) \right) v M dv
		\\
		& + \int_{\mathbb{R}^3}
		g_\eps \gamma_\eps (E_\eps+v\wedge B_\eps) \cdot \left( \nabla_v \varphi \right) \chi\left( {|v|^2\over K_\eps}\right) M dv,
	\end{aligned}
\end{equation}
for the acceleration terms, and
\begin{equation}\label{D-def}
	% \begin{aligned}
		D_\eps (\varphi) =
		\frac1{\eps^3} \int_{\mathbb{R}^3}
		\hat \gamma_\eps \cQ(G_\eps, G_\eps) \varphi \chi\left( {|v|^2\over K_\eps}\right) M dv,
		% \\
		% & =
		% \frac1{\eps^3} \int_{\mathbb{R}^3}
		% \hat \gamma_\eps Q(f_\eps, f_\eps) \varphi \chi\left( {|v|^2\over K_\eps}\right) dv,
	% \end{aligned}
\end{equation}
for the corresponding conservation defects.

By describing the asymptotic behavior of $F_\eps (\varphi)$, $A_\eps(\varphi)$ and $D_\eps (\varphi)$, we will prove the following consistency result (compare with the formal macroscopic conservation laws \eqref{moment-eps}).

\begin{prop}\label{approx1-prop}
	Let $\left(f_\eps, E_\eps, B_\eps\right)$ be the sequence of renormalized solutions to the scaled one species Vlasov-Maxwell-Boltzmann system \eqref{VMB1} considered in Theorem \ref{NS-WEAKCV} and denote by $\tilde \rho_\eps$, $\tilde u_\eps$ and $\tilde \theta_\eps$ the density, bulk velocity and temperature associated with the renormalized fluctuations $g_\eps \gamma_\eps \chi\left( {|v|^2\over K_\eps} \right)$.
	
	Then, one has the approximate conservation laws
	\begin{equation*}%\label{approximate1}
		\begin{aligned}
			\d_t \tilde \rho_\eps + \frac1\eps \nabla_x\cdot \tilde u_\eps
			& = R_{\eps,1},
			\\
			\d_t \tilde u_\eps
			+ \nabla_x\cdot \left( \tilde u_\eps \otimes \tilde u_\eps -\frac{|\tilde u_\eps|^2}{3} \operatorname{Id} - \int_{\mathbb{R}^3\times\mathbb{R}^3\times\mathbb{S}^2} \hat q_\eps \tilde \phi MM_* dvdv_*d\sigma \right) \hspace{-50mm} &
			\\
			& = - \frac 1\eps \nabla_x\left(\tilde \rho_\eps+\tilde \theta_\eps\right)
			+ \frac1\eps E_\eps
			+\tilde \rho_\eps E_\eps
			+ \tilde u_\eps \wedge B_\eps + R_{\eps,2},
			\\
			\d_t \left(\frac 32\tilde \theta_\eps-\tilde \rho_\eps\right) + \nabla_x \cdot \left( \frac52 \tilde u_\eps \tilde \theta_\eps - \int_{\mathbb{R}^3\times\mathbb{R}^3\times\mathbb{S}^2} \hat q_\eps \tilde \psi MM_* dvdv_*d\sigma\right) \hspace{-50mm} &
			\\
			& = \tilde u_\eps \cdot E_\eps+ R_{\eps,3},
		\end{aligned}
	\end{equation*}
	where $\tilde\phi$ and $\tilde\psi$ are defined by \eqref{phi-psi-def} and \eqref{phi-psi-def inverses}, and the remainders $R_{\eps,i}$, $i=1,2,3$, converge to $0$ in $L^1_\mathrm{loc}\left(dt;W^{-1,1}_\mathrm{loc}\left(dx\right)\right)$.
\end{prop}

The proof of Proposition \ref{approx1-prop} consists in three steps respectively devoted to the study of conservation defects, fluxes and acceleration terms in \eqref{conservation FAD}. It does not present any particular difficulty and relies on refined decompositions of the different terms in the same spirit as the proof of Proposition \ref{solenoidalOhm}.

For the sake of clarity, these three steps are respectively detailed in Sections \ref{section step 1}, \ref{section step 2} and \ref{section step 3}, below.

More precisely, Proposition \ref{approx1-prop} will clearly follow from the combination of the approximate conservation laws \eqref{conservation FAD} with Lemma \ref{defect1-lem}, which handles the vanishing of conservation defects $D_\eps(\varphi)$, for any collision invariant $\varphi$, Lemma \ref{flux1}, which establishes the asymptotic behavior of the fluxes $F_\eps(v)$ and $F_\eps\left(\frac{|v|^2}{2}-\frac 52\right)$, and Lemma \ref{acceleration1}, which characterizes the acceleration terms $A_\eps(1)$, $A_\eps(v)$ and $A_\eps\left(\frac{|v|^2}{2}-\frac 52\right)$ as $\eps\to 0$.

In order to easily extend, later on in Section \ref{conservation defects 2 species}, the arguments from the present section to the case of two species, we are going to carefully keep track of and emphasize the different points where the equi-integrability property from Lemma \ref{x-compactness1 0} is used.

\subsection{Conservation defects}\label{section step 1}

The first step of the proof is to establish the vanishing of conservation defects.

\begin{lem}\label{defect1-lem}
	The conservation defects defined by \eqref{D-def} converge to zero. More precisely, for any collision invariant $\varphi$,
	\begin{equation*}
		D_\eps(\varphi) \to 0
		\text{ in }L^1_\mathrm{loc}(dtdx)
		\text{ as }\eps \to 0.
	\end{equation*}
	%
	% Moreover, we have the global bound
	% \begin{equation*}
	% 	\eps D_\eps(\varphi) =
	% 	O(1)_{L^2(dtdx)}
	% 	+ O(\eps)_{L^2\left(dt;L^1(dx)\right)}
	% 	+ O\left(\eps^2 K_\eps^\frac{3}{2}\right)_{L^1(dtdx)}.
	% \end{equation*}
\end{lem}

\begin{proof}
	Following the strategy of proof of Proposition \ref{high weak-comp}, we introduce a convenient decomposition of $D_\eps(\varphi)$, for any collision invariant $\varphi$, and then estimate the different terms using the uniform bounds from Lemmas \ref{L2-lem} and \ref{L2-qlem} (provided by the relative entropy and entropy dissipation), the relaxation estimate \eqref{relaxation-est} from Lemma \ref{relaxation-control}, as well as the equi-integrability coming from Lemma \ref{x-compactness1 0}.

	Thus, using \eqref{integrands-decomposition}, we decompose $D_\eps(\varphi)$, taking advantage of collisional symmetries~:
	\begin{equation}\label{defect-decomposition}
		\begin{aligned}
			D_\eps(\varphi) & =
			\frac{\eps}{4}
			\int_{\mathbb{R}^3\times\mathbb{R}^3\times\mathbb{S}^2}
			\hat \gamma_\eps \varphi \chi\left( {|v|^2\over K_\eps} \right)
			\hat q_\eps^2 MM_*dvdv_*d\sigma
			\\
			& - \frac{1}{\eps}\int_{\mathbb{R}^3\times\mathbb{R}^3\times\mathbb{S}^2}
			\hat \gamma_\eps \varphi \left(1-\chi\left( {|v|^2\over K_\eps} \right)\right)
			\hat q_\eps
			\sqrt{G_\eps G_{\eps *}} MM_* dvdv_*d\sigma
			\\
			& + \frac{1}{\eps}\int_{\mathbb{R}^3\times\mathbb{R}^3\times\mathbb{S}^2}
			\hat \gamma_\eps \left(1-\hat\gamma_{\eps *}\right) \varphi
			\hat q_\eps
			\sqrt{G_\eps G_{\eps *}} MM_* dvdv_*d\sigma
			\\
			& + \frac{1}{\eps}\int_{\mathbb{R}^3\times\mathbb{R}^3\times\mathbb{S}^2}
			\hat \gamma_\eps \hat\gamma_{\eps *} \left(1-\hat\gamma_{\eps}'\hat\gamma_{\eps *}'\right) \varphi
			\hat q_\eps
			\sqrt{G_\eps G_{\eps *}} MM_* dvdv_*d\sigma
			\\
			& - \frac{\eps}{4}\int_{\mathbb{R}^3\times\mathbb{R}^3\times\mathbb{S}^2}
			\hat \gamma_\eps \hat\gamma_{\eps *} \hat\gamma_{\eps}'\hat\gamma_{\eps *}'
			\varphi
			\hat q_\eps^2 MM_* dvdv_*d\sigma
			\\
			& \eqdefa D^1_\eps(\varphi)+D^2_\eps(\varphi)+D^3_\eps(\varphi)+D^4_\eps(\varphi)+D^5_\eps(\varphi),
		\end{aligned}
	\end{equation}
	where we have used that $\varphi$ is a collision invariant to symmetrize the last term.

	Now, we show that each term $D^i_\eps(\varphi)$, $i=1,\ldots,5$, vanishes separately.

	\noindent$\bullet$ The vanishing of the first term $D^1_\eps(\varphi)$, for any function $\varphi(v)$ growing at most quadratically at infinity, easily follows, using Lemma \ref{L2-qlem}, from the estimate
	\begin{equation*}
		\begin{aligned}
			\left\|D_\eps^1(\varphi)\right\|_{L^1(dtdx)} & \leq
			\frac{\eps}{4}\left\| \hat q_\eps \right\|^2_{L^2\left(MM_* dtdxdvdv_* d\sigma\right)}
			\left\| \hat \gamma_\eps \right\|_{L^\infty}
			\left\| \chi \left( {|v|^2\over K_\eps} \right)\varphi \right\|_{L^\infty} \\
			& \leq C \eps K_\eps = CK \eps|\log\eps|.
		\end{aligned}
	\end{equation*}

	\noindent$\bullet$ The second term $D^2_\eps(\varphi)$ is controlled by estimate \eqref{gaussian-decay 0} on the tails of Gaussian distributions. Using the bound from Lemma \ref{L2-qlem} and the pointwise boundedness of $\Gamma'(z)\sqrt{z}$, we get indeed, for all $\varphi(v)$ growing at most quadratically at infinity,
	\begin{equation*}
		\begin{aligned}
			\left\|D_\eps^2(\varphi)\right\|_{L^1_\mathrm{loc}(dtdx)}
			&\leq
			\frac C\eps \left\| \hat q_\eps \right\|_{L^2\left(MM_*dtdxdvdv_*d\sigma\right)}
			\left\| \hat \gamma_\eps \sqrt{G_\eps} \right\|_{L^\infty}
			\\
			&\times
			\left\| \sqrt{G_\eps} \right\|_{ L^2_\mathrm{loc}\left(dtdx;L^2\left(Mdv\right)\right)}
			\left\| \mathds{1}_{\left\{|v|^2 \geq K_\eps\right\}} \varphi \right\|_{L^2\left(Mdv\right)}\\
			& \leq C\eps^{\frac{K}{4}-1}\left|\log\eps\right|^\frac{5}{4},
		\end{aligned}
	\end{equation*}
	which tends to zero as soon as $K>4$.

	\noindent$\bullet$ The last term $D^5_\eps (\varphi)$ is mastered using the same tools.
	% For high energies, i.e.\ when $|v|^2+|v_*|^2 > K|\log \eps|$, we use a higher-dimensional variant of the estimate \eqref{gaussian-decay} on the tails of Gaussian distributions, namely
	% \begin{equation*}
	% 	\int_{\left\{|v|^2+|v_*|^2 >R\right\}} \left(|v|^2+|v_*|^2\right)^\frac{p}{2} M(v)M(v_*) dvdv_*
	% 	\sim \frac{1}{8} R^\frac{p+4}{2} e^{-\frac{R}{2}},
	% \end{equation*}
	% valid for any $p\in\mathbb{R}$, as $R\rightarrow \infty$, to obtain
	% $$
	% \begin{array}{l}
	% \ds|D_\eps^{5>}(\varphi)|\\
	% \eqdefa \ds\left| \frac1{2\eps^3}\iiint  \hat
	% \gamma_\eps  \hat \gamma_\eps' \hat \gamma_{\eps *} \hat \gamma_{\eps *}'\indc_{|v|^2+|v_*|^2 >K_\eps} (\sqrt{f'_\eps f'_{\eps *}}-\sqrt{f_\eps f_{\eps *}})^2(\varphi+\varphi_*)  dvdv_*d\sigma \right|\\
	% \leq \ds  \frac2{\eps^3}\left\| G_\eps\hat \gamma_\eps \right\|_{L^\infty}^2 \left\| \hat \gamma_\eps \right\|_{L^\infty}^2  \|\indc_{|v|^2+|v_*|^2 >K_\eps} (\varphi+\varphi_*) \|_{L^1(MM_*dvdv_*d\sigma)}
	% \end{array}
	% $$
	For high energies, i.e.\ when $|v|^2 \geq K|\log \eps|$, we obtain
	\begin{equation*}
		\begin{aligned}
			& D_\eps^{5>}(\varphi) \\
			& \eqdefa
			\frac{\eps}{4}\int_{\mathbb{R}^3\times\mathbb{R}^3\times\mathbb{S}^2}
			\hat \gamma_\eps \hat\gamma_{\eps *} \hat\gamma_{\eps}'\hat\gamma_{\eps *}'
			\varphi
			\mathds{1}_{\left\{|v|^2\geq K_\eps\right\}}
			\hat q_\eps^2 MM_* dvdv_*d\sigma
			\\
			& \leq \frac 1\eps
			\left\| \hat \gamma_\eps \sqrt{G_\eps} \right\|_{L^\infty}^2
			\left\| \hat \gamma_\eps \right\|_{L^\infty}^2
			\left\| \hat q_\eps \right\|_{L^2\left(MM_*dvdv_*d\sigma\right)}
			\left\| \varphi \mathds{1}_{\left\{|v|^2 \geq K_\eps\right\}} \right\|_{L^2\left(MM_*dvdv_*d\sigma\right)},
		\end{aligned}
	\end{equation*}
	so that, using the estimate \eqref{gaussian-decay 0} on the tails of Gaussian distributions and the bound on $\hat q_\eps$ from Lemma \ref{L2-qlem},
	\begin{equation*}
		D_\eps^{5>}(\varphi)=O\left(\eps^{\frac{K}{4}-1}\left|\log\eps\right|^\frac{5}{4}\right)_{L^2\left(dtdx\right)},
	\end{equation*}
	which tends to zero as soon as $K>4$.
	
	For moderate energies, i.e.\ when $|v|^2 <  K|\log \eps|$, we easily find
	\begin{equation*}
		\begin{aligned}
			D_\eps^{5<}(\varphi)
			& \eqdefa
			\frac{\eps}{4}\int_{\mathbb{R}^3\times\mathbb{R}^3\times\mathbb{S}^2}
			\hat \gamma_\eps \hat\gamma_{\eps *} \hat\gamma_{\eps}'\hat\gamma_{\eps *}'
			\varphi
			\mathds{1}_{\left\{|v|^2 < K_\eps\right\}}
			\hat q_\eps^2 MM_* dvdv_*d\sigma \\
			& \leq CK\eps\left|\log\eps\right|\left\|\hat q_\eps\right\|_{L^2\left(MM_*dvdv_*d\sigma\right)}^2,
		\end{aligned}
	\end{equation*}
	so that the entropy dissipation bound from Lemma \ref{L2-qlem} provides
	\begin{equation*}
		D_\eps^{5<}(\varphi) = O\left(\eps\left|\log\eps\right|\right)_{L^1\left(dtdx\right)}.
	\end{equation*}

	\noindent$\bullet$ The handling of $D_\eps^3(\varphi)$ requires the  equi-integrability coming from Lemma \ref{x-compactness1 0}. First, one has, by the Cauchy-Schwarz inequality,
	\begin{equation*}
		\begin{aligned}
			& \left|D_\eps^3(\varphi)\right|
			\\
			& \leq
			\frac 1\eps \left\| \hat q_\eps \right\|_{L^2\left(MM_*dvdv_*d\sigma\right)}
			\left\| \hat \gamma_\eps \left(1-\hat\gamma_{\eps *}\right) \varphi
			\sqrt{G_\eps G_{\eps *}}\right\|_{L^2\left(MM_*dvdv_*d\sigma\right)}
			\\
			& \leq C
			\left\| \hat q_\eps \right\|_{L^2\left(MM_*dvdv_*d\sigma\right)}
			\left\| \hat \gamma_\eps \sqrt{G_\eps} \right\|_{L^\infty}
			\left\| \frac{1}{\eps}\left(1-\hat\gamma_\eps\right)\sqrt{G_\eps} \right\|_{ L^2\left(Mdv\right)}
			\left\| \varphi \right\|_{L^2\left(Mdv\right)}
			\\
			& \leq C
			\left\| \hat q_\eps \right\|_{L^2\left(MM_*dvdv_*d\sigma\right)}
			\left\| \hat \gamma_\eps \sqrt{G_\eps} \right\|_{L^\infty}
			\left\| \frac{1}{\eps}\left(1-\hat\gamma_\eps\right)\right\|_{ L^2\left(Mdv\right)}
			\left\| \varphi \right\|_{L^2\left(Mdv\right)}
			\\
			& +C
			\left\| \hat q_\eps \right\|_{L^2\left(MM_*dvdv_*d\sigma\right)}
			\left\| \hat \gamma_\eps \sqrt{G_\eps} \right\|_{L^\infty}
			\left\| \left(1-\hat\gamma_\eps\right) \hat g_\eps \right\|_{ L^2\left(Mdv\right)}
			\left\| \varphi \right\|_{L^2\left(Mdv\right)}.
		\end{aligned}
	\end{equation*}
	As the support of $\Gamma'(z)-1=\gamma(z)-1+(z-1)\gamma'(z)$ is a subset of  $[2,\infty)$ and since $G_\eps\geq 2$ implies that $\eps\hat g_\eps\geq 2(\sqrt 2 - 1)$, we infer
	\begin{equation}\label{gamma-cv 0}
		\left\|D_\eps^3(\varphi)\right\|_{L^1_\mathrm{loc}\left(dtdx\right)}
		\leq C
		\left\| \left(1-\hat\gamma_\eps\right)\hat g_\eps \right\|_{ L^2_\mathrm{loc}\left(dtdx;L^2\left(Mdv\right)\right)}.
	\end{equation}
	
	Next, from the equi-integrability of $\hat g_\eps^2$ (see Lemma \ref{x-compactness1 0}) and the fact that $1-\hat \gamma_\eps$ is uniformly bounded in $L^\infty$ and converges almost everywhere to zero (possibly up to extraction of a subsequence), we deduce by the Product Limit Theorem that
	\begin{equation}\label{gamma-cv}
		\left(1-\hat \gamma_\eps\right) \hat g_\eps
		\rightarrow 0 \text{ in } L^2_\mathrm{loc}\left(dtdx;L^2\left(Mdv\right)\right).
	\end{equation}
	Thus, we conclude that
	\begin{equation*}
		D_\eps^3(\varphi)\rightarrow 0
		\text{ in } L^1_\mathrm{loc}\left(dtdx\right).
	\end{equation*}

	\noindent$\bullet$ A similar argument provides the convergence of the remaining term $D_\eps^4(\varphi)$.
	Thus, one has by the Cauchy-Schwarz inequality, for any $2<p<\infty$,
	\begin{equation*}% \label{eq-use2}
		\begin{aligned}
			& \left|D_\eps^4(\varphi)\right|
			\\
			& \leq
			\frac 1\eps \left\| \hat q_\eps \right\|_{L^2\left(MM_*dvdv_*d\sigma\right)}
			\left\| \hat \gamma_\eps \hat\gamma_{\eps *}\left(1-\hat\gamma_{\eps}^{\prime}\hat\gamma_{\eps*}^{\prime}\right) \varphi
			\sqrt{G_\eps G_{\eps *}}\right\|_{L^2\left(MM_*dvdv_*d\sigma\right)}
			\\
			& \leq C
			\left\| \hat q_\eps \right\|_{L^2\left(MM_*dvdv_*d\sigma\right)}
			\left\| \hat \gamma_\eps \sqrt{G_\eps} \right\|_{L^\infty}^2
			\left\| \frac 1\eps\left(1-\hat\gamma_{\eps}^{\prime}\hat\gamma_{\eps*}^{\prime}\right) \varphi \right\|_{L^2\left(MM_*dvdv_*d\sigma\right)}
			\\
			& \leq C_p
			\left\| \hat q_\eps \right\|_{L^2\left(MM_*dvdv_*d\sigma\right)}
			\left\| \frac 1\eps \left(1-\hat\gamma_{\eps}\right) \right\|_{L^p\left(Mdv\right)}.
		\end{aligned}
	\end{equation*}
	Therefore, thanks to the bound on $\hat q_\eps$ from Lemma \ref{L2-qlem}, we infer, for any $2<p<\infty$,
	\begin{equation}\label{1-gamma 0}
		\left\|D_\eps^4(\varphi)\right\|_{L^1_\mathrm{loc}\left(dtdx\right)}
		\leq C
		\left\| \frac 1\eps\left(1-\hat\gamma_\eps\right) \right\|_{ L^2_\mathrm{loc}\left(dtdx;L^p\left(Mdv\right)\right)}.
	\end{equation}

	Next, the hypotheses \eqref{gamma-def} on $\gamma(z)$ imply that
	\begin{equation*}
		\left|\frac1\eps \left(1-\hat \gamma_\eps\right) \right|
		\leq \frac 1{2\left(\sqrt 2 - 1\right)} \left|1-\hat \gamma_\eps\right| \left|\hat g_\eps\right|
		\leq \frac 1{2\left(\sqrt 2 - 1\right)} \left|1-\hat \gamma_\eps\right| \left( \left|\Pi\hat g_\eps\right|
		+ \left|\hat g_\eps - \Pi\hat g_\eps\right| \right) ,
	\end{equation*}
	whence
	\begin{equation*}
		\begin{aligned}
			\left|\frac1{\eps^2} \left(1-\hat \gamma_\eps\right)\right|
			& \leq \frac 1{2\left(\sqrt 2 - 1\right)\eps} \left|1-\hat \gamma_\eps\right|
			\left(\left|\Pi\hat g_\eps\right|
			+\left|\hat g_\eps - \Pi\hat g_\eps\right|\right) \\
			& \leq \frac 1{4\left(\sqrt 2 - 1\right)^2} \left|1-\hat \gamma_\eps\right| \left|\hat g_\eps \Pi\hat g_\eps\right|
			+\frac 1{2\left(\sqrt 2 - 1\right)} \left|1-\hat \gamma_\eps\right| \frac{1}{\eps}\left|\hat g_\eps - \Pi\hat g_\eps\right|,
		\end{aligned}
	\end{equation*}
	which, with the relaxation estimate \eqref{relaxation-est} from Lemma \ref{relaxation-control}, shows that, for all $1\leq r<2$,
	\begin{equation*}
		\frac1{\eps^2} \left(1-\hat \gamma_\eps\right)=
		O(1)_{L^1_\mathrm{loc}\left(dtdx ; L^r\left(Mdv\right)\right)}.
	\end{equation*}
	Therefore, for every $2\leq p<4$,
	\begin{equation}\label{special bound}
		\frac1{\eps} \left(1-\hat \gamma_\eps\right)=
		O(1)_{L^2_\mathrm{loc}\left(dtdx ; L^p\left(Mdv\right)\right)}.
	\end{equation}
	
	Moreover, from the equi-integrability of $\hat g_\eps^2$ and the fact that $1-\hat \gamma_\eps$ is uniformly bounded in $L^\infty$ and converges almost everywhere to zero (possibly up to extraction of a subsequence), we deduce by the Product Limit Theorem that
	\begin{equation*}
		\frac 1\eps\left(1-\hat \gamma_\eps\right)
		\rightarrow 0 \text{ in } L^2_\mathrm{loc}\left(dtdx;L^2\left(Mdv\right)\right).
	\end{equation*}
	Therefore, by interpolation, we obtain that, for every $2\leq p<4$,
	\begin{equation}\label{1-gamma}
		\frac 1\eps\left(1-\hat \gamma_\eps\right)
		\rightarrow 0 \text{ in } L^2_\mathrm{loc}\left(dtdx;L^p\left(Mdv\right)\right).
	\end{equation}
	Thus, we conclude that
	\begin{equation*}
		D_\eps^4(\varphi)\rightarrow 0
		\text{ in } L^1_\mathrm{loc}\left(dtdx\right).
	\end{equation*}
	
	On the whole, we have shown that each term from \eqref{defect-decomposition} vanishes as $\eps\rightarrow 0$ in $L^1_\mathrm{loc}(dtdx)$, which leads to the expected convergence and concludes the proof of the lemma.
\end{proof}

\subsection{Decomposition of flux terms}\label{section step 2}

We characterize now the asymptotic behavior of the flux terms.

\begin{lem}\label{flux1}
	The flux terms defined by \eqref{F-def} satisfy
	\begin{equation*}
		\begin{aligned}
			F_\eps(v)
			- \frac1\eps\left(\tilde \rho_\eps+\tilde \theta_\eps\right)\operatorname{Id}
			- \left(\tilde u_\eps\otimes \tilde u_\eps -\frac{|\tilde u_\eps|^2}{3} \operatorname{Id}\right)
			+ \int_{\mathbb{R}^3\times\mathbb{R}^3\times\mathbb{S}^2}
			\hat q_\eps \tilde \phi
			MM_*dvdv_*d\sigma
			& \rightarrow 0,
			\\
			F_\eps\left( \frac{|v|^2}{2}-\frac 52 \right)
			- \frac52 \tilde u_\eps \tilde \theta_\eps
			+ \int_{\mathbb{R}^3\times\mathbb{R}^3\times\mathbb{S}^2}
			\hat q_\eps \tilde \psi
			MM_*dvdv_*d\sigma
			& \rightarrow 0,
		\end{aligned}
	\end{equation*}
	in $L^1_\mathrm{loc}(dtdx)$ as $\eps \to 0$, where $\tilde \phi, \tilde \psi \in L^2\left(Mdv\right)$ are the kinetic momentum and energy fluxes defined by \eqref{phi-psi-def} and \eqref{phi-psi-def inverses}.
	%
	% Moreover, for any collision invariant $\varphi$, we have the global bound
	% \begin{equation*}
	% 	\eps F_\eps (\varphi) = O(1)_{L^\infty\left(dt;L^2\left(dx\right)\right)}.
	% \end{equation*}
\end{lem}

\begin{proof}
	% The global bound simply follows from the Cauchy-Schwarz inequality and the uniform control provided by Lemma \ref{L2-lem}, noticing $g_\eps=\frac 12 \hat g_\eps\left(1+\sqrt{G_\eps}\right)$ so that $\left|g_\eps\gamma_\eps\right|\leq C\left|\hat g_\eps\right|$,
	% \begin{equation*}
	% 	\begin{aligned}
	% 		\left\| \eps F_\eps (\varphi) \right\| _{L^\infty\left(dt;L^2(dx)\right)}
	% 		& = \left\| \int_{\mathbb{R}^3}
	% 		g_\eps \gamma_\eps\varphi \chi\left( {|v|^2\over K_\eps}\right) v M dv \right\| _{L^\infty\left(dt;L^2(dx)\right)}\\
	% 		&\leq \| \hat g_\eps \|_{L^\infty\left(dt;L^2(Mdxdv)\right)}
	% 		\|v\varphi\|_{L^2(Mdv)} =O(1).
	% 	\end{aligned}
	% \end{equation*}

	In order to characterize the asymptotic behavior of fluxes, we use, following the strategy of proof of Proposition \ref{high weak-comp}, the linearized version of the Chapman-Enskog decomposition
	\begin{equation*}
		\hat g_\eps = \Pi \hat g_\eps + \left(\hat g_\eps -\Pi \hat g_\eps\right),
	\end{equation*}
	where $\Pi $ is the orthogonal projection onto $\Ker \cL$ in $L^2\left(Mdv\right)$ and $\hat g_\eps$ is the renormalized fluctuation. Note, however, that we need here a more refined decomposition than the one used in the proof of Proposition \ref{high weak-comp} as we consider now a more singular limit.

	\noindent$\bullet$ Notice that, modulo the diagonal term in the momentum flux
	\begin{equation*}
		\frac1{\eps} \int_{\mathbb{R}^3}
		g_\eps \gamma_\eps \chi\left( {|v|^2\over K_\eps}\right)
		\frac{|v|^2}{3}
		Mdv = \frac1\eps (\tilde \rho_\eps +\tilde \theta_\eps),
	\end{equation*}
	the flux terms have the following structure
	\begin{equation*}
		\tilde F_\eps (\zeta)=\frac1\eps \int_{\mathbb{R}^3}
		g_\eps \gamma_\eps \zeta \chi\left( {|v|^2\over K_\eps}\right) M dv,
	\end{equation*}
	where $\zeta \in \Ker (\cL)^\perp\subset L^2\left(Mdv\right)$. Indeed, it is readily seen that the kinetic fluxes $\phi(v)$ and $\psi(v)$, defined by \eqref{phi-psi-def}, are orthogonal to collision invariants.

	Furthermore, using the identity \eqref{fluct-decomposition}, the fluxes can be rewritten in the following form
	\begin{equation}\label{tilde flux 1}
		\begin{aligned}
			\tilde F_\eps (\zeta)
			& =
			\frac14 \int_{\mathbb{R}^3}
			\hat g_\eps^2 \gamma_\eps \zeta \chi\left( {|v|^2\over K_\eps}\right) M dv
			+
			\frac1\eps \int_{\mathbb{R}^3}
			\hat g_\eps \gamma_\eps \zeta \chi\left( {|v|^2\over K_\eps}\right) M dv
			\\
			& = \frac14\int_{\mathbb{R}^3} (\Pi \hat g_\eps)^2 \zeta M dv
			+ \frac1\eps \int_{\mathbb{R}^3} \hat g_\eps \zeta M dv
			\\
			& + \frac14 \int_{\mathbb{R}^3} \left(\hat g_\eps^2 -\left(\Pi \hat g_\eps\right)^2\right) \gamma_\eps \chi\left( {|v|^2\over K_\eps} \right)\zeta M dv
			\\
			& + \frac14 \int_{\mathbb{R}^3} \left(\gamma_\eps \chi\left( {|v|^2\over K_\eps}\right)-1\right)\left(\Pi \hat g_\eps\right)^2 \zeta M dv
			\\
			& + \frac1\eps \int_{\mathbb{R}^3} \hat g_\eps \left(\gamma_\eps \chi\left( {|v|^2\over K_\eps}\right)-1 \right) \zeta M dv
			\\
			& \eqdefa \frac14\int_{\mathbb{R}^3} \left(\Pi \hat g_\eps\right)^2 \zeta M dv
			+ \frac1\eps \int_{\mathbb{R}^3} \hat g_\eps \zeta M dv
			+ F_\eps^1(\zeta)+F_\eps^2(\zeta)+F_\eps^3(\zeta).
		\end{aligned}
	\end{equation}

	Now, by \eqref{g-L2}, \eqref{relaxation-est} and Lemma \ref{x-compactness1 0}, the remainder terms $F_\eps^1(\zeta)$, $F_\eps^2(\zeta)$ and $F_\eps^3(\zeta)$ will all be shown below to converge to $0$ in $L^1_\mathrm{loc}(dtdx)$ as $\eps\rightarrow 0$. Furthermore, explicit computations will identify the asymptotic behavior of the first term in the above right-hand side. However, there still remains to handle the second term in the right-hand side above, for the limit of this singular expression is not apparent yet (even formally). It is precisely for this term that we have to employ the crucial fact that $\zeta$ belongs to $\Ker (\cL)^\perp$, as we now explain.

	Indeed, note first that the properties of the linearized Boltzmann operator $\mathcal{L}$ stated in Propositions \ref{hilbert-prop} and \ref{coercivity} combined with the Fredholm alternative imply that $\mathcal{L}$ is self-adjoint and Fredholm of index zero on $L^2(Mdv)$. Therefore, its range is exactly the orthogonal complement of its kernel. It follows that any $\zeta \in \Ker (\cL)^\perp\subset L^2\left(Mdv\right)$ belongs to the range of $\mathcal{L}$ and, thus, that there is an inverse $\tilde\zeta\in L^2(Mdv)$ such that
	\begin{equation*}
		\zeta =\cL \tilde \zeta,
	\end{equation*}
	uniquely determined by the fact that it is orthogonal to the kernel of $\mathcal{L}$ (i.e.\ to the collision invariants).

	Then, making use of the simple identity
	\begin{equation*}% \label{bilinear}
		\cL \hat g_\eps = \frac{\eps}2 \cQ\left(\hat g_\eps, \hat g_\eps\right) - \eps \int_{\mathbb{R}^3\times\mathbb{S}^2} \hat q_\eps M_* dv_*d\sigma,
	\end{equation*}
	one has therefore
	\begin{equation}\label{tilde flux 2}
		\begin{aligned}
			\frac1\eps \int_{\mathbb{R}^3} \hat g_\eps \zeta M dv
			& = \frac1\eps \int_{\mathbb{R}^3} \hat g_\eps \mathcal{L}\tilde\zeta M dv
			= \frac1\eps \int_{\mathbb{R}^3} \mathcal{L}\hat g_\eps \tilde\zeta M dv
			\\
			& = \frac12 \int_{\mathbb{R}^3} \cQ\left(\hat g_\eps, \hat g_\eps\right) \tilde \zeta M dv
			- \int_{\mathbb{R}^3\times\mathbb{R}^3\times\mathbb{S}^2}
			\hat q_\eps \tilde \zeta MM_*dvdv_*d\sigma
			\\
			& = \frac12 \int_{\mathbb{R}^3} \cQ\left(\Pi\hat g_\eps, \Pi\hat g_\eps\right) \tilde \zeta M dv
			- \int_{\mathbb{R}^3\times\mathbb{R}^3\times\mathbb{S}^2}
			\hat q_\eps \tilde \zeta MM_*dvdv_*d\sigma
			+ F_\eps^4(\zeta),
		\end{aligned}
	\end{equation}
	where
	\begin{equation*}
		F_\eps^4(\zeta) \eqdefa
		\frac 14
		\int_{\mathbb{R}^3}
		\mathcal{Q}\left(\hat g_\eps - \Pi \hat g_\eps, \hat g_\eps + \Pi \hat g_\eps \right)
		\tilde \zeta M dv
		+
		\frac 14
		\int_{\mathbb{R}^3}
		\mathcal{Q}\left(\hat g_\eps + \Pi \hat g_\eps, \hat g_\eps - \Pi \hat g_\eps \right)
		\tilde \zeta M dv.
	\end{equation*}

	Now, combining \eqref{tilde flux 1} with \eqref{tilde flux 2} and using the identity %(see \cite{BGL2}, for instance)
	\begin{equation*}% \label{L-identity}
		\cQ\left(\Pi \hat g_\eps,\Pi \hat g_\eps\right)=\frac12 \cL \left(\left(\Pi \hat g_\eps\right)^2\right),
	\end{equation*}
	which straightforwardly follows from the following computation, valid for any collision invariant $\varphi$,
	\begin{equation*}
		\varphi'\varphi_*'-\varphi\varphi_*
		=
		\frac 12
		\underbrace{\left(\left(\varphi'+\varphi_*'\right)^2-\left(\varphi+\varphi_*\right)^2\right)}_{=0}
		+
		\frac 12 \left(\varphi^2+\varphi_*^2-\varphi^{\prime 2}-\varphi_*^{\prime 2}\right),
	\end{equation*}
	we deduce that
	\begin{equation}\label{flux1-decomposition1}
		\begin{aligned}
			\tilde F_\eps(\zeta)
			-\frac12\int_{\mathbb{R}^3}
			&
			\left(\Pi \hat g_\eps\right)^2 \zeta Mdv
			+ \int_{\mathbb{R}^3\times\mathbb{R}^3\times\mathbb{S}^2}
			\hat q_\eps \tilde \zeta MM_* dvdv_*d\sigma
			\\
			& =
			\tilde F_\eps(\zeta)
			-\frac14\int_{\mathbb{R}^3} \left(\Pi \hat g_\eps\right)^2 \zeta Mdv
			\\
			& -\frac14\int_{\mathbb{R}^3} \mathcal{L}\left(\left(\Pi \hat g_\eps\right)^2\right) \tilde \zeta Mdv
			+ \int_{\mathbb{R}^3\times\mathbb{R}^3\times\mathbb{S}^2}
			\hat q_\eps \tilde \zeta MM_* dvdv_*d\sigma
			\\
			& =F_\eps^1(\zeta)+F_\eps^2(\zeta)+F_\eps^3(\zeta)+F_\eps^4(\zeta).
		\end{aligned}
	\end{equation}

	Explicit computations show that the advection terms can be conveniently expressed with the moments of $\Pi \hat g_\eps$ (which are equal, by definition, to those of $\hat g_\eps$). Indeed, decomposing
	\begin{equation*}
		\begin{aligned}
			\left(\Pi\hat g_\eps\right)^2
			& =
			\left(\hat\rho_\eps + \hat u_\eps\cdot v + \hat \theta_\eps\left(\frac{|v|^2}{2}-\frac 32\right)\right)^2
			\\
			& =
			\underbrace{
			\hat\rho_\eps^2-3\hat\rho_\eps\hat\theta_\eps-\frac 32\hat\theta_\eps^2
			+ 2\left(\hat\rho_\eps+\hat\theta_\eps\right)\hat u_\eps\cdot v
			+ \left(\frac{\left|\hat u_\eps\right|^2}{3}+\hat\rho_\eps\hat\theta_\eps + \hat \theta_\eps^2\right)|v|^2
			}_{\in\Ker\mathcal{L}}
			\\
			& +
			\underbrace{
			\hat u_\eps \otimes \hat u_\eps : \phi + 2\hat\theta_\eps\hat u_\eps\cdot\psi + \hat\theta_\eps^2\left(\frac{|v|^4}{4}-\frac{5|v|^2}{2}+\frac{15}{4}\right)
			}_{\perp\Ker\mathcal{L}},
		\end{aligned}
	\end{equation*}
	where $\phi$ and $\psi$ are defined in \eqref{phi-psi-def} and $\hat \rho_\eps$, $\hat u_\eps$ and $\hat\theta_\eps$ are, respectively, the density, bulk velocity and temperature associated with $\hat g_\eps$, we find that %(see \cite{BGL2} p. 71 for instance) (there is no p.71 in this reference...)
	\begin{equation}\label{square infty maxwell}
		\begin{aligned}
			\frac12\int_{\mathbb{R}^3}\left(\Pi \hat g_\eps\right)^2 \phi Mdv
			& = \hat u_\eps\otimes\hat u_\eps -\frac{\left|\hat u_\eps\right|^2}{3}\operatorname{Id},
			\\
			\frac12\int_{\mathbb{R}^3}\left(\Pi \hat g_\eps\right)^2 \psi Mdv
			& = \frac 52 \hat u_\eps \hat \theta_\eps.
		\end{aligned}
	\end{equation}
	In particular, it follows from \eqref{flux1-decomposition1} that
	\begin{equation}\label{flux1-decomposition2}
		\begin{aligned}
			\tilde F_\eps(\phi)
			-\left(\hat u_\eps\otimes\hat u_\eps -\frac{\left|\hat u_\eps\right|^2}{3}\operatorname{Id}\right)
			+ \int_{\mathbb{R}^3\times\mathbb{R}^3\times\mathbb{S}^2}
			& \hat q_\eps \tilde \phi MM_* dvdv_*d\sigma
			\\
			& =F_\eps^1(\phi)+F_\eps^2(\phi)+F_\eps^3(\phi)+F_\eps^4(\phi),
			\\
			\tilde F_\eps(\psi)
			-\frac 52 \hat u_\eps \hat \theta_\eps
			+ \int_{\mathbb{R}^3\times\mathbb{R}^3\times\mathbb{S}^2}
			\hat q_\eps \tilde \psi MM_* dvdv_*d\sigma
			& =F_\eps^1(\psi)+F_\eps^2(\psi)+F_\eps^3(\psi)+F_\eps^4(\psi).
		\end{aligned}
	\end{equation}

	Next, writing
	\begin{equation*}
		g_\eps \gamma_\eps \chi\left( {|v|^2\over K_\eps}\right) -\hat g_\eps = \frac12 \hat g_\eps \left(\gamma_\eps \chi\left( {|v|^2\over K_\eps}\right)\left(\sqrt{G_\eps} +1\right)-2\right),
	\end{equation*}
	using the equi-integrability of $\hat g_\eps^2$ from Lemma \ref{x-compactness1 0}, the fact that the second factor $\gamma_\eps \chi\left( {|v|^2\over K_\eps}\right)\left(\sqrt{G_\eps} +1\right)-2$ is uniformly bounded in $L^\infty$ and converges almost everywhere to $0$, observe that, by the Product Limit Theorem,
	\begin{equation}\label{gamma-comparison}
		g_\eps \gamma_\eps \chi\left( {|v|^2\over K_\eps}\right)-\hat g_\eps \to 0 \quad \text{in }L^2_\mathrm{loc}\left(dtdx ; L^2(Mdv)\right).
	\end{equation}
	In particular
	\begin{equation}\label{moments-comparison}
		\tilde \rho_\eps-\hat \rho_\eps \to 0,
		\quad \tilde u_\eps -\hat u_\eps \to 0
		\quad\text{and}\quad
		\tilde \theta_\eps -\hat \theta_\eps \to 0
		\quad\text{in }L^2_\mathrm{loc}(dtdx) \text{ as }\eps \to 0.
	\end{equation}

	Therefore, on the whole, combining \eqref{flux1-decomposition2} with \eqref{moments-comparison}, we see that proving Lemma \ref{flux1} comes down to establishing the vanishing of the four remainder terms $F_\eps^1(\zeta)$, $F_\eps^2(\zeta)$, $F_\eps^3(\zeta)$ and $F_\eps^4(\zeta)$, for any $\zeta=O\left(|v|^3\right)$ as $|v|\rightarrow\infty$.

	\noindent$\bullet$ The first term
	\begin{equation*}
		F_\eps^1(\zeta)
		=
		\frac14 \int_{\mathbb{R}^3} \left(\hat g_\eps^2 -\left(\Pi \hat g_\eps\right)^2\right) \gamma_\eps \chi\left( {|v|^2\over K_\eps} \right)\zeta M dv,
	\end{equation*}
	requires a careful treatment because of the growth of $\zeta(v)=O(|v|^3)$ for large velocities. By the Cauchy-Schwarz inequality, it holds that
	\begin{equation}\label{eq-use3} 
		\begin{aligned}
			\| F_\eps^1(\zeta) \|_{L^1_\mathrm{loc}(dtdx)}
			& \leq \left\| ( \hat g_\eps +\Pi \hat g_\eps ) \gamma_\eps \chi\left( {|v|^2\over K_\eps}\right) \zeta \right\|_{L^2_\mathrm{loc}\left(dtdx;L^2(Mdv)\right)}
			\\
			& \times \left\| \hat g_\eps -\Pi \hat g_\eps \right\|_{L^2_\mathrm{loc}\left(dtdx;L^2(Mdv)\right)}.
		\end{aligned}
	\end{equation}
	
	We already know from Lemma \ref{L2 relaxation} that
	\begin{equation}\label{eq-relaxation}
		\left\|\hat g_\eps-\Pi\hat g_\eps\right\|_{L^2_\mathrm{loc}\left(dtdx;L^2(Mdv)\right)}
		\rightarrow 0
		\quad\text{as }\eps\rightarrow 0.
	\end{equation}
	It remains then to bound the first term in the right-hand side of \eqref{eq-use3} by obtaining a suitable control of large velocities. This follows from Lemma \ref{v2-int} and the definition of $\Pi$, which yields, for all $2\leq p<4$,
	\begin{equation}\label{eq-use5}
		\left|(\hat g_\eps +\Pi \hat g_\eps)\gamma_\eps\right|
		\leq
		\frac{C\left|\hat g_\eps\right|}{1+\sqrt{G_\eps}}
		+
		C\left|\Pi \hat g_\eps\right|
		=O(1)_{L^2_\mathrm{loc}\left(dtdx;L^p(Mdv)\right)}.
	\end{equation}

	Hence, incorporating this last estimate in \eqref{eq-use3} leads to, in view of \eqref{eq-relaxation},
	\begin{equation}\label{flux remainder 1}
		F_\eps^1(\zeta) \to 0 \quad \text{in } L^1_\mathrm{loc}(dtdx) \text{ as } \eps \to 0.
	\end{equation}

	\noindent$\bullet$ The term
	\begin{equation*}
		F_\eps^2(\zeta)
		=
		\frac14 \int_{\mathbb{R}^3} \left(\gamma_\eps \chi\left( {|v|^2\over K_\eps}\right)-1\right)\left(\Pi \hat g_\eps\right)^2 \zeta M dv,
	\end{equation*}
	is easily disposed of, using the equi-integrability of $\hat g_\eps^2$ from Lemma \ref{x-compactness1 0} which implies in particular that
	\begin{equation*}
		\left(\Pi \hat g_\eps\right)^2 \left(1+|v|^p\right) M \text{ is uniformly integrable on } [0,T]\times K\times \mathbb{R}^3,
	\end{equation*}
	for each $T>0$, each compact $K\subset \mathbb{R}^3$ and each $p\in\mathbb{R}$. Indeed, by the Product Limit Theorem, as $\left(\gamma_\eps\chi \left( {|v|^2\over K_\eps}\right)-1\right)$ is bounded in $L^\infty$ and converges almost everywhere to zero, we obtain, for any $p\in\mathbb{R}$,
	\begin{equation}\label{eq-use4}
		\left(\Pi \hat g_\eps\right)^2
		\left( \gamma_\eps \chi\left( {|v|^2\over K_\eps}\right)-1\right)
		\to 0 \quad \text{in } L^1_\mathrm{loc}\left(dtdx;L^1\left(\left(1+|v|^p\right)Mdv\right)\right).
	\end{equation}
	In particular, it follows that
	\begin{equation}\label{flux remainder 2}
		F_\eps^2(\zeta) \to 0 \quad \text{in } L^1_\mathrm{loc}(dtdx) \text{ as } \eps \to 0.
	\end{equation}

	\noindent$\bullet$ In order to get the convergence of
	\begin{equation*}
		F_\eps^3(\zeta)
		=
		\frac1\eps \int_{\mathbb{R}^3} \hat g_\eps \left(\gamma_\eps \chi\left( {|v|^2\over K_\eps}\right)-1 \right) \zeta M dv,
	\end{equation*}
	we use both the estimate \eqref{gaussian-decay 0} on the tails of Gaussian distributions and the convergence \eqref{1-gamma} previously obtained in the proof of Lemma \ref{defect1-lem}.

	Since $\zeta^2(v)=O\left(|v|^6\right)$ as $|v|\to \infty$, one has first, by \eqref{gaussian-decay 0}, that
	\begin{equation}\label{eq-use7}
		\begin{aligned}
			& \left \| \frac1\eps \int_{\mathbb{R}^3}
			\hat g_\eps \gamma_\eps \left(\chi\left( {|v|^2\over K_\eps}\right)-1\right)\zeta M dv\right\|_{L^\infty\left(dt;L^2(dx)\right)}
			\\
			& \leq \frac1{\eps} \| \gamma_\eps\|_{L^\infty} \| \hat g_\eps\|_{L^\infty\left(dt;L^2(Mdxdv)\right)}
			\left(\int_{\mathbb{R}^3} \mathds{1}_{\left\{|v|^2> K_\eps\right\}} \zeta^2 M dv\right)^{1/2}\leq C\eps^{\frac K4-1} |\log \eps|^\frac{7}{4},
		\end{aligned}
	\end{equation}
	which vanishes as soon as $K>4$. Furthermore, by \eqref{1-gamma}, we find
	\begin{equation}\label{eq-use8}
		\begin{aligned}
			& \left \| \frac1\eps \int_{\mathbb{R}^3} \hat g_\eps (\gamma_\eps -1) \zeta M dv\right\|_{L^2_\mathrm{loc}\left(dt;L^1_\mathrm{loc}(dx)\right)}
			\\
			& \leq \left\| \zeta {\gamma_\eps-1\over \eps}\right\|_{L^2_\mathrm{loc}\left(dtdx;L^2(Mdv)\right)}
			\left\| \hat g_\eps \right\|_{L^\infty\left(dt;L^2(Mdxdv)\right)} = o(1).
		\end{aligned}
	\end{equation}
	Thus, combining the preceding estimates yields
	\begin{equation}\label{flux remainder 3}
		F_\eps^3(\zeta) \to 0 \quad \text{in } L^1_\mathrm{loc}(dtdx) \text{ as } \eps \to 0.
	\end{equation}

	\noindent$\bullet$ Finally, the continuity of $\cQ$
	\begin{equation*}
		\left\| \cQ(g,g)\right\|_{L^2(Mdv)} \leq C\| g \|_{L^2(M dv)}^2,
	\end{equation*}
	easily implies that
	\begin{equation}\label{eq-use6}
		\begin{aligned}
		\|F_\eps^4(\zeta)\|_{L^1_\mathrm{loc}\left(dtdx\right)}
		& \leq C \| \tilde \zeta \|_{L^2( M dv)}
		\left\| \hat g_\eps -\Pi \hat g_\eps \right\|_{L^2_\mathrm{loc}\left(dtdx ; L^2( M dv)\right)}
		\\
		& \times \left\| \hat g_\eps +\Pi \hat g_\eps \right\|_{L^2_\mathrm{loc}\left(dtdx; L^2( M dv)\right)}
		\\
		& \leq C \| \tilde \zeta \|_{L^2( M dv)}
		\left\| \hat g_\eps -\Pi \hat g_\eps \right\|_{L^2_\mathrm{loc}\left(dtdx ; L^2( M dv)\right)}
		\\
		& \times \left\| \hat g_\eps \right\|_{L^2_\mathrm{loc}\left(dtdx; L^2( M dv)\right)},
		\end{aligned}
	\end{equation}
	whence, in view of \eqref{eq-relaxation},
	\begin{equation}\label{flux remainder 4}
		F_\eps^4(\zeta) \to 0 \quad \text{in } L^1_\mathrm{loc}(dtdx) \text{ as } \eps \to 0.
	\end{equation}

	On the whole, combining estimates \eqref{flux remainder 1}, \eqref{flux remainder 2}, \eqref{flux remainder 3} and \eqref{flux remainder 4} leads to the expected vanishing of flux remainders which concludes the proof of the lemma.
\end{proof}

\subsection{Decomposition of acceleration terms}\label{section step 3}

It only remains to deal with the acceleration terms.

\begin{lem}\label{acceleration1}
	The acceleration terms defined by \eqref{A-def} satisfy
	\begin{equation*}
		\begin{aligned}
			A_\eps(1) & \to 0,
			\\
			A_\eps(v) -\frac1\eps E_\eps - \left(\tilde \rho_\eps E_\eps +\tilde u_\eps \wedge B_\eps\right) & \to 0,
			\\
			A_\eps\left( \frac{|v|^2}{2}-\frac 52 \right)-\tilde u_\eps \cdot E_\eps & \to 0,
		\end{aligned}
	\end{equation*}
	in $L^1_\mathrm{loc}(dtdx)$ as $\eps \to 0$.
	%
	% Moreover, for any collision invariant $\varphi$, we have the global bound
	% \begin{equation*}
	% 	\eps A_\eps (\varphi) = O(1)_{L^\infty\left(dt;L^2\left(dx\right)\right)}
	% 	+ O(\eps)_{L^\infty\left(dt;L^1\left(dx\right)\right)}.
	% \end{equation*}
\end{lem}

\begin{proof}
	% First, a rapid inspection at the acceleration terms in \eqref{A-def} shows that, since $\left|g_\eps\gamma_\eps\right|,\left|g_\eps\hat\gamma_\eps\right|\leq C\left|\hat g_\eps\right|$,
	% \begin{equation*}
	% 	\begin{aligned}
	% 		& \left\| \eps A_\eps (\varphi)
	% 		-
	% 		E_\eps \cdot
	% 		\int_{\mathbb{R}^3} \hat \gamma_\eps \varphi\chi\left( {|v|^2\over K_\eps}\right) vM dv
	% 		\right\|_{L^1(dx)}
	% 		\\
	% 		& \hspace{10mm} \leq
	% 		C\eps\left\|E_\eps\right\|_{L^2(dx)}
	% 		\left\|\hat g_\eps\right\|_{L^2\left(Mdxdv\right)}
	% 		\left(\left\| v\varphi \right\|_{L^2(Mdv)}+\left\| \nabla_v \varphi \right\|_{L^2(Mdv)}\right)
	% 		\\
	% 		& \hspace{10mm} +
	% 		C\eps\left\|B_\eps\right\|_{L^2(dx)}
	% 		\left\|\hat g_\eps\right\|_{L^2\left(Mdxdv\right)}
	% 		\left\| v\nabla_v \varphi \right\|_{L^2(Mdv)},
	% 	\end{aligned}
	% \end{equation*}
	% which easily provides the expected global bound.

	By definition of the acceleration terms, one has the decomposition
	\begin{equation}\label{acceleration decomposition}
		\begin{aligned}
			A_\eps (\varphi)-\frac1\eps
			E_\eps \cdot \int_{\mathbb{R}^3}
			\varphi v M dv
			- \int_{\mathbb{R}^3}
			g_\eps \gamma_\eps \left(E_\eps+v\wedge B_\eps\right)
			& \cdot \left(\nabla _v \varphi\right)
			\chi\left( {|v|^2\over K_\eps}\right)
			M dv
			\\
			&= A_\eps^1\left(\varphi \right) + A_\eps^2(\varphi) + A_\eps^3(\varphi),
		\end{aligned}
	\end{equation}
	with
	\begin{equation*}
		\begin{aligned}
			A_\eps^1(\varphi)
			& = -\frac1\eps E_\eps \cdot \int_{\mathbb{R}^3}
			\varphi(v)(1-\chi)\left( {|v|^2\over K_\eps}\right)vMdv,
			\\
			A_\eps^2(\varphi)
			& = E_\eps \cdot \int_{\mathbb{R}^3}
			\left( g_\eps \hat \gamma_\eps
			- {1-\hat \gamma_\eps \over \eps}
			- g_\eps \gamma_\eps\right) \varphi(v) \chi\left( {|v|^2\over K_\eps}\right)
			v M dv,
			\\
			A_\eps^3 (\varphi)
			& ={2\over K_\eps} E_\eps \cdot
			\int_{\mathbb{R}^3}
			g_\eps \gamma_\eps
			\varphi(v) \chi' \left( {|v|^2\over K_\eps}\right)
			v M dv.
		\end{aligned}
	\end{equation*}
	As previously, describing the convergence requires a careful treatment.
	
	\noindent$\bullet$ By the Gaussian decay estimate \eqref{gaussian-decay 0} and the uniform $L^2$ bound on $E_\eps$ inherited from the entropy inequality \eqref{entropy1}, we get, for all $v\varphi(v)=O\left(|v|^3\right)$ as $|v|\to \infty$,
	\begin{equation}\label{acceleration remainder 1}
		A_\eps^1(\varphi)=O\left(\eps^{\frac K2-1}|\log \eps|^2\right)_{L^2_\mathrm{loc}(dtdx)},
	\end{equation}
	which tends to $0$ as soon as $K>2$.

	\noindent$\bullet$ For the second term, recalling $\hat \gamma_\eps=\gamma_\eps + \eps g_\eps\gamma'\left(G_\eps\right)$ and writing $g_\eps=\frac 12 \hat g_\eps\left(1+\sqrt{G_\eps}\right)$, an easy computation provides
	\begin{equation*}
		\begin{aligned}
			A_\eps^2(\varphi)
			& =E_\eps \cdot \int_{\mathbb{R}^3}
			\left( \eps g_\eps^2 \gamma '(G_\eps)
			- {1-\hat \gamma_\eps \over \eps} \right)
			\varphi(v) \chi\left( {|v|^2\over K_\eps}\right)
			vM dv
			\\
			& =E_\eps \cdot \int_{\mathbb{R}^3}
			\left( \frac 12 \hat g_\eps \left(1+\sqrt{G_\eps}\right) \left(G_\eps-1\right) \gamma '(G_\eps)
			- {1-\hat \gamma_\eps \over \eps} \right)
			\varphi(v) \chi\left( {|v|^2\over K_\eps}\right)
			vM dv.
		\end{aligned}
	\end{equation*}

	By \eqref{1-gamma}, we have
	\begin{equation*}
		\int_{\mathbb{R}^3}
		{1-\hat \gamma_\eps \over \eps} \varphi(v) \chi\left( {|v|^2\over K_\eps}\right) vM dv
		\to 0 \text{ in } L^2_\mathrm{loc} (dtdx).
	\end{equation*}
	Similarly, from the equi-integrability of $\hat g_\eps^2$ (see Lemma \ref{x-compactness1 0}) and the fact that, by the hypotheses \eqref{gamma-def} on $\gamma(z)$, $\left(1+\sqrt{G_\eps}\right) \left(G_\eps-1\right) \gamma '(G_\eps)$ is uniformly bounded in $L^\infty$ and converges almost everywhere to zero (possibly up to extraction of a subsequence), we deduce by the Product Limit Theorem that
	\begin{equation}\label{gamma-cv 2}
		\hat g_\eps \left(1+\sqrt{G_\eps}\right) \left(G_\eps-1\right) \gamma '(G_\eps)
		\rightarrow 0 \text{ in } L^2_\mathrm{loc}\left(dtdx;L^2\left(Mdv\right)\right).
	\end{equation}
	Therefore, it follows that
	\begin{equation*}
		\int_{\mathbb{R}^3}
		\left( \frac 12 \hat g_\eps \left(1+\sqrt{G_\eps}\right) \left(G_\eps-1\right) \gamma '(G_\eps)
		- {1-\hat \gamma_\eps \over \eps} \right)
		\varphi(v) \chi\left( {|v|^2\over K_\eps}\right)
		vM dv \to 0,
	\end{equation*}
	in $L^2_\mathrm{loc}(dtdx)$, which, when combined with the uniform $L^2$ bound on $E_\eps$, implies that
	\begin{equation}\label{acceleration remainder 2}
		A_\eps^2 (\varphi) \to 0 \quad \text{in } L^1_\mathrm{loc} (dtdx).
	\end{equation}

	\noindent$\bullet$ The last remainder term is easy to control. From the uniform $L^2$ estimates on $E_\eps$ and $\left|g_\eps \gamma_\eps\right|=\frac 12\left|\hat g_\eps \left(1+\sqrt{G_\eps}\right) \gamma\left(G_\eps\right)\right|\leq C\left|\hat g_\eps\right|$, and the fact that
	\begin{equation*}
		{1\over K^2_\eps }\int_{\mathbb{R}^3} \varphi^2 \chi' \left( {|v|^2\over K_\eps}\right)^2 |v|^2 M dv \to 0,
	\end{equation*}
	we deduce that
	\begin{equation}\label{acceleration remainder 3}
		A_\eps^3 (\varphi) \to 0 \quad \text{in } L^1_\mathrm{loc}(dtdx).
	\end{equation}

	Thus, on the whole, incorporating the convergences of remainder terms \eqref{acceleration remainder 1}, \eqref{acceleration remainder 2} and \eqref{acceleration remainder 3} into the decomposition \eqref{acceleration decomposition} and performing direct computations of $\int_{\mathbb{R}^3} \varphi v M dv$ and $ \int_{\mathbb{R}^3} g_\eps \gamma_\eps (E_\eps+v\wedge B_\eps)\cdot \left(\nabla_v  \varphi\right) \chi\left( {|v|^2\over K_\eps}\right) M dv$ leads then to the expected convergences and concludes the proof of the lemma.
	\end{proof}

% \section{Conservation of  densities, total momentum and total energy for \eqref{VMB2}}
\section[Approximate conservation of mass, momentum and energy\ldots]{Approximate conservation of mass, momentum and energy for two species}\label{conservation defects 2 species}

In a way very similar to the one species case from Section \ref{conservation defects 1 species}, we can write approximate conservation laws for the two species Vlasov-Maxwell-Boltzmann system \eqref{VMB2}. However, there are two main differences. The first one is that we do not expect the momentum and energy of each species to be conserved separately, for the mixed collision operators in \eqref{VMB2} do not vanish (even formally) when integrated against collision invariants (except constants) unless they are added together. The second one is that the perturbation in \eqref{VMB2} is more singular so that we do not expect to be able to establish a weak compactness statement such as Lemma \ref{x-compactness1 0}~: the remainders will therefore be controlled by a modulated entropy, which will yield the convergence of remainders to zero at the very end of the proof using Gr\"onwall's lemma (see Chapter \ref{entropy method}).

Having in mind to establish some loop estimates with Gr\"onwall's lemma, which are characteristic of modulated energy (or relative entropy) methods, we impose now some bound from below on the renormalizations. More precisely, we consider here an admissible nonlinearity $\Gamma(z)$ defined by
\begin{equation*}
	\Gamma(z) - 1 = (z-1) \gamma(z),
\end{equation*}
where $\gamma\in C^1\left([0,\infty);\mathbb{R}\right)$ satisfies the following assumptions, for some given $C_1,C_2>0$~:
\begin{equation}\label{gamma-assumption}
	\begin{aligned}
		\gamma(z) & \equiv 1, && \hbox{for all } z \in [0, 2],\\
		\gamma(z) & \leq 1, && \hbox{for all } z \in [0,\infty),\\
		\gamma(z) & \rightarrow 0, && \hbox{as } z\rightarrow\infty,\\
		\gamma(z) & \geq {C_1 \over \left(1+ z\right)^\frac{1}{2}}, && \hbox{for all } z \in [0,\infty),\\
		\left|\gamma'(z)\right| & \leq {C_2 \over \left(1+ z\right)^\frac{3}{2}}, && \hbox{for all } z \in [0,\infty).\\
		% z^{\frac 32}\gamma'(z) & \rightarrow c\in\mathbb{R}\setminus\{0\}, && \hbox{as } z\rightarrow\infty.
	\end{aligned}
\end{equation}
The above hypotheses on $\Gamma(z)$ are clearly more restrictive than the corresponding assumptions \eqref{gamma-def} in the one species case. Note that necessarily $\left|\gamma(z)\right| \leq {2C_2 \over \left(1+ z\right)^\frac{1}{2}}$ and
\begin{equation}\label{gamma below}
	C_1\left(\sqrt z - 1\right)^2 \leq \left(z-1\right)^2\gamma(z)^2 = \left(\Gamma(z)-1\right)^2.
\end{equation}

% \begin{equation*}
% 	(1+z)^\frac 12\Gamma'(z) = (1+z)^\frac 12\gamma(z)+(1+z)^\frac 12 (z-1) \gamma'(z)\geq C_1-C_2>0,
% \end{equation*}
% whence $\sqrt{1+\Gamma^{-1}(y)}$ is uniformly Lipschitz, for one easily computes that
% \begin{equation*}
% 	\frac{d}{dy}\sqrt{1+\Gamma^{-1}(y)}=\frac{1}{2\left(1+\Gamma^{-1}(y)\right)^\frac 12\Gamma'\left(\Gamma^{-1}(y)\right)}\leq\frac{1}{2(C_1-C_2)}.
% \end{equation*}

With the notation $\gamma_\eps^\pm$ for $\gamma\left(G_\eps^\pm\right)$ and $\hat \gamma_\eps^\pm\left(G_\eps\right)$ for $\Gamma'\left(G_\eps^\pm\right)$, the scaled Vlasov-Boltzmann equation in \eqref{VMB2} renormalized relatively to the Maxwellian $M$ with the admissible nonlinearity $\Gamma(z)$ reads
\begin{equation}\label{moment2 0}
	\begin{aligned}
		\d_t \left(g_\eps^\pm\gamma_\eps^\pm\right) + \frac 1\eps v \cdot \nabla_x \left(g_\eps^\pm\gamma_\eps^\pm\right)
		\pm \frac \delta\eps
		& (\eps E_\eps+ v\wedge B_\eps)
		\cdot \nabla_v \left(g_\eps^\pm\gamma_\eps^\pm\right)
		\mp \frac \delta \eps E_\eps \cdot v G_\eps^\pm \hat \gamma_\eps^\pm
		\\
		& =
		\frac 1{\eps^3}\hat \gamma_\eps^\pm\mathcal{Q}\left(G_\eps^\pm,G_\eps^\pm\right)
		+
		\frac {\delta^2}{\eps^3}\hat \gamma_\eps^\pm\mathcal{Q}\left(G_\eps^\pm,G_\eps^\mp\right).
	\end{aligned}
\end{equation}

Following the strategy of Section \ref{conservation defects 1 species}, we also introduce a truncation of large velocities $\chi\left( {|v|^2 \over K_\eps} \right)$, with $K_\eps =K|\log \eps|$, for some large $K>0$ to be fixed later on, and $\chi\in C_c^\infty\left([0,\infty)\right)$ a smooth compactly supported function such that $\mathds{1}_{[0,1]}\leq \chi \leq \mathds{1}_{[0,2]}$.

Thus, multiplying each side of the above equation by $\varphi(v) \chi\left( {|v|^2 \over K_\eps} \right)$, where $\varphi$ is a collision invariant, and averaging with respect to $Mdv$ leads to the moment equations
\begin{equation}\label{moment2}
	\d_t \int_{\mathbb{R}^3} g_\eps^\pm \gamma_\eps^\pm \varphi \chi\left( {|v|^2\over K_\eps} \right) Mdv
	+ \nabla_x\cdot F_\eps^\pm(\varphi) = \pm A_\eps^\pm(\varphi)+ D_\eps^\pm(\varphi) + \Delta_\eps^\pm(\varphi),
\end{equation}
% \begin{equation}
% \begin{aligned}
% \ds \d_t \int Mg_\eps^\pm \gamma_\eps^\pm&\ds \!\!\! \varphi \chi\left( {|v|^2\over K_\eps}\right)dv + \frac1\eps \nabla_x  \cdot \int
% Mg_\eps^\pm \gamma_\eps^\pm \chi\left( {|v|^2\over K_\eps}\right)\varphi vdv\\
% &\ds =\pm \frac\delta\eps E_\eps \cdot \int Mv (1+\eps g_\eps^\pm) \hat \gamma_\eps^\pm\varphi\chi\left( {|v|^2\over K_\eps}\right)dv \\
% &\ds  \pm \frac\delta\eps \int g^\pm_\eps \gamma^\pm_\eps (\eps E_\eps+v\wedge B_\eps) \cdot \nabla_v \left(M \varphi \chi\left( {|v|^2\over K_\eps}\right)\right) dv\\
% &\ds  +\frac1{\eps^3}\int  \hat
% \gamma_\eps^\pm \cQ(G_\eps^\pm, G_\eps^\pm) M\varphi \chi\left( {|v|^2\over K_\eps}\right) dv\\
% &\ds +\frac{\delta^2}{\eps^3}\int \hat \gamma_\eps^\pm  \cQ(G_\eps^\mp, G_\eps^\pm) M\varphi \chi\left( {|v|^2\over K_\eps}\right) dv.
% \end{aligned}
% \end{equation}
with the notations
\begin{equation}\label{F-def 2}
	F_\eps^\pm (\varphi)=\frac1\eps \int_{\mathbb{R}^3}
	g_\eps^\pm \gamma_\eps^\pm \varphi \chi\left( {|v|^2\over K_\eps}\right) v M dv,
\end{equation}
for the fluxes, 
\begin{equation}\label{A-def 2}
	\begin{aligned}
		A_\eps^\pm (\varphi)
		& = \frac\delta\eps E_\eps \cdot
		\int_{\mathbb{R}^3} (1+\eps g_\eps^\pm) \hat \gamma_\eps^\pm \varphi\chi\left( {|v|^2\over K_\eps}\right) vM dv \\
		& + \frac \delta\eps \int_{\mathbb{R}^3}
		g_\eps^\pm \gamma_\eps^\pm (\eps E_\eps+v\wedge B_\eps) \cdot \nabla_v \left( \varphi \chi\left( {|v|^2\over K_\eps}\right) M \right) dv
		\\
		& = \frac\delta\eps E_\eps \cdot
		\int_{\mathbb{R}^3} \hat \gamma_\eps^\pm \varphi\chi\left( {|v|^2\over K_\eps}\right) vM dv
		+ \delta E_\eps \cdot
		\int_{\mathbb{R}^3} g_\eps^\pm \hat \gamma_\eps^\pm \varphi\chi\left( {|v|^2\over K_\eps}\right) vM dv
		\\
		& + \delta E_\eps\cdot\int_{\mathbb{R}^3}
		g_\eps^\pm \gamma_\eps^\pm \varphi \left( \frac{2}{K_\eps}\chi'\left( {|v|^2\over K_\eps}\right) - \chi\left( {|v|^2\over K_\eps}\right) \right) v M dv
		\\
		& + \frac\delta\eps \int_{\mathbb{R}^3}
		g_\eps^\pm \gamma_\eps^\pm (\eps E_\eps+v\wedge B_\eps) \cdot \left( \nabla_v \varphi \right) \chi\left( {|v|^2\over K_\eps}\right) M dv,
	\end{aligned}
\end{equation}
for the acceleration terms, and
\begin{equation}\label{D-def 2}
	\begin{aligned}
		D_\eps^\pm (\varphi)
		& =
		\frac1{\eps^3} \int_{\mathbb{R}^3}
		\hat \gamma_\eps^\pm \cQ\left(G_\eps^\pm, G_\eps^\pm\right) \varphi \chi\left( {|v|^2\over K_\eps}\right) M dv,
		\\
		\Delta_\eps^\pm(\varphi)
		& =
		\frac{\delta^2}{\eps^3} \int_{\mathbb{R}^3}
		\hat \gamma_\eps^\pm \cQ\left(G_\eps^\pm, G_\eps^\mp\right) \varphi \chi\left( {|v|^2\over K_\eps}\right) M dv,
	\end{aligned}
\end{equation}
for the corresponding conservation defects.

By describing the asymptotic behavior of $F_\eps^\pm (\varphi)$, $A_\eps^\pm(\varphi)$, $D_\eps^\pm (\varphi)$ and $\Delta_\eps^\pm(\varphi)$, we will prove the following consistency result (compare with the formal macroscopic conservation laws \eqref{moment-eps two species} and \eqref{difference ohm 2} by setting $\alpha=\delta\eps$, $\beta=\delta$ and $\gamma=1$ therein).

\begin{prop}\label{approx2-prop}
	Let $\left(f_\eps^\pm, E_\eps, B_\eps\right)$ be the sequence of renormalized solutions to the scaled two species Vlasov-Maxwell-Boltzmann system \eqref{VMB2} considered in Theorem \ref{CV-OMHD} for weak interspecies interactions, i.e.\ $\delta=o(1)$ and $\frac\delta\eps$ unbounded, or in Theorem \ref{CV-OMHDSTRONG} for strong interspecies interactions, i.e.\ $\delta=1$, and denote by $\tilde \rho_\eps^\pm$, $\tilde u_\eps^\pm$ and $\tilde \theta_\eps^\pm$ the density, bulk velocity and temperature associated with the renormalized fluctuations $g_\eps^\pm \gamma_\eps^\pm \chi\left( {|v|^2\over K_\eps} \right)$. Further define the hydrodynamic variables
	\begin{equation*}
		\tilde \rho_\eps = \frac{\tilde \rho^+_\eps+\tilde \rho^-_\eps}{2}, \qquad
		\tilde u_\eps = \frac{\tilde u^+_\eps + \tilde u^-_\eps}{2}, \qquad
		\tilde \theta_\eps = \frac{\tilde \theta^+_\eps+\tilde \theta_\eps^-}{2},
	\end{equation*}
	and electrodynamic variables
	\begin{equation*}
		\tilde n_\eps = \tilde \rho^+_\eps - \tilde \rho^-_\eps, \qquad
		\tilde j_\eps = \frac\delta\eps\left(\tilde u^+_\eps - \tilde u^-_\eps\right), \qquad
		\tilde w_\eps = \frac\delta\eps\left(\tilde \theta^+_\eps -\tilde \theta^-_\eps\right).
	\end{equation*}

	Then, one has the approximate hydrodynamic conservation laws
	\begin{equation*}%\label{approximate2}
		\begin{aligned}
			\d_t \tilde \rho_\eps + \frac1\eps \nabla_x\cdot \tilde u_\eps
			& = R_{\eps,1},
			\\
			\d_t \tilde u_\eps
			+ \nabla_x\cdot \left( \tilde u_\eps \otimes \tilde u_\eps -\frac{\left|\tilde u_\eps\right|^2}{3} \operatorname{Id}
			- \int_{\mathbb{R}^3\times\mathbb{R}^3\times\mathbb{S}^2} \frac{\hat q_\eps^++\hat q_\eps^-}{2} \tilde \phi MM_* dvdv_*d\sigma \right)
			\hspace{-60mm} &
			\\
			& = - \frac 1\eps \nabla_x\left(\tilde \rho_\eps+\tilde \theta_\eps\right)
			+ \frac 12\left( \delta \tilde n_\eps E_\eps + \tilde j_\eps \wedge B_\eps \right)
			+ R_{\eps,2},
			\\
			\d_t \left(\frac 32\tilde \theta_\eps-\tilde \rho_\eps\right) + \nabla_x \cdot \left( \frac52 \tilde u_\eps \tilde \theta_\eps
			- \int_{\mathbb{R}^3\times\mathbb{R}^3\times\mathbb{S}^2} \frac{\hat q_\eps^++\hat q_\eps^-}{2} \tilde \psi MM_* dvdv_*d\sigma\right)
			\hspace{-60mm} &
			\\
			& = R_{\eps,3},
		\end{aligned}
	\end{equation*}
	where~:
	\begin{itemize}
		\item
		$\tilde\phi$ and $\tilde\psi$ are defined by \eqref{phi-psi-def} and \eqref{phi-psi-def inverses},
		
		\item
		and the remainders $R_{\eps,i}$, $i=1,2,3$, satisfy
		\begin{equation}\label{approximate remainders}
			\begin{aligned}
				& \left\|R_{\eps,i}\right\|_{W^{-1,1}_\mathrm{loc}\left(dx\right)}
				\\
				& \leq C\delta
				\left\|E_\eps-\bar E\right\|_{L^2(dx)}
				\left\| \left(g_\eps ^+ \gamma_\eps^+ \chi \left({|v|^2\over K_\eps}\right) -\bar g^+,
				g_\eps ^- \gamma_\eps^- \chi \left({|v|^2\over K_\eps}\right) -\bar g^-\right)\right\|_{L^2(Mdxdv)}
				\\
				& + C
				\left\| \left(g_\eps ^+ \gamma_\eps^+ \chi \left({|v|^2\over K_\eps}\right) -\bar g^+,
				g_\eps ^- \gamma_\eps^- \chi \left({|v|^2\over K_\eps}\right) -\bar g^-\right)\right\|_{L^2(Mdxdv)}^2
				\\
				& + C
				\left\|\left(\hat q_\eps^+ - \bar q^+,\hat q_\eps^- - \bar q^-,\hat q_\eps^{+,-} - \bar q^{+,-}, \hat q_\eps^{-,+} - \bar q^{-,+}\right)\right\|_{L^2\left(MM_*dxdvdv_*d\sigma\right)}
				\\
				& \times \left\| \left(g_\eps ^+ \gamma_\eps^+ \chi \left({|v|^2\over K_\eps}\right) -\bar g^+,
				g_\eps ^- \gamma_\eps^- \chi \left({|v|^2\over K_\eps}\right) -\bar g^-\right)\right\|_{L^2(Mdxdv)}
				\\
				& +o(1)_{L^1_\mathrm{loc}(dt)},
			\end{aligned}
		\end{equation}
		for any two given infinitesimal Maxwellians, which differ only by their densities,
		\begin{equation*}
			\bar g^\pm=\bar \rho^\pm+\bar u\cdot v + \bar \theta \left(\frac{|v|^2}{2}-\frac 32\right),
		\end{equation*}
		with $\bar \rho^\pm,\bar u,\bar \theta\in L^\infty(dtdx)\cap L^\infty\left(dt;L^2(dx)\right)$, any collision integrands $\bar q^\pm,\bar q^{\pm,\mp}\in L^\infty\left(dtdx;L^2\left(MM_*dvdv_*d\sigma\right)\right)\cap L^2\left(MM_*dtdxdvdv_*d\sigma\right)$ and any electric field $\bar E \in L^\infty(dtdx)\cap L^\infty\left(dt;L^2(dx)\right)$.
	\end{itemize}
	One also has the approximate electrodynamic conservation laws
	\begin{equation*}%\label{approximate2}
		\begin{aligned}
			\d_t \tilde n_\eps + \frac1\delta \nabla_x\cdot \tilde j_\eps
			& = R_{\eps,4},
			\\
			\frac{\eps^2}{\delta^2}\d_t \tilde j_\eps
			+
			\frac 1\delta \nabla_x\left(\tilde n_\eps+\frac{\eps}{\delta}\tilde w_\eps\right)
			& =
			2\left( E_\eps + \tilde u_\eps \wedge B_\eps \right)
			\\
			& +
			\int_{\mathbb{R}^3\times\mathbb{R}^3\times\mathbb{S}^2}
			\left(\hat q_\eps^{+,-}-\hat q_\eps^{-,+}\right) v MM_* dvdv_*d\sigma
			+ R_{\eps,5},
			\\
			\d_t \left(\frac 32\frac{\eps^2}{\delta^2} \tilde w_\eps-\frac\eps \delta \tilde n_\eps\right)
			& =
			\int_{\mathbb{R}^3\times\mathbb{R}^3\times\mathbb{S}^2}
			\left(\hat q_\eps^{+,-}-\hat q_\eps^{-,+}\right) \left(\frac{|v|^2}{2}-\frac 52\right) MM_* dvdv_*d\sigma
			\\
			& + R_{\eps,6},
		\end{aligned}
	\end{equation*}
	where~:
	\begin{itemize}
		\item
		the remainder $R_{\eps,4}$ also satisfies \eqref{approximate remainders},
		\item
		and the remainders $R_{\eps,i}$, $i=5,6$, converge to $0$ in $L^1_\mathrm{loc}\left(dt;W^{-1,1}_\mathrm{loc}\left(dx\right)\right)$.
	\end{itemize}
\end{prop}

Just like in the proof of proposition \ref{approx1-prop}, the proof of Proposition \ref{approx2-prop} consists in three steps respectively devoted to the study of conservation defects, fluxes and acceleration terms in \eqref{moment2}.

For the sake of clarity, these three steps are respectively detailed in Sections \ref{section step 1 two species}, \ref{section step 2 two species} and \ref{section step 3 two species}, below.

More precisely, Proposition \ref{approx2-prop} will clearly follow from the combination of the approximate conservation laws \eqref{moment2} with Lemma \ref{defect2-lem}, which handles the conservation defects $D_\eps^\pm(\varphi)$ and $\Delta_\eps^\pm(\varphi)$, for any collision invariant $\varphi$, Lemma \ref{flux2}, which establishes the asymptotic behavior of the fluxes $F_\eps^\pm(v)$ and $F_\eps^\pm\left(\frac{|v|^2}{2}-\frac 52\right)$, Lemma \ref{acceleration2}, which characterizes the acceleration terms $A_\eps^\pm(1)$, $A_\eps^\pm(v)$ and $A_\eps^\pm\left(\frac{|v|^2}{2}-\frac 52\right)$ as $\eps\to 0$, and with the following simple estimates of nonlinear terms~:
\begin{equation}\label{nonlinear modulation}
	\begin{aligned}
		& \left|
		\frac 12\left(
		\tilde u_\eps^+\otimes \tilde u_\eps^+
		+
		\tilde u_\eps^-\otimes \tilde u_\eps^-
		\right)
		-\tilde u_\eps\otimes \tilde u_\eps
		\right|
		\\
		& =
		\frac 14
		\big|
		\left(\tilde u_\eps^+ - \bar u\right)
		\otimes \left(\tilde u_\eps^+ - \bar u\right)
		+
		\left(\tilde u_\eps^- - \bar u\right)
		\otimes \left(\tilde u_\eps^- - \bar u\right)
		\\
		& -
		\left(\tilde u_\eps^+ - \bar u\right)
		\otimes \left(\tilde u_\eps^- - \bar u\right)
		-
		\left(\tilde u_\eps^- - \bar u\right)
		\otimes \left(\tilde u_\eps^+ - \bar u\right)
		\big|
		\\
		& \leq C
		\left\| \left(g_\eps ^+ \gamma_\eps^+ \chi \left({|v|^2\over K_\eps}\right) -\bar g^+,
		g_\eps ^- \gamma_\eps^- \chi \left({|v|^2\over K_\eps}\right) -\bar g^-\right)\right\|_{L^2(Mdv)}^2,
		\\
		& \left|
		\frac 12\left(
		\tilde u_\eps^+\tilde \theta_\eps^+
		+
		\tilde u_\eps^-\tilde \theta_\eps^-
		\right)
		-\tilde u_\eps\tilde \theta_\eps
		\right|
		\\
		& =
		\frac 14
		\big|
		\left(\tilde u_\eps^+-\bar u\right)\left(\tilde \theta_\eps^+ - \bar \theta\right)
		+
		\left(\tilde u_\eps^--\bar u\right)\left(\tilde \theta_\eps^- - \bar \theta\right)
		\\
		& -
		\left(\tilde u_\eps^+-\bar u\right)\left(\tilde \theta_\eps^- - \bar \theta\right)
		-
		\left(\tilde u_\eps^--\bar u\right)\left(\tilde \theta_\eps^+ - \bar \theta\right)
		\big|
		\\
		& \leq C
		\left\| \left(g_\eps ^+ \gamma_\eps^+ \chi \left({|v|^2\over K_\eps}\right) -\bar g^+,
		g_\eps ^- \gamma_\eps^- \chi \left({|v|^2\over K_\eps}\right) -\bar g^-\right)\right\|_{L^2(Mdv)}^2.
	\end{aligned}
\end{equation}

\bigskip

	As it turns out, the macroscopic conservation laws provided by Proposition \ref{approx2-prop} will not be sufficient to complete the renormalized relative entropy method in Chapter \ref{entropy method}, for the renormalized electric current $\tilde j_\eps$ in the approximate conservation of momentum of Proposition \ref{approx2-prop} is not controlled by the entropy dissipation. This difficulty will be bypassed by expressing the Lorentz force with the Poynting vector $E_\eps\wedge B_\eps$ (as performed in Section \ref{macroscopic defects}), which will consequently require the handling of the defect measures $m_\eps$ and $a_\eps$, introduced in Section \ref{macroscopic defects}, stemming from the terms $\int_{\mathbb{R}^3}\left(f_\eps^++f_\eps^-\right)v\otimes v dv$ and
	$\begin{pmatrix}
		E_\eps\\B_\eps
	\end{pmatrix}
	\otimes
	\begin{pmatrix}
		E_\eps\\B_\eps
	\end{pmatrix}$, respectively. Fortunately, the defects $m_\eps$ and $a_\eps$ are naturally controlled by the scaled entropy inequality \eqref{entropy2}.

	The following proposition appropriately provides an alternate approximate conservation of momentum law based on the Poynting vector, which will be crucial for the renormalized relative entropy method detailed in Chapter \ref{entropy method}. For convenience, the proof of this proposition is deferred to Section \ref{proofofprop} below.

	\begin{prop}\label{approx3-prop}
		Let $\left(f_\eps^\pm, E_\eps, B_\eps\right)$ be the sequence of renormalized solutions to the scaled two species Vlasov-Maxwell-Boltzmann system \eqref{VMB2} considered in Theorem \ref{CV-OMHD} for weak interspecies interactions, i.e.\ $\delta=o(1)$ and $\frac\delta\eps$ unbounded, or in Theorem \ref{CV-OMHDSTRONG} for strong interspecies interactions, i.e.\ $\delta=1$, and denote by $\tilde \rho_\eps^\pm$, $\tilde u_\eps^\pm$ and $\tilde \theta_\eps^\pm$ the density, bulk velocity and temperature associated with the renormalized fluctuations $g_\eps^\pm \gamma_\eps^\pm \chi\left( {|v|^2\over K_\eps} \right)$. Further define the hydrodynamic variables
		\begin{equation*}
			\tilde \rho_\eps = \frac{\tilde \rho^+_\eps+\tilde \rho^-_\eps}{2}, \qquad
			\tilde u_\eps = \frac{\tilde u^+_\eps + \tilde u^-_\eps}{2}, \qquad
			\tilde \theta_\eps = \frac{\tilde \theta^+_\eps+\tilde \theta_\eps^-}{2}.
		\end{equation*}

		Then, one has the approximate conservation of momentum law
		\begin{equation*}
			\begin{aligned}
				\partial_t & \left(\tilde u_\eps
				+\frac 12 E_\eps\wedge B_\eps
				+\frac 12
				\begin{pmatrix}
					a_{\eps 26}-a_{\eps 35}\\a_{\eps 34}-a_{\eps 16}\\a_{\eps 15}-a_{\eps 24}
				\end{pmatrix}
				\right)
				\\
				& +
				\nabla_x \cdot \left( \tilde u_\eps \otimes \tilde u_\eps
				-\frac{\left|\tilde u_\eps\right|^2}{3} \operatorname{Id}
				+ \frac 1{2\eps^2}m_\eps
				-
				\int_{\mathbb{R}^3\times\mathbb{R}^3\times\mathbb{S}^2} \frac{\hat q_\eps^++\hat q_\eps^-}{2} \tilde \phi MM_* dvdv_*d\sigma
				\right)
				\\
				& -\frac 12 \nabla_x\cdot\left(
				E_\eps\otimes E_\eps + e_\eps + B_\eps\otimes B_\eps + b_\eps
				\right)
				+\nabla_x\left(\frac{|E_\eps|^2+|B_\eps|^2+\operatorname{Tr}a_\eps}{4}\right)
				\\ & = - \frac 1\eps \nabla_x\left(\tilde \rho_\eps+\tilde \theta_\eps\right)+\partial_t \left(o(1)_{L^\infty\left(dt;L^1_{\mathrm{loc}}(dx)\right)}\right) + R_{\eps,7},
			\end{aligned}
		\end{equation*}
		where~:
		\begin{itemize}
			\item
			$\tilde\phi$ is defined by \eqref{phi-psi-def} and \eqref{phi-psi-def inverses},
			
			\item
			the remainder $R_{\eps,7}$ satisfies
			\begin{equation*}
				\begin{aligned}
					\left\| R_{\eps,7}\right\|_{W^{-1,1}_\mathrm{loc}(dx)}
					& \leq
					C_1\int_{\mathbb{R}^3\times\mathbb{R}^3} \left(\frac 1{\eps^2}{h\left(\eps g_\eps^+\right)}
					- \frac{1}{2} \left( g_\eps^+ \gamma_\eps^+\chi \left({|v|^2\over K_\eps}\right) \right)^2 \right)Mdxdv
					\\
					& +
					C_1\int_{\mathbb{R}^3\times\mathbb{R}^3} \left(\frac 1{\eps^2}{h\left(\eps g_\eps^-\right)}
					- \frac{1}{2} \left( g_\eps^- \gamma_\eps^-\chi \left({|v|^2\over K_\eps}\right) \right)^2 \right)Mdxdv
					\\
					& + C_2
					\left\| \left(g_\eps ^+ \gamma_\eps^+\chi \left({|v|^2\over K_\eps}\right) -\bar g^+,
					g_\eps ^- \gamma_\eps^-\chi \left({|v|^2\over K_\eps}\right) -\bar g^-\right)\right\|_{L^2(Mdxdv)}^2
					\\
					& + o(1)_{L^1_{\mathrm{loc}}(dt)},
				\end{aligned}
			\end{equation*}
			for any two given infinitesimal Maxwellians, which differ only by their densities,
			\begin{equation*}
				\bar g^\pm=\bar \rho^\pm+\bar u\cdot v + \bar \theta \left(\frac{|v|^2}{2}-\frac 32\right),
			\end{equation*}
			with $\bar \rho^\pm,\bar u,\bar \theta\in L^\infty(dtdx)\cap L^\infty\left(dt;L^2(dx)\right)$,
			
			\item
			and the symmetric positive definite matrix measures $m_\eps$ and $a_\eps$ are the defects introduced in Section \ref{macroscopic defects} stemming from the terms $\int_{\mathbb{R}^3}\left(f_\eps^++f_\eps^-\right)v\otimes v dv$ and
		$\begin{pmatrix}
			E_\eps\\B_\eps
		\end{pmatrix}
		\otimes
		\begin{pmatrix}
			E_\eps\\B_\eps
		\end{pmatrix}$, respectively, with the notation $e_\eps=\left(a_{\eps ij}\right)_{1\leq i,j\leq 3}$ and $b_\eps=\left(a_{\eps (i+3)(j+3)}\right)_{1\leq i,j\leq 3}$.
		\end{itemize}
	\end{prop}

In the limit $\eps \to 0$ and for well-prepared initial data, we expect that $\left\| g_\eps^\pm \gamma_\eps^\pm \chi \left({|v|^2\over K_\eps}\right) - \bar g^\pm\right\|^2_{L^2(Mdv)}$ should converge strongly to zero for a suitable choice of $\bar g^\pm$. Propositions \ref{approx2-prop} and \ref{approx3-prop} provide then the expected consistency.

A close inspection of \eqref{moment1} and \eqref{moment2 0} shows that the main specificities of the two species case handled here, by comparison with the one species case treated in Section \ref{conservation defects 1 species}, are the following~:
\begin{itemize}
	\item Mixed collision terms do not have all the usual microscopic symmetries, so that we cannot expect macroscopic momentum and energy conservation to hold for each species separately. In other words, there is an exchange of momentum and energy (but not mass) between cations and anions. Symmetries and conservation laws are retrieved by considering the total momentum and total energy.
	
	\item The magnetic force is stronger, so that its contribution to the acceleration terms has to be studied carefully.
	
	\item The assumptions \eqref{gamma-assumption} on the renormalization $\Gamma(z)$ are more restrictive than \eqref{gamma-def}. Whereas \eqref{gamma-def} permitted us to consider a uniformly bounded renormalization if necessary, \eqref{gamma-assumption} requires $\Gamma(z)$ to behave like $\sqrt{z}$ for large values of $z$. Thus, we can no longer have an $L^\infty$ bound on the renormalized fluctuations. However, it still holds true that $\sqrt{G_\eps^\pm} \gamma_\eps^\pm$ and $\sqrt{G_\eps^\pm} \hat\gamma_\eps^\pm$ are uniformly bounded pointwise, which is the only property of $\Gamma(z)$ that we have actually used in Section \ref{conservation defects 1 species}.
	
	The precise usefulness of hypotheses \eqref{gamma-assumption} will become apparent in the proof of Lemma \ref{trunc-lem} below, where the growth properties of $\Gamma$ are employed to compare $\hat g_\eps^\pm$ with $g_\eps^\pm\gamma_\eps$. (Note that this is an instance of the importance of having a theory of renormalized solutions valid for square root renormalizations.) It would be possible to consider here more general renormalizations by working with auxiliary renormalizations when controlling the remainders in Proposition \ref{approx2-prop}. However, this would only add useless technical cumbersomeness to the estimates, which we prefer to avoid by imposing the more restrictive assumptions \eqref{gamma-assumption} on $\Gamma$.
	
	\item The equi-integrability of $|\hat g_\eps|^2$ stated in Lemma \ref{x-compactness1 0} is no longer valid here (only Lemma \ref{x-compactness2 0} holds here) and we have to substitute compactness estimates by the consistency estimates provided by Lemmas \ref{trunc-lem} and \ref{trunc-lem 2} below.
\end{itemize}

To be precise, in Section \ref{conservation defects 1 species}, the equi-integrability of $|\hat g_\eps|^2$ has been used to control $D_\eps^3$ and $D_\eps^4$ in the conservation defects, $F_\eps^1$, $F_\eps^2$, $F_\eps^3$ and $F_\eps^4$ in the fluxes, as well as $A_\eps^2$ in the acceleration terms.

In order to circumvent this lack of compactness, we need to understand how to substitute the convergences \eqref{gamma-cv}, \eqref{1-gamma}, \eqref{gamma-comparison}, \eqref{eq-relaxation}, \eqref{eq-use4} and \eqref{gamma-cv 2} by bounds which will be absorbed through appropriate loop estimates later on (using Gr\"onwall's lemma). This is precisely the goal of the following lemmas, whose technical proofs are postponed to Section \ref{proof of lemma} below, for clarity.

	\begin{lem}\label{trunc-lem}
		For any $2\leq p <4$ and $1\leq q<\infty$, and denoting, for convenience,
		\begin{equation*}
			\left[ g_\eps^\pm \gamma_\eps^\pm -\bar g^\pm \right]
			=
			\left\|g_\eps^\pm \gamma_\eps^\pm -\bar g^\pm\right\|_{L^2\left(Mdv\right)}
			+o(1)_{L^2_{\mathrm{loc}}(dtdx)},
		\end{equation*}
		one has the following consistency estimates
		% \begin{equation*}% \label{trunc-est}
			\begin{align}
				\label{bound 1}
				\left\|\mathds{1}_{\left\{G_\eps^\pm\geq 2\right\}}\hat g_\eps^\pm\right\|_{L^2(Mdv)}
				& \leq C
				\left[ g_\eps^\pm \gamma_\eps^\pm -\bar g^\pm \right],
				\\
				\label{bound 2}
				\left\|g_\eps^\pm\gamma_\eps^\pm\chi\left(\frac{|v|^2}{K_\eps}\right)-\hat g_\eps^\pm\right\|_{L^2(Mdv)}
				& \leq C
				\left[ g_\eps^\pm \gamma_\eps^\pm -\bar g^\pm \right],
				\\
				\label{bound 3}
				\left\|g_\eps^\pm\gamma_\eps^\pm-\hat g_\eps^\pm\right\|_{L^2(Mdv)}
				& \leq C
				\left[ g_\eps^\pm \gamma_\eps^\pm -\bar g^\pm \right],
				\\
				\label{bound 5}
				\left\| \hat g_\eps^\pm -\Pi \hat g_\eps^\pm \right\|_{L^2(Mdv)}
				& \leq C\left[ g_\eps^\pm \gamma_\eps^\pm -\bar g^\pm \right],
				\\
				\label{bound 6}
				\left\|\left(\Pi \hat g_\eps^\pm\right)^2\left(\gamma_\eps^\pm\chi\left(\frac{|v|^2}{K_\eps}\right)-1\right)\right\|_{L^q(Mdv)}
				& \leq C
				\left[ g_\eps^\pm \gamma_\eps^\pm -\bar g^\pm \right]^2,
				\\
				\label{bound 4}
				\left\|\frac 1\eps\mathds{1}_{\left\{G_\eps^\pm\geq 2\right\}}\right\|_{L^p(Mdv)}
				& \leq C
				\left[ g_\eps^\pm \gamma_\eps^\pm -\bar g^\pm \right].
				% \\
				% 				\left\|{1-\chi_\eps^\pm \over \eps }\right\|_{L^2(dx)}
				% 				& \leq C
				% 				\left\|g_\eps^\pm \gamma_\eps^\pm -\bar g^\pm\right\|_{L^2\left(Mdxdv\right)},
				% \\
% 				\label{bound 7}
% 				\left\|\hat g_\eps^\pm\left(1+\sqrt{G_\eps^\pm}\right)\left(G_\eps^\pm -1\right)\gamma'\left(G_\eps^\pm\right)
% 				\right\|_{L^2(Mdxdv)}
% 				& \leq C
% 				R\left(g_\eps^\pm \gamma_\eps^\pm -g^\pm\right).
				% \\
				% \label{bound 8}
				% \left\| \frac1\eps  (1- \tilde\gamma (\eps g_\eps))\sqrt{G_\eps }\chi\left( {|v|^2\over K_\eps}\right) \right\|_{L^2(M dxdv)}
				% & \leq C\| g_\eps \gamma_\eps -g\|_{L^2(Mdxdv)} +O(\eps K_\eps^{p/2})_{L^2(dt)},
				% \\
				% \label{bound 9}
				% \left\| \frac1\eps  (1- \tilde \gamma (\eps g_\eps))\chi\left( {|v|^2\over K_\eps}\right) \right\|_{L^2(M(1+|v|^p)dx dv)}
				% & \leq C\| g_\eps \gamma_\eps -g\|_{L^2(Mdxdv)} +O(\eps K_\eps^{p/2})_{L^2(dt)}.
				% \\
				% \label{bound 100}
				% \left\|\mathds{1}_{\left\{G_\eps^\pm\geq 2\right\}}\hat g_\eps^\pm\right\|_{L^2\left((1+|v|)^2 Mdv\right)}^2
				% & \leq
				% C\int_{\mathbb{R}^3} \left(\frac 1{\eps^2}{h\left(\eps g_\eps^\pm\right)}
				% - \frac{1}{2} \left(\hat g_\eps^\pm\right)^2 \right)Mdv
				% \\
				% \notag
				% & +
				% C\left[ g_\eps^\pm \gamma_\eps^\pm -\bar g^\pm \right]^2.
			\end{align}
		% \end{equation*}
	\end{lem}

	\begin{lem}\label{trunc-lem 2}
		For any $2\leq p <4$ and $1\leq q<2$, one has the following consistency estimates
			\begin{align}
				\label{bound 7}
				\left\|\mathds{1}_{\left\{G_\eps^\pm\geq 2\right\}}\hat g_\eps^\pm\right\|_{L^2(Mdv)}
				& = o(1)_{L^q_\mathrm{loc}(dtdx)},
				\\
				\label{bound 8}
				\left\|g_\eps^\pm\gamma_\eps^\pm\chi\left(\frac{|v|^2}{K_\eps}\right)-\hat g_\eps^\pm\right\|_{L^2(Mdv)}
				& = o(1)_{L^q_\mathrm{loc}(dtdx)},
				\\
				\label{bound 9}
				\left\|g_\eps^\pm\gamma_\eps^\pm-\hat g_\eps^\pm\right\|_{L^2(Mdv)}
				& = o(1)_{L^q_\mathrm{loc}(dtdx)},
				\\
				\label{bound 11}
				\left\| \hat g_\eps^\pm -\Pi \hat g_\eps^\pm \right\|_{L^2(Mdv)}
				& = o(1)_{L^q_\mathrm{loc}(dtdx)},
				\\
				\label{bound 10}
				\left\|\frac 1\eps\mathds{1}_{\left\{G_\eps^\pm\geq 2\right\}}\right\|_{L^p(Mdv)}
				& = o(1)_{L^q_\mathrm{loc}(dtdx)}.
			\end{align}
	\end{lem}

	The following lemma provides a refinement, displaying improved velocity integrability, of the bound \eqref{bound 1} from Lemma \ref{trunc-lem}. It is based on the method of proof of Lemma \ref{v2-int} and is crucial in the demonstration of Proposition \ref{approx3-prop}.

	\begin{lem}\label{trunc-lem 3}
		One has the following consistency estimates
		\begin{equation*}% \label{bound 100}
			\begin{aligned}
				\left\|\mathds{1}_{\left\{G_\eps^\pm\geq 2\right\}}\hat g_\eps^\pm\right\|_{L^2\left((1+|v|)^2 Mdv\right)}^2
				& \leq
				C_1\int_{\mathbb{R}^3} \left(\frac 1{\eps^2}{h\left(\eps g_\eps^\pm\right)}
				- \frac{1}{2} \left( g_\eps^\pm \gamma_\eps^\pm \right)^2 \right)Mdv
				\\
				& + C_2\left\|g_\eps^\pm \gamma_\eps^\pm -\bar g^\pm\right\|_{L^2\left(Mdv\right)}^2
				+o(1)_{L^1_{\mathrm{loc}}(dtdx)},
				% \\
				% & +
				% C\left[ g_\eps^\pm \gamma_\eps^\pm -\bar g^\pm \right]^2.
			\end{aligned}
		\end{equation*}
		and
		\begin{equation*}
			\begin{aligned}
				\left\|\mathds{1}_{\left\{G_\eps^\pm\geq 2\right\}}\hat g_\eps^\pm\right\|_{L^2_\mathrm{loc}\left(dx;L^2\left((1+|v|)^2 Mdv\right)\right)}^2 \hspace{-30mm} &
				\\
				& \leq
				C_1\int_{\mathbb{R}^3\times\mathbb{R}^3} \left(\frac 1{\eps^2}{h\left(\eps g_\eps^\pm\right)}
				- \frac{1}{2} \left( g_\eps^\pm \gamma_\eps^\pm\chi \left({|v|^2\over K_\eps}\right) \right)^2 \right)Mdxdv
				\\
				& + C_2\left\|g_\eps^\pm \gamma_\eps^\pm\chi \left({|v|^2\over K_\eps}\right) -\bar g^\pm\right\|_{L^2\left(Mdxdv\right)}^2
				+o(1)_{L^1_{\mathrm{loc}}(dt)}.
			\end{aligned}
		\end{equation*}
	\end{lem}

	The next result comprises yet another important consistency estimate following from the preceding lemma. This estimate is not used in the present chapter, we only record it here for later reference in the proof of Theorem \ref{CV-OMHDSTRONG} in Chapter \ref{entropy method} for strong interspecies interactions.

	\begin{lem}\label{trunc-lem 4}
		One has the following consistency estimates
		\begin{equation*}
			\begin{aligned}
				\left\|\left(\hat g_\eps^+-\hat g_\eps^--\hat n_\eps\right)\hat g_\eps^\pm\right\|_{L^1\left((1+|v|)^2 Mdv\right)}
				& \leq
				C_1\int_{\mathbb{R}^3} \left(\frac 1{\eps^2}{h\left(\eps g_\eps^\pm\right)}
				- \frac{1}{2} \left( g_\eps^\pm \gamma_\eps^\pm \right)^2 \right)Mdv
				\\
				& + C_2\left\|g_\eps^\pm \gamma_\eps^\pm -\bar g^\pm\right\|_{L^2\left(Mdv\right)}^2
				+o(1)_{L^1_{\mathrm{loc}}(dtdx)},
			\end{aligned}
		\end{equation*}
		and
		\begin{equation*}
			\begin{aligned}
				\left\|\left(\hat g_\eps^+-\hat g_\eps^--\hat n_\eps\right)\hat g_\eps^\pm\right\|_{L^1_\mathrm{loc}\left(dx;L^1\left((1+|v|)^2 Mdv\right)\right)} \hspace{-30mm} &
				\\
				& \leq
				C_1\int_{\mathbb{R}^3\times\mathbb{R}^3} \left(\frac 1{\eps^2}{h\left(\eps g_\eps^\pm\right)}
				- \frac{1}{2} \left( g_\eps^\pm \gamma_\eps^\pm\chi \left({|v|^2\over K_\eps}\right) \right)^2 \right)Mdxdv
				\\
				& + C_2\left\|g_\eps^\pm \gamma_\eps^\pm\chi \left({|v|^2\over K_\eps}\right) -\bar g^\pm\right\|_{L^2\left(Mdxdv\right)}^2
				+o(1)_{L^1_{\mathrm{loc}}(dt)},
			\end{aligned}
		\end{equation*}
		where $\hat n_\eps$ is the charge density associated with $\hat g_\eps^\pm$.
	\end{lem}

\subsection{Conservation defects}\label{section step 1 two species}

The first step of the proof of Proposition \ref{approx2-prop} is to establish the control of conservation defects.

\begin{lem}\label{defect2-lem}
	The conservation defects defined by \eqref{D-def 2} satisfy the controls, for any collision invariant $\varphi$,
	\begin{equation*}
		\begin{aligned}
			\left|D_\eps^{\pm}(\varphi)\right|
			& \leq C
			\left\| \hat q_\eps^\pm - \bar q^\pm \right\|_{L^2\left(MM_*dvdv_*d\sigma\right)}
			\left\| g_\eps^\pm \gamma_\eps^\pm - \bar g^\pm \right\|_{ L^2\left(Mdv\right)}
			\\
			& + o(1)_{L^1_\mathrm{loc}(dtdx)},
			\\
			\left|\Delta_\eps^+(\varphi)+\Delta_\eps^-(\varphi)\right|
			& \leq
			C\delta
			\left\| \left( \hat q_\eps^{+,-} - \bar q^{+,-},
			\hat q_\eps^{-,+} - \bar q^{-,+} \right) \right\|_{L^2\left(MM_*dvdv_*d\sigma\right)}
			\\
			& \times
			\left\|\left(g_\eps^+\gamma_\eps^+ - \bar g^+,
			g_\eps^-\gamma_\eps^- - \bar g^-\right) \right\|_{L^2\left(Mdv\right)}
			+o(1)_{L^1_\mathrm{loc}(dtdx)},
			\\
			\left|\Delta_\eps^\pm(1)\right|
			& \leq
			C\delta
			\left\| \hat q_\eps^{\pm,\mp} - \bar q^{\pm,\mp} \right\|_{L^2\left(MM_*dvdv_*d\sigma\right)}
			\\
			& \times
			\left\|\left(g_\eps^+\gamma_\eps^+ - \bar g^+,
			g_\eps^-\gamma_\eps^- - \bar g^-\right) \right\|_{L^2\left(Mdv\right)}
			+o(1)_{L^1_\mathrm{loc}(dtdx)},
			\\
			\frac\eps\delta\Delta_\eps^\pm(\varphi)
			& =
			\int_{\mathbb{R}^3\times\mathbb{R}^3\times\mathbb{S}^2}
			\hat q_\eps^{\pm,\mp}
			\varphi MM_*dvdv_*d\sigma
			+
			o(1)_{L^1_\mathrm{loc}(dtdx)}.
		\end{aligned}
	\end{equation*}
\end{lem}

\begin{proof}
	We follow the proof of Lemma \ref{defect1-lem} in the one species case. Thus, we first note that $D_\eps^\pm(\varphi)$ can be decomposed exactly as in \eqref{defect-decomposition}, which yields
	\begin{equation}\label{defect-decomposition-two}
		\begin{aligned}
			D_\eps^\pm(\varphi) & =
			\frac{\eps}{4}
			\int_{\mathbb{R}^3\times\mathbb{R}^3\times\mathbb{S}^2}
			\hat \gamma_\eps^\pm \varphi \chi\left( {|v|^2\over K_\eps} \right)
			\hat q_\eps^{\pm 2} MM_*dvdv_*d\sigma
			\\
			& - \frac{1}{\eps}\int_{\mathbb{R}^3\times\mathbb{R}^3\times\mathbb{S}^2}
			\hat \gamma_\eps^\pm \varphi \left(1-\chi\left( {|v|^2\over K_\eps} \right)\right)
			\hat q_\eps^\pm
			\sqrt{G_\eps^\pm G_{\eps *}^\pm} MM_* dvdv_*d\sigma
			\\
			& + \frac{1}{\eps}\int_{\mathbb{R}^3\times\mathbb{R}^3\times\mathbb{S}^2}
			\hat \gamma_\eps^\pm \left(1-\hat\gamma_{\eps *}^\pm\right) \varphi
			\hat q_\eps^\pm
			\sqrt{G_\eps^\pm G_{\eps *}^\pm} MM_* dvdv_*d\sigma
			\\
			& + \frac{1}{\eps}\int_{\mathbb{R}^3\times\mathbb{R}^3\times\mathbb{S}^2}
			\hat \gamma_\eps^\pm \hat\gamma_{\eps *}^\pm \left(1-\hat\gamma_{\eps}^{\pm\prime}\hat\gamma_{\eps *}^{\pm\prime}\right) \varphi
			\hat q_\eps^\pm
			\sqrt{G_\eps^\pm G_{\eps *}^\pm} MM_* dvdv_*d\sigma
			\\
			& - \frac{\eps}{4}\int_{\mathbb{R}^3\times\mathbb{R}^3\times\mathbb{S}^2}
			\hat \gamma_\eps^\pm \hat\gamma_{\eps *}^\pm \hat\gamma_{\eps}^{\pm\prime}\hat\gamma_{\eps *}^{\pm\prime}
			\varphi
			\hat q_\eps^{\pm 2} MM_* dvdv_*d\sigma
			\\
			& \eqdefa D^{\pm 1}_\eps(\varphi)+D^{\pm 2}_\eps(\varphi)+D^{\pm 3}_\eps(\varphi)+D^{\pm 4}_\eps(\varphi)+D^{\pm 5}_\eps(\varphi),
		\end{aligned}
	\end{equation}
	where we have used that $\varphi$ is a collision invariant to symmetrize the last term.
	
	Then, we estimate the defects $D^{\pm 1}_\eps(\varphi)$, $D^{\pm 2}_\eps(\varphi)$ and $D^{\pm 5}_\eps(\varphi)$ exactly as $D^{1}_\eps(\varphi)$, $D^{2}_\eps(\varphi)$ and $D^{5}_\eps(\varphi)$ in the one species case. Indeed, the control of these terms only depends on the bounds provided by the relative entropy and entropy dissipation through Lemmas \ref{L2-lem} and \ref{L2-qlem} and, therefore, holds in both the one species and two species cases. Thus, we have that
	\begin{equation*}
		D^{\pm 1}_\eps(\varphi), D^{\pm 2}_\eps(\varphi), D^{\pm 5}_\eps(\varphi) \to 0
		\text{ in }L^1_\mathrm{loc}(dtdx)
		\text{ as }\eps \to 0.
	\end{equation*}
	
	The remaining terms cannot be handled as in Lemma \ref{defect1-lem} and do not necessarily vanish, because of the lack of equi-integrability of $\left|\hat g_\eps^\pm\right|^2$. Note, however, that the estimates \eqref{gamma-cv 0} and \eqref{1-gamma 0} can be reproduced here without difficulty, which yields, for any $2<p<\infty$,
	\begin{equation*}
		\begin{aligned}
		\left|D_\eps^{\pm 3}(\varphi)\right|
		& \leq C
		\left\| \hat q_\eps^\pm \right\|_{L^2\left(MM_*dvdv_*d\sigma\right)}
		\left\| \left(1-\hat\gamma_\eps^\pm\right)\hat g_\eps^\pm \right\|_{ L^2\left(Mdv\right)}
		\\
		& \leq C
		\left\| \hat q_\eps^\pm - \bar q^\pm \right\|_{L^2\left(MM_*dvdv_*d\sigma\right)}
		\left\| \left(1-\hat\gamma_\eps^\pm\right)\hat g_\eps^\pm \right\|_{ L^2\left(Mdv\right)}
		\\
		& + C
		\left\| \bar q^\pm \right\|_{L^\infty\left(dtdx;L^2\left(MM_*dvdv_*d\sigma\right)\right)}
		\left\| \left(1-\hat\gamma_\eps^\pm\right)\hat g_\eps^\pm \right\|_{ L^2\left(Mdv\right)},
		\\
		\left|D_\eps^{\pm 4}(\varphi)\right|
		& \leq C
		\left\| \hat q_\eps^\pm \right\|_{L^2\left(MM_*dvdv_*d\sigma\right)}
		\left\| \frac 1\eps\left(1-\hat\gamma_\eps^\pm\right) \right\|_{ L^p\left(Mdv\right)}
		\\
		& \leq C
		\left\| \hat q_\eps^\pm - \bar q^\pm \right\|_{L^2\left(MM_*dvdv_*d\sigma\right)}
		\left\| \frac 1\eps\left(1-\hat\gamma_\eps^\pm\right) \right\|_{ L^p\left(Mdv\right)}
		\\
		& + C
		\left\| \bar q^\pm \right\|_{L^\infty\left(dtdx;L^2\left(MM_*dvdv_*d\sigma\right)\right)}
		\left\| \frac 1\eps\left(1-\hat\gamma_\eps^\pm\right) \right\|_{ L^p\left(Mdv\right)}.
		\end{aligned}
	\end{equation*}
	Then, instead of using the convergences \eqref{gamma-cv} and \eqref{1-gamma} (which are not valid here), we employ the pairs of controls \eqref{bound 1}-\eqref{bound 7} and \eqref{bound 4}-\eqref{bound 10}, respectively, provided by Lemmas \ref{trunc-lem} and \ref{trunc-lem 2}, which yields
	\begin{equation*}
		\begin{aligned}
		\left|D_\eps^{\pm 3}(\varphi)\right|
		& \leq C
		\left\| \hat q_\eps^\pm - \bar q^\pm \right\|_{L^2\left(MM_*dvdv_*d\sigma\right)}
		\left\| g_\eps^\pm \gamma_\eps^\pm - \bar g^\pm \right\|_{ L^2\left(Mdv\right)} + o(1)_{L^1_\mathrm{loc}(dtdx)},
		\\
		\left|D_\eps^{\pm 4}(\varphi)\right|
		& \leq C
		\left\| \hat q_\eps^\pm - \bar q^\pm \right\|_{L^2\left(MM_*dvdv_*d\sigma\right)}
		\left\| g_\eps^\pm \gamma_\eps^\pm - \bar g^\pm \right\|_{ L^2\left(Mdv\right)} + o(1)_{L^1_\mathrm{loc}(dtdx)}.
		\end{aligned}
	\end{equation*}
	
	On the whole, combining the preceding estimates clearly concludes the proof of the control of $D^\pm_\eps(\varphi)$.

	We turn to the analysis of the mixed terms $\Delta_\eps^\pm(\varphi)$, which are handled in a very similar fashion. Let us just recall that we do not expect the conservation of momentum and energy to hold for each species separately, so that only the total mixed conservation defects $\Delta_\eps^+(\varphi)+\Delta_\eps^-(\varphi)$ are expected to vanish in the limit, in general.
	
	First, we decompose
	\begin{equation*}
		\Delta_\eps^+(\varphi)+\Delta_\eps^-(\varphi) = \Delta^{1}_\eps(\varphi)+\Delta^{2}_\eps(\varphi)+\Delta^{3}_\eps(\varphi)+\Delta^{4}_\eps(\varphi)+\Delta^{5}_\eps(\varphi),
	\end{equation*}
	where we define
	\begin{equation}\label{defect-decomposition-mixed}
		\begin{aligned}
			\Delta^{1}_\eps(\varphi)
			& =
			\frac{\eps}{4}
			\int_{\mathbb{R}^3\times\mathbb{R}^3\times\mathbb{S}^2}
			\hat \gamma_\eps^+
			\varphi \chi\left( {|v|^2\over K_\eps} \right)
			\left(\hat q_\eps^{+,-}\right)^2
			MM_*dvdv_*d\sigma
			\\
			& +
			\frac{\eps}{4}
			\int_{\mathbb{R}^3\times\mathbb{R}^3\times\mathbb{S}^2}
			\hat \gamma_\eps^-
			\varphi \chi\left( {|v|^2\over K_\eps} \right)
			\left(\hat q_\eps^{-,+}\right)^2
			MM_*dvdv_*d\sigma,
			\\
			\Delta^2_\eps(\varphi)
			& =
			- \frac{\delta}{\eps}\int_{\mathbb{R}^3\times\mathbb{R}^3\times\mathbb{S}^2}
			\hat \gamma_\eps^+
			\varphi \left(1-\chi\left( {|v|^2\over K_\eps} \right)\right)
			\hat q_\eps^{+,-}
			\sqrt{G_\eps^+ G_{\eps *}^-}
			MM_* dvdv_*d\sigma
			\\
			&
			- \frac{\delta}{\eps}\int_{\mathbb{R}^3\times\mathbb{R}^3\times\mathbb{S}^2}
			\hat \gamma_\eps^-
			\varphi \left(1-\chi\left( {|v|^2\over K_\eps} \right)\right)
			\hat q_\eps^{-,+}
			\sqrt{G_\eps^- G_{\eps *}^+}
			MM_* dvdv_*d\sigma,
			\\
			\Delta^3_\eps(\varphi)
			& =
			\frac{\delta}{\eps}\int_{\mathbb{R}^3\times\mathbb{R}^3\times\mathbb{S}^2}
			\hat \gamma_\eps^+ \left(1-\hat\gamma_{\eps *}^-\right) \varphi
			\hat q_\eps^{+,-}
			\sqrt{G_\eps^+ G_{\eps *}^-} MM_* dvdv_*d\sigma
			\\
			& +
			\frac{\delta}{\eps}\int_{\mathbb{R}^3\times\mathbb{R}^3\times\mathbb{S}^2}
			\hat \gamma_\eps^- \left(1-\hat\gamma_{\eps *}^+\right) \varphi
			\hat q_\eps^{-,+}
			\sqrt{G_\eps^- G_{\eps *}^+} MM_* dvdv_*d\sigma,
			\\
			\Delta^4_\eps(\varphi)
			& =
			\frac{\delta}{\eps}\int_{\mathbb{R}^3\times\mathbb{R}^3\times\mathbb{S}^2}
			\hat \gamma_\eps^+ \hat\gamma_{\eps *}^- \left(1-\hat\gamma_{\eps}^{+\prime}\hat\gamma_{\eps *}^{-\prime}\right) \varphi
			\hat q_\eps^{+,-}
			\sqrt{G_\eps^+ G_{\eps *}^-} MM_* dvdv_*d\sigma
			\\
			& +
			\frac{\delta}{\eps}\int_{\mathbb{R}^3\times\mathbb{R}^3\times\mathbb{S}^2}
			\hat \gamma_\eps^- \hat\gamma_{\eps *}^+ \left(1-\hat\gamma_{\eps}^{-\prime}\hat\gamma_{\eps *}^{+\prime}\right) \varphi
			\hat q_\eps^{-,+}
			\sqrt{G_\eps^- G_{\eps *}^+} MM_* dvdv_*d\sigma,
			\\
			\Delta^5_\eps(\varphi)
			& =
			- \frac{\eps}{4}\int_{\mathbb{R}^3\times\mathbb{R}^3\times\mathbb{S}^2}
			\hat \gamma_\eps^+ \hat\gamma_{\eps *}^- \hat\gamma_{\eps}^{+\prime}\hat\gamma_{\eps *}^{-\prime}
			\varphi
			\left(\hat q_\eps^{+,-}\right)^2 MM_* dvdv_*d\sigma
			\\
			& - \frac{\eps}{4}\int_{\mathbb{R}^3\times\mathbb{R}^3\times\mathbb{S}^2}
			\hat \gamma_\eps^- \hat\gamma_{\eps *}^+ \hat\gamma_{\eps}^{-\prime}\hat\gamma_{\eps *}^{+\prime}
			\varphi
			\left(\hat q_\eps^{-,+}\right)^2 MM_* dvdv_*d\sigma.
		\end{aligned}
	\end{equation}
	Note that we use the fact that $\varphi$ is a collision invariant, i.e.\ that $\varphi+\varphi_* = \varphi'+\varphi_*'$, only to symmetrize $\Delta_\eps^5(\varphi)$~:
	\begin{equation*}
		\begin{aligned}
			\Delta^5_\eps(\varphi)
			& =
			\frac{\delta}{\eps}\int_{\mathbb{R}^3\times\mathbb{R}^3\times\mathbb{S}^2}
			\hat \gamma_\eps^+ \hat\gamma_{\eps *}^- \hat\gamma_{\eps}^{+\prime}\hat\gamma_{\eps *}^{-\prime}
			\varphi
			\hat q_\eps^{+,-}
			\sqrt{G_\eps^+ G_{\eps *}^-} MM_* dvdv_*d\sigma
			\\
			& +
			\frac{\delta}{\eps}\int_{\mathbb{R}^3\times\mathbb{R}^3\times\mathbb{S}^2}
			\hat \gamma_\eps^- \hat\gamma_{\eps *}^+ \hat\gamma_{\eps}^{-\prime}\hat\gamma_{\eps *}^{+\prime}
			\varphi
			\hat q_\eps^{-,+}
			\sqrt{G_\eps^- G_{\eps *}^+} MM_* dvdv_*d\sigma
			\\
			& =
			\frac{\delta}{\eps}\int_{\mathbb{R}^3\times\mathbb{R}^3\times\mathbb{S}^2}
			\hat \gamma_\eps^+ \hat\gamma_{\eps *}^- \hat\gamma_{\eps}^{+\prime}\hat\gamma_{\eps *}^{-\prime}
			\left(\varphi+\varphi_*\right)
			\hat q_\eps^{+,-}
			\sqrt{G_\eps^+ G_{\eps *}^-} MM_* dvdv_*d\sigma
			\\
			& =
			\frac{\delta}{2\eps}\int_{\mathbb{R}^3\times\mathbb{R}^3\times\mathbb{S}^2}
			\hat \gamma_\eps^+ \hat\gamma_{\eps *}^- \hat\gamma_{\eps}^{+\prime}\hat\gamma_{\eps *}^{-\prime}
			\left(\varphi+\varphi_*\right)
			\hat q_\eps^{+,-}
			\\
			& \times \left(\sqrt{G_\eps^+ G_{\eps *}^-}-\sqrt{G_\eps^{+\prime} G_{\eps *}^{-\prime}}\right) MM_* dvdv_*d\sigma
			\\
			& =
			- \frac{\eps}{4}\int_{\mathbb{R}^3\times\mathbb{R}^3\times\mathbb{S}^2}
			\hat \gamma_\eps^+ \hat\gamma_{\eps *}^- \hat\gamma_{\eps}^{+\prime}\hat\gamma_{\eps *}^{-\prime}
			\left(\varphi+\varphi_*\right)
			\left(\hat q_\eps^{+,-}\right)^2
			MM_* dvdv_*d\sigma.
		\end{aligned}
	\end{equation*}
	This is precisely the point where we need to consider the sum of the mixed collision integrands over both species. Note that, if $\varphi \equiv 1$, then we have $\varphi = \varphi_*$, so that the conservation defects can be dealt with separately.
	
	Anyway, the terms in \eqref{defect-decomposition-mixed} are all similar to those in \eqref{defect-decomposition-two}. We even have an additional factor $\delta$ in the terms $\Delta^2_\eps(\varphi)$, $\Delta^3_\eps(\varphi)$ and $\Delta^4_\eps(\varphi)$. Therefore, with the exact same arguments used to treat the conservation defects $D^\pm_\eps(\varphi)$, we conclude the proof of the controls over $\Delta^+_\eps(\varphi)+\Delta^-_\eps(\varphi)$ and $\Delta^\pm_\eps(1)$.
	
	Finally, in order to derive the control of $\frac\eps\delta\Delta_\eps^\pm(\varphi)$, we consider the simple decomposition, writing $\sqrt{G_{*\eps}^\mp}=1+\frac\eps 2\hat g_{*\eps}^\mp$,
	\begin{equation*}
		\begin{aligned}
			\frac\eps\delta\Delta_\eps^\pm(\varphi)
			& -
			\int_{\mathbb{R}^3\times\mathbb{R}^3\times\mathbb{S}^2}
			\hat q_\eps^{\pm,\mp}
			\varphi MM_*dvdv_*d\sigma
			\\
			& =
			\frac\eps 2
			\int_{\mathbb{R}^3\times\mathbb{R}^3\times\mathbb{S}^2} \hat\gamma_\eps^\pm\sqrt{G_\eps^\pm}
			\chi\left( {|v|^2\over K_\eps}\right)\hat g_{\eps *}^\mp\hat q_\eps^{\pm,\mp}\varphi MM_*dvdv_*d\sigma
			\\
			& +
			\int_{\mathbb{R}^3\times\mathbb{R}^3\times\mathbb{S}^2} \left(\hat\gamma_\eps^\pm\sqrt{G_\eps^\pm}
			\chi\left( {|v|^2\over K_\eps}\right)-1\right)\hat q_\eps^{\pm,\mp}\varphi MM_*dvdv_*d\sigma
			\\
			& +
			\frac{\eps^2}{4\delta}
			\int_{\mathbb{R}^3\times\mathbb{R}^3\times\mathbb{S}^2} \hat\gamma_\eps^\pm \left(\hat q_\eps^{\pm,\mp}\right)^2
			\varphi\chi\left( {|v|^2\over K_\eps}\right)MM_*dvdv_*d\sigma.
		\end{aligned}
	\end{equation*}
	
	Then, since $\frac\eps\delta$ vanishes, $\varphi\chi\left( {|v|^2\over K_\eps}\right)$ is bounded pointwise by a constant multiple of $\left|\log \eps\right|$, the renormalized fluctuations $\hat g_{\eps *}^\mp$ are uniformly bounded in $L^\infty\left(dt ; L^2\left(M_* dxdv_*\right)\right)$ and the collision integrands $\hat q_\eps^{\pm,\mp}$ are uniformly bounded in $L^2\left(MM_*dtdxdvdv_*d\sigma\right)$, we find that the first and third terms from the right-hand side above vanish in $L^1_{\mathrm{loc}}\left(dtdx\right)$.
	
	Further noticing that $\left(\hat\gamma_\eps^\pm\sqrt{G_\eps^\pm}\chi\left( {|v|^2\over K_\eps}\right)-1\right)\varphi$ is dominated by $|\varphi|$ and converges almost everywhere to $0$, it is easily shown that the second term in the above right-hand side vanishes in $L^1_\mathrm{loc}(dtdx)$, as well.
	
	The proof of the lemma is now complete.
\end{proof}

\subsection{Decomposition of flux terms}\label{section step 2 two species}

We characterize now the asymptotic behavior of the flux terms.

\begin{lem}\label{flux2}
	The flux terms defined by \eqref{F-def 2} satisfy
	\begin{equation*}
		\begin{aligned}
			& \left| F_\eps^\pm (v) - \frac 1\eps\left(\tilde \rho_\eps^\pm+\tilde \theta_\eps^\pm\right)\operatorname{Id}
			-
			\tilde u_\eps^\pm\otimes \tilde u_\eps^\pm +\frac{\left| \tilde u_\eps^\pm\right|^2}{3}\operatorname{Id}
			+ \int_{\mathbb{R}^3\times\mathbb{R}^3\times\mathbb{S}^2}
			\hat q_\eps^\pm \tilde \phi MM_*dvdv_*d\sigma\right|
			\\
			& \leq C
			\left\| g_\eps^\pm \gamma_\eps^\pm - \bar g^\pm \right\|_{ L^2\left(Mdv\right)}^2
			+ o(1)_{L^1_\mathrm{loc}(dtdx)},
		\end{aligned}
	\end{equation*}
	and
	\begin{equation*}
		\begin{aligned}
			& \left| F_\eps^\pm \left(\frac{|v|^2}{2}-\frac 52\right)
			-
			\frac 52 \tilde u_\eps^\pm \tilde \theta_\eps^\pm
			+ \int_{\mathbb{R}^3\times\mathbb{R}^3\times\mathbb{S}^2}
			\hat q_\eps^\pm \tilde \psi MM_*dvdv_*d\sigma\right|
			\\
			& \leq C
			\left\| g_\eps^\pm \gamma_\eps^\pm - \bar g^\pm \right\|_{ L^2\left(Mdv\right)}^2
			+ o(1)_{L^1_\mathrm{loc}(dtdx)},
		\end{aligned}
	\end{equation*}
	where $\tilde \phi, \tilde \psi \in L^2\left(Mdv\right)$ are the kinetic momentum and energy fluxes defined by \eqref{phi-psi-def} and \eqref{phi-psi-def inverses}.
\end{lem}

\begin{proof}
	Flux terms are strictly identical to those handled in Lemma \ref{flux1} for the one species case, so that we can reproduce essentially the same arguments.
	
	Thus, we notice first that, modulo the diagonal term in the momentum flux
	\begin{equation*}
		\frac1{\eps} \int_{\mathbb{R}^3}
		g_\eps^\pm \gamma_\eps^\pm \chi\left( {|v|^2\over K_\eps}\right)
		\frac{|v|^2}{3}
		Mdv = \frac1\eps (\tilde \rho_\eps^\pm +\tilde \theta_\eps^\pm),
	\end{equation*}
	the flux terms have the following structure
	\begin{equation*}
		\tilde F_\eps^\pm (\zeta)=\frac1\eps \int_{\mathbb{R}^3}
		g_\eps^\pm \gamma_\eps^\pm \zeta \chi\left( {|v|^2\over K_\eps}\right) M dv,
	\end{equation*}
	where $\zeta \in \Ker (\cL)^\perp\subset L^2\left(Mdv\right)$. Indeed, it is readily seen that the kinetic fluxes $\phi(v)$ and $\psi(v)$, defined by \eqref{phi-psi-def}, are orthogonal to collision invariants.

	Then, reproducing the decomposition \eqref{flux1-decomposition1} from the proof of Lemma \ref{flux1}, we find
	\begin{equation}\label{decomposition fluxes two species}
		\begin{aligned}
			\tilde F_\eps^\pm (\zeta)
			-
			\frac12\int_{\mathbb{R}^3} \left(\Pi \hat g_\eps^\pm\right)^2 \zeta M dv
			& + \int_{\mathbb{R}^3\times\mathbb{R}^3\times\mathbb{S}^2}
			\hat q_\eps^\pm \tilde \zeta MM_*dvdv_*d\sigma
			\\
			& = F_\eps^{\pm 1}(\zeta)+F_\eps^{\pm 2}(\zeta)+F_\eps^{\pm 3}(\zeta)+F_\eps^{\pm 4}(\zeta),
		\end{aligned}
	\end{equation}
	where $\zeta=\mathcal{L}\tilde\zeta$ and
	\begin{equation*}
		\begin{aligned}
			F_\eps^{\pm 1}(\zeta) & =
			\frac14 \int_{\mathbb{R}^3} \left(\hat g_\eps^{\pm 2} -\left(\Pi \hat g_\eps^\pm\right)^2\right) \gamma_\eps^\pm \chi\left( {|v|^2\over K_\eps} \right)\zeta M dv,
			\\
			F_\eps^{\pm 2}(\zeta) & =
			\frac14 \int_{\mathbb{R}^3} \left(\gamma_\eps^\pm \chi\left( {|v|^2\over K_\eps}\right)-1\right)\left(\Pi \hat g_\eps^\pm\right)^2 \zeta M dv,
			\\
			F_\eps^{\pm 3}(\zeta) & =
			\frac1\eps \int_{\mathbb{R}^3} \hat g_\eps^\pm \left(\gamma_\eps^\pm \chi\left( {|v|^2\over K_\eps}\right)-1 \right) \zeta M dv,
			\\
			F_\eps^{\pm 4}(\zeta) & =
			\frac 14
			\int_{\mathbb{R}^3}
			\mathcal{Q}\left(\hat g_\eps^\pm - \Pi \hat g_\eps^\pm, \hat g_\eps^\pm + \Pi \hat g_\eps^\pm \right)
			\tilde \zeta M dv
			\\
			& +
			\frac 14
			\int_{\mathbb{R}^3}
			\mathcal{Q}\left(\hat g_\eps^\pm + \Pi \hat g_\eps^\pm, \hat g_\eps^\pm - \Pi \hat g_\eps^\pm \right)
			\tilde \zeta M dv.
		\end{aligned}
	\end{equation*}

	The remainder terms $F_\eps^{\pm 1}(\zeta)$, $F_\eps^{\pm 2}(\zeta)$, $F_\eps^{\pm 3}(\zeta)$ and $F_\eps^{\pm 4}(\zeta)$ cannot be handled here as in Lemma \ref{flux1} and do not necessarily vanish, because of the lack of equi-integrability of $\left|\hat g_\eps^\pm\right|^2$. Note, however, that the estimates \eqref{eq-use3}, \eqref{eq-use5}, \eqref{eq-use7}, \eqref{eq-use8} and \eqref{eq-use6} can be reproduced here without difficulty, which yields, for any $2<p<4$ and $1<q<\infty$,
	\begin{equation*}
		\begin{aligned}
			\left| F_\eps^{\pm 1}(\zeta)\right|
			& \leq C
			\left(\left\| \hat g_\eps^\pm \right\|_{L^2\left(Mdv\right)}+\left\| \hat q_\eps^\pm \right\|_{L^2\left(MM_*dvdv_*d\sigma\right)}^\frac 12\right)
			\left\| \hat g_\eps^\pm -\Pi \hat g_\eps^\pm \right\|_{ L^2\left(Mdv\right)}
			\\
			& \leq C
			\left\| \hat g_\eps^\pm - \bar g^\pm \right\|_{L^2\left(Mdv\right)}
			\left\| \hat g_\eps^\pm -\Pi \hat g_\eps^\pm \right\|_{ L^2\left(Mdv\right)}
			\\
			& +C
			\left(\left\| \bar g^\pm \right\|_{L^2\left(Mdv\right)}
			+\left\| \hat q_\eps^\pm \right\|_{L^2\left(MM_*dvdv_*d\sigma\right)}^\frac 12
			\right)
			\left\| \hat g_\eps^\pm -\Pi \hat g_\eps^\pm \right\|_{ L^2\left(Mdv\right)},
			\\
			\left| F_\eps^{\pm 2}(\zeta)\right|
			& \leq C
			\left\|\left(\Pi \hat g_\eps^\pm\right)^2\left(\gamma_\eps^\pm \chi\left( {|v|^2\over K_\eps}\right)-1\right)\right\|_{L^q(Mdv)},
			\\
			\left| F_\eps^{\pm 3}(\zeta)\right|
			& \leq C
			\left\| \hat g_\eps^\pm \right\|_{L^2\left(Mdv\right)}
			\left\| \frac{1-\gamma_\eps^\pm}{\eps} \right\|_{ L^p\left(Mdv\right)}
			+C\eps^{\frac K4-1} |\log \eps|^\frac{7}{4}\left\| \hat g_\eps^\pm \right\|_{L^2\left(Mdv\right)}
			\\
			& \leq C
			\left\| \hat g_\eps^\pm -\bar g^\pm \right\|_{L^2\left(Mdv\right)}
			\left\| \frac{1-\gamma_\eps^\pm}{\eps} \right\|_{ L^p\left(Mdv\right)}
			\\
			& + C
			\left\| \bar g^\pm \right\|_{L^2\left(Mdv\right)}
			\left\| \frac{1-\gamma_\eps^\pm}{\eps} \right\|_{ L^p\left(Mdv\right)}
			+C\eps^{\frac K4-1} |\log \eps|^\frac{7}{4}\left\| \hat g_\eps^\pm \right\|_{L^2\left(Mdv\right)},
			\\
			\left| F_\eps^{\pm 4}(\zeta)\right|
			& \leq C
			\left\| \hat g_\eps^\pm \right\|_{L^2\left(Mdv\right)}
			\left\| \hat g_\eps^\pm -\Pi \hat g_\eps^\pm \right\|_{ L^2\left(Mdv\right)}
			\\
			& \leq C
			\left\| \hat g_\eps^\pm - \bar g^\pm \right\|_{L^2\left(Mdv\right)}
			\left\| \hat g_\eps^\pm -\Pi \hat g_\eps^\pm \right\|_{ L^2\left(Mdv\right)}
			\\
			& + C
			\left\| \bar g^\pm \right\|_{L^2\left(Mdv\right)}
			\left\| \hat g_\eps^\pm -\Pi \hat g_\eps^\pm \right\|_{ L^2\left(Mdv\right)}.
		\end{aligned}
	\end{equation*}
	Then, instead of using the convergences \eqref{1-gamma}, \eqref{eq-relaxation} and \eqref{eq-use4} (which are not valid here), we employ the combinations of controls \eqref{bound 4}-\eqref{bound 10}, \eqref{bound 5}-\eqref{bound 11} and \eqref{bound 6}, respectively, provided by Lemmas \ref{trunc-lem} and \ref{trunc-lem 2}, which yields
	\begin{equation*}
		\begin{aligned}
			\left| F_\eps^{\pm 1}(\zeta)\right|
			& \leq C
			\left\| \hat g_\eps^\pm - \bar g^\pm \right\|_{L^2\left(Mdv\right)}
			\left\| g_\eps^\pm \gamma_\eps^\pm - \bar g^\pm \right\|_{ L^2\left(Mdv\right)}
			+ o(1)_{L^1_\mathrm{loc}(dtdx)},
			\\
			\left| F_\eps^{\pm 2}(\zeta)\right|
			& \leq C
			\left\| g_\eps^\pm \gamma_\eps^\pm - \bar g^\pm \right\|_{ L^2\left(Mdv\right)}^2 + o(1)_{L^1_\mathrm{loc}(dtdx)},
			\\
			\left| F_\eps^{\pm 3}(\zeta)\right|
			& \leq C
			\left\| \hat g_\eps^\pm - \bar g^\pm \right\|_{L^2\left(Mdv\right)}
			\left\| g_\eps^\pm \gamma_\eps^\pm - \bar g^\pm \right\|_{ L^2\left(Mdv\right)} + o(1)_{L^1_\mathrm{loc}(dtdx)},
			\\
			\left| F_\eps^{\pm 4}(\zeta)\right|
			& \leq C
			\left\| \hat g_\eps^\pm - \bar g^\pm \right\|_{L^2\left(Mdv\right)}
			\left\| g_\eps^\pm \gamma_\eps^\pm - \bar g^\pm \right\|_{ L^2\left(Mdv\right)} + o(1)_{L^1_\mathrm{loc}(dtdx)}.
		\end{aligned}
	\end{equation*}
	
	On the whole, using \eqref{bound 3} and then incorporating the preceding estimates into \eqref{decomposition fluxes two species}, we obtain
	\begin{equation*}
		\begin{aligned}
			\left|\tilde F_\eps^\pm (\zeta)
			-
			\frac12\int_{\mathbb{R}^3} \left(\Pi \hat g_\eps^\pm\right)^2 \zeta M dv
			+ \int_{\mathbb{R}^3\times\mathbb{R}^3\times\mathbb{S}^2}
			\hat q_\eps^\pm \tilde \zeta MM_*dvdv_*d\sigma\right| &
			\\
			\leq C
			\left\| g_\eps^\pm \gamma_\eps^\pm - \bar g^\pm \right\|_{ L^2\left(Mdv\right)}^2
			& + o(1)_{L^1_\mathrm{loc}(dtdx)},
		\end{aligned}
	\end{equation*}
	which, when further combined with the direct computation \eqref{square infty maxwell}, yields
	\begin{equation*}
		\begin{aligned}
			\left|\tilde F_\eps^\pm (\phi)
			-
			\left(\hat u_\eps^\pm\otimes\hat u_\eps^\pm -\frac{\left|\hat u_\eps^\pm\right|^2}{3}\operatorname{Id}\right)
			+ \int_{\mathbb{R}^3\times\mathbb{R}^3\times\mathbb{S}^2}
			\hat q_\eps^\pm \tilde \phi MM_*dvdv_*d\sigma\right| &
			\\
			\leq C
			\left\| g_\eps^\pm \gamma_\eps^\pm - \bar g^\pm \right\|_{ L^2\left(Mdv\right)}^2
			& + o(1)_{L^1_\mathrm{loc}(dtdx)},
			\\
			\left|\tilde F_\eps^\pm (\psi)
			-
			\frac 52 \hat u_\eps^\pm \hat \theta_\eps^\pm
			+ \int_{\mathbb{R}^3\times\mathbb{R}^3\times\mathbb{S}^2}
			\hat q_\eps^\pm \tilde \psi MM_*dvdv_*d\sigma\right| &
			\\
			\leq C
			\left\| g_\eps^\pm \gamma_\eps^\pm - \bar g^\pm \right\|_{ L^2\left(Mdv\right)}^2
			& + o(1)_{L^1_\mathrm{loc}(dtdx)},
		\end{aligned}
	\end{equation*}
	where $\hat u_\eps^\pm$ and $\hat\theta_\eps^\pm$ are, respectively, the bulk velocity and temperature associated with $\hat g_\eps^\pm$.
	
	Finally, employing \eqref{bound 2}, \eqref{bound 3} and \eqref{bound 8}, we easily obtain that
	\begin{equation*}
		\begin{aligned}
			& \left|\left(\tilde u_\eps^\pm\otimes \tilde u_\eps^\pm -\frac{\left|\tilde u_\eps^\pm\right|^2}{3}\operatorname{Id}\right)
			-\left(\hat u_\eps^\pm\otimes\hat u_\eps^\pm -\frac{\left|\hat u_\eps^\pm\right|^2}{3}\operatorname{Id}\right)\right|
			\\
			& \hspace{20mm} \leq C\left|\tilde u_\eps^\pm-\hat u_\eps^\pm\right| \left|\tilde u_\eps^\pm- \bar u\right|
			+C\left|\hat u_\eps^\pm-\bar u\right| \left|\tilde u_\eps^\pm- \hat u_\eps^\pm\right|
			+ C\left|\tilde u_\eps^\pm-\hat u_\eps^\pm\right| \left|\bar u\right|
			\\
			& \hspace{20mm} \leq C
			\left\| g_\eps^\pm \gamma_\eps^\pm - \bar g^\pm \right\|_{ L^2\left(Mdv\right)}^2
			+ o(1)_{L^1_\mathrm{loc}(dtdx)},
			\\
			& \left|\frac 52 \tilde u_\eps^\pm \tilde \theta_\eps^\pm
			-\frac 52 \hat u_\eps^\pm \hat \theta_\eps^\pm\right|
			\\
			& \hspace{20mm} \leq C\left|\tilde u_\eps^\pm-\hat u_\eps^\pm\right| \left|\tilde \theta_\eps^\pm- \bar \theta\right|
			+C\left|\hat u_\eps^\pm-\bar u\right| \left|\tilde \theta_\eps^\pm- \hat \theta_\eps^\pm\right|
			\\
			& \hspace{20mm} + C\left|\tilde u_\eps^\pm-\hat u_\eps^\pm\right| \left|\bar \theta\right|
			+C\left|\bar u\right| \left|\tilde \theta_\eps^\pm- \hat \theta_\eps^\pm\right|
			\\
			& \hspace{20mm} \leq C
			\left\| g_\eps^\pm \gamma_\eps^\pm - \bar g^\pm \right\|_{ L^2\left(Mdv\right)}^2
			+ o(1)_{L^1_\mathrm{loc}(dtdx)}.
		\end{aligned}
	\end{equation*}
	Combining the preceding estimates concludes the proof of the lemma.
\end{proof}

\subsection{Decomposition of acceleration terms}\label{section step 3 two species}

It only remains to deal with the acceleration terms.

\begin{lem}\label{acceleration2}
	The acceleration terms defined by \eqref{A-def 2} satisfy
	\begin{equation*}
		\begin{aligned}
			\left|A_\eps^{\pm }(1)\right|
			& \leq C\delta \left|E_\eps-\bar E\right|
			\left\|g_\eps^\pm \gamma_\eps^\pm -\bar g^\pm\right\|_{L^2\left(Mdv\right)}
			\\
			& +o(1)_{L^1_{\mathrm{loc}}(dtdx)},
			\\
			\left|A_\eps^{\pm }(v) -\frac\delta\eps E_\eps - \delta \tilde \rho_\eps^\pm E_\eps - \frac\delta\eps \tilde u_\eps^\pm \wedge B_\eps\right|
			& \leq C\delta \left|E_\eps-\bar E\right|
			\left\|g_\eps^\pm \gamma_\eps^\pm -\bar g^\pm\right\|_{L^2\left(Mdv\right)}
			\\
			& +o(1)_{L^1_{\mathrm{loc}}(dtdx)},
			\\
			\left|A_\eps^{\pm }\left( \frac{|v|^2}{2}-\frac 52 \right)-\delta \tilde u_\eps^\pm \cdot E_\eps\right|
			& \leq C\delta \left|E_\eps-\bar E\right|
			\left\|g_\eps^\pm \gamma_\eps^\pm -\bar g^\pm\right\|_{L^2\left(Mdv\right)}
			\\
			& +o(1)_{L^1_{\mathrm{loc}}(dtdx)}.
		\end{aligned}
	\end{equation*}
\end{lem}

\begin{proof}
	We follow a strategy similar to the proof of Lemma \ref{acceleration1} in the one species case.

	Thus, by definition of the acceleration terms, we consider the decomposition
	\begin{equation}\label{acceleration decomposition two species}
		\begin{aligned}
			A_\eps^\pm (\varphi)-\frac\delta\eps
			E_\eps \cdot \int_{\mathbb{R}^3}
			\varphi v M dv
			- \frac\delta\eps \int_{\mathbb{R}^3}
			g_\eps^\pm \gamma_\eps^\pm & \left(\eps E_\eps+v\wedge B_\eps\right)
			\cdot \left(\nabla _v \varphi\right)
			\chi\left( {|v|^2\over K_\eps}\right)
			M dv
			\\
			& = A_\eps^{\pm 1}\left(\varphi \right) + A_\eps^{\pm 2}(\varphi) + A_\eps^{\pm 3}(\varphi),
		\end{aligned}
	\end{equation}
	with
	\begin{equation*}
		\begin{aligned}
			A_\eps^{\pm 1}(\varphi)
			& = -\frac\delta\eps E_\eps \cdot \int_{\mathbb{R}^3}
			\varphi(v)(1-\chi)\left( {|v|^2\over K_\eps}\right)vMdv,
			\\
			A_\eps^{\pm 2}(\varphi)
			& = \delta E_\eps \cdot \int_{\mathbb{R}^3}
			\left( g_\eps^\pm \hat \gamma_\eps^\pm
			- {1-\hat \gamma_\eps^\pm \over \eps}
			- g_\eps^\pm \gamma_\eps^\pm\right) \varphi(v) \chi\left( {|v|^2\over K_\eps}\right)
			v M dv,
			\\
			A_\eps^{\pm 3} (\varphi)
			& ={2\delta \over K_\eps} E_\eps \cdot
			\int_{\mathbb{R}^3}
			g_\eps^\pm \gamma_\eps^\pm
			\varphi(v) \chi' \left( {|v|^2\over K_\eps}\right)
			v M dv.
		\end{aligned}
	\end{equation*}

	Then, we estimate the remainders $A^{\pm 1}_\eps(\varphi)$ and $A^{\pm 3}_\eps(\varphi)$ exactly as $A^{1}_\eps(\varphi)$ and $A^{3}_\eps(\varphi)$ in the one species case. It yields that
	\begin{equation*}
		A^{\pm 1}_\eps(\varphi), A^{\pm 3}_\eps(\varphi) \to 0
		\text{ in }L^1_\mathrm{loc}(dtdx)
		\text{ as }\eps \to 0.
	\end{equation*}
	
	The remaining term $A_\eps^{\pm 2}(\varphi)$ cannot be handled as in Lemma \ref{acceleration1}. % and, in fact, does not necessarily vanish, because of the lack of equi-integrability of $\left|\hat g_\eps^\pm\right|^2$.
	Note, however, that, writing $g_\eps^\pm=\frac 12 \hat g_\eps^\pm\left(1+\sqrt{G_\eps^\pm}\right),$ an easy computation provides
		\begin{equation*}
			\begin{aligned}
				A_\eps^{\pm 2}(\varphi)
				& =\delta E_\eps \cdot \int_{\mathbb{R}^3}
				\left( \eps g_\eps^{\pm 2} \gamma'(G_\eps^\pm)
				- {1-\hat \gamma_\eps^\pm \over \eps} \right)
				\varphi(v) \chi\left( {|v|^2\over K_\eps}\right)
				vM dv
				\\
				& =\delta E_\eps \cdot \int_{\mathbb{R}^3}
				\frac 12 \hat g_\eps^\pm \left(1+\sqrt{G_\eps^\pm}\right) \left(G_\eps^\pm-1\right) \gamma '(G_\eps^\pm)
				\varphi(v) \chi\left( {|v|^2\over K_\eps}\right)
				vM dv
				\\
				& -\delta E_\eps \cdot \int_{\mathbb{R}^3}
				{1-\hat \gamma_\eps^\pm \over \eps}
				\varphi(v) \chi\left( {|v|^2\over K_\eps}\right)
				vM dv.
			\end{aligned}
		\end{equation*}
	
	Then, simply noticing, in view of the hypotheses \eqref{gamma-assumption} on the renormalization, that $(1+\sqrt z)(z-1)\gamma'(z)$ is bounded pointwise and supported on values $z\geq 2$, we deduce
	\begin{equation*}
		\begin{aligned}
			\left|A_\eps^{\pm 2}(\varphi)\right|
			& \leq
			C\delta\left|E_\eps\right|
			\left\|\hat g_\eps^\pm\left(1+\sqrt{G_\eps^\pm}\right)\left(G_\eps^\pm -1\right)\gamma'\left(G_\eps^\pm\right)
			\right\|_{L^2(Mdv)}
			\\
			& +
			C\delta\left|E_\eps\right|
			\left\|{1-\hat \gamma_\eps^\pm \over \eps}
			\right\|_{L^2(Mdv)}
			\\
			& \leq
			C\delta\left|E_\eps\right|
			\left\|\mathds{1}_{\left\{G_\eps^\pm\geq 2\right\}}\hat g_\eps^\pm
			\right\|_{L^2(Mdv)}
			\\
			& \leq
			C\delta\left|E_\eps-\bar E\right|
			\left\|\mathds{1}_{\left\{G_\eps^\pm\geq 2\right\}}\hat g_\eps^\pm
			\right\|_{L^2(Mdv)}
			+
			C\delta\left|\bar E\right|
			\left\|\mathds{1}_{\left\{G_\eps^\pm\geq 2\right\}}\hat g_\eps^\pm
			\right\|_{L^2(Mdv)},
		\end{aligned}
	\end{equation*}
	so that we easily obtain from \eqref{bound 1} in Lemma \ref{trunc-lem} and \eqref{bound 7} in Lemma \ref{trunc-lem 2} that
	\begin{equation*}
		\left|A_\eps^{\pm 2}(\varphi)\right|
		\leq C\delta \left|E_\eps - \bar E\right|
		\left\|g_\eps^\pm \gamma_\eps^\pm -\bar g^\pm\right\|_{L^2\left(Mdv\right)}
		+o(1)_{L^1_{\mathrm{loc}}(dtdx)}.
	\end{equation*}

	Finally, incorporating the preceding remainder estimates into \eqref{acceleration decomposition two species} and performing direct computations of $\int_{\mathbb{R}^3} \varphi v M dv$ and $ \int_{\mathbb{R}^3} g_\eps^\pm \gamma_\eps^\pm (\eps E_\eps+v\wedge B_\eps)\cdot \left(\nabla_v  \varphi\right) \chi\left( {|v|^2\over K_\eps}\right) M dv$ leads then to the expected controls of acceleration terms and concludes the proof of the lemma.
\end{proof}

\subsection{Proof of Proposition \ref{approx3-prop}}\label{proofofprop}

We justify here the validity of the approximate conservation of momentum law provided by Proposition \ref{approx3-prop}.

\begin{proof}[Proof of Proposition \ref{approx3-prop}]
	According to Section \ref{macroscopic defects}, renormalized solutions of the two species Vlasov-Maxwell-Boltzmann system \eqref{VMB2} satisfy the conservation of momentum
	\begin{equation}\label{poynting momentum}
		\begin{aligned}
			\partial_t & \left(\int_{\mathbb{R}^3}\left(g_\eps^++g_\eps^-\right) vMdv
			+E_\eps\wedge B_\eps
			+
			\begin{pmatrix}
				a_{\eps 26}-a_{\eps 35}\\a_{\eps 34}-a_{\eps 16}\\a_{\eps 15}-a_{\eps 24}
			\end{pmatrix}
			\right)
			\\
			& +
			\nabla_x \cdot \left(
			\frac 1\eps\int_{\mathbb{R}^3}\left(g_\eps^++g_\eps^-\right)
			v\otimes v
			M dv + \frac 1{\eps^2}m_\eps
			-E_\eps\otimes E_\eps - e_\eps - B_\eps\otimes B_\eps - b_\eps
			\right)
			\\
			& +\nabla_x\left(\frac{|E_\eps|^2+|B_\eps|^2+\operatorname{Tr}a_\eps}{2}\right)
			= 0.
		\end{aligned}
	\end{equation}

	Next, we decompose
	\begin{equation}\label{decomposition 2}
		\begin{aligned}
			g_\eps^\pm & =
			g_\eps^\pm\gamma_\eps^\pm\chi\left(\frac{|v|^2}{K_\eps}\right)
			+
			g_\eps^\pm\gamma_\eps^\pm\left(1-\chi\left(\frac{|v|^2}{K_\eps}\right)\right)
			+
			g_\eps^\pm\left(1-\gamma_\eps^\pm\right)
			\\
			& =
			g_\eps^\pm\gamma_\eps^\pm\chi\left(\frac{|v|^2}{K_\eps}\right)
			+
			g_\eps^\pm\gamma_\eps^\pm\left(1-\chi\left(\frac{|v|^2}{K_\eps}\right)\right)
			+
			\hat g_\eps^\pm\left(1-\gamma_\eps^\pm\right)
			+
			\frac\eps 4
			\hat g_\eps^{\pm 2}\left(1-\gamma_\eps^\pm\right).
		\end{aligned}
	\end{equation}
	Then, using that $g_\eps^\pm\gamma_\eps^\pm$ is dominated by $\left|\hat g_\eps^\pm\right|$ with the uniform bounds from Lemma \ref{L2-lem} and the control of Gaussian tails \eqref{gaussian-decay 0}, it holds that, for any $p\in\mathbb{R}$,
	\begin{equation*}
		\begin{aligned}
			\int_{\mathbb{R}^3} \left|g_\eps^\pm\gamma_\eps^\pm\left(1-\chi\left(\frac{|v|^2}{K_\eps}\right)\right)\right||v|^\frac p2 Mdv
			& \leq C
			\int_{\left\{|v|^2\geq K_\eps\right\}} \left|\hat g_\eps^\pm\right||v|^\frac p2 Mdv
			\\
			& \leq C
			\left\|\hat g_\eps^\pm\right\|_{L^2\left(Mdv\right)}
			\left(\int_{\left\{|v|^2\geq K_\eps\right\}} |v|^p Mdv\right)^\frac 12
			\\
			& \leq C \left(K|\log \eps|\right)^\frac {p+1}4\eps^{\frac K4}
			\left\|\hat g_\eps^\pm\right\|_{L^2\left(Mdv\right)}
			\\
			& =o(\eps)_{L^\infty\left(dt;L^2\left(dx\right)\right)},
		\end{aligned}
	\end{equation*}
	as soon as $K>4$.
	
	Moreover, since $G_\eps^\pm\geq 2$ implies $\eps\hat g_\eps^\pm\geq 2(\sqrt 2 - 1)$, whence
	\begin{equation*}
		\left|g_\eps^\pm\right|=\left|\hat g_\eps^\pm+\frac\eps 4 \hat g_\eps^{\pm 2}\right|\leq C \eps \hat g_\eps^{\pm 2},
	\end{equation*}
	we find, employing the uniform bounds from Lemmas \ref{L1-lem} and \ref{L2-lem}, that
	\begin{equation*}
		\begin{aligned}
			\int_{\mathbb{R}^3}
			\left|g_\eps^\pm\left(1-\gamma_\eps^\pm\right)\right|
			\left(1+|v|\right)Mdv
			& \leq C\eps^\frac 12
			\int_{\mathbb{R}^3}
			\left|g_\eps^\pm\right|^\frac 12\left|\hat g_\eps^\pm\right|
			\left(1+|v|\right)Mdv
			\\
			& \leq
			C\eps^\frac 12 \left\|g_\eps^\pm\right\|_{L^1\left(\left(1+|v|^2\right)Mdv\right)}^\frac 12
			\left\|\hat g_\eps^\pm\right\|_{L^2\left(Mdv\right)}
			\\
			& =
			o(1)_{L^\infty\left(dt;L^1_{\mathrm{loc}}\left(dx\right)\right)}.
		\end{aligned}
	\end{equation*}
	Alternately, using Lemma \ref{v2-int}, we obtain
	\begin{equation*}
		\begin{aligned}
			\int_{\mathbb{R}^3}
			\left|
			\hat g_\eps^\pm\left(1-\gamma_\eps^\pm\right)
			+
			\frac\eps 4
			\hat g_\eps^{\pm 2}\left(1-\gamma_\eps^\pm\right)
			\right|
			\left(1+|v|\right)Mdv
			& \leq
			C\eps\left\|\hat g_\eps^\pm\right\|_{L^2\left((1+|v|)Mdv\right)}^2
			\\
			& =
			o(1)_{L^1_{\mathrm{loc}}\left(dtdx\right)}.
		\end{aligned}
	\end{equation*}
	If, instead of Lemma \ref{v2-int}, one applies Lemma \ref{trunc-lem 3}, then one finds that
	\begin{equation*}
		\begin{aligned}
			\frac 1\eps & \left\|\int_{\mathbb{R}^3}
			\left|
			\hat g_\eps^\pm\left(1-\gamma_\eps^\pm\right)
			+
			\frac\eps 4
			\hat g_\eps^{\pm 2}\left(1-\gamma_\eps^\pm\right)
			\right|
			\left(1+|v|^2\right)Mdv\right\|_{L^1_\mathrm{loc}(dx)}
			\\
			& \leq
			C \left\|\mathds{1}_{\left\{G_\eps^\pm\geq 2\right\}}\hat g_\eps^\pm\right\|_{L^2_\mathrm{loc}\left(dx;L^2\left((1+|v|^2)Mdv\right)\right)}^2
			\\
			& \leq
			C_1\int_{\mathbb{R}^3\times\mathbb{R}^3} \left(\frac 1{\eps^2}{h\left(\eps g_\eps^\pm\right)}
			- \frac{1}{2} \left( g_\eps^\pm \gamma_\eps^\pm\chi \left({|v|^2\over K_\eps}\right) \right)^2 \right)Mdxdv
			\\
			& + C_2\left\|g_\eps^\pm \gamma_\eps^\pm\chi \left({|v|^2\over K_\eps}\right) -\bar g^\pm\right\|_{L^2\left(Mdxdv\right)}^2
			+ o(1)_{L^1_{\mathrm{loc}}\left(dt\right)}.
		\end{aligned}
	\end{equation*}
	
	Thus, combining the preceding estimates with the decomposition \eqref{decomposition 2}, we arrive at
	\begin{equation}\label{densities approx}
		\int_{\mathbb{R}^3}g_\eps^\pm \varphi(v) Mdv =
		\int_{\mathbb{R}^3}
		g_\eps^\pm\gamma_\eps^\pm\chi\left(\frac{|v|^2}{K_\eps}\right) \varphi(v) Mdv
		+
		o(1)_{L^\infty\left(dt;L^1_{\mathrm{loc}}(dx)\right)},
	\end{equation}
	for all $\varphi(v)$ such that $\frac{\varphi(v)}{1+|v|}\in L^\infty(dv)$, and
	\begin{equation*}
		\begin{aligned}
			\frac 1\eps & \left\|\int_{\mathbb{R}^3}
			\left(g_\eps^\pm
			-g_\eps^\pm\gamma_\eps^\pm\chi\left(\frac{|v|^2}{K_\eps}\right)\right) \varphi(v)
			Mdv\right\|_{L^1_\mathrm{loc}(dx)}
			\\
			& \leq
			C_1\int_{\mathbb{R}^3\times\mathbb{R}^3} \left(\frac 1{\eps^2}{h\left(\eps g_\eps^\pm\right)}
			- \frac{1}{2} \left( g_\eps^\pm \gamma_\eps^\pm\chi \left({|v|^2\over K_\eps}\right) \right)^2 \right)Mdxdv
			\\
			& + C_2\left\|g_\eps^\pm \gamma_\eps^\pm\chi \left({|v|^2\over K_\eps}\right) -\bar g^\pm\right\|_{L^2\left(Mdxdv\right)}^2
			+
			o(1)_{L^1_{\mathrm{loc}}(dt)},
		\end{aligned}
	\end{equation*}
	for all $\varphi(v)$ such that $\frac{\varphi(v)}{1+|v|^2}\in L^\infty(dv)$, which, when incorporated into \eqref{poynting momentum}, yields the approximate conservation law
	\begin{equation*}
		\begin{aligned}
			\partial_t & \left(\int_{\mathbb{R}^3}\left(g_\eps^+\gamma_\eps^++g_\eps^-\gamma_\eps^-\right)
			\chi\left(\frac{|v|^2}{K_\eps}\right) vMdv
			+E_\eps\wedge B_\eps
			+
			\begin{pmatrix}
				a_{\eps 26}-a_{\eps 35}\\a_{\eps 34}-a_{\eps 16}\\a_{\eps 15}-a_{\eps 24}
			\end{pmatrix}
			\right)
			\\
			& +
			\nabla_x \cdot \left(
			\frac 1\eps\int_{\mathbb{R}^3}\left(g_\eps^+\gamma_\eps^++g_\eps^-\gamma_\eps^-\right)
			\chi\left(\frac{|v|^2}{K_\eps}\right)
			v\otimes v
			M dv + \frac 1{\eps^2}m_\eps
			\right)
			\\
			& -
			\nabla_x \cdot \left(
			E_\eps\otimes E_\eps + e_\eps + B_\eps\otimes B_\eps + b_\eps
			\right)
			+\nabla_x\left(\frac{|E_\eps|^2+|B_\eps|^2+\operatorname{Tr}a_\eps}{2}\right)
			\\
			& = \partial_t \left(o(1)_{L^\infty\left(dt;L^1_{\mathrm{loc}}(dx)\right)}\right) + \widetilde R_{\eps},
		\end{aligned}
	\end{equation*}
	where the remainder $\widetilde R_{\eps}$ satisfies
	\begin{equation}\label{tilde remainder}
		\begin{aligned}
			\left\|\widetilde R_{\eps}\right\|_{W^{-1,1}_\mathrm{loc}(dx)}
			& \leq
			C_1\int_{\mathbb{R}^3\times\mathbb{R}^3} \left(\frac 1{\eps^2}{h\left(\eps g_\eps^+\right)}
			- \frac{1}{2} \left( g_\eps^+ \gamma_\eps^+\chi \left({|v|^2\over K_\eps}\right) \right)^2 \right)Mdxdv
			\\
			& +
			C_1\int_{\mathbb{R}^3\times\mathbb{R}^3} \left(\frac 1{\eps^2}{h\left(\eps g_\eps^-\right)}
			- \frac{1}{2} \left( g_\eps^- \gamma_\eps^-\chi \left({|v|^2\over K_\eps}\right) \right)^2 \right)Mdxdv
			\\
			& + C_2
			\left\| \left(g_\eps ^+ \gamma_\eps^+\chi \left({|v|^2\over K_\eps}\right) -\bar g^+,
			g_\eps ^- \gamma_\eps^-\chi \left({|v|^2\over K_\eps}\right) -\bar g^-\right)\right\|_{L^2(Mdxdv)}^2
			\\
			& + o(1)_{L^1_{\mathrm{loc}}(dt)}.
		\end{aligned}
	\end{equation}

	Then, expressing the flux terms above with Lemma \ref{flux2}, we find that
	\begin{equation*}
		\begin{aligned}
			\partial_t & \left(2\tilde u_\eps
			+E_\eps\wedge B_\eps
			+
			\begin{pmatrix}
				a_{\eps 26}-a_{\eps 35}\\a_{\eps 34}-a_{\eps 16}\\a_{\eps 15}-a_{\eps 24}
			\end{pmatrix}
			\right)
			\\
			& +
			\nabla_x \cdot \left( \tilde u_\eps^+ \otimes \tilde u_\eps^+ + \tilde u_\eps^- \otimes \tilde u_\eps^-
			-\frac{\left|\tilde u_\eps^+\right|^2+\left|\tilde u_\eps^-\right|^2}{3} \operatorname{Id}
			+ \frac 1{\eps^2}m_\eps
			\right)
			\\
			& -\nabla_x\cdot\left( \int_{\mathbb{R}^3\times\mathbb{R}^3\times\mathbb{S}^2} \left(\hat q_\eps^++\hat q_\eps^-\right) \tilde \phi MM_* dvdv_*d\sigma
			\right)
			\\
			& -\nabla_x\cdot\left(
			E_\eps\otimes E_\eps + e_\eps + B_\eps\otimes B_\eps + b_\eps
			\right)
			+\nabla_x\left(\frac{|E_\eps|^2+|B_\eps|^2+\operatorname{Tr}a_\eps}{2}\right)
			\\
			& = - \frac 2\eps \nabla_x\left(\tilde \rho_\eps+\tilde \theta_\eps\right)+\partial_t \left(o(1)_{L^\infty\left(dt;L^1_{\mathrm{loc}}(dx)\right)}\right)
			+ \bar R_{\eps},
		\end{aligned}
	\end{equation*}
	where the remainder $\bar R_{\eps}$ also satisfies \eqref{tilde remainder}.
	
	Finally, an application of the estimate \eqref{nonlinear modulation} concludes the proof of the proposition.
\end{proof}

\subsection{Proofs of Lemmas \ref{trunc-lem}, \ref{trunc-lem 2}, \ref{trunc-lem 3} and \ref{trunc-lem 4}}\label{proof of lemma}

At last, we provide a complete justification of Lemmas \ref{trunc-lem}, \ref{trunc-lem 2}, \ref{trunc-lem 3} and \ref{trunc-lem 4}.

	\begin{proof}[Proof of Lemma \ref{trunc-lem}]
		This lemma hinges upon the simple fact that the renormalization $\Gamma(z)$ enjoys the suitable bound from below \eqref{gamma below}. In terms of the renormalized fluctuations, this bound implies that
		\begin{equation*}
			\frac 1C\left|\hat g_\eps^\pm\right|\leq \left|g_\eps^\pm\gamma_\eps^\pm\right|\leq C \left|\hat g_\eps^\pm\right|,
		\end{equation*}
		for some $C>1$, which will be used repeatedly throughout the present proof.
		
		\noindent$\bullet$ In order to establish the first bound \eqref{bound 1}, notice that, since $G_\eps^\pm\geq 2$ implies $\eps\hat g_\eps^\pm\geq 2(\sqrt 2 - 1)$,
		\begin{equation*}
			\begin{aligned}
				\left|\mathds{1}_{\left\{G_\eps^\pm\geq 2\right\}}\hat g_\eps^\pm\right|
				& \leq
				C\left|\mathds{1}_{\left\{G_\eps^\pm\geq 2\right\}}\gamma_\eps^\pm g_\eps^\pm\right|
				\\
				& \leq
				C\left|g_\eps^\pm \gamma_\eps^\pm -\bar g^\pm\right|
				+
				C\left|\mathds{1}_{\left\{G_\eps^\pm\geq 2\right\}} \bar g^\pm\right|
				\\
				& \leq
				C\left|g_\eps^\pm \gamma_\eps^\pm -\bar g^\pm\right|
				+
				C\eps\left| \hat g_\eps^\pm\right|\left| \bar g^\pm\right|,
			\end{aligned}
		\end{equation*}
		whence
		\begin{equation*}
			\begin{aligned}
				\left\|\mathds{1}_{\left\{G_\eps^\pm\geq 2\right\}}\hat g_\eps^\pm\chi\left(\frac{|v|^2}{K_\eps}\right)\right\|_{L^2\left(Mdv\right)}
				& \leq
				C\left\|g_\eps^\pm \gamma_\eps^\pm -\bar g^\pm\right\|_{L^2\left(Mdv\right)}
				\\
				& +
				C\eps K_\eps \left\| \hat g_\eps^\pm\right\|_{L^2\left(Mdv\right)}
				\left\| \frac{\bar g^\pm}{1+|v|^2}\right\|_{L^\infty\left(dtdxdv\right)}.
			\end{aligned}
		\end{equation*}
		Moreover, it is readily seen that
		\begin{equation*}
			\left\|\mathds{1}_{\left\{G_\eps^\pm\geq 2\right\}}\hat g_\eps^\pm\left(1-\chi\left(\frac{|v|^2}{K_\eps}\right)\right)\right\|_{L^2\left(Mdv\right)}
			\leq \frac{C}{K_\eps^\frac 12}
			\left\|\hat g_\eps^\pm\right\|_{L^2\left((1+|v|^2)Mdv\right)}.
		\end{equation*}

		Therefore, it follows that, combining the preceding estimates and considering the uniform bounds from Lemmas \ref{L2-lem} and \ref{v2-int},
		\begin{equation*}
			\begin{aligned}
				\left\|\mathds{1}_{\left\{G_\eps^\pm\geq 2\right\}}\hat g_\eps^\pm\right\|_{L^2\left(Mdv\right)}
				& \leq
				C\left\|g_\eps^\pm \gamma_\eps^\pm -\bar g^\pm\right\|_{L^2\left(Mdv\right)}
				\\
				& +
				O\left(\eps |\log\eps|\right)_{L^\infty\left(dt;L^2(dx)\right)}
				+ O\left(\frac{1}{|\log\eps|^\frac 12}\right)_{L^2_\mathrm{loc}(dtdx)},
			\end{aligned}
		\end{equation*}
		which concludes the proof of \eqref{bound 1}.

		\noindent$\bullet$ To deduce the second bound \eqref{bound 2}, we decompose, writing $g_\eps^\pm=\hat g_\eps^\pm\left(1+\frac \eps 4\hat g_\eps^\pm\right)$ and using that $\left(1+\frac \eps 4\hat g_\eps^\pm\right)\gamma_\eps^\pm$ is uniformly bounded pointwise,
		\begin{equation*}
			\begin{aligned}
				\left|g_\eps^\pm\gamma_\eps^\pm\chi\left(\frac{|v|^2}{K_\eps}\right)-\hat g_\eps^\pm\right|
				& =
				\left|\hat g_\eps^\pm\left(\left(1+\frac \eps 4\hat g_\eps^\pm\right)\gamma_\eps^\pm\chi\left(\frac{|v|^2}{K_\eps}\right)-1\right)\right|
				\\
				& \leq C
				\left|g_\eps^\pm\gamma_\eps^\pm\left(\left(1+\frac \eps 4\hat g_\eps^\pm\right)\gamma_\eps^\pm\chi\left(\frac{|v|^2}{K_\eps}\right)-1\right)\right|
				\\
				& \leq
				C\left|g_\eps^\pm\gamma_\eps^\pm - \bar g^\pm\right|
				\\
				& + C\left|\bar g^\pm\left(\left(1+\frac \eps 4\hat g_\eps^\pm\right)\gamma_\eps^\pm\chi\left(\frac{|v|^2}{K_\eps}\right)-1\right)\right|
				\\
				& \leq
				C\left|g_\eps^\pm\gamma_\eps^\pm - \bar g^\pm\right|
				+ C\left|\bar g^\pm\left(1-\chi\left(\frac{|v|^2}{K_\eps}\right)\right)\right|
				\\
				& + C\eps\left|\bar g^\pm\hat g_\eps^\pm \gamma_\eps^\pm\chi\left(\frac{|v|^2}{K_\eps}\right)\right|
				+ C\left|\bar g^\pm\left(\gamma_\eps^\pm-1\right)\chi\left(\frac{|v|^2}{K_\eps}\right)\right|,
			\end{aligned}
		\end{equation*}
		which implies
		\begin{equation*}
			\begin{aligned}
				& \left\|g_\eps^\pm\gamma_\eps^\pm\chi\left(\frac{|v|^2}{K_\eps}\right)-\hat g_\eps^\pm\right\|_{L^2(Mdv)}
				\\
				& \hspace{10mm} \leq
				C\left\|g_\eps^\pm\gamma_\eps^\pm - \bar g^\pm\right\|_{L^2(Mdv)}
				+ C\eps K_\eps\left\|\hat g_\eps^\pm\right\|_{L^2(Mdv)}\left\| \frac{\bar g^\pm}{1+|v|^2}\right\|_{L^\infty\left(dtdxdv\right)}
				\\
				& \hspace{10mm} + C\left\|(1+|v|^2)\left(1-\chi\left(\frac{|v|^2}{K_\eps}\right)\right)\right\|_{L^2(Mdv)}\left\| \frac{\bar g^\pm}{1+|v|^2}\right\|_{L^\infty\left(dtdxdv\right)}.
			\end{aligned}
		\end{equation*}

		Then, employing the control of Gaussian tails \eqref{gaussian-decay 0} and the uniform bound from Lemma \ref{L2-lem}, we infer that
		\begin{equation*}
			\begin{aligned}
				& \left\|g_\eps^\pm\gamma_\eps^\pm\chi\left(\frac{|v|^2}{K_\eps}\right)-\hat g_\eps^\pm\right\|_{L^2(Mdv)}
				\\
				& \hspace{10mm} \leq
				C\left\|g_\eps^\pm\gamma_\eps^\pm - \bar g^\pm\right\|_{L^2(Mdv)}
				+ O\left(\eps |\log\eps|\right)_{L^\infty\left(dt;L^2(dx)\right)}
				+ C\eps^\frac K4|\log\eps|^\frac 54,
			\end{aligned}
		\end{equation*}
		which establishes \eqref{bound 2}.

		\noindent$\bullet$ The third bound \eqref{bound 3} easily follows from the estimate
		\begin{equation*}
			\begin{aligned}
				& \left\|g_\eps^\pm\gamma_\eps^\pm-\hat g_\eps^\pm\right\|_{L^2(Mdv)}
				\\
				& \leq
				\left\|g_\eps^\pm\gamma_\eps^\pm\chi\left(\frac{|v|^2}{K_\eps}\right)-\hat g_\eps^\pm\right\|_{L^2(Mdv)}
				+
				\left\| \hat g_\eps^\pm\left(1-\chi\left(\frac{|v|^2}{K_\eps}\right)\right)\right\|_{L^2(Mdv)}
				\\
				& \leq
				\left\|g_\eps^\pm\gamma_\eps^\pm\chi\left(\frac{|v|^2}{K_\eps}\right)-\hat g_\eps^\pm\right\|_{L^2(Mdv)}
				+\frac{C}{K_\eps^\frac 12}
				\left\| \hat g_\eps^\pm \right\|_{L^2\left(\left(1+|v|^2\right)Mdv\right)},
			\end{aligned}
		\end{equation*}
		which, when combined with the second bound \eqref{bound 2}, concludes its justification.

		\noindent$\bullet$ The justification of \eqref{bound 5} is simple. Since $\Pi \bar g^\pm=\bar g^\pm$, we easily estimate
		\begin{equation*}
			\begin{aligned}
				\left\| \hat g_\eps^\pm -\Pi \hat g_\eps^\pm \right\|_{L^2( Mdv)}
				& \leq C
				\left\| \hat g_\eps^\pm - \bar g^\pm\right\|_{L^2( Mdv)}
				+C
				\left\| \Pi \bar g^\pm-\Pi \hat g_\eps^\pm \right\|_{L^2( Mdv)}
				\\
				& \leq C
				\left\| \hat g_\eps^\pm - \bar g^\pm\right\|_{L^2( Mdv)}.
			\end{aligned}
		\end{equation*}
		Therefore, the bound \eqref{bound 5} is obtained by combining the preceding control with \eqref{bound 3}.

		\noindent$\bullet$ We focus now on \eqref{bound 6}. We first easily find that
		\begin{equation*}
			\begin{aligned}
				& \left\|\left(\Pi \hat g_\eps^\pm\right)^2\left(\gamma_\eps^\pm\chi\left(\frac{|v|^2}{K_\eps}\right)-1\right)\right\|_{L^q(Mdv)}
				\\
				& \leq
				C\left\|\Pi \left(\hat g_\eps^\pm-\bar g^\pm\right)\right\|_{L^{2q}(Mdv)}^2
				+
				C\left\|\left(\bar g^\pm\right)^2 \left(\gamma_\eps^\pm\chi\left(\frac{|v|^2}{K_\eps}\right)-1\right)\right\|_{L^q(Mdv)}
				\\
				& \leq
				C\left\|\hat g_\eps^\pm-\bar g^\pm\right\|_{L^{2}(Mdv)}^2
				\\
				& +
				C\left\| \frac{\bar g^\pm}{1+|v|^2} \right\|_{L^\infty(dtdxdv)}^2
				\left\|\left(1+|v|^4\right)\chi\left(\frac{|v|^2}{K_\eps}\right)\left(\gamma_\eps^\pm-1\right)\right\|_{L^{q}(Mdv)}
				\\
				& +
				C\left\| \frac{\bar g^\pm}{1+|v|^2} \right\|_{L^\infty(dtdxdv)}^2
				\left\|\left(1+|v|^4\right)\left(\chi\left(\frac{|v|^2}{K_\eps}\right)-1\right)\right\|_{L^{q}(Mdv)}.
			\end{aligned}
		\end{equation*}
		Therefore, utilizing the control of Gaussian tails \eqref{gaussian-decay 0} and the fact that $G_\eps^\pm\geq 2$ on the support of $\gamma_\eps^\pm-1$, we deduce that
		\begin{equation*}
			\begin{aligned}
				& \left\|\left(\Pi \hat g_\eps^\pm\right)^2\left(\gamma_\eps^\pm\chi\left(\frac{|v|^2}{K_\eps}\right)-1\right)\right\|_{L^q(Mdv)}
				\\
				& \leq
				C\left\|\hat g_\eps^\pm-\bar g^\pm\right\|_{L^{2}(Mdv)}^2
				+
				CK_\eps^2\eps^\frac 2q\left\|\hat g_\eps^\pm\right\|_{L^2(Mdv)}^{\frac 2q}
				+
				C |\log\eps|^{2+\frac 1{2q}}\eps^{\frac{K}{2q}},
			\end{aligned}
		\end{equation*}
		which, when combined with \eqref{bound 3}, concludes the proof of \eqref{bound 6}.

		\noindent$\bullet$ Next, we establish the last bound \eqref{bound 4}. Note first that the case $p=2$ is easily deduced from \eqref{bound 1}, using again that $G_\eps^\pm\geq 2$ implies $\eps\hat g_\eps^\pm\geq 2(\sqrt 2 - 1)$. Thus, the difficulty here lies in obtaining a gain of velocity integrability.

		To this end, we introduce the macroscopic truncation
		\begin{equation*}% \label{macro cutoff}
			\chi_\eps^\pm = \mathds{1}_{\left\{ \eps \left\| g_\eps^\pm \gamma_\eps^\pm - \bar g^\pm \right \|_{L^2(Mdv)} \leq 1 \right\}}.
		\end{equation*}
		Then, we have
		\begin{equation*}
			\left(1-\chi_\eps^\pm\right)\left\|\frac 1\eps\mathds{1}_{\left\{G_\eps^\pm\geq 2\right\}}\right\|_{L^p(Mdv)}
			\leq C
			\left\| g_\eps^\pm \gamma_\eps - \bar g^\pm \right \|_{L^2(Mdv)}.
		\end{equation*}
		Moreover, controlling Gaussian tails with \eqref{gaussian-decay 0}, it clearly holds that
		\begin{equation*}
			\chi_\eps^\pm \left\|\frac 1\eps\mathds{1}_{\left\{G_\eps^\pm\geq 2\right\}}\left(1-\chi\left(\frac{|v|^2}{K_\eps}\right)\right)\right\|_{L^p(Mdv)}
			\leq C \eps^{\frac{K}{2p}-1}\left|\log\eps\right|^{\frac 1{2p}},
		\end{equation*}
		which is small as soon as $K\geq 8>2p$, so that we only have to control the size of $\frac 1\eps\mathds{1}_{\left\{G_\eps^\pm\geq 2\right\}}\chi\left(\frac{|v|^2}{K_\eps}\right)$ on the support of $\chi_\eps^\pm$.
		
		Thus, employing the decomposition
		\begin{equation*}
			\frac1\eps \mathds{1}_{\left\{G_\eps^\pm\geq 2\right\}}
			\leq \frac 32 \mathds{1}_{\left\{G_\eps^\pm\geq 2\right\}} \left|\hat g_\eps^\pm\right|
			\leq \frac 32 \mathds{1}_{\left\{G_\eps^\pm\geq 2\right\}}\left( \left|\Pi\hat g_\eps^\pm\right|
			+ \left|\hat g_\eps^\pm - \Pi\hat g_\eps^\pm\right| \right) ,
		\end{equation*}
		we find, for any $1<r<2$, the interpolation estimate
		\begin{equation*}
			\begin{aligned}
				\frac1{\eps^2} \mathds{1}_{\left\{G_\eps^\pm\geq 2\right\}} & \chi\left(\frac{|v|^2}{K_\eps}\right)
				\\
				& \leq \frac C{\eps^{2-\frac 2 r}} \mathds{1}_{\left\{G_\eps^\pm\geq 2\right\}}
				\chi\left(\frac{|v|^2}{K_\eps}\right)
				\left(\left|\Pi\hat g_\eps^\pm\right|
				+\left|\hat g_\eps^\pm - \Pi\hat g_\eps^\pm\right|\right)^{\frac 2r}
				\\
				& \leq C\mathds{1}_{\left\{G_\eps^\pm\geq 2\right\}} \left|\hat g_\eps^\pm\right|^{2-\frac 2r} \left|\Pi\hat g_\eps^\pm\right|^{\frac 2r}
				\chi\left(\frac{|v|^2}{K_\eps}\right)
				+ \frac{C}{\eps^{2-\frac 2r}}\left|\hat g_\eps^\pm - \Pi\hat g_\eps^\pm\right|^{\frac 2r}
				\\
				& \leq C\mathds{1}_{\left\{G_\eps^\pm\geq 2\right\}} \left|\hat g_\eps^\pm\right|^{2-\frac 2r} \left|\Pi\left(\hat g_\eps^\pm-\bar g^\pm\right)\right|^{\frac 2r}
				\\ & +
				C\mathds{1}_{\left\{G_\eps^\pm\geq 2\right\}} \left|\hat g_\eps^\pm\right|^{2-\frac 2r} \left|\bar g^\pm\right|^{\frac 2r}
				\chi\left(\frac{|v|^2}{K_\eps}\right)
				+\frac{C}{\eps^{2-\frac 2r}}\left|\hat g_\eps^\pm - \Pi\hat g_\eps^\pm\right|^{\frac 2r}
				\\
				& \leq C\mathds{1}_{\left\{G_\eps^\pm\geq 2\right\}} \left|\hat g_\eps^\pm\right|^{2-\frac 2r} \left|\Pi\left(\hat g_\eps^\pm-\bar g^\pm\right)\right|^{\frac 2r}
				\\ & +
				CK_\eps^\frac 2r \eps^{\frac 4r - 2} \left|\hat g_\eps^\pm\right|^{\frac 2r}
				\left\| \frac{\bar g^\pm}{1+|v|^2}\right\|_{L^\infty\left(dtdxdv\right)}^\frac 2r
				+\frac{C}{\eps^{2-\frac 2r}}\left|\hat g_\eps^\pm - \Pi\hat g_\eps^\pm\right|^{\frac 2r}.
			\end{aligned}
		\end{equation*}
		Therefore, combining the preceding estimate with the relaxation \eqref{relaxation-est} from Lemma \ref{relaxation-control}, we deduce that
		\begin{equation*}
			\begin{aligned}
				& \left\|\frac1{\eps^2} \mathds{1}_{\left\{G_\eps^\pm\geq 2\right\}}\chi\left(\frac{|v|^2}{K_\eps}\right)\right\|_{L^r\left(Mdv\right)}
% 				& \leq
% 				C\left\|\mathds{1}_{\left\{G_\eps^\pm\geq 2\right\}} \left|\hat g_\eps^\pm\right|^{2-\frac 2r} \left|\Pi\left(\hat g_\eps^\pm-g^\pm\right)\right|^{\frac 2r}\right\|_{L^r(Mdv)}
% 				\\
% 				& +
% 				C\left|\log\eps\right|^\frac 2r \eps^{\frac 4r - 2}\left\| \hat g_\eps^\pm \right\|_{L^2(Mdv)}^\frac 2r
% 				+\frac{C}{\eps^{2-\frac 2r}}\left\|\hat g_\eps^\pm - \Pi\hat g_\eps^\pm\right\|_{L^2(Mdv)}^{\frac 2r}
% 				\\
				\\
				& \leq
				C\left\|\mathds{1}_{\left\{G_\eps^\pm\geq 2\right\}} \hat g_\eps^\pm \right\|_{L^2(Mdv)}^{2-\frac{2}{r}}
				\left\| \Pi\left(\hat g_\eps^\pm-\bar g^\pm\right)\right\|_{L^\frac{2}{2-r}(Mdv)}^\frac{2}{r}
				\\
				& +
				C\left|\log\eps\right|^\frac 2r \eps^{\frac 4r - 2}\left\| \hat g_\eps^\pm \right\|_{L^2(Mdv)}^\frac 2r
				+\frac{C}{\eps^{2-\frac 2r}}\left\|\hat g_\eps^\pm - \Pi\hat g_\eps^\pm\right\|_{L^2(Mdv)}^{\frac 2r}
				\\
				& \leq
				C\left\|\mathds{1}_{\left\{G_\eps^\pm\geq 2\right\}} \hat g_\eps^\pm \right\|_{L^2(Mdv)}^{2-\frac{2}{r}}
				\left\| \hat g_\eps^\pm-\bar g^\pm \right\|_{L^2(Mdv)}^\frac{2}{r}
				\\
				& +
				C\left|\log\eps\right|^\frac 2r \eps^{\frac 4r - 2}\left\| \hat g_\eps^\pm \right\|_{L^2(Mdv)}^\frac 2r
				+\frac{C}{\eps^{2-\frac 2r}}\left(\eps\left\|\hat g_\eps^\pm\right\|_{L^2(Mdv)}^2+O(\eps)_{L^2(dtdx)}\right)^{\frac 2r}
				\\
				& \leq
				C\left\|\mathds{1}_{\left\{G_\eps^\pm\geq 2\right\}} \hat g_\eps^\pm \right\|_{L^2(Mdv)}^2
				+C
				\left\| \hat g_\eps^\pm-\bar g^\pm \right\|_{L^2(Mdv)}^2
				\\
				& +
				C\left|\log\eps\right|^\frac 2r \eps^{\frac 4r - 2}\left\| \hat g_\eps^\pm \right\|_{L^2(Mdv)}^\frac 2r
				+
				C \eps^{\frac 4r - 2}\left\| g_\eps^\pm \gamma_\eps^\pm \right\|_{L^2(Mdv)}^\frac 4r
				+O\left(\eps^{\frac 4r - 2}\right)_{L^r(dtdx)}
				\\
				& \leq
				C\left\|\mathds{1}_{\left\{G_\eps^\pm\geq 2\right\}} \hat g_\eps^\pm \right\|_{L^2(Mdv)}^2
				+C
				\left\| \hat g_\eps^\pm-\bar g^\pm \right\|_{L^2(Mdv)}^2
				\\
				& +
				C\left|\log\eps\right|^\frac 2r \eps^{\frac 4r - 2}\left\| \hat g_\eps^\pm \right\|_{L^2(Mdv)}^\frac 2r
				+C \eps^{\frac 4r - 2}\left\| g_\eps^\pm \gamma_\eps^\pm - \bar g^\pm \right\|_{L^2(Mdv)}^\frac 4r
				\\
				& +C \eps^{\frac 4r - 2}\left\| \bar g^\pm \right\|_{L^2(Mdv)}^\frac 4r
				+O\left(\eps^{\frac 4r - 2}\right)_{L^r(dtdx)}.
			\end{aligned}
		\end{equation*}
		Note that these controls do not yield vanishing remainders in the endpoint case $r=2$, which explains the use of the interpolation parameter $1<r<2$.
		
		Then, recalling that the preceding estimate only needs to be performed on the support of $\chi_\eps^\pm$ and denoting $p=2r$, we infer
		\begin{equation*}
			\begin{aligned}
				& \chi_\eps^\pm\left\|\frac1{\eps} \mathds{1}_{\left\{G_\eps^\pm\geq 2\right\}}\chi\left(\frac{|v|^2}{K_\eps}\right)\right\|_{L^p\left(Mdv\right)}
				\\
				& \leq
				C\left\|\mathds{1}_{\left\{G_\eps^\pm\geq 2\right\}} \hat g_\eps^\pm \right\|_{L^2(Mdv)}
				+C
				\left\| \hat g_\eps^\pm-\bar g^\pm \right\|_{L^2(Mdv)}
				+C \left\| g_\eps^\pm \gamma_\eps^\pm - \bar g^\pm \right\|_{L^2(Mdv)}
				\\
				& +
				C\left|\log\eps\right|^\frac 2p \eps^{\frac 4p - 1}\left\| \hat g_\eps^\pm \right\|_{L^2(Mdv)}^\frac 2p
				+C \eps^{\frac 4p - 1}\left\| \bar g^\pm \right\|_{L^2(Mdv)}^\frac 4p
				+O\left(\eps^{\frac 4p - 1}\right)_{L^p(dtdx)},
			\end{aligned}
		\end{equation*}
		which, when combined with the bounds \eqref{bound 1} and \eqref{bound 3}, concludes the proof of \eqref{bound 4}.

		The proof of the lemma is now complete.
	\end{proof}

	\begin{proof}[Proof of Lemma \ref{trunc-lem 2}]
		This lemma is a simple consequence of the relaxation estimate provided by Lemma \ref{relaxation-control}~:
		\begin{equation*}
			\hat g_\eps^\pm - \Pi \hat g_\eps^\pm = O(\eps)_{L^1_\mathrm{loc}\left(dtdx ; L^2(Mdv)\right)}.
		\end{equation*}
		
		\noindent$\bullet$ In order to establish the first bound \eqref{bound 7}, notice that, since $G_\eps^\pm\geq 2$ implies $\eps\hat g_\eps^\pm\geq 2(\sqrt 2 - 1)$,
		\begin{equation*}
			\begin{aligned}
				\left|\mathds{1}_{\left\{G_\eps^\pm\geq 2\right\}}\hat g_\eps^\pm\right|
				& \leq
				\left|\hat g_\eps^\pm - \Pi \hat g_\eps^\pm\right|
				+
				\left|\mathds{1}_{\left\{G_\eps^\pm\geq 2\right\}} \Pi \hat g_\eps^\pm\right|
				\\
				& \leq
				\left|\hat g_\eps^\pm - \Pi \hat g_\eps^\pm\right|
				+
				C \left|\mathds{1}_{\left\{G_\eps^\pm\geq 2\right\}}\left(1+|v|^2\right)\right| \left\| \hat g_\eps^\pm\right\|_{L^2(Mdv)},
			\end{aligned}
		\end{equation*}
		whence
		\begin{equation*}
			\left\|\mathds{1}_{\left\{G_\eps^\pm\geq 2\right\}}\hat g_\eps^\pm\chi\left(\frac{|v|^2}{K_\eps}\right)\right\|_{L^2\left(Mdv\right)}
			\leq
			\left\|\hat g_\eps^\pm - \Pi \hat g_\eps^\pm\right\|_{L^2\left(Mdv\right)}
			+
			C\eps K_\eps \left\| \hat g_\eps^\pm\right\|_{L^2\left(Mdv\right)}^2.
		\end{equation*}
		Moreover, it is readily seen that
		\begin{equation*}
			\left\|\mathds{1}_{\left\{G_\eps^\pm\geq 2\right\}}\hat g_\eps^\pm\left(1-\chi\left(\frac{|v|^2}{K_\eps}\right)\right)\right\|_{L^2\left(Mdv\right)}
			\leq \frac{C}{K_\eps^\frac 12}
			\left\|\hat g_\eps^\pm\right\|_{L^2\left((1+|v|^2)Mdv\right)}.
		\end{equation*}

		Therefore, it follows that, combining the preceding estimates and considering the uniform bounds from Lemmas \ref{L2-lem} and \ref{v2-int},
		\begin{equation*}
				\left\|\mathds{1}_{\left\{G_\eps^\pm\geq 2\right\}}\hat g_\eps^\pm\right\|_{L^2\left(Mdv\right)}
				= o(1)_{L^1_\mathrm{loc}(dtdx)},
		\end{equation*}
		which concludes the proof of \eqref{bound 7} by interpolation.

		\noindent$\bullet$ To deduce the second bound \eqref{bound 8}, we decompose, writing $g_\eps^\pm=\hat g_\eps^\pm\left(1+\frac \eps 4\hat g_\eps^\pm\right)$ and using that $\left(1+\frac \eps 4\hat g_\eps^\pm\right)\gamma_\eps^\pm$ is uniformly bounded pointwise,
		\begin{equation*}
			\begin{aligned}
				\left|g_\eps^\pm\gamma_\eps^\pm\chi\left(\frac{|v|^2}{K_\eps}\right)-\hat g_\eps^\pm\right|
				& =
				\left|\hat g_\eps^\pm\left(\left(1+\frac \eps 4\hat g_\eps^\pm\right)\gamma_\eps^\pm\chi\left(\frac{|v|^2}{K_\eps}\right)-1\right)\right|
				\\
				& \leq
				C\left|\hat g_\eps^\pm - \Pi \hat g_\eps^\pm\right|
				+ C\left|\Pi \hat g_\eps^\pm\left(\left(1+\frac \eps 4\hat g_\eps^\pm\right)\gamma_\eps^\pm\chi\left(\frac{|v|^2}{K_\eps}\right)-1\right)\right|
				\\
				& \leq
				C\left|\hat g_\eps^\pm - \Pi \hat g_\eps^\pm\right|
				+ C\left|\Pi \hat g_\eps^\pm\left(1-\chi\left(\frac{|v|^2}{K_\eps}\right)\right)\right|
				\\
				& + C\eps\left|\Pi \hat g_\eps^\pm\hat g_\eps^\pm \gamma_\eps^\pm\chi\left(\frac{|v|^2}{K_\eps}\right)\right|
				+ C\left|\Pi \hat g_\eps^\pm\left(\gamma_\eps^\pm-1\right)\chi\left(\frac{|v|^2}{K_\eps}\right)\right|,
			\end{aligned}
		\end{equation*}
		which implies, since $\left\|\frac{\Pi\hat g_\eps}{1+|v|^2}\right\|_{L^\infty(Mdv)}\leq C\left\|\hat g_\eps\right\|_{L^2(Mdv)}$,
		\begin{equation*}
			\begin{aligned}
				& \left\|g_\eps^\pm\gamma_\eps^\pm\chi\left(\frac{|v|^2}{K_\eps}\right)-\hat g_\eps^\pm\right\|_{L^2(Mdv)}
				\\
				& \hspace{10mm} \leq
				C\left\|\hat g_\eps^\pm - \Pi \hat g_\eps^\pm\right\|_{L^2(Mdv)}
				+ C\eps K_\eps\left\|\hat g_\eps^\pm\right\|_{L^2(Mdv)}^2
				\\
				& \hspace{10mm} + C\left\|(1+|v|^2)\left(1-\chi\left(\frac{|v|^2}{K_\eps}\right)\right)\right\|_{L^2(Mdv)}\left\| \hat g_\eps^\pm \right\|_{L^2\left(Mdv\right)}.
			\end{aligned}
		\end{equation*}

		Then, employing the control of Gaussian tails \eqref{gaussian-decay 0} and the uniform bound from Lemma \ref{L2-lem}, we infer that
		\begin{equation*}
			\begin{aligned}
				\left\|g_\eps^\pm\gamma_\eps^\pm\chi\left(\frac{|v|^2}{K_\eps}\right)-\hat g_\eps^\pm\right\|_{L^2(Mdv)}
				& \leq
				C\left\|\hat g_\eps^\pm - \Pi \hat g_\eps^\pm\right\|_{L^2(Mdv)}
				\\
				& + O\left(\eps |\log\eps|\right)_{L^\infty\left(dt;L^1(dx)\right)}
				\\
				& + O\left(\eps^\frac K4|\log\eps|^\frac 54\right)_{L^\infty\left(dt;L^2(dx)\right)},
			\end{aligned}
		\end{equation*}
		which, with an interpolation argument, establishes \eqref{bound 8}.

		\noindent$\bullet$ The third bound \eqref{bound 9} easily follows from the estimate
		\begin{equation*}
			\begin{aligned}
				& \left\|g_\eps^\pm\gamma_\eps^\pm-\hat g_\eps^\pm\right\|_{L^2(Mdv)}
				\\
				& \leq
				\left\|g_\eps^\pm\gamma_\eps^\pm\chi\left(\frac{|v|^2}{K_\eps}\right)-\hat g_\eps^\pm\right\|_{L^2(Mdv)}
				+
				\left\| \hat g_\eps^\pm\left(1-\chi\left(\frac{|v|^2}{K_\eps}\right)\right)\right\|_{L^2(Mdv)}
				\\
				& \leq
				\left\|g_\eps^\pm\gamma_\eps^\pm\chi\left(\frac{|v|^2}{K_\eps}\right)-\hat g_\eps^\pm\right\|_{L^2(Mdv)}
				+\frac{C}{K_\eps^\frac 12}
				\left\| \hat g_\eps^\pm \right\|_{L^2\left(\left(1+|v|^2\right)Mdv\right)},
			\end{aligned}
		\end{equation*}
		which, when combined with the second bound \eqref{bound 8}, concludes its justification.

		\noindent$\bullet$ Next, simply notice that the fourth bound \eqref{bound 11} is a reformulation of Lemma \ref{relaxation-control} with an interpolation argument, which we have incorporated here for mere convenience.

		\noindent$\bullet$ Finally, we easily establish the last bound \eqref{bound 10}. To this end, note first that the case $p=2$ is easily deduced from \eqref{bound 7}, using again that $G_\eps^\pm\geq 2$ implies $\eps\hat g_\eps^\pm\geq 2(\sqrt 2 - 1)$. Furthermore, repeating the estimate leading to the bound \eqref{special bound} yields here that, for every $2\leq p<4$,
		\begin{equation*}
			\frac1{\eps} \mathds{1}_{\left\{G_\eps^\pm\geq 2\right\}}=
			O(1)_{L^2_\mathrm{loc}\left(dtdx ; L^p\left(Mdv\right)\right)}.
		\end{equation*}
		Therefore, the bound \eqref{bound 10} is obtained by interpolation.
		
		The proof of the lemma is thus complete.
	\end{proof}

\begin{proof}[Proof of Lemma \ref{trunc-lem 3}]
	First, we utilize \eqref{v2 trick}, with some fixed $\gamma>4$ therein to be determined later on, to estimate
	\begin{equation*}
		\begin{aligned}
			(1+|v|)^2 \left|\hat g_\eps^\pm  \right|^2
			& \leq
			{2\gamma^\frac 12 \over \eps} \sqrt{h\left(\eps g_\eps^\pm\right)}(1+|v|) \left|\hat g_\eps^\pm \right|
			+
			\frac{2\gamma^\frac 12}{\eps}e^{\frac{(1+|v|)^2}{2\gamma}}(1+|v|) \left|\hat g_\eps^\pm \right|
			\\
			& \leq
			{2\gamma \over \eps^2} {h\left(\eps g_\eps^\pm\right)}+\frac 12 (1+|v|)^2 \left|\hat g_\eps^\pm  \right|^2
			+
			\frac{2\gamma^\frac 12}{\eps}e^{\frac{(1+|v|)^2}{2\gamma}}(1+|v|) \left|\hat g_\eps^\pm  \right|,
		\end{aligned}
	\end{equation*}
	whence
	\begin{equation*}
		(1+|v|)^2 \left|\hat g_\eps^\pm  \right|^2
		\leq
		{4\gamma \over \eps^2} {h\left(\eps g_\eps^\pm\right)}
		+
		\frac{4\gamma^\frac 12}{\eps}e^{\frac{(1+|v|)^2}{2\gamma}}(1+|v|) \left|\hat g_\eps^\pm  \right|.
	\end{equation*}
	It follows that, for any $2<p<4$,
	\begin{equation}\label{previous bound}
		\begin{aligned} 
			& \left\|\mathds{1}_{\left\{G_\eps^\pm\geq 2\right\}}\hat g_\eps^\pm\right\|_{L^2\left((1+|v|)^2 Mdv\right)}^2
			\\
			& \leq
			{4\gamma \over \eps^2} \int_{\mathbb{R}^3} \mathds{1}_{\left\{G_\eps^\pm\geq 2\right\}} {h\left(\eps g_\eps^\pm\right)}Mdv
			\\
			& + 2 \left\|\mathds{1}_{\left\{G_\eps^\pm\geq 2\right\}}\hat g_\eps^\pm\right\|_{L^2\left(Mdv\right)}^2
			+2 \left\|\frac {\gamma^\frac 12}\eps\mathds{1}_{\left\{G_\eps^\pm\geq 2\right\}}
			e^{\frac{(1+|v|)^2}{2\gamma}}(1+|v|) \right\|_{L^2\left(Mdv\right)}^2
			\\
			& \leq
			{4\gamma \over \eps^2} \int_{\mathbb{R}^3} \mathds{1}_{\left\{G_\eps^\pm\geq 2\right\}} {h\left(\eps g_\eps^\pm\right)}Mdv
			+ 2 \left\|\mathds{1}_{\left\{G_\eps^\pm\geq 2\right\}}\hat g_\eps^\pm\right\|_{L^2\left(Mdv\right)}^2
			\\
			& +2\gamma \left\|\frac {1}\eps\mathds{1}_{\left\{G_\eps^\pm\geq 2\right\}}\right\|_{L^p\left(Mdv\right)}^2
			\left\|
			e^{\frac{(1+|v|)^2}{2\gamma}}(1+|v|) \right\|_{L^{\frac{2p}{p-2}}\left(Mdv\right)}^2.
		\end{aligned}
	\end{equation}
	Here, we need to set the parameter $\gamma$ so large that $\frac{p}{p-2}<\frac\gamma 2$ in order to yield a finite constant in the last term above.
	
	Then, further using that $h(z)=\frac 12 \left(z\gamma(1+z)\right)^2+O(z^3)$, we deduce
	\begin{equation*}
		\begin{aligned}
			& \left\|\mathds{1}_{\left\{G_\eps^\pm\geq 2\right\}}\hat g_\eps^\pm\right\|_{L^2\left((1+|v|)^2 Mdv\right)}^2
			\\
			& \leq
			{C_1 \over \eps^2} \int_{\mathbb{R}^3} \mathds{1}_{\left\{G_\eps^\pm\geq 2\right\}} \left({h\left(\eps g_\eps^\pm\right)}
			- \frac{\eps^2}{2} \left(g_\eps^\pm \gamma_\eps^\pm\right)^2 \right)Mdv
			\\
			& + C_2 \left\|\mathds{1}_{\left\{G_\eps^\pm\geq 2\right\}}\hat g_\eps^\pm\right\|_{L^2\left(Mdv\right)}^2
			+C_2\left\|\frac {1}\eps\mathds{1}_{\left\{G_\eps^\pm\geq 2\right\}}\right\|_{L^p\left(Mdv\right)}^2
			\\
			& \leq
			C_1\int_{\mathbb{R}^3} \left(\frac 1{\eps^2}{h\left(\eps g_\eps^\pm\right)}
			- \frac{1}{2} \left(g_\eps^\pm\gamma_\eps^\pm\right)^2 \right)Mdv
			\\
			& + C_2 \left\|\mathds{1}_{\left\{G_\eps^\pm\geq 2\right\}}\hat g_\eps^\pm\right\|_{L^2\left(Mdv\right)}^2
			+C_2 \left\|\frac {1}\eps\mathds{1}_{\left\{G_\eps^\pm\geq 2\right\}}\right\|_{L^p\left(Mdv\right)}^2
			\\
			& + \frac {C_2}{\eps^2}\int_{\mathbb{R}^3} \mathds{1}_{\left\{G_\eps^\pm< 2\right\}} \left|{h\left(\eps g_\eps^\pm\right)}
			- \frac{\eps^2}{2}\left( g_\eps^\pm\gamma_\eps^\pm\right)^2 \right|Mdv
			\\
			& \leq
			C_1\int_{\mathbb{R}^3} \left(\frac 1{\eps^2}{h\left(\eps g_\eps^\pm\right)}
			- \frac{1}{2} \left( g_\eps^\pm \gamma_\eps^\pm\right)^2 \right)Mdv
			\\
			& + C_2 \left\|\mathds{1}_{\left\{G_\eps^\pm\geq 2\right\}}\hat g_\eps^\pm\right\|_{L^2\left(Mdv\right)}^2
			+C_2\left\|\frac {1}\eps\mathds{1}_{\left\{G_\eps^\pm\geq 2\right\}}\right\|_{L^p\left(Mdv\right)}^2
			\\
			& + C_2 \int_{\mathbb{R}^3} \mathds{1}_{\left\{\left|\eps g_\eps^\pm\right|< 1\right\}}
			\left|\eps g_\eps^\pm\right|\left|g_\eps^\pm\gamma_\eps^\pm\right|^2Mdv,
		\end{aligned}
	\end{equation*}
	where, $C_1>0$ and $C_2>0$ denote diverse constants which only depend on fixed parameters and which we do not distiguinsh for simplicity.
	
	Finally, combining the preceding estimate with the bounds \eqref{bound 1} and \eqref{bound 4}, modulating the last term $\left|g_\eps^\pm\gamma_\eps^\pm\right|^2 = \left|g_\eps^\pm\gamma_\eps^\pm-\bar g^\pm\right|^2+2g_\eps^\pm\gamma_\eps^\pm\bar g^\pm - \left|\bar g^\pm\right|^2$, and using that $\mathds{1}_{\left\{\left|\eps g_\eps^\pm\right|< 1\right\}}\left|\eps g_\eps^\pm\right|$ is bounded pointwise and converges almost everywhere to zero (possibly up to extraction of subsequences) with the Product Limit Theorem, we deduce the first estimate of the lemma.

	The remaining estimate requires some care, for the function $h(z)-\frac 12\left(z\gamma(1+z)\right)^2$ can take negative values. Integrating locally in $x$ the previous bound \eqref{previous bound}, we first observe that
	\begin{equation*}
		\begin{aligned} 
			& \left\|\mathds{1}_{\left\{G_\eps^\pm\geq 2\right\}}\hat g_\eps^\pm\right\|_{L^2_\mathrm{loc}\left(dx;L^2\left((1+|v|)^2 Mdv\right)\right)}^2
			\\
			& \leq
			{C_1 \over \eps^2} \int_{\mathbb{R}^3\times\mathbb{R}^3} \mathds{1}_{\left\{
			{h\left(\eps g_\eps^\pm\right)}
			\geq \frac{\eps^2}{2} \left(g_\eps^\pm \gamma_\eps^\pm\chi \left({|v|^2\over K_\eps}\right)\right)^2
			\right\}}
			\\
			& \times \left({h\left(\eps g_\eps^\pm\right)}
			- \frac{\eps^2}{2} \left(g_\eps^\pm \gamma_\eps^\pm\chi \left({|v|^2\over K_\eps}\right)\right)^2 \right)Mdxdv
			\\
			& + C_2 \left\|\mathds{1}_{\left\{G_\eps^\pm\geq 2\right\}}\hat g_\eps^\pm\right\|_{L^2_\mathrm{loc}\left(dx;L^2\left(Mdv\right)\right)}^2
			+C_2 \left\|\frac {1}\eps\mathds{1}_{\left\{G_\eps^\pm\geq 2\right\}}\right\|_{L^2_\mathrm{loc}\left(dx;L^p\left(Mdv\right)\right)}^2
			\\
			& \leq
			C_1\int_{\mathbb{R}^3\times\mathbb{R}^3} \left(\frac 1{\eps^2}{h\left(\eps g_\eps^\pm\right)}
			- \frac{1}{2} \left(g_\eps^\pm\gamma_\eps^\pm\chi \left({|v|^2\over K_\eps}\right)\right)^2 \right)Mdxdv
			\\
			& + C_2 \left\|\mathds{1}_{\left\{G_\eps^\pm\geq 2\right\}}\hat g_\eps^\pm\right\|_{L^2_\mathrm{loc}\left(dx;L^2\left(Mdv\right)\right)}^2
			+C_2\left\|\frac {1}\eps\mathds{1}_{\left\{G_\eps^\pm\geq 2\right\}}\right\|_{L^2_\mathrm{loc}\left(dx;L^p\left(Mdv\right)\right)}^2
			\\
			& + \frac {C_2}{\eps^2}\int_{\mathbb{R}^3\times\mathbb{R}^3} \mathds{1}_{\left\{
			{h\left(\eps g_\eps^\pm\right)}
			< \frac{\eps^2}{2} \left(g_\eps^\pm \gamma_\eps^\pm\chi \left({|v|^2\over K_\eps}\right)\right)^2
			\right\}}
			\\
			&\times \left(
			\frac{\eps^2}{2}\left( g_\eps^\pm\gamma_\eps^\pm\chi \left({|v|^2\over K_\eps}\right)\right)^2
			-{h\left(\eps g_\eps^\pm\right)}
			\right)Mdxdv,
		\end{aligned}
	\end{equation*}
	whence
	\begin{equation*}
		\begin{aligned} 
			& \left\|\mathds{1}_{\left\{G_\eps^\pm\geq 2\right\}}\hat g_\eps^\pm\right\|_{L^2_\mathrm{loc}\left(dx;L^2\left((1+|v|)^2 Mdv\right)\right)}^2
			\\
			& \leq
			C_1\int_{\mathbb{R}^3\times\mathbb{R}^3} \left(\frac 1{\eps^2}{h\left(\eps g_\eps^\pm\right)}
			- \frac{1}{2} \left(g_\eps^\pm\gamma_\eps^\pm\chi \left({|v|^2\over K_\eps}\right)\right)^2 \right)Mdxdv
			\\
			& + C_2 \left\|\mathds{1}_{\left\{G_\eps^\pm\geq 2\right\}}\hat g_\eps^\pm\right\|_{L^2_\mathrm{loc}\left(dx;L^2\left(Mdv\right)\right)}^2
			+C_2\left\|\frac {1}\eps\mathds{1}_{\left\{G_\eps^\pm\geq 2\right\}}\right\|_{L^2_\mathrm{loc}\left(dx;L^p\left(Mdv\right)\right)}^2
			\\
			& + \frac {C_2}{\eps^2}\int_{\mathbb{R}^3\times\mathbb{R}^3} \mathds{1}_{\left\{
			{h\left(\eps g_\eps^\pm\right)}
			< \frac{\eps^2}{2} \left(g_\eps^\pm \gamma_\eps^\pm\right)^2
			\right\}}
			\left(
			\frac{\eps^2}{2}\left( g_\eps^\pm\gamma_\eps^\pm\right)^2
			-{h\left(\eps g_\eps^\pm\right)}
			\right)\chi \left({|v|^2\over K_\eps}\right)^2Mdxdv.
		\end{aligned}
	\end{equation*}

	Next, considering $N>0$ so large that $h(z) < \frac 12\left(z\gamma(1+z)\right)^2$ implies $|z|\leq N$, for any $z\in [-1,\infty)$, which is always possible in view of the assumptions \eqref{gamma-assumption} on $\gamma(z)$, and using that $h(z)=\frac 12 \left(z\gamma(1+z)\right)^2+O(z^3)$ again, we infer that
	\begin{equation*}
		\begin{aligned}
			& \left\|\mathds{1}_{\left\{G_\eps^\pm\geq 2\right\}}\hat g_\eps^\pm\right\|_{L^2_\mathrm{loc}\left(dx;L^2\left((1+|v|)^2 Mdv\right)\right)}^2
			\\
			& \leq
			C_1\int_{\mathbb{R}^3\times\mathbb{R}^3} \left(\frac 1{\eps^2}{h\left(\eps g_\eps^\pm\right)}
			- \frac{1}{2} \left( g_\eps^\pm \gamma_\eps^\pm\chi \left({|v|^2\over K_\eps}\right)\right)^2 \right)Mdxdv
			\\
			& + C_2 \left\|\mathds{1}_{\left\{G_\eps^\pm\geq 2\right\}}\hat g_\eps^\pm\right\|_{L^2_\mathrm{loc}\left(dx;L^2\left(Mdv\right)\right)}^2
			+C_2\left\|\frac {1}\eps\mathds{1}_{\left\{G_\eps^\pm\geq 2\right\}}\right\|_{L^2_\mathrm{loc}\left(dx;L^p\left(Mdv\right)\right)}^2
			\\
			& + C_2 \int_{\mathbb{R}^3\times\mathbb{R}^3} \mathds{1}_{\left\{\left|\eps g_\eps^\pm\right|\leq N\right\}}
			\left|\eps g_\eps^\pm\right|\left|g_\eps^\pm\gamma_\eps^\pm\chi \left({|v|^2\over K_\eps}\right)\right|^2Mdxdv.
		\end{aligned}
	\end{equation*}
	
	Then, as before, combining the preceding estimate with the bounds \eqref{bound 1} and \eqref{bound 4}, and modulating the last term $\frac 12\left|g_\eps^\pm\gamma_\eps^\pm\chi \left({|v|^2\over K_\eps}\right)\right|^2 \leq \left|g_\eps^\pm\gamma_\eps^\pm\chi \left({|v|^2\over K_\eps}\right)-\bar g^\pm\right|^2+ \left|\bar g^\pm\right|^2$, we arrive at
	\begin{equation*}
		\begin{aligned} 
			& \left\|\mathds{1}_{\left\{G_\eps^\pm\geq 2\right\}}\hat g_\eps^\pm\right\|_{L^2_\mathrm{loc}\left(dx;L^2\left((1+|v|)^2 Mdv\right)\right)}^2
			\\
			& \leq
			C_1\int_{\mathbb{R}^3\times\mathbb{R}^3} \left(\frac 1{\eps^2}{h\left(\eps g_\eps^\pm\right)}
			- \frac{1}{2} \left( g_\eps^\pm \gamma_\eps^\pm\chi \left({|v|^2\over K_\eps}\right)\right)^2 \right)Mdxdv
			\\
			& + C_2 \left\|g_\eps^\pm\gamma_\eps^\pm\chi \left({|v|^2\over K_\eps}\right)-\bar g^\pm\right\|_{L^2\left(Mdxdv\right)}^2 + o(1)_{L^1_\mathrm{loc}(dt)}
			\\
			& + C_2 \int_{\mathbb{R}^3\times\mathbb{R}^3} \mathds{1}_{\left\{\left|\eps g_\eps^\pm\right|\leq N\right\}}
			\left|\eps g_\eps^\pm\right|\left|\bar g^\pm\right|^2Mdxdv.
		\end{aligned}
	\end{equation*}
	Finally, since $\bar g^\pm$ belongs to $L^\infty\left(dt;L^2\left(Mdxdv\right)\right)$ and $\mathds{1}_{\left\{\left|\eps g_\eps^\pm\right|< 1\right\}}\left|\eps g_\eps^\pm\right|$ is bounded pointwise and converges almost everywhere to zero, we deduce, through a straightforward application of Egorov's theorem, that the last term above vanishes locally in $L^1(dt)$, which concludes the proof of the lemma.
\end{proof}

\begin{proof}[Proof of Lemma \ref{trunc-lem 4}]
	We begin by estimating, using the relaxation estimate \eqref{relaxation estimate} from Lemma \ref{relaxation2-control},
	\begin{equation*}
		\begin{aligned}
			& \left\|\left(\hat g_\eps^+-\hat g_\eps^--\hat n_\eps\right)\hat g_\eps^\pm\mathds{1}_{\left\{G_\eps^\pm\geq 2\right\}}\right\|_{L^1\left((1+|v|)^2 Mdv\right)}
			\\
			& \leq
			\left\|\hat g_\eps^\pm\mathds{1}_{\left\{G_\eps^\pm\geq 2\right\}}\right\|_{L^2\left((1+|v|)^2 Mdv\right)}^2
			+\frac 14
			\left\|\left(\hat g_\eps^+-\hat g_\eps^--\hat n_\eps\right)\right\|_{L^2\left((1+|v|)^2 Mdv\right)}^2
			\\
			& \leq
			\left\|\hat g_\eps^\pm\mathds{1}_{\left\{G_\eps^\pm\geq 2\right\}}\right\|_{L^2\left((1+|v|)^2 Mdv\right)}^2
			+\frac 14
			\left\|\left(\hat g_\eps^+-\hat g_\eps^--\hat n_\eps\right)\left(\left|\hat g_\eps^+\right|+\left|\hat g_\eps^-\right|\right)\right\|_{L^1\left((1+|v|)^2 Mdv\right)}
			\\
			& +\frac 14
			\left\|\hat g_\eps^+-\hat g_\eps^--\hat n_\eps\right\|_{L^2\left( Mdv\right)}
			\left\|\hat n_\eps\right\|_{L^2\left((1+|v|)^4 Mdv\right)}
			\\
			& \leq
			\left\|\hat g_\eps^\pm\mathds{1}_{\left\{G_\eps^\pm\geq 2\right\}}\right\|_{L^2\left((1+|v|)^2 Mdv\right)}^2
			+\frac 14
			\left\|\left(\hat g_\eps^+-\hat g_\eps^--\hat n_\eps\right)\left(\left|\hat g_\eps^+\right|+\left|\hat g_\eps^-\right|\right)\right\|_{L^1\left((1+|v|)^2 Mdv\right)}
			\\
			& +C
			\left\|\hat g_\eps^+-\hat g_\eps^--\hat n_\eps\right\|_{L^2\left( Mdv\right)}
			\left\|\hat g_\eps^+-\hat g_\eps^-\right\|_{L^2\left(Mdv\right)}
			\\
			& \leq
			\left\|\hat g_\eps^\pm\mathds{1}_{\left\{G_\eps^\pm\geq 2\right\}}\right\|_{L^2\left((1+|v|)^2 Mdv\right)}^2
			+\frac 14\sum_{\pm}
			\left\|\left(\hat g_\eps^+-\hat g_\eps^--\hat n_\eps\right)\hat g_\eps^\pm \right\|_{L^1\left((1+|v|)^2 Mdv\right)}
			\\
			& +C\sum_{\pm}
			\left\|\hat g_\eps^\pm-\bar g^\pm\right\|_{L^2\left( Mdv\right)}^2 + o(1)_{L^1_\mathrm{loc}(dtdx)}.
		\end{aligned}
	\end{equation*}
	It follows that
	\begin{equation*}
		\begin{aligned}
			\sum_{\pm} & \left\|\left(\hat g_\eps^+-\hat g_\eps^--\hat n_\eps\right)\hat g_\eps^\pm\right\|_{L^1\left((1+|v|)^2 Mdv\right)}
			\\
			& \leq
			\sum_{\pm}\left\|\left(\hat g_\eps^+-\hat g_\eps^--\hat n_\eps\right)\hat g_\eps^\pm\mathds{1}_{\left\{G_\eps^\pm< 2\right\}}\right\|_{L^1\left((1+|v|)^2 Mdv\right)}
			\\
			& +
			\sum_{\pm}\left\|\left(\hat g_\eps^+-\hat g_\eps^--\hat n_\eps\right)\hat g_\eps^\pm\mathds{1}_{\left\{G_\eps^\pm\geq 2\right\}}\right\|_{L^1\left((1+|v|)^2 Mdv\right)}
			\\
			& \leq
			\sum_{\pm}\left\|\left(\hat g_\eps^+-\hat g_\eps^--\hat n_\eps\right)\hat g_\eps^\pm\mathds{1}_{\left\{G_\eps^\pm< 2\right\}}\right\|_{L^1\left((1+|v|)^2 Mdv\right)}
			\\
			& +
			\sum_{\pm}\left\|\hat g_\eps^\pm\mathds{1}_{\left\{G_\eps^\pm\geq 2\right\}}\right\|_{L^2\left((1+|v|)^2 Mdv\right)}^2
			+\frac 12\sum_{\pm}
			\left\|\left(\hat g_\eps^+-\hat g_\eps^--\hat n_\eps\right)\hat g_\eps^\pm \right\|_{L^1\left((1+|v|)^2 Mdv\right)}
			\\
			& +C\sum_{\pm}
			\left\|\hat g_\eps^\pm-\bar g^\pm\right\|_{L^2\left( Mdv\right)}^2 + o(1)_{L^1_\mathrm{loc}(dtdx)},
		\end{aligned}
	\end{equation*}
	whence
	\begin{equation*}
		\begin{aligned}
			& \sum_{\pm} \left\|\left(\hat g_\eps^+-\hat g_\eps^--\hat n_\eps\right)\hat g_\eps^\pm\right\|_{L^1\left((1+|v|)^2 Mdv\right)}
			\\
			& \leq
			2 \sum_{\pm}\left\|\left(\hat g_\eps^+-\hat g_\eps^--\hat n_\eps\right)\hat g_\eps^\pm\mathds{1}_{\left\{G_\eps^\pm< 2\right\}}\right\|_{L^1\left((1+|v|)^2 Mdv\right)}
			\\
			& +
			2 \sum_{\pm}\left\|\hat g_\eps^\pm\mathds{1}_{\left\{G_\eps^\pm\geq 2\right\}}\right\|_{L^2\left((1+|v|)^2 Mdv\right)}^2
			+C\sum_{\pm}
			\left\|\hat g_\eps^\pm-\bar g^\pm\right\|_{L^2\left( Mdv\right)}^2 + o(1)_{L^1_\mathrm{loc}(dtdx)}.
		\end{aligned}
	\end{equation*}
	
	Thus, in view of the estimates \eqref{bound 2} and \eqref{bound 3} from Lemma \ref{trunc-lem} and utilizing Lemma \ref{trunc-lem 3}, we deduce that, in order to conclude the proof of the lemma, it is sufficient to establish that
	\begin{equation}\label{trunc-lem 4 claim}
		\begin{aligned}
			\left\|\left(\hat g_\eps^+-\hat g_\eps^--\hat n_\eps\right)\hat g_\eps^\pm\mathds{1}_{\left\{G_\eps^\pm< 2\right\}}\right\|_{L^1\left((1+|v|)^2 Mdv\right)}
			& \leq C\sum_{\pm}
			\left\|\hat g_\eps^\pm-\bar g^\pm\right\|_{L^2\left( Mdv\right)}^2
			\\
			& + o(1)_{L^1_\mathrm{loc}(dtdx)}.
		\end{aligned}
	\end{equation}
	To this end, employing the estimate \eqref{relaxation2-control 2} from Lemma \ref{relaxation2-control}, we first obtain that
	\begin{equation*}
		\begin{aligned}
			& \left\|\left(\hat g_\eps^+-\hat g_\eps^--\hat n_\eps\right)\hat g_\eps^\pm\mathds{1}_{\left\{G_\eps^\pm< 2\right\}}\right\|_{L^1\left((1+|v|)^2 Mdv\right)}
			\\
			& \leq
			\left\|\left(\hat g_\eps^+-\hat g_\eps^--\hat n_\eps
			-\frac\eps 2\hat n_\eps\left(\hat g_\eps^\pm-\hat \rho_\eps^\pm\right)
			\right)\hat g_\eps^\pm\mathds{1}_{\left\{G_\eps^\pm< 2\right\}}\right\|_{L^1\left((1+|v|)^2 Mdv\right)}
			\\
			& +\frac 12 \left|\hat n_\eps\right|
			\left\|\left(\hat g_\eps^\pm-\hat \rho_\eps^\pm\right)\eps\hat g_\eps^\pm\mathds{1}_{\left\{G_\eps^\pm< 2\right\}}\right\|_{L^1\left((1+|v|)^2 Mdv\right)}
			\\
			& \leq
			\left\|\hat g_\eps^+-\hat g_\eps^--\hat n_\eps
			-\frac\eps 2\hat n_\eps\left(\hat g_\eps^\pm-\hat \rho_\eps^\pm\right)
			\right\|_{L^2\left(Mdv\right)}
			\left\|\hat g_\eps^\pm\mathds{1}_{\left\{G_\eps^\pm< 2\right\}}\right\|_{L^2\left((1+|v|)^4 Mdv\right)}
			\\
			& +C\left(\sum_\pm \left\|\hat g_\eps^\pm\right\|_{L^2\left(Mdv\right)}^2\right)
			\left\|\eps\hat g_\eps^\pm\mathds{1}_{\left\{G_\eps^\pm< 2\right\}}\right\|_{L^2\left((1+|v|)^4 Mdv\right)}
			\\
			& \leq
			C\left(\sum_\pm \left\|\hat g_\eps^\pm\right\|_{L^2\left(Mdv\right)}^2\right)
			\left\|\eps\hat g_\eps^\pm\mathds{1}_{\left\{G_\eps^\pm< 2\right\}}\right\|_{L^2\left((1+|v|)^4 Mdv\right)}
			+o(1)_{L^1_\mathrm{loc}(dtdx)}.
		\end{aligned}
	\end{equation*}
	Then, noticing that, in view of Lemma \ref{v2-int},
	\begin{equation*}
		\begin{aligned}
			\left\|\eps\hat g_\eps^\pm\mathds{1}_{\left\{G_\eps^\pm< 2\right\}}\right\|_{L^2\left((1+|v|)^4 Mdv\right)}
			& =O(1)_{L^\infty(dtdx)},
			\\
			\left\|\eps\hat g_\eps^\pm\mathds{1}_{\left\{G_\eps^\pm< 2\right\}}\right\|_{L^2\left((1+|v|)^4 Mdv\right)}
			& =O(\eps)_{L^2_\mathrm{loc}(dtdx)},
		\end{aligned}
	\end{equation*}
	we conclude that \eqref{trunc-lem 4 claim} holds, which completes the proof of the lemma.
\end{proof}

%% file: t-oscillations0.tex
% \chapter{Time oscillations and compensated compactness}
\chapter{Acoustic and electromagnetic waves}
\label{oscillations}

In Chapter \ref{weak bounds}, we conducted a rather extensive study of the scaled relative entropy and entropy dissipation bounds. These yielded controls on the fluctuations in all variables $t$, $x$ and $v$ in appropriate function spaces and, thus, allowed us to establish essential weak compactness estimates on the fluctuations. Moreover, relaxation estimates were also obtained therein, showing that fluctuations remain close to their hydrodynamic projection, which implied improved controls in the $v$ variable.

Then, in Chapter \ref{hypoellipticity}, we showed that the control of the behavior of fluctuations in $v$ could be improved to strong compactness estimates in $v$, which could then be transfered --~exploiting the hypoelliptic phenomenon in kinetic transport equations~-- to the $x$ variable to deduce strong compactness estimates in both $x$ and $v$.

Thus, we know so far that there are no oscillations in $x$ and $v$ in the fluctuations as the Knudsen number tends to zero. Note, however, that nothing is claimed about the control of oscillations in the $t$ variable in the fluctuations and the control of oscillations in $t$ and $x$ in the electromagnetic fields.

In fact, because of the scaling of the transport operator $\eps \d_t + v\cdot \nabla_x$, we do not expect to obtain additional regularity or compactness with respect to time on the fluctuations~: the natural variable is indeed the fast time $\frac t\eps$ (see discussion in Section \ref{space compactness}). We are however interested in the slow macroscopic dynamics. Since there is nothing to prevent an oscillatory behavior in $t$, we need to further describe the dependence of fluctuations with respect to time and filter the fast oscillations.

There may also be persistence of fast oscillations in both $t$ and $x$ in the electromagnetic fields (and electrodynamic macroscopic variables, such as the electric current), which we do not expect to control due to the hyperbolic nature of Maxwell's equations.

It turns out that oscillations in fluctuations and electromagnetic fields are sometimes coupled. We will therefore need to treat and filter them simultaneously.

In the context of the viscous incompressible hydrodynamic limit of the Boltzmann equation, the filtering of acoustic waves was first understood by Lions and Masmoudi in \cite{lions4}.

\bigskip

In the present chapter, we are going to focus exclusively on the one species setting treated in Theorem \ref{NS-WEAKCV}, i.e.\ on the regime leading to the incompressible quasi-static Navier-Stokes-Fourier-Maxwell-Poisson system. The proof of this result is based on weak compactness methods which require the handling of possible time and electromagnetic oscillations. In this case, the available strong compactness with respect to spatial variables is good enough and we are actually able to get here a rough description of oscillations, which will be sufficient to derive the weak stability and convergence of the Vlasov-Maxwell-Boltzmann system \eqref{VMB1} as $\eps \to 0$.

As for Theorems \ref{CV-OMHD} and \ref{CV-OMHDSTRONG} concerning the two species setting, i.e.\ in the regime leading to the two-fluid incompressible Navier-Stokes-Fourier-Maxwell system with (solenoidal) Ohm's law, the previous filtering method cannot be applied, and --~as already mentioned~-- there is no asymptotic weak stability of the Vlasov-Maxwell-Boltzmann system \eqref{VMB2} (nor existence of weak solutions to the corresponding limiting model). In order to bypass this difficulty, the idea in this setting is then to compare the actual solutions of the scaled Vlasov-Maxwell-Boltzmann system to some approximate solutions (known a priori to be regular in $t$ and $x$) capturing the fast oscillations. This method of proof, detailed in Chapter \ref{entropy method} later on, is the so-called renormalized modulated entropy method, which is only performed in this work in the case of well-prepared initial data, for the sake of simplicity. The oscillations are therefore automatically filtered out by the method and we do not need to further describe the time dependence of fluctuations. Of course, the case of ill-prepared initial data for two species is interesting and should be addressed. Nevertheless, this issue only seems to present difficulties somewhat similar to those encountered in the handling of initial data in the asymptotic problems considered in \cite{PSR} and \cite[Chapter 5]{SR}, for instance.

\section{Formal filtering of oscillations}

Now, as mentionned above, let us focus exclusively, for the remainder of the present chapter, on the regime of Theorem \ref{NS-WEAKCV} (with one species of particles only) leading to the incompressible quasi-static Navier-Stokes-Fourier-Maxwell-Poisson system \eqref{NSFMP 2}.

On the one hand, going back to the corresponding formal analysis from Chapter \ref{formal-chap}, we expect that the fast time oscillations are governed by the following singular linear system given by \eqref{wave} and \eqref{wave operator 3}~:
\begin{equation}\label{linear structure}
\d_t
\begin{pmatrix}
	\rho_\eps \\ u_\eps\\ \sqrt{3\over 2}\theta_\eps \\ E_\eps \\ B_\eps
\end{pmatrix}
+\frac1\eps
W
% \begin{pmatrix}
% 	0 &\DIV&0&0&0
% 	\\
% 	\nabla_x &0&\sqrt{2\over 3} \nabla_x & - \operatorname{Id} & 0
% 	\\
% 	0 & \sqrt{2\over 3}\DIV &0&0&0
% 	\\
% 	0& \operatorname{Id} &0&0&- \ROT
% 	\\
% 	0&0&0& \ROT &0
% \end{pmatrix}
\begin{pmatrix}
	\rho_\eps\\ u_\eps\\ \sqrt{3\over 2}\theta_\eps \\ E_\eps \\ B_\eps
\end{pmatrix}
=O\left(1\right),
\end{equation}
where the antisymmetric differential operator $W:L^2(dx)\to H^{-1}(dx)$ --~the wave operator~-- is defined by
\begin{equation}\label{W-def}
	W=
	\begin{pmatrix}
		0 &\DIV&0&0&0
		\\
		\nabla_x &0&\sqrt{2\over 3} \nabla_x & - \operatorname{Id} & 0
		\\
		0 & \sqrt{2\over 3}\DIV &0&0&0
		\\
		0& \operatorname{Id} &0&0&- \ROT
		\\
		0&0&0& \ROT &0
	\end{pmatrix}.
\end{equation}

On the other hand, looking back at the formal macroscopic nonlinear system \eqref{moment-eps3}, we see that, in order to derive the limiting system \eqref{NSFMP 2}, we will eventually need to pass to the limit in the nonlinear terms
\begin{equation}\label{nonlinear structure}
		P\left( \nabla_x\cdot \left(u_\eps\otimes u_\eps\right) - \rho_\eps E_\eps - u_\eps\wedge B_\eps \right)
		\quad\text{and}\quad
		\frac 52 \nabla_x\cdot\left(u_\eps\theta_\eps\right) - u_\eps \cdot E_\eps,
\end{equation}
where $P:L^2(dx)\to L^2(dx)$ denotes the Leray projector onto solenoidal vector fields, and establish their weak stability. Since there are oscillations, this will only be possible if one can show that the linear structure \eqref{linear structure} is somehow ``compatible'' with the quadratic forms defined by \eqref{nonlinear structure}.

We explain now why such a ``compatibility'' between the structures of \eqref{linear structure} and \eqref{nonlinear structure} is to be expected, at least formally. First, we decompose the nonlinear terms \eqref{nonlinear structure} into
\begin{equation}\label{compensated 1}
	\begin{aligned}
		P & \left( \nabla_x\cdot \left(u_\eps\otimes u_\eps\right) - \rho_\eps E_\eps - u_\eps\wedge B_\eps \right)
		\\
		& =
		P\left( u_\eps \Div u_\eps + \frac 12 \nabla_x \left|u_\eps\right|^2 - u_\eps\wedge\rot u_\eps - \rho_\eps E_\eps - u_\eps\wedge B_\eps \right)
		\\
		& =
		P\left( u_\eps \left(\eps\partial_t\rho_\eps + \Div u_\eps\right)
		+ \rho_\eps\left(\eps\partial_t u_\eps + \nabla_x\left(\rho_\eps+\theta_\eps\right) - E_\eps \right) \right)
		\\
		& - P\left( \rho_\eps\nabla_x\left(\rho_\eps+\theta_\eps\right)
		+ u_\eps\wedge\left(\rot u_\eps + B_\eps\right) + \eps\partial_t\left(\rho_\eps u_\eps\right) \right)
		\\
		& =
		P\left( u_\eps \left(\eps\partial_t\rho_\eps + \Div u_\eps\right)
		+ \rho_\eps\left(\eps\partial_t u_\eps + \nabla_x\left(\rho_\eps+\theta_\eps\right) - E_\eps \right) \right)
		\\
		& - P\left( \frac{2\rho_\eps-3\theta_\eps}{5}\nabla_x\left(\rho_\eps+\theta_\eps\right)
		+ \frac 3{10}\nabla_x\left(\rho_\eps+\theta_\eps\right)^2\right)
		\\
		& - P\left( u_\eps\wedge\left(\rot u_\eps + B_\eps\right) + \eps\partial_t\left(\rho_\eps u_\eps\right) \right)
		\\
		& =
		P\left( u_\eps \left(\eps\partial_t\rho_\eps + \Div u_\eps\right)
		+ \rho_\eps\left(\eps\partial_t u_\eps + \nabla_x\left(\rho_\eps+\theta_\eps\right) - E_\eps \right) \right)
		\\
		& - P\left( \frac{2\rho_\eps-3\theta_\eps}{5}\nabla_x\left(\rho_\eps+\theta_\eps\right)
		+ u_\eps\wedge\left(\rot u_\eps + B_\eps\right) + \eps\partial_t\left(\rho_\eps u_\eps\right) \right),
	\end{aligned}
\end{equation}
where we used that $P\frac 12\nabla_x \left|u_\eps\right|^2=P\frac 3{10}\nabla_x\left(\rho_\eps+\theta_\eps\right)^2=0$, and
\begin{equation}\label{compensated 2}
	\begin{aligned}
		\frac 52 & \nabla_x\cdot\left(u_\eps\theta_\eps\right) - u_\eps \cdot E_\eps
		\\
		& = \frac 52 \theta_\eps \Div u_\eps + \frac 52 u_\eps\cdot\nabla_x\theta_\eps - u_\eps \cdot E_\eps
		\\
		& = \frac 52 \theta_\eps \left(\frac 32\eps\partial_t\theta_\eps + \Div u_\eps\right) + u_\eps\cdot\left(\eps\partial_t u_\eps + \nabla_x\left(\rho_\eps+\theta_\eps\right) - E_\eps\right)
		\\
		& + u_\eps\cdot\nabla_x\left(\frac 32\theta_\eps-\rho_\eps\right) - \frac{15}{8}\eps\partial_t\theta_\eps^2 - \frac 12 \eps\partial_t u_\eps^2,
	\end{aligned}
\end{equation}
which formally implies, using the first three equations from \eqref{linear structure}, that
\begin{equation*}
	\begin{aligned}
		P ( \nabla_x\cdot \left(u_\eps\otimes u_\eps\right) - \rho_\eps & E_\eps - u_\eps\wedge B_\eps )
		\\
		& = - P\left( \frac{2\rho_\eps-3\theta_\eps}{5}\nabla_x\left(\rho_\eps+\theta_\eps\right)
		+ u_\eps\wedge\left(\rot u_\eps + B_\eps\right) \right) + O(\eps),
		\\
		\frac 52 \nabla_x\cdot\left(u_\eps\theta_\eps\right) - u_\eps \cdot E_\eps
		& = u_\eps\cdot\nabla_x\left(\frac 32\theta_\eps-\rho_\eps\right) + O(\eps).
	\end{aligned}
\end{equation*}

Thus, this decomposition is sufficient to deduce the weak stability of the nonlinear terms \eqref{nonlinear structure} provided the oscillating part of $\left(\rho_\eps, u_\eps, \sqrt{3\over 2}\theta_\eps, E_\eps, B_\eps\right)$ can be restricted to the constraints $3\theta_\eps-2\rho_\eps=0$ and $\rot u_\eps + B_\eps=0$, i.e.\ provided we can find a decomposition
\begin{equation}\label{oscillations decomposition}
	\begin{pmatrix}
		\rho_\eps
		\\
		u_\eps
		\\
		\sqrt{3\over 2} \theta_\eps
		\\
		E_\eps
		\\
		B_\eps
	\end{pmatrix}
	=
	\begin{pmatrix}
		\bar \rho_\eps
		\\
		\bar u_\eps
		\\
		\sqrt{3\over 2} \bar \theta_\eps
		\\
		\bar E_\eps
		\\
		\bar B_\eps
	\end{pmatrix}
	+
	\begin{pmatrix}
		\tilde \rho_\eps
		\\
		\tilde u_\eps
		\\
		\sqrt{3\over 2} \tilde \theta_\eps
		\\
		\tilde E_\eps
		\\
		\tilde B_\eps
	\end{pmatrix},
\end{equation}
such that $\left(\bar \rho_\eps, \bar u_\eps, \sqrt{3\over 2}\bar \theta_\eps, \bar E_\eps, \bar B_\eps\right)$ is relatively compact in the strong topology of $L^2_\mathrm{loc}(dtdx)$, whereas $\left(\tilde \rho_\eps, \tilde u_\eps, \sqrt{3\over 2}\tilde \theta_\eps, \tilde E_\eps, \tilde B_\eps\right)\rightharpoonup 0$ in $L^2_\mathrm{loc}(dtdx)$, with $3\tilde \theta_\eps-2\tilde\rho_\eps=0$ and $\rot \tilde u_\eps + \tilde B_\eps=0$.

In order to obtain such a decomposition, it is very natural to orthogonally project $\left(\rho_\eps, u_\eps, \sqrt{3\over 2}\theta_\eps, E_\eps, B_\eps\right)$ onto the kernel of $W$
\begin{equation*}
		\Ker W
		=
		\set{
		\left(\rho, u, \sqrt{3\over 2}\theta, E, B\right)
% 			\begin{pmatrix}
% 			\rho \\ u \\ \sqrt{3\over 2}\theta \\ E \\ B
% 		\end{pmatrix}
		\in L^2(dx)
		}{
		E=\nabla_x\left(\rho+\theta\right)
		\text{ and }
		u=\rot B
		},
\end{equation*}
% subspace
% \begin{equation*}
% 		H
% 		=
% 		\set{
% 		\left(\rho, u, \sqrt{3\over 2}\theta, E, B\right)
% 		\in L^2(dx)
% 		}{
% 		3\theta-2\rho=0
% 		\text{ and }
% 		\rot u+B=0
% 		}.
% \end{equation*}
and on its orthogonal complement
\begin{equation*}
		\Ker W^\perp
		=
		\set{
		\left(\rho, u, \sqrt{3\over 2}\theta, E, B\right)
% 			\begin{pmatrix}
% 			\rho \\ u \\ \sqrt{3\over 2}\theta \\ E \\ B
% 		\end{pmatrix}
		\in L^2(dx)
		}{
		\frac 32 \theta=\rho=\Div E
		\text{ and }
		B = -\rot u
		}.
\end{equation*}
More precisely, we define
\begin{equation*}%\label{PW-def}
	\begin{pmatrix}
		\bar \rho_\eps
		\\
		\bar u_\eps
		\\
		\sqrt{3\over 2} \bar \theta_\eps
		\\
		\bar E_\eps
		\\
		\bar B_\eps
	\end{pmatrix}
	=
	P_W
	\begin{pmatrix} \rho \\ u \\ \sqrt{3\over 2}\theta \\ E \\ B \end{pmatrix}
	=
	\begin{pmatrix}
		\frac{3}{3-5\Delta_x}\left(\rho-\Div E\right)
		+
		\frac{\Delta_x}{3-5\Delta_x}\left(3\theta-2\rho\right)
		\\
		\frac{\ROT}{1-\Delta_x}\left(\ROT u+  B\right)
		\\
		\frac{2}{3-5\Delta_x}\sqrt{\frac 32}\left(\rho-\Div E\right)
		+
		\frac{1-\Delta_x}{3-5\Delta_x}\sqrt{\frac 32}\left(3\theta-2\rho\right)
		\\
		\frac{5}{3-5\Delta_x}\nabla_x\left(\rho-\Div E\right)
		+
		\frac{1}{3-5\Delta_x}\nabla_x\left(3\theta-2\rho\right)
		\\
		\frac{1}{1-\Delta_x}\left( \ROT u+ B - \nabla_x\Div B\right)
	\end{pmatrix},
\end{equation*}
where $P_W:L^2(dx)\to L^2(dx)$ is the orthogonal projection onto $\Ker W$, and
\begin{equation*}
	\begin{pmatrix}
		\tilde \rho_\eps
		\\
		\tilde u_\eps
		\\
		\sqrt{3\over 2} \tilde \theta_\eps
		\\
		\tilde E_\eps
		\\
		\tilde B_\eps
	\end{pmatrix}
	=
	P_W^\perp
	\begin{pmatrix} \rho \\ u \\ \sqrt{3\over 2}\theta \\ E \\ B \end{pmatrix}
	=
	\begin{pmatrix}
		\frac{3\Div}{3-5\Delta_x}\left(E-\nabla_x\left(\rho+\theta\right)\right)
		\\
		\frac{1}{1-\Delta_x}\left(u-\ROT B - \nabla_x\Div u \right)
		\\
		\frac{2\Div}{3-5\Delta_x}\sqrt{\frac 32}\left(E - \nabla_x\left(\rho+\theta\right)\right)
		\\
		\frac{5\rot\rot}{3-5\Delta_x} E
		+
		\frac{3}{3-5\Delta_x}\left(E - \nabla_x\left(\rho+\theta\right)\right)
		\\
		\frac{\rot}{1-\Delta_x}\left( \rot B - u \right)
	\end{pmatrix},
\end{equation*}
where $P_W^\perp:L^2(dx)\to L^2(dx)$ is the orthogonal projection onto $\Ker W^\perp$. Note that these projections can also be computed explicitly using basic linear algebra in Fourier variables.

Then, recalling that $W$ is antisymmetric so that its range is orthogonal to its kernel, it holds that $P_W W = 0$, whence
\begin{equation*}
	\partial_t
	\begin{pmatrix}
		\bar \rho_\eps
		\\
		\bar u_\eps
		\\
		\sqrt{3\over 2} \bar \theta_\eps
		\\
		\bar E_\eps
		\\
		\bar B_\eps
	\end{pmatrix}
	=O\left(1\right),
\end{equation*}
which implies that $\left(\bar \rho_\eps, \bar u_\eps, \sqrt{3\over 2}\bar \theta_\eps, \bar E_\eps, \bar B_\eps\right)$ is expected to be compact in $t$. Moreover, since $\left(\rho_\eps, u_\eps, \sqrt{3\over 2}\theta_\eps, E_\eps, B_\eps\right)$ belongs to $\Ker W$ in the limit $\eps\to 0$, it is naturally expected that $\left(\tilde \rho_\eps, \tilde u_\eps, \sqrt{3\over 2}\tilde \theta_\eps, \tilde E_\eps, \tilde B_\eps\right)\rightharpoonup 0$. Finally, the constraints $3\tilde \theta_\eps-2\tilde\rho_\eps=0$ and $\rot \tilde u_\eps + \tilde B_\eps=0$ are implied by the fact that $\left(\tilde \rho_\eps, \tilde u_\eps, \sqrt{3\over 2}\tilde \theta_\eps, \tilde E_\eps, \tilde B_\eps\right)$ belongs to $\Ker W^\perp$.

On the whole, we have shown the formal existence of a decomposition \eqref{oscillations decomposition}, which explains why the nonlinear terms \eqref{nonlinear structure} are expected to be weakly stable as $\eps\to 0$, at least formally.

Generally speaking, such methods are called ``compensated compactness'' (following Murat and Tartar \cite{murat, murat2, tartar}), for they consist in compensating the lack of strong compactness in the quadratic terms \eqref{nonlinear structure} by carefully studying the linear structure \eqref{linear structure} of oscillations.

\section{Rigorous filtering of oscillations}

This chapter aims at rendering the preceding developments rigorous. % by establishing a decomposition \eqref{oscillations decomposition} and characterizing its compactness properties.
Thus, the main result concerning the filtering of acoustic and electromagnetic waves in the nonlinear terms through the method of compensated compactness is described in the following proposition.

\begin{prop}\label{time weak stability}
	Let $\left(f_\eps, E_\eps, B_\eps\right)$ be the sequence of renormalized solutions to the scaled one species Vlasov-Maxwell-Boltzmann system \eqref{VMB1} considered in Theorem \ref{NS-WEAKCV} and denote by $\tilde \rho_\eps$, $\tilde u_\eps$ and $\tilde \theta_\eps$ the density, bulk velocity and temperature associated with the renormalized fluctuations $g_\eps \gamma_\eps \chi\left( {|v|^2\over K_\eps} \right)$ employed in Proposition \ref{approx1-prop}. In accordance with Lemma \ref{L2-lem}, denote by
		\begin{equation*}
				\rho,u,\theta,E,B\in L^\infty\left(dt;L^2\left(dx\right)\right),
		\end{equation*}
	any joint limit points of the families $\tilde \rho_\eps$, $\tilde u_\eps$, $\tilde \theta_\eps$, $E_\eps$ and $B_\eps$, respectively.
	
	Then, as $\eps\to 0$, one has the weak stability of nonlinear terms
	\begin{equation}\label{nonlinear convergence}
		\begin{aligned}
			P\left( \nabla_x\cdot \left(\tilde u_\eps\otimes \tilde u_\eps\right) - \tilde \rho_\eps E_\eps - \tilde u_\eps\wedge B_\eps \right)
			& \rightharpoonup
			P\left( \nabla_x\cdot \left(u\otimes u\right) - \rho E - u\wedge B \right),
			\\
			\frac 52 \nabla_x\cdot\left(\tilde u_\eps\tilde \theta_\eps\right) - \tilde u_\eps \cdot E_\eps
			& \rightharpoonup
			\frac 52 \nabla_x\cdot\left(u\theta\right) - u \cdot E,
		\end{aligned}
	\end{equation}
	in the sense of distributions (where we only consider smooth compactly supported solenoidal test functions).
\end{prop}

\begin{proof}
	First of all, it is to be emplasized that compactness in $x$ is not an issue at all here, for none of the nonlinear terms in \eqref{nonlinear convergence} involves a product of the electromagnetic fields $E_\eps$ and $B_\eps$ only.
	
	Indeed, from the strong compactness \eqref{xreg-hatmoments} obtained in Chapter \ref{hypoellipticity}
	\begin{equation*}
		\begin{aligned}
			\lim_{|h|\rightarrow 0}\sup_{\eps>0} \left\| \hat \rho_\eps(t,x+h) - \hat \rho_\eps(t,x)\right\|_{L^2_\mathrm{loc}(dtdx)} & = 0, \\
			\lim_{|h|\rightarrow 0}\sup_{\eps>0} \left\| \hat u_\eps(t,x+h) - \hat u_\eps(t,x)\right\|_{L^2_\mathrm{loc}(dtdx)} & = 0, \\
			\lim_{|h|\rightarrow 0}\sup_{\eps>0} \left\| \hat \theta_\eps(t,x+h) - \hat \theta_\eps(t,x)\right\|_{L^2_\mathrm{loc}(dtdx)} & = 0,
		\end{aligned}
	\end{equation*}
	where $\hat \rho_\eps$, $\hat u_\eps$ and $\hat \theta_\eps$ respectively denote the density, bulk velocity and temperature of the renormalized fluctuations $\hat g_\eps$ defined by \eqref{hatg}, and the comparison \eqref{moments-comparison} between $\hat g_\eps$ and $g_\eps \gamma_\eps \chi\left( {|v|^2\over K_\eps}\right)$ established in Chapter \ref{conservation0-chap}
	\begin{equation*}
		\tilde \rho_\eps-\hat \rho_\eps \to 0,
		\quad \tilde u_\eps -\hat u_\eps \to 0
		\quad\text{and}\quad
		\tilde \theta_\eps -\hat \theta_\eps \to 0
		\quad\text{in }L^2_\mathrm{loc}(dtdx) \text{ as }\eps \to 0,
	\end{equation*}
	we deduce that the following local spatial compactness property holds~:
	\begin{equation*}%\label{xreg-moments}
		\begin{aligned}
			\lim_{|h|\rightarrow 0}\sup_{\eps>0} \left\| \tilde \rho_\eps(t,x+h) - \tilde \rho_\eps(t,x)\right\|_{L^2_\mathrm{loc}(dtdx)} & = 0, \\
			\lim_{|h|\rightarrow 0}\sup_{\eps>0} \left\| \tilde u_\eps(t,x+h) - \tilde u_\eps(t,x)\right\|_{L^2_\mathrm{loc}(dtdx)} & = 0, \\
			\lim_{|h|\rightarrow 0}\sup_{\eps>0} \left\| \tilde \theta_\eps(t,x+h) - \tilde \theta_\eps(t,x)\right\|_{L^2_\mathrm{loc}(dtdx)} & = 0.
		\end{aligned}
	\end{equation*}
	
	In particular, denoting by $\tilde \rho_\eps^\delta$, $\tilde u_\eps^\delta$, $\tilde \theta_\eps^\delta$, $\rho^\delta$, $u^\delta$ and $\theta^\delta$ the respective spatial convolutions of $\tilde \rho_\eps$, $\tilde u_\eps$, $\tilde \theta_\eps$, $\rho$, $u$ and $\theta$ with a smooth compactly supported mollifier $\frac{1}{\delta^3}\chi\left(\frac{x}{\delta}\right)$, where $\delta>0$ and $\chi\in C_c^\infty\left(\mathbb{R}^3\right)$, with $\int_{\mathbb{R}^3}\chi(x)dx = 1$, we see that it is possible to replace each $\tilde \rho_\eps$, $\tilde u_\eps$ and $\tilde \theta_\eps$ in \eqref{nonlinear convergence} by $\tilde \rho_\eps^\delta$, $\tilde u_\eps^\delta$ and $\tilde \theta_\eps^\delta$, respectively, producing remainders that are uniformly small in $L^1_\mathrm{loc}\left(dt;W^{-1,1}_\mathrm{loc}\left(dx\right)\right)$ as $\delta\to 0$. This reduces \eqref{nonlinear convergence} to showing the nonlinear convergence
	\begin{equation*}
		\begin{aligned}
			P\left( \nabla_x\cdot \left(\tilde u_\eps^\delta\otimes \tilde u_\eps^\delta\right) - \tilde \rho_\eps^\delta E_\eps - \tilde u_\eps^\delta\wedge B_\eps \right)
			& \rightharpoonup
			P\left( \nabla_x\cdot \left(u^\delta\otimes u^\delta\right) - \rho^\delta E - u^\delta\wedge B \right),
			\\
			\frac 52 \nabla_x\cdot\left(\tilde u_\eps^\delta\tilde \theta_\eps^\delta\right) - \tilde u_\eps^\delta \cdot E_\eps
			& \rightharpoonup
			\frac 52 \nabla_x\cdot\left(u^\delta\theta^\delta\right) - u^\delta \cdot E,
		\end{aligned}
	\end{equation*}
	in the sense of distributions.
	
	Then, denoting by $E_\eps^\delta$, $B_\eps^\delta$, $E^\delta$ and $B^\delta$ the respective spatial convolutions of $E_\eps$, $B_\eps$, $E$ and $B$ with the mollifier $\frac{1}{\delta^3}\chi\left(\frac{x}{\delta}\right)$, we notice, since we are only seeking to establish a convergence in the sence of distributions, that we may also replace $E_\eps$ and $B_\eps$ by $E_\eps^\delta$ and $B_\eps^\delta$, respectively, thus further reducing the proof of the present proposition to establishing the nonlinear convergence, for any fixed $\delta>0$,
	\begin{equation*}
		\begin{aligned}
			P\left( \nabla_x\cdot \left(\tilde u_\eps^\delta\otimes \tilde u_\eps^\delta\right) - \tilde \rho_\eps^\delta E_\eps^\delta - \tilde u_\eps^\delta\wedge B_\eps^\delta \right)
			& \rightharpoonup
			P\left( \nabla_x\cdot \left(u^\delta\otimes u^\delta\right) - \rho^\delta E^\delta - u^\delta\wedge B^\delta \right),
			\\
			\frac 52 \nabla_x\cdot\left(\tilde u_\eps^\delta\tilde \theta_\eps^\delta\right) - \tilde u_\eps^\delta \cdot E_\eps^\delta
			& \rightharpoonup
			\frac 52 \nabla_x\cdot\left(u^\delta\theta^\delta\right) - u^\delta \cdot E^\delta,
		\end{aligned}
	\end{equation*}
	in the sense of distributions.

	Now, according to Proposition \ref{approx1-prop}, coupling the linear part of the macroscopic equations derived therein with Maxwell's equations, one has the following acoustic-electromagnetic wave system, for any fixed $\delta>0$~:
	\begin{equation}\label{wave system}
		\d_t
		\begin{pmatrix}
			\tilde \rho_\eps^\delta \\ \tilde u_\eps^\delta\\ \sqrt{3\over 2}\tilde \theta_\eps^\delta \\ E_\eps^\delta \\ B_\eps^\delta
		\end{pmatrix}
		+\frac1\eps
		W
		\begin{pmatrix}
			\tilde \rho_\eps^\delta\\ \tilde u_\eps^\delta\\ \sqrt{3\over 2}\tilde \theta_\eps^\delta \\ E_\eps^\delta \\ B_\eps^\delta
		\end{pmatrix}
		=O\left(1\right)_{L^1_\mathrm{loc}\left(dt;L^\infty_\mathrm{loc}\left(dx\right)\right)},
	\end{equation}
	where the wave operator $W$ is defined in \eqref{W-def}. In particular, it holds that
	\begin{equation*}
		\begin{aligned}
			\partial_t\left(2\tilde \rho_\eps^\delta-3\tilde \theta_\eps^\delta\right) & =O\left(1\right)_{L^1_\mathrm{loc}\left(dtdx\right)},
			\\
			\partial_t\left(\rot \tilde u_\eps^\delta + B_\eps^\delta\right) & =O\left(1\right)_{L^1_\mathrm{loc}\left(dtdx\right)},
		\end{aligned}
	\end{equation*}
	whence $2\tilde \rho_\eps^\delta-3\tilde \theta_\eps^\delta$ and $\rot \tilde u_\eps^\delta + B_\eps^\delta$ are relatively compact in the strong topology of $L^2_\mathrm{loc}(dtdx)$ (in both variables $t$ and $x$). It follows that
	\begin{equation}\label{compensated strong convergence}
		\begin{aligned}
			2\tilde \rho_\eps^\delta-3\tilde \theta_\eps^\delta
			& \to
			2\rho^\delta-3\theta^\delta,
			\\
			\rot \tilde u_\eps^\delta + B_\eps^\delta
			& \to
			\rot u^\delta + B^\delta,
		\end{aligned}
	\end{equation}
	in $L^2_\mathrm{loc}(dtdx)$.

	Next, we reproduce here rigorously the formal identities \eqref{compensated 1} and \eqref{compensated 2}, which yields (for fixed $\eps$, notice that $\tilde \rho_\eps^\delta$, $\tilde u_\eps^\delta$ and $\tilde \theta_\eps^\delta$ are now differentiable once in $t$ with a derivative lying in $L^1_\mathrm{loc}\left(dt;L^\infty_\mathrm{loc}\left(dx\right)\right)$, due to \eqref{wave system})
	\begin{equation*}
		\begin{aligned}
			P & \left( \nabla_x\cdot \left(\tilde u_\eps^\delta \otimes \tilde u_\eps^\delta \right) - \tilde \rho_\eps^\delta  E_\eps^\delta - \tilde u_\eps^\delta \wedge B_\eps^\delta \right)
			\\
			& =
			P\left( \tilde u_\eps^\delta \left(\eps\partial_t\tilde \rho_\eps^\delta + \Div \tilde u_\eps^\delta\right)
			+ \tilde \rho_\eps^\delta\left(\eps\partial_t \tilde u_\eps^\delta + \nabla_x\left(\tilde \rho_\eps^\delta+\tilde \theta_\eps^\delta\right) - E_\eps^\delta \right) \right)
			\\
			& - P\left( \frac{2\tilde \rho_\eps^\delta-3\tilde \theta_\eps^\delta}{5}\nabla_x\left(\tilde \rho_\eps^\delta+\tilde \theta_\eps^\delta\right)
			+ \tilde u_\eps^\delta\wedge\left(\rot \tilde u_\eps^\delta + B_\eps^\delta\right) + \eps\partial_t\left(\tilde \rho_\eps^\delta \tilde u_\eps^\delta\right) \right),
		\end{aligned}
	\end{equation*}
	and
	\begin{equation*}
		\begin{aligned}
			\frac 52 & \nabla_x\cdot\left(\tilde u_\eps^\delta\tilde \theta_\eps^\delta\right) - \tilde u_\eps^\delta \cdot E_\eps^\delta
			\\
			& = \frac 52 \tilde \theta_\eps^\delta \left(\frac 32\eps\partial_t\tilde \theta_\eps^\delta + \Div \tilde u_\eps^\delta\right) + \tilde u_\eps^\delta\cdot\left(\eps\partial_t \tilde u_\eps^\delta + \nabla_x\left(\tilde \rho_\eps^\delta+\tilde \theta_\eps^\delta\right) - E_\eps^\delta\right)
			\\
			& + \tilde u_\eps^\delta\cdot\nabla_x\left(\frac 32\tilde \theta_\eps^\delta-\tilde \rho_\eps^\delta\right) - \frac{15}{8}\eps\partial_t\left(\tilde \theta_\eps^\delta\right)^2 - \frac 12 \eps\partial_t \left|\tilde u_\eps^\delta\right|^2.
		\end{aligned}
	\end{equation*}
	Consequently, since \eqref{wave system} implies that, for fixed $\delta>0$,
	\begin{equation*}
		\begin{aligned}
			\eps\partial_t\tilde \rho_\eps^\delta + \Div \tilde u_\eps^\delta & = O(\eps)_{L^1_\mathrm{loc}\left(dt;L^\infty_\mathrm{loc}(dx)\right)},
			\\
			\eps\partial_t \tilde u_\eps^\delta + \nabla_x\left(\tilde \rho_\eps^\delta+\tilde \theta_\eps^\delta\right) - E_\eps^\delta & = O(\eps)_{L^1_\mathrm{loc}\left(dt;L^\infty_\mathrm{loc}(dx)\right)},
			\\
			\frac 32\eps\partial_t\tilde \theta_\eps^\delta + \Div \tilde u_\eps^\delta & = O(\eps)_{L^1_\mathrm{loc}\left(dt;L^\infty_\mathrm{loc}(dx)\right)},
		\end{aligned}
	\end{equation*}
	we deduce, in view of the strong convergences \eqref{compensated strong convergence}, that
	\begin{equation*}
		\begin{aligned}
			P\big( \nabla_x\cdot \left(\tilde u_\eps^\delta\otimes \tilde u_\eps^\delta\right) - \tilde \rho_\eps^\delta E_\eps^\delta & - \tilde u_\eps^\delta\wedge B_\eps^\delta \big)
			\\
			& \rightharpoonup
			- P\left( \frac{2\rho^\delta-3\theta^\delta}{5}\nabla_x\left(\rho^\delta+\theta^\delta\right)
			+ u^\delta\wedge\left(\rot u^\delta + B^\delta\right) \right),
			% \\
% 			& =
% 			P\left( \nabla_x\cdot \left(u^\delta\otimes u^\delta\right) - \rho^\delta E^\delta - u^\delta\wedge B^\delta \right),
			\\
			\frac 52 \nabla_x\cdot\left(\tilde u_\eps^\delta\tilde \theta_\eps^\delta\right) - \tilde u_\eps^\delta \cdot E_\eps^\delta
			& \rightharpoonup
			u^\delta\cdot\nabla_x\left(\frac 32\theta^\delta-\rho^\delta\right),
			% \\
% 			& =
% 			\frac 52 \nabla_x\cdot\left(u^\delta\theta^\delta\right) - u^\delta \cdot E^\delta,
		\end{aligned}
	\end{equation*}
	in the sense of distributions. Finally, using from Proposition \ref{weak-comp} that $\Div u^\delta = 0$ and $E^\delta=\nabla_x\left(\rho^\delta+\theta^\delta\right)$ to deduce
	\begin{equation*}
		\begin{aligned}
			\frac{2\rho^\delta-3\theta^\delta}{5}\nabla_x\left(\rho^\delta+\theta^\delta\right)
			& = \rho^\delta\nabla_x\left(\rho^\delta+\theta^\delta\right) - \frac 3{10}\nabla_x \left(\rho^\delta+\theta^\delta\right)^2
			\\
			& = \rho^\delta E^\delta - \frac 3{10}\nabla_x \left(\rho^\delta+\theta^\delta\right)^2,
			\\
			u^\delta\wedge\rot u^\delta
			& = -\nabla_x\cdot\left(u^\delta\otimes u^\delta\right) + u^\delta\Div u^\delta + \frac 12 \nabla_x\left|u^\delta\right|^2
			\\
			& = -\nabla_x\cdot\left(u^\delta\otimes u^\delta\right) + \frac 12 \nabla_x\left|u^\delta\right|^2,
			\\
			u^\delta\cdot\nabla_x\left(\frac 32\theta^\delta-\rho^\delta\right)
			& =
			\frac 52 \nabla_x \left(u^\delta\theta^\delta\right)
			-u^\delta\cdot\nabla_x\left(\rho^\delta+\theta^\delta\right)
			- \frac 52\theta^\delta\Div u^\delta
			\\
			& =
			\frac 52 \nabla_x\cdot\left(u^\delta\theta^\delta\right) - u^\delta \cdot E^\delta,
		\end{aligned}
	\end{equation*}
	concludes the proof of the proposition.
\end{proof}

%% file: Grad0.tex
\chapter{Grad's moment method}\label{grad}

We are now in a position to proceed to the demonstration of Theorem \ref{NS-WEAKCV}. Generally speaking, the formal approach to this proof follows the method of Grad from \cite{grad, grad2}, which consists in studying the moments of Boltzmann equations as the densities remain close to statistical equilibrium through the formal Hilbert's expansions from \cite{hilbert}. In our fully rigorous setting, since we are considering renormalized solutions of the Vlasov-Maxwell-Boltzmann system \eqref{VMB1} (which, we recall, are not known to exist in general), our method of proof proceeds through the asymptotic analysis of renormalized moments satisfying approximate macroscopic conservation laws leading to the incompressible quasi-static Navier-Stokes-Fourier-Maxwell-Poisson system \eqref{NSFMP 2}.

We insist on the fact that the result we are about to establish holds globally in time and does not require any additional assumption on the initial data, neither on the initial velocity profile, nor on the initial thermodynamic fields, nor on the corresponding solution to the limiting system.

\section{Proof of Theorem \ref{NS-WEAKCV}}

Most of the difficult steps of this proof have been performed in the preceding chapters. We therefore only have to appropriately gather previous results together.

% \section{From consistency to convergence}
%
%
% Gathering all the previous results together, we are now able to prove the convergence to the electrostatic Navier-Stokes-Maxwell equations (\eqref{NSFMP 2}) starting from the scaled Vlasov-Maxwell-Boltzmann system (\eqref{VMB1}).
%
% Due to the scaling of the transport in this regime, we can indeed discard fast spatial oscillations using the hypoellipticity mechanism presented in Chapter 7.
% Fast time oscillations will be handled in the present chapter using compensated compactness. We therefore get the stability of nonlinear herms under the only size condition on the initial relative entropy.
%
% In particular, the convergence result
%
% \section{Consistency}

\subsection{Weak convergence of fluctuations, collision integrands and electromagnetic fields}

Thus, we are considering here a family of renormalized solutions $(f_\eps,E_\eps,B_\eps)$ to the scaled one species Vlasov-Maxwell-Boltzmann system \eqref{VMB1} satisfying the scaled entropy inequality \eqref{entropy1}.

According to Lemmas \ref{L1-lem} and \ref{L2-lem}, the corresponding families of fluctuations $g_\eps$ and renormalized fluctuations $\hat g_\eps$ are weakly compact in $L^1_\mathrm{loc}\left(dtdx;L^1\left(\left(1+|v|^2\right)Mdv\right)\right)$ and $L^2_\mathrm{loc}\left(dt;L^2\left(Mdxdv\right)\right)$, respectively, while, in view of Lemma \ref{L2-qlem}, the corresponding collision integrands $\hat q_\eps$ are weakly compact in $L^2\left(MM_*dtdxdvdv_*d\sigma\right)$. Thus, using Lemma \ref{L1-lem} again and the decomposition \eqref{fluct-decomposition}, we know that there exist $g\in L^\infty\left(dt;L^2\left(Mdxdv\right)\right)$, $(E,B)\in L^\infty\left(dt;L^2(dx)\right)$ and $q\in L^2\left(MM_*dtdxdvdv_*d\sigma\right)$, such that, up to extraction of subsequences,
\begin{equation*}
	\begin{aligned}
		g_\eps & \rightharpoonup g & & \text{in }L^1_\mathrm{loc}\left(dtdx;L^1\left(\left(1+|v|^2\right)Mdv\right)\right),
		\\
		\hat g_\eps & \stackrel{*}{\rightharpoonup} g & & \text{in }L^\infty\left(dt;L^2\left(Mdxdv\right)\right),
		\\
		\left(E_\eps,B_\eps\right) & \stackrel{*}{\rightharpoonup} \left(E,B\right) & & \text{in }L^\infty\left(dt;L^2\left(dx\right)\right),
		\\
		\hat q_\eps & \rightharpoonup q & & \text{in }L^2\left(MM_*dtdxdvdv_*d\sigma\right),
	\end{aligned}
\end{equation*}
as $\eps\to 0$. Therefore, one also has the weak convergence of the density $\rho_\eps$, bulk velocity $u_\eps$ and temperature $\theta_\eps$ corresponding to $g_\eps$~:
\begin{equation*}
	\rho_\eps \rightharpoonup \rho,
	\quad u_\eps \rightharpoonup u
	\quad \text{and}\quad
	\theta_\eps \rightharpoonup \theta
	\quad\text{in }L^1_\mathrm{loc}(dtdx) \text{ as }\eps \to 0,
\end{equation*}
where $\rho, u, \theta \in L^\infty\left(dt;L^2(dx)\right)$ are, respectively, the density, bulk velocity and temperature corresponding to $g$. In fact, Lemma \ref{relaxation-control} implies that
\begin{equation}\label{infinitesimal maxwellian proof}
	g=\Pi g = \rho + u\cdot v + \theta\left(\frac{|v|^2}{2}-\frac 32\right).
\end{equation}

Based on the uniform initial bound \eqref{init-fluctuation}, a very slight modification of Lemma \ref{L1-lem} allows us to deduce similar weak compactness properties on the initial data. Thus, the initial fluctuations $g_\eps^\mathrm{in}$ are weakly relatively compact in $L^1_\mathrm{loc}\left(dx;L^1\left(\left(1+|v|^2\right)Mdv\right)\right)$ and so, up to further extraction of subsequences, we may also assume that there are $g^\mathrm{in}_0\in L^2\left(Mdxdv\right)$ and $(E^\mathrm{in}_0,B^\mathrm{in}_0)\in L^2(dx)$, such that, up to extraction of subsequences,
\begin{equation*}
	\begin{aligned}
		g_\eps^\mathrm{in} & \rightharpoonup g^\mathrm{in}_0 & & \text{in }L^1_\mathrm{loc}\left(dx;L^1\left(\left(1+|v|^2\right)Mdv\right)\right),
		\\
		\left(E_\eps^\mathrm{in},B_\eps^\mathrm{in}\right) & \rightharpoonup \left(E_0^\mathrm{in},B_0^\mathrm{in}\right) & & \text{in }L^2\left(dx\right),
	\end{aligned}
\end{equation*}
as $\eps\to 0$. Therefore, one also has the weak convergence of the initial density $\rho_\eps^\mathrm{in}$, bulk velocity $u_\eps^\mathrm{in}$ and temperature $\theta_\eps^\mathrm{in}$ corresponding to $g_\eps^\mathrm{in}$~:
\begin{equation*}
	\rho_\eps^\mathrm{in} \rightharpoonup \rho^\mathrm{in}_0,
	\quad u_\eps^\mathrm{in} \rightharpoonup u^\mathrm{in}_0
	\quad \text{and}\quad
	\theta_\eps^\mathrm{in} \rightharpoonup \theta^\mathrm{in}_0
	\quad\text{in }L^1_\mathrm{loc}(dx) \text{ as }\eps \to 0,
\end{equation*}
where $\rho^\mathrm{in}_0, u^\mathrm{in}_0, \theta^\mathrm{in}_0 \in L^2(dx)$ are, respectively, the initial density, bulk velocity and temperature corresponding to $g^\mathrm{in}_0$. Note that the infinitesimal Maxwellian form \eqref{infinitesimal maxwellian proof} does not necessarily hold for the initial data $g_0^\mathrm{in}$.

\subsection{Constraint equations, Maxwell's system and energy inequality}

In view of Proposition \ref{weak-comp}, we already know that the limiting thermodynamic fields $\rho$, $u$ and $\theta$ satisfy the incompressibility and Boussinesq relations
\begin{equation}\label{w-constraint1}
	\Div u = 0, \qquad \nabla_x\left(\rho+\theta\right)-E = 0.
\end{equation}
Furthermore, the discussion in Section \ref{limit maxwell} shows that the limiting electromagnetic field satisfies the electrostatic approximation of Maxwell's equations~:
\begin{equation}\label{w-constraint2}
	\ROT E = 0, \qquad \Div E = \rho, \qquad  \ROT B = u, \qquad \Div B = 0.
\end{equation}
By passing to the weak limit in the initial Gauss' laws \eqref{initial gauss}, one also has initially that
\begin{equation*}
	\Div E_0^\mathrm{in}=\rho_0^\mathrm{in},\qquad \Div B_0^\mathrm{in} = 0.
\end{equation*}

As for the energy bound, Proposition \ref{energy ineq 1} states that, for almost every $t\geq 0$,
\begin{equation*}
	\begin{aligned}
		& \frac 12\left(\left\|\rho\right\|_{L^2_x}^2 + \left\|u\right\|_{L^2_x}^2
		+ \frac 32\left\|\theta\right\|_{L^2_x}^2 + \left\|E\right\|_{L^2_x}^2
		+ \left\|B\right\|_{L^2_x}^2 \right)(t)
		\\
		& \hspace{30mm} +
		\int_0^t \left(\mu
		\left\|\nabla_x u\right\|_{L^2_x}^2
		+ \frac 52\kappa
		\left\|\nabla_x\theta\right\|_{L^2_x}^2\right)(s) ds
		\leq C^\mathrm{in},
	\end{aligned}
\end{equation*}
where the viscosity $\mu>0$ and thermal conductivity $\kappa>0$ are defined by \eqref{mu kappa}. In particular, it holds that
\begin{equation*}
		\begin{aligned}
			\left(\rho, u, \theta, B\right) & \in L^\infty \left( [0,\infty), dt ; L^2\left(\mathbb{R}^3, dx\right)\right), \\
			\left(u,\theta\right) & \in L^2\left([0,\infty), dt ; \dot H^1\left(\mathbb{R}^3, dx\right)\right).
		\end{aligned}
\end{equation*}

This energy bound can be improved to the actual energy inequality \eqref{energy} provided some well-preparedness of the initial data is assumed. This is discussed in the few remarks following the statement of Theorem \ref{NS-WEAKCV}.

% We have seen in Chapter 4 that the entropy inequality provides uniform bounds (and thus weak compactness) on the sequences of renormalized fluctuations
% $$ g_\eps \gamma_\eps= \frac1\eps  (\Gamma (G_\eps)-1)$$
% and corresponding thermodynamic fields
% $$\rho_\eps=\int Mg_\eps \gamma_\eps dv , \qquad u_\eps=\int Mg_\eps \gamma_\eps vdv,\qquad \theta_\eps =\frac12 \int Mg_\eps \gamma_\eps (|v|^2-3) dv $$
% as well as on the sequences of electric and magnetic fields
% $E_\eps, B_\eps$.
%
% \medskip
% These variables satisfy approximately the following linear constraint relations or dynamical equations
% $$
% \begin{aligned}
% \Div u_\eps = O(\eps), \quad \nabla(\rho_\eps+\theta_\eps)-E_\eps = O(\eps)\,,\\
% \ROT E_\eps =O(\eps) ,\quad  \Div E_\eps = \rho_\eps +O(\eps),\\
% -\ROT B_\eps =u_\eps +O(\eps), \quad \Div B_\eps = 0\,.
% \end{aligned}
% $$
% which are stable under weak convergence.

% Therefore any joint limit point $(\rho, u,\theta,E,B)$ satisfies the {\bf incompressibility and Boussinesq relations}
% \begin{equation}
% \label{w-constraint1}
% \Div u = 0, \quad \nabla(\rho+\theta)-E = 0,
% \end{equation}
% together with the {\bf electrostatic Maxwell equations}
% \begin{equation}
% \label{w-constraint2}
% \begin{aligned}
%  \ROT E = 0, \quad \Div E = \rho, \\
%   -\ROT B = u \hbox{ and } \Div B = 0\,.
%   \end{aligned}
% \end{equation}

\subsection{Evolution equations}
% \subsection{Approximate conservation laws}

We move on now to the rigorous derivation of the asymptotic macroscopic evolution equations.

We know from Chapter \ref{conservation0-chap} that some approximate macroscopic evolution equations, which look like the Navier-Stokes-Fourier system with electromagnetic forces, are satisfied up to a remainder which is small in some distribution space. More precisely, according to Proposition \ref{approx1-prop}, defining the macroscopic variables $\tilde\rho_\eps$, $\tilde u_\eps$ and $\tilde \theta_\eps$ as the density, bulk velocity and temperature, respectively, corresponding to the renormalized fluctuations $g_\eps\gamma_\eps\chi\left(\frac{|v|^2}{K_\eps}\right)$ used therein, it holds that
\begin{equation}\label{approximate1}
	\left\{
	\begin{aligned}
		\d_t \tilde\rho_\eps + \frac1\eps \nabla_x\cdot\tilde u_\eps
		& = R_{\eps,1},
		\\
		\d_t\tilde u_\eps
		+ \nabla_x \cdot \left(\tilde u_\eps \otimes\tilde u_\eps -\frac{|\tilde u_\eps|^2}{3} \operatorname{Id} - \int_{\mathbb{R}^3\times\mathbb{R}^3\times\mathbb{S}^2} \hat q_\eps \tilde \phi MM_* dvdv_*d\sigma \right) \hspace{-60mm} &
		\\
		& = - \frac 1\eps \nabla_x\left(\tilde\rho_\eps+\tilde\theta_\eps\right)
		+ \frac1\eps E_\eps
		+\tilde\rho_\eps E_\eps
		+\tilde u_\eps \wedge B_\eps + R_{\eps,2},
		\\
		\d_t \left(\frac 32\tilde\theta_\eps-\tilde\rho_\eps\right) + \nabla_x \cdot \left( \frac52\tilde u_\eps\tilde \theta_\eps - \int_{\mathbb{R}^3\times\mathbb{R}^3\times\mathbb{S}^2} \hat q_\eps \tilde \psi MM_* dvdv_*d\sigma\right) \hspace{-60mm} &
		\\
		& =\tilde u_\eps \cdot E_\eps+ R_{\eps,3},
	\end{aligned}
	\right.
\end{equation}
where $\tilde\phi$ and $\tilde\psi$ are defined by \eqref{phi-psi-def} and \eqref{phi-psi-def inverses}, and the remainders $R_{\eps,i}$, $i=1,2,3$, converge to $0$ in $L^1_\mathrm{loc}\left(dt;W^{-1,1}_\mathrm{loc}\left(dx\right)\right)$.

Since, up to further extraction of subsequences, $\gamma_\eps \chi\left( {|v|^2\over K_\eps}\right)$ converges almost everywhere towards $1$, $g_\eps$ is weakly compact in $L^1_\mathrm{loc}\left(dtdx;L^1\left(\left(1+|v|^2\right)Mdv\right)\right)$ and $g_\eps\gamma_\eps$ is uniformly bounded in $L^\infty\left(dt;L^2\left(Mdxdv\right)\right)$, we deduce, by the Product Limit Theorem, that
\begin{equation*}
	g_\eps \gamma_\eps \chi\left( {|v|^2\over K_\eps}\right)
	\stackrel{*}{\rightharpoonup} g \quad \text{in }L^\infty\left(dt;L^2\left(Mdxdv\right)\right).
\end{equation*}
In particular, one also has the convergence of the renormalized moments
\begin{equation*}
	\tilde \rho_\eps \stackrel{*}{\rightharpoonup} \rho,
	\quad \tilde u_\eps \stackrel{*}{\rightharpoonup} u
	\quad\text{and}\quad
	\tilde \theta_\eps \stackrel{*}{\rightharpoonup} \theta
	\quad\text{in }L^\infty\left(dt;L^2(dx)\right),
\end{equation*}
and the same argument yields the convergence of the initial renormalized moments
\begin{equation*}
	\tilde \rho_\eps^\mathrm{in} \rightharpoonup \rho_0^\mathrm{in},
	\quad \tilde u_\eps^\mathrm{in} \rightharpoonup u_0^\mathrm{in}
	\quad\text{and}\quad
	\tilde \theta_\eps^\mathrm{in} \rightharpoonup \theta_0^\mathrm{in}
	\quad\text{in }L^2(dx),
\end{equation*}
where $\tilde \rho_\eps^\mathrm{in}=\tilde\rho_\eps(t=0)$, $\tilde u_\eps^\mathrm{in}=\tilde u_\eps(t=0)$ and $\tilde \theta_\eps^\mathrm{in}=\tilde\theta_\eps(t=0)$.

Next we consider the magnetic potentials $A_\eps, A \in L^\infty\left(dt;\dot H^1\left(dx\right)\right)$ and $A_\eps^\mathrm{in}, A_0^\mathrm{in} \in \dot H^1\left(dx\right)$ in the Coulomb gauge defined by
\begin{equation*}
	A_\eps = \frac{\rot}{-\Delta_x}B_\eps,
	\qquad
	A = \frac{\rot}{-\Delta_x}B
	\qquad\text{and}\qquad
	A_\eps^\mathrm{in} = \frac{\rot}{-\Delta_x}B_\eps^\mathrm{in},
	\qquad
	A_0^\mathrm{in} = \frac{\rot}{-\Delta_x}B_0^\mathrm{in},
\end{equation*}
so that
\begin{equation*}
	\begin{aligned}
		B_\eps & = \rot A_\eps, & \Div A_\eps & = 0,
		\\
		B & = \rot A, & \Div A & = 0,
		\\
		B_\eps^\mathrm{in} & = \rot A_\eps^\mathrm{in}, & \Div A_\eps^\mathrm{in} & = 0,
		\\
		B_0^\mathrm{in} & = \rot A_\eps^\mathrm{in}, & \Div A_0^\mathrm{in} & = 0.
	\end{aligned}
\end{equation*}
Faraday's equation from \eqref{VMB1} can then be recast as
\begin{equation}\label{faraday potential}
	\eps\partial_t A_\eps + \frac{\nabla_x \Div}{-\Delta_x} E_\eps + E_\eps = 0.
\end{equation}
Further note that
\begin{equation*}
	A_\eps\rightharpoonup A
	\qquad\text{and}\qquad
	A_\eps^\mathrm{in}\rightharpoonup A_0^\mathrm{in},
\end{equation*}
in the sense of distributions.

Now, incorporating the preceding relation \eqref{faraday potential} into the evolution equation for $\tilde u_\eps$ in \eqref{approximate1}, we obtain the following system of evolution equations~:
\begin{equation*}
	\begin{aligned}
		\d_t\left(\tilde u_\eps+ A_\eps\right)
		+ \nabla_x \cdot \left(\tilde u_\eps \otimes\tilde u_\eps -\frac{|\tilde u_\eps|^2}{3} \operatorname{Id} - \int_{\mathbb{R}^3\times\mathbb{R}^3\times\mathbb{S}^2} \hat q_\eps \tilde \phi MM_* dvdv_*d\sigma \right) \hspace{-70mm} &
		\\
		& = - \frac 1\eps \nabla_x\left(\tilde\rho_\eps+\tilde\theta_\eps + \frac{\Div}{-\Delta_x}E_\eps\right)
		+\tilde\rho_\eps E_\eps
		+\tilde u_\eps \wedge B_\eps + R_{\eps,2},
		\\
		\d_t \left(\frac 32\tilde\theta_\eps-\tilde\rho_\eps\right) + \nabla_x \cdot \left( \frac52\tilde u_\eps\tilde \theta_\eps - \int_{\mathbb{R}^3\times\mathbb{R}^3\times\mathbb{S}^2} \hat q_\eps \tilde \psi MM_* dvdv_*d\sigma\right) \hspace{-70mm} &
		\\
		& =\tilde u_\eps \cdot E_\eps+ R_{\eps,3},
	\end{aligned}
\end{equation*}
whose weak formulation is given by
\begin{equation*}
	\begin{aligned}
		-\int_{\mathbb{R}^3}\left(\tilde u_\eps^\mathrm{in}+ A_\eps^\mathrm{in}\right)\cdot\varphi(t=0)dx
		-\int_{[0,\infty)\times\mathbb{R}^3}\left(\tilde u_\eps+ A_\eps\right)\cdot \partial_t\varphi dtdx
		\hspace{-42mm}&
		\\
		- \int_{[0,\infty)\times\mathbb{R}^3} \left(\tilde u_\eps \otimes\tilde u_\eps - \int_{\mathbb{R}^3\times\mathbb{R}^3\times\mathbb{S}^2} \hat q_\eps \tilde \phi MM_* dvdv_*d\sigma \right)
		:\nabla_x\varphi dtdx
		\hspace{-42mm}&
		\\
		& =\int_{[0,\infty)\times\mathbb{R}^3} \left(\tilde\rho_\eps E_\eps+\tilde u_\eps \wedge B_\eps\right) \cdot \varphi dtdx + o(1),
		\\
		-\int_{\mathbb{R}^3}\left(\frac 32\tilde\theta_\eps^\mathrm{in}-\tilde\rho_\eps^\mathrm{in}\right)\chi(t=0) dx
		-\int_{[0,\infty)\times\mathbb{R}^3}\left(\frac 32\tilde\theta_\eps-\tilde\rho_\eps\right)\partial_t\chi dtdx
		\hspace{-42mm}&
		\\
		-\int_{[0,\infty)\times\mathbb{R}^3}
		\left( \frac52\tilde u_\eps\tilde \theta_\eps - \int_{\mathbb{R}^3\times\mathbb{R}^3\times\mathbb{S}^2} \hat q_\eps \tilde \psi MM_* dvdv_*d\sigma\right)
		\cdot\nabla_x\chi dtdx
		\hspace{-42mm}&
		\\
		& =\int_{[0,\infty)\times\mathbb{R}^3}\tilde u_\eps \cdot E_\eps \chi dtdx + o(1),
	\end{aligned}
\end{equation*}
where $\varphi(t,x)\in C_c^\infty\left([0,\infty)\times\mathbb{R}^3;\mathbb{R}^3\right)$ and $\chi(t,x)\in C_c^\infty\left([0,\infty)\times\mathbb{R}^3;\mathbb{R}\right)$ are test functions such that $\Div \varphi =0$.

By the weak stability result stated in Proposition \ref{time weak stability}, we can then pass to limit $\eps\to 0$ in the above weak formulation, thus yielding the following asymptotic system~:
\begin{equation}\label{weak formulation limiting system}
	\begin{aligned}
		-\int_{\mathbb{R}^3}\left(u_0^\mathrm{in}+ A_0^\mathrm{in}\right)\cdot\varphi(t=0)dx
		-\int_{[0,\infty)\times\mathbb{R}^3}\left(u+ A\right)\cdot \partial_t\varphi dtdx
		\hspace{-30mm}&
		\\
		- \int_{[0,\infty)\times\mathbb{R}^3} \left(u \otimes u - \int_{\mathbb{R}^3\times\mathbb{R}^3\times\mathbb{S}^2} q \tilde \phi MM_* dvdv_*d\sigma \right)
		:\nabla_x\varphi dtdx
		\hspace{-30mm}&
		\\
		& =\int_{[0,\infty)\times\mathbb{R}^3} \left(\rho E + u \wedge B\right) \cdot \varphi dtdx,
		\\
		-\int_{\mathbb{R}^3}\left(\frac 32 \theta_0^\mathrm{in}-\rho_0^\mathrm{in}\right)\chi(t=0) dx
		-\int_{[0,\infty)\times\mathbb{R}^3}\left(\frac 32\theta-\rho\right)\partial_t\chi dtdx
		\hspace{-30mm}&
		\\
		-\int_{[0,\infty)\times\mathbb{R}^3}
		\left( \frac52 u\theta - \int_{\mathbb{R}^3\times\mathbb{R}^3\times\mathbb{S}^2} q \tilde \psi MM_* dvdv_*d\sigma\right)
		\cdot\nabla_x\chi dtdx
		\hspace{-30mm}&
		\\
		& =\int_{[0,\infty)\times\mathbb{R}^3}u \cdot E \chi dtdx,
	\end{aligned}
\end{equation}
which is precisely the weak formulation of the system
\begin{equation}\label{limiting system 1}
	\begin{aligned}
		\d_t\left( u + A \right)
		+ \nabla_x \cdot \left(u \otimes u - \int_{\mathbb{R}^3\times\mathbb{R}^3\times\mathbb{S}^2} q \tilde \phi MM_* dvdv_*d\sigma \right)
		& = - \nabla_x p
		+\rho E
		+u \wedge B,
		\\
		\d_t \left(\frac 32 \theta-\rho\right) + \nabla_x \cdot \left( \frac52 u \theta - \int_{\mathbb{R}^3\times\mathbb{R}^3\times\mathbb{S}^2} q \tilde \psi MM_* dvdv_*d\sigma\right)
		& = u \cdot E,
	\end{aligned}
\end{equation}
with initial data
\begin{equation*}
	\left(u+A\right)(t=0)=u_0^\mathrm{in}+A_0^\mathrm{in}
	\qquad\text{and}\qquad
	\left(\frac 32\theta-\rho\right)(t=0)=\frac 32\theta_0^\mathrm{in}-\rho_0^\mathrm{in}.
\end{equation*}

By Proposition \ref{weak-comp}, we can further identify the diffusion terms involving the limiting collision integrand $q$. Indeed, utilizing identity \eqref{q phi psi} with formulas \eqref{delta identities}, we obtain
\begin{equation*}
	\begin{aligned}
		\int_{\mathbb{R}^3\times\mathbb{R}^3\times\mathbb{S}^2} q \tilde \phi MM_* dv dv_*d\sigma
		& =\int_{\mathbb{R}^3}\phi:\nabla_x u \tilde\phi Mdv
		=
		\mu\left(\nabla_x u + \nabla_x^t u - \frac 23 \Div u \operatorname{Id}\right),
		\\
		\int_{\mathbb{R}^3\times\mathbb{R}^3\times\mathbb{S}^2} q \tilde \psi MM_* dv dv_*d\sigma
		& =\int_{\mathbb{R}^3}\psi\cdot\nabla_x \theta \tilde\psi Mdv
		=\frac 52\kappa\nabla_x\theta,
	\end{aligned}
\end{equation*}
where the constants $\mu,\kappa>0$ are defined in \eqref{mu kappa} and $\phi$, $\psi$ are the kinetic fluxes defined by \eqref{phi-psi-def}. Incorporating the above relations into \eqref{limiting system 1} and recalling that $u$ is a solenoidal vector field, we finally find the evolution system
\begin{equation*}
	\begin{aligned}
		\d_t\left( u + A \right)
		+ \nabla_x \cdot \left(u \otimes u\right) - \mu\Delta_x u
		& = - \nabla_x p
		+\rho E
		+u \wedge B,
		\\
		\d_t \left(\frac 32 \theta-\rho\right) + \nabla_x \cdot \left( \frac52 u \theta\right)
		-\frac 52\kappa\Delta_x\theta
		& = u \cdot E.
	\end{aligned}
\end{equation*}

Then, defining the adjusted electric field by
\begin{equation*}
	\tilde E= P\tilde E + P^\perp \tilde E = -\partial_t A + \nabla_x\left(\rho+\theta\right),
\end{equation*}
the above evolutions system, when combined with the constraint equations \eqref{w-constraint1} and \eqref{w-constraint2}, can be recast as
\begin{equation*}
	\begin{cases}
		\begin{aligned}
			\d_t u
			+
			u\cdot\nabla_x u - \mu\Delta_x u
			& = -\nabla_x p + \tilde E
			+ \rho \nabla_x\theta + u \wedge B , \hspace{-20mm}&& \\
			&& \Div u & = 0,\\
			\d_t \left(\frac32\theta-\rho\right)
			+
			u\cdot\nabla_x\left(\frac32\theta-\rho\right)
			- \frac 52 \kappa \Delta_x\theta
			& = 0,
			& \Delta_x(\rho+\theta) & =\rho, \\
			\ROT B & = u, & \Div \tilde E & = \rho , \\
			\partial_t B + \rot \tilde E  & = 0, & \Div B & = 0,
		\end{aligned}
	\end{cases}
\end{equation*}
which is precisely the incompressible quasi-static Navier-Stokes-Fourier-Maxwell-Poisson system \eqref{NSFMP 2}.

\subsection{Temporal continuity, initial data and conclusion of proof}

There only remains to establish the weak temporal continuity of the observables~:
\begin{equation}\label{temporal}
	\left(\rho, u, \theta, B\right) \in C\left([0,\infty) ; \textit{w-}L^2\left(\mathbb{R}^3,dx\right)\right),
\end{equation}
and to identify their respective initial data.

For the moment, we only know from the weak formulation \eqref{weak formulation limiting system} that, for any solenoidal $\varphi(x)\in C_c^\infty\left(\mathbb{R}^3;\mathbb{R}^3\right)$ and $\chi(x)\in C_c^\infty\left(\mathbb{R}^3;\mathbb{R}\right)$,
\begin{equation*}
	\begin{aligned}
		\int_{\mathbb{R}^3}\left(u+A\right)(t,x)\varphi(x) dx
		& \in C\left([0,\infty);\mathbb{R}\right),
		\\
		\int_{\mathbb{R}^3}\left(\frac 32\theta-\rho\right)(t,x)\chi(x)dx
		& \in C\left([0,\infty);\mathbb{R}\right),
	\end{aligned}
\end{equation*}
and
\begin{equation*}
	\begin{aligned}
		\left(u+A\right)(0,x)
		& =
		\left(u_0^\mathrm{in}+A_0^\mathrm{in}\right)(x),
		\\
		\left(\frac 32\theta-\rho\right)(0,x)
		& =
		\left(\frac 32\theta_0^\mathrm{in}-\rho_0^\mathrm{in}\right)(x).
	\end{aligned}
\end{equation*}
Notice, replacing the test function $\varphi$ by $\rot\varphi$, that one also has
\begin{equation*}
	\int_{\mathbb{R}^3}\left(\rot u+B\right)(t,x)\varphi(x) dx
	\in C\left([0,\infty);\mathbb{R}\right),
\end{equation*}
for any $\varphi(x)\in C_c^\infty\left(\mathbb{R}^3;\mathbb{R}^3\right)$, and
\begin{equation*}
	\left(\rot u+B\right)(0,x)
	=
	\left(\rot u_0^\mathrm{in}+B_0^\mathrm{in}\right)(x).
\end{equation*}
In particular, a straightforward density argument yields that
\begin{equation*}
	\rot u+B, \frac 32\theta-\rho \in C\left([0,\infty) ; \textit{w-}L^2\left(\mathbb{R}^3,dx\right)\right).
\end{equation*}

Finally, using the relations $\rot B = u$ and $\Delta_x\left(\rho+\theta\right)=\rho$, it is easy to express each observable $\rho$, $u$, $\theta$ and $B$ in terms of $\rot u+B$ and $\frac 32\theta-\rho$, only~:
\begin{equation*}
	\begin{aligned}
		\rho & =\frac{3\Delta_x\left(\rho+\theta\right)-5\Delta_x\rho}{3-5\Delta_x}
		=\frac{2\Delta_x}{3-5\Delta_x}\left(\frac 32\theta-\rho\right),
		\\
		u & = \frac{\rot B - \Delta_x u}{1-\Delta_x}= \frac{\rot}{1-\Delta_x}\left(\rot u+B\right),
		\\
		\theta & =\frac{1-\Delta_x}{\Delta_x}\rho
		=\frac{2-2\Delta_x}{3-5\Delta_x}\left(\frac 32\theta-\rho\right),
		\\
		B & = \frac{\rot}{-\Delta_x}u = \frac{1}{1-\Delta_x}\left(\rot u+B\right).
	\end{aligned}
\end{equation*}
It follows that \eqref{temporal} holds true and that the initial data is provided by
\begin{equation*}
	\begin{aligned}
		\rho(t=0) & = \frac{\Delta_x}{3-5\Delta_x}\left(3\theta^\mathrm{in}_0-2\rho^\mathrm{in}_0\right),
		&
		u(t=0) & = \frac{\ROT}{1-\Delta_x}\left(\ROT u^\mathrm{in}_0+B^\mathrm{in}_0\right),
		\\
		\theta(t=0) & = \frac{1-\Delta_x}{3-5\Delta_x}\left(3\theta^\mathrm{in}_0-2\rho^\mathrm{in}_0\right),
		&
		B(t=0) & = \frac{1}{1-\Delta_x}\left(\ROT u^\mathrm{in}_0+B^\mathrm{in}_0\right),
	\end{aligned}
\end{equation*}
which, at last, concludes the proof of Theorem \ref{NS-WEAKCV}.\qed

%% file: rel-entropy0.tex
\chapter{The renormalized relative entropy method}\label{entropy method}

%%%%%%%%%%%%%%%%%%%%%%%%%%%%%%%%%%%%%%%%%%%%%%%%%

We are now going to investigate the more singular asymptotics leading to the two-fluid incompressible Navier-Stokes-Fourier-Maxwell systems with (solenoidal) Ohm's law. As explained in Section \ref{stability existence 2}, the limiting models obtained in these regimes are not weakly stable and, thus, are not known to have global solutions (except under suitable regularity and smallness assumptions on the initial data).

However, from the physical point of view, these asymptotic regimes are important insofar as they justify Ohm's laws, which are fundamental in plasma physics.

From the mathematical point of view, the Navier-Stokes-Fourier-Maxwell systems obtained in the limit share many features with the three-dimensional incompressible Euler equations. Proving some convergence results requires then methods which are different from the weak compactness techniques used in the proof of Theorem \ref{NS-WEAKCV} in Chapter \ref{grad} and which are typically based on weak-strong stability principles and dissipative solutions (see Section \ref{laure-diogo}). The main novelty here, compared to the convergence results from the Boltzmann equation to the incompressible Euler equations (see \cite[Chapter 5]{SR}), is to use renormalization techniques together with the relative entropy method.

%%%%%%%%%%%%%%%%%%%%%%%%%%%%%%%%%%%%%%%%%%%%%%%%%

\section{The relative entropy method~: old and new}

The principle of the relative entropy method is to compare the distribution with its formal asymptotics in some appropriate metrics~:
\begin{itemize}
	\item
	The idea of using the relative entropy $H\left(f_\eps^\pm|M\right)$ to build such metrics goes back to Yau \cite{yau} in the framework of the asymptotic study of Ginzburg-Landau's equation, then to Golse \cite{golse0} for the hydrodynamic limits of the Boltzmann equation. The important points of this method are the fact that the scaled relative entropy is a Lyapunov functional for the Boltzmann equation and that it controls the size of the fluctuations.
	
	\item
	An approximate solution is obtained by formal expansions (the so-called Hilbert expansions~; see \cite{hilbert}), which consist in seeking a formal solution to the scaled system in the form
	\begin{equation*}
		f^\pm _\mathrm{app}= M\left( 1 + \eps g_0^\pm + \eps^2 g_1^\pm + \ldots \right).
	\end{equation*}
	Note that the successive approximations $g_n^\pm$ should depend a priori both on macroscopic variables $t$, $x$ and on fast variables $\frac{t}{\eps}$, $\frac{t}{\eps^2}$, \ldots, $\frac{x}{\eps}$, $\frac{x}{\eps^2}$, \ldots .
	
	For well-prepared initial data, that is for data satisfying some profile condition (thermodynamic equilibrium) as well as macroscopic linear constraints (incompressibility and Boussinesq relations, for instance), there is neither kinetic initial layer nor fast oscillating waves, so that $g_0^\pm$ reduces actually to the solution of the limiting system. Note that, in the cases considered here, the nonlinear constraints (Ohm's laws) have a different status~: solutions of the limiting models are well-defined even though these constraints are not defined initially. This is similar to the existence theory for parabolic equations with initial data which are not in the domain of the diffusion operator.
	
	\item
	The core of the proof consists then in getting some stability inequality for the scaled modulated entropy
	\begin{equation*}
		\sum_\pm \frac1{\eps^2} H\left(f_\eps^\pm|f_\mathrm{app}^\pm\right)
		=\sum_\pm \frac1{\eps^2}
		\int_{\mathbb{R}^3\times\mathbb{R}^3} \left( f_\eps ^\pm \log {f^\pm_\eps \over f^\pm_\mathrm{app}} - f^\pm_\eps + f^\pm_\mathrm{app} \right) dxdv,
	\end{equation*}
	which measures in some sense the distance between the fluctuations $g_\eps^\pm$ and their expected limits $g_0^\pm$. The convergence relies then on some technical computations and Gr\"onwall's lemma.
\end{itemize}

The stability inequality we expect to obtain should be reminiscent of the inequality defining the corresponding dissipative solutions of the limiting systems (see Section \ref{laure-diogo}). Thus, it should be based solely on the decay of the entropy and on local conservation laws. In particular, there is no need for a priori strong compactness~: nonlinear terms should be controlled by a loop estimate using Gr\"onwall's lemma.

Unfortunately, this simple strategy fails, in general~:
\begin{itemize}
	\item even for weak solutions in the sense of distributions (not renormalized), provided they exist, we have no control on large velocities~;
	\item for renormalized solutions in the sense of DiPerna and Lions, provided they exist, local conservation laws are not known to hold.
	% \item for renormalized solutions with a defect measure in the sense of Alexandre and Villani, even after renormalization, the kinetic equation is relaxed into an inequality.
\end{itemize}

The main novelty here is to use {\bf renormalization techniques combined with the relative entropy method}. More precisely, we will not use the usual modulated entropy inequality for renormalized solutions to the kinetic equations. Rather, we will modulate a renormalized version of the entropy inequality, which requires much less a priori information on the solutions.

\section{Proof of Theorem \ref{CV-OMHD} on weak interactions}\label{proof of theorem weak}

Several steps of this demonstration have been performed in the preceding chapters. We therefore begin our proof by appropriately gathering previous results together.

\subsection{Weak convergence of fluctuations, collision integrands and electromagnetic fields}

Thus, we are considering here a family of renormalized solutions $(f_\eps^\pm,E_\eps,B_\eps)$ to the scaled two species Vlasov-Maxwell-Boltzmann system \eqref{VMB2}, in the regime of weak interspecies interactions, i.e.\ $\delta=o(1)$ and $\frac{\delta}{\eps}$ is unbounded, satisfying the scaled entropy inequality \eqref{entropy2}.

According to Lemmas \ref{L1-lem} and \ref{L2-lem}, the corresponding families of fluctuations $g_\eps^\pm$ and renormalized fluctuations $\hat g_\eps^\pm$ are weakly compact in $L^1_\mathrm{loc}\left(dtdx;L^1\left(\left(1+|v|^2\right)Mdv\right)\right)$ and $L^2_\mathrm{loc}\left(dt;L^2\left(Mdxdv\right)\right)$, respectively, while, in view of Lemma \ref{L2-qlem}, the corresponding collision integrands $\hat q_\eps^\pm$ and $\hat q_\eps^{\pm,\mp}$ are weakly compact in $L^2\left(MM_*dtdxdvdv_*d\sigma\right)$. Thus, using Lemma \ref{L1-lem} again and the decomposition \eqref{fluct-decomposition}, we know that there exist $g^\pm\in L^\infty\left(dt;L^2\left(Mdxdv\right)\right)$, $(E,B)\in L^\infty\left(dt;L^2(dx)\right)$ and $q^\pm,q^{\pm,\mp}\in L^2\left(MM_*dtdxdvdv_*d\sigma\right)$, such that, up to extraction of subsequences,
\begin{equation}\label{weak conv fluctuations}
	\begin{aligned}
		g_\eps^\pm & \rightharpoonup g^\pm & & \text{in }L^1_\mathrm{loc}\left(dtdx;L^1\left(\left(1+|v|^2\right)Mdv\right)\right),
		\\
		\hat g_\eps^\pm & \stackrel{*}{\rightharpoonup} g^\pm & & \text{in }L^\infty\left(dt;L^2\left(Mdxdv\right)\right),
		\\
		\left(E_\eps,B_\eps\right) & \stackrel{*}{\rightharpoonup} \left(E,B\right) & & \text{in }L^\infty\left(dt;L^2\left(dx\right)\right),
		\\
		\hat q_\eps^\pm & \rightharpoonup q^\pm & & \text{in }L^2\left(MM_*dtdxdvdv_*d\sigma\right),
		\\
		\hat q_\eps^{\pm,\mp} & \rightharpoonup q^{\pm,\mp} & & \text{in }L^2\left(MM_*dtdxdvdv_*d\sigma\right),
	\end{aligned}
\end{equation}
as $\eps\to 0$. Therefore, one also has the weak convergence of the densities $\rho_\eps^\pm$, bulk velocities $u_\eps^\pm$ and temperatures $\theta_\eps^\pm$ corresponding to $g_\eps^\pm$~:
\begin{equation*}
	\rho_\eps^\pm \rightharpoonup \rho^\pm,
	\quad u_\eps^\pm \rightharpoonup u^\pm
	\quad \text{and}\quad
	\theta_\eps^\pm \rightharpoonup \theta^\pm
	\quad\text{in }L^1_\mathrm{loc}(dtdx) \text{ as }\eps \to 0,
\end{equation*}
where $\rho^\pm, u^\pm, \theta^\pm \in L^\infty\left(dt;L^2(dx)\right)$ are, respectively, the densities, bulk velocities and temperatures corresponding to $g^\pm$. In fact, Lemma \ref{relaxation-control} implies that
\begin{equation*}
	g^\pm=\Pi g^\pm = \rho^\pm + u^\pm\cdot v + \theta^\pm\left(\frac{|v|^2}{2}-\frac 32\right).
\end{equation*}

Next, we further introduce the scaled fluctuations
\begin{equation*}
	h_\eps = \frac{\delta}{\eps}\left[\left(g_\eps^+-g_\eps^-\right) - n_\eps\right],
\end{equation*}
where $n_\eps=\rho_\eps^+-\rho_\eps^-$ is the charge density, and the electrodynamic variables
\begin{equation*}
	j_\eps=\frac \delta\eps\left(u_\eps^+-u_\eps^-\right),\qquad w_\eps=\frac\delta\eps\left(\theta_\eps^+-\theta_\eps^-\right),
\end{equation*}
which are precisely the bulk velocity and temperature associated with the scaled fluctuations $h_\eps$. In view of Lemma \ref{bound hjw}, the electric current $j_\eps$ and the internal electric energy $w_\eps$ are uniformly bounded in $L^1_\mathrm{loc}\left(dtdx\right)$, which necessarily implies, letting $\frac\eps\delta\rightarrow 0$, that $u^+=u^-$ and $\theta^+=\theta^-$. Furthermore, Proposition \ref{weak-comp2} asserts that $\rho^+=\rho^-$, as well. Therefore, we appropriately rename the limiting macroscopic variables
\begin{equation*}
	\rho=\rho^+=\rho^-, \qquad u=u^+=u^-, \qquad \theta=\theta^+=\theta^-,
\end{equation*}
and the limiting fluctuation
\begin{equation*}
	g=g^+=g^-=\rho + u\cdot v + \theta\left(\frac{|v|^2}{2}-\frac 32\right).
\end{equation*}

Now, according to Lemma \ref{weak compactness h}, it holds that $h_\eps$ is weakly compact in $L^1_\mathrm{loc}\left(dtdx;L^1\left(\left(1+|v|^2\right)Mdv\right)\right)$ and that $j_\eps$ and $w_\eps$ are weakly compact in $L^1_\mathrm{loc}(dtdx)$, so that, up to extraction of subsequences, there are $h\in L^1_\mathrm{loc}\left(dtdx;L^1\left(\left(1+|v|^2\right)Mdv\right)\right)$ and $j,w\in L^1_\mathrm{loc}(dtdx)$ such that
\begin{equation*}
	\begin{aligned}
		h_\eps & \rightharpoonup h & & \text{in }L^1_\mathrm{loc}\left(dtdx;L^1\left(\left(1+|v|^2\right)Mdv\right)\right),
		\\
		j_\eps & \rightharpoonup j & & \text{in }L^1_\mathrm{loc}\left(dtdx\right),
		\\
		w_\eps & \rightharpoonup w & & \text{in }L^1_\mathrm{loc}\left(dtdx\right),
	\end{aligned}
\end{equation*}
as $\eps\to 0$. Moreover, by Proposition \ref{weak-comp3}, one has the infinitesimal Maxwellian form
\begin{equation*}
	h=j\cdot v + w\left(\frac{|v|^2}{2}-\frac 32\right).
\end{equation*}

\subsection{Constraint equations, Maxwell's system and energy inequality}

In view of Proposition \ref{weak-comp2}, we already know that the limiting thermodynamic fields $\rho$, $u$ and $\theta$ satisfy the incompressibility and Boussinesq relations
\begin{equation}\label{constraints}
	\Div u = 0, \qquad \rho+\theta = 0.
\end{equation}

Moreover, the discussion in Section \ref{limit maxwell} shows that the limiting electromagnetic field satisfies the following form of Maxwell's equations~:
\begin{equation*}
	\begin{cases}
		\begin{aligned}
			\d_t E - \ROT B &= -j,
			\\
			\d_t B + \ROT E& = 0,
			\\
			\DIV E &=0,
			\\
			\DIV B &=0.
		\end{aligned}
	\end{cases}
\end{equation*}
Note that, taking the divergence of the Amp\`ere equation above, necessarily $\Div j = 0$.

Finally, Proposition \ref{solenoidalOhm} further establishes that the electrodynamic variables $j$ and $w$ satisfy the solenoidal Ohm's law and the internal electric energy equilibrium relation
\begin{equation*}
	j = \sigma\left(-\nabla_x \bar p + E + u\wedge B \right),
	\qquad
	w = 0,
\end{equation*}
where the electric conductivity $\sigma>0$ is defined by \eqref{sigma} and the pressure gradient $-\nabla_x \bar p$ is the Lagrange multiplier associated with the solenoidal constraint $\Div j =0$.

As for the energy bound, Proposition \ref{energy inequality weak interactions} states that, for almost every $t\geq 0$,
\begin{equation*}
	\begin{aligned}
		& \frac 12\left( 2 \left\|u\right\|_{L^2_x}^2
		+ 5 \left\|\theta\right\|_{L^2_x}^2 + \left\|E\right\|_{L^2_x}^2
		+ \left\|B\right\|_{L^2_x}^2 \right)(t)
		\\
		& \hspace{10mm} +
		\int_0^t \left(2\mu
		\left\|\nabla_x u\right\|_{L^2_x}^2
		+ 5 \kappa
		\left\|\nabla_x\theta\right\|_{L^2_x}^2
		+\frac 1{\sigma}
		\left\|j\right\|_{L^2_x}^2
		\right)(s) ds \leq C^\mathrm{in},
	\end{aligned}
\end{equation*}
where the viscosity $\mu>0$, thermal conductivity $\kappa>0$ and electric conductivity $\sigma>0$ are respectively defined by \eqref{mu kappa 2} and \eqref{sigma}. In particular, it holds that
\begin{equation*}
		\begin{aligned}
			\left(u, \theta, E, B\right) & \in L^\infty \left( [0,\infty), dt ; L^2\left(\mathbb{R}^3, dx\right)\right), \\
			\left(u,\theta\right) & \in L^2\left([0,\infty), dt ; \dot H^1\left(\mathbb{R}^3, dx\right)\right), \\
			j & \in L^2\left([0,\infty)\times\mathbb{R}^3, dtdx\right).
		\end{aligned}
\end{equation*}
This energy bound can be improved to the actual energy inequality
\begin{equation*}
	\begin{aligned}
		& \frac 12\left( 2 \left\|u\right\|_{L^2_x}^2
		+ 5 \left\|\theta\right\|_{L^2_x}^2 + \left\|E\right\|_{L^2_x}^2
		+ \left\|B\right\|_{L^2_x}^2 \right)(t)
		\\
		& \hspace{10mm} +
		\int_0^t \left(2\mu
		\left\|\nabla_x u\right\|_{L^2_x}^2
		+ 5 \kappa
		\left\|\nabla_x\theta\right\|_{L^2_x}^2
		+\frac 1{\sigma}
		\left\|j\right\|_{L^2_x}^2
		\right)(s) ds
		\\
		& \hspace{10mm} \leq \frac 12\left( 2 \left\|u^\mathrm{in}\right\|_{L^2_x}^2
		+ 5 \left\|\theta^\mathrm{in}\right\|_{L^2_x}^2 + \left\|E^\mathrm{in}\right\|_{L^2_x}^2
		+ \left\|B^\mathrm{in}\right\|_{L^2_x}^2 \right),
	\end{aligned}
\end{equation*}
using the well-preparedness of the initial data \eqref{well-prepared init data}.

\subsection{The renormalized modulated entropy inequality}\label{stability weak}
% \subsection{Evolution equations and stability inequality}

We move on now to the rigorous derivation of a stability inequality encoding the asymptotic macroscopic evolution equations for $u$ and $\theta$ in the spirit of the weak-strong stability inequalities used in Section \ref{laure-diogo} to define dissipative solutions for some Navier-Stokes-Maxwell systems. Recall that, as explained therein, such systems are in general not known to display weak stability so that their weak solutions in the energy space are not known to exist.

To this end, as in Section \ref{conservation defects 2 species}, we define the renormalized fluctuations $g_\eps^\pm \gamma_\eps^\pm\chi\left(\frac{|v|^2}{K_\eps}\right)$, with $K_\eps =K|\log \eps|$, for some large $K>0$, and $\chi\in C_c^\infty\left([0,\infty)\right)$ a smooth compactly supported function such that $\mathds{1}_{[0,1]}\leq \chi \leq \mathds{1}_{[0,2]}$, and where $\gamma_\eps^\pm=\gamma\left(G_\eps^\pm\right)$ for some renormalization $\gamma\in C^1\left([0,\infty);\mathbb{R}\right)$ satisfying \eqref{gamma-assumption}.

Since, up to further extraction of subsequences, $\gamma_\eps^\pm \chi\left( {|v|^2\over K_\eps}\right)$ converges almost everywhere towards $1$, $g_\eps^\pm$ is weakly compact in $L^1_\mathrm{loc}\left(dtdx;L^1\left(\left(1+|v|^2\right)Mdv\right)\right)$ and $g_\eps^\pm\gamma_\eps^\pm$ is uniformly bounded in $L^\infty\left(dt;L^2\left(Mdxdv\right)\right)$, we deduce, by the Product Limit Theorem, that
\begin{equation*}
	g_\eps^\pm \gamma_\eps^\pm \chi\left( {|v|^2\over K_\eps}\right)
	\stackrel{*}{\rightharpoonup} g \quad \text{in }L^\infty\left(dt;L^2\left(Mdxdv\right)\right).
\end{equation*}
Therefore, one has the weak convergence of the densities $\tilde\rho_\eps^\pm$, bulk velocities $\tilde u_\eps^\pm$ and temperatures $\tilde\theta_\eps^\pm$ corresponding to $g_\eps^\pm\gamma_\eps^\pm \chi\left( {|v|^2\over K_\eps}\right)$~:
\begin{equation*}
	\tilde\rho_\eps^\pm \stackrel{*}{\rightharpoonup} \rho,
	\quad \tilde u_\eps^\pm \stackrel{*}{\rightharpoonup} u
	\quad \text{and}\quad
	\tilde\theta_\eps^\pm \stackrel{*}{\rightharpoonup} \theta
	\quad\text{in }L^\infty\left(dt;L^2(dx)\right) \text{ as }\eps \to 0.
\end{equation*}
In particular, the hydrodynamic variables $\tilde\rho_\eps=\frac{\tilde\rho_\eps^++\tilde\rho_\eps^-}{2}$, $\tilde u_\eps=\frac{\tilde u_\eps^++\tilde u_\eps^-}{2}$ and $\tilde\theta_\eps=\frac{\tilde\theta_\eps^++\tilde\theta_\eps^-}{2}$ also obviously verify
\begin{equation}\label{weak limit observables}
	\tilde\rho_\eps \stackrel{*}{\rightharpoonup} \rho,
	\quad \tilde u_\eps \stackrel{*}{\rightharpoonup} u
	\quad \text{and}\quad
	\tilde\theta_\eps \stackrel{*}{\rightharpoonup} \theta
	\quad\text{in }L^\infty\left(dt;L^2(dx)\right) \text{ as }\eps \to 0.
\end{equation}
It follows that, since $u$ is solenoidal,
\begin{equation}\label{limit solenoidal}
	P^\perp \tilde u_\eps \stackrel{*}{\rightharpoonup} 0
	\quad\text{in }L^\infty\left(dt;L^2(dx)\right) \text{ as }\eps \to 0,
\end{equation}
and, in view of the limiting Boussinesq relation,
\begin{equation}\label{limit boussinesq}
	\tilde \rho_\eps+\tilde\theta_\eps \stackrel{*}{\rightharpoonup} 0
	\quad\text{in }L^\infty\left(dt;L^2(dx)\right) \text{ as }\eps \to 0.
\end{equation}

We establish now the convergence of the electric current $\tilde j_\eps =\frac \delta\eps\left(\tilde u_\eps^+-\tilde u_\eps^-\right)$. Since $\eps \hat g_\eps^\pm\geq 2\left(\sqrt 2 -1 \right)$ on the support of $1-\gamma_\eps^\pm$, we easily estimate, using the uniform bound from Lemma \ref{v2-int}, that
\begin{equation*}
	\begin{aligned}
		\left|\frac\delta\eps g_\eps^\pm\left(1-\gamma_\eps^\pm\right)\right|
		& =
		\left|\frac\delta\eps\left( \hat g_\eps^\pm + \frac\eps 4 \hat g_\eps^{\pm 2}\right)\left(1-\gamma_\eps^\pm\right)\right|
		\\
		& \leq
		\left|\frac\delta\eps \hat g_\eps^\pm\left(1-\gamma_\eps^\pm\right)\right|
		+\frac\delta 4 \hat g_\eps^{\pm 2}
		\\
		& \leq C\delta \hat g_\eps^{\pm 2} = O(\delta)_{L^1_\mathrm{loc}\left(dtdx;L^1\left(\left(1+|v|^2\right)Mdv\right)\right)},
	\end{aligned}
\end{equation*}
whereas, using the Gaussian decay \eqref{gaussian-decay 0}, we also obtain, provided $K>4$,
\begin{equation*}
	\begin{aligned}
		\left|\frac\delta\eps g_\eps^\pm\left(1-\chi\left(\frac{|v|^2}{K_\eps}\right)\right)\right|
		& =
		\left|\frac\delta\eps\left( \hat g_\eps^\pm + \frac\eps 4 \hat g_\eps^{\pm 2}\right)\left(1-\chi\left(\frac{|v|^2}{K_\eps}\right)\right)\right|
		\\
		& \leq
		\left|\frac\delta\eps \hat g_\eps^\pm\left(1-\chi\left(\frac{|v|^2}{K_\eps}\right)\right)\right|
		+\frac\delta 4 \hat g_\eps^{\pm 2}
		\\
		& \leq C\frac{\delta}{\eps^2}\left(1-\chi\left(\frac{|v|^2}{K_\eps}\right)\right)^2 + C\delta \hat g_\eps^{\pm 2}
		\\
		& = O(\delta)_{L^1_\mathrm{loc}\left(dtdx;L^1\left(\left(1+|v|^2\right)Mdv\right)\right)}.
	\end{aligned}
\end{equation*}
Thus, we infer that
\begin{equation}\label{tilde u approx u}
	\frac\delta\eps\left(\tilde u_\eps^\pm-u_\eps^\pm\right)\rightarrow 0
	\qquad\text{in }L^1_\mathrm{loc}\left(dtdx\right) \text{ as }\eps \to 0,
\end{equation}
whence
\begin{equation*}
	\tilde j_\eps\rightharpoonup j
	\qquad\text{in }L^1_\mathrm{loc}\left(dtdx\right) \text{ as }\eps \to 0.
\end{equation*}

\bigskip

Now, the $L^2\left(Mdxdv\right)$ norm of $g_\eps^\pm \gamma_\eps^\pm\chi\left(\frac{|v|^2}{K_\eps}\right)$ is not a Lyapunov functional but it is nevertheless controlled by the relative entropy
\begin{equation}\label{renormalized energy entropy control}
	\frac 12 \left\|g_\eps^\pm \gamma_\eps^\pm\chi\left(\frac{|v|^2}{K_\eps}\right)\right\|_{L^2\left(Mdxdv\right)}^2
	\leq \frac{C}{\eps^2}H\left(f_\eps^\pm\right),
\end{equation}
for some $C>1$, and therefore by the initial data \eqref{init-fluctuation 2}. One may therefore try, in a preliminary attempt to show an asymptotic stability inequality, to modulate the approximate energy associated with $g_\eps^\pm \gamma_\eps^\pm\chi\left(\frac{|v|^2}{K_\eps}\right)$, i.e.\ its $L^2\left(Mdxdv\right)$ norm, by introducing a test function $\bar g$ in infinitesimal Maxwellian form~:
\begin{equation*}
	\bar g= \bar u\cdot v + \bar \theta \left(\frac{|v|^2}{2}-\frac 52\right),
\end{equation*}
where
\begin{equation*}
	\bar u(t,x),\bar \theta(t,x)\in C^\infty_c\left([0,\infty)\times\mathbb{R}^3\right)
	\quad\text{with } \Div\bar u = 0,
\end{equation*}
and then establishing a stability inequality for the modulated energies
\begin{equation}\label{modulated energy}
	\frac 12 \left\|g_\eps^\pm \gamma_\eps^\pm\chi\left(\frac{|v|^2}{K_\eps}\right)-\bar g\right\|_{L^2\left(Mdxdv\right)}^2.
\end{equation}
Notice that it holds, utilizing the elementary identity $a^2+\frac 32 b^2=\frac 35\left(a+b\right)^2+\frac 52\left(\frac{3b-2a}{5}\right)^2$, for any $a,b\in\mathbb{R}$,
\begin{equation*}
	\begin{aligned}
		& \left\|g_\eps^\pm \gamma_\eps^\pm\chi\left(\frac{|v|^2}{K_\eps}\right)-\bar g\right\|_{L^2\left(Mdxdv\right)}^2
		\\ & \geq
		\left\|\Pi\left(g_\eps^\pm \gamma_\eps^\pm\chi\left(\frac{|v|^2}{K_\eps}\right)\right)-\bar g\right\|_{L^2\left(Mdxdv\right)}^2
		\\ & =
		\left(\left\|\tilde\rho_\eps^\pm+\bar\theta\right\|_{L^2\left(dx\right)}^2
		+\left\|\tilde u_\eps^\pm-\bar u\right\|_{L^2\left(dx\right)}^2
		+\frac 32 \left\|\tilde\theta_\eps^\pm-\bar \theta\right\|_{L^2\left(dx\right)}^2\right)
		\\ & =
		\left(\frac 35\left\|\tilde\rho_\eps^\pm+\tilde\theta_\eps^\pm\right\|_{L^2\left(dx\right)}^2
		+\left\|\tilde u_\eps^\pm-\bar u\right\|_{L^2\left(dx\right)}^2
		+\frac 52 \left\|\frac{3\tilde\theta_\eps^\pm-2\tilde\rho_\eps^\pm}{5}-\bar \theta\right\|_{L^2\left(dx\right)}^2\right).
	\end{aligned}
\end{equation*}

It turns out that this approach is not quite suitable for our purpose because, even though, for any $0\leq t_1<t_2$ (see the proof of Lemma \ref{L1-lem}),
\begin{equation}\label{limit entropy}
		\int_{t_1}^{t_2}
		\frac 1 {2}
		\left\|g\right\|_{L^2\left(Mdxdv\right)}^2 dt
		\leq
		\liminf_{\eps\rightarrow 0}\int_{t_1}^{t_2}\frac 1{\eps^2}H\left(f_\eps^\pm\right)dt,
\end{equation}
it is not possible to set $C=1$ in \eqref{renormalized energy entropy control}. Indeed, the first term in the polynomial expansion of the function $h(z)=(1+z)\log(1+z)-z$ defining the entropy is $\frac 12 z^2$, but the second term is $-\frac 16 z^3$ and may be negative.

Some entropy (or energy) is therefore lost by considering the modulated energies \eqref{modulated energy}. These considerations lead us to introduce a more precise modulated functional in replacement of \eqref{modulated energy} capturing more information on the fluctuations. To be precise, instead of \eqref{modulated energy}, we consider now the {\bf renormalized modulated entropies}
\begin{equation}\label{modulated entropy}
	\frac 1{\eps^2}H\left(f_\eps^\pm\right)
	-
	\int_{\mathbb{R}^3\times\mathbb{R}^3}
	g_\eps^\pm \gamma_\eps^\pm\chi\left(\frac{|v|^2}{K_\eps}\right)\bar g
	Mdxdv
	+
	\frac 1 {2}
	\left\|\bar g\right\|_{L^2\left(Mdxdv\right)}^2.
\end{equation}
Note that the above functional may be negative for fixed $\eps>0$. However, in view of \eqref{limit entropy}, it recovers asymptotically a non-negative quantity, which is precisely the asymptotic modulated energy~:
\begin{equation}\label{modulated entropy liminf}
	\begin{aligned}
		& \int_{t_1}^{t_2}
		\frac 12 \left(\left\|u-\bar u\right\|_{L^2(dx)}^2+\frac 52\left\|\theta-\bar \theta\right\|_{L^2(dx)}^2\right)
		dt
		\\
		& = \int_{t_1}^{t_2}
		\frac 1 {2}
		\left\|g - \bar g\right\|_{L^2\left(Mdxdv\right)}^2 dt
		\\
		& \leq \liminf_{\eps\rightarrow 0}
		\int_{t_1}^{t_2}\left(\frac 1{\eps^2}H\left(f_\eps^\pm\right)
		-
		\int_{\mathbb{R}^3\times\mathbb{R}^3}
		g_\eps^\pm \gamma_\eps^\pm\chi\left(\frac{|v|^2}{K_\eps}\right)\bar g
		Mdxdv
		+
		\frac 1 {2}
		\left\|\bar g\right\|_{L^2\left(Mdxdv\right)}^2\right)dt,
	\end{aligned}
\end{equation}
for all $0\leq t_1<t_2$.

The first term in \eqref{modulated entropy} is precisely the entropy of $f_\eps^\pm$ and will be controlled by the scaled entropy inequality \eqref{entropy2}, whereas the last term in \eqref{modulated entropy} only involves smooth quantities and will therefore be controlled directly. As for the middle term in the modulated entropy \eqref{modulated entropy}, its time derivative will involve the {\bf approximate macroscopic conservation laws} for $ g_\eps^\pm \gamma_\eps ^\pm \chi\left(\frac{|v|^2}{K_\eps}\right) $.

\bigskip

Recall that a major difficulty in the relative entropy methods developed for the hydrodynamic limit of the Boltzmann equation towards the incompressible Euler equations (see \cite{SR2, SR, SR3}) pertains to the handling of large velocities. Here, large velocities are no longer a problem, for we deal now with conservations laws of renormalized fluctuations $ g_\eps^\pm \gamma_\eps ^\pm \chi\left(\frac{|v|^2}{K_\eps}\right) $ whose defects are well controlled. Thus, the present method is more robust than the usual relative entropy method which cannot deal with fluctuations of temperature.

Furthermore, thanks to the flatness assumption \eqref{gamma-assumption} on $\gamma(z)$ near $z=1$, the conservation defects are expected to vanish in the limit (at least formally). In fact, as shown in Section \ref{conservation defects 2 species}, they can be estimated in terms of the modulated entropy and entropy dissipation. The convergence will then be obtained through a loop estimate based on an appropriate use of Gr\"onwall's lemma.

\bigskip

Now, in order to establish the renormalized modulated entropy inequality leading to the convergence stated in Theorem \ref{CV-OMHD}, we introduce further test functions
\begin{equation*}
	\bar E(t,x),\bar B(t,x), \bar j(t,x)\in C^\infty_c\left([0,\infty)\times\mathbb{R}^3\right)
	\quad\text{with }\Div\bar E = \Div\bar B = \Div\bar j = 0,
\end{equation*}
and we define the {\bf renormalized modulated entropy}
\begin{equation*}
	\begin{aligned}
		\delta\mathcal{H}_\eps(t)
		& =
		\frac 1{\eps^2}H\left(f_\eps^+\right) + \frac 1{\eps^2}H\left(f_\eps^-\right)
		\\
		& -
		\int_{\mathbb{R}^3\times\mathbb{R}^3}
		\left(g_\eps^+ \gamma_\eps^++g_\eps^- \gamma_\eps^-\right)\chi\left(\frac{|v|^2}{K_\eps}\right)\bar g
		Mdxdv
		+
		\left\|\bar g\right\|_{L^2\left(Mdxdv\right)}^2
		\\
		& +
		\frac 12\left\|E_\eps - \bar E\right\|_{L^2(dx)}^2
		+
		\frac 12\left\|B_\eps - \bar B\right\|_{L^2(dx)}^2
		+ \frac 1{2}\int_{\mathbb{R}^3}\left(\frac 1{\eps^2}\operatorname{Tr}m_\eps+\operatorname{Tr}a_\eps\right) dx
		\\
		& - \int_{\mathbb{R}^3}
		\left(\left(E_\eps-\bar E\right)\wedge\left(B_\eps-\bar B\right)
		+
		\begin{pmatrix}
			a_{\eps 26}-a_{\eps 35}\\a_{\eps 34}-a_{\eps 16}\\a_{\eps 15}-a_{\eps 24}
		\end{pmatrix}
		\right)\cdot\bar u
		dx,
	\end{aligned}
\end{equation*}
where the matrix measures $m_\eps$ and $a_\eps$ are the defects introduced in Section \ref{macroscopic defects} and controlled by the scaled entropy inequality \eqref{entropy2}.

We also define the {\bf renormalized modulated energy}
\begin{equation*}
	\begin{aligned}
		\delta\mathcal{E}_\eps(t)
		& =
		\frac 12\left\|g_\eps^+ \gamma_\eps^+ \chi\left(\frac{|v|^2}{K_\eps}\right)-\bar g\right\|_{L^2(Mdxdv)}^2
		+
		\frac 12\left\|g_\eps^- \gamma_\eps^- \chi\left(\frac{|v|^2}{K_\eps}\right)-\bar g\right\|_{L^2(Mdxdv)}^2
		\\
		& +
		\frac 12\left\|E_\eps - \bar E\right\|_{L^2(dx)}^2
		+
		\frac 12\left\|B_\eps - \bar B\right\|_{L^2(dx)}^2
		+ \frac 1{2}\int_{\mathbb{R}^3}\left(\frac 1{\eps^2}\operatorname{Tr}m_\eps+\operatorname{Tr}a_\eps\right) dx
		\\
		& - \int_{\mathbb{R}^3}
		\left(\left(E_\eps-\bar E\right)\wedge\left(B_\eps-\bar B\right)
		+
		\begin{pmatrix}
			a_{\eps 26}-a_{\eps 35}\\a_{\eps 34}-a_{\eps 16}\\a_{\eps 15}-a_{\eps 24}
		\end{pmatrix}
		\right)\cdot\bar u
		dx,
	\end{aligned}
\end{equation*}
which is asymptotically equivalent to $\delta\mathcal{H}_\eps(t)$, at least formally. Note that $\delta\mathcal{H}_\eps(t)$ controls more accurately the large values of the fluctuations $g_\eps^\pm$ than $\delta\mathcal{E}_\eps(t)$. Lemma \ref{entropy energy} below shows how the modulated entropy $\delta\mathcal{H}_\eps(t)$ controls the modulated energy $\delta\mathcal{E}_\eps(t)$.

Finally, we introduce the {\bf renormalized modulated entropy dissipation}
\begin{equation*}
	\begin{aligned}
		\delta\mathcal{D}_\eps(t) & =
		\frac 14
		\left\|\hat q^+_\eps - \bar q^+\right\|_{L^2\left(MM_*dxdvdv_*d\sigma\right)}^2
		+
		\frac 14
		\left\|\hat q^-_\eps - \bar q^-\right\|_{L^2\left(MM_*dxdvdv_*d\sigma\right)}^2
		\\
		& + \frac 14
		\left\|\hat q^{+,-}_\eps - \bar q^{+,-}\right\|_{L^2\left(MM_*dxdvdv_*d\sigma\right)}^2
		+
		\frac 14
		\left\|\hat q^{-,+}_\eps - \bar q^{-,+}\right\|_{L^2\left(MM_*dxdvdv_*d\sigma\right)}^2,
	\end{aligned}
\end{equation*}
where
\begin{equation}\label{test integrand}
	\begin{aligned}
		\bar q^\pm
		& =
		\nabla_x \bar u:\left(\tilde \phi + \tilde \phi_*
		-\tilde\phi'-\tilde\phi_*'\right)
		+\nabla_x \bar \theta\cdot\left(\tilde \psi + \tilde \psi_*
		-\tilde\psi'-\tilde\psi_*'\right),
		\\
		\bar q^{\pm,\mp}
		& = \mp\frac 12 \bar j\cdot\left(v-v_*-v'+v_*'\right),
	\end{aligned}
\end{equation}
so that
\begin{equation*}
	\begin{aligned}
		\int_{\mathbb{R}^3\times\mathbb{S}^2} \bar q^\pm M_*dv_*d\sigma
		& = \nabla_x \bar u:\mathcal{L}\tilde\phi + \nabla_x \bar \theta\cdot \mathcal{L}\tilde\psi
		= \nabla_x \bar u:\phi + \nabla_x \bar \theta\cdot \psi,
		\\
		\int_{\mathbb{R}^3\times\mathbb{S}^2} \bar q^{\pm,\mp} M_* dv_*d\sigma
		& =
		\mp\frac 12 \bar j\cdot \mathfrak{L} \left(v\right),
	\end{aligned}
\end{equation*}
with $\phi$, $\psi$, $\tilde\phi$ and $\tilde \psi$ defined by \eqref{phi-psi-def} and \eqref{phi-psi-def inverses}.

Then, assuming from now on that $\left\|\bar u\right\|_{L^\infty(dtdx)}<1$ and using the lower weak sequential semi-continuity of the entropies \eqref{modulated entropy liminf} and of the electromagnetic energy \eqref{convex wedge} together with Lemma \ref{vector defect}, we find that, for all $0\leq t_1<t_2$,
\begin{equation}\label{energy liminf}
	\int_{t_1}^{t_2}
	\delta\mathcal{E}(t) dt
	\leq \liminf_{\eps\rightarrow 0}
	\min\left\{
	\int_{t_1}^{t_2}\delta\mathcal{H}_\eps(t) dt,
	\int_{t_1}^{t_2}\delta\mathcal{E}_\eps(t) dt
	\right\},
\end{equation}
where
\begin{equation*}
	\begin{aligned}
		\delta\mathcal{E}(t)
		& =
		\left\|g - \bar g\right\|_{L^2\left(Mdxdv\right)}^2
		+\frac 12\left\|E-\bar E\right\|_{L^2(dx)}^2
		+\frac 12\left\|B-\bar B\right\|_{L^2(dx)}^2
		\\
		\\
		& - \int_{\mathbb{R}^3}
		\left(\left(E_\eps-\bar E\right)\wedge\left(B_\eps-\bar B\right)
		\right)\cdot\bar u
		dx
		\\
		& =
		\left\|u-\bar u\right\|_{L^2(dx)}^2+\frac 52\left\|\theta-\bar \theta\right\|_{L^2(dx)}^2
		+\frac 12\left\|E-\bar E\right\|_{L^2(dx)}^2
		+\frac 12\left\|B-\bar B\right\|_{L^2(dx)}^2
		\\
		\\
		& - \int_{\mathbb{R}^3}
		\left(\left(E_\eps-\bar E\right)\wedge\left(B_\eps-\bar B\right)
		\right)\cdot\bar u
		dx,
	\end{aligned}
\end{equation*}
while, repeating mutatis mutandis the computations leading to \eqref{entropy2 1 limit} and \eqref{entropy2 1 limit 2} in the proof of Proposition \ref{energy inequality weak interactions}, we obtain, for all $0\leq t_1<t_2$,
\begin{equation}\label{dissipation liminf}
	\int_{t_1}^{t_2}
	\delta\mathcal{D}(t)dt
	\leq \liminf_{\eps\rightarrow 0}
	\int_{t_1}^{t_2}\delta\mathcal{D}_\eps(t) dt,
\end{equation}
where
\begin{equation*}
	\begin{aligned}
		\delta\mathcal{D}(t)
		& =
		2\mu
		\left\|\nabla_x \left(u-\bar u\right)\right\|_{L^2_x}^2
		+ 5\kappa
		\left\|\nabla_x\left(\theta-\bar\theta\right)\right\|_{L^2_x}^2
		+ \frac 1\sigma
		\left\|j-\bar j\right\|_{L^2_x}^2
		\\
		& \leq
		\frac 14
		\left\|q^+ - \bar q^+\right\|_{L^2\left(MM_*dxdvdv_*d\sigma\right)}^2
		+
		\frac 14
		\left\|q^- - \bar q^-\right\|_{L^2\left(MM_*dxdvdv_*d\sigma\right)}^2
		\\
		& + \frac 14
		\left\|q^{+,-} - \bar q^{+,-}\right\|_{L^2\left(MM_*dxdvdv_*d\sigma\right)}^2
		+
		\frac 14
		\left\|q^{-,+} - \bar q^{-,+}\right\|_{L^2\left(MM_*dxdvdv_*d\sigma\right)}^2.
	\end{aligned}
\end{equation*}

The following lemma shows how the modulated entropy $\delta\mathcal{H}_\eps$ controls the modulated energy $\delta\mathcal{E}_\eps$ up to a small remainder.

\begin{lem}\label{entropy energy}
	It holds that
	\begin{equation*}
		\delta\mathcal{E}_\eps(t)
		\leq
		C \delta\mathcal{H}_\eps(t)
		+o(1)_{L^\infty(dt)},
	\end{equation*}
	for some fixed constant $C>1$.
\end{lem}

\begin{proof}
	Recall first, in view of the hypotheses on the renormalization \eqref{gamma-assumption}, that
	\begin{equation*}
		\frac 1C\left|\hat g_\eps^\pm\right|\leq \left|g_\eps^\pm\gamma_\eps^\pm\right|\leq C \left|\hat g_\eps^\pm\right|,
	\end{equation*}
	for some $C>1$, and that the elementary inequality \eqref{sqrt entropy} implies
	\begin{equation*}
		\frac 14 \hat g_\eps^{\pm 2}\leq \frac 1{\eps^2}h\left(\eps g_\eps^\pm\right).
	\end{equation*}
	
	We proceed then with the straightforward estimate, where $C_0>1$ is a large constant and $\lambda>0$ is a small parameter to be determined later on,
	\begin{equation}\label{control EH 1}
		\begin{aligned}
			\frac 12 & \left(g_\eps^\pm \gamma_\eps^\pm \chi\left(\frac{|v|^2}{K_\eps}\right)-\bar g\right)^2
			\mathds{1}_{\left\{\left|\eps g_\eps^\pm \gamma_\eps^\pm\right| > \lambda\right\}}
			\\
			& \leq
			\left(\frac {C_0}{\eps^2}h\left(\eps g_\eps^\pm\right)
			-
			g_\eps^\pm \gamma_\eps^\pm \chi\left(\frac{|v|^2}{K_\eps}\right)\bar g
			+
			\frac 12\bar g^2\right)
			\mathds{1}_{\left\{\left|\eps g_\eps^\pm \gamma_\eps^\pm\right| > \lambda\right\}}
			\\
			& \leq C_0
			\left(\frac 1{\eps^2}h\left(\eps g_\eps^\pm\right)
			-
			g_\eps^\pm \gamma_\eps^\pm \chi\left(\frac{|v|^2}{K_\eps}\right)\bar g
			+
			\frac 12\bar g^2\right)
			\mathds{1}_{\left\{\left|\eps g_\eps^\pm \gamma_\eps^\pm\right| > \lambda\right\}}
			\\
			& + (C_0-1)
			g_\eps^\pm \gamma_\eps^\pm \chi\left(\frac{|v|^2}{K_\eps}\right)\bar g
			\mathds{1}_{\left\{\left|\eps g_\eps^\pm \gamma_\eps^\pm\right| > \lambda\right\}}
			\\
			& \leq C_0
			\left(\frac 1{\eps^2}h\left(\eps g_\eps^\pm\right)
			-
			g_\eps^\pm \gamma_\eps^\pm \chi\left(\frac{|v|^2}{K_\eps}\right)\bar g
			+
			\frac 12\bar g^2\right)
			\mathds{1}_{\left\{\left|\eps g_\eps^\pm \gamma_\eps^\pm\right| > \lambda\right\}}
			+ \frac {CK\eps\left|\log \eps\right|}\lambda
			\hat g_\eps^{\pm 2}.
		\end{aligned}
	\end{equation}
	Furthermore, utilizing the simple inequality \eqref{second order h expansion ineq}, we deduce that
	\begin{equation*}
		\begin{aligned}
			\frac 12 & \left(g_\eps^\pm \gamma_\eps^\pm \chi\left(\frac{|v|^2}{K_\eps}\right)-\bar g\right)^2
			\mathds{1}_{\left\{\left|\eps g_\eps^\pm \gamma_\eps^\pm\right| \leq \lambda\right\}}
			\\
			& \leq
			\left(\frac {1}{\eps^2}h\left(\eps g_\eps^\pm \gamma_\eps^\pm \chi\left(\frac{|v|^2}{K_\eps}\right)\right)
			-
			g_\eps^\pm \gamma_\eps^\pm \chi\left(\frac{|v|^2}{K_\eps}\right)\bar g
			+
			\frac 12\bar g^2\right)
			\mathds{1}_{\left\{\left|\eps g_\eps^\pm \gamma_\eps^\pm\right| \leq \lambda\right\}}
			\\
			& + \frac \eps 6
			\left(g_\eps^\pm \gamma_\eps^\pm \chi\left(\frac{|v|^2}{K_\eps}\right)\right)^3
			\mathds{1}_{\left\{\left|\eps g_\eps^\pm \gamma_\eps^\pm\right| \leq \lambda\right\}}
			\\
			& \leq
			\left(\frac {1}{\eps^2}h\left(\eps g_\eps^\pm\right)
			-
			g_\eps^\pm \gamma_\eps^\pm \chi\left(\frac{|v|^2}{K_\eps}\right)\bar g
			+
			\frac 12\bar g^2\right)
			\mathds{1}_{\left\{\left|\eps g_\eps^\pm \gamma_\eps^\pm\right| \leq \lambda\right\}}
			\\
			& + \frac 1 6
			\left(g_\eps^\pm \gamma_\eps^\pm \chi\left(\frac{|v|^2}{K_\eps}\right) - \bar g\right)^2
			\left|\eps g_\eps^\pm \gamma_\eps^\pm\right|
			\mathds{1}_{\left\{\left|\eps g_\eps^\pm \gamma_\eps^\pm\right| \leq \lambda\right\}}
			\\
			& + \frac 1 3
			g_\eps^\pm \gamma_\eps^\pm \chi\left(\frac{|v|^2}{K_\eps}\right)\bar g
			\left|\eps g_\eps^\pm \gamma_\eps^\pm\right|
			\mathds{1}_{\left\{\left|\eps g_\eps^\pm \gamma_\eps^\pm\right| \leq \lambda\right\}}
			\\
			& \leq
			\left(\frac {1}{\eps^2}h\left(\eps g_\eps^\pm\right)
			-
			g_\eps^\pm \gamma_\eps^\pm \chi\left(\frac{|v|^2}{K_\eps}\right)\bar g
			+
			\frac 12\bar g^2\right)
			\mathds{1}_{\left\{\left|\eps g_\eps^\pm \gamma_\eps^\pm\right| \leq \lambda\right\}}
			\\
			& + \frac \lambda 6
			\left(g_\eps^\pm \gamma_\eps^\pm \chi\left(\frac{|v|^2}{K_\eps}\right) - \bar g\right)^2
			\mathds{1}_{\left\{\left|\eps g_\eps^\pm \gamma_\eps^\pm\right| \leq \lambda\right\}}
			+ CK\eps\left|\log \eps\right|
			\hat g_\eps^{\pm 2},
		\end{aligned}
	\end{equation*}
	whence, fixing $0<\lambda<3$,
	\begin{equation}\label{control EH 2}
		\begin{aligned}
			\left(\frac 12-\frac \lambda 6\right) & \left(g_\eps^\pm \gamma_\eps^\pm \chi\left(\frac{|v|^2}{K_\eps}\right)-\bar g\right)^2
			\mathds{1}_{\left\{\left|\eps g_\eps^\pm \gamma_\eps^\pm\right| \leq \lambda\right\}}
			\\
			\leq
			&
			\left(\frac {1}{\eps^2}h\left(\eps g_\eps^\pm\right)
			-
			g_\eps^\pm \gamma_\eps^\pm \chi\left(\frac{|v|^2}{K_\eps}\right)\bar g
			+
			\frac 12\bar g^2\right)
			\mathds{1}_{\left\{\left|\eps g_\eps^\pm \gamma_\eps^\pm\right| \leq \lambda\right\}}
			+ CK\eps\left|\log \eps\right|
			\hat g_\eps^{\pm 2}.
		\end{aligned}
	\end{equation}
	
	Therefore, combining estimates \eqref{control EH 1} and \eqref{control EH 2} and using the uniform bound on $\hat g_\eps^\pm$ from Lemma \ref{L2-lem}, we obtain that
	\begin{equation*}
		\begin{aligned}
			\min\left\{
			\frac{1}{C_0},
			\left(1-\frac \lambda 3\right)
			\right\}\frac 12
			& \left(g_\eps^\pm \gamma_\eps^\pm \chi\left(\frac{|v|^2}{K_\eps}\right)-\bar g\right)^2
			\\
			& \leq
			\left(\frac {1}{\eps^2}h\left(\eps g_\eps^\pm\right)
			-
			g_\eps^\pm \gamma_\eps^\pm \chi\left(\frac{|v|^2}{K_\eps}\right)\bar g
			+
			\frac 12\bar g^2\right)
			\\
			& + O\left(\eps\left|\log \eps\right|\right)_{L^\infty\left(dt;L^1\left(Mdxdv\right)\right)},
		\end{aligned}
	\end{equation*}
	which, upon integrating against $Mdxdv$ and adding the contributions of the electromagnetic field and the defect measures to the energy, concludes the proof of the lemma.
\end{proof}

The following result establishes the renormalized modulated entropy inequality at the order $\eps$, which will eventually allow us to deduce the crucial weak-strong stability of the limiting thermodynamic fields, thus defining dissipative solutions.

\begin{prop}\label{stab-prop}
	One has the stability inequality
	\begin{equation}\label{stab ineq eps}
		\begin{aligned}
			\delta\mathcal{H}_\eps(t) & + \frac 12 \int_0^t \delta\mathcal{D}_\eps(s) e^{\int_s^t\lambda(\sigma)d\sigma}ds
			\\
			& \leq \delta\mathcal{H}_\eps(0) e^{\int_0^t\lambda(s)ds}
			\\
			& +\int_0^t
			% \left[
			\int_{\mathbb{R}^3} \mathbf{A}\cdot
			\begin{pmatrix}
				\tilde u_\eps-\bar u \\ \frac 32\tilde \theta_\eps-\tilde\rho_\eps-\frac 52\bar\theta \\
				-\frac\sigma 2\int_{\mathbb{R}^3\times\mathbb{R}^3\times\mathbb{S}^2} \left(\hat q_\eps^{+,-}-\hat q_\eps^{-,+}\right) v MM_* dvdv_*d\sigma
				-\bar j
				\\ E_\eps-\bar E +\bar u \wedge \left(B_\eps-\bar B\right) \\ B_\eps-\bar B +\left(E_\eps-\bar E\right)\wedge\bar u
			\end{pmatrix}(s)
			dx
			% \right](s)
			\\
			& \times
			e^{\int_s^t\lambda(\sigma)d\sigma}ds
			+ o(1)_{L^\infty_\mathrm{loc}(dt)},
		\end{aligned}
	\end{equation}
	where the {\bf acceleration operator} is defined by
	\begin{equation*}
		\mathbf{A} \left( \bar u, \bar \theta, \bar j, \bar E, \bar B\right)
		=
		\begin{pmatrix}
			\mathbf{A}_1
			\\
			\mathbf{A}_2
			\\
			\mathbf{A}_3
			\\
			\mathbf{A}_4
			\\
			\mathbf{A}_5
		\end{pmatrix}
		=
		\begin{pmatrix}
			-2\left(\d_t \bar u +
			P\left(\bar u\cdot\nabla_x \bar u\right) - \mu\Delta_x \bar u\right)
			+ P \left(\bar j \wedge \bar B\right)
			\\
			-2\left(\partial_t\bar\theta + \bar u \cdot\nabla_x\bar \theta - \kappa\Delta_x\bar \theta\right)
			\\
			- \frac 1{\sigma}\bar j + P\left(\bar E + \bar u\wedge \bar B\right)
			\\
			-\left(\partial_t\bar E - \rot\bar B + \bar j\right)
			\\
			-\left(\partial_t\bar B + \rot\bar E\right)
		\end{pmatrix},
	\end{equation*}
	and the {\bf growth rate} is given by
	\begin{equation*}
		\begin{aligned}
			\lambda(t) =
			C\Bigg(\frac{\left\|\bar u(t)\right\|_{W^{1,\infty}\left(dx\right)}+\left\|\partial_t\bar u(t)\right\|_{L^\infty(dx)}+\left\|\bar j(t)\right\|_{L^\infty(dx)}}
			{1-\left\|\bar u(t)\right\|_{L^\infty(dx)}} &
			\\
			+\left\|\bar\theta(t)\right\|_{W^{1,\infty}(dx)} & +\left\|\bar\theta(t)\right\|_{W^{1,\infty}(dx)}^2\Bigg),
		\end{aligned}
	\end{equation*}
	with a constant $C>0$ independent of test functions and $\eps$.
\end{prop}

\begin{proof}
	The main ingredients of the proof of this stability inequality are~:
	\begin{itemize}
		
		\item
		The scaled entropy inequality \eqref{entropy2}
		\begin{equation}\label{ent}
			\begin{aligned}
				\frac1{\eps^2} & H\left(f_\eps^{+}\right)
				+ \frac1{\eps^2} H\left(f_\eps^{-}\right)
				+ \frac 1{2\eps^2}\int_{\mathbb{R}^3}\operatorname{Tr}m_\eps dx
				+ \frac 1{2} \int_{\mathbb{R}^3} \left(|E_\eps|^2+ |B_\eps|^2 + \operatorname{Tr}a_\eps\right) dx
				\\
				& +\frac{1}{4}\int_0^t\int_{\mathbb{R}^3\times\mathbb{R}^3\times\mathbb{R}^3\times\mathbb{S}^2}
				\left(
				\left(\hat q_\eps^{+}\right)^2
				+ \left(\hat q_\eps^{-}\right)^2
				+ \left(\hat q_\eps^{+,-}\right)^2
				+ \left(\hat q_\eps^{-,+}\right)^2\right)
				MM_*dxdvdv_*d\sigma ds
				% \\
				% & +\frac{1}{\epsilon^4}\int_0^t\int_{\mathbb{R}^3}\left(D\left(f_\eps^+\right)+D\left(f_\eps^-\right)
				% + \delta^2 D\left(f_\eps^+,f_\eps^-\right)\right)(s) dx ds
				\\
				& \leq
				\frac1{\eps^2} H\left(f_\eps^{+\mathrm{in}}\right)
				+ \frac1{\eps^2} H\left(f_\eps^{-\mathrm{in}}\right)
				+ \frac1{2}\int_{\mathbb{R}^3} \left(|E_\eps^{\rm in}|^2+ |B_\eps^{\rm in}|^2\right) dx,
			\end{aligned}
		\end{equation}
		which is naturally satisfied by renormalized solutions of the scaled two species Vlasov-Maxwell-Boltzmann system \eqref{VMB2} (provided they exist) and where we have used the inequality \eqref{q-est} from Lemma \ref{L2-qlem} in order to conveniently simplify the dissipation terms.
		
		\item
		The approximate conservation of energy obtained in Proposition \ref{approx2-prop}
		\begin{equation}\label{conservation-laws}
% 			\left\{
% 			\begin{aligned}
% 				\d_t \tilde \rho_\eps + \frac1\eps \nabla_x\cdot \tilde u_\eps
% 				& = R_{\eps,1},
% 				\\
				\d_t \left(\frac 32\tilde \theta_\eps-\tilde \rho_\eps\right) + \nabla_x \cdot \left( \frac52 \tilde u_\eps \tilde \theta_\eps
				- \int_{\mathbb{R}^3\times\mathbb{R}^3\times\mathbb{S}^2} \frac{\hat q_\eps^++\hat q_\eps^-}{2} \tilde \psi MM_* dvdv_*d\sigma\right)
				= R_{\eps,1},
			% \end{aligned}
			% \right.
		\end{equation}
		% \begin{equation}\label{conservation-laws}%\label{approximate2}
		% 	\left\{
		% 	\begin{aligned}
		% 		\d_t \tilde \rho_\eps + \frac1\eps \nabla_x\cdot \tilde u_\eps
		% 		& = R_{\eps,1},
		% 		\\
		% 		\d_t \tilde u_\eps
		% 		+ \nabla_x\cdot \left( \tilde u_\eps \otimes \tilde u_\eps -\frac{\left|\tilde u_\eps\right|^2}{3} \operatorname{Id}
		% 		- \int_{\mathbb{R}^3\times\mathbb{R}^3\times\mathbb{S}^2} \frac{\hat q_\eps^++\hat q_\eps^-}{2} \tilde \phi MM_* dvdv_*d\sigma \right)
		% 		\hspace{-65mm} &
		% 		\\
		% 		& = - \frac 1\eps \nabla_x\left(\tilde \rho_\eps+\tilde \theta_\eps\right)
		% 		+ \frac 12\left(\delta\tilde n_\eps E_\eps+\tilde j_\eps \wedge B_\eps\right)
		% 		+ R_{\eps,2},
		% 		\\
		% 		\d_t \left(\frac 32\tilde \theta_\eps-\tilde \rho_\eps\right) + \nabla_x \cdot \left( \frac52 \tilde u_\eps \tilde \theta_\eps
		% 		- \int_{\mathbb{R}^3\times\mathbb{R}^3\times\mathbb{S}^2} \frac{\hat q_\eps^++\hat q_\eps^-}{2} \tilde \psi MM_* dvdv_*d\sigma\right)
		% 		\hspace{-65mm} &
		% 		\\
		% 		& = R_{\eps,3},
		% 	\end{aligned}
		% 	\right.
		% \end{equation}
		where the remainder $R_{\eps,1}$ satisfies
		\begin{equation}\label{remainder1}
			\left\|R_{\eps,1}\right\|_{W^{-1,1}_\mathrm{loc}\left(dx\right)}
			\leq C\delta \mathcal{E}_\eps(t) + C\left(\delta\mathcal{E}_\eps(t)\delta\mathcal{D}_\eps(t)\right)^\frac{1}{2} + o(1)_{L^1_\mathrm{loc}(dt)},
		\end{equation}
		for some $C>0$. Note that we do not employ the approximate conservation of momentum from Proposition \ref{approx2-prop}.
		
		\item The approximate conservation of momentum law from Proposition \ref{approx3-prop}
		\begin{equation}\label{conservation-laws2}
			\begin{aligned}
				\partial_t & \left(\tilde u_\eps
				+\frac 12 E_\eps\wedge B_\eps
				+\frac 12
				\begin{pmatrix}
					a_{\eps 26}-a_{\eps 35}\\a_{\eps 34}-a_{\eps 16}\\a_{\eps 15}-a_{\eps 24}
				\end{pmatrix}
				\right)
				\\
				& +
				\nabla_x \cdot \left( \tilde u_\eps \otimes \tilde u_\eps
				-\frac{\left|\tilde u_\eps\right|^2}{3} \operatorname{Id}
				+ \frac 1{2\eps^2}m_\eps
				-
				\int_{\mathbb{R}^3\times\mathbb{R}^3\times\mathbb{S}^2} \frac{\hat q_\eps^++\hat q_\eps^-}{2} \tilde \phi MM_* dvdv_*d\sigma
				\right)
				\\
				& -\frac 12 \nabla_x\cdot\left(
				E_\eps\otimes E_\eps + e_\eps + B_\eps\otimes B_\eps + b_\eps
				\right)
				+\nabla_x\left(\frac{|E_\eps|^2+|B_\eps|^2+\operatorname{Tr}a_\eps}{4}\right)
				\\ & = - \frac 1\eps \nabla_x\left(\tilde \rho_\eps+\tilde \theta_\eps\right)+\partial_t \left(R_{\eps,2}\right) + R_{\eps,3},
			\end{aligned}
		\end{equation}
		where the remainders $R_{\eps,2}$ and $R_{\eps,3}$ satisfy
		\begin{equation}\label{remainder2}
			\begin{aligned}
				R_{\eps,2} & = o(1)_{L^\infty\left(dt;L^1_{\mathrm{loc}}(dx)\right)}
				\\
				\left\|R_{\eps,3}\right\|_{W^{-1,1}_\mathrm{loc}\left(dx\right)}
				& \leq C_1\delta\mathcal{H}_\eps(t) + C_2\delta\mathcal{E}_\eps(t) + o(1)_{L^1_\mathrm{loc}(dt)},
			\end{aligned}
		\end{equation}
		for some $C_1,C_2>0$.
		
		\item The approximate Ohm's law
		\begin{equation}\label{Ohm consistent}
			\pm
			\int_{\mathbb{R}^3\times\mathbb{R}^3\times\mathbb{S}^2}
			\hat q_\eps^{\pm,\mp}
			v MM_*dvdv_*d\sigma = \nabla_x\bar p_\eps -\left(E_\eps + \tilde u_\eps\wedge B_\eps\right)
			+R_{\eps,4}+R_{\eps,5},
		\end{equation}
		where the remainder $R_{\eps,4}$ vanishes weakly
		\begin{equation}\label{remainder4}
			R_{\eps,4}=o(1)_{\textit{w-}L^1_\mathrm{loc}(dtdx)},
		\end{equation}
		whereas $R_{\eps,5}$ satisfies
		\begin{equation}\label{remainder5}
			\left\|R_{\eps,5}\right\|_{L^1(dx)}
			\leq \frac{C\delta \mathcal{E}_\eps(t)}{1-\left\|\bar u\right\|_{L^\infty\left(dx\right)}}
			+o(1)_{\textit{w-}L_\mathrm{loc}^1(dt)}.
		\end{equation}
		
		This approximate law is obtained directly from the limiting laws derived in Proposition \ref{high weak-comp}. Indeed, it is easily deduced from \eqref{pre Ohm} that \eqref{Ohm consistent} holds with the remainders
		\begin{equation*}
			\begin{aligned}
				R_{\eps,4}
% 				& =
% 				\left(E_\eps-E
% 				\right)
% 				\pm
% 				\int_{\mathbb{R}^3\times\mathbb{R}^3\times\mathbb{S}^2}
% 				\left(\hat q_\eps^{\pm,\mp}-q^{\pm,\mp}\right)
% 				v MM_*dvdv_*d\sigma
% 				\\
% 				& + \tilde u_\eps\wedge B_\eps
% 				-u\wedge B
% 				\\
				& =
				\left(E_\eps-E
				\right)
				\pm
				\int_{\mathbb{R}^3\times\mathbb{R}^3\times\mathbb{S}^2}
				\left(\hat q_\eps^{\pm,\mp}-q^{\pm,\mp}\right)
				v MM_*dvdv_*d\sigma
				\\
				& + \left(\tilde u_\eps-u\right)\wedge B + u \wedge \left(B_\eps-B\right),
				\\
				R_{\eps,5} & =
				\left(\tilde u_\eps-u\right)\wedge\left(B_\eps - B\right)
				=
				\frac{O\left(\delta\mathcal{E}_\eps(t)+\delta\mathcal{E}(t)\right)_{L^1\left(dx\right)}}{1-\left\|\bar u\right\|_{L^\infty\left(dx\right)}},
			\end{aligned}
		\end{equation*}
		where we have used that
		\begin{equation}\label{test u estimate}
			\begin{aligned}
				& \frac{1-\left\|\bar u\right\|_{L^\infty(dx)}}{2} \left(\left\|E_\eps-\bar E\right\|_{L^2(dx)}^2+\left\|B_\eps-\bar B\right\|_{L^2(dx)}^2\right)
				\\
				& \leq\frac 12 \left(\left\|E_\eps-\bar E\right\|_{L^2(dx)}^2+\left\|B_\eps-\bar B\right\|_{L^2(dx)}^2\right)
				-\int_{\mathbb{R}^3}\left(\left(E_\eps-\bar E\right)\wedge\left(B_\eps-\bar B\right)\right)\cdot\bar u dx.
			\end{aligned}
		\end{equation}
		The above estimate on $R_{\eps,5}$ is then readily improved to \eqref{remainder5} upon noticing from \eqref{energy liminf} that
		\begin{equation*}
			\delta\mathcal{E}(t)\leq \delta\mathcal{E}_0(t),
		\end{equation*}
		where $\delta\mathcal{E}_0(t)$ is the limit, up to extraction of subsequences, of $\delta\mathcal{E}_\eps(t)$ in $\textit{w$^*$-}L^\infty(dt)$, and then writing
		\begin{equation*}
			\begin{aligned}
				\delta\mathcal{E}_\eps(t)+\delta\mathcal{E}(t)
				& \leq
				2\delta \mathcal{E}_\eps(t)+\delta\mathcal{E}_0(t)-\delta\mathcal{E}_\eps(t)
				\\
				& =2\delta \mathcal{E}_\eps(t)+o(1)_{\textit{w$^*$-}L^\infty(dt)}.
			\end{aligned}
		\end{equation*}
		
		% , for every non-negative $\varphi\in C_c\left([0,\infty)\right)$,
		% \begin{equation*}
		% 	\lim_{\eps\rightarrow 0}\left(\int_{0}^{\infty}\varphi(t)
		% 	\left(\delta\mathcal{E}(t)-\delta\mathcal{E}_\eps(t)\right)
		% 	dt\right)_+ = 0.
		% \end{equation*}
		
		It is to be emphasized that it would be possible to derive the above approximate Ohm's law directly from Proposition \ref{approx2-prop}. However, the method used here is more robust and we find it much more satisfiying to derive the consistency of the approximate law from the knowledge of the limiting equation. Indeed, morally, it is much more difficult to the derive limiting equations, which require some kind of weak stability property, than to merely establish the consistency of approximate equations.

		\item Maxwell's equations
		\begin{equation}\label{Maxwell consistent}
			\begin{cases}
				\begin{aligned}
					\d_t E_\eps - \ROT B_\eps &= -j_\eps = -\tilde j_\eps +R_{\eps,6},
					\\
					\d_t B_\eps + \ROT E_\eps & = 0,
					\\
					\Div E_\eps & = \delta n_\eps,
					\\
					\Div B_\eps & = 0,
				\end{aligned}
			\end{cases}
		\end{equation}
		where the remainder $R_{\eps,6}=\tilde j_\eps - j_\eps$ satisfies, thanks to the convergence \eqref{tilde u approx u},
		\begin{equation}\label{remainder6}
			R_{\eps,6}=o(1)_{L^1_\mathrm{loc}(dtdx)}.
		\end{equation}
		
		Notice that we cannot rigorously write the identities \eqref{maxwell energy two species} and \eqref{maxwell poynting two species} for the above system, because the source terms $j_\eps$ and $n_\eps$ do not belong to $L^2_\mathrm{loc}(dtdx)$ a priori. Nevertheless, one has the following modulated identities~:
		\begin{equation}\label{modulated maxwell energy}
			\begin{aligned}
				\partial_t\left(E_\eps\cdot \bar E+B_\eps\cdot\bar B\right)
				& +
				\nabla_x\cdot\left(E_\eps\wedge \bar B + \bar E\wedge B_\eps\right)
				\\
				& = -\left(\tilde j_\eps -R_{\eps,6}\right)\cdot\bar E
				- \left(\bar j +\mathbf{A}_4\right)\cdot E_\eps
				- \mathbf{A}_5\cdot B_\eps,
			\end{aligned}
		\end{equation}
		and
		\begin{equation}\label{modulated poynting}
			\begin{aligned}
				\partial_t & \left(\left(E_\eps-\bar E\right)\wedge \left(B_\eps-\bar B\right) \right)
				+\frac 12\nabla_x\left(\left|E_\eps-\bar E\right|^2 + \left|B_\eps-\bar B\right|^2 \right)
				\\
				& -\nabla_x\cdot\left(\left(E_\eps-\bar E\right)\otimes \left(E_\eps-\bar E\right)+\left(B_\eps-\bar B\right)\otimes \left(B_\eps-\bar B\right)\right)
				\\
				& =
				\partial_t\left(E_\eps\wedge B_\eps \right)
				+\frac 12\nabla_x\left(\left|E_\eps\right|^2 + \left|B_\eps\right|^2 \right)
				-\nabla_x\cdot\left(E_\eps\otimes E_\eps+B_\eps\otimes B_\eps\right)
				\\
				& + \left(\bar j+\mathbf{A}_4\right)\wedge \left(B_\eps-\bar B\right) +\left(E_\eps-\bar E\right)\wedge \mathbf{A}_5
				\\
				& +\left(\tilde j_\eps-R_{\eps, 6}\right)\wedge \bar B + \delta n_\eps\bar E.
			\end{aligned}
		\end{equation}

		\item The asymptotic characterization \eqref{q phi psi 2} of the limiting collision integrands $q^\pm$ from Proposition \ref{weak-comp2}, which implies that
		\begin{equation}\label{asymptotic mu kappa}
			\begin{aligned}
				\left(\int_{\mathbb{R}^3\times\mathbb{R}^3\times\mathbb{S}^2} \frac{\hat q_\eps^++\hat q_\eps^-}{2} \tilde \phi MM_* dvdv_*d\sigma\right)
				-\mu\left(\nabla_x\tilde u_\eps+\nabla_x^t\tilde u_\eps\right)
				& \to 0,
				\\
				\left(\int_{\mathbb{R}^3\times\mathbb{R}^3\times\mathbb{S}^2} \frac{\hat q_\eps^++\hat q_\eps^-}{2} \tilde \psi MM_* dvdv_*d\sigma\right)
				-\frac 52\kappa\nabla_x\tilde\theta_\eps
				& \to 0,
			\end{aligned}
		\end{equation}
		in the sense of distributions, where $\mu,\kappa>0$ are defined by \eqref{mu kappa 2}.

		\item The asymptotic characterization \eqref{mixed q phi psi} of the limiting collision integrands $q^{\pm,\mp}$ from Proposition \ref{weak-comp3}, which implies that
		\begin{equation}\label{asymptotic sigma}
			\begin{aligned}
				\sigma\left(\int_{\mathbb{R}^3\times\mathbb{R}^3\times\mathbb{S}^2} \frac{\hat q_\eps^{+,-}-\hat q_\eps^{-,+}}{2} v MM_* dvdv_*d\sigma\right)+\tilde j_\eps
				\to 0,
			\end{aligned}
		\end{equation}
		in the sense of distributions, where $\sigma>0$ is defined by \eqref{sigma}. Moreover, since $j$ is solenoidal, it holds that
		\begin{equation}\label{limit j solenoidal}
			P^\perp \left(\int_{\mathbb{R}^3\times\mathbb{R}^3\times\mathbb{S}^2} \frac{\hat q_\eps^{+,-}-\hat q_\eps^{-,+}}{2} v MM_* dvdv_*d\sigma\right)
			\rightharpoonup 0
			\quad\text{in }L^2(dtdx).
		\end{equation}
	\end{itemize}

	Now, by definition of the acceleration operator $\mathbf{A}$, straightforward energy computations, similar to those performed in the proof of Proposition \ref{energy estimate 2}, applied to the test functions $\left( \bar u, \bar \theta, \bar j, \bar E, \bar B\right)$, show that the following energy identity holds~:
	\begin{equation}\label{test energy}
		\frac{d}{dt}\bar\CE(t)+\bar\CD(t)= - \int_{\mathbb{R}^3} \mathbf{A}\cdot
		\begin{pmatrix}
			\bar u \\ \frac 52 \bar \theta \\ \bar j \\ \bar E \\ \bar B
		\end{pmatrix}
		dx,
	\end{equation}
	where the energy $\bar\CE$ and energy dissipation $\bar\CD$ are defined by
	\begin{equation*}
		\begin{aligned}
			\bar{\mathcal{E}}(t)
			& =
			\left\|\bar g\right\|_{L^2\left(Mdxdv\right)}^2
			+\frac 12\left\|\bar E\right\|_{L^2(dx)}^2
			+\frac 12\left\|\bar B\right\|_{L^2(dx)}^2
			\\
			& =
			\left\|\bar u\right\|_{L^2(dx)}^2+\frac 52\left\|\bar \theta\right\|_{L^2(dx)}^2
			+\frac 12\left\|\bar E\right\|_{L^2(dx)}^2
			+\frac 12\left\|\bar B\right\|_{L^2(dx)}^2,
		\end{aligned}
	\end{equation*}
	and
	\begin{equation*}
		\begin{aligned}
			\bar{\mathcal{D}}(t)
			& =
			2\mu
			\left\|\nabla_x \bar u\right\|_{L^2_x}^2
			+ 5\kappa
			\left\|\nabla_x \bar\theta\right\|_{L^2_x}^2
			+ \frac 1\sigma
			\left\|\bar j\right\|_{L^2_x}^2
			\\
			& =
			\frac 14
			\left\|\bar q^+\right\|_{L^2\left(MM_*dxdvdv_*d\sigma\right)}^2
			+
			\frac 14
			\left\|\bar q^-\right\|_{L^2\left(MM_*dxdvdv_*d\sigma\right)}^2
			\\
			& + \frac 14
			\left\|\bar q^{+,-}\right\|_{L^2\left(MM_*dxdvdv_*d\sigma\right)}^2
			+
			\frac 14
			\left\|\bar q^{-,+}\right\|_{L^2\left(MM_*dxdvdv_*d\sigma\right)}^2.
		\end{aligned}
	\end{equation*}

	Next, similar duality computations applied to the approximate conservation of energy \eqref{conservation-laws} yield that
	\begin{equation*}
		\begin{aligned}
			\frac{d}{dt} & \int_{\mathbb{R}^3} \left(\frac 32 \tilde\theta_\eps-\tilde\rho_\eps\right)\cdot\bar\theta dx
			+
			\int_{\mathbb{R}^3}
			\left(\bar u\cdot\nabla_x\bar\theta\left(\frac 32 \tilde\theta_\eps-\tilde\rho_\eps\right)
			-\kappa\Delta_x \bar\theta\left(\frac 32 \tilde\theta_\eps-\tilde\rho_\eps\right)\right)
			dx
			\\
			& -
			\int_{\mathbb{R}^3}
			\left( \frac52 \tilde u_\eps \tilde \theta_\eps
			- \int_{\mathbb{R}^3\times\mathbb{R}^3\times\mathbb{S}^2} \frac{\hat q_\eps^++\hat q_\eps^-}{2} \tilde \psi MM_* dvdv_*d\sigma\right)
			\cdot\nabla_x\bar\theta
			dx
			\\
			& = \int_{\mathbb{R}^3}R_{\eps,1}\cdot \bar\theta
			-\frac{1}{2}\mathbf{A}_2\left(\frac 32 \tilde\theta_\eps-\tilde\rho_\eps\right)
			dx.
		\end{aligned}
	\end{equation*}
	Further reorganizing the preceding equation so that all remainder terms appear on its right-hand side, we find that
	\begin{equation*}
		\begin{aligned}
			\frac{d}{dt} & \int_{\mathbb{R}^3} \left(\frac 32 \tilde\theta_\eps-\tilde\rho_\eps\right)\cdot\bar\theta dx
			\\
			& +
			\int_{\mathbb{R}^3}
			\left( \int_{\mathbb{R}^3\times\mathbb{R}^3\times\mathbb{S}^2} \left(\hat q_\eps^++\hat q_\eps^-\right) \tilde \psi MM_* dvdv_*d\sigma\right)
			\cdot\nabla_x\bar\theta
			dx
			\\
			& = \int_{\mathbb{R}^3}R_{\eps,1}\cdot \bar\theta
			-\frac{1}{2}\mathbf{A}_2\left(\frac 32 \tilde\theta_\eps-\tilde\rho_\eps\right)
			+\frac 52
			\left(\tilde u_\eps-\bar u\right)\cdot\nabla_x\bar\theta\left( \tilde\theta_\eps-\bar\theta\right)
			dx
			\\
			& +
			\int_{\mathbb{R}^3}
			\frac 54 P^\perp\tilde u_\eps \cdot\nabla_x \left(\bar\theta^2 \right)
			+
			\left(\bar u \cdot\nabla_x\bar\theta-\kappa\Delta_x\bar\theta\right) \left( \tilde\theta_\eps+\tilde\rho_\eps\right)
			dx
			\\
			& +
			\int_{\mathbb{R}^3}
			\frac 52\kappa\tilde\theta_\eps\Delta_x \bar\theta
			+
			\left( \int_{\mathbb{R}^3\times\mathbb{R}^3\times\mathbb{S}^2} \frac{\hat q_\eps^++\hat q_\eps^-}{2} \tilde \psi MM_* dvdv_*d\sigma\right)
			\cdot\nabla_x\bar\theta
			dx.
		\end{aligned}
	\end{equation*}
	It then follows, using the convergences \eqref{limit solenoidal}, \eqref{limit boussinesq}, \eqref{asymptotic mu kappa}, the estimate \eqref{remainder1} and Lemma \ref{entropy energy} (allowing to control the energy by the entropy), that
	\begin{equation}\label{conservation-laws4}
		\begin{aligned}
			\frac{d}{dt} \int_{\mathbb{R}^3} \left(\frac 32 \tilde\theta_\eps-\tilde\rho_\eps\right)\cdot\bar\theta dx
			& +
			\int_{\mathbb{R}^3}
			\left( \int_{\mathbb{R}^3\times\mathbb{R}^3\times\mathbb{S}^2} \left(\hat q_\eps^++\hat q_\eps^-\right) \tilde \psi MM_* dvdv_*d\sigma\right)
			\cdot\nabla_x\bar\theta
			dx
			\\
			& \geq -C\left\|\bar\theta\right\|_{W^{1,\infty}(dx)}\left(\delta\mathcal{E}_\eps(t)+\left(\delta\mathcal{E}_\eps(t)\delta\mathcal{D}_\eps(t)\right)^\frac{1}{2}\right)
			\\
			& - \frac{1}{2}\int_{\mathbb{R}^3}
			\mathbf{A}_2\left(\frac 32 \tilde\theta_\eps-\tilde\rho_\eps\right)
			dx
			+o(1)_{\textit{w-}L^1_\mathrm{loc}(dt)}
			\\
			& \geq -C\left(\left\|\bar\theta\right\|_{W^{1,\infty}(dx)}
			+\left\|\bar\theta\right\|_{W^{1,\infty}(dx)}^2\right)
			\delta\mathcal{H}_\eps(t)
			-\frac 14 \delta\mathcal{D}_\eps(t)
			\\
			& - \frac{1}{2}\int_{\mathbb{R}^3}
			\mathbf{A}_2\left(\frac 32 \tilde\theta_\eps-\tilde\rho_\eps\right)
			dx
			+o(1)_{\textit{w-}L^1_\mathrm{loc}(dt)}.
		\end{aligned}
	\end{equation}

	Likewise, using the solenoidal property $\Div \bar u = 0$, analogous duality computations applied to the approximate conservation of momentum \eqref{conservation-laws2} yield that
	\begin{equation*}
		\begin{aligned}
			\frac{d}{dt} & \int_{\mathbb{R}^3}\left(\tilde u_\eps\cdot\bar u+\frac 12 \left(E_\eps\wedge B_\eps\right)\cdot \bar u
			+\frac 12
			\begin{pmatrix}
				a_{\eps 26}-a_{\eps 35}\\a_{\eps 34}-a_{\eps 16}\\a_{\eps 15}-a_{\eps 24}
			\end{pmatrix}\cdot\bar u - R_{\eps,2}\cdot\bar u\right) dx
			\\
			& - \int_{\mathbb{R}^3}\left(\frac 12 \left(E_\eps\wedge B_\eps\right)\cdot \partial_t\bar u
			+\frac 12
			\begin{pmatrix}
				a_{\eps 26}-a_{\eps 35}\\a_{\eps 34}-a_{\eps 16}\\a_{\eps 15}-a_{\eps 24}
			\end{pmatrix}\cdot\partial_t\bar u - R_{\eps,2}\cdot\partial_t\bar u\right) dx
			\\
			& +
			\int_{\mathbb{R}^3}
			\left(\left(P\tilde u_\eps\right)\otimes \bar u- \tilde u_\eps\otimes\tilde u_\eps - \frac{1}{2\eps^2}m_\eps\right) :\nabla_x\bar u
			dx
			\\
			& +
			\int_{\mathbb{R}^3}
			\left(\int_{\mathbb{R}^3\times\mathbb{R}^3\times\mathbb{S}^2} \frac{\hat q_\eps^++\hat q_\eps^-}{2} \tilde \phi MM_* dvdv_*d\sigma\right)
			:\nabla_x\bar u
			-\mu \Delta_x\bar u \cdot\tilde u_\eps dx
			\\
			& +\frac 12
			\int_{\mathbb{R}^3}
			\left(E_\eps\otimes E_\eps + e_\eps + B_\eps\otimes B_\eps + b_\eps\right) :\nabla_x\bar u
			dx
			\\
			& =
			\int_{\mathbb{R}^3}
			R_{\eps,3}\cdot\bar u + \frac 12P\left(\bar j\wedge\bar B\right)\cdot\tilde u_\eps
			-\frac{1}{2}\mathbf{A}_1\cdot \tilde u_\eps
			dx,
		\end{aligned}
	\end{equation*}
	whence, reorganizing some terms so that remainders are moved to the right-hand side,
	\begin{equation*}
		\begin{aligned}
			\frac{d}{dt} & \int_{\mathbb{R}^3}\left(\tilde u_\eps\cdot\bar u+\frac 12 \left(E_\eps\wedge B_\eps\right)\cdot \bar u
			+\frac 12
			\begin{pmatrix}
				a_{\eps 26}-a_{\eps 35}\\a_{\eps 34}-a_{\eps 16}\\a_{\eps 15}-a_{\eps 24}
			\end{pmatrix}\cdot\bar u \right) dx
			\\
			& +
			\int_{\mathbb{R}^3}
			\left(\int_{\mathbb{R}^3\times\mathbb{R}^3\times\mathbb{S}^2} \left(\hat q_\eps^++\hat q_\eps^-\right) \tilde \phi MM_* dvdv_*d\sigma\right)
			:\nabla_x\bar u
			dx
			\\
			& =
			\int_{\mathbb{R}^3}
			R_{\eps,3}\cdot\bar u + \frac 12P\left(\bar j\wedge\bar B\right)\cdot\tilde u_\eps
			-\frac{1}{2}\mathbf{A}_1\cdot \tilde u_\eps
			dx
			\\
			& +
			\int_{\mathbb{R}^3}
			\left(\bar u\otimes \left(P^\perp \tilde u_\eps\right) + \left(P^\perp\tilde u_\eps\right)\otimes\bar u +\left(\tilde u_\eps-\bar u\right)\otimes \left(\tilde u_\eps-\bar u\right)
			+ \frac{1}{2\eps^2}m_\eps\right) :\nabla_x\bar u
			dx
			\\
			& +
			\int_{\mathbb{R}^3}
			\mu\tilde u_\eps\cdot \Delta_x\bar u
			+
			\left(\int_{\mathbb{R}^3\times\mathbb{R}^3\times\mathbb{S}^2} \frac{\hat q_\eps^++\hat q_\eps^-}{2} \tilde \phi MM_* dvdv_*d\sigma\right)
			:\nabla_x\bar u dx
			\\
			& + \int_{\mathbb{R}^3}\left(\frac 12 \left(E_\eps\wedge B_\eps\right)\cdot \partial_t\bar u
			+\frac 12
			\begin{pmatrix}
				a_{\eps 26}-a_{\eps 35}\\a_{\eps 34}-a_{\eps 16}\\a_{\eps 15}-a_{\eps 24}
			\end{pmatrix}\cdot\partial_t\bar u - R_{\eps,2}\cdot\partial_t\bar u\right) dx
			\\
			& -\frac 12
			\int_{\mathbb{R}^3}
			\left(E_\eps\otimes E_\eps + e_\eps + B_\eps\otimes B_\eps + b_\eps\right) :\nabla_x\bar u
			dx
			+\frac {d}{dt}\int_{\mathbb{R}^3}R_{\eps,2}\cdot\bar udx.
		\end{aligned}
	\end{equation*}
	Then, using the convergences \eqref{limit solenoidal}, \eqref{asymptotic mu kappa}, the estimates \eqref{measureCS1}, \eqref{measureCS2}, \eqref{remainder2} and Lemmas \ref{vector defect} and \ref{entropy energy} (allowing to control the energy by the entropy), we arrive at
	\begin{equation*}
		\begin{aligned}
			\frac{d}{dt} & \int_{\mathbb{R}^3}\left(\tilde u_\eps\cdot\bar u+\frac 12 \left(E_\eps\wedge B_\eps\right)\cdot \bar u
			+\frac 12
			\begin{pmatrix}
				a_{\eps 26}-a_{\eps 35}\\a_{\eps 34}-a_{\eps 16}\\a_{\eps 15}-a_{\eps 24}
			\end{pmatrix}\cdot\bar u \right) dx
			\\
			& +
			\int_{\mathbb{R}^3}
			\left(\int_{\mathbb{R}^3\times\mathbb{R}^3\times\mathbb{S}^2} \left(\hat q_\eps^++\hat q_\eps^-\right) \tilde \phi MM_* dvdv_*d\sigma\right)
			:\nabla_x\bar u
			dx
			\\
			& \geq
			- C\left(\left\|\bar u\right\|_{W^{1,\infty}\left(dx\right)}+\frac{\left\|\partial_t\bar u\right\|_{L^\infty(dx)}}{1-\left\|\bar u\right\|_{L^\infty(dx)}}\right)\delta\mathcal{H}_\eps(t)
			-\frac{1}{2}\int_{\mathbb{R}^3}\mathbf{A}_1\cdot \tilde u_\eps dx
			\\
			& + o(1)_{\textit{w-}L^1_\mathrm{loc}(dt)}
			+\frac{d}{dt}\left(o(1)_{L^\infty\left(dt\right)}\right)
			+\frac 12 \int_{\mathbb{R}^3}
			P\left(\bar j\wedge\bar B\right)\cdot\tilde u_\eps
			dx
			\\
			& + \frac 12\int_{\mathbb{R}^3}\left(E_\eps\wedge B_\eps\right)\cdot \partial_t\bar u
			-
			\left(E_\eps\otimes E_\eps + e_\eps + B_\eps\otimes B_\eps + b_\eps\right) :\nabla_x\bar u
			dx.
		\end{aligned}
	\end{equation*}

	The next step consists in combining the preceding inequality with the identity \eqref{modulated poynting} in order to modulate the Poynting vector $E_\eps\wedge B_\eps$. This yields
	\begin{equation*}
		\begin{aligned}
			\frac{d}{dt} & \int_{\mathbb{R}^3}\left(\tilde u_\eps\cdot\bar u+\frac 12 \left(\left(E_\eps-\bar E\right)\wedge \left(B_\eps-\bar B\right)\right)\cdot \bar u
			+\frac 12
			\begin{pmatrix}
				a_{\eps 26}-a_{\eps 35}\\a_{\eps 34}-a_{\eps 16}\\a_{\eps 15}-a_{\eps 24}
			\end{pmatrix}\cdot\bar u \right) dx
			\\
			& +
			\int_{\mathbb{R}^3}
			\left(\int_{\mathbb{R}^3\times\mathbb{R}^3\times\mathbb{S}^2} \left(\hat q_\eps^++\hat q_\eps^-\right) \tilde \phi MM_* dvdv_*d\sigma\right)
			:\nabla_x\bar u
			dx
			\\
			& \geq
			- C\left(\left\|\bar u\right\|_{W^{1,\infty}\left(dx\right)}+\frac{\left\|\partial_t\bar u\right\|_{L^\infty(dx)}}{1-\left\|\bar u\right\|_{L^\infty(dx)}}\right)\delta\mathcal{H}_\eps(t)
			\\
			& - \frac{1}{2}\int_{\mathbb{R}^3}
			\mathbf{A}_1\cdot \tilde u_\eps
			dx
			+\frac{1}{2}\int_{\mathbb{R}^3}\left(\mathbf{A}_4\wedge \left(B_\eps-\bar B\right) +\left(E_\eps-\bar E\right)\wedge \mathbf{A}_5\right)\cdot\bar u dx
			\\
			& + o(1)_{\textit{w-}L^1_\mathrm{loc}(dt)}
			+\frac{d}{dt}\left(o(1)_{L^\infty\left(dt\right)}\right)
			\\
			& +\frac 12 \int_{\mathbb{R}^3}\left( \delta n_\eps\bar E - R_{\eps, 6}\wedge \bar B \right)\cdot\bar u
			-\left(\tilde u_\eps\wedge B_\eps\right)\cdot\bar j
			- \left(\bar u\wedge \bar B\right)\cdot\tilde j_\eps
			- \left(\bar j\wedge\bar B\right)\cdot P^\perp\tilde u_\eps
			dx
			\\
			& + \frac 12\int_{\mathbb{R}^3}
			\left(\bar j\wedge \left(B_\eps-\bar B\right)\right)\cdot\left(\bar u-\tilde u_\eps\right)
			+
			\left(\left(E_\eps-\bar E\right)\wedge \left(B_\eps-\bar B\right)\right)\cdot \partial_t\bar udx
			\\
			& - \frac 12\int_{\mathbb{R}^3}
			\left(\left(E_\eps-\bar E\right)\otimes \left(E_\eps-\bar E\right) + e_\eps + \left(B_\eps-\bar B\right)\otimes \left(B_\eps-\bar B\right) + b_\eps\right) :\nabla_x\bar u
			dx.
		\end{aligned}
	\end{equation*}
	It then follows, using the convergence \eqref{limit solenoidal}, the estimates \eqref{measureCS2}, \eqref{remainder6} and Lemmas \ref{vector defect} and \ref{entropy energy} (allowing to control the energy by the entropy), that
	\begin{equation*}
		\begin{aligned}
			\frac{d}{dt} & \int_{\mathbb{R}^3}\left(\tilde u_\eps\cdot\bar u+\frac 12 \left(\left(E_\eps-\bar E\right)\wedge \left(B_\eps-\bar B\right)\right)\cdot \bar u
			+\frac 12
			\begin{pmatrix}
				a_{\eps 26}-a_{\eps 35}\\a_{\eps 34}-a_{\eps 16}\\a_{\eps 15}-a_{\eps 24}
			\end{pmatrix}\cdot\bar u \right) dx
			\\
			& +
			\int_{\mathbb{R}^3}
			\left(\int_{\mathbb{R}^3\times\mathbb{R}^3\times\mathbb{S}^2} \left(\hat q_\eps^++\hat q_\eps^-\right) \tilde \phi MM_* dvdv_*d\sigma\right)
			:\nabla_x\bar u
			dx
			\\
			& \geq
			- C\frac{\left\|\bar u\right\|_{W^{1,\infty}\left(dx\right)}+\left\|\partial_t\bar u\right\|_{L^\infty(dx)}+\left\|\bar j\right\|_{L^\infty(dx)}}{1-\left\|\bar u\right\|_{L^\infty(dx)}} \delta\mathcal{H}_\eps(t)
			\\
			& - \frac{1}{2}\int_{\mathbb{R}^3}
			\mathbf{A}_1\cdot \tilde u_\eps
			dx
			+ \frac{1}{2}\int_{\mathbb{R}^3}\left(\mathbf{A}_4\wedge \left(B_\eps-\bar B\right) +\left(E_\eps-\bar E\right)\wedge \mathbf{A}_5\right)\cdot\bar u dx
			\\
			& + o(1)_{\textit{w-}L^1_\mathrm{loc}(dt)}
			+\frac{d}{dt}\left(o(1)_{L^\infty\left(dt\right)}\right)
			- \frac 12 \int_{\mathbb{R}^3}
			\left(\tilde u_\eps\wedge B_\eps\right)\cdot\bar j
			+ \left(\bar u\wedge \bar B\right)\cdot\tilde j_\eps
			dx.
		\end{aligned}
	\end{equation*}

	Now, for mere convenience of notation, we introduce the following integrand~:
	\begin{equation*}
		\begin{aligned}
			\mathcal{I}
			& =
			\tilde u_\eps\cdot\bar u
			+\frac 12 \left(E_\eps\cdot \bar E+B_\eps\cdot\bar B\right)
			\\
			& +\frac 12 \left(\left(E_\eps-\bar E\right)\wedge \left(B_\eps-\bar B\right)\right)\cdot \bar u
			+\frac 12
			\begin{pmatrix}
				a_{\eps 26}-a_{\eps 35}\\a_{\eps 34}-a_{\eps 16}\\a_{\eps 15}-a_{\eps 24}
			\end{pmatrix}\cdot\bar u.
		\end{aligned}
	\end{equation*}
	Thus, further employing the identity \eqref{modulated maxwell energy}, we find that
	\begin{equation*}
		\begin{aligned}
			% \frac{d}{dt} & \int_{\mathbb{R}^3}
			% \left(\tilde u_\eps\cdot\bar u
			% +\frac 12 \left(E_\eps\cdot \bar E+B_\eps\cdot\bar B\right)
			% +\frac 12 \left(\left(E_\eps-\bar E\right)\wedge \left(B_\eps-\bar B\right)\right)\cdot \bar u
			% \vphantom{+\frac 12
			% \begin{pmatrix}
			% 	a_{\eps 26}-a_{\eps 35}\\a_{\eps 34}-a_{\eps 16}\\a_{\eps 15}-a_{\eps 24}
			% \end{pmatrix}\cdot\bar u}
			% \right.
			% \\
			% & +\left.
			% \vphantom{\tilde u_\eps\cdot\bar u
			% +\frac 12 \left(E_\eps\cdot \bar E+B_\eps\cdot\bar B\right)
			% +\frac 12 \left(\left(E_\eps-\bar E\right)\wedge \left(B_\eps-\bar B\right)\right)\cdot \bar u}
			% \frac 12
			% \begin{pmatrix}
			% 	a_{\eps 26}-a_{\eps 35}\\a_{\eps 34}-a_{\eps 16}\\a_{\eps 15}-a_{\eps 24}
			% \end{pmatrix}\cdot\bar u \right)
			% dx
			% \\
			\frac{d}{dt} & \int_{\mathbb{R}^3}\mathcal{I}dx
			+
			\int_{\mathbb{R}^3}
			\left(\int_{\mathbb{R}^3\times\mathbb{R}^3\times\mathbb{S}^2} \left(\hat q_\eps^++\hat q_\eps^-\right) \tilde \phi MM_* dvdv_*d\sigma\right)
			:\nabla_x\bar u
			dx
			\\
			& \geq
			- C\frac{\left\|\bar u\right\|_{W^{1,\infty}\left(dx\right)}+\left\|\partial_t\bar u\right\|_{L^\infty(dx)}+\left\|\bar j\right\|_{L^\infty(dx)}}{1-\left\|\bar u\right\|_{L^\infty(dx)}} \delta\mathcal{H}_\eps(t)
			\\
			& - \frac{1}{2} \int_{\mathbb{R}^3}
			\mathbf{A}_1\cdot \tilde u_\eps
			+
			\mathbf{A}_4\cdot\left(E_\eps+\bar u \wedge \left(B_\eps-\bar B\right)\right)
			+\mathbf{A}_5\cdot\left(B_\eps+\left(E_\eps-\bar E\right)\wedge\bar u\right) dx
			\\
			& + o(1)_{\textit{w-}L^1_\mathrm{loc}(dt)}
			+\frac{d}{dt}\left(o(1)_{L^\infty\left(dt\right)}\right)
			\\
			&
			+\frac 12
			\int_{\mathbb{R}^3}R_{\eps,6}\cdot\bar E
			-\left(E_\eps+\tilde u_\eps\wedge B_\eps\right)\cdot\bar j
			-\left(\bar E + \bar u\wedge \bar B\right)\cdot\tilde j_\eps
			dx,
		\end{aligned}
	\end{equation*}
	whence, in view of the estimate \eqref{remainder6},
	\begin{equation*}
		\begin{aligned}
			% \frac{d}{dt} & \int_{\mathbb{R}^3}
			% \left(\tilde u_\eps\cdot\bar u
			% +\frac 12 \left(E_\eps\cdot \bar E+B_\eps\cdot\bar B\right)
			% +\frac 12 \left(\left(E_\eps-\bar E\right)\wedge \left(B_\eps-\bar B\right)\right)\cdot \bar u
			% \vphantom{+\frac 12
			% \begin{pmatrix}
			% 	a_{\eps 26}-a_{\eps 35}\\a_{\eps 34}-a_{\eps 16}\\a_{\eps 15}-a_{\eps 24}
			% \end{pmatrix}\cdot\bar u}
			% \right.
			% \\
			% & + \left.
			% \vphantom{\tilde u_\eps\cdot\bar u
			% +\frac 12 \left(E_\eps\cdot \bar E+B_\eps\cdot\bar B\right)
			% +\frac 12 \left(\left(E_\eps-\bar E\right)\wedge \left(B_\eps-\bar B\right)\right)\cdot \bar u}
			% \frac 12
			% \begin{pmatrix}
			% 	a_{\eps 26}-a_{\eps 35}\\a_{\eps 34}-a_{\eps 16}\\a_{\eps 15}-a_{\eps 24}
			% \end{pmatrix}\cdot\bar u \right)
			% dx
			% \\
			\frac{d}{dt} & \int_{\mathbb{R}^3}\mathcal{I}dx
			+
			\int_{\mathbb{R}^3}
			\left(\int_{\mathbb{R}^3\times\mathbb{R}^3\times\mathbb{S}^2} \left(\hat q_\eps^++\hat q_\eps^-\right) \tilde \phi MM_* dvdv_*d\sigma\right)
			:\nabla_x\bar u
			dx
			\\
			& \geq
			- C\frac{\left\|\bar u\right\|_{W^{1,\infty}\left(dx\right)}+\left\|\partial_t\bar u\right\|_{L^\infty(dx)}+\left\|\bar j\right\|_{L^\infty(dx)}}{1-\left\|\bar u\right\|_{L^\infty(dx)}} \delta\mathcal{H}_\eps(t)
			\\
			& - \frac{1}{2} \int_{\mathbb{R}^3}
			\mathbf{A}_1\cdot \tilde u_\eps
			+
			\mathbf{A}_4\cdot\left(E_\eps+\bar u \wedge \left(B_\eps-\bar B\right)\right)
			+\mathbf{A}_5\cdot\left(B_\eps+\left(E_\eps-\bar E\right)\wedge\bar u\right) dx
			\\
			& + o(1)_{\textit{w-}L^1_\mathrm{loc}(dt)}
			+\frac{d}{dt}\left(o(1)_{L^\infty\left(dt\right)}\right)
			\\
			& -\frac 12
			\int_{\mathbb{R}^3}
			\left(E_\eps+\tilde u_\eps\wedge B_\eps\right)\cdot\bar j
			+\left(\bar E + \bar u\wedge \bar B\right)\cdot\tilde j_\eps
			dx.
		\end{aligned}
	\end{equation*}

	Using then the approximate Ohm's law \eqref{Ohm consistent} and reorganizing the resulting inequality so that all remainder terms appear on its right-hand side, we obtain
	\begin{equation*}
		\begin{aligned}
			% \frac{d}{dt} & \int_{\mathbb{R}^3}
			% \left(\tilde u_\eps\cdot\bar u
			% +\frac 12 \left(E_\eps\cdot \bar E+B_\eps\cdot\bar B\right)
			% +\frac 12 \left(\left(E_\eps-\bar E\right)\wedge \left(B_\eps-\bar B\right)\right)\cdot \bar u
			% \vphantom{+\frac 12
			% \begin{pmatrix}
			% 	a_{\eps 26}-a_{\eps 35}\\a_{\eps 34}-a_{\eps 16}\\a_{\eps 15}-a_{\eps 24}
			% \end{pmatrix}\cdot\bar u}
			% \right.
			% \\
			% & + \left.
			% \vphantom{\tilde u_\eps\cdot\bar u
			% +\frac 12 \left(E_\eps\cdot \bar E+B_\eps\cdot\bar B\right)
			% +\frac 12 \left(\left(E_\eps-\bar E\right)\wedge \left(B_\eps-\bar B\right)\right)\cdot \bar u}
			% \frac 12
			% \begin{pmatrix}
			% 	a_{\eps 26}-a_{\eps 35}\\a_{\eps 34}-a_{\eps 16}\\a_{\eps 15}-a_{\eps 24}
			% \end{pmatrix}\cdot\bar u \right)
			% dx
			% \\
			\frac{d}{dt} & \int_{\mathbb{R}^3}\mathcal{I}dx
			+
			\int_{\mathbb{R}^3}
			\left(\int_{\mathbb{R}^3\times\mathbb{R}^3\times\mathbb{S}^2} \left(\hat q_\eps^++\hat q_\eps^-\right) \tilde \phi MM_* dvdv_*d\sigma\right)
			:\nabla_x\bar u
			dx
			\\
			& - \int_{\mathbb{R}^3}
			\left(\int_{\mathbb{R}^3\times\mathbb{R}^3\times\mathbb{S}^2} \frac{\hat q_\eps^{+,-}-\hat q_\eps^{-,+}}{2} v MM_* dvdv_*d\sigma\right)
			\cdot\bar j
			dx
			\\
			& \geq
			- C\frac{\left\|\bar u\right\|_{W^{1,\infty}\left(dx\right)}+\left\|\partial_t\bar u\right\|_{L^\infty(dx)}+\left\|\bar j\right\|_{L^\infty(dx)}}{1-\left\|\bar u\right\|_{L^\infty(dx)}} \delta\mathcal{H}_\eps(t)
			\\
			& - \frac{1}{2} \int_{\mathbb{R}^3}
			\mathbf{A}_1\cdot \tilde u_\eps
			+
			\mathbf{A}_4\cdot\left(E_\eps+\bar u \wedge \left(B_\eps-\bar B\right)\right)
			+ \mathbf{A}_5\cdot\left(B_\eps+\left(E_\eps-\bar E\right)\wedge\bar u\right) dx
			\\
			&
			+\frac{\sigma}{2}
			\int_{\mathbb{R}^3}
			\mathbf{A}_3
			\cdot
			\left(\int_{\mathbb{R}^3\times\mathbb{R}^3\times\mathbb{S}^2} \frac{\hat q_\eps^{+,-}-\hat q_\eps^{-,+}}{2} v MM_* dvdv_*d\sigma\right)
			dx
			\\
			& + o(1)_{\textit{w-}L^1_\mathrm{loc}(dt)}
			+\frac{d}{dt}\left(o(1)_{L^\infty\left(dt\right)}\right)
			-\frac 12
			\int_{\mathbb{R}^3}
			\left(R_{\eps,4}+R_{\eps,5}\right)\cdot\bar j dx
			\\
			&
			+\frac \sigma 2
			\int_{\mathbb{R}^3}
			\left(\bar E + \bar u\wedge \bar B\right)
			\cdot
			P^\perp\left(\int_{\mathbb{R}^3\times\mathbb{R}^3\times\mathbb{S}^2} \frac{\hat q_\eps^{+,-}-\hat q_\eps^{-,+}}{2} v MM_* dvdv_*d\sigma\right)
			dx
			\\
			& -\frac 12 \int_{\mathbb{R}^3}
			\left(\sigma\left(\int_{\mathbb{R}^3\times\mathbb{R}^3\times\mathbb{S}^2} \frac{\hat q_\eps^{+,-}-\hat q_\eps^{-,+}}{2} v MM_* dvdv_*d\sigma\right)
			+ \tilde j_\eps
			\right)
			\cdot \left(\bar E + \bar u\wedge \bar B\right)
			dx.
		\end{aligned}
	\end{equation*}
	Thus, in view of the convergences \eqref{asymptotic sigma}, \eqref{limit j solenoidal}, the estimates \eqref{remainder4}, \eqref{remainder5} and Lemma \ref{entropy energy} (allowing to control the energy by the entropy), we infer that
	\begin{equation}\label{conservation-laws3}
		\begin{aligned}
			% \frac{d}{dt} & \int_{\mathbb{R}^3}
			% \left(\tilde u_\eps\cdot\bar u
			% +\frac 12 \left(E_\eps\cdot \bar E+B_\eps\cdot\bar B\right)
			% +\frac 12 \left(\left(E_\eps-\bar E\right)\wedge \left(B_\eps-\bar B\right)\right)\cdot \bar u
			% \vphantom{+\frac 12
			% \begin{pmatrix}
			% 	a_{\eps 26}-a_{\eps 35}\\a_{\eps 34}-a_{\eps 16}\\a_{\eps 15}-a_{\eps 24}
			% \end{pmatrix}\cdot\bar u}
			% \right.
			% \\
			% & + \left.
			% \vphantom{\tilde u_\eps\cdot\bar u
			% +\frac 12 \left(E_\eps\cdot \bar E+B_\eps\cdot\bar B\right)
			% +\frac 12 \left(\left(E_\eps-\bar E\right)\wedge \left(B_\eps-\bar B\right)\right)\cdot \bar u}
			% \frac 12
			% \begin{pmatrix}
			% 	a_{\eps 26}-a_{\eps 35}\\a_{\eps 34}-a_{\eps 16}\\a_{\eps 15}-a_{\eps 24}
			% \end{pmatrix}\cdot\bar u \right)
			% dx
			% \\
			\frac{d}{dt} & \int_{\mathbb{R}^3}\mathcal{I}dx
			+
			\int_{\mathbb{R}^3}
			\left(\int_{\mathbb{R}^3\times\mathbb{R}^3\times\mathbb{S}^2} \left(\hat q_\eps^++\hat q_\eps^-\right) \tilde \phi MM_* dvdv_*d\sigma\right)
			:\nabla_x\bar u
			dx
			\\
			& - \int_{\mathbb{R}^3}
			\left(\int_{\mathbb{R}^3\times\mathbb{R}^3\times\mathbb{S}^2} \frac{\hat q_\eps^{+,-}-\hat q_\eps^{-,+}}{2} v MM_* dvdv_*d\sigma\right)
			\cdot\bar j
			dx
			\\
			& \geq
			- C\frac{\left\|\bar u\right\|_{W^{1,\infty}\left(dx\right)}+\left\|\partial_t\bar u\right\|_{L^\infty(dx)}+\left\|\bar j\right\|_{L^\infty(dx)}}{1-\left\|\bar u\right\|_{L^\infty(dx)}} \delta\mathcal{H}_\eps(t)
			\\
			& - \frac{1}{2}\int_{\mathbb{R}^3}
			\mathbf{A}_1\cdot \tilde u_\eps
			+
			\mathbf{A}_4\cdot\left(E_\eps+\bar u \wedge \left(B_\eps-\bar B\right)\right)
			+ \mathbf{A}_5\cdot\left(B_\eps+\left(E_\eps-\bar E\right)\wedge\bar u\right) dx
			\\
			&
			+\frac{\sigma}{2}
			\int_{\mathbb{R}^3}
			\mathbf{A}_3
			\cdot
			\left(\int_{\mathbb{R}^3\times\mathbb{R}^3\times\mathbb{S}^2} \frac{\hat q_\eps^{+,-}-\hat q_\eps^{-,+}}{2} v MM_* dvdv_*d\sigma\right)
			dx
			\\
			& + o(1)_{\textit{w-}L^1_\mathrm{loc}(dt)}
			+\frac{d}{dt}\left(o(1)_{L^\infty\left(dt\right)}\right).
		\end{aligned}
	\end{equation}

	At last, we may now combine the inequalities \eqref{conservation-laws4} and \eqref{conservation-laws3} to deduce, employing the symmetries of collision integrands and \eqref{test integrand} to rewrite dissipation terms, that
	\begin{equation*}
		\begin{aligned}
			& \frac{d}{dt} \int_{\mathbb{R}^3}
			\left(
			\left(g_\eps^+ \gamma_\eps^++g_\eps^- \gamma_\eps^-\right) \chi\left(\frac{|v|^2}{K_\eps}\right)
			\bar g
			+E_\eps\cdot \bar E+B_\eps\cdot\bar B
			\vphantom{
			+ \left(\left(E_\eps-\bar E\right)\wedge \left(B_\eps-\bar B\right)
			+
			\begin{pmatrix}
				a_{\eps 26}-a_{\eps 35}\\a_{\eps 34}-a_{\eps 16}\\a_{\eps 15}-a_{\eps 24}
			\end{pmatrix}
			\right)\cdot \bar u
			}
			\right.
			\\
			& + \left.
			\vphantom{
			\left(g_\eps^+ \gamma_\eps^++g_\eps^- \gamma_\eps^-\right) \chi\left(\frac{|v|^2}{K_\eps}\right)
			\bar g
			+E_\eps\cdot \bar E+B_\eps\cdot\bar B
			}
			\left(\left(E_\eps-\bar E\right)\wedge \left(B_\eps-\bar B\right)
			+
			\begin{pmatrix}
				a_{\eps 26}-a_{\eps 35}\\a_{\eps 34}-a_{\eps 16}\\a_{\eps 15}-a_{\eps 24}
			\end{pmatrix}
			\right)\cdot \bar u
			\right)
			dx
			\\
			& +\frac 12
			\int_{\mathbb{R}^3\times\mathbb{R}^3\times\mathbb{R}^3\times\mathbb{S}^2}
			\left(\hat q_\eps^+\bar q^++\hat q_\eps^-\bar q^-
			+\hat q_\eps^{+,-}\bar q^{+,-}+\hat q_\eps^{-,+}\bar q^{-,+}
			\right)
			MM_* dxdvdv_*d\sigma
			\\
			& \geq
			- \lambda(t) \delta\mathcal{H}_\eps(t)
			- \int_{\mathbb{R}^3} \mathbf{A}\cdot
			\begin{pmatrix}
				\tilde u_\eps \\ \frac 32\tilde \theta_\eps-\tilde\rho_\eps \\
				-\frac\sigma 2\int_{\mathbb{R}^3\times\mathbb{R}^3\times\mathbb{S}^2} \left(\hat q_\eps^{+,-}-\hat q_\eps^{-,+}\right) v MM_* dvdv_*d\sigma
				\\ E_\eps+\bar u \wedge \left(B_\eps-\bar B\right) \\ B_\eps+\left(E_\eps-\bar E\right)\wedge\bar u
			\end{pmatrix}
			dx
			\\
			&
			+ o(1)_{\textit{w-}L^1_\mathrm{loc}(dt)}
			+\frac{d}{dt}\left(o(1)_{L^\infty\left(dt\right)}\right)-\frac 12 \delta\mathcal{D}_\eps(t).
		\end{aligned}
	\end{equation*}

	Next, assembling the preceding inequality with the scaled entropy inequality \eqref{ent} and the energy estimate \eqref{test energy}, we finally obtain
	\begin{equation*}
		\begin{aligned}
			& \frac{d}{dt}\delta\mathcal{H}_\eps(t)+\delta\mathcal{D}_\eps(t)
			\\
			& \leq
			\lambda(t) \delta\mathcal{H}_\eps(t)
			+\int_{\mathbb{R}^3} \mathbf{A}\cdot
			\begin{pmatrix}
				\tilde u_\eps-\bar u \\ \frac 32\tilde \theta_\eps-\tilde\rho_\eps-\frac 52\bar\theta \\
				-\frac\sigma 2\int_{\mathbb{R}^3\times\mathbb{R}^3\times\mathbb{S}^2} \left(\hat q_\eps^{+,-}-\hat q_\eps^{-,+}\right) v MM_* dvdv_*d\sigma
				-\bar j
				\\ E_\eps-\bar E +\bar u \wedge \left(B_\eps-\bar B\right) \\ B_\eps-\bar B +\left(E_\eps-\bar E\right)\wedge\bar u
			\end{pmatrix}
			dx
			\\
			&
			+ o(1)_{\textit{w-}L^1_\mathrm{loc}(dt)}
			+\frac{d}{dt}\left(o(1)_{L^\infty\left(dt\right)}\right)+\frac 12 \delta\mathcal{D}_\eps(t),
		\end{aligned}
	\end{equation*}
	which, with a straightforward application of Gr\"onwall's lemma (carefully note that this is valid even though $\delta\mathcal{H}_\eps(t)$ may be negative), concludes the proof of the proposition.
\end{proof}

\begin{rem}
	The proof of Proposition \ref{stab-prop} is based on the construction of the stability inequality \eqref{stability 4} from Proposition \ref{stability poynting 2} for the two-fluid incompressible Navier-Stokes-Maxwell system with solenoidal Ohm's law \eqref{TFINSMSO}. This approach has the great advantage of using the approximate macroscopic conservation of momentum established in Proposition \ref{approx3-prop} rather than the one from Proposition \ref{approx2-prop}.
	
	Indeed, if we were to use the latter approximate conservation law from Proposition \ref{approx2-prop}, we would have to modulate the nonlinear term $\tilde j_\eps\wedge B_\eps$ into $\left(\tilde j_\eps-\bar j\right)\wedge \left(B_\eps-\bar B\right)$ (much like in the proof of Proposition \ref{modulated energy estimate 2}~; see \eqref{bad modulation}). The term $\left(B_\eps-\bar B\right)$ would then have to be absorbed (through Gr\"onwall's lemma) by a renormalized modulated energy (or entropy), whereas $\left(\tilde j_\eps-\bar j\right)$ would need to be controlled by a renormalized modulated entropy dissipation provided $\tilde j_\eps$ is replaced by the collision integrands $-\sigma\left(\int_{\mathbb{R}^3\times\mathbb{R}^3\times\mathbb{S}^2} \frac{\hat q_\eps^{+,-}-\hat q_\eps^{-,+}}{2} v MM_* dvdv_*d\sigma\right)$. However, this last step produces remainders which may not belong to $L^2(dtdx)$ and, therefore, cannot multiply $B_\eps$. Thus, this procedure would fail.
	
	It is therefore not possible (at least, we do not know how to make it work) to establish a similar renormalized relative entropy inequality for renormalized solutions of the scaled two species Vlasov-Maxwell-Boltzmann system \eqref{VMB2} based on the construction of the stability inequality \eqref{stability 2} from Proposition \ref{modulated energy estimate 2}.
	
	Using the strategy of Proposition \ref{stability poynting 2} removes this difficulty altogether by expressing the Lorentz force $\tilde j_\eps\wedge B_\eps$ with the Poynting vector $E_\eps\wedge B_\eps$ (and some other terms). However, the drawback of this approach resides in the necessity of the restriction $\left\|\bar u\right\|_{L^\infty_{t,x}}<1$. Recall, nevertheless, that this restriction is physically relevant, since it merely entails that the modulus of the velocity $\bar u$ be less than the speed of light (see comments after the proofs of Propositions \ref{stability poynting} and \ref{stability poynting 2}).
\end{rem}

\subsection{Convergence and conclusion of proof}

We may now pass to the limit in the approximate stability inequality \eqref{stab ineq eps} and, thus, derive the crucial modulated energy inequality for the limiting system \eqref{TFINSFMSO 2}. To this end, we simply integrate \eqref{stab ineq eps} in time against non-negative test functions and then let $\eps\to 0$, which yields, in view of the well-preparedness of the initial data \eqref{well-prepared init data}, the weak convergences \eqref{weak conv fluctuations}, \eqref{weak limit observables} and the lower semi-continuities \eqref{energy liminf}, \eqref{dissipation liminf}, that
\begin{equation*}
	\begin{aligned}
		\delta\mathcal{E}(t) & + \frac 12 \int_0^t \delta\mathcal{D}(s) e^{\int_s^t\lambda(\sigma)d\sigma}ds
		\\
		& \leq \delta\mathcal{E}(0) e^{\int_0^t\lambda(s)ds}
		\\
		& +\int_0^t
		\int_{\mathbb{R}^3} \mathbf{A}\cdot
		\begin{pmatrix}
			u-\bar u \\ \frac 32 \theta-\rho-\frac 52\bar\theta \\
			-\frac\sigma 2\int_{\mathbb{R}^3\times\mathbb{R}^3\times\mathbb{S}^2} \left(q^{+,-}-q^{-,+}\right) v MM_* dvdv_*d\sigma
			-\bar j
			\\ E-\bar E +\bar u \wedge \left(B-\bar B\right) \\ B-\bar B +\left(E-\bar E\right)\wedge\bar u
		\end{pmatrix}(s)
		dx
		\\
		& \times
		e^{\int_s^t\lambda(\sigma)d\sigma}ds.
	\end{aligned}
\end{equation*}
Finally, using \eqref{constraints} and the characterization \eqref{mixed q phi psi} of the limiting collision integrands $q^{\pm,\mp}$ from Proposition \ref{weak-comp3}, we deduce that
\begin{equation*}
	\begin{aligned}
		& \delta\mathcal{E}(t) + \frac 12 \int_0^t \delta\mathcal{D}(s) e^{\int_s^t\lambda(\sigma)d\sigma}ds
		\\
		& \leq \delta\mathcal{E}(0) e^{\int_0^t\lambda(s)ds}
		+\int_0^t
		\int_{\mathbb{R}^3} \mathbf{A}\cdot
		\begin{pmatrix}
			u-\bar u \\ \frac 52\left(\theta-\bar\theta\right) \\
			j-\bar j
			\\ E-\bar E +\bar u \wedge \left(B-\bar B\right) \\ B-\bar B +\left(E-\bar E\right)\wedge\bar u
		\end{pmatrix}(s)
		dx
		e^{\int_s^t\lambda(\sigma)d\sigma}ds,
	\end{aligned}
\end{equation*}
which is precisely the stability inequality we were after.

As for the temporal continuity of $\left(u,\frac 52\theta,E,B\right)$, it is readily seen from the approximate macroscopic conservation laws from Proposition \ref{approx2-prop} and Maxwell's equations \eqref{Maxwell consistent} that $\partial_t P\tilde u_\eps$, $\partial_t \left(\frac 32\tilde\theta_\eps-\tilde\rho_\eps\right)$, $\partial_t E_\eps$ and $\partial_t B_\eps$ are uniformly bounded, in $L^1_\mathrm{loc}$ in time and in some negative index Sobolev space in $x$. It is therefore possible to show (see \cite[Appendix C]{lions7}) that $\left(P\tilde u_\eps,\frac 32\tilde\theta_\eps-\tilde\rho_\eps,E_\eps,B_\eps\right)$ converges to $\left(u,\frac 52\theta,E,B\right)\in C\left([0,\infty);\textit{w-}L^2\left(\mathbb{R}^3\right)\right)$ weakly in $L^2(dx)$ uniformly locally in time.

At last, the proof of Theorem \ref{CV-OMHD} is complete.\qed

% By an easy density argument, we then obtain the same inequality for any $( \bar u,\bar \theta)\in L^\infty_t( L^2_x) \cap L^1_t(W^{1,\infty} _x)  $ and $(\bar E, \bar B) \in L^\infty_t (L^2_x)$ such that
% $$\bar  j = \P (\bar E+\bar u\wedge \bar B) \in L^\infty_{t,x}\,.$$
%
% Together with the constraints obtained in the first section of this chapter, this implies that the weak limit $(u,\theta, j,E,B)$ is a dissipative solution to the Navier-Stokes-Maxwell system with solenoidal Ohm's law.

\section{Proof of Theorem \ref{CV-OMHDSTRONG} on strong interactions}\label{proof of theorem strong}

This demonstration closely follows the method of proof of Theorem \ref{CV-OMHD} presented in the preceding section. However, the asymptotic limit treated in Theorem \ref{CV-OMHDSTRONG} is more singular than the one from Theorem \ref{CV-OMHD}. Some steps in the coming proof will therefore require some greater care than their counterparts from the previous section.

As before, we begin our proof by appropriately gathering previous results together.

\subsection{Weak convergence of fluctuations, collision integrands and electromagnetic fields}

Thus, we are considering here a family of renormalized solutions $(f_\eps^\pm,E_\eps,B_\eps)$ to the scaled two species Vlasov-Maxwell-Boltzmann system \eqref{VMB2}, in the regime of strong interspecies interactions, i.e.\ $\delta=1$, satisfying the scaled entropy inequality \eqref{entropy2}.

According to Lemmas \ref{L1-lem} and \ref{L2-lem}, the corresponding families of fluctuations $g_\eps^\pm$ and renormalized fluctuations $\hat g_\eps^\pm$ are weakly compact in $L^1_\mathrm{loc}\left(dtdx;L^1\left(\left(1+|v|^2\right)Mdv\right)\right)$ and $L^2_\mathrm{loc}\left(dt;L^2\left(Mdxdv\right)\right)$, respectively, while, in view of Lemma \ref{L2-qlem}, the corresponding collision integrands $\hat q_\eps^\pm$ and $\hat q_\eps^{\pm,\mp}$ are weakly compact in $L^2\left(MM_*dtdxdvdv_*d\sigma\right)$. Thus, using Lemma \ref{L1-lem} again and the decomposition \eqref{fluct-decomposition}, we know that there exist $g^\pm\in L^\infty\left(dt;L^2\left(Mdxdv\right)\right)$, $(E,B)\in L^\infty\left(dt;L^2(dx)\right)$ and $q^\pm,q^{\pm,\mp}\in L^2\left(MM_*dtdxdvdv_*d\sigma\right)$, such that, up to extraction of subsequences,
\begin{equation}\label{weak conv fluctuations strong}
	\begin{aligned}
		g_\eps^\pm & \rightharpoonup g^\pm & & \text{in }L^1_\mathrm{loc}\left(dtdx;L^1\left(\left(1+|v|^2\right)Mdv\right)\right),
		\\
		\hat g_\eps^\pm & \stackrel{*}{\rightharpoonup} g^\pm & & \text{in }L^\infty\left(dt;L^2\left(Mdxdv\right)\right),
		\\
		\left(E_\eps,B_\eps\right) & \stackrel{*}{\rightharpoonup} \left(E,B\right) & & \text{in }L^\infty\left(dt;L^2\left(dx\right)\right),
		\\
		\hat q_\eps^\pm & \rightharpoonup q^\pm & & \text{in }L^2\left(MM_*dtdxdvdv_*d\sigma\right),
		\\
		\hat q_\eps^{\pm,\mp} & \rightharpoonup q^{\pm,\mp} & & \text{in }L^2\left(MM_*dtdxdvdv_*d\sigma\right),
	\end{aligned}
\end{equation}
as $\eps\to 0$. Therefore, one also has the weak convergence of the densities $\rho_\eps^\pm$, bulk velocities $u_\eps^\pm$ and temperatures $\theta_\eps^\pm$ corresponding to $g_\eps^\pm$~:
\begin{equation*}% \label{hydro variables conv strong}
	\rho_\eps^\pm \rightharpoonup \rho^\pm,
	\quad u_\eps^\pm \rightharpoonup u^\pm
	\quad \text{and}\quad
	\theta_\eps^\pm \rightharpoonup \theta^\pm
	\quad\text{in }L^1_\mathrm{loc}(dtdx) \text{ as }\eps \to 0,
\end{equation*}
where $\rho^\pm, u^\pm, \theta^\pm \in L^\infty\left(dt;L^2(dx)\right)$ are, respectively, the densities, bulk velocities and temperatures corresponding to $g^\pm$. In fact, Lemma \ref{relaxation-control} implies that
\begin{equation*}
	g^\pm=\Pi g^\pm = \rho^\pm + u^\pm\cdot v + \theta^\pm\left(\frac{|v|^2}{2}-\frac 32\right).
\end{equation*}

Next, we further introduce the scaled fluctuations
\begin{equation*}
	h_\eps = \frac{1}{\eps}\left[\left(g_\eps^+-g_\eps^-\right) - n_\eps\right],
\end{equation*}
where $n_\eps=\rho_\eps^+-\rho_\eps^-$ is the charge density, and the electrodynamic variables
\begin{equation*}
	j_\eps=\frac 1\eps\left(u_\eps^+-u_\eps^-\right),\qquad w_\eps=\frac 1\eps\left(\theta_\eps^+-\theta_\eps^-\right),
\end{equation*}
which are precisely the bulk velocity and temperature associated with the scaled fluctuations $h_\eps$. In view of Lemma \ref{bound hjw}, the electric current $j_\eps$ and the internal electric energy $w_\eps$ are uniformly bounded in $L^1_\mathrm{loc}\left(dtdx\right)$, which necessarily implies, letting $\eps\rightarrow 0$, that $u^+=u^-$ and $\theta^+=\theta^-$. Carefully note, though, that the limiting densities $\rho^+$ and $\rho^-$ may be distinct here. Therefore, we appropriately rename the limiting macroscopic variables
\begin{equation*}
	\rho=\frac{\rho^++\rho^-}{2}, \qquad n=\rho^+-\rho^-, \qquad u=u^+=u^-, \qquad \theta=\theta^+=\theta^-,
\end{equation*}
whence
\begin{equation*}
	\begin{aligned}
		g^\pm & =\rho^\pm + u\cdot v + \theta\left(\frac{|v|^2}{2}-\frac 32\right),
		\\
		\frac{g^++g^-}{2} & =\rho + u\cdot v + \theta\left(\frac{|v|^2}{2}-\frac 32\right),
		\\
		g^+-g^- & = n.
	\end{aligned}
\end{equation*}

Now, according to Lemma \ref{weak compactness h}, it holds that $\frac{h_\eps}{1+\left\|\hat g_\eps^+-\hat g_\eps^-\right\|_{L^2\left(Mdv\right)}}$ is weakly compact in $L^1_\mathrm{loc}\left(dtdx;L^1\left(\left(1+|v|\right)Mdv\right)\right)$ and that $\frac{j_\eps}{1+\left\|\hat g_\eps^+-\hat g_\eps^-\right\|_{L^2\left(Mdv\right)}}$ is weakly compact in $L^2_\mathrm{loc}(dtdx)$. Moreover, Lemma \ref{strong n} indicates that, up to extraction of subsequences, there is $h\in L^1_\mathrm{loc}\left(dtdx;L^1\left(\left(1+|v|^2\right)Mdv\right)\right)$ such that
\begin{equation*}
	\frac{h_\eps}{1+\left\|\hat g_\eps^+-\hat g_\eps^-\right\|_{L^2\left(Mdv\right)}} \rightharpoonup \frac{h}{1+\left|n\right|}
	\quad \text{in }L^1_\mathrm{loc}\left(dtdx;L^1\left(\left(1+|v|\right)Mdv\right)\right),
\end{equation*}
as $\eps\to 0$. Note, however, that $h$ is not characterized by an infinitesimal Maxwellian form. Here, we only have that
\begin{equation*}
	\Pi h=j\cdot v + w\left(\frac{|v|^2}{2}-\frac 32\right),
\end{equation*}
where the electric current $j$ and the internal electric energy $w$ are defined by
\begin{equation*}
	j=\int_{\mathbb{R}^3}hvMdv,\qquad
	w=\int_{\mathbb{R}^3}h\left(\frac{|v|^2}{3}-1\right)Mdv.
\end{equation*}
In particular, it holds that
\begin{equation}\label{convergence current strong}
	\frac{j_\eps}{1+\left\|\hat g_\eps^+-\hat g_\eps^-\right\|_{L^2\left(Mdv\right)}} \rightharpoonup \frac{j}{1+\left|n\right|}
	\quad \text{in }L^2_\mathrm{loc}\left(dtdx\right).
\end{equation}

Finally, since $n\in L^\infty\left(dt;L^2\left(dx\right)\right)$, setting
\begin{equation*}
	r_\eps=\frac{1+|n|}{1+\left\|\hat g_\eps^+-\hat g_\eps^-\right\|_{L^2\left(Mdv\right)}}\in L^\infty\left(dt;L^2_\mathrm{loc}\left(dx\right)\right),
\end{equation*}
which, according to Lemma \ref{strong n} and up to extraction of subsequences, converges almost everywhere towards the constant function $1$, we see that
\begin{equation*}
	r_\eps j_\eps \rightharpoonup j
	\quad \text{in }L^1_\mathrm{loc}\left(dtdx\right).
\end{equation*}
Similarly, since $\frac{h_\eps}{1+\left\|\hat g_\eps^+-\hat g_\eps^-\right\|_{L^2\left(Mdv\right)}}$ is bounded in $L^2_\mathrm{loc}\left(dtdx;L^1\left(\left(1+|v|\right)Mdv\right)\right)$ uniformly, by virtue of Lemma \ref{weak compactness h}, it is possible to show that
\begin{equation*}
	r_\eps h_\eps \rightharpoonup h
	\quad \text{in }L^1_\mathrm{loc}\left(dtdx;L^1\left(\left(1+|v|\right)Mdv\right)\right).
\end{equation*}

\subsection{Constraint equations, Maxwell's system and energy inequality}

In view of Proposition \ref{high weak-comp2}, we already know that the limiting thermodynamic fields $\rho$, $u$ and $\theta$ satisfy the incompressibility and Boussinesq relations
\begin{equation}\label{constraints strong}
	\Div u = 0, \qquad \rho+\theta = 0.
\end{equation}

Moreover, the discussion in Section \ref{limit maxwell} shows that the limiting electromagnetic field satisfies the Faraday equation and Gauss' laws~:
\begin{equation*}
	\begin{cases}
		\begin{aligned}
			\d_t B + \ROT E& = 0,
			\\
			\DIV E &=n,
			\\
			\DIV B &=0.
		\end{aligned}
	\end{cases}
\end{equation*}
Recall, however, that we do not know from Section \ref{limit maxwell} whether Amp\`ere's equation is necessarily satisfied in the sense of distributions in the limit.

Finally, Proposition \ref{strongOhm} further establishes that the electrodynamic variables $j$ and $w$ satisfy Ohm's law and the internal electric energy equilibrium relation
\begin{equation*}
	j-nu = \sigma\left(-\frac 12 \nabla_x n + E + u\wedge B \right),
	\qquad
	w = n\theta,
\end{equation*}
where the electric conductivity $\sigma>0$ is defined by \eqref{sigma 3}.

As for the energy bound, Proposition \ref{energy inequality strong interactions} states that, for almost every $t\geq 0$,
\begin{equation*}
	\begin{aligned}
		& \frac 12\left(\frac 12\left\|n\right\|_{L^2_x}^2+2\left\|u\right\|_{L^2_x}^2
		+ 5\left\|\theta\right\|_{L^2_x}^2 + \left\|E\right\|_{L^2_x}^2
		+ \left\|B\right\|_{L^2_x}^2 \right)(t)
		\\
		& +
		\int_0^t \left(2\mu
		\left\|\nabla_x u\right\|_{L^2_x}^2
		+ 5\kappa
		\left\|\nabla_x\theta\right\|_{L^2_x}^2 +
		\frac 1\sigma\left\|j-nu\right\|_{L^2_x}^2 \right)(s) ds
		\leq C^\mathrm{in},
	\end{aligned}
\end{equation*}
where the viscosity $\mu>0$, thermal conductivity $\kappa>0$ and electric conductivity $\sigma>0$ are respectively defined by \eqref{mu kappa 2} and \eqref{sigma 3}. In particular, it holds that
\begin{equation*}
		\begin{aligned}
			\left(n,u, \theta, E, B\right) & \in L^\infty \left( [0,\infty), dt ; L^2\left(\mathbb{R}^3, dx\right)\right), \\
			\left(u,\theta\right) & \in L^2\left([0,\infty), dt ; \dot H^1\left(\mathbb{R}^3, dx\right)\right), \\
			j-nu & \in L^2\left([0,\infty)\times\mathbb{R}^3, dtdx\right).
		\end{aligned}
\end{equation*}
This energy bound can be improved to the actual energy inequality
\begin{equation*}
	\begin{aligned}
		& \frac 12\left(\frac 12\left\|n\right\|_{L^2_x}^2+2\left\|u\right\|_{L^2_x}^2
		+ 5\left\|\theta\right\|_{L^2_x}^2 + \left\|E\right\|_{L^2_x}^2
		+ \left\|B\right\|_{L^2_x}^2 \right)(t)
		\\
		& \hspace{10mm} +
		\int_0^t \left(2\mu
		\left\|\nabla_x u\right\|_{L^2_x}^2
		+ 5\kappa
		\left\|\nabla_x\theta\right\|_{L^2_x}^2 +
		\frac 1\sigma\left\|j-nu\right\|_{L^2_x}^2 \right)(s) ds
		\\
		& \hspace{10mm} \leq
		\frac 12\left(\frac 12\left\|n^\mathrm{in}\right\|_{L^2_x}^2+2\left\|u^\mathrm{in}\right\|_{L^2_x}^2
		+ 5\left\|\theta^\mathrm{in}\right\|_{L^2_x}^2 + \left\|E^\mathrm{in}\right\|_{L^2_x}^2
		+ \left\|B^\mathrm{in}\right\|_{L^2_x}^2 \right),
	\end{aligned}
\end{equation*}
using the well-preparedness of the initial data \eqref{well-prepared init data strong}.

\subsection{The renormalized modulated entropy inequality}
% \subsection{Evolution equations and stability inequality}

We move on now to the rigorous derivation of a stability inequality encoding the asymptotic macroscopic evolution equations for $u$ and $\theta$ and the Amp\`ere equation in the spirit of the weak-strong stability inequalities used in Section \ref{laure-diogo} to define dissipative solutions for some Navier-Stokes-Maxwell systems. Recall that, as explained therein, such systems are in general not known to display weak stability so that their weak solutions in the energy space are not known to exist.

The strategy used here closely follows the method employed in the case of weak interactions detailed in Section \ref{stability weak}.

Thus, as in Section \ref{conservation defects 2 species}, we define the renormalized fluctuations $g_\eps^\pm \gamma_\eps^\pm\chi\left(\frac{|v|^2}{K_\eps}\right)$, with $K_\eps =K|\log \eps|$, for some large $K>0$, and $\chi\in C_c^\infty\left([0,\infty)\right)$ a smooth compactly supported function such that $\mathds{1}_{[0,1]}\leq \chi \leq \mathds{1}_{[0,2]}$, and where $\gamma_\eps^\pm=\gamma\left(G_\eps^\pm\right)$ for some renormalization $\gamma\in C^1\left([0,\infty);\mathbb{R}\right)$ satisfying \eqref{gamma-assumption}.

Since, up to further extraction of subsequences, $\gamma_\eps^\pm \chi\left( {|v|^2\over K_\eps}\right)$ converges almost everywhere towards $1$, $g_\eps^\pm$ is weakly compact in $L^1_\mathrm{loc}\left(dtdx;L^1\left(\left(1+|v|^2\right)Mdv\right)\right)$ and $g_\eps^\pm\gamma_\eps^\pm$ is uniformly bounded in $L^\infty\left(dt;L^2\left(Mdxdv\right)\right)$, we deduce, by the Product Limit Theorem, that
\begin{equation*}
	g_\eps^\pm \gamma_\eps^\pm \chi\left( {|v|^2\over K_\eps}\right)
	\stackrel{*}{\rightharpoonup} g^\pm \quad \text{in }L^\infty\left(dt;L^2\left(Mdxdv\right)\right).
\end{equation*}
Therefore, one has the weak convergence of the densities $\tilde\rho_\eps^\pm$, bulk velocities $\tilde u_\eps^\pm$ and temperatures $\tilde\theta_\eps^\pm$ corresponding to $g_\eps^\pm\gamma_\eps^\pm \chi\left( {|v|^2\over K_\eps}\right)$~:
\begin{equation*}
	\tilde\rho_\eps^\pm \stackrel{*}{\rightharpoonup} \rho^\pm,
	\quad \tilde u_\eps^\pm \stackrel{*}{\rightharpoonup} u
	\quad \text{and}\quad
	\tilde\theta_\eps^\pm \stackrel{*}{\rightharpoonup} \theta
	\quad\text{in }L^\infty\left(dt;L^2(dx)\right) \text{ as }\eps \to 0.
\end{equation*}
In particular, the hydrodynamic variables $\tilde\rho_\eps=\frac{\tilde\rho_\eps^++\tilde\rho_\eps^-}{2}$, $\tilde u_\eps=\frac{\tilde u_\eps^++\tilde u_\eps^-}{2}$ and $\tilde\theta_\eps=\frac{\tilde\theta_\eps^++\tilde\theta_\eps^-}{2}$ also obviously verify
\begin{equation}\label{weak limit observables strong}
	\tilde\rho_\eps \stackrel{*}{\rightharpoonup} \rho,
	\quad \tilde u_\eps \stackrel{*}{\rightharpoonup} u
	\quad \text{and}\quad
	\tilde\theta_\eps \stackrel{*}{\rightharpoonup} \theta
	\quad\text{in }L^\infty\left(dt;L^2(dx)\right) \text{ as }\eps \to 0,
\end{equation}
while the charge density $\tilde n_\eps=\tilde\rho^+_\eps-\tilde\rho_\eps^-$ satisfies
\begin{equation}\label{weak limit charges strong}
	\tilde n_\eps \stackrel{*}{\rightharpoonup} n
	\quad\text{in }L^\infty\left(dt;L^2(dx)\right) \text{ as }\eps \to 0.
\end{equation}
It follows that, since $u$ is solenoidal,
\begin{equation}\label{limit solenoidal strong}
	P^\perp \tilde u_\eps \stackrel{*}{\rightharpoonup} 0
	\quad\text{in }L^\infty\left(dt;L^2(dx)\right) \text{ as }\eps \to 0,
\end{equation}
and, in view of the limiting Boussinesq relation,
\begin{equation}\label{limit boussinesq strong}
	\tilde \rho_\eps+\tilde\theta_\eps \stackrel{*}{\rightharpoonup} 0
	\quad\text{in }L^\infty\left(dt;L^2(dx)\right) \text{ as }\eps \to 0.
\end{equation}

Here, in constrast with the convergence properties of the electric current established in Section \ref{stability weak} for weak interactions, we cannot show the convergence of the electric current $\tilde j_\eps =\frac 1\eps\left(\tilde u_\eps^+-\tilde u_\eps^-\right)$ towards $j$ unless we renormalize it as in \eqref{convergence current strong}. Instead, we establish below in \eqref{tilde u approx u strong} a useful consistency relation for $\tilde j_\eps$ by suitably controlling remainders in the spirit of Section \ref{conservation defects 2 species}.

\bigskip

Now, just as in the case of weak interspecies interactions, the $L^2\left(Mdxdv\right)$ norm of $g_\eps^\pm \gamma_\eps^\pm\chi\left(\frac{|v|^2}{K_\eps}\right)$ is not a Lyapunov functional but it is nevertheless controlled by the relative entropy
\begin{equation}\label{renormalized energy entropy control strong}
	\frac 12 \left\|g_\eps^\pm \gamma_\eps^\pm\chi\left(\frac{|v|^2}{K_\eps}\right)\right\|_{L^2\left(Mdxdv\right)}^2
	\leq \frac{C}{\eps^2}H\left(f_\eps^\pm\right),
\end{equation}
for some $C>1$, and therefore by the initial data \eqref{init-fluctuation 2 strong}. One may therefore try, in a preliminary attempt to show an asymptotic stability inequality, to modulate the approximate energy associated with $g_\eps^\pm \gamma_\eps^\pm\chi\left(\frac{|v|^2}{K_\eps}\right)$, i.e.\ its $L^2\left(Mdxdv\right)$ norm, by introducing a test functions $\bar g^\pm$ in infinitesimal Maxwellian form~:
\begin{equation*}
	\bar g^\pm= \bar\rho^\pm + \bar u\cdot v + \bar \theta \left(\frac{|v|^2}{2}-\frac 32\right),
\end{equation*}
where
\begin{equation*}
	\bar\rho^\pm(t,x),\bar u(t,x),\bar \theta(t,x)\in C^\infty_c\left([0,\infty)\times\mathbb{R}^3\right)
	\quad\text{with } \Div\bar u = 0,\ \frac{\bar\rho^++\bar\rho^-}{2}+\bar\theta = 0,
\end{equation*}
and then establishing a stability inequality for the modulated energies
\begin{equation}\label{modulated energy strong}
	\frac 12 \left\|g_\eps^\pm \gamma_\eps^\pm\chi\left(\frac{|v|^2}{K_\eps}\right)-\bar g^\pm\right\|_{L^2\left(Mdxdv\right)}^2.
\end{equation}
Notice that it holds, utilizing the elementary identity $a^2+\frac 32 b^2=\frac 35\left(a+b\right)^2+\frac 52\left(\frac{3b-2a}{5}\right)^2$, for any $a,b\in\mathbb{R}$,
\begin{equation*}
	\begin{aligned}
		\sum_{\pm}& \left\|g_\eps^\pm \gamma_\eps^\pm\chi\left(\frac{|v|^2}{K_\eps}\right)-\bar g^\pm\right\|_{L^2\left(Mdxdv\right)}^2
		\\ & \geq
		\sum_{\pm} \left\|\Pi\left(g_\eps^\pm \gamma_\eps^\pm\chi\left(\frac{|v|^2}{K_\eps}\right)\right)-\bar g^\pm\right\|_{L^2\left(Mdxdv\right)}^2
		\\ & =
		\sum_{\pm}\left(\left\|\tilde\rho_\eps^\pm-\bar\rho^\pm\right\|_{L^2\left(dx\right)}^2
		+\left\|\tilde u_\eps^\pm-\bar u\right\|_{L^2\left(dx\right)}^2
		+\frac 32 \left\|\tilde\theta_\eps^\pm-\bar \theta\right\|_{L^2\left(dx\right)}^2\right)
		\\ & =
		2\left\|\tilde\rho_\eps-\bar\rho\right\|_{L^2\left(dx\right)}^2
		+\frac 12 \left\|\tilde n_\eps-\bar n\right\|_{L^2\left(dx\right)}^2
		\\
		& +\sum_{\pm}\left(
		\left\|\tilde u_\eps^\pm-\bar u\right\|_{L^2\left(dx\right)}^2
		+\frac 32 \left\|\tilde\theta_\eps^\pm-\bar \theta\right\|_{L^2\left(dx\right)}^2\right)
		\\ & =
		\frac 12 \left\|\tilde n_\eps-\bar n\right\|_{L^2\left(dx\right)}^2
		\\
		& +\sum_{\pm}\left(
		\frac 35 \left\|\tilde\rho_\eps+\tilde\theta_\eps^\pm\right\|_{L^2\left(dx\right)}^2
		+
		\left\|\tilde u_\eps^\pm-\bar u\right\|_{L^2\left(dx\right)}^2
		+\frac 52 \left\|\frac{3\tilde\theta_\eps^\pm-2\tilde\rho_\eps}{5}-\bar\theta\right\|_{L^2\left(dx\right)}^2\right),
	\end{aligned}
\end{equation*}
where we have denoted $\bar \rho=\frac{\bar\rho^++\bar\rho^-}{2}$ and $\bar n = \bar\rho^+-\bar\rho^-$.

As before, it turns out that this approach is not quite suitable for our purpose because, even though, for any $0\leq t_1<t_2$ (see the proof of Lemma \ref{L1-lem}),
\begin{equation}\label{limit entropy strong}
		\int_{t_1}^{t_2}
		\frac 1 {2}
		\left\|g^\pm\right\|_{L^2\left(Mdxdv\right)}^2 dt
		\leq
		\liminf_{\eps\rightarrow 0}\int_{t_1}^{t_2}\frac 1{\eps^2}H\left(f_\eps^\pm\right)dt,
\end{equation}
it is not possible to set $C=1$ in \eqref{renormalized energy entropy control strong}. Indeed, the first term in the polynomial expansion of the function $h(z)=(1+z)\log(1+z)-z$ defining the entropy is $\frac 12 z^2$, but the second term is $-\frac 16 z^3$ and may be negative.

Some entropy (or energy) is therefore lost by considering the modulated energies \eqref{modulated energy strong}. These considerations lead us to introduce a more precise modulated functional in replacement of \eqref{modulated energy strong} capturing more information on the fluctuations. To be precise, instead of \eqref{modulated energy strong}, we consider now the {\bf renormalized modulated entropies}
\begin{equation}\label{modulated entropy strong}
	\frac 1{\eps^2}H\left(f_\eps^\pm\right)
	-
	\int_{\mathbb{R}^3\times\mathbb{R}^3}
	g_\eps^\pm \gamma_\eps^\pm\chi\left(\frac{|v|^2}{K_\eps}\right)\bar g^\pm
	Mdxdv
	+
	\frac 1 {2}
	\left\|\bar g^\pm\right\|_{L^2\left(Mdxdv\right)}^2.
\end{equation}
Note that the above functional may be negative for fixed $\eps>0$. However, in view of \eqref{limit entropy strong}, it recovers asymptotically a non-negative quantity, which is precisely the asymptotic modulated energy~:
\begin{equation}\label{modulated entropy liminf strong}
	\begin{aligned}
		\int_{t_1}^{t_2}
		\frac 12 & \left(\left\|\rho^\pm-\bar \rho^\pm\right\|_{L^2(dx)}^2+\left\|u-\bar u\right\|_{L^2(dx)}^2+\frac 32\left\|\theta-\bar \theta\right\|_{L^2(dx)}^2\right)
		dt
		\\
		& = \int_{t_1}^{t_2}
		\frac 1 {2}
		\left\|g^\pm - \bar g^\pm\right\|_{L^2\left(Mdxdv\right)}^2 dt
		\\
		& \leq \liminf_{\eps\rightarrow 0}
		\int_{t_1}^{t_2}\bigg(\frac 1{\eps^2}H\left(f_\eps^\pm\right)
		\\
		& -
		\int_{\mathbb{R}^3\times\mathbb{R}^3}
		g_\eps^\pm \gamma_\eps^\pm\chi\left(\frac{|v|^2}{K_\eps}\right)\bar g^\pm
		Mdxdv
		+
		\frac 1 {2}
		\left\|\bar g^\pm\right\|_{L^2\left(Mdxdv\right)}^2\bigg)dt,
	\end{aligned}
\end{equation}
for all $0\leq t_1<t_2$.

The first term in \eqref{modulated entropy strong} is precisely the entropy of $f_\eps^\pm$ and will be controlled by the scaled entropy inequality \eqref{entropy2}, whereas the last term in \eqref{modulated entropy strong} only involves smooth quantities and will therefore be controlled directly. As for the middle term in the modulated entropy \eqref{modulated entropy strong}, its time derivative will involve the {\bf approximate macroscopic conservation laws} for $ g_\eps^\pm \gamma_\eps ^\pm \chi\left(\frac{|v|^2}{K_\eps}\right) $.

\bigskip

Now, in order to establish the renormalized modulated entropy inequality leading to the convergence stated in Theorem \ref{CV-OMHDSTRONG}, we introduce further test functions
\begin{equation*}
	\bar E(t,x),\bar B(t,x), \bar j(t,x)\in C^\infty_c\left([0,\infty)\times\mathbb{R}^3\right)
	\quad\text{with }\Div\bar E = \bar n,\ \Div\bar B = 0,
\end{equation*}
and we define the {\bf renormalized modulated entropy}
\begin{equation*}
	\begin{aligned}
		\delta\mathcal{H}_\eps(t)
		& =
		\frac 1{\eps^2}H\left(f_\eps^+\right) + \frac 1{\eps^2}H\left(f_\eps^-\right)
		\\
		& -
		\int_{\mathbb{R}^3\times\mathbb{R}^3}
		\left(g_\eps^+ \gamma_\eps^+\chi\left(\frac{|v|^2}{K_\eps}\right)\bar g^++g_\eps^- \gamma_\eps^-\chi\left(\frac{|v|^2}{K_\eps}\right)\bar g^-\right)
		Mdxdv
		\\
		& +
		\frac 12\left\|\bar g^+\right\|_{L^2\left(Mdxdv\right)}^2 + \frac 12\left\|\bar g^-\right\|_{L^2\left(Mdxdv\right)}^2
		\\
		& +
		\frac 12\left\|E_\eps - \bar E\right\|_{L^2(dx)}^2
		+
		\frac 12\left\|B_\eps - \bar B\right\|_{L^2(dx)}^2
		+ \frac 1{2}\int_{\mathbb{R}^3}\left(\frac 1{\eps^2}\operatorname{Tr}m_\eps+\operatorname{Tr}a_\eps\right) dx
		\\
		& - \int_{\mathbb{R}^3}
		\left(\left(E_\eps-\bar E\right)\wedge\left(B_\eps-\bar B\right)
		+
		\begin{pmatrix}
			a_{\eps 26}-a_{\eps 35}\\a_{\eps 34}-a_{\eps 16}\\a_{\eps 15}-a_{\eps 24}
		\end{pmatrix}
		\right)\cdot\bar u
		dx,
	\end{aligned}
\end{equation*}
where the matrix measures $m_\eps$ and $a_\eps$ are the defects introduced in Section \ref{macroscopic defects} and controlled by the scaled entropy inequality \eqref{entropy2}.

We also define the {\bf renormalized modulated energy}
\begin{equation*}
	\begin{aligned}
		\delta\mathcal{E}_\eps(t)
		& =
		\frac 12\left\|g_\eps^+ \gamma_\eps^+ \chi\left(\frac{|v|^2}{K_\eps}\right)-\bar g^+\right\|_{L^2(Mdxdv)}^2
		+
		\frac 12\left\|g_\eps^- \gamma_\eps^- \chi\left(\frac{|v|^2}{K_\eps}\right)-\bar g^-\right\|_{L^2(Mdxdv)}^2
		\\
		& +
		\frac 12\left\|E_\eps - \bar E\right\|_{L^2(dx)}^2
		+
		\frac 12\left\|B_\eps - \bar B\right\|_{L^2(dx)}^2
		+ \frac 1{2}\int_{\mathbb{R}^3}\left(\frac 1{\eps^2}\operatorname{Tr}m_\eps+\operatorname{Tr}a_\eps\right) dx
		\\
		& - \int_{\mathbb{R}^3}
		\left(\left(E_\eps-\bar E\right)\wedge\left(B_\eps-\bar B\right)
		+
		\begin{pmatrix}
			a_{\eps 26}-a_{\eps 35}\\a_{\eps 34}-a_{\eps 16}\\a_{\eps 15}-a_{\eps 24}
		\end{pmatrix}
		\right)\cdot\bar u
		dx,
	\end{aligned}
\end{equation*}
which is asymptotically equivalent to $\delta\mathcal{H}_\eps(t)$, at least formally. Note that $\delta\mathcal{H}_\eps(t)$ controls more accurately the large values of the fluctuations $g_\eps^\pm$ than $\delta\mathcal{E}_\eps(t)$. Lemma \ref{entropy energy strong} below shows how the modulated entropy $\delta\mathcal{H}_\eps(t)$ controls the modulated energy $\delta\mathcal{E}_\eps(t)$.

Finally, we introduce the {\bf renormalized modulated entropy dissipation}
\begin{equation*}
	\begin{aligned}
		\delta\mathcal{D}_\eps(t) & =
		\frac 14
		\left\|\hat q^+_\eps - \bar q^+\right\|_{L^2\left(MM_*dxdvdv_*d\sigma\right)}^2
		+
		\frac 14
		\left\|\hat q^-_\eps - \bar q^-\right\|_{L^2\left(MM_*dxdvdv_*d\sigma\right)}^2
		\\
		& + \frac 14
		\left\|\hat q^{+,-}_\eps - \bar q^{+,-}\right\|_{L^2\left(MM_*dxdvdv_*d\sigma\right)}^2
		+
		\frac 14
		\left\|\hat q^{-,+}_\eps - \bar q^{-,+}\right\|_{L^2\left(MM_*dxdvdv_*d\sigma\right)}^2,
	\end{aligned}
\end{equation*}
where
\begin{equation}\label{test integrand strong}
	\begin{aligned}
		\bar q^\pm
		& =
		\frac 12\nabla_x \bar u:\left(\tilde \phi + \tilde \phi_*
		-\tilde\phi'-\tilde\phi_*'\right)
		+\frac 12 \nabla_x \bar \theta\cdot\left(\tilde \psi + \tilde \psi_*
		-\tilde\psi'-\tilde\psi_*'\right)
		\\
		& \mp\frac 1{\sigma}
		\left(\bar j-\bar n\bar u\right)\cdot\left(\tilde\Phi+\tilde\Phi_*-\tilde\Phi'-\tilde\Phi_*'\right),
		\\
		\bar q^{\pm,\mp}
		& =
		\frac 12\nabla_x \bar u:\left(\tilde \phi + \tilde \phi_*
		-\tilde\phi'-\tilde\phi_*'\right)
		+\frac 12 \nabla_x \bar \theta\cdot\left(\tilde \psi + \tilde \psi_*
		-\tilde\psi'-\tilde\psi_*'\right)
		\\
		& \mp\frac 1{\sigma}
		\left(\bar j-\bar n\bar u\right)\cdot\left(\tilde\Phi-\tilde\Phi_*-\tilde\Phi'+\tilde\Phi_*'\right),
	\end{aligned}
\end{equation}
so that
\begin{equation*}
	\begin{aligned}
		\frac 12\int_{\mathbb{R}^3\times\mathbb{S}^2} \left(\bar q^++\bar q^-+\bar q^{+,-}+\bar q^{-,+}\right) M_*dv_*d\sigma
		& = \nabla_x \bar u:\mathcal{L}\tilde\phi + \nabla_x \bar \theta\cdot \mathcal{L}\tilde\psi
		\\
		& = \nabla_x \bar u:\phi + \nabla_x \bar \theta\cdot \psi,
		\\
		\int_{\mathbb{R}^3\times\mathbb{S}^2} \left(\bar q^+-\bar q^-\right) M_*dv_*d\sigma
		& = -\frac{2}{\sigma}\left(\bar j - \bar n\bar u\right)\cdot\mathcal{L}\left(\tilde\Phi\right),
		\\
		\int_{\mathbb{R}^3\times\mathbb{S}^2} \left(\bar q^{+,-}-\bar q^{-,+}\right) M_*dv_*d\sigma
		& = -\frac{2}{\sigma}\left(\bar j - \bar n\bar u\right)\cdot\mathfrak{L}\left(\tilde\Phi\right),
		\\
		\bar q^{+}+\bar q^{-}-\bar q^{+,-}-\bar q^{-,+} & =0,
	\end{aligned}
\end{equation*}
with $\phi$, $\psi$, $\tilde\phi$ and $\tilde \psi$ defined by \eqref{phi-psi-def} and \eqref{phi-psi-def inverses} and $\tilde\Phi$ defined by \eqref{phi-psi-def inverses two species}.

Then, assuming from now on that $\left\|\bar u\right\|_{L^\infty(dtdx)}<1$ and using the lower weak sequential semi-continuity of the entropies \eqref{modulated entropy liminf strong} and of the electromagnetic energy \eqref{convex wedge} together with Lemma \ref{vector defect}, we find that, for all $0\leq t_1<t_2$,
\begin{equation}\label{energy liminf strong}
	\int_{t_1}^{t_2}
	\delta\mathcal{E}(t) dt
	\leq \liminf_{\eps\rightarrow 0}
	\min\left\{
	\int_{t_1}^{t_2}\delta\mathcal{H}_\eps(t) dt,
	\int_{t_1}^{t_2}\delta\mathcal{E}_\eps(t) dt
	\right\},
\end{equation}
where
\begin{equation*}
	\begin{aligned}
		\delta\mathcal{E}(t)
		& =
		\frac 12\left\|g^+ - \bar g^+\right\|_{L^2\left(Mdxdv\right)}^2+\frac 12\left\|g^- - \bar g^-\right\|_{L^2\left(Mdxdv\right)}^2
		\\
		& +\frac 12\left\|E-\bar E\right\|_{L^2(dx)}^2
		+\frac 12\left\|B-\bar B\right\|_{L^2(dx)}^2
		- \int_{\mathbb{R}^3}
		\left(\left(E_\eps-\bar E\right)\wedge\left(B_\eps-\bar B\right)
		\right)\cdot\bar u
		dx
		\\
		& =
		\left\|\rho-\bar \rho\right\|_{L^2(dx)}^2+\frac 14\left\|n-\bar n\right\|_{L^2(dx)}^2
		+ \left\|u-\bar u\right\|_{L^2(dx)}^2+\frac 32\left\|\theta-\bar \theta\right\|_{L^2(dx)}^2
		\\
		& +\frac 12\left\|E-\bar E\right\|_{L^2(dx)}^2
		+\frac 12\left\|B-\bar B\right\|_{L^2(dx)}^2
		- \int_{\mathbb{R}^3}
		\left(\left(E_\eps-\bar E\right)\wedge\left(B_\eps-\bar B\right)
		\right)\cdot\bar u
		dx
		\\
		& =
		\frac 14\left\|n-\bar n\right\|_{L^2(dx)}^2
		+ \left\|u-\bar u\right\|_{L^2(dx)}^2+\frac 52\left\|\theta-\bar \theta\right\|_{L^2(dx)}^2
		\\
		& +\frac 12\left\|E-\bar E\right\|_{L^2(dx)}^2
		+\frac 12\left\|B-\bar B\right\|_{L^2(dx)}^2
		- \int_{\mathbb{R}^3}
		\left(\left(E_\eps-\bar E\right)\wedge\left(B_\eps-\bar B\right)
		\right)\cdot\bar u
		dx,
	\end{aligned}
\end{equation*}
while, repeating mutatis mutandis the computations leading to \eqref{entropy3 1 limit} and \eqref{entropy3 1 limit 2} in the proof of Proposition \ref{energy inequality strong interactions}, we obtain, for all $0\leq t_1<t_2$,
\begin{equation}\label{dissipation liminf strong}
	\int_{t_1}^{t_2}
	\delta\mathcal{D}(t)dt
	\leq \liminf_{\eps\rightarrow 0}
	\int_{t_1}^{t_2}\delta\mathcal{D}_\eps(t) dt,
\end{equation}
where
\begin{equation*}
	\begin{aligned}
		\delta\mathcal{D}(t)
		& =
		2\mu
		\left\|\nabla_x \left(u-\bar u\right)\right\|_{L^2_x}^2
		+ 5\kappa
		\left\|\nabla_x\left(\theta-\bar\theta\right)\right\|_{L^2_x}^2
		+ \frac 1\sigma
		\left\|\left(j-nu\right)-\left(\bar j-\bar n\bar u\right)\right\|_{L^2_x}^2
		\\
		& \leq
		\frac 14
		\left\|q^+ - \bar q^+\right\|_{L^2\left(MM_*dxdvdv_*d\sigma\right)}^2
		+
		\frac 14
		\left\|q^- - \bar q^-\right\|_{L^2\left(MM_*dxdvdv_*d\sigma\right)}^2
		\\
		& + \frac 14
		\left\|q^{+,-} - \bar q^{+,-}\right\|_{L^2\left(MM_*dxdvdv_*d\sigma\right)}^2
		+
		\frac 14
		\left\|q^{-,+} - \bar q^{-,+}\right\|_{L^2\left(MM_*dxdvdv_*d\sigma\right)}^2.
	\end{aligned}
\end{equation*}

The following lemma shows how the modulated entropy $\delta\mathcal{H}_\eps$ controls the modulated energy $\delta\mathcal{E}_\eps$ up to a small remainder. It is obtained by repeating the proof of Lemma \ref{entropy energy} and, thus, we skip the details of its demonstration.

\begin{lem}\label{entropy energy strong}
	It holds that
	\begin{equation*}
		\delta\mathcal{E}_\eps(t)
		\leq
		C \delta\mathcal{H}_\eps(t)
		+o(1)_{L^\infty(dt)},
	\end{equation*}
	for some fixed constant $C>1$.
\end{lem}

The following result establishes the renormalized modulated entropy inequality at the order $\eps$, which will eventually allow us to deduce the crucial weak-strong stability of the limiting thermodynamic fields, thus defining dissipative solutions.

\begin{prop}\label{stab-prop strong}
	One has the stability inequality
	\begin{equation}\label{stab ineq eps strong}
		\begin{aligned}
			& \delta\mathcal{H}_\eps(t) + \frac 12 \int_0^t \delta\mathcal{D}_\eps(s) e^{\int_s^t\lambda(\sigma)d\sigma}ds
			\\
			& \leq \delta\mathcal{H}_\eps(0) e^{\int_0^t\lambda(s)ds}
			\\
			& +\int_0^t
			\int_{\mathbb{R}^3} \mathbf{A}\cdot
			\begin{pmatrix}
				\tilde u_\eps-\bar u \\ \frac 32\tilde \theta_\eps-\tilde\rho_\eps-\frac 52\bar\theta \\
				- \int_{\mathbb{R}^3\times\mathbb{R}^3\times\mathbb{S}^2}
				\sum_{\pm}\left(\pm\hat q_\eps^\pm\pm\hat q_\eps^{\pm,\mp}\right)
				% \left(\hat q_\eps^+ - \hat q_\eps^- + \hat q_\eps^{+,-}-\hat q_\eps^{-,+}\right)
				\tilde\Phi MM_* dvdv_*d\sigma
				-\left(\bar j-\bar n\bar u\right)
				\\ E_\eps-\bar E +\bar u \wedge \left(B_\eps-\bar B\right) -\frac 12 \nabla_x\left(\tilde n_\eps-\bar n\right) \\ B_\eps-\bar B +\left(E_\eps-\bar E\right)\wedge\bar u
			\end{pmatrix}(s)
			dx
			\\
			& \times
			e^{\int_s^t\lambda(\sigma)d\sigma}ds
			+ o(1)_{L^\infty_\mathrm{loc}(dt)},
		\end{aligned}
	\end{equation}
	where the {\bf acceleration operator} is defined by
	\begin{equation*}
		\begin{aligned}
			\mathbf{A} \left( \bar u, \bar \theta, \bar n, \bar j, \bar E, \bar B\right)
			& =
			\begin{pmatrix}
				\mathbf{A}_1
				\\
				\mathbf{A}_2
				\\
				\mathbf{A}_3
				\\
				\mathbf{A}_4
				\\
				\mathbf{A}_5
			\end{pmatrix}
			\\
			& =
			\begin{pmatrix}
				-2\left(\d_t \bar u +
				P\left(\bar u\cdot\nabla_x \bar u\right) - \mu\Delta_x \bar u\right)
				+ P \left(\bar n\bar E+\bar j \wedge \bar B\right)
				\\
				-2\left(\partial_t\bar\theta + \bar u \cdot\nabla_x\bar \theta - \kappa\Delta_x\bar \theta\right)
				\\
				- \frac 1{\sigma}\left(\bar j-\bar n\bar u\right) - \frac 12\nabla_x\bar n + \bar E + \bar u\wedge \bar B
				\\
				-\left(\partial_t\bar E - \rot\bar B + \bar j\right)
				\\
				-\left(\partial_t\bar B + \rot\bar E\right)
			\end{pmatrix},
		\end{aligned}
	\end{equation*}
	and the {\bf growth rate} is given by
	\begin{equation*}
		\begin{aligned}
			& \lambda(t) =
			\\
			& C\Bigg(\frac{\left\|\bar u(t)\right\|_{W^{1,\infty}\left(dx\right)}+\left\|\partial_t\bar u(t)\right\|_{L^\infty(dx)}
			+\left\|\bar\theta(t)\right\|_{W^{1,\infty}(dx)}
			+\left\|\left(\bar j-\bar n\bar u\right)(t)\right\|_{L^\infty(dx)}}
			{1-\left\|\bar u(t)\right\|_{L^\infty(dx)}}
			\\
			& \hspace{35mm} +\left\|\bar\theta(t)\right\|_{W^{1,\infty}(dx)}^2
			+\left\|\left(\frac 12 \nabla_x\bar n - \bar E - \bar u\wedge\bar B\right)(t)\right\|_{L^\infty(dx)} \Bigg),
		\end{aligned}
	\end{equation*}
	with a constant $C>0$ independent of test functions and $\eps$.
\end{prop}

\begin{proof}
	The main ingredients of the proof of this stability inequality are~:
	\begin{itemize}
		
		\item
		The scaled entropy inequality \eqref{entropy2}
		\begin{equation}\label{ent strong}
			\begin{aligned}
				\frac1{\eps^2} & H\left(f_\eps^{+}\right)
				+ \frac1{\eps^2} H\left(f_\eps^{-}\right)
				+ \frac 1{2\eps^2}\int_{\mathbb{R}^3}\operatorname{Tr}m_\eps dx
				+ \frac 1{2} \int_{\mathbb{R}^3} \left(|E_\eps|^2+ |B_\eps|^2 + \operatorname{Tr}a_\eps\right) dx
				\\
				& +\frac{1}{4}\int_0^t\int_{\mathbb{R}^3\times\mathbb{R}^3\times\mathbb{R}^3\times\mathbb{S}^2}
				\left(
				\left(\hat q_\eps^{+}\right)^2
				+ \left(\hat q_\eps^{-}\right)^2
				+ \left(\hat q_\eps^{+,-}\right)^2
				+ \left(\hat q_\eps^{-,+}\right)^2\right)
				MM_*dxdvdv_*d\sigma ds
				% \\
				% & +\frac{1}{\epsilon^4}\int_0^t\int_{\mathbb{R}^3}\left(D\left(f_\eps^+\right)+D\left(f_\eps^-\right)
				% + \delta^2 D\left(f_\eps^+,f_\eps^-\right)\right)(s) dx ds
				\\
				& \leq
				\frac1{\eps^2} H\left(f_\eps^{+\mathrm{in}}\right)
				+ \frac1{\eps^2} H\left(f_\eps^{-\mathrm{in}}\right)
				+ \frac1{2}\int_{\mathbb{R}^3} \left(|E_\eps^{\rm in}|^2+ |B_\eps^{\rm in}|^2\right) dx,
			\end{aligned}
		\end{equation}
		which is naturally satisfied by renormalized solutions of the scaled two species Vlasov-Maxwell-Boltzmann system \eqref{VMB2} (provided they exist) and where we have used the inequality \eqref{q-est} from Lemma \ref{L2-qlem} in order to conveniently simplify the dissipation terms.
		
		\item
		The approximate conservation of energy obtained in Proposition \ref{approx2-prop}
		\begin{equation}\label{conservation-laws strong}
			\d_t \left(\frac 32\tilde \theta_\eps-\tilde \rho_\eps\right) + \nabla_x \cdot \left( \frac52 \tilde u_\eps \tilde \theta_\eps
			- \int_{\mathbb{R}^3\times\mathbb{R}^3\times\mathbb{S}^2} \frac{\hat q_\eps^++\hat q_\eps^-}{2} \tilde \psi MM_* dvdv_*d\sigma\right)
		\end{equation}
		where the remainder $R_{\eps,1}$ satisfies
		\begin{equation}\label{remainder1 strong}
			\left\|R_{\eps,1}\right\|_{W^{-1,1}_\mathrm{loc}\left(dx\right)}
			\leq \frac{C\delta \mathcal{E}_\eps(t)}{1-\left\|\bar u\right\|_{L^\infty\left(dx\right)}}
			+ C\left(\delta\mathcal{E}_\eps(t)\delta\mathcal{D}_\eps(t)\right)^\frac{1}{2} + o(1)_{L^1_\mathrm{loc}(dt)},
		\end{equation}
		for some $C>0$, where we have used \eqref{test u estimate}. Note that we do not employ the approximate conservation of momentum from Proposition \ref{approx2-prop}.
		
		\item The approximate conservation of momentum law from Proposition \ref{approx3-prop}
		\begin{equation}\label{conservation-laws2 strong}
			\begin{aligned}
				\partial_t & \left(\tilde u_\eps
				+\frac 12 E_\eps\wedge B_\eps
				+\frac 12
				\begin{pmatrix}
					a_{\eps 26}-a_{\eps 35}\\a_{\eps 34}-a_{\eps 16}\\a_{\eps 15}-a_{\eps 24}
				\end{pmatrix}
				\right)
				\\
				& +
				\nabla_x \cdot \left( \tilde u_\eps \otimes \tilde u_\eps
				-\frac{\left|\tilde u_\eps\right|^2}{3} \operatorname{Id}
				+ \frac 1{2\eps^2}m_\eps
				-
				\int_{\mathbb{R}^3\times\mathbb{R}^3\times\mathbb{S}^2} \frac{\hat q_\eps^++\hat q_\eps^-}{2} \tilde \phi MM_* dvdv_*d\sigma
				\right)
				\\
				& -\frac 12 \nabla_x\cdot\left(
				E_\eps\otimes E_\eps + e_\eps + B_\eps\otimes B_\eps + b_\eps
				\right)
				+\nabla_x\left(\frac{|E_\eps|^2+|B_\eps|^2+\operatorname{Tr}a_\eps}{4}\right)
				\\ & = - \frac 1\eps \nabla_x\left(\tilde \rho_\eps+\tilde \theta_\eps\right)+\partial_t \left(R_{\eps,2}\right) + R_{\eps,3},
			\end{aligned}
		\end{equation}
		where the remainders $R_{\eps,2}$ and $R_{\eps,3}$ satisfy
		\begin{equation}\label{remainder2 strong}
			\begin{aligned}
				R_{\eps,2} & = o(1)_{L^\infty\left(dt;L^1_{\mathrm{loc}}(dx)\right)}
				\\
				\left\|R_{\eps,3}\right\|_{W^{-1,1}_\mathrm{loc}\left(dx\right)}
				& \leq C_1\delta\mathcal{H}_\eps(t) + C_2\delta\mathcal{E}_\eps(t) + o(1)_{L^1_\mathrm{loc}(dt)},
			\end{aligned}
		\end{equation}
		for some $C_1,C_2>0$.
		
		\item The approximate Ohm's law
		\begin{equation}\label{Ohm consistent strong}
			\begin{aligned}
				\frac 1\sigma
				\int_{\mathbb{R}^3\times\mathbb{R}^3\times\mathbb{S}^2}
				\left(
				\hat q_\eps^+ - \hat q_\eps^-
				+ \hat q_\eps^{+,-}-\hat q_\eps^{-,+}\right)
				\tilde\Phi MM_*dvdv_*d\sigma
				\hspace{-40mm} &
				\\
				& = \frac 12 \nabla_x\tilde n_\eps -\left(E_\eps + \tilde u_\eps\wedge B_\eps\right)
				+R_{\eps,4}+\nabla_x R_{\eps,5}+R_{\eps,6},
			\end{aligned}
		\end{equation}
		where $\sigma>0$ is defined by \eqref{sigma 3} and the remainders $R_{\eps,4}$ and $R_{\eps,5}$ vanish weakly
		\begin{equation}\label{remainder4 strong}
			R_{\eps,4}=o(1)_{\textit{w-}L^1_\mathrm{loc}(dtdx)}
			\quad\text{and}\quad
			R_{\eps,5}=o(1)_{\textit{w-}L^1_\mathrm{loc}(dtdx)},
		\end{equation}
		whereas $R_{\eps,6}$ satisfies
		\begin{equation}\label{remainder5 strong}
			\left\|R_{\eps,6}\right\|_{L^1(dx)}
			\leq \frac{C\delta \mathcal{E}_\eps(t)}{1-\left\|\bar u\right\|_{L^\infty\left(dx\right)}}
			+o(1)_{\textit{w-}L_\mathrm{loc}^1(dt)}.
		\end{equation}
		
		This approximate law is obtained directly from the limiting laws derived in Proposition \ref{high weak-comp2}. Indeed, it is easily deduced from \eqref{pre Ohm strong} that \eqref{Ohm consistent strong} holds with the remainders
		\begin{equation*}
			\begin{aligned}
				R_{\eps,4}
				& =
				\frac 1\sigma
				\int_{\mathbb{R}^3\times\mathbb{R}^3\times\mathbb{S}^2}
				\left(
				\hat q_\eps^+ - \hat q_\eps^-
				+\hat q_\eps^{+,-}-\hat q_\eps^{-,+}\right)
				\tilde \Phi MM_*dvdv_*d\sigma
				\\
				& - \frac 1\sigma \int_{\mathbb{R}^3\times\mathbb{R}^3\times\mathbb{S}^2} \left(q^+-q^-+q^{+,-}-q^{-,+}\right)
				\tilde \Phi MM_*dvdv_*d\sigma
				\\
				& + \left(E_\eps-E\right)
				+ \left(\tilde u_\eps-u\right)\wedge B + u \wedge \left(B_\eps-B\right),
				\\
				R_{\eps,5} & =-\frac 12\left(\tilde n_\eps-n\right),
				\\
				R_{\eps,6} & =
				\left(\tilde u_\eps-u\right)\wedge\left(B_\eps - B\right)
				=
				\frac{O\left(\delta\mathcal{E}_\eps(t)+\delta\mathcal{E}(t)\right)_{L^1\left(dx\right)}}{1-\left\|\bar u\right\|_{L^\infty\left(dx\right)}},
			\end{aligned}
		\end{equation*}
		where we have used \eqref{test u estimate}. The above estimate on $R_{\eps,6}$ is then readily improved to \eqref{remainder5 strong} upon noticing from \eqref{energy liminf strong} that
		\begin{equation*}
			\delta\mathcal{E}(t)\leq \delta\mathcal{E}_0(t),
		\end{equation*}
		where $\delta\mathcal{E}_0(t)$ is the limit, up to extraction of subsequences, of $\delta\mathcal{E}_\eps(t)$ in $\textit{w$^*$-}L^\infty(dt)$, and then writing
		\begin{equation*}
			\begin{aligned}
				\delta\mathcal{E}_\eps(t)+\delta\mathcal{E}(t)
				& \leq
				2\delta \mathcal{E}_\eps(t)+\delta\mathcal{E}_0(t)-\delta\mathcal{E}_\eps(t)
				\\
				& =2\delta \mathcal{E}_\eps(t)+o(1)_{\textit{w$^*$-}L^\infty(dt)}.
			\end{aligned}
		\end{equation*}
		
		As in the case of the approximate solenoidal Ohm's law \eqref{Ohm consistent} for weak interspecies interactions, it would be possible to derive the above approximate Ohm's law employing the methods of proof of Proposition \ref{approx2-prop}. Nevertheless, the method presented here is more robust.

		\item Maxwell's equations
		\begin{equation}\label{Maxwell consistent strong}
			\begin{cases}
				\begin{aligned}
					\d_t E_\eps - \ROT B_\eps &= -j_\eps = -\tilde j_\eps +R_{\eps,7},
					\\
					\d_t B_\eps + \ROT E_\eps & = 0,
					\\
					\Div E_\eps & = n_\eps = \tilde n_\eps - R_{\eps,8},
					\\
					\Div B_\eps & = 0,
				\end{aligned}
			\end{cases}
		\end{equation}
		where the remainders $R_{\eps,7}=\tilde j_\eps - j_\eps$ and $R_{\eps,8}=\tilde n_\eps - n_\eps$ satisfy
		\begin{equation}\label{remainder6 strong}
			\left\| R_{\eps,7} \right\|_{L^1_\mathrm{loc}\left(dx\right)}
			\leq C\delta\mathcal{H}_\eps(t) + o(1)_{L^1_{\mathrm{loc}}\left(dt\right)},
		\end{equation}
		and
		\begin{equation}\label{remainder7 strong}
			R_{\eps,8}=o(1)_{L^\infty\left(dt;L^1_\mathrm{loc}(dx)\right)}.
		\end{equation}
		
		The convergence \eqref{remainder7 strong} straightforwardly follows from \eqref{densities approx}. As for the control \eqref{remainder6 strong}, it is obtained through the following estimate. First, since $G^\pm_\eps\geq 2$ and $\eps \hat g_\eps^\pm\geq 2\left(\sqrt 2 -1 \right)$ on the support of $1-\gamma_\eps^\pm$, we easily deduce, using Lemma \ref{trunc-lem 3}, that
		\begin{equation*}
			\begin{aligned}
				\left\|\frac 1\eps g_\eps^\pm\left(1-\gamma_\eps^\pm\right)\right\|_{L^1_\mathrm{loc}\left(dx;L^1\left((1+|v|)^2 Mdv\right)\right)}
				\hspace{-30mm}&
				\\
				& =
				\left\|\frac 1\eps\left( \hat g_\eps^\pm + \frac\eps 4 \hat g_\eps^{\pm 2}\right)\left(1-\gamma_\eps^\pm\right)\right\|_{L^1_\mathrm{loc}\left(dx;L^1\left((1+|v|)^2 Mdv\right)\right)}
				\\
				& \leq C \left\|\mathds{1}_{\left\{G^\pm_\eps\geq 2\right\}}\hat g_\eps^{\pm}\right\|_{L^2_\mathrm{loc}\left(dx;L^2\left((1+|v|)^2 Mdv\right)\right)}^2
				\\
				& \leq
				C_1\int_{\mathbb{R}^3\times\mathbb{R}^3} \left(\frac 1{\eps^2}{h\left(\eps g_\eps^\pm\right)}
				- \frac{1}{2} \left( g_\eps^\pm \gamma_\eps^\pm\chi \left({|v|^2\over K_\eps}\right) \right)^2 \right)Mdxdv
				\\
				& + C_2\left\|g_\eps^\pm \gamma_\eps^\pm\chi \left({|v|^2\over K_\eps}\right) -\bar g^\pm\right\|_{L^2\left(Mdxdv\right)}^2
				+ o(1)_{L^1_{\mathrm{loc}}\left(dt\right)}
				\\
				& \leq C_1\delta\mathcal{H}_\eps(t)+C_2\delta\mathcal{E}_\eps(t)+ o(1)_{L^1_{\mathrm{loc}}\left(dt\right)}.
			\end{aligned}
		\end{equation*}
		Furthermore, using the Gaussian decay \eqref{gaussian-decay 0} and that $g_\eps^\pm\gamma_\eps^\pm$ is comparable to $\hat g_\eps$, we also obtain
		\begin{equation*}
			\begin{aligned}
				\left\|\frac 1\eps g_\eps^\pm\gamma_\eps^\pm\left(1-\chi\left(\frac{|v|^2}{K_\eps}\right)\right)\right\|_{L^1\left((1+|v|)^2 Mdv\right)}
				\hspace{-30mm}&
				\\
				& \leq C
				\left\|\frac 1\eps \hat g_\eps^\pm\left(1-\chi\left(\frac{|v|^2}{K_\eps}\right)\right)\right\|_{L^1\left((1+|v|)^2 Mdv\right)}
				\\
				& \leq C
				\left\|\hat g_\eps^\pm\right\|_{L^2\left(Mdv\right)}
				\left\|\frac 1\eps \left(1-\chi\left(\frac{|v|^2}{K_\eps}\right)\right)\right\|_{L^2\left((1+|v|)^4 Mdv\right)}
				\\
				& = O\left(K^\frac 54\left|\log\eps\right|^\frac 54\eps^{\frac K4-1}\right)_{L^\infty\left(dt;L^2(dx)\right)}.
			\end{aligned}
		\end{equation*}
		Thus, further using Lemma \ref{entropy energy strong}, we infer, provided $K>4$, that
		\begin{equation*}
			\frac 1\eps \left\| g_\eps^\pm - g_\eps^\pm\gamma_\eps^\pm\chi\left(\frac{|v|^2}{K_\eps}\right)\right\|_{L^1_\mathrm{loc}\left(dx;L^1\left((1+|v|)^2 Mdv\right)\right)}
			\leq C\delta\mathcal{H}_\eps(t) + o(1)_{L^1_{\mathrm{loc}}\left(dt\right)},
		\end{equation*}
		whence
		\begin{equation}\label{tilde u approx u strong}
			\left\| \tilde j_\eps-j_\eps \right\|_{L^1_\mathrm{loc}\left(dx\right)}
			\leq C\delta\mathcal{H}_\eps(t) + o(1)_{L^1_{\mathrm{loc}}\left(dt\right)},
		\end{equation}
		which establishes \eqref{remainder6 strong}.

		Notice that we cannot rigorously write the identities \eqref{maxwell energy two species} and \eqref{maxwell poynting two species} for the above system, because the source terms $j_\eps$ and $n_\eps$ do not belong to $L^2_\mathrm{loc}(dtdx)$ a priori. Nevertheless, one has the following modulated identities~:
		\begin{equation}\label{modulated maxwell energy strong}
			\begin{aligned}
				\partial_t\left(E_\eps\cdot \bar E+B_\eps\cdot\bar B\right)
				& +
				\nabla_x\cdot\left(E_\eps\wedge \bar B + \bar E\wedge B_\eps\right)
				\\
				& = -\left(\tilde j_\eps -R_{\eps,7}\right)\cdot\bar E
				- \left(\bar j +\mathbf{A}_4\right)\cdot E_\eps
				- \mathbf{A}_5\cdot B_\eps,
			\end{aligned}
		\end{equation}
		and
		\begin{equation}\label{modulated poynting strong}
			\begin{aligned}
				\partial_t & \left(\left(E_\eps-\bar E\right)\wedge \left(B_\eps-\bar B\right) \right)
				+\frac 12\nabla_x\left(\left|E_\eps-\bar E\right|^2 + \left|B_\eps-\bar B\right|^2 \right)
				\\
				& -\nabla_x\cdot\left(\left(E_\eps-\bar E\right)\otimes \left(E_\eps-\bar E\right)+\left(B_\eps-\bar B\right)\otimes \left(B_\eps-\bar B\right)\right)
				\\
				& =
				\partial_t\left(E_\eps\wedge B_\eps \right)
				+\frac 12\nabla_x\left(\left|E_\eps\right|^2 + \left|B_\eps\right|^2 \right)
				-\nabla_x\cdot\left(E_\eps\otimes E_\eps+B_\eps\otimes B_\eps\right)
				\\
				& + \left(\bar j+\mathbf{A}_4\right)\wedge \left(B_\eps-\bar B\right) +\left(E_\eps-\bar E\right)\wedge \mathbf{A}_5
				\\
				& +\left(\tilde j_\eps-R_{\eps, 7}\right)\wedge \bar B + \left( \tilde n_\eps - R_{\eps,8} \right)\bar E
				+\bar n\left(E_\eps-\bar E\right).
			\end{aligned}
		\end{equation}

		Finally, taking the divergence of the approximate Amp\`ere equation from \eqref{Maxwell consistent strong}, we obtain the approximate conservation of charge (or approximate continuity equation)
		\begin{equation}\label{conservation charge strong}
			\partial_t \tilde n_\eps + \nabla_x\cdot \tilde j_\eps
			=
			\partial_t R_{\eps,8} + \nabla_x\cdot R_{\eps,7}.
		\end{equation}
		Note that we could just as well use the approximate conservation of charge from Proposition \ref{approx2-prop}.

		\item The asymptotic characterization \eqref{q phi psi 3} of the limiting collision integrands from Proposition \ref{high weak-comp2} combined with \eqref{null dissipation} from Proposition \ref{high weak-comp3}, which implies that
		\begin{equation}\label{asymptotic mu kappa strong}
			\begin{aligned}
				\left(\int_{\mathbb{R}^3\times\mathbb{R}^3\times\mathbb{S}^2} \frac{\hat q_\eps^++\hat q_\eps^-}{2} \tilde \phi MM_* dvdv_*d\sigma\right)
				-\mu\left(\nabla_x\tilde u_\eps+\nabla_x^t\tilde u_\eps\right)
				& \to 0,
				\\
				\left(\int_{\mathbb{R}^3\times\mathbb{R}^3\times\mathbb{S}^2} \frac{\hat q_\eps^++\hat q_\eps^-}{2} \tilde \psi MM_* dvdv_*d\sigma\right)
				-\frac 52\kappa\nabla_x\tilde\theta_\eps
				& \to 0,
				\\
				\hat q_\eps^++\hat q_\eps^--\hat q_\eps^{+,-}-\hat q_\eps^{-,+}
				& \to 0,
			\end{aligned}
		\end{equation}
		in the sense of distributions, where $\mu,\kappa>0$ are defined by \eqref{mu kappa 2}.

		\item The asymptotic characterizations \eqref{mixed q phi psi 2} and \eqref{mixed q phi psi 3} of the limiting collision integrands from Proposition \ref{high weak-comp3}, whose proofs imply that
		\begin{equation}\label{asymptotic sigma strong}
			\int_{\mathbb{R}^3\times\mathbb{R}^3\times\mathbb{S}^2} \left(\hat q_\eps^{+}-\hat q_\eps^{-}+\hat q_\eps^{+,-}-\hat q_\eps^{-,+}\right) \tilde\Phi M M_* dvdv_*d\sigma
			+\tilde j_\eps - \tilde n_\eps \tilde u_\eps
			=R_{\eps,9},
		\end{equation}
		where the remainder $R_{\eps,9}$ satisfies
		\begin{equation}\label{remainder8 strong}
			\left\| R_{\eps,9} \right\|_{L^1_\mathrm{loc}\left(dx\right)}
			\leq C\delta\mathcal{H}_\eps(t) + o(1)_{L^1_{\mathrm{loc}}\left(dt\right)}.
		\end{equation}
		
		Indeed, we first obtain from \eqref{asymptotic sigma strong 2}, using Lemmas \ref{trunc-lem 4} and \ref{entropy energy strong}, that
		\begin{equation}\label{asymptotic sigma strong 5}
			\begin{aligned}
				\left\|h_\eps- \hat h_\eps - \frac 14 \hat n_\eps\left(\hat g_\eps^+-\hat \rho_\eps^+ + \hat g_\eps^--\hat \rho_\eps^-\right)\right\|_{L^1_\mathrm{loc}\left(dx;L^1\left(\left(1+|v|^2\right)Mdv\right)\right)} \hspace{-20mm} &
				\\
				& \leq C\delta\mathcal{H}_\eps(t)+o(1)_{L^1_\mathrm{loc}(dt)},
			\end{aligned}
		\end{equation}
		where $\hat h_\eps = \frac{1}{\eps}\left[\left(\hat g_\eps^+-\hat g_\eps^-\right) - \hat n_\eps\right]$, $\hat n_\eps=\hat \rho_\eps^+-\hat\rho_\eps^-$ and $\hat\rho_\eps^\pm$ are the densities associated with the fluctuations $\hat g_\eps^\pm$. Next, combining \eqref{asymptotic sigma strong 3} with \eqref{asymptotic sigma strong 4}, straightforward computations yield that
		\begin{equation*}
			\begin{aligned}
				& \left(\mathcal{L}+\mathfrak{L}\right) \left(h_\eps
				- \frac 12\hat n_\eps \left(\hat g_\eps^+ + \hat g_\eps^-\right)\right)
				+\int_{\mathbb{R}^3\times\mathbb{S}^2} \left(\hat q_\eps^{+}-\hat q_\eps^{-}+\hat q_\eps^{+,-}-\hat q_\eps^{-,+}\right) M_* dv_*d\sigma
				\\
				& =
				\frac 1 2 \cQ\left(\hat g_\eps^+-\hat g_\eps^- - \hat n_\eps ,\hat g_\eps^+ + \hat g_\eps^- \right)
				+\left(\mathcal{L}+\mathfrak{L}\right)\left(h_\eps- \hat h_\eps - \frac 14 \hat n_\eps\left(\hat g_\eps^+-\hat \rho_\eps^+ + \hat g_\eps^--\hat \rho_\eps^-\right)\right)
				\\
				& - \frac 12\hat n_\eps \mathcal{L}\left(\hat g_\eps^+ + \hat g_\eps^--\Pi\left(\hat g_\eps^+ + \hat g_\eps^-\right)\right),
			\end{aligned}
		\end{equation*}
		which implies, using \eqref{phi-psi-def inverses two species} and the self-adjointness of $\mathcal{L}+\mathfrak{L}$ and then employing Lemmas \ref{trunc-lem}, \ref{trunc-lem 2} (on consistency estimates) and \ref{entropy energy strong} (allowing to control the energy by the entropy) with the estimates \eqref{relaxation estimate} and \eqref{asymptotic sigma strong 5}, that
		\begin{equation*}
			\begin{aligned}
				\left\| j_\eps - \hat n_\eps \hat u_\eps
				+\int_{\mathbb{R}^3\times\mathbb{R}^3\times\mathbb{S}^2} \left(\hat q_\eps^{+}-\hat q_\eps^{-}+\hat q_\eps^{+,-}-\hat q_\eps^{-,+}\right) \tilde\Phi M M_* dvdv_*d\sigma
				\right\|_{L^1_\mathrm{loc}(dx)} \hspace{-30mm} &
				\\
				& \leq C\delta\mathcal{H}_\eps(t) + o(1)_{L^1_\mathrm{loc}(dt)},
			\end{aligned}
		\end{equation*}
		where $\hat u_\eps=\frac 12 \int_{\mathbb{R}^3}\left(\hat g_\eps^+ + \hat g_\eps^-\right)vMdv$. Finally, utilizing the control \eqref{tilde u approx u strong} with yet another application of Lemma \ref{trunc-lem} allows us to deduce the validity of \eqref{asymptotic sigma strong} from the preceding estimate.
	\end{itemize}

	Now, by definition of the acceleration operator $\mathbf{A}$, straightforward energy computations, similar to those performed in the proof of Proposition \ref{energy estimate 2}, applied to the test functions $\left( \bar u, \bar \theta, \bar n, \bar j, \bar E, \bar B\right)$, show that the following energy identity holds~:
	\begin{equation}\label{test energy strong}
		\frac{d}{dt}\bar\CE(t)+\bar\CD(t)= - \int_{\mathbb{R}^3} \mathbf{A}\cdot
		\begin{pmatrix}
			\bar u \\ \frac 52 \bar \theta \\ \bar j - \bar n\bar u \\ \bar E - \frac 12 \nabla_x\bar n \\ \bar B
		\end{pmatrix}
		dx,
	\end{equation}
	where the energy $\bar\CE$ and energy dissipation $\bar\CD$ are defined by
	\begin{equation*}
		\begin{aligned}
			\bar{\mathcal{E}}(t)
			& =
			\frac 12\left\|\bar g^+\right\|_{L^2\left(Mdxdv\right)}^2
			+\frac 12\left\|\bar g^-\right\|_{L^2\left(Mdxdv\right)}^2
			+\frac 12\left\|\bar E\right\|_{L^2(dx)}^2
			+\frac 12\left\|\bar B\right\|_{L^2(dx)}^2
			\\
			& =
			\left\|\bar \rho\right\|_{L^2(dx)}^2 + \frac 14\left\|\bar n\right\|_{L^2(dx)}^2
			+\left\|\bar u\right\|_{L^2(dx)}^2+\frac 32\left\|\bar \theta\right\|_{L^2(dx)}^2
			\\
			& +\frac 12\left\|\bar E\right\|_{L^2(dx)}^2
			+\frac 12\left\|\bar B\right\|_{L^2(dx)}^2,
			\\
			& =
			\frac 14\left\|\bar n\right\|_{L^2(dx)}^2
			+\left\|\bar u\right\|_{L^2(dx)}^2+\frac 52\left\|\bar \theta\right\|_{L^2(dx)}^2
			+\frac 12\left\|\bar E\right\|_{L^2(dx)}^2
			+\frac 12\left\|\bar B\right\|_{L^2(dx)}^2,
		\end{aligned}
	\end{equation*}
	and
	\begin{equation*}
		\begin{aligned}
			\bar{\mathcal{D}}(t)
			& =
			2\mu
			\left\|\nabla_x \bar u\right\|_{L^2_x}^2
			+ 5\kappa
			\left\|\nabla_x \bar\theta\right\|_{L^2_x}^2
			+ \frac 1\sigma
			\left\|\bar j - \bar n\bar u\right\|_{L^2_x}^2
			\\
			& =
			\frac 1{16}
			\left\|\bar q^+ + \bar q^- + \bar q^{+,-} + \bar q^{-,+}\right\|_{L^2\left(MM_*dxdvdv_*d\sigma\right)}^2
			\\
			& + \frac 18
			\left\|\bar q^+ - \bar q^-\right\|_{L^2\left(MM_*dxdvdv_*d\sigma\right)}^2
			+
			\frac 18
			\left\|\bar q^{+,-} - \bar q^{-,+}\right\|_{L^2\left(MM_*dxdvdv_*d\sigma\right)}^2
			\\
			& =
			\frac 14
			\left\|\bar q^+\right\|_{L^2\left(MM_*dxdvdv_*d\sigma\right)}^2
			+
			\frac 14
			\left\|\bar q^-\right\|_{L^2\left(MM_*dxdvdv_*d\sigma\right)}^2
			\\
			& + \frac 14
			\left\|\bar q^{+,-}\right\|_{L^2\left(MM_*dxdvdv_*d\sigma\right)}^2
			+
			\frac 14
			\left\|\bar q^{-,+}\right\|_{L^2\left(MM_*dxdvdv_*d\sigma\right)}^2.
		\end{aligned}
	\end{equation*}

	Next, notice that a slight variant of the estimate \eqref{conservation-laws4} derived in the proof of Proposition \ref{stab-prop} on weak interactions is also valid here in the case of strong interactions. Indeed, reproducing the very same duality computations preceding \eqref{conservation-laws4} onto the approximate conservation of energy \eqref{conservation-laws strong} and, then, using the convergences \eqref{limit solenoidal strong}, \eqref{limit boussinesq strong}, \eqref{asymptotic mu kappa strong}, the estimate \eqref{remainder1 strong} and Lemma \ref{entropy energy strong} (allowing to control the energy by the entropy), yields that
	\begin{equation}\label{conservation-laws4 strong}
		\begin{aligned}
			\frac{d}{dt} \int_{\mathbb{R}^3} & \left(\frac 32 \tilde\theta_\eps-\tilde\rho_\eps\right)\cdot\bar\theta dx
			\\
			& +\frac 12
			\int_{\mathbb{R}^3}
			\left( \int_{\mathbb{R}^3\times\mathbb{R}^3\times\mathbb{S}^2} \left(\hat q_\eps^++\hat q_\eps^- + \hat q_\eps^{+,-}+\hat q_\eps^{-,+}\right) \tilde \psi MM_* dvdv_*d\sigma\right)
			\cdot\nabla_x\bar\theta
			dx
			\\
			& \geq
			-C\left\|\bar\theta\right\|_{W^{1,\infty}(dx)}
			\left(\frac{\delta\mathcal{E}_\eps(t)}{1-\left\|\bar u\right\|_{L^\infty(dx)}}
			+\left(\delta\mathcal{E}_\eps(t)\delta\mathcal{D}_\eps(t)\right)^\frac{1}{2}\right)
			\\
			& - \frac{1}{2}\int_{\mathbb{R}^3}
			\mathbf{A}_2\left(\frac 32 \tilde\theta_\eps-\tilde\rho_\eps\right)
			dx
			+o(1)_{\textit{w-}L^1_\mathrm{loc}(dt)}
			\\
			& \geq -C\left(\frac{\left\|\bar\theta\right\|_{W^{1,\infty}(dx)}}{1-\left\|\bar u\right\|_{L^\infty(dx)}}
			+\left\|\bar\theta\right\|_{W^{1,\infty}(dx)}^2\right)
			\delta\mathcal{H}_\eps(t)
			-\frac 14 \delta\mathcal{D}_\eps(t)
			\\
			& - \frac{1}{2}\int_{\mathbb{R}^3}
			\mathbf{A}_2\left(\frac 32 \tilde\theta_\eps-\tilde\rho_\eps\right)
			dx
			+o(1)_{\textit{w-}L^1_\mathrm{loc}(dt)}.
		\end{aligned}
	\end{equation}

	Likewise, following the proof of Proposition \ref{stab-prop}, using the solenoidal property $\Div \bar u = 0$, analogous duality computations applied to the approximate conservations of momentum \eqref{conservation-laws2 strong} and charge \eqref{conservation charge strong} yield that
	\begin{equation*}
		\begin{aligned}
			\frac{d}{dt} & \int_{\mathbb{R}^3}
			\left(
			\frac 14 \tilde n_\eps\bar n +
			\tilde u_\eps\cdot\bar u+\frac 12 \left(E_\eps\wedge B_\eps\right)\cdot \bar u
			\vphantom{+\frac 12
			\begin{pmatrix}
				a_{\eps 26}-a_{\eps 35}\\a_{\eps 34}-a_{\eps 16}\\a_{\eps 15}-a_{\eps 24}
			\end{pmatrix}\cdot\bar u - R_{\eps,2}\cdot\bar u
			-\frac 14 R_{\eps,8}\bar n}
			\right.
			\\
			& 
			\left.
			\vphantom{\frac 14 \tilde n_\eps\bar n +
			\tilde u_\eps\cdot\bar u+\frac 12 \left(E_\eps\wedge B_\eps\right)\cdot \bar u}
			+\frac 12
			\begin{pmatrix}
				a_{\eps 26}-a_{\eps 35}\\a_{\eps 34}-a_{\eps 16}\\a_{\eps 15}-a_{\eps 24}
			\end{pmatrix}\cdot\bar u - R_{\eps,2}\cdot\bar u
			-\frac 14 R_{\eps,8}\bar n
			\right)
			dx
			\\
			& - \int_{\mathbb{R}^3}\left(\frac 12 \left(E_\eps\wedge B_\eps\right)\cdot \partial_t\bar u
			+\frac 12
			\begin{pmatrix}
				a_{\eps 26}-a_{\eps 35}\\a_{\eps 34}-a_{\eps 16}\\a_{\eps 15}-a_{\eps 24}
			\end{pmatrix}\cdot\partial_t\bar u - R_{\eps,2}\cdot\partial_t\bar u
			-\frac 14 R_{\eps,8}\partial_t\bar n\right) dx
			\\
			& +
			\int_{\mathbb{R}^3}
			\frac 14 \left(\tilde n_\eps\nabla_x\cdot \bar j - \tilde j_\eps\nabla_x\bar n\right)
			+
			\left(\left(P\tilde u_\eps\right)\otimes \bar u- \tilde u_\eps\otimes\tilde u_\eps - \frac{1}{2\eps^2}m_\eps\right) :\nabla_x\bar u
			dx
			\\
			& +
			\int_{\mathbb{R}^3}
			\left(\int_{\mathbb{R}^3\times\mathbb{R}^3\times\mathbb{S}^2} \frac{\hat q_\eps^++\hat q_\eps^-}{2} \tilde \phi MM_* dvdv_*d\sigma\right)
			:\nabla_x\bar u
			-\mu \Delta_x\bar u \cdot\tilde u_\eps dx
			\\
			& +\frac 12
			\int_{\mathbb{R}^3}
			\left(E_\eps\otimes E_\eps + e_\eps + B_\eps\otimes B_\eps + b_\eps\right) :\nabla_x\bar u
			dx
			\\
			& =
			\int_{\mathbb{R}^3}
			R_{\eps,3}\cdot\bar u - \frac 14 R_{\eps,7}\cdot \nabla_x\bar n dx
			\\
			& + \int_{\mathbb{R}^3} \frac 12P\left(\bar n\bar E+\bar j\wedge\bar B\right)\cdot\tilde u_\eps
			-\frac{1}{2}\mathbf{A}_1\cdot \tilde u_\eps - \frac 14 \tilde n_\eps\nabla_x\cdot \mathbf{A}_4
			dx,
		\end{aligned}
	\end{equation*}
	whence, reorganizing some terms so that remainders are moved to the right-hand side,
	\begin{equation*}
		\begin{aligned}
			\frac{d}{dt} & \int_{\mathbb{R}^3}\left(
			\frac 14 \tilde n_\eps \bar n +
			\tilde u_\eps\cdot\bar u+\frac 12 \left(E_\eps\wedge B_\eps\right)\cdot \bar u
			+\frac 12
			\begin{pmatrix}
				a_{\eps 26}-a_{\eps 35}\\a_{\eps 34}-a_{\eps 16}\\a_{\eps 15}-a_{\eps 24}
			\end{pmatrix}\cdot\bar u \right) dx
			\\
			& +
			\int_{\mathbb{R}^3}
			\left(\int_{\mathbb{R}^3\times\mathbb{R}^3\times\mathbb{S}^2} \left(\hat q_\eps^++\hat q_\eps^-\right) \tilde \phi MM_* dvdv_*d\sigma\right)
			:\nabla_x\bar u
			dx
			\\
			& =
			\int_{\mathbb{R}^3}
			R_{\eps,3}\cdot\bar u - \frac 14 R_{\eps,7}\cdot\nabla_x\bar n
			-\frac{1}{2}\mathbf{A}_1\cdot \tilde u_\eps
			-\frac 14\tilde n_\eps\nabla_x\cdot\mathbf{A}_4
			dx
			\\
			& + \int_{\mathbb{R}^3}
			\frac 12P\left(\bar n\bar E+\bar j\wedge\bar B\right)\cdot\tilde u_\eps
			+
			\frac 14 \left(\tilde j_\eps\nabla_x\bar n - \tilde n_\eps\nabla_x\cdot \bar j\right)
			dx
			\\
			& +
			\int_{\mathbb{R}^3}
			\left(\bar u\otimes \left(P^\perp \tilde u_\eps\right) + \left(P^\perp\tilde u_\eps\right)\otimes\bar u +\left(\tilde u_\eps-\bar u\right)\otimes \left(\tilde u_\eps-\bar u\right)
			+ \frac{1}{2\eps^2}m_\eps\right) :\nabla_x\bar u
			dx
			\\
			& +
			\int_{\mathbb{R}^3}
			\mu\tilde u_\eps\cdot \Delta_x\bar u
			+
			\left(\int_{\mathbb{R}^3\times\mathbb{R}^3\times\mathbb{S}^2} \frac{\hat q_\eps^++\hat q_\eps^-}{2} \tilde \phi MM_* dvdv_*d\sigma\right)
			:\nabla_x\bar u dx
			\\
			& + \int_{\mathbb{R}^3}\left(\frac 12 \left(E_\eps\wedge B_\eps\right)\cdot \partial_t\bar u
			+\frac 12
			\begin{pmatrix}
				a_{\eps 26}-a_{\eps 35}\\a_{\eps 34}-a_{\eps 16}\\a_{\eps 15}-a_{\eps 24}
			\end{pmatrix}\cdot\partial_t\bar u - R_{\eps,2}\cdot\partial_t\bar u
			-\frac 14 R_{\eps,8}\partial_t\bar n \right)
			dx
			\\
			& -\frac 12
			\int_{\mathbb{R}^3}
			\left(E_\eps\otimes E_\eps + e_\eps + B_\eps\otimes B_\eps + b_\eps\right) :\nabla_x\bar u
			dx
			+\frac {d}{dt}\int_{\mathbb{R}^3}R_{\eps,2}\cdot\bar u
			+\frac 14 R_{\eps,8}\bar ndx.
		\end{aligned}
	\end{equation*}
	Then, using the convergences \eqref{limit solenoidal strong}, \eqref{asymptotic mu kappa strong}, the estimates \eqref{measureCS1}, \eqref{measureCS2}, \eqref{remainder2 strong}, \eqref{remainder7 strong} and Lemmas \ref{vector defect} and \ref{entropy energy strong} (allowing to control the energy by the entropy), we arrive at
	\begin{equation*}
		\begin{aligned}
			\frac{d}{dt} & \int_{\mathbb{R}^3}\left(
			\frac 14 \tilde n_\eps\bar n
			+ \tilde u_\eps\cdot\bar u+\frac 12 \left(E_\eps\wedge B_\eps\right)\cdot \bar u
			+\frac 12
			\begin{pmatrix}
				a_{\eps 26}-a_{\eps 35}\\a_{\eps 34}-a_{\eps 16}\\a_{\eps 15}-a_{\eps 24}
			\end{pmatrix}\cdot\bar u \right) dx
			\\
			& +
			\int_{\mathbb{R}^3}
			\left(\int_{\mathbb{R}^3\times\mathbb{R}^3\times\mathbb{S}^2} \left(\hat q_\eps^++\hat q_\eps^-\right) \tilde \phi MM_* dvdv_*d\sigma\right)
			:\nabla_x\bar u
			dx
			\\
			& \geq
			- C\left(\left\|\bar u\right\|_{W^{1,\infty}\left(dx\right)}+\frac{\left\|\partial_t\bar u\right\|_{L^\infty(dx)}}{1-\left\|\bar u\right\|_{L^\infty(dx)}}\right)\delta\mathcal{H}_\eps(t)
			-\frac 14\int_{\mathbb{R}^3}R_{\eps,7}\cdot\nabla_x\bar n dx
			\\
			& -\int_{\mathbb{R}^3}\frac 12 \mathbf{A}_1\cdot \tilde u_\eps
			+\frac 14\tilde n_\eps\nabla_x\cdot\mathbf{A}_4
			dx
			+ o(1)_{\textit{w-}L^1_\mathrm{loc}(dt)}
			+\frac{d}{dt}\left(o(1)_{L^\infty\left(dt\right)}\right)
			\\
			& + \int_{\mathbb{R}^3}
			\frac 12P\left(\bar n\bar E+\bar j\wedge\bar B\right)\cdot\tilde u_\eps
			+
			\frac 14 \left(\tilde j_\eps\nabla_x\bar n - \tilde n_\eps\nabla_x\cdot \bar j\right)
			dx
			\\
			& + \frac 12\int_{\mathbb{R}^3}\left(E_\eps\wedge B_\eps\right)\cdot \partial_t\bar u
			-
			\left(E_\eps\otimes E_\eps + e_\eps + B_\eps\otimes B_\eps + b_\eps\right) :\nabla_x\bar u
			dx.
		\end{aligned}
	\end{equation*}

	The next step consists in combining the preceding inequality with the identity \eqref{modulated poynting strong} in order to modulate the Poynting vector $E_\eps\wedge B_\eps$. This yields
	\begin{equation*}
		\begin{aligned}
			\frac{d}{dt} & \int_{\mathbb{R}^3}\left(
			\frac 14\tilde n_\eps\bar n +
			\tilde u_\eps\cdot\bar u+\frac 12 \left(\left(E_\eps-\bar E\right)\wedge \left(B_\eps-\bar B\right)\right)\cdot \bar u
			+\frac 12
			\begin{pmatrix}
				a_{\eps 26}-a_{\eps 35}\\a_{\eps 34}-a_{\eps 16}\\a_{\eps 15}-a_{\eps 24}
			\end{pmatrix}\cdot\bar u \right) dx
			\\
			& +
			\int_{\mathbb{R}^3}
			\left(\int_{\mathbb{R}^3\times\mathbb{R}^3\times\mathbb{S}^2} \left(\hat q_\eps^++\hat q_\eps^-\right) \tilde \phi MM_* dvdv_*d\sigma\right)
			:\nabla_x\bar u
			dx
			\\
			& \geq
			- C\left(\left\|\bar u\right\|_{W^{1,\infty}\left(dx\right)}+\frac{\left\|\partial_t\bar u\right\|_{L^\infty(dx)}}{1-\left\|\bar u\right\|_{L^\infty(dx)}}\right)\delta\mathcal{H}_\eps(t)
			\\
			& + \int_{\mathbb{R}^3}
			\frac{1}{2}\left(\mathbf{A}_4\wedge \left(B_\eps-\bar B\right) +\left(E_\eps-\bar E\right)\wedge \mathbf{A}_5\right)\cdot\bar u
			-\frac 12\mathbf{A}_1\cdot \tilde u_\eps-\frac 14 \tilde n_\eps\nabla_x\cdot\mathbf{A}_4
			dx
			\\
			& + o(1)_{\textit{w-}L^1_\mathrm{loc}(dt)}
			+\frac{d}{dt}\left(o(1)_{L^\infty\left(dt\right)}\right)
			\\
			& +\frac 12 \int_{\mathbb{R}^3} \left( \left(\bar u\wedge \bar B-\frac 12 \nabla_x\bar n\right)\cdot R_{\eps, 7}
			- R_{\eps,8}\bar E\cdot\bar u \right)
			+ \left(\bar n\bar E + \bar j\wedge\bar B\right)\cdot P^\perp\tilde u_\eps
			dx
			\\
			& +\frac 12 \int_{\mathbb{R}^3}
			\bar n\bar E\cdot \tilde u_\eps + \tilde n_\eps\bar E \cdot\bar u
			+\bar n\left(E_\eps-\bar E\right)\cdot\bar u
			dx
			\\
			& +\frac 12\int_{\mathbb{R}^3}
			\left(\frac 12\nabla_x\tilde n_\eps - \tilde u_\eps\wedge B_\eps\right)\cdot\bar j
			+ \left(\frac 12\nabla_x\bar n - \bar u\wedge \bar B\right)\cdot\tilde j_\eps
			dx
			\\
			& + \frac 12\int_{\mathbb{R}^3}
			\left(\bar j\wedge \left(B_\eps-\bar B\right)\right)\cdot\left(\bar u-\tilde u_\eps\right)
			+
			\left(\left(E_\eps-\bar E\right)\wedge \left(B_\eps-\bar B\right)\right)\cdot \partial_t\bar udx
			\\
			& - \frac 12\int_{\mathbb{R}^3}
			\left(\left(E_\eps-\bar E\right)\otimes \left(E_\eps-\bar E\right) + e_\eps + \left(B_\eps-\bar B\right)\otimes \left(B_\eps-\bar B\right) + b_\eps\right) :\nabla_x\bar u
			dx.
		\end{aligned}
	\end{equation*}
	It then follows, using the convergence \eqref{limit solenoidal strong}, the estimates \eqref{measureCS2}, \eqref{remainder7 strong} and Lemma \ref{vector defect}, that
	\begin{equation*}
		\begin{aligned}
			\frac{d}{dt} & \int_{\mathbb{R}^3}\left(
			\frac 14\tilde n_\eps\bar n +
			\tilde u_\eps\cdot\bar u+\frac 12 \left(\left(E_\eps-\bar E\right)\wedge \left(B_\eps-\bar B\right)\right)\cdot \bar u
			+\frac 12
			\begin{pmatrix}
				a_{\eps 26}-a_{\eps 35}\\a_{\eps 34}-a_{\eps 16}\\a_{\eps 15}-a_{\eps 24}
			\end{pmatrix}\cdot\bar u \right) dx
			\\
			& +
			\int_{\mathbb{R}^3}
			\left(\int_{\mathbb{R}^3\times\mathbb{R}^3\times\mathbb{S}^2} \left(\hat q_\eps^++\hat q_\eps^-\right) \tilde \phi MM_* dvdv_*d\sigma\right)
			:\nabla_x\bar u
			dx
			\\
			& \geq
			- C\left(\frac{\left\|\bar u\right\|_{W^{1,\infty}\left(dx\right)}+\left\|\partial_t\bar u\right\|_{L^\infty(dx)}}{1-\left\|\bar u\right\|_{L^\infty(dx)}}
			\right)
			\delta\mathcal{H}_\eps(t)
			\\
			& + \int_{\mathbb{R}^3}
			\frac{1}{2}\left(\mathbf{A}_4\wedge \left(B_\eps-\bar B\right) +\left(E_\eps-\bar E\right)\wedge \mathbf{A}_5\right)\cdot\bar u
			-\frac 12\mathbf{A}_1\cdot \tilde u_\eps-\frac 14 \tilde n_\eps\nabla_x\cdot\mathbf{A}_4
			dx
			\\
			& + o(1)_{\textit{w-}L^1_\mathrm{loc}(dt)}
			+\frac{d}{dt}\left(o(1)_{L^\infty\left(dt\right)}\right)
			+\frac 12 \int_{\mathbb{R}^3} \left(\bar u\wedge \bar B-\frac 12 \nabla_x\bar n\right)\cdot R_{\eps, 7}
			dx
			\\
			& +\frac 12 \int_{\mathbb{R}^3}
			\bar n\bar E\cdot \tilde u_\eps + \tilde n_\eps\bar E \cdot\bar u
			+\bar n\left(E_\eps-\bar E\right)\cdot\bar u
			dx
			\\
			& +\frac 12\int_{\mathbb{R}^3}
			\left(\frac 12\nabla_x\tilde n_\eps - \tilde u_\eps\wedge B_\eps\right)\cdot\bar j
			+ \left(\frac 12\nabla_x\bar n - \bar u\wedge \bar B\right)\cdot\tilde j_\eps
			dx
			\\
			& + \frac 12\int_{\mathbb{R}^3}
			\left(\bar j\wedge \left(B_\eps-\bar B\right)\right)\cdot\left(\bar u-\tilde u_\eps\right)
			dx.
		\end{aligned}
	\end{equation*}

	Now, for mere convenience of notation, we introduce the following integrand~:
	\begin{equation*}
		\begin{aligned}
			\mathcal{I}
			& =
			\frac 14\tilde n_\eps\bar n +
			\tilde u_\eps\cdot\bar u
			+\frac 12 \left(E_\eps\cdot \bar E+B_\eps\cdot\bar B\right)
			\\
			& +\frac 12 \left(\left(E_\eps-\bar E\right)\wedge \left(B_\eps-\bar B\right)\right)\cdot \bar u
			+\frac 12
			\begin{pmatrix}
				a_{\eps 26}-a_{\eps 35}\\a_{\eps 34}-a_{\eps 16}\\a_{\eps 15}-a_{\eps 24}
			\end{pmatrix}\cdot\bar u.
		\end{aligned}
	\end{equation*}
	Thus, further employing the identity \eqref{modulated maxwell energy strong}, we find that
	\begin{equation*}
		\begin{aligned}
			% & \frac{d}{dt} \int_{\mathbb{R}^3}
			% \left(
			% \frac 14\tilde n_\eps\bar n +
			% \tilde u_\eps\cdot\bar u
			% +\frac 12 \left(E_\eps\cdot \bar E+B_\eps\cdot\bar B\right)
			% +\frac 12 \left(\left(E_\eps-\bar E\right)\wedge \left(B_\eps-\bar B\right)\right)\cdot \bar u
			% \vphantom{+\frac 12
			% \begin{pmatrix}
			% 	a_{\eps 26}-a_{\eps 35}\\a_{\eps 34}-a_{\eps 16}\\a_{\eps 15}-a_{\eps 24}
			% \end{pmatrix}\cdot\bar u}
			% \right.
			% \\
			% & +\left.
			% \vphantom{\tilde u_\eps\cdot\bar u
			% +\frac 12 \left(E_\eps\cdot \bar E+B_\eps\cdot\bar B\right)
			% +\frac 12 \left(\left(E_\eps-\bar E\right)\wedge \left(B_\eps-\bar B\right)\right)\cdot \bar u}
			% \frac 12
			% \begin{pmatrix}
			% 	a_{\eps 26}-a_{\eps 35}\\a_{\eps 34}-a_{\eps 16}\\a_{\eps 15}-a_{\eps 24}
			% \end{pmatrix}\cdot\bar u \right)
			% dx
			% \\
			& \frac{d}{dt} \int_{\mathbb{R}^3}\mathcal{I}dx
			+
			\int_{\mathbb{R}^3}
			\left(\int_{\mathbb{R}^3\times\mathbb{R}^3\times\mathbb{S}^2} \sum_\pm \hat q_\eps^\pm \tilde \phi MM_* dvdv_*d\sigma\right)
			:\nabla_x\bar u
			dx
			\\
			& \geq
			- C\left(\frac{\left\|\bar u\right\|_{W^{1,\infty}\left(dx\right)}+\left\|\partial_t\bar u\right\|_{L^\infty(dx)}}{1-\left\|\bar u\right\|_{L^\infty(dx)}}
			\right)\delta\mathcal{H}_\eps(t)
			-\frac 12 \int_{\mathbb{R}^3} \mathbf{A}_1\cdot \tilde u_\eps dx
			\\
			& - \frac{1}{2} \int_{\mathbb{R}^3}
			\mathbf{A}_4\cdot\left(E_\eps+\bar u \wedge \left(B_\eps-\bar B\right)-\frac 12 \nabla_x\tilde n_\eps\right)
			+\mathbf{A}_5\cdot\left(B_\eps+\left(E_\eps-\bar E\right)\wedge\bar u\right) dx
			\\
			& + o(1)_{\textit{w-}L^1_\mathrm{loc}(dt)}
			+\frac{d}{dt}\left(o(1)_{L^\infty\left(dt\right)}\right)
			+\frac 12\int_{\mathbb{R}^3}
			\left(\frac 12\nabla_x\tilde n_\eps - E_\eps - \tilde u_\eps\wedge B_\eps\right)\cdot\left(\bar j-\bar n \bar u\right)
			dx
			\\
			& +\frac 12\int_{\mathbb{R}^3}
			\left(\frac 12\nabla_x\bar n - \bar E - \bar u\wedge \bar B\right)\cdot\left(\tilde j_\eps-\tilde n_\eps\tilde u_\eps\right)
			dx
			\\
			& + \frac 12\int_{\mathbb{R}^3}
			\left(\left(\tilde u_\eps-\bar u\right) \wedge \left(B_\eps-\bar B\right)\right)\cdot \left(\bar j-\bar n\bar u\right)
			+\frac 14P^\perp \tilde u_\eps\cdot \nabla_x\left(\bar n^2\right)
			dx
			\\
			& +\frac 12 \int_{\mathbb{R}^3}\left(\left(n-\bar n\right)\left(u-\bar u\right)-R_{\eps,7}\right)\cdot\left(\frac 12\nabla_x\bar n-\bar E-\bar u\wedge\bar B\right)
			dx,
		\end{aligned}
	\end{equation*}
	whence, in view of the convergence \eqref{limit solenoidal strong}, the estimate \eqref{remainder6 strong} and Lemma \ref{entropy energy strong} (allowing to control the energy by the entropy),
	\begin{equation*}
		\begin{aligned}
			% & \frac{d}{dt} \int_{\mathbb{R}^3}
			% \left(
			% \frac 14\tilde n_\eps\bar n +
			% \tilde u_\eps\cdot\bar u
			% +\frac 12 \left(E_\eps\cdot \bar E+B_\eps\cdot\bar B\right)
			% +\frac 12 \left(\left(E_\eps-\bar E\right)\wedge \left(B_\eps-\bar B\right)\right)\cdot \bar u
			% \vphantom{+\frac 12
			% \begin{pmatrix}
			% 	a_{\eps 26}-a_{\eps 35}\\a_{\eps 34}-a_{\eps 16}\\a_{\eps 15}-a_{\eps 24}
			% \end{pmatrix}\cdot\bar u}
			% \right.
			% \\
			% & + \left.
			% \vphantom{\tilde u_\eps\cdot\bar u
			% +\frac 12 \left(E_\eps\cdot \bar E+B_\eps\cdot\bar B\right)
			% +\frac 12 \left(\left(E_\eps-\bar E\right)\wedge \left(B_\eps-\bar B\right)\right)\cdot \bar u}
			% \frac 12
			% \begin{pmatrix}
			% 	a_{\eps 26}-a_{\eps 35}\\a_{\eps 34}-a_{\eps 16}\\a_{\eps 15}-a_{\eps 24}
			% \end{pmatrix}\cdot\bar u \right)
			% dx
			% +
			% \\
			& \frac{d}{dt} \int_{\mathbb{R}^3}\mathcal{I}dx
			+
			\int_{\mathbb{R}^3}
			\left(\int_{\mathbb{R}^3\times\mathbb{R}^3\times\mathbb{S}^2} \sum_\pm \hat q_\eps^\pm \tilde \phi MM_* dvdv_*d\sigma\right)
			:\nabla_x\bar u
			dx
			\\
			& \geq
			- C\Bigg(\frac{\left\|\bar u\right\|_{W^{1,\infty}\left(dx\right)}+\left\|\partial_t\bar u\right\|_{L^\infty(dx)}+\left\|\bar j-\bar n \bar u\right\|_{L^\infty(dx)}}{1-\left\|\bar u\right\|_{L^\infty(dx)}}
			\\
			& + \left\|\frac 12\nabla_x\bar n-\bar E-\bar u\wedge\bar B\right\|_{L^\infty(dx)}
			\Bigg)\delta\mathcal{H}_\eps(t)
			-\frac 12 \int_{\mathbb{R}^3} \mathbf{A}_1\cdot \tilde u_\eps dx
			\\
			& - \frac{1}{2} \int_{\mathbb{R}^3}
			\mathbf{A}_4\cdot\left(E_\eps+\bar u \wedge \left(B_\eps-\bar B\right)-\frac 12 \nabla_x\tilde n_\eps\right)
			+\mathbf{A}_5\cdot\left(B_\eps+\left(E_\eps-\bar E\right)\wedge\bar u\right) dx
			\\
			& + o(1)_{\textit{w-}L^1_\mathrm{loc}(dt)}
			+\frac{d}{dt}\left(o(1)_{L^\infty\left(dt\right)}\right)
			+\frac 12\int_{\mathbb{R}^3}
			\left(\frac 12\nabla_x\tilde n_\eps - E_\eps - \tilde u_\eps\wedge B_\eps\right)\cdot\left(\bar j-\bar n \bar u\right)
			dx
			\\
			& +\frac 12\int_{\mathbb{R}^3}
			\left(\frac 12\nabla_x\bar n - \bar E - \bar u\wedge \bar B\right)\cdot\left(\tilde j_\eps-\tilde n_\eps\tilde u_\eps\right)
			dx.
		\end{aligned}
	\end{equation*}

	Using then the approximate Ohm's law \eqref{Ohm consistent strong} with the control \eqref{asymptotic sigma strong} and reorganizing the resulting inequality so that all remainder terms appear on its right-hand side, we obtain
	\begin{equation*}
		\begin{aligned}
			% \frac{d}{dt} & \int_{\mathbb{R}^3}
			% \left(
			% \frac 14\tilde n_\eps\bar n +
			% \tilde u_\eps\cdot\bar u
			% +\frac 12 \left(E_\eps\cdot \bar E+B_\eps\cdot\bar B\right)
			% +\frac 12 \left(\left(E_\eps-\bar E\right)\wedge \left(B_\eps-\bar B\right)\right)\cdot \bar u
			% \vphantom{+\frac 12
			% \begin{pmatrix}
			% 	a_{\eps 26}-a_{\eps 35}\\a_{\eps 34}-a_{\eps 16}\\a_{\eps 15}-a_{\eps 24}
			% \end{pmatrix}\cdot\bar u}
			% \right.
			% \\
			% & + \left.
			% \vphantom{\tilde u_\eps\cdot\bar u
			% +\frac 12 \left(E_\eps\cdot \bar E+B_\eps\cdot\bar B\right)
			% +\frac 12 \left(\left(E_\eps-\bar E\right)\wedge \left(B_\eps-\bar B\right)\right)\cdot \bar u}
			% \frac 12
			% \begin{pmatrix}
			% 	a_{\eps 26}-a_{\eps 35}\\a_{\eps 34}-a_{\eps 16}\\a_{\eps 15}-a_{\eps 24}
			% \end{pmatrix}\cdot\bar u \right)
			% dx
			% +
			% \\
			\frac{d}{dt} & \int_{\mathbb{R}^3}\mathcal{I}dx
			+
			\int_{\mathbb{R}^3}
			\left(\int_{\mathbb{R}^3\times\mathbb{R}^3\times\mathbb{S}^2} \sum_\pm\hat q_\eps^\pm \tilde \phi MM_* dvdv_*d\sigma\right)
			:\nabla_x\bar u
			dx
			\\
			& - \frac 1\sigma \int_{\mathbb{R}^3}
			\left(\int_{\mathbb{R}^3\times\mathbb{R}^3\times\mathbb{S}^2}
			\sum_\pm\left(
			\pm\hat q_\eps^\pm
			\pm \hat q_\eps^{\pm,\mp}
			\right)
			\tilde\Phi MM_* dvdv_*d\sigma\right)
			\cdot\left(\bar j-\bar n\bar u\right)
			dx
			\\
			& \geq
			- C\Bigg(\frac{\left\|\bar u\right\|_{W^{1,\infty}\left(dx\right)}+\left\|\partial_t\bar u\right\|_{L^\infty(dx)}+\left\|\bar j-\bar n\bar u\right\|_{L^\infty(dx)}}{1-\left\|\bar u\right\|_{L^\infty(dx)}}
			\\
			&
			+\left\|\frac 12 \nabla_x\bar n - \bar E - \bar u\wedge\bar B\right\|_{L^\infty(dx)}\Bigg)\delta\mathcal{H}_\eps(t)
			+ o(1)_{\textit{w-}L^1_\mathrm{loc}(dt)}
			+\frac{d}{dt}\left(o(1)_{L^\infty\left(dt\right)}\right)
			\\
			& - \frac{1}{2} \int_{\mathbb{R}^3}
			\mathbf{A}_4\cdot\left(E_\eps+\bar u \wedge \left(B_\eps-\bar B\right)-\frac 12 \nabla_x\tilde n_\eps\right)
			+ \mathbf{A}_5\cdot\left(B_\eps+\left(E_\eps-\bar E\right)\wedge\bar u\right) dx
			\\
			&
			+\frac{1}{2}
			\int_{\mathbb{R}^3}
			\mathbf{A}_3
			\cdot
			\left(
			\int_{\mathbb{R}^3\times\mathbb{R}^3\times\mathbb{S}^2}
			\sum_\pm\left(
			\pm\hat q_\eps^\pm
			\pm \hat q_\eps^{\pm,\mp}\right)
			\tilde\Phi MM_*dvdv_*d\sigma
			\right)
			-\mathbf{A}_1\cdot\tilde u_\eps
			dx
			\\
			&
			-\frac 12
			\int_{\mathbb{R}^3}
			\left(R_{\eps,4}+\nabla_x R_{\eps,5}+R_{\eps,6}\right)\cdot\left(\bar j-\bar n\bar u\right)
			-\left(\frac 12\nabla_x\bar n - \bar E - \bar u\wedge \bar B\right)
			\cdot
			R_{\eps,9}
			dx.
		\end{aligned}
	\end{equation*}
	Thus, in view of the estimates \eqref{remainder4 strong}, \eqref{remainder5 strong}, \eqref{asymptotic mu kappa strong}, \eqref{remainder8 strong} and Lemma \ref{entropy energy strong} (allowing to control the energy by the entropy), we infer that
	\begin{equation}\label{conservation-laws3 strong}
		\begin{aligned}
			% \frac{d}{dt} & \int_{\mathbb{R}^3}
			% \left(
			% \frac 14\tilde n_\eps\bar n +
			% \tilde u_\eps\cdot\bar u
			% +\frac 12 \left(E_\eps\cdot \bar E+B_\eps\cdot\bar B\right)
			% +\frac 12 \left(\left(E_\eps-\bar E\right)\wedge \left(B_\eps-\bar B\right)\right)\cdot \bar u
			% \vphantom{+\frac 12
			% \begin{pmatrix}
			% 	a_{\eps 26}-a_{\eps 35}\\a_{\eps 34}-a_{\eps 16}\\a_{\eps 15}-a_{\eps 24}
			% \end{pmatrix}\cdot\bar u}
			% \right.
			% \\
			% & + \left.
			% \vphantom{\tilde u_\eps\cdot\bar u
			% +\frac 12 \left(E_\eps\cdot \bar E+B_\eps\cdot\bar B\right)
			% +\frac 12 \left(\left(E_\eps-\bar E\right)\wedge \left(B_\eps-\bar B\right)\right)\cdot \bar u}
			% \frac 12
			% \begin{pmatrix}
			% 	a_{\eps 26}-a_{\eps 35}\\a_{\eps 34}-a_{\eps 16}\\a_{\eps 15}-a_{\eps 24}
			% \end{pmatrix}\cdot\bar u \right)
			% dx
			% +
			% \\
			\frac{d}{dt} & \int_{\mathbb{R}^3}\mathcal{I}dx
			+
			\frac 12\int_{\mathbb{R}^3}
			\left(\int_{\mathbb{R}^3\times\mathbb{R}^3\times\mathbb{S}^2} \sum_\pm \left(\hat q_\eps^\pm+\hat q_\eps^{\pm,\mp}\right) \tilde \phi MM_* dvdv_*d\sigma\right)
			:\nabla_x\bar u
			dx
			\\
			& - \frac 1\sigma \int_{\mathbb{R}^3}
			\left(\int_{\mathbb{R}^3\times\mathbb{R}^3\times\mathbb{S}^2}
			\sum_\pm\left(
			\pm\hat q_\eps^\pm
			\pm \hat q_\eps^{\pm,\mp}
			\right)
			\tilde\Phi MM_* dvdv_*d\sigma\right)
			\cdot\left(\bar j-\bar n\bar u\right)
			dx
			\\
			& \geq
			- C\Bigg(\frac{\left\|\bar u\right\|_{W^{1,\infty}\left(dx\right)}+\left\|\partial_t\bar u\right\|_{L^\infty(dx)}+\left\|\bar j-\bar n\bar u\right\|_{L^\infty(dx)}}{1-\left\|\bar u\right\|_{L^\infty(dx)}}
			\\
			&
			+\left\|\frac 12 \nabla_x\bar n - \bar E - \bar u\wedge\bar B\right\|_{L^\infty(dx)}\Bigg)\delta\mathcal{H}_\eps(t)
			+ o(1)_{\textit{w-}L^1_\mathrm{loc}(dt)}
			+\frac{d}{dt}\left(o(1)_{L^\infty\left(dt\right)}\right)
			\\
			& - \frac{1}{2} \int_{\mathbb{R}^3}
			\mathbf{A}_4\cdot\left(E_\eps+\bar u \wedge \left(B_\eps-\bar B\right)-\frac 12 \nabla_x\tilde n_\eps\right)
			+ \mathbf{A}_5\cdot\left(B_\eps+\left(E_\eps-\bar E\right)\wedge\bar u\right) dx
			\\
			&
			+\frac{1}{2}
			\int_{\mathbb{R}^3}
			\mathbf{A}_3
			\cdot
			\left(
			\int_{\mathbb{R}^3\times\mathbb{R}^3\times\mathbb{S}^2}
			\sum_\pm\left(
			\pm\hat q_\eps^\pm
			\pm \hat q_\eps^{\pm,\mp}\right)
			\tilde\Phi MM_*dvdv_*d\sigma
			\right)
			-\mathbf{A}_1\cdot \tilde u_\eps
			dx.
		\end{aligned}
	\end{equation}

	At last, we may now combine the inequalities \eqref{conservation-laws4 strong} and \eqref{conservation-laws3 strong} to deduce, employing the symmetries of collision integrands and \eqref{test integrand strong} to rewrite dissipation terms, that
	\begin{equation*}
		\begin{aligned}
			& \frac{d}{dt} \int_{\mathbb{R}^3}
			\left(
			\left(g_\eps^+ \gamma_\eps^+ \chi\left(\frac{|v|^2}{K_\eps}\right) \bar g^+ +g_\eps^- \gamma_\eps^- \chi\left(\frac{|v|^2}{K_\eps}\right) \bar g^- \right)
			+E_\eps\cdot \bar E+B_\eps\cdot\bar B
			\vphantom{
			+ \left(\left(E_\eps-\bar E\right)\wedge \left(B_\eps-\bar B\right)
			+
			\begin{pmatrix}
				a_{\eps 26}-a_{\eps 35}\\a_{\eps 34}-a_{\eps 16}\\a_{\eps 15}-a_{\eps 24}
			\end{pmatrix}
			\right)\cdot \bar u
			}
			\right.
			\\
			& + \left.
			\vphantom{
			\left(g_\eps^+ \gamma_\eps^++g_\eps^- \gamma_\eps^-\right) \chi\left(\frac{|v|^2}{K_\eps}\right)
			\bar g
			+E_\eps\cdot \bar E+B_\eps\cdot\bar B
			}
			\left(\left(E_\eps-\bar E\right)\wedge \left(B_\eps-\bar B\right)
			+
			\begin{pmatrix}
				a_{\eps 26}-a_{\eps 35}\\a_{\eps 34}-a_{\eps 16}\\a_{\eps 15}-a_{\eps 24}
			\end{pmatrix}
			\right)\cdot \bar u
			\right)
			dx
			\\
			& +\frac 12
			\int_{\mathbb{R}^3\times\mathbb{R}^3\times\mathbb{R}^3\times\mathbb{S}^2}
			\left(\hat q_\eps^+\bar q^++\hat q_\eps^-\bar q^-
			+\hat q_\eps^{+,-}\bar q^{+,-}+\hat q_\eps^{-,+}\bar q^{-,+}
			\right)
			MM_* dxdvdv_*d\sigma
			\\
			& \geq
			- \lambda(t) \delta\mathcal{H}_\eps(t)
			+ o(1)_{\textit{w-}L^1_\mathrm{loc}(dt)}
			+\frac{d}{dt}\left(o(1)_{L^\infty\left(dt\right)}\right)-\frac 12 \delta\mathcal{D}_\eps(t)
			\\
			& - \int_{\mathbb{R}^3} \mathbf{A}\cdot
			\begin{pmatrix}
				\tilde u_\eps \\ \frac 32\tilde \theta_\eps-\tilde\rho_\eps \\
				-\int_{\mathbb{R}^3\times\mathbb{R}^3\times\mathbb{S}^2} \left(\hat q_\eps^+ - \hat q_\eps^- + \hat q_\eps^{+,-}-\hat q_\eps^{-,+}\right) \tilde\Phi MM_* dvdv_*d\sigma
				\\ E_\eps+\bar u \wedge \left(B_\eps-\bar B\right) - \frac 12 \nabla_x\tilde n_\eps \\ B_\eps+\left(E_\eps-\bar E\right)\wedge\bar u
			\end{pmatrix}
			dx.
		\end{aligned}
	\end{equation*}

	Next, assembling the preceding inequality with the scaled entropy inequality \eqref{ent strong} and the energy estimate \eqref{test energy strong}, we finally obtain
	\begin{equation*}
		\begin{aligned}
			& \frac{d}{dt}\delta\mathcal{H}_\eps(t)+\delta\mathcal{D}_\eps(t)
			\\
			& \leq
			\lambda(t) \delta\mathcal{H}_\eps(t)
			+ o(1)_{\textit{w-}L^1_\mathrm{loc}(dt)}
			+\frac{d}{dt}\left(o(1)_{L^\infty\left(dt\right)}\right)+\frac 12 \delta\mathcal{D}_\eps(t)
			\\
			& +\int_{\mathbb{R}^3} \mathbf{A}\cdot
			\begin{pmatrix}
				\tilde u_\eps-\bar u \\ \frac 32\tilde \theta_\eps-\tilde\rho_\eps-\frac 52\bar\theta \\
				- \int_{\mathbb{R}^3\times\mathbb{R}^3\times\mathbb{S}^2}
				\left(\hat q_\eps^+ - \hat q_\eps^- + \hat q_\eps^{+,-}-\hat q_\eps^{-,+}\right) \tilde\Phi MM_* dvdv_*d\sigma
				-\left(\bar j-\bar n\bar u\right)
				\\ E_\eps-\bar E +\bar u \wedge \left(B_\eps-\bar B\right) -\frac 12 \nabla_x\left(\tilde n_\eps-\bar n\right) \\ B_\eps-\bar B +\left(E_\eps-\bar E\right)\wedge\bar u
			\end{pmatrix}
			dx,
		\end{aligned}
	\end{equation*}
	which, with a straightforward application of Gr\"onwall's lemma (carefully note that this is valid even though $\delta\mathcal{H}_\eps(t)$ may be negative), concludes the proof of the proposition.
\end{proof}

\begin{rem}
	The proof of Proposition \ref{stab-prop strong} is based on the construction of the stability inequality \eqref{stability 3} from Proposition \ref{stability poynting} for the two-fluid incompressible Navier-Stokes-Maxwell system with Ohm's law \eqref{TFINSMO}. As in the proof of Proposition \ref{stab-prop}, this approach has the great advantage of using the approximate macroscopic conservation of momentum established in Proposition \ref{approx3-prop} rather than the one from Proposition \ref{approx2-prop} and, thus, removes the difficulties associated with the nonlinear Lorentz force $\tilde n_\eps E_\eps + \tilde j_\eps\wedge B_\eps$ by expressing it with the Poynting vector $E_\eps\wedge B_\eps$ (and some other terms).
	
	However, the drawback of this approach resides in the necessity of the restriction $\left\|\bar u\right\|_{L^\infty_{t,x}}<1$. Recall, nevertheless, that this restriction is physically relevant, since it merely entails that the modulus of the velocity $\bar u$ be less than the speed of light (see comments after the proofs of Propositions \ref{stability poynting} and \ref{stability poynting 2}).
	
	Note finally that it is not possible (at least, we do not know how to make it work) to establish a similar renormalized relative entropy inequality for renormalized solutions of the scaled two species Vlasov-Maxwell-Boltzmann system \eqref{VMB2} based on the construction of the stability inequality \eqref{stability 1} from Proposition \ref{modulated energy estimate 1} (see the remark following the proof of Proposition \ref{stab-prop}).
\end{rem}

\subsection{Convergence and conclusion of proof}

We may now pass to the limit in the approximate stability inequality \eqref{stab ineq eps strong} and, thus, derive the crucial modulated energy inequality for the limiting system \eqref{TFINSFMO 2}. To this end, we simply integrate \eqref{stab ineq eps strong} in time against non-negative test functions and then let $\eps\to 0$, which yields, in view of the well-preparedness of the initial data \eqref{well-prepared init data strong}, the weak convergences \eqref{weak conv fluctuations strong}, \eqref{weak limit observables strong}, \eqref{weak limit charges strong} and the lower semi-continuities \eqref{energy liminf strong}, \eqref{dissipation liminf strong}, that
\begin{equation*}
	\begin{aligned}
		& \delta\mathcal{E}(t) + \frac 12 \int_0^t \delta\mathcal{D}(s) e^{\int_s^t\lambda(\sigma)d\sigma}ds
		\\
		& \leq \delta\mathcal{E}(0) e^{\int_0^t\lambda(s)ds}
		\\
		& +\int_0^t
		\int_{\mathbb{R}^3} \mathbf{A}\cdot
		\begin{pmatrix}
			u-\bar u \\ \frac 32\theta-\rho-\frac 52\bar\theta \\
			- \int_{\mathbb{R}^3\times\mathbb{R}^3\times\mathbb{S}^2}
			\sum_{\pm}\left(\pm q^\pm\pm q^{\pm,\mp}\right)
			\tilde\Phi MM_* dvdv_*d\sigma
			-\left(\bar j-\bar n\bar u\right)
			\\ E-\bar E +\bar u \wedge \left(B-\bar B\right) -\frac 12 \nabla_x\left( n-\bar n\right) \\ B-\bar B +\left(E-\bar E\right)\wedge\bar u
		\end{pmatrix}(s)
		dx
		\\
		& \times
		e^{\int_s^t\lambda(\sigma)d\sigma}ds.
	\end{aligned}
\end{equation*}
Finally, using \eqref{constraints strong} and the characterizations \eqref{mixed q phi psi 2}, \eqref{mixed q phi psi 3} of the limiting collision integrands $q^{\pm}$, $q^{\pm,\mp}$ from Proposition \ref{high weak-comp3}, we deduce that (see also the proof of Proposition \ref{strongOhm} for more detailed computations yielding the term $j-nu$)
\begin{equation*}
	\begin{aligned}
		\delta\mathcal{E}(t) & + \frac 12 \int_0^t \delta\mathcal{D}(s) e^{\int_s^t\lambda(\sigma)d\sigma}ds
		\\
		& \leq \delta\mathcal{E}(0) e^{\int_0^t\lambda(s)ds}
		\\
		& +\int_0^t
		\int_{\mathbb{R}^3} \mathbf{A}\cdot
		\begin{pmatrix}
			u-\bar u \\ \frac 52\left(\theta-\bar\theta\right) \\
			j-nu-\left(\bar j-\bar n\bar u\right)
			\\ E-\bar E +\bar u \wedge \left(B-\bar B\right) - \frac 12 \nabla_x\left(n-\bar n\right) \\ B-\bar B +\left(E-\bar E\right)\wedge\bar u
		\end{pmatrix}(s)
		dx
		e^{\int_s^t\lambda(\sigma)d\sigma}ds,
	\end{aligned}
\end{equation*}
which is precisely the stability inequality we were after.

As for the temporal continuity of $\left(u,n,\frac 52\theta,E,B\right)$, it is readily seen from the approximate macroscopic conservation laws from Proposition \ref{approx2-prop} and Maxwell's equations \eqref{Maxwell consistent strong} that $\partial_t P\tilde u_\eps$, $\partial_t \tilde n_\eps$, $\partial_t \left(\frac 32\tilde\theta_\eps-\tilde\rho_\eps\right)$, $\partial_t E_\eps$ and $\partial_t B_\eps$ are uniformly bounded, in $L^1_\mathrm{loc}$ in time and in some negative index Sobolev space in $x$. It is therefore possible to show (see \cite[Appendix C]{lions7}) that $\left(P\tilde u_\eps, \tilde n_\eps,\frac 32\tilde\theta_\eps-\tilde\rho_\eps,E_\eps,B_\eps\right)$ converges to $\left(u, n,\frac 52\theta,E,B\right)\in C\left([0,\infty);\textit{w-}L^2\left(\mathbb{R}^3\right)\right)$ weakly in $L^2(dx)$ uniformly locally in time.

At last, the proof of Theorem \ref{CV-OMHDSTRONG} is complete.\qed

%% file: appendixtransfer.tex
\chapter{Cross-section for momentum and energy transfer}

The \textbf{cross-section for momentum and energy transfer} $m(z)=m(|z|)\in L^1_\mathrm{loc}\left(\mathbb{R}^3\right)$, such that $m(z)\geq 0$, is defined by
\begin{equation*}
	\begin{aligned}
		\int_{\mathbb{S}^2}\left(v-v'\right) b(v-v_*,\sigma) d\sigma & = m(v-v_*)\left(v-v_*\right), \\
		\int_{\mathbb{S}^2}\left(\frac{|v|^2}2-\frac{|v'|^2}{2}\right) b(v-v_*,\sigma) d\sigma & = m(v-v_*)\left(\frac{|v|^2}2-\frac{|v_*|^2}2\right).
	\end{aligned}
\end{equation*}
Clearly, it is defined as the average transfer of momentum and energy in any collision between any two particles having pre-collisional velocities $v\in\mathbb{R}^3$ and $v_*\in\mathbb{R}^3$. The following proposition guarantees that $m(z)$ is well-defined by the relations above.

\begin{prop}\label{cross section transfer}
	Let
	\begin{equation*}
		m(v-v_*)=m(|v-v_*|)=\frac 12 \int_{\mathbb{S}^2}\left(1-\cos\theta\right) b\left(|v-v_*|,\cos\theta\right)d\sigma,
	\end{equation*}
	with $\cos\theta=\frac{v-v_*}{|v-v_*|}\cdot\sigma$.
	
	It holds that
	\begin{equation*}
		\begin{aligned}
			\int_{\mathbb{S}^2}\left(v-v'\right) b(v-v_*,\sigma) d\sigma & =
			m(v-v_*) \left(v-v_*\right), \\
			\int_{\mathbb{S}^2}\left|v-v'\right|^2 b(v-v_*,\sigma) d\sigma & =
			m(v-v_*) \left|v-v_*\right|^2, \\
			\int_{\mathbb{S}^2}\left(|v|^2-{|v'|}^2\right) b(v-v_*,\sigma) d\sigma & =
			m(v-v_*) \left(|v|^2-|v_*|^2\right),
		\end{aligned}
	\end{equation*}
	and
	\begin{equation*}
		\begin{aligned}
			\int_{\mathbb{R}^3\times\mathbb{R}^3\times\mathbb{S}^2} \left(v-v'\right)\otimes\left(v-v'\right)
			& b(v-v_*,\sigma)MM_*dvdv_* d\sigma \\
			& =\frac 13
			\int_{\mathbb{R}^3\times\mathbb{R}^3}m(v-v_*) |v-v_*|^2 MM_* dvdv_*\operatorname{Id},
			\\
			\int_{\mathbb{R}^3\times\mathbb{R}^3\times\mathbb{S}^2} \left(v-v'\right)\left(|v|^2-{|v'|}^2\right) & b(v-v_*,\sigma)MM_*dvdv_* d\sigma
			\\
			& = 0,
			\\
			\int_{\mathbb{R}^3\times\mathbb{R}^3\times\mathbb{S}^2} \left(|v|^2-{|v'|}^2\right)^2 & b(v-v_*,\sigma)MM_*dvdv_* d\sigma \\
			& =
			\int_{\mathbb{R}^3\times\mathbb{R}^3}m(v-v_*) \left(|v|^2-|v_*|^2\right)^2 MM_* dvdv_*.
		\end{aligned}
	\end{equation*}
\end{prop}

\begin{proof}
	Note first that, using the spherical symmetries of the cross-section $b$,
	\begin{equation*}
		b(v-v_*,\sigma) = b(v-v_*,\tilde\sigma),
	\end{equation*}
	where
	\begin{equation*}
		\tilde\sigma = 2\left(\sigma\cdot\frac{v-v_*}{|v-v_*|}\right) \frac{v-v_*}{|v-v_*|} - \sigma \in\mathbb{S}^2.
	\end{equation*}
	Hence, we compute
	\begin{equation*}
		\begin{aligned}
			\int_{\mathbb{S}^2}\left(v-v'\right) b(v-v_*,\sigma) d\sigma & =
			\int_{\mathbb{S}^2}\left(\left(v-v'\right)\cdot\frac{v-v_*}{|v-v_*|}\right) \frac{v-v_*}{|v-v_*|} b(v-v_*,\sigma) d\sigma \\
			& =
			\frac 12\int_{\mathbb{S}^2}\left(1-\cos\theta\right) (v-v_*) b(v-v_*,\sigma) d\sigma.
		\end{aligned}
	\end{equation*}
	Next, since $|v-v'|^2=\frac 12(1-\cos\theta)|v-v_*|^2$, we easily find that
	\begin{equation*}
		\int_{\mathbb{S}^2}|v-v'|^2 b(v-v_*,\sigma) d\sigma =
		\frac 12\int_{\mathbb{S}^2}\left(1-\cos\theta\right) |v-v_*|^2 b(v-v_*,\sigma) d\sigma.
	\end{equation*}
	Further straightforward computations yield, employing the previous identities, that
	\begin{equation*}
		\begin{aligned}
			\int_{\mathbb{S}^2}\left(|v|^2-|v'|^2\right) b(v-v_*,\sigma) d\sigma
			& =
			\int_{\mathbb{S}^2}\left(2 v\cdot(v-v') - |v-v'|^2\right) b(v-v_*,\sigma) d\sigma \\
			& =
			m(v-v_*)
			\left(2 v\cdot(v-v_*) - |v-v_*|^2\right) \\
			& =
			m(v-v_*)
			\left(|v|^2 - |v_*|^2\right).
		\end{aligned}
	\end{equation*}
	Finally, we obtain, using the pre-post-collisional change of variables and the previous identities,
	\begin{equation*}
		\begin{aligned}
			\int_{\mathbb{R}^3\times\mathbb{R}^3\times\mathbb{S}^2} & \left(v-v'\right)\otimes\left(v-v'\right) b(v-v_*,\sigma)MM_*dvdv_* d\sigma \\
			& =
			2\int_{\mathbb{R}^3\times\mathbb{R}^3\times\mathbb{S}^2} \left(v-v'\right)\otimes v b(v-v_*,\sigma)MM_*dvdv_* d\sigma \\
			& =2
			\int_{\mathbb{R}^3\times\mathbb{R}^3}m(v-v_*) \left(v-v_*\right) \otimes v MM_* dvdv_* \\
			& =
			\int_{\mathbb{R}^3\times\mathbb{R}^3}m(v-v_*) (v-v_*)\otimes(v-v_*) MM_* dvdv_* \\
			& =\frac 13
			\int_{\mathbb{R}^3\times\mathbb{R}^3}m(v-v_*) |v-v_*|^2 MM_* dvdv_* \operatorname{Id},
		\end{aligned}
	\end{equation*}
	\begin{equation*}
		\begin{aligned}
			\int_{\mathbb{R}^3\times\mathbb{R}^3\times\mathbb{S}^2} & \left(v-v'\right)\left(|v|^2-{|v'|}^2\right) b(v-v_*,\sigma)MM_*dvdv_* d\sigma \\
			& =
			2\int_{\mathbb{R}^3\times\mathbb{R}^3\times\mathbb{S}^2} \left(v-v'\right)|v|^2 b(v-v_*,\sigma)MM_*dvdv_* d\sigma \\
			& =2
			\int_{\mathbb{R}^3\times\mathbb{R}^3}m(v-v_*) \left(v-v_*\right) |v|^2 MM_* dvdv_* \\
			& =
			\int_{\mathbb{R}^3\times\mathbb{R}^3}m(v-v_*) (v-v_*)\left(|v|^2-{|v_*|}^2\right) MM_* dvdv_* \\
			& =0,
		\end{aligned}
	\end{equation*}
	and
	\begin{equation*}
		\begin{aligned}
			\int_{\mathbb{R}^3\times\mathbb{R}^3\times\mathbb{S}^2} & \left(|v|^2-{|v'|}^2\right)^2 b(v-v_*,\sigma)MM_*dvdv_* d\sigma \\
			& =
			2\int_{\mathbb{R}^3\times\mathbb{R}^3\times\mathbb{S}^2} \left(|v|^2-{|v'|}^2\right)|v|^2 b(v-v_*,\sigma)MM_*dvdv_* d\sigma \\
			& =2
			\int_{\mathbb{R}^3\times\mathbb{R}^3}m(v-v_*) \left(|v|^2-|v_*|^2\right) |v|^2 MM_* dvdv_* \\
			& =
			\int_{\mathbb{R}^3\times\mathbb{R}^3}m(v-v_*) \left(|v|^2-|v_*|^2\right)^2 MM_* dvdv_*,
		\end{aligned}
	\end{equation*}
	which concludes the justification of the proposition.
\end{proof}

%% file: appendixyoung.tex
\chapter{Young inequalitites}\label{inequalities}

The use of generalized Young inequalities has been ubiquitous in the theory of hydrodynamic limits of the Boltzmann equation since its early treatment in \cite{BGL2}. In its most general version, \textbf{Young's inequality} (also known as Fenchel-Young inequality or Fenchel's inequality, in this case) states that
\begin{equation}\label{fenchel inequality}
	\ip{z,y}_{E,E^*}\leq f(z) + f^*(y),
\end{equation}
where $E$ is a real vector space, $E^*$ is its algebraic dual space, $f(z)$ is a real-valued functional defined on a domain $D\subset E$ and $f^*(y)$ is its \textbf{Legendre transform} (or Legendre-Fenchel transform) defined by
\begin{equation}\label{fenchel transform}
	f^*(y)=\sup_{z\in D}\left(\ip{z,y}_{E,E^*}-f(z)\right),
\end{equation}
on the dual domain
\begin{equation*}
	D^*=\set{y\in E^*}{\sup_{z\in D}\left(\ip{z,y}_{E,E^*}-f(z)\right)<\infty}.
\end{equation*}
Note that $D^*$ is convex and that $f^*$ is lower semi-continuous and convex, for it is defined as the supremum of affine functions. Thus, the transform $f^*$ is also called the convex conjugate of $f$. Clearly, the inequality \eqref{fenchel inequality} is an obvious consequence of the definition \eqref{fenchel transform}.

The Young inequalities are fundamental in extracting useful information and controls from the entropy and the entropy dissipation bounds in \eqref{scaled entropy} and \eqref{scaled entropy two species}. To this end, following \cite{BGL2}, we introduce the non-negative convex functions
\begin{equation*}
	\begin{aligned}
		h(z) & =(1+z)\log(1+z)-z, \\
		r(z) & =z\log(1+z), \\
	\end{aligned}
\end{equation*}
defined over $z>-1$. Notice that $h(z)\leq r(z)$. In this notation, we may then recast the entropy as
\begin{equation*}
	\frac{1}{\eps^2}H(f_\eps)=\frac{1}{\eps^2}H(f_\eps|M)=
	\int_{\mathbb{R}^3\times\mathbb{R}^3} \frac{1}{\eps^2}h\left(\eps g_\eps\right) M dxdv,
\end{equation*}
where $f_\eps=M(1+\eps g_\eps)$, and the entropy dissipation as
\begin{equation*}
	\begin{aligned}
		\frac{1}{\eps^4}D(f_\eps)
		& = \frac{1}{4\eps^4} \int _{\mathbb{R}^3 \times \mathbb{R}^3 \times \mathbb{S}^2}
		r\left(\frac{f'_\eps f'_{\eps *}}{f_\eps f_{\eps *}}-1\right) f_\eps f_{\eps *} b(v-v_*,\sigma)
		dvdv_* d\sigma.
	\end{aligned}
\end{equation*}

We recall now some useful properties, which are already found in \cite{BGL2}, of the convex functions $h(z)$ and $r(z)$. Thus, we consider the Legendre transformations $h^*(y)$ and $r^*(y)$ of $h(z)$ and $r(z)$, respectively, well-defined for any $y\in\mathbb{R}$ by
\begin{equation*}
	h^*(y)=e^y-1-y,
\end{equation*}
and
\begin{equation*}
	r^*(y)=\frac{z_0^2}{1+z_0},
\end{equation*}
where $z_0> -1$ is the unique solution to $y=\log\left(1+z_0\right)+\frac{z_0}{1+z_0}$. Since $h(z)\leq r(z)$, notice that $h^*(y)\geq r^*(y)$. Then, of course, for any $z>-1$ and $y\in\mathbb{R}$, the Young inequalities hold~:
\begin{equation*}
	\begin{aligned}
		zy & \leq h(z)+h^*(y) = \left[(1+z)\log(1+z)-z\right] + \left[e^y-1-y\right], \\
		zy & \leq r(z)+r^*(y) = \left[z\log(1+z)\right] + r^*(y).
	\end{aligned}
\end{equation*}

It is also possible to show that $r(z)$ and $h(z)$ satisfy the reflection inequalities, for any $z>-1$,
\begin{equation*}
	\begin{aligned}
		h\left(|z|\right) & \leq h(z),\\
		r\left(|z|\right) & \leq r(z),
	\end{aligned}
\end{equation*}
and that $h^*(y)$ and $r^*(y)$ have the following exponential growth control, for any $y\geq 0$,
\begin{equation*}
	\begin{aligned}
		h^*\left( y\right) & \leq e^y, \\
		r^*\left( y\right) & \leq e^y,
	\end{aligned}
\end{equation*}
and the superquadratic homogeneity, for any $y\geq 0$ and $0\leq\lambda\leq 1$,
\begin{equation*}
	\begin{aligned}
		h^*\left(\lambda y\right) & \leq \lambda^2 h^*(y), \\
		r^*\left(\lambda y\right) & \leq \lambda^2 r^*(y),
	\end{aligned}
\end{equation*}
which is easily obtained by proving that $\lambda\mapsto\frac{h^*\left(\lambda y\right)}{\lambda^2}$ and $\lambda\mapsto\frac{r^*\left(\lambda y\right)}{\lambda^2}$ are increasing functions.

Thus, combining the above properties we arrive at the most useful inequalities
\begin{equation}
	\begin{gathered}
		\label{Young h}
		\left|zy\beta\right|
		\leq
		\frac{\alpha}{\eps^2}\left(h\left(\eps |z|\right) + h^*\left(\frac{\eps|\beta y|}{\alpha}\right)\right)
		\leq
		\frac{\alpha}{\eps^2}h(\eps z) + \frac{\beta^2}{\alpha} e^{|y|},\\
		\text{for any}\quad z>-\frac{1}{\eps},\quad y\in\mathbb{R},\quad \alpha\geq\eps|\beta|>0,
	\end{gathered}
\end{equation}
and
\begin{equation}
	\begin{gathered}
		\label{Young r}
		\left|zy\beta\right|
		\leq
		\frac{\alpha}{\eps^4}\left(
		r\left(\eps^2 |z|\right)
		+
		r^*\left(\frac{\eps^2|\beta y|}{\alpha}\right)\right)
		\leq
		\frac{\alpha}{\eps^4}r\left(\eps^2 z\right)
		+
		\frac{\beta^2}{\alpha} e^{|y|},\\
		\text{for any}\quad z>-\frac{1}{\eps^2},\quad y\in\mathbb{R},\quad \alpha\geq\eps^2|\beta|>0.
	\end{gathered}
\end{equation}

The above Young inequalities \eqref{Young h} and \eqref{Young r} are intensively used throughout this work to extract bounds and compactness properties from the various entropy and entropy dissipation bounds.

Finally, for reference, we list some elementary inequalities in connection with the convex functions $h(z)$ and $r(z)$ and their corresponding Young inequalities~:
\begin{itemize}
	
	\item For every $z>-1$, it holds that
	\begin{equation}\label{sqrt entropy}
		\left(\sqrt{1+z}-1\right)^2\leq h(z),
	\end{equation}
	which follows straightforwardly from the identity, for all $z>-1$,
	\begin{equation*}
		\left(\sqrt{1+z}-1\right)^2
		+2h\left(\sqrt{1+z}-1\right)\sqrt{1+z}
		=h(z),
	\end{equation*}
	or by showing that the function defined on $z> -1$ by
	\begin{equation*}
		f(z)=
		\begin{cases}
			\frac{h(z)}{\left(\sqrt{1+z}-1\right)^2} & \text{if }z\neq 0,
			\\
			2 & \text{if }z=0,
		\end{cases}
	\end{equation*}
	is increasing and reaches the value $1$ as $z\rightarrow -1$.
	
	\item For every $z>-1$, it holds that
	\begin{equation}\label{second order h expansion ineq}
		\frac 12 z^2 \leq h(z) + \frac 16 z^3,
	\end{equation}
	which is obtained by integrating twice the elementary inequality
	\begin{equation*}
		1\leq \frac 1{1+z} + z = h''(z)+z,
	\end{equation*}
	valid for every $z>-1$.
	
	\item For every $z>-1$, it holds that
	\begin{equation}\label{sqrt abs}
		\left(\sqrt{1+z}-1\right)^2\leq |z|,
	\end{equation}
	which is a consequence of the direct computation
	\begin{equation*}
		|z|-\left(\sqrt{1+z}-1\right)^2 =
		\begin{cases}
			2\left(\sqrt{1+z}-1\right)
			&\text{if } z\geq 0,
			\\
			2\sqrt{1+z}\left(1-\sqrt{1+z}\right)
			&\text{if } z\leq 0.
		\end{cases}
	\end{equation*}
	
	\item For every $z>-1$, it holds that
	\begin{equation}\label{sqrt dissip}
		\left(\sqrt{1+z}-1\right)^2\leq\frac{1}{4}r(z),
	\end{equation}
	which follows straightforwardly by integrating twice the inequality
	\begin{equation*}
		\frac{1}{(1+z)^\frac{3}{2}}\leq \frac 12\left(\frac 1{1+z}+\frac{1}{(1+z)^2}\right).
	\end{equation*}
	
	\item For every $z>-1$ and $y\in\mathbb{R}$, it holds that
	\begin{equation}\label{weaker young r}
		\begin{aligned}
			zy & = zy \mathds{1}_{\left\{y\geq 0, y\leq \log(1+z)\right\}\cup\left\{y<0,y>\log(1+z)\right\}}\\
			& + zy \mathds{1}_{\left\{y\geq 0, z< e^y - 1\right\}\cup \left\{y<0, z\geq e^y - 1\right\}} \\
			& \leq z\log(1+z) + (e^y - 1)y = r(z) + (e^y - 1)y.
		\end{aligned}
	\end{equation}
	This implies, in particular, that $r^*(y)\leq (e^y - 1)y$, for every $y\in\mathbb{R}$.
	
\end{itemize}

%% file: appendixhypoellipticity.tex
\chapter[End of proof of Lemma \ref{HYPO-THM 3} on hypoelliptic transfer\ldots]{End of proof of Lemma \ref{HYPO-THM 3} on hypoelliptic transfer of compactness}\label{endofproof}

The justification of Lemma \ref{HYPO-THM 3} has not been fully completed in Chapter \ref{hypoellipticity} lest it become unclear and tedious. Instead, we complete it now with a full justification of Lemma \ref{HYPO-THM 3}, when the equi-integrability of $\left\{\phi_\eps(t,x,v)\right\}_{\eps>0}$ is not known a priori and when each $\phi_\eps$ is non-negative, based on the proof of Theorem 2.4 from \cite{arsenio}. Recall that, according to the partial proof following the statement of Lemma \ref{HYPO-THM 3}, it is sufficient in this case to establish the equi-integrability of $\left\{\phi_\eps(t,x,v)\right\}_{\eps>0}$ in all variables.

For convenience of the reader, we first recall the precise result which we are about to justify.

\begin{lem}
	Let the bounded family of non-negative functions
	\begin{equation*}
		\left\{\phi_\eps(t,x,v)\right\}_{\eps>0}
		\subset L^1\left(\mathbb{R}_t\times\mathbb{R}^3_x; L^r\left(\mathbb{R}^3_v\right)\right),
	\end{equation*}
	for some $1<r<\infty$, be locally relatively compact in $v$ and such that
	\begin{equation*}
		\left(\eps\d_t + v\cdot\nabla_x\right) \phi_\eps =
		\left(1-\Delta_{x}\right)^\frac{\beta}{2}\left(1-\Delta_v\right)^\frac{\alpha}{2}S_\eps,
	\end{equation*}
	for all $\eps>0$ and for some bounded family
	\begin{equation*}
		\left\{S_\eps(t,x,v)\right\}_{\eps>0}
		\subset L^{1}\left(\mathbb{R}_t\times\mathbb{R}^3_x; L^r(\mathbb{R}^3_v)\right),
	\end{equation*}
	where $\alpha\geq 0$ and $0\leq \beta<1$. We further assume that, for any compact set $K\subset\mathbb{R}^3\times\mathbb{R}^3$,
	\begin{equation*}
		\left\{\int_{K}\phi_\eps(t,x,v) dxdv\right\}_{\eps>0}\quad \text{is equi-integrable (in $t$).}
	\end{equation*}
	
	Then, $\left\{\phi_\eps(t,x,v)\right\}_{\eps>0}$ is equi-integrable (in all variables).
\end{lem}

We advise the reader of the difficulty to grasp the full content of the proof below without any prior knowledge of the work from \cite{arsenio}. Therefore, we suggest that this appendix be read in parallel with the article \cite{arsenio}.

\begin{proof}
	We first notice, repeating the proof of Lemma 3.1 from \cite{arsenio}, that we have the following interpolation formula (compare with (5.17) in \cite{arsenio}), for any $R>0$ and $\delta >0$, denoting by $(\tau,\eta,\xi)$ the Fourier variables of $(t,x,v)$ and by $\mathcal{F}$ the Fourier transforms~:
	\begin{equation}\label{interpolation non stationary}
		\begin{aligned}
			& \left[ \mathcal{F}_{x,v}^{-1} \left(1-\chi\right)\left(\frac{\eta}{R}\right)
			\chi\left(\frac{\xi}{\lambda}\right)\mathcal{F}_{x,v} \right] \phi_\eps (t,x,v) \\
			&\hspace{20mm} = \frac{1}{\left(2\pi\right)^3} \int_{\mathbb{R}\times\mathbb{R}^3\times\mathbb{R}^3} \lambda^3 \hat \chi \left(\lambda(u-v)\right)
			K_1(t-s,x-y,u)\\
			&\hspace{20mm} \left[ \mathcal{F}_{x,v}^{-1}
			\left(1-\chi\right)\left(\frac{\eta}{R}\right)\psi \left( {\xi\over \lambda}\right) \mathcal{F}_{x,v}\right] \phi_\eps (s,y,u)dsdydu \\
			&\hspace{20mm} - \frac{2\delta i}{\left(2\pi\right)^3}
			\int_{\mathbb{R}\times\mathbb{R}^3\times\mathbb{R}^3} \lambda^3\hat \chi \left(\lambda(u-v)\right)
			K_2(t-s,x-y,u)\\
			&\hspace{20mm} \left[ \mathcal{F}_{x,v}^{-1}
			\left(1-\chi\right)\left(\frac{\eta}{R}\right)R^{1-\beta}{\ip{\eta}^{\beta}\over |\eta|}
			\chi\left( {\xi\over 10 \lambda}\right) {\ip{\xi}^\alpha\over \lambda^\alpha} \mathcal{F}_{x,v} \right] S_\eps (s,y,u)dsdydu,
		\end{aligned}
	\end{equation}
	% \begin{equation*}
	% 	\begin{aligned}
	% 		& \left[ \mathcal{F}^{-1} \left(1-\chi\right)\left(\frac{\eta}{R}\right)
	% 		\chi\left(\frac{\xi}{K}\right)\mathcal{F} \right] \phi_\eps (t,x,v) \\
	% 		&= \frac{1}{\left(2\pi\right)^3} \int_{\mathbb{R}^3}K^3 \hat \chi \left(K(u-v)\right)\\
	% 		& \left[ \mathcal{F}^{-1} \rho\left(\frac{2K}{|\eta|}\left(\eta \cdot u+\eps\tau\right)\right)
	% 		\left(1-\chi\right)\left(\frac{\eta}{R}\right)
	% 		\psi \left( {\xi\over K}\right) \mathcal{F}\right] \phi_\eps (t,x,u)du \\
	% 		& - \frac{2 i K^{1+\alpha}}{\left(2\pi\right)^3 R^{1-\beta'}}
	% 		\int_{\mathbb{R}^3} K^3\hat \chi \left(K(u-v)\right)\\
	% 		& \left[ \mathcal{F}^{-1} \gamma\left(\frac{2K}{|\eta|}\left(\eta \cdot u+\eps\tau\right)\right)
	% 		\left(1-\chi\right)\left(\frac{\eta}{R}\right)R^{1-\beta'}{\ip{\eta}^{\beta'}\ip{\eta}^{\beta-\beta'}\over |\eta|}
	% 		\chi\left( {\xi\over 10 K}\right) {\ip{\xi}^\alpha\over K^\alpha} \mathcal{F} \right] S_\eps (t,x,u)du,
	% 	\end{aligned}
	% \end{equation*}
	where $\lambda=\delta^{\frac{1}{1+\alpha}}R^{\frac{1-\beta}{1+\alpha}}$, $\chi,\psi\in C_c^\infty\left(\mathbb{R}^3\right)$ are cutoff functions such that $\mathds{1}_{\left\{|r|\leq \frac 12\right\}} \leq \chi(r) \leq \mathds{1}_{\left\{ |r| \leq 1 \right\}}$ and $\mathds{1}_{\left\{1\leq|r|\leq 5\right\}}\leq \psi(r) \leq\mathds{1}_{\left\{\frac 12\leq|r|\leq \frac{11}{2}\right\}}$, and the singular kernels $K_i$, $i=1,2$, are defined by
	\begin{equation*}
		\begin{aligned}
			K_1(t,x,u)&=\mathcal{F}_{t,x}^{-1} \left[\rho\left(\frac{2\lambda}{|\eta|}\left(\eta \cdot u+\eps\tau\right)\right)
			\right],\\
			K_2(t,x,u)&=\mathcal{F}_{t,x}^{-1} \left[\gamma\left(\frac{2\lambda}{|\eta|}\left(\eta \cdot u+\eps\tau\right)\right)
			\right],
		\end{aligned}
	\end{equation*}
	where $\rho\in\mathcal{S}\left(\mathbb{R}\right)$ is such that $\rho(0)=\frac{1}{2\pi}\int_{\mathbb{R}}\hat\rho(\sigma)d\sigma=1$ and $\operatorname{supp}\hat\rho\subset\left\{1\leq |\sigma|\leq 2\right\}$, and $\gamma(r)=\frac{1-\rho(r)}{r}$ for all $r\in\mathbb{R}$.

	The rather deep meaning of the above interpolation formula \eqref{interpolation non stationary} resides in the presence of the frequency cutoff functions $\psi\left(\frac \xi\lambda\right)$ and $\chi\left(\frac\xi{10\lambda}\right)$ in its right-hand side, which is only possible through a precise analysis of dispersive and hypoelliptic phenomena, i.e.\ the transport of frequencies, in the kinetic transport equation. This requires that the support of $\hat\rho(\sigma)$ be restricted to $\left\{1\leq |\sigma|\leq 2\right\}$, though.

	Now, we claim that the mappings
	\begin{equation}\label{mappings}
		f(t,x,v)\mapsto\frac{1}{\left(2\pi\right)^3} \int_{\mathbb{R}\times\mathbb{R}^3\times\mathbb{R}^3} \lambda^3 \hat \chi \left(\lambda(u-v)\right)
		K_i(t-s,x-y,u) f (s,y,u)dsdydu,
	\end{equation}
	with $i=1,2$, have bounded extensions between
	\begin{equation}\label{between}
		L^1\left(\mathbb{R}_t\times\mathbb{R}^3_x;L^r\left(\mathbb{R}^3_v\right)\right)
		\longrightarrow
		L^r\left(\mathbb{R}^3_v;L^{1,\infty}\left(\mathbb{R}_t\times\mathbb{R}^3_x\right)\right).
	\end{equation}
	As shown in \cite{arsenio}, this boundedness follows from an application of Theorem 5.2 therein (or a very slight variant of it allowing different dimensions for different variables, i.e.\ $(t,x)\in\mathbb{R}^4$ and $v\in\mathbb{R}^3$) provided we establish that
	\begin{equation}\label{CZ operator}
		\begin{aligned}
			\sup_{v\in\mathbb{R}^3}
			\Bigg\| & \sup_{(s,y)\neq 0}\int_{\left\{|(t,x)|\geq 2|(s,y)|\right\}}
			\\
			& \left|A_v\sharp K_i\left(t-s,x-y,u\right)-A_v\sharp K_i \left(t,x,u\right)\right| dtdx \Bigg\|_{L^{r'}\left(\left|\lambda^3\hat\chi\left(\lambda(u-v)\right)\right|du\right)}
			<\infty,
		\end{aligned}
	\end{equation}
	where $A_v\sharp K_i(t,x,u)=\left|\det A_v\right|K_i\left(A_v(t,x),u\right)$, for some family of automorphisms $A_v$ of $\mathbb{R}\times\mathbb{R}^3$ fixing the origin (the remaining hypotheses of Theorem 5.2 from \cite{arsenio} being easily verified through standard arguments from the analysis of Fourier multipliers). Here, we will consider, for each $v\in\mathbb{R}^3$, the automorphism $A_v$ of $\mathbb{R}\times\mathbb{R}^3$ defined by
	\begin{equation*}
		A_v(t,x)=\left(\lambda\eps t, x+\lambda tv\right),
	\end{equation*}
	so that
	\begin{equation*}
		\begin{aligned}
			A_v\sharp K_1(t,x,u)&=\mathcal{F}_{t,x}^{-1} \left[\rho\left(2\lambda\frac{\eta\cdot(u-v)}{|\eta|}+2\frac{\tau}{|\eta|}\right)
			\right],\\
			A_v\sharp K_2(t,x,u)&=\mathcal{F}_{t,x}^{-1} \left[\gamma\left(2\lambda\frac{\eta\cdot(u-v)}{|\eta|}+2\frac{\tau}{|\eta|}\right)
			\right].
		\end{aligned}
	\end{equation*}

	In turn, the estimate \eqref{CZ operator} is established employing classical methods from harmonic analysis (in \cite{arsenio}, the corresponding step is performed in Lemmas 5.3 and 5.4). Thus, it can be shown that \eqref{CZ operator} will hold as soon as the following H\"ormander-Mikhlin condition for homogeneous Fourier multipliers is verified~:
	\begin{equation}\label{HM}
		\int_{1 < \left|(\tau,\eta)\right| < 2 }
		\left| \partial_{\tau,\eta}^\nu\mathcal{F}_{t,x}\left(A_v\sharp K_i\right)(\tau,\eta,u) \right|^2d\tau d\eta \leq C\left(1+\lambda|u-v|\right)^N,
	\end{equation}
	for some independent constant $C>0$ and some possibly very large $N\in\mathbb{N}$, and for any multi-index $\nu\in\mathbb{N}^4$ such that $\nu_1\in\left\{0,1\right\}$ and $0\leq \nu_2+\nu_3+\nu_4 \leq 2$.

	The preceding control is easily verified for $K_1$ through a straightforward calculation using that $\rho(r)$ decays rapidly. As for $K_2$, this step requires some greater care because $\gamma^{(n)}(r)$ only decays as $\frac 1{|r|^{n+1}}$ for large values of $r$, for all $n\in\mathbb{N}$. Nevertheless, it is also straightforwardly verified upon noticing, for any multi-index $\nu$ as before, that
	\begin{equation*}
		\begin{aligned}
			\left|(\tau,\eta)\right|^{|\nu|}
			&
			\left|\partial_{\tau,\eta}^\nu
			\left[\gamma\left(2\lambda\frac{\eta\cdot(u-v)}{|\eta|}+2\frac{\tau}{|\eta|}\right) \right]\right|
			\\
			& \leq
			C \sum_{n=0}^3
			\left(1+\lambda|u-v|\right)^n
			\left(1+\frac{|\tau|}{|\eta|}\right)^{n+2}
			\left|\gamma^{(n)}\left(2\lambda\frac{\eta\cdot(u-v)}{|\eta|}+2\frac{\tau}{|\eta|}\right)\right|
			\\
			& \leq
			C \left(1+\frac{|\tau|}{|\eta|}\right) \sum_{n=0}^3
			\left(1+\lambda|u-v|\right)^{2n+1}
			\\
			& \left|
			\left(1+\left|\lambda\frac{\eta\cdot(u-v)}{|\eta|}+\frac{\tau}{|\eta|}\right|\right)^{n+1}
			\gamma^{(n)}\left(2\lambda\frac{\eta\cdot(u-v)}{|\eta|}+2\frac{\tau}{|\eta|}\right)\right|
			\\
			& \leq
			C \left(1+\frac{|\tau|}{|\eta|}\right)
			\left(1+\lambda|u-v|\right)^{7}
			\leq
			\frac{C}{|\eta|}
			\left(1+\lambda|u-v|\right)^{7},
		\end{aligned}
	\end{equation*}
	which implies that \eqref{HM} holds with $N=14$, for $\frac{1}{|\eta|^2}$ is locally integrable in $\mathbb{R}^3$.

	So far, we have thus established the weak type boundedness on the spaces \eqref{between} of the mappings \eqref{mappings}. Therefore, we conclude from the interpolation formula \eqref{interpolation non stationary} that
	\begin{equation*}
		\begin{aligned}
			\left\|
			\left[ \mathcal{F}_{x}^{-1} \left(1-\chi\right)\left(\frac{\eta}{R}\right)
			\mathcal{F}_{x} \right] \phi_\eps
			\right\|_{L^r_vL^{1,\infty}_{t,x}}
			& \leq C \left\|
			\left[ \mathcal{F}_{v}^{-1}
			\left(1-\chi\right)\left(\frac{\xi}{\lambda}\right)\mathcal{F}_{v} \right] \phi_\eps
			\right\|_{L^{1}_{t,x}L^r_v}
			\\
			& + C \left\|
			\left[ \mathcal{F}_{v}^{-1} \psi\left(\frac{\xi}{\lambda}\right)\mathcal{F}_{v} \right] \phi_\eps
			\right\|_{L^{1}_{t,x}L^r_v}
			\\
			& + C \delta \left\|
			S_\eps
			\right\|_{L^1_{t,x}L^r_v},
		\end{aligned}
	\end{equation*}
	whence, since the $\phi_\eps$'s are relatively compact in $v$, since $\lim_{R\rightarrow\infty}\lambda=\infty$ and by the arbitrariness of $\delta>0$,
	\begin{equation}\label{compactness1}
		\lim_{R\rightarrow\infty}\sup_{\eps>0}\left\|
		\left[ \mathcal{F}_{x}^{-1} \left(1-\chi\right)\left(\frac{\eta}{R}\right)
		\mathcal{F}_{x} \right] \phi_\eps
		\right\|_{L^r_vL^{1,\infty}_{t,x}} = 0.
	\end{equation}

	Next, applying spatial Riesz transforms $f\mapsto \mathcal{F}^{-1}_x\frac{\eta_i}{|\eta|}\mathcal{F}_x f$, $i=1,2,3$, to the identity \eqref{interpolation non stationary} and repeating the preceding arguments, we deduce that it also holds
	\begin{equation}\label{compactness2}
		\lim_{R\rightarrow\infty}\sup_{\eps>0}\left\|
		\left[ \mathcal{F}_{x}^{-1} \left(1-\chi\right)\left(\frac{\eta}{R}\right)\frac{\eta_i}{|\eta|}
		\mathcal{F}_{x} \right] \phi_\eps
		\right\|_{L^r_vL^{1,\infty}_{t,x}} = 0.
	\end{equation}
	Consequently, we obtain that the $\phi_\eps$'s are relatively compact in $x$ in the quasi-Banach space defined by the quasi-norm
	\begin{equation*}
		\left\|f\right\|_{L^r_vL^{1,\infty}_{t,x}}+\sum_{i=1}^3 \left\|\mathcal{F}^{-1}_x\frac{\eta_i}{|\eta|}\mathcal{F}_x f\right\|_{L^r_vL^{1,\infty}_{t,x}},
	\end{equation*}
	which can be shown, repeating standard (but difficult) arguments from the analysis of weak Hardy spaces (see \cite{arsenio} and the references therein for details), to be equivalent to the quasi-norm
	\begin{equation}\label{quasinorm}
		\left\|\sup_{s>0}\left|\varphi_s *_x f\right|\right\|_{L^r_vL^{1,\infty}_{t,x}},
	\end{equation}
	where $\varphi(x)\in C_c^\infty\left(\mathbb{R}^3\right)$ is such that $\int_{\mathbb{R}^3}\varphi(x)dx\neq 0$ and $\varphi_s(x)=\frac{1}{s^3}\varphi\left(\frac{x}{s}\right)$.

	Finally, following \cite{arsenio}, we explain how the ensuing spatial compactness in the topology given by the quasi-norm \eqref{quasinorm} is sufficient to entail the equi-integrability of the $\phi_\eps$'s, provided they are non-negative. To this end, we note first, for every $a,b>0$ and any compact subset $K\subset \mathbb{R}\times\mathbb{R}^3\times\mathbb{R}^3$, that
	\begin{equation*}
		\begin{aligned}
			& \int_{K\cap\left\{\left|\phi_\eps\right|>a \right\}} \left|\phi_\eps (t,x,v)\right| dtdxdv\\
			& \leq
			\int_{K\cap\left\{\left|\phi_\eps\right|>a,\ \left\|f_\lambda\right\|_{L^1_x}\leq b \right\}} \left|\phi_\eps (t,x,v)\right| dtdxdv
			+
			\int_{\left\{\left\|\phi_\eps\right\|_{L^1_x}> b \right\}} \left\|\mathds{1}_{K}\phi_\eps (t,x,v)\right\|_{L^1_x}dtdv.
		\end{aligned}
	\end{equation*}
	Clearly, in view of the equi-integrability hypotheses, the family $\left\|\phi_\eps (t,x,v)\right\|_{L^1_x}$ is locally equi-integrable in $t$ and $v$, so that the last term above can be made uniformly small by choosing $b>0$ arbitrarily large.

	Thus, the proof will be finished upon showing that, for each fixed $b>0$ and every compact subset $K\subset \mathbb{R}\times\mathbb{R}^3\times\mathbb{R}^3$,
	\begin{equation*}
		\lim_{a\rightarrow\infty}\sup_{\eps>0}
		\int_{\set{(t,x,v)\in K}{\left|\phi_\eps\right|>a,\ \left\|\phi_\eps\right\|_{L^1_x}\leq b}} \left|\phi_\eps (t,x,v)\right| dtdxdv=0,
	\end{equation*}
	which, following the arguments from the end of the proof of Theorem 2.4 in \cite{arsenio} based on Proposition 5.5 therein (which is only valid for non-negative functions), is a direct consequence of the compactness estimates \eqref{compactness1} and \eqref{compactness2}.

	The justification of the lemma is thus complete.
\end{proof}